\documentclass[twoside,reqno,openright,a4paper,11pt]{report}
\usepackage[top=3cm, bottom=3cm, left=2.75cm, right=2.75cm]{geometry}
\usepackage{fancyhdr}
\usepackage{anyfontsize}
\usepackage{tikz}
\usetikzlibrary{shapes, arrows, arrows.meta, decorations.markings, babel}
\usepackage[european, nooldvoltagedirection]{circuitikz}
\usepackage{transparent}
\usepackage{array}
\usepackage{booktabs}
\usepackage{longtable}
\usepackage{multirow}
\usepackage{amsmath}
\usepackage{amsthm}
\usepackage{mathtools}
\usepackage[integrals]{wasysym}
\usepackage{stmaryrd}
\usepackage{mathdots}

\usepackage{bbm}
\usepackage[mathcal]{euscript}

\usepackage[utf8]{inputenc}
\usepackage[T1]{fontenc}

\usepackage[brazil, french, english]{babel}

\usepackage[noadjust]{cite}
\usepackage{mathscinet}

\usepackage{kpfonts}

\usepackage[chapter]{algorithm}
\usepackage{algorithmic}
\usepackage{animate}
\usepackage{caption}
\usepackage{url}
\usepackage{epigraph}
\usepackage[inline]{enumitem}
\usepackage{pdfpages}
\usepackage{siunitx}
\usepackage[normalem]{ulem}
\usepackage{appendix}
\usepackage{etoolbox}
\usepackage{lipsum}
\usepackage{import}
\usepackage{ifthen}

\newtoggle{print}
\togglefalse{print}

\iftoggle{print}{%
	\usepackage[hidelinks, bookmarks, bookmarksnumbered, pdfstartview={XYZ null null 1.00}]{hyperref}%
}{%
  \usepackage[colorlinks, linkcolor=red!50!black, citecolor=blue!75!black, urlcolor=violet, bookmarks, bookmarksnumbered, pdfstartview={XYZ null null 1.00}]{hyperref}%
}
\usepackage{bookmark}

\newlength{\medianiz}
\iftoggle{print}{%
  \setlength{\medianiz}{0.25cm}%
}{%
  \setlength{\medianiz}{0cm}%
}
\theoremstyle{plain}
\newtheorem{theorem}{Theorem}[chapter]
\newtheorem{lemma}[theorem]{Lemma}
\newtheorem{corollary}[theorem]{Corollary}
\newtheorem{proposition}[theorem]{Proposition}
\newtheorem{conjecture}[theorem]{Conjecture}

\newtheorem*{claim*}{Claim}
\theoremstyle{definition}
\newtheorem{definition}[theorem]{Definition}
\newtheorem{hypotheses}[theorem]{Hypotheses}
\newtheorem{remark}[theorem]{Remark}
\newtheorem{example}[theorem]{Example}

\makeatletter
\def\step#1{\def\tempa{#1}\futurelet\next\step@i} 
\def\step@i{\ifx\next\bgroup\expandafter\step@ii\else\expandafter\step@end\fi} 
\def\step@ii#1{ 
  \par
  \addvspace{\medskipamount}%
  \noindent\emph{Step~\tempa: #1.}\par\nobreak\smallskip
  \@afterheading%
} 
\def\step@end{ 
  \par
  \addvspace{\medskipamount}%
  \noindent\emph{Step~\tempa.}%
}
\makeatother

\makeatletter
\def\case#1{\def\tempa{#1}\futurelet\next\case@i} 
\def\case@i{\ifx\next\bgroup\expandafter\case@ii\else\expandafter\case@end\fi} 
\def\case@ii#1{ 
  \par
  \addvspace{\medskipamount}%
  \noindent\emph{Case~\tempa: #1.}\par\nobreak\smallskip
  \@afterheading%
} 
\def\case@end{ 
  \par
  \addvspace{\medskipamount}%
  \noindent\emph{Case~\tempa.}%
}
\makeatother

\DeclareMathOperator{\sign}{sign}

\DeclareMathOperator{\Ker}{Ker}
\DeclareMathOperator{\Real}{Re}
\DeclareMathOperator{\Imag}{Im}
\DeclareMathOperator{\Tr}{Tr}
\DeclareMathOperator{\spr}{spr}

\DeclareMathOperator{\rank}{rk}
\DeclareMathOperator{\Ran}{Ran}
\DeclareMathOperator{\Adj}{Adj}

\DeclareMathOperator{\id}{id}

\DeclareMathOperator{\Span}{Span}

\DeclareMathOperator{\supp}{supp}

\DeclareMathOperator{\V}{V}

\DeclarePairedDelimiter{\norm}{\lVert}{\rVert}
\DeclarePairedDelimiter{\abs}{\lvert}{\rvert}
\DeclarePairedDelimiter{\floor}{\lfloor}{\rfloor}
\DeclarePairedDelimiter{\ceil}{\lceil}{\rceil}

\DeclarePairedDelimiter{\scalprod}{\langle}{\rangle}

\newcommand{\1}{\mathbbm{1}}
\newcommand{\suchthat}{\ifnum\currentgrouptype=16 \mathrel{}\middle\vert\mathrel{}\else\mid\fi}
\newcommand{\diff}{\,\mathrm{d}}
\newcommand{\df}{\mathrm{d}}
\newcommand{\transp}{{\mathrm{T}}}
\newcommand{\loc}{\mathrm{loc}}
\newcommand{\Id}{\mathrm{Id}}

\DeclareMathOperator{\NBV}{NBV}

\DeclareMathOperator{\WellPosed}{\mathsf{WP}}
\DeclareMathOperator{\NoSingular}{\mathsf{NS}}
\DeclareMathOperator{\LowerBound}{\mathsf{LB}}
\DeclareMathOperator{\SDG}{\mathsf{SDG}}

\fancypagestyle{plain}{%
\fancyhf{}%
\fancyhf[ROF, LEF]{\it \thepage}%
\renewcommand{\headrulewidth}{0pt}%
\renewcommand{\footrulewidth}{0pt}}

\makeatletter
\def\cleardoublepage{\clearpage\if@twoside \ifodd\c@page\else
\hbox{}
\thispagestyle{plain}
\newpage
\if@twocolumn\hbox{}\newpage\fi\fi\fi}
\makeatother

\makeatletter
\newcommand{\customlabel}[2]{%
   \protected@write \@auxout {}{\string \newlabel {#1}{{#2}{\thepage}{#2}{#1}{}} }%
   \hypertarget{#1}{#2}
}
\makeatother

\begin{document}

\setlength{\parskip}{1pt plus 1pt minus 1pt}

\setlist[enumerate, 1]{label={\textnormal{(\alph*)}}, ref={(\alph*)}, leftmargin=0pt, itemindent=*}
\setlist[enumerate, 2]{label={\textnormal{(\roman*)}}, ref={(\roman*)}, itemindent=\parindent}
\setlist[description, 1]{leftmargin=0pt, itemindent=*}
\setlist[itemize, 1]{label={\textbullet}, leftmargin=0pt, itemindent=*}

\newlist{hypothesis}{enumerate}{1}
\setlist[hypothesis]{label={\textup{(H\arabic*)}}, ref={(H\arabic*)}, leftmargin=*, widest*=10}
\newlist{hypolist}{enumerate}{1}
\setlist[hypolist]{label={\textup{(H\arabic*)}}, ref={(H\arabic*)}, leftmargin=*, widest*=10, resume=hypolist}

\makeatletter
\renewcommand{\@pnumwidth}{2em}
\renewcommand{\@tocrmarg}{3em}
\renewcommand{\@dotsep}{4}
\makeatother

\newcommand{\ThesisTitle}{Control, stability, and stabilization of some classes of linear time-delay systems}

\hypersetup{pdftitle={\ThesisTitle}, pdfauthor={Guilherme Mazanti}, pdfkeywords={Habilitation à diriger des recherches}}%

\pagenumbering{Alph}

%
%

\includepdf[pages=1]{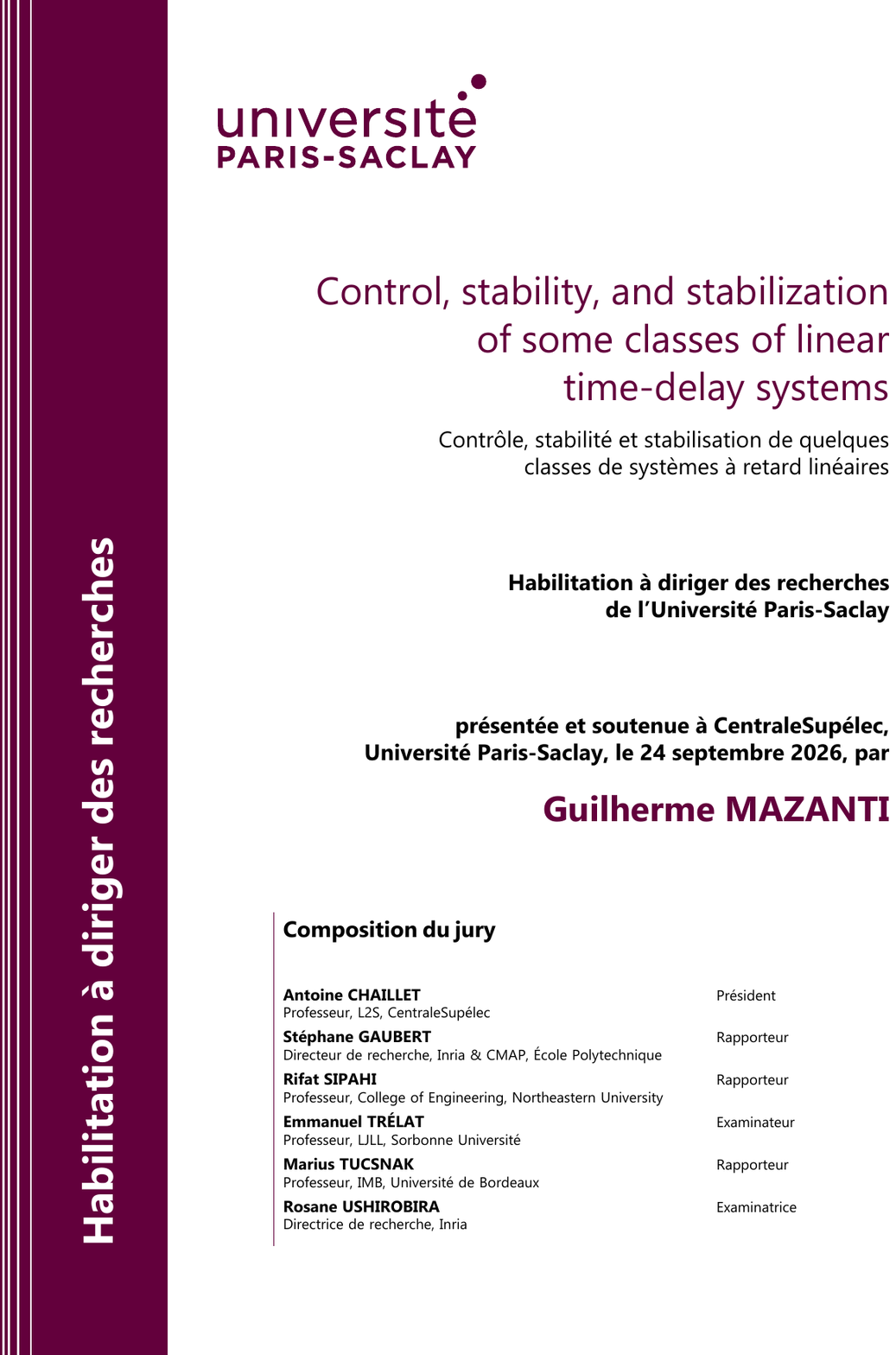}

\thispagestyle{empty}
\mbox{}
\newpage

\pagestyle{fancy}
\fancyhf{}
\fancyhf[ROF, LEF]{\it \thepage}
\renewcommand{\headrulewidth}{0pt}
\renewcommand{\footrulewidth}{0pt}
\pagenumbering{roman}

\chapter*{Agradecimentos, remerciements \& acknowledgments}

I would like to start by thanking Antoine Chaillet, Stéphane Gaubert, Rifat Sipahi, Emmanuel Trélat, Marius Tucsnak, and Rosane Ushirobira for agreeing to take part in the evaluation of this work, with special thanks to Stéphane Gaubert, Rifat Sipahi, and Marius Tucsnak for their careful reading of the (rather long) manuscript in order to write their reports.

The work presented in this manuscript is the product of many collaborations over the years, and I would thus like to take this opportunity to thank all my collaborators, 
Gonzalo Arias,
Saeed Sa\-de\-ghi Arj\-mand,
Jean Auriol,
Is\-maï\-la Ba\-lo\-goun,
Ami\-na Benarab,
Cathe\-rine Bon\-net,
Is\-lam Boussaada,
Adam Braun,
Hugo Ca\-va\-lera,
Ya\-cine Chi\-tour,
Anh-Tuan Cla\-baut,
Fritz Co\-lo\-nius,
Ro\-main Du\-casse,
Samer Dweik,
Sé\-bas\-tien Fueyo,
Ayr\-ton Ham\-mou\-mou,
Ju\-lien Huynh,
Adri\-en Le\-clerc,
Xu-Guang Li,
Pierre Lissy,
João Mi\-guel Ma\-cha\-do,
Swann Marx,
Fré\-dé\-ric Ma\-zenc,
Cé\-sar Fer\-nan\-do Mén\-dez Bar\-rios,
Ja\-que\-li\-ne Go\-doy Mes\-qui\-ta,
Wim Mi\-chiels,
Ti\-touan Mil\-let,
Pierre Mon\-mar\-ché,
Thi\-bault Mo\-quet,
Fe\-li\-pe Gon\-çal\-ves Net\-to,
Sil\-viu Ni\-cu\-les\-cu,
Max Per\-rau\-din,
Lau\-rent Pfeiffer,
Jay\-vir Raj,
Fi\-lippo San\-tam\-bro\-gio,
Ma\-rio Si\-ga\-lotti,
Franck Sim,
Ying Tang,
Mat\-thieu Thomas,
Karim Tra\-bel\-si,
Ni\-co\-las Van\-spranghe, and
To\-máš Vyhlídal.
Be they established, longstanding collaborations or occasional or more recent ones, I learn a lot from each collaboration, and they all contributed to make me the researcher I am today. More than colleagues, most of my long-term collaborators have become friends, and working with friends is always a great pleasure. And I would also like to thank all my collaborators for putting up with me, even when I suggest that we should change the notations in the whole paper when we are about to finish it, even when I use those horrible \(\Xi\) as variable names, even when I decide that \(\mathsf X\) looks better than \(X\) if \(X\) is a Banach space and I change everything in the paper, even when I insist that it is very important that we should decide once and for all if \(0 \in \mathbb N\) and \(0 \in \mathbb R_+\) or not in the context of the paper and make it clear in the notations, or, more importantly, even when I start my \LaTeX{} rants about hyphens, en dashes, and em dashes; about unbreakable spaces and manual spacing after periods when they do not mark the end of the sentence; about why \verb@\colon@ is better than~\verb@:@ for defining functions; and about why one should \emph{never} use \verb@eqnarray@ \cite{Madsen2006Avoid}.

\begin{otherlanguage}{french}
Je tiens à remercier l'équipe Inria DISCO de m'avoir accueilli en tant que post-doc en 2019 et de m'avoir fortement soutenu et encouragé lors de ma candidature au sein de l'équipe en 2020. Merci aussi, Frédéric B., Catherine, Islam, Ziad, Frédéric M., Silviu, Laurent et Giorgio, pour être à mes côtés en tant qu'équipe, pour toutes nos discussions, scientifiques ou pas, et pour tous vos conseils. Je tiens aussi à remercier nos successives assistantes d'équipe ces dernières années, Bamissa Sangare, Mélanie Da Silva et Irina Lahaye, ainsi qu'à tout le personnel du SAER qui nous a accompagnés pendant les périodes où DISCO n'avait pas d'assistant propre. Votre travail nous aide beaucoup dans nos activités de chercheurs et votre grande efficacité nous permet de nous concentrer sur ce qui est le plus important pour nous, la recherche. Je remercie également le personnel technique et administratif de l'Inria Saclay, notamment Sébastien Gilles et Christian Poli pour leurs contributions dans nos développements logiciels, Estelle Fiévé et Nathalie Pinto pour leur aide dans nos dépôts de logiciels et collaborations industrielles, et Delphine Charbonneau pour son accompagnement dans nos dossiers RH.

Notre équipe DISCO n'est pas isolée, mais fait partie du Laboratoire des Signaux et Systèmes (L2S), que je tiens aussi à remercier. Outre la qualité scientifique du laboratoire et de notre pôle Automatique et Systèmes, qui me permet d'avoir d'échanges scientifiques du plus haut niveau avec mes collègues, je m'estime chanceux d'être arrivé, dès mon recrutement en 2020, dans un laboratoire avec une très bonne ambiance et des collègues qui ont pu devenir des amis au fil des années. Merci à toutes et à tous qui y contribuent quotidiennement et qui font de ce laboratoire un endroit convivial et propice à faire avancer la science. Je tiens à remercier en particulier tous les membres de notre pôle Automatique et Systèmes, Nina, Jean, Riccardo, Frédéric B., Catherine, Benoît, Islam, Lucas, Antoine C., Yacine, Mladen, Sette, Didier, Stéphane, Antoine G., Luca, Alessio, Ziad, Françoise, Hugo, Antonio, Maria, Paolo, Frédéric M., Hugues, Florentina, Silviu, Dorothée, Sorin, Elena, William, Laurent (même si tu nous a abandonnés il y a quelques jours), Dario, Pedro, Guillaume, Houria, Cristina S., Sihem, Aarsh, Sami, Giorgio et Cristina V., ainsi qu'à tous les doctorants et post-docs qui font vivre notre laboratoire. Notre vie de chercheur est aussi rendue beaucoup plus simple grâce aux services techniques et administratifs du laboratoire, et je remercie ainsi Christian David et Stéphane Gérard pour leur soutien informatique, Audrey Bertinet, Stéphanie Douesnard et Catherine Kilic pour leur soutien de gestion et leur grande disponibilité pour répondre à nos questions, y compris les plus bêtes ou naïves, et à José Fonseca, Anne Kohl, et toute l'équipe de gestion financière du laboratoire, toujours à nos côtés pour nous aider à accueillir nos visiteurs et à partir en voyage professionnel.

À CentraleSupélec, nous avons aussi la chance d'avoir notre Fédération de Mathématiques (FdM), soutenue aussi par le CNRS, et qui regroupe les chercheurs de plusieurs laboratoires de CentraleSupélec dont la recherche se trouve à l'interface entre les thématiques de leur laboratoire et les mathématiques. C'est une structure parfaite pour quelqu'un comme moi, souvent considéré comme trop mathématicien par les automaticiens et trop automaticien par les mathématiciens. Merci donc à Pauline Lafitte de m'avoir accueilli au sein de la FdM, à Alexandre Richard pour l'excellent travail de direction de la FdM ces dernières années, aux organisateurs actuels et précédents des séminaires de la fédération et des ``explorations mathématiques'' pour contribuer à l'animation scientifique en mathématiques à CentraleSupélec, et à tous les collègues de la fédération, qui contribuent à faire vivre notre discipline au sein de CentraleSupélec dans une très bonne ambiance.

L'enseignement joue un rôle important dans ma carrière, et je tiens dans ce sens à remercier Guillaume Sandou, qui m'a attribué mes premiers cours lors du début de mon post-doc au L2S en 2019~; Charles Soussen, qui m'a aidé dans mes premiers pas d'enseignant à CentraleSupélec dans son cours d'Optimisation pour les élèves en dual diplôme~; Pierre-Yves Richard et Stéphane Vialle, qui successivement, pendant 6 années, m'ont confié l'enseignement de l'Optimisation aux élèves en alternance de CentraleSupélec~; John Cagnol, qui m'a confié un groupe de CM et de TD pour que j'aie l'honneur d'enseigner les équations aux dérivées partielles aux élèves de première année de CentraleSupélec, ainsi que toute l'équipe enseignante d'EDP, qui contribue à faire de l'enseignement de ce cours une tâche extremement agréable~; Erick Herbin, qui m'a soutenu dans ma démarche pour devenir professeur attaché au département de Mathématiques de CentraleSupélec~; et à nouveau Erick ainsi que Pauline Lafitte, pour m'avoir confié les étudiants de 3ème année de CentraleSupélec pour que je leur enseigne la théorie des jeux et le contrôle.
\end{otherlanguage}

\begin{otherlanguage}{brazil}
Se cheguei ao ponto da minha carreira acadêmica de defender a minha habilitação a dirigir pesquisas, é porque sempre tive ao meu lado uma família que me incentivava a estudar (apesar das piadas de falar que era para eu parar de estudar e sair para a rua) (ou ao menos eu interpretava como piada), com mãe, avós e tios que ficavam orgulhosos de mim a cada nota boa na escola, a cada aprovação em vestibular, e a cada sucesso acadêmico, mesmo eu tendo decidido morar do outro lado do oceano Atlântico. Não deve ser fácil para uma mãe ouvir o filho falar que foi aprovado para estudar na École Polytechnique, a 9.503~km de casa, e ter que dar os parabéns ao mesmo tempo em que pergunta se ele realmente vai se mudar para o outro lado do planeta. Acho que ela já desconfiava que a minha estadia na França ia durar um pouco mais do que os 2 anos e meio que eu ficaria na Polytechnique. Por isso tudo, e por muito mais, meu muito obrigado!

Além da família, tive a sorte de encontrar, em meu percurso escolar e acadêmico, diversos professores que me incentivaram. Seria impossível listar os nomes de todos aqueles que, de uma forma ou de outra, contribuíram com a minha formação, mas quero deixar agradecimentos especiais à Maria Luzia, que despertou em mim o gosto pelas ciências e me ensinou que estudar vai muito além de decorar, ao `Seu' Pedro (\emph{in memoriam}) e ao Estéberson, que, muito além de me ensinar matemática, me deram o gosto pela disciplina, e ao Hélio, que me incentivou a sonhar mais alto e buscar as melhores formações. A lista de professores da USP que me incentivaram nos meus estudos é também longa demais para ser detalhada aqui, mas não posso deixar de mencionar o professor Hildebrando Munhoz Rodrigues (\emph{in memoriam}) que, além de me ensinar EDOs e sistemas dinâmicos em cursos da graduação e do doutorado, foi meu orientador de iniciação científica por 2 anos e meio, um longo período no qual, muito além de análise, integração de Lebesgue e álgebra linear, eu aprendi o que é a pesquisa em matemática e o que é a carreira de pesquisador. O professor Hildebrando também me ensinou diversas outras coisas que tiveram um impacto fundamental na minha carreira, como a redação de textos matemáticos (e o uso de \LaTeX{}), a arte de se apresentar resultados matemáticos (inclusive em inglês!), e outras coisas bem mais básicas como o que é a École Polytechnique e como me inscrever no concurso internacional.

Os incentivos para chegar onde cheguei também vieram, direta ou indiretamente, dos meus amigos, aos quais agradeço profundamente. Obrigado Stênio, David, José, Ana, Yara e Raíssa por terem feito dos meus anos de escola uma época inesquecível e por todos os bons momentos que passamos juntos a cada reencontro; obrigado John, PHP, Fábio e Tenshi por tudo que passamos juntos na nossa época de graduação na USP, incluindo os nossos Nerd Xurras; obrigado Ricardo, Chico, George, Levi, Soneca, Thaenan, Hayssam, Ferraro, Natal e Mineiro por termos formado esse excelente grupo de brasileiros X2008; e obrigado também aos meus amigos brasileiros na França, minha segunda família, Ricardo (de novo!), Laís, Léo, Bianca, Bernardo, Rodrigo, Tomer, Dudu, Andrews e Carol. São tão segunda família que alguns falaram que vão até tirar folga para ir assistir à minha defesa de HDR, e nem é só pela comida!

Por fim, quero agradecer muito especificamente à direção do Colégio Santa Clara por ter organizado, em 11 de dezembro de 2005, um evento reunindo alunos do 3º colegial das escolas de Assis e Cândido Mota. Pois foi neste evento que conheci aquela que hoje é minha companheira para todos os momentos, Ariana, minha esposa há mais de 12 anos. Ela poderia ter escolhido alguém mais bonito, alguém mais rico, alguém com piadas mais engraçadas e com um senso de humor menos matemático, alguém que tivesse planos para trabalhar ``de verdade'' em vez de fazer um doutorado e seguir carreira acadêmica, alguém que quisesse se mudar para um país cuja língua ela já falava, mas por algum motivo ela me escolheu, e agradeço a ela todos os dias por isto. E também pela incrível vida familiar que construímos, fonte da minha felicidade cotidiana, e pela nossa maravilhosa filha Liliana, a quem dedico este trabalho.
\end{otherlanguage}
\cleardoublepage

\pdfbookmark[1]{Résumé}{abstract-fr}

\begin{otherlanguage}{french}

\begin{center}
\textbf{\Large Contrôle, stabilité et stabilisation de quelques classes de systèmes à retard linéaires}
\end{center}

\noindent Ce document présente une sélection des travaux récents de l'auteur sur le contrôle, la stabilité et la stabilisation de certaines classes de systèmes à retard linéaires. Après une introduction détaillée aux sujets sur les systèmes à retard abordés dans ce travail, ses contenus principaux sont présentés dans une série de chapitres majoritairement indépendants.

Trois chapitres traitent de divers aspects des équations aux différences linéaires. Le premier s'intéresse à des questions de contrôlabilité pour une équation aux différences linéaire avec une quantité finie de retards ponctuels, se servant de la théorie de la réalisation pour obtenir une condition nécessaire et suffisante de contrôlabilité approchée de type Hautus dans $L^p$, $p \in [1, +\infty)$, ainsi qu'une condition nécessaire de contrôlabilité exacte dans $L^1$. Le deuxième étudie la robustesse de la stabilité exponentielle par rapport à des perturbations de paramètres pour une équation aux différences linéaire plus générale, étendant en particulier le critère de Hale--Silkowski, développé à l'origine pour des équations avec une quantité finie de retards ponctuels, à des équations comportant également des retards distribués. Enfin, le troisième chapitre considère une équation aux différences linéaire avec un seul retard variant en temps, étudiant l'existence et l'unicité de ses solutions ainsi que son comportement asymptotique dans quelques cadres fonctionnels, soulignant l'importance de l'espace fonctionnel dans une telle analyse.

Le chapitre suivant présente plusieurs résultats sur la propriété de dominance induite par multiplicité (MID, d'après l'expression en anglais \foreignlanguage{english}{\emph{mul\-ti\-pli\-city-in\-duced dominance}}) pour des équations différentielles à retard de types retardé et neutre, principalement dans le cas d'un seul retard. On présente d'abord un résultat général établissant la propriété MID pour des racines réelles de multiplicité maximale, avant une série de discussions et résultats dans le cas de racines réelles de multiplicité suffisamment grande. Les cas de la propriété MID pour des racines complexes conjuguées ou pour des systèmes à plus d'un retard sont aussi analysés pour certaines équations de petit ordre. Le chapitre se conclut par une série d'exemples et d'applications illustrant les usages et limitations de la propriété MID en théorie du contrôle.

Deux boîtes à outils logicielles pour l'analyse et le contrôle de systèmes à retard, YALTAPy et P3$\delta$, sont présentées dans le chapitre suivant, avec une courte description de leurs principales fonctionnalités illustrée par quelques exemples. Enfin, ce document se conclut par un chapitre présentant plusieurs problèmes ouverts sur les sujets des chapitres précédents, détaillant les difficultés et perspectives pour chacun d'entre eux.

\vfill

\noindent \textbf{Mots-clés.} Systèmes à retard, équations aux différences, équations différentielles à retard, contrôlabilité, stabilité forte, comportement asymptotique, stabilité exponentielle, stabilisation, théorie de la réalisation, méthodes spectrales, placement de pôles, dominance induite par multiplicité.

\medskip

\noindent \textbf{2020 Mathematics Subject Classification.} Principal: 93-02, 39-02, 39A06, 45A05, 93C43, 93B05, 93D23; Secondaire: 93B55, 39A30, 93D09, 33C15, 41A21, 35Q49, 34K43.
\end{otherlanguage}

\pdfbookmark[1]{Abstract}{abstract-en}

\begin{center}
\textbf{\Large\ThesisTitle}
\end{center}

\noindent This monograph presents a selection of the author's recent works on the control, stability, and stabilization of some classes of linear time-delay systems. After a detailed introduction to topics related to time-delay systems covered in the monograph, its main contents are presented in a series of mostly independent chapters.

Three chapters cover different aspects of linear difference equations. The first one addresses controllability issues for a linear difference equation with finitely many pointwise delays, making use of realization theory to deduce a Hautus-type necessary and sufficient approximate controllability condition in $L^p$, $p \in [1, +\infty)$, as well as a necessary condition for exact controllability in $L^1$. The second one studies the robustness of exponential stability with respect to perturbations of parameters for a more general linear difference equation, extending in particular the Hale--Silkowski criterion, originally developed for equations with finitely many pointwise delays, to equations also involving distributed delays. Finally, the third one considers a linear difference equation with a single time-varying delay, studying existence and uniqueness of its solutions as well as its asymptotic behavior in several functional settings and highlighting the importance of the functional space in such an analysis.

The next chapter presents several results on the mul\-ti\-pli\-city-in\-duced dominance (MID) property for delay-differential equations of retarded and neutral types, mostly in the case of a single delay. We first provide a general result establishing the MID property for real roots of maximal multiplicity, before a series of discussions and results on the case of real roots with large enough multiplicity. The cases of the MID property for complex-conjugate roots or for systems with more than one delay are also examined for some low-order equations. The chapter is concluded by a series of examples and applications illustrating the uses and limitations of the MID property in control theory.

Two software packages for the analysis and control of time-delay systems, YALTAPy and P3$\delta$, are presented in the subsequent chapter, with a brief description of their main functionalities illustrated by some examples. Finally, the monograph is concluded by a chapter presenting several open problems on the topics from the previous chapters, discussing difficulties and perspectives for each of them.

\vfill

\noindent \textbf{Keywords.} Time-delay systems, difference equations, delay-differential equations, controllability, strong stability, asymptotic behavior, exponential stability, stabilization, realization theory, spectral methods, pole placement, multiplicity-induced dominance.

\medskip

\noindent \textbf{2020 Mathematics Subject Classification.} Primary: 93-02, 39-02, 39A06, 45A05, 93C43, 93B05, 93D23; Secondary: 93B55, 39A30, 93D09, 33C15, 41A21, 35Q49, 34K43.
\pdfbookmark[1]{Resumo}{abstract-pt-br}

\begin{otherlanguage}{brazil}

\begin{center}
\textbf{\Large Controle, estabilidade e estabilização de algumas classes de sistemas lineares com atraso}
\end{center}

\noindent Este documento contém uma seleção dos trabalhos recentes do autor sobre controle, estabilidade e estabilização de algumas classes de sistemas lineares com atraso. Após uma introdução detalhada aos tópicos relacionados a sistemas com atraso abordados neste trabalho, seu conteúdo principal é apresentado em uma série de capítulos majoritariamente independentes.

Três capítulos abordam diferentes aspectos de equações de diferenças lineares. O primeiro trata de questões de controlabilidade para uma equação de diferenças linear com uma quantidade finita de atrasos pontuais, usando a teoria da realização para obter condições necessárias e suficientes de controlabilidade aproximada do tipo Hautus em $L^p$, $p \in [1, +\infty)$, bem como uma condição necessária de controlabilidade exata em $L^1$. O segundo estuda a robustez da estabilidade exponencial em relação a perturbações dos parâmetros de uma equação de diferenças linear mais geral, generalizando em particular o critério de Hale--Silkowski, desenvolvido originalmente para equações com uma quantidade finita de atrasos pontuais, para equações que também envolvem atrasos distribuídos. Finalmente, o terceiro capítulo considera uma equação de diferenças linear com um único atraso variante no tempo, estudando existência e unicidade de soluções bem como seu comportamento assintótico em vários espaços funcionais, destacando a importância do espaço funcional para este tipo de análise.

O capítulo seguinte apresenta diversos resultados sobre a propriedade de dominância induzida por multiplicidade (MID, do inglês \foreignlanguage{english}{\emph{mul\-ti\-pli\-city-in\-duced dominance}}) para equações diferenciais com atraso dos tipos retardado ou neutro, principalmente no caso de um único atraso. Apresentamos primeiro um resultado geral demonstrando a propriedade MID para raízes reais de multiplicidade máxima, seguido por uma série de discussões e resultados para o caso de raízes reais com multiplicidade suficientemente grande. Os casos da propriedade MID para raízes complexas conjugadas ou para sistemas com mais de um atraso também são examinados para algumas equações de ordem baixa. O capítulo é concluído por uma série de exemplos e aplicações que ilustram os usos e as limitações da propriedade MID em teoria de controle.

Dois softwares para a análise e o controle de sistemas com atraso, YALTAPy e P3$\delta$, são apresentados no capítulo seguinte, com uma breve descrição de suas principais funcionalidades ilustrada por exemplos. Finalmente, este trabalho é concluído com um capítulo apresentando diversos problemas abertos sobre os tópicos dos capítulos anteriores, discutindo dificuldades e perspectivas para cada um deles.

\vfill

\noindent \textbf{Palavras-chave.} Sistemas com atraso, equações de diferenças, equações diferenciais com atraso, controlabilidade, estabilidade forte, comportamento assintótico, estabilidade exponencial, estabilização, teoria da realização, métodos espectrais, alocação de polos, dominância induzida por multiplicidade.

\medskip

\noindent \textbf{2020 Mathematics Subject Classification.} Principais: 93-02, 39-02, 39A06, 45A05, 93C43, 93B05, 93D23; Secundários: 93B55, 39A30, 93D09, 33C15, 41A21, 35Q49, 34K43.
\end{otherlanguage}

\pdfbookmark[1]{\contentsname}{table-of-contents}
\tableofcontents
\cleardoublepage

\phantomsection
\addcontentsline{toc}{chapter}{\listfigurename}
\listoffigures
\cleardoublepage

\phantomsection
\addcontentsline{toc}{chapter}{\listtablename}
\listoftables
\cleardoublepage

\newenvironment{notationtable}{\renewcommand{\arraystretch}{1.3}
\begin{longtable}{@{} >{\raggedright} p{0.18\textwidth} @{\hspace*{0.02\textwidth}} p{0.8\textwidth} @{}}
\toprule
\textbf{Notation} & \textbf{Meaning} \tabularnewline
\midrule}{\bottomrule
\end{longtable}}

\chapter*{Notation}
\addcontentsline{toc}{chapter}{Notation}
\markboth{Notation}{Notation}

\label{notation}

\section*{Sets of numbers}
\addcontentsline{toc}{section}{Sets of numbers}
\markright{Sets of numbers}

\begin{notationtable}
$\mathbb N$, $\mathbb Z$, $\mathbb Q$, $\mathbb R$, $\mathbb C$ & Sets of nonnegative integers, of integers, of rational numbers, of real numbers, and of complex numbers, respectively. \tabularnewline
$\mathbb A_+$, $\mathbb A_-$ & For $\mathbb A \in \{\mathbb N, \mathbb Z, \mathbb Q, \mathbb R\}$, $\mathbb A_+$ denotes the set $\mathbb A \cap [0, +\infty)$ and $\mathbb A_-$ denotes the set $\mathbb A \cap (-\infty, 0]$. In particular, $\mathbb Z_+ = \mathbb N_+ = \mathbb N$ and $\mathbb N_- = \{0\}$. \tabularnewline
$\mathbb A^\ast$ & For $\mathbb A \in \{\mathbb N, \mathbb N_+, \mathbb N_-, \mathbb Z, \mathbb Z_+, \mathbb Z_-, \mathbb Q, \mathbb Q_+, \mathbb Q_-, \mathbb R, \mathbb R_+, \mathbb R_-, \mathbb C\}$, $\mathbb A^\ast$ denotes the set $\mathbb A \setminus \{0\}$. In particular, $\mathbb R_-$, $\mathbb R_-^\ast$, $\mathbb R_+$, and $\mathbb R_+^\ast$ are the intervals $(-\infty, 0]$, $(-\infty, 0)$, $[0, +\infty)$, $(0, +\infty)$, respectively. \tabularnewline
$\llbracket a, b\rrbracket$ & $[a, b] \cap \mathbb Z$, where $a, b \in \mathbb R$ and $[a, b] = \emptyset$ if $a > b$. \tabularnewline
$\lvert a, b]$, $\lvert a, b)$, $[a, b\rvert$, $(a, b\rvert$, $\lvert a, b\rvert$ & For $a, b \in \overline{\mathbb R}$, $\lvert a, b]$ is the interval $[a, b]$ if $a > -\infty$ and $(a, b]$ if $a = -\infty$. The other intervals are defined in a similar manner. \tabularnewline
\end{notationtable}

\section*{Elementary functions}
\addcontentsline{toc}{section}{Elementary functions}
\markright{Elementary functions}

\begin{notationtable}
$\id_X$ & Identity map in a set $X$. The subscript $X$ is omitted when it is clear from the context. \tabularnewline
$\floor{x},\, \ceil{x}$ & Floor and ceiling of $x \in \mathbb R$, denoting the unique integers satisfying $x - 1 < \floor{x} \leq x$ and $x \leq \ceil{x} < x + 1$. This notation is extended componentwise to vectors and matrices. \tabularnewline
$\sign(x)$ & Sign of $x$, defined as $1$ if $x > 0$, $-1$ if $x < 0$, or $0$ if $x = 0$. \tabularnewline
$x_+$, $x_-$ & Positive and negative parts of $x$, defined, respectively, as $\max(x, 0)$ and $\max(-x, 0)$. \tabularnewline
$\1_A$ & Indicator function of the set $A$, equal to $1$ in $A$ and $0$ elsewhere. \tabularnewline
$\log_b a$ & Logarithm of $a > 0$ in base $b > 1$. The base $b$ is omitted when $b = 10$. \tabularnewline
$\ln a$ & Natural logarithm of $a > 0$. \tabularnewline
$\Real z,\, \Imag z$ & Real and imaginary parts of $z \in \mathbb C$, respectively. \tabularnewline
$\overline z$ & Complex conjugate of $z \in \mathbb C$. \tabularnewline
$\binom{n}{k}$ & Binomial coefficient, defined by $\binom{n}{k} = \frac{n!}{k! (n-k)!}$ if $k$ and $n$ integers with $k \in \llbracket 0, n\rrbracket$, and by $\binom{n}{k} = 0$ if $n < 0$, $k < 0$, or $k > n$.\tabularnewline
$(\alpha)_k$ & Pochhammer symbol for the ascending factorial for $\alpha \in \mathbb C$ and $k \in \mathbb N$, defined inductively by $(\alpha)_0 = 1$ and $(\alpha)_{k+1} = (\alpha+k) (\alpha)_k$. \tabularnewline
\end{notationtable}

\section*{Function notation and operations with functions}
\addcontentsline{toc}{section}{Function notation and operations with functions}
\markright{Function notation and operations with functions}

\begin{notationtable}
$f(a^+)$, $f(a^-)$ & $\displaystyle \lim_{x \to a^+} f(x)$ and $\displaystyle \lim_{x \to a^-} f(x)$, respectively. \tabularnewline
$f(A)$ & Image under the function $f\colon X \to Y$ of the set $A \subset X$. \tabularnewline
$f\rvert_{C}$ & Restriction of the function $f\colon A \to B$ to the set $C \subset A$, defined as $f\rvert_{C}\colon C \to B$ by setting $f\rvert_{C}(x) = f(x)$ for every $x \in C$. \tabularnewline
$f^{\prime}$, $f^{\prime\prime}$, $f^{(k)}$ & Derivatives of the function $f$ of a single variable of order $1$, $2$, and $k$, respectively, for $k \in \mathbb N$, with the convention $f^{(0)} = f$. \tabularnewline
$\partial^k_{x_1 \dotso x_k} f$ & Partial derivative of the function $f$ of order $k \in \mathbb N^\ast$ with respect to the variables named $x_k, \dotsc, x_1$, sometimes also denoted as $\frac{\partial^k f}{\partial x_1 \dotsm \partial x_k}$. The index $k$ is omitted when $k = 1$. \tabularnewline
$\partial^k_{i_1 \dots i_k} f$ & If $f$ is a function of $d$ variables named $x_1, \dotsc, x_d$ and $(i_1, \dotsc, i_k) \in \llbracket 1, d\rrbracket^k$, then $\partial^k_{i_1 \dots i_k} f$ is a shorthand notation for $\partial^k_{x_{i_1} \dots x_{i_k}} f$. The index $k$ is omitted when $k = 1$. \tabularnewline
$f * g$ & Convolution product between $f$ and $g$. \tabularnewline
$\V_X(f)$ & Total variation over the set $X$ of the function $f$. \tabularnewline
$\V_a^b(f)$ & Shorthand notation for $\V_{[a, b]}(f)$. \tabularnewline
$\supp f$ & Support of the function $f$, defined as $\supp f = \overline{\{x \suchthat f(x) \neq 0\}}$. \tabularnewline
$\deg P$ & Degree of the polynomial function $P$. \tabularnewline
$O(\cdot)$, $o(\cdot)$ & Landau notation: $f(x) = O(g(x))$ as $x \to x_0$ if there exist $C > 0$ and a neighborhood of $x_0$ in which $\abs{f(x)} \leq C \abs{g(x)}$, and $f(x) = o(g(x))$ as $x \to x_0$ if $\lim_{x \to x_0} \frac{f(x)}{g(x)} = 0$. \tabularnewline
\end{notationtable}

\section*{Topology and metric spaces}
\addcontentsline{toc}{section}{Topology and metric spaces}
\markright{Topology and metric spaces}

\begin{notationtable}
$\partial A$ & Boundary of the set $A$. \tabularnewline
$\overline A$ & Closure of the set $A$. \tabularnewline
$B_{\mathsf X}(x, r)$, $\overline B_{\mathsf X}(x, r)$ & Open and closed balls, respectively, in the metric space $\mathsf X$ centered at $x$ and with radius $r$. If $\mathsf X = \mathbb R^d$ with the usual Euclidean distance, we use simply $d$ as index. The index is omitted when clear from the context. \tabularnewline
\end{notationtable}

\clearpage

\section*{Matrices, vector spaces, and linear operators}
\addcontentsline{toc}{section}{Matrices, vector spaces, and linear operators}
\markright{Matrices, vector spaces, and linear operators}

\begin{notationtable}
$e_1, \dotsc, e_d$ & Canonical basis of $\mathbb R^d$ or $\mathbb C^d$. \tabularnewline
$\mathcal M_{d, m}(\mathbb K)$ & Set of $d \times m$ matrices with coefficients in $\mathbb K \in \{\mathbb R, \mathbb C\}$. \tabularnewline
$\mathcal M_d(\mathbb K)$ & $\mathcal M_{d, d}(\mathbb K)$. \tabularnewline
$\Id_d$ & Identity matrix in $\mathcal M_d(\mathbb R)$ or $\mathcal M_d(\mathbb C)$. The subscript $d$ is omitted when clear from the context. \tabularnewline
$0_{d, m}$ & Zero matrix in $\mathcal M_{d, m}(\mathbb R)$ or $\mathcal M_{d, m}(\mathbb C)$, denoted simply by $0_d$ if $m = d$ or by $0$ when $d$ and $m$ are clear from the context. \tabularnewline
$\begin{pmatrix} A & B \end{pmatrix}$ & Matrix in $\mathcal M_{d, m_1 + m_2}(\mathbb K)$ obtained by juxtaposition of the matrices $A \in \mathcal M_{d, m_1}(\mathbb K)$ and $B \in \mathcal M_{d, m_2}(\mathbb K)$. \tabularnewline
$A^\transp$, $A^\ast$ & Transpose and Hermitian transpose of the matrix $A \in \mathcal M_{d, m}(\mathbb K)$, respectively.\tabularnewline
$A^{-1}$ & Inverse of the matrix $A \in \mathcal M_d(\mathbb K)$, when it exists. \tabularnewline
$\Tr(A)$, $\det(A)$, $\Adj(A)$ & Trace, determinant, and adjugate of the matrix $A \in \mathcal M_d(\mathbb K)$, respectively. \tabularnewline
$\Span(V)$ & Vector space spanned by the vectors of the set $V$. \tabularnewline
$\Ran A$, $\rank A$ & Range and rank of a matrix $A \in \mathcal M_{d, m}(\mathbb K)$, respectively, or of a linear operator $A$. \tabularnewline
$\Ker A$ & Kernel of a matrix $A \in \mathcal M_{d, m}(\mathbb K)$ or of a linear operator $A$. \tabularnewline
$\sigma(A)$ & Spectrum of the matrix $A \in \mathcal M_d(\mathbb C)$ or of the operator $A$ in the Banach space $\mathsf X$. \tabularnewline
$\spr(A)$ & Spectral radius of the matrix $A \in \mathcal M_d(\mathbb C)$ or of the bounded operator $A$ in the Banach space $\mathsf X$, defined as $\sup\{\abs{\lambda} \suchthat \lambda \in \sigma(A)\}$. \tabularnewline
$\abs{\cdot}$ & An arbitrary norm in $\mathbb R^d$ or $\mathbb C^d$ and the corresponding induced matrix norm in $\mathcal M_{d, m}(\mathbb R)$ or $\mathcal M_{d, m}(\mathbb C)$. \tabularnewline
$\abs{\cdot}_p$ & $\ell_p$ norm of a vector in $\mathbb R^d$ or $\mathbb C^d$ and the induced matrix norm in $\mathcal M_{d, m}(\mathbb R)$ or $\mathcal M_{d, m}(\mathbb C)$. \tabularnewline
$\scalprod{\cdot, \cdot}_{\mathsf H}$ & Inner product in the (real or complex) Hilbert space $\mathsf H$. When $\mathsf H$ is complex, $\scalprod{\cdot, \cdot}_{\mathsf H}$ is assumed to be anti-linear in the first variable and linear in the second one. The index $\mathsf H$ is omitted when clear from the context. \tabularnewline
$\scalprod{\cdot, \cdot}_{\mathsf X^\prime, \mathsf X}$ & Duality product in the topological vector space $\mathsf X$ and its topological dual $\mathsf X^\prime$. The index $\mathsf X^\prime, \mathsf X$ is omitted when clear from the context. \tabularnewline
$x \cdot y$ & $\scalprod{x, y}_{\mathbb R^d}$ for $x, y \in \mathbb R^d$. \tabularnewline
\end{notationtable}

\clearpage

\section*{Measure theory}
\addcontentsline{toc}{section}{Measure theory}
\markright{Measure theory}

\begin{notationtable}
$\mathfrak B$, $\mathfrak L$ & Borel and Lebesgue $\sigma$-algebras, respectively, in $\mathbb R^d$ or in any of its (measurable) subsets. \tabularnewline
$\mathcal L(A)$, $\abs{A}$ & Lebesgue measure of the measurable set $A \subset \mathbb R^d$. \tabularnewline
$\V_X(\mu)$ & Total variation over the set $X$ of the measure $\mu$. \tabularnewline
$\V_a^b(\mu)$ & Shorthand notation for $\V_{[a, b]}(\mu)$. \tabularnewline
$f_{\#} \mu$ & Pushforward of the measure $\mu$ by the map $f$. \tabularnewline
$\supp \mu$ & Support of the Borel measure $\mu$, defined as $\supp\mu = \{x \suchthat \mu(N) > 0 \text{ for any open } N \ni x\}$. \tabularnewline
\end{notationtable}

\section*{Functional spaces}
\addcontentsline{toc}{section}{Functional spaces}
\markright{Functional spaces}

\begin{notationtable}
$L^p(X, Y)$ & Usual Lebesgue space of $p$-integrable functions defined in $X$ and taking values in $Y$. In this notation and the others from this table, the set $Y$ is omitted from the notation if $Y = \mathbb R$, and we write simply $L^p$ if $X$ and $Y$ are clear from the context. \tabularnewline
$L^p_{\loc}(X, Y)$ & Set of all locally $p$-integrable functions (with respect to the Lebesgue measure) defined in $X$ and taking values in $Y$. \tabularnewline
$W^{k, p}(X, Y)$,\newline $H^k(X, Y)$ & Usual Sobolev space of $k$-times weakly differentiable functions with de\-ri\-va\-tives in $L^p$ defined in $X$ and taking values in $Y$, the notation $H^k$ being used when $p = 2$. \tabularnewline
$\mathcal C(X, Y)$ & Set of continuous functions from the topological space $X$ to the topological space $Y$. \tabularnewline
$\mathcal C^k(A, Y)$ & Set of $k$ times differentiable $Y$-valued functions defined on the set $A \subset \mathbb R^d$, for $k \in \mathbb N$ and $Y$ a Banach space with its strong topology. \tabularnewline
$\mathcal C^k_{\mathrm{c}}(A, Y)$ & Subset of $\mathcal C^k(A, Y)$ of all the compactly supported functions. \tabularnewline
\end{notationtable}
\label{end:notation}

\chapter*{Acronyms}
\addcontentsline{toc}{chapter}{Acronyms}
\markboth{Acronyms}{Acronyms}

\noindent
{\renewcommand{\arraystretch}{1.3}
\begin{longtable}{@{} m{0.18\textwidth} @{\hspace*{0.02\textwidth}} m{0.68\textwidth} @{\hspace*{0.02\textwidth}} m{0.1\textwidth}@{}}
\toprule
\textbf{Acronym} & \textbf{Meaning} & \textbf{Page} \tabularnewline
\midrule
CRRID & Coexisting real roots-induced dominance & \pageref{CRRID} \tabularnewline
DAE & Differential-algebraic equation & \pageref{DAE} \tabularnewline
DDAE & Delay differential-algebraic equation & \pageref{DDAE} \tabularnewline
DDE & Delay-differential equation & \pageref{DDE} \tabularnewline
FDE & Functional differential equation & \pageref{FDE} \tabularnewline
LMI & Linear matrix inequality & \pageref{LMI} \tabularnewline
MID & Multiplicity-induced-dominancy & \pageref{MID} \tabularnewline
ODE & Ordinary differential equation & \pageref{ODE} \tabularnewline
PDE & Partial differential equation & \pageref{PDE} \tabularnewline
TDS & Time-delay system & \pageref{TDS} \tabularnewline
\bottomrule
\end{longtable}}

\chapter*{Preface}
\addcontentsline{toc}{chapter}{Preface}
\markboth{Preface}{Preface}

\epigraph{\foreignlanguage{brazil}{Ai! que preguiça!}}{\foreignlanguage{brazil}{\emph{Macunaíma}, Mario de Andrade}}

{\setlength{\parskip}{1pt plus 5pt minus 1pt}

This monograph provides a presentation of some of the author's recent contributions to the control, stability, and stabilization of some classes of linear time-delay systems. The contents presented here cover a broad spectrum of the author's research between 2020 and 2025 and, more precisely, the monograph reports on the results of the works \cite{Chitour2023Hautus, GoncalvesNetto2025Strong, ChitourStrong, MazantiWellPosedness, Boussaada2025Decay, Fueyo2023Pole, Boussaada2022Generic, Boussaada2022Some, Mazanti2021Multiplicity, Boussaada2025OverOrder, Boussaada2022MID, Boussaada2022Pade, Fueyo2021Insights, Niculescu2021Stability, Mazanti2021Effects, Benarab2020MID, Mazanti2020Spectral, Mazanti2020Qualitative, Cavalera2022YALTAPy, Boussaada2022New, Boussaada2020Partial, Boussaada2021New}. Not included in this monograph are the author's older works \cite{Mazanti2013Stabilization}, \cite{Mazanti2014Stabilization}, and \cite{Chitour2020Approximate, Colonius2019Decay, Mazanti2017Relative, Chitour2017Persistently, Chitour2016Stability}, whose contents were reported, respectively, in the Master's theses \cite{Mazanti2011Stabilisation} and \cite{Mazanti2013Condition} and in the PhD thesis \cite{Mazanti2016Stability}. Also not included are works on other of the author's research topics, such as switched systems \cite{Mazanti2011Sistemas, Chitour2013Stabilization, Chitour2023Gap, Chitour2021Gap}, hyperbolic systems \cite{Balogoun2025Novel, Chitour2024Approximate, Chitour2021One, Tang2017Stability, Arias2023Frequency}, and mean field games and optimization \cite{MachadoNash, Mazanti2024Note, Sadeghi2022Nonsmooth, Sadeghi2022Multi, Ducasse2022Second, Dweik2020Sharp, Mazanti2019Minimal, Sadeghi2021Characterization, Mazanti2025Nonsmooth}, as well as the recent results on time-delay systems from \cite{Michiels2025Relations, Boussaada2025Iterated}, which were obtained after the corresponding parts of this monograph were concluded. Writing of this monograph started in December~2022.

The monograph starts in Chapter~\ref{chap:introduction} by a general introduction to time-delay systems and the main topics addressed in the sequel. We introduce systems of retarded and neutral types in a rather general setting, as well as continuous-time difference equations and delay differential-algebraic equations, which are the subject of some of the subsequent chapters. The main topics of the subsequent chapters, such as controllability of difference equations, spectral methods for linear time-delay systems, and the multiplicity-induced dominance property are successively introduced in the remainder of Chapter~\ref{chap:introduction}. While part of the material of Chapter~\ref{chap:introduction} comes from the introductions of works on which this monograph is based, these contents were expanded and put into a wider context, and additional comments, examples, and references were provided whenever suitable in order to provide a comprehensive and self-contained introduction to the topics of this monograph.

Chapters~\ref{chap:controllability}--\ref{chap:diff-eqn-varying-delay} deal with some topics in the control and stability of continuous-time difference equations, i.e., equations in which the instantaneous state of the system at time $t$ is a function of the same state at previous times and possibly of some control terms (see Section~\ref{sec:intro-cont-time-diff-eqn}). Chapter~\ref{chap:controllability} is based on \cite{Chitour2023Hautus} and reports on controllability criteria obtained in that reference for linear difference equations with (finitely many) constant delays in finite time in $L^p$ spaces, $p \in [1, +\infty)$. The main technical tool used for that purpose is the realization theory developed by Y.~Yamamoto for general infinite-dimensional dynamical systems (see Section~\ref{sec:Realization}). Approximate controllability is characterized in terms of a Hautus-type controllability criterion, and we also provide a necessary condition for $L^1$ exact controllability, which can be seen as the closure of the $L^1$ approximate controllability criterion. Furthermore, we give an explicit upper bound on the minimal times of approximate and exact controllability. The material from Chapter~\ref{chap:controllability} follows very closely the article \cite{Chitour2023Hautus}, one of the main novelties from this monograph being a more pedagogical introduction to the required tools from distribution theory in Section~\ref{sec:distributions} with respect to \cite[Section~2]{Chitour2023Hautus}. Readers interested only in the main results on this topic may restrict their attention to Section~\ref{sec:statements}, which provides the statements of the controllability criteria in Theorems~\ref{main_result1} and \ref{main_result2}.

As discussed in Section~\ref{sec:sensitivity}, the exponential stability of a linear difference equation with finitely many pointwise delays is very sensitive to perturbations in its delays and, for this reason, stronger notions of stability were introduced in the literature. The aim of Chapter~\ref{chapter:strong-stability} is to extend this analysis to a more general class of linear difference equations, namely \eqref{sys_general}, reporting on the results from the article \cite{GoncalvesNetto2025Strong} and the submitted manuscript \cite{ChitourStrong}. Most of the material of Chapter~\ref{chapter:strong-stability} comes from these two references, the main novelty of the chapter being a unified presentation of both. A particular effort has been put into clarifying the roles played by the different functional spaces $\WellPosed$, $\NoSingular$, $\LowerBound$, and $\SDG$ introduced in Sections~\ref{sec:well-posed-strong-stability} and \ref{sec:definition-spaces} (whose notations are modified with respect to \cite{ChitourStrong} for improving clarity). Some results that are only implicitly stated in \cite{GoncalvesNetto2025Strong, ChitourStrong}, or that follow easily from the discussion in those references, were made explicit in Chapter~\ref{chapter:strong-stability}, such as Proposition~\ref{prop:TV-less-1-implies-global-strong-stab} and Corollary~\ref{coro:lss-implies-rho-leq-1}. The definitions of strong stability from Definitions~\ref{def_strong-stab} and \ref{def_stron-stab-UR} were modified with respect to the corresponding ones in \cite{GoncalvesNetto2025Strong, ChitourStrong}, which allows one to formulate our main results and conjectures in a clearer way.

The difference equations considered in Chapters~\ref{chap:controllability} and \ref{chapter:strong-stability} only deal with the case of constant delays, but, as illustrated by Example~\ref{expl:AvellarMarconato} in Section~\ref{sec:intro-varying}, time-varying delays introduce many subtleties in the asymptotic behavior of difference equations. Chapter~\ref{chap:diff-eqn-varying-delay} aims to illustrate many other subtleties induced by time-varying delays by considering a linear difference equation with a single time-varying delay and thoroughly studying questions on the existence and uniqueness of its solutions and their asymptotic behavior. The contents of Chapter~\ref{chap:diff-eqn-varying-delay} come from the submitted manuscript \cite{MazantiWellPosedness}, and, as the conclusion of the writing of the chapter precedes that of \cite{MazantiWellPosedness}, the contents of the chapter present no major novelty with respect to the latter. Chapter~\ref{chap:diff-eqn-varying-delay} provides several results on the existence and uniqueness of solutions of the linear difference equation with a time-varying delay of interest, which are different according to which kind of solution one is interested in (pointwise, continuous, regulated, measurable, almost everywhere, or $L^p$ solutions). It then proceeds to characterize asymptotic and exponential stability for these kinds of solutions. The results from Chapter~\ref{chap:diff-eqn-varying-delay} show in particular that the same system can be asymptotically stable for one family of solutions but unstable for others, due to the variations in time of the delay. It is the author's belief that the proofs presented in \cite{MazantiWellPosedness} have rather simple and instructive underlying ideas and that their details actually provide much insight into the phenomena induced by time variations in the delays. For this reason, complete proofs are provided all along Chapter~\ref{chap:diff-eqn-varying-delay}.

Moving away from difference equations, Chapter~\ref{chap:MID} considers linear time-invariant scalar delay-differential equations of retarded or neutral types with a single delay (with the exception of Section~\ref{sec:MID-two-delays}, which considers a simple equation with two delays) and is motivated by a stabilizability problem: when some of the coefficients of the equation can be chosen, how to choose them in order to render its zero solution exponentially stable? Chapter~\ref{chap:MID} focuses on a property that allows one to address this question, the multiplicity-induced domi\-nance (MID\label{MID}) property, a fact previously observed in the literature stating that a real spectral value with the maximal possible multiplicity tends to be the rightmost element in the spectrum of a time-delay system, determining thus its asymptotic behavior. After some preliminary discussions on Padé approximations and Kummer and Whittaker functions, Chapter~\ref{chap:MID} shows the validity of the MID property for linear scalar delay-differential equations of arbitrary order with a single delay and of retarded or neutral types, discussing in the sequel recent results on related questions, including dominance for spectral values with large multiplicity, but smaller than the maximal, dominance of pairs of complex conjugate roots with maximal multiplicity for low-order systems, and dominance in the case of systems with more than one delay. The chapter also presents an extensive list of examples.

Much effort was devoted to the presentation of Chapter~\ref{chap:MID}, which is based on the journal articles \cite{Boussaada2025Decay, Fueyo2023Pole, Boussaada2022Generic, Boussaada2022Some, Mazanti2021Multiplicity, Boussaada2025OverOrder} and the communications \cite{Boussaada2022MID, Boussaada2022Pade, Fueyo2021Insights, Niculescu2021Stability, Mazanti2021Effects, Benarab2020MID, Mazanti2020Spectral, Mazanti2020Qualitative}. Chapter~\ref{chap:MID} provides what the author believes to be a clear and complete presentation of the author's most important recent results on the MID property. Some proofs were rewritten with respect to the published versions in order to make them clearer, while others were adapted or simplified whenever possible, and additional discussion was included whenever necessary, so that Chapter~\ref{chap:MID} can serve as a useful reference for anyone wishing to understand what is the MID property and what is currently known about it.

More precisely, the presentation of Padé approximants from Section~\ref{sec:Pade} was expanded with respect to the corresponding one in \cite{Boussaada2022Pade}. The preliminary material from Section~\ref{sec:Kummer-Whittaker} comes from \cite{Boussaada2022Some}, but the proof of Proposition~\ref{PropTsvetkoffCorrected} was rewritten in order to properly take into account the singularity at $0$ by a limiting argument, which was not detailed in \cite{Boussaada2022Some}. As for Section~\ref{sec:MID-maximal}, its main result, Theorem~\ref{thm:main-MID}, was first proved in \cite{Mazanti2020Qualitative} in the case $n = 2$, $m = 1$, then generalized to the case of any $n \in \mathbb N^\ast$ and $m = n - 1$ in \cite{Mazanti2021Multiplicity}, and a final generalization was obtained in \cite{Boussaada2022Generic} in order to treat the case of any nonnegative integers $n$ and $m$. However, all those proofs used rather long computations based on properties of binomial coefficients. It was remarked in \cite{Boussaada2022Pade} that such computations can be avoided if one notices the link between the characteristic function of a delay-differential equation with a root of maximal multiplicity and the Padé approximant of the exponential function. The proof of Theorem~\ref{thm:main-MID} presented in Chapter~\ref{chap:MID} is thus based on that of \cite{Boussaada2022Pade} and uses this idea (which also allows one to obtain Theorem~\ref{thm:main-MID}\ref{item:Pade}, which is absent in the earlier works). The proof presented in Chapter~\ref{chap:MID} is an expanded version of that of \cite{Boussaada2022Pade}, taking additional care to handle the normalization condition that guarantees uniqueness of the Padé approximant. We also note that the material from Section~\ref{sec:MID-DDAE} is new: while an application of the MID property to delay-differential algebraic equations was previously proposed in \cite{Mazanti2021Effects}, that reference only considered the case $n = m = 1$ of \eqref{eq:DDAE-MID}, and the contribution of Section~\ref{sec:MID-DDAE} is the generalization to any nonnegative integers $n$ and $m$.

The contents of Section~\ref{sec:MID-intermediate} come from \cite{Boussaada2025OverOrder, Boussaada2025Decay, Boussaada2022MID}, and mostly from \cite{Boussaada2025OverOrder}, as the other two references deal only with particular cases. The presentation was modified with respect to those references in order to make it clearer and more integrated in the chapter and, in particular, Examples~\ref{expl:not-dominant} and \ref{expl:not-dominant-parametric} are new. With respect to the proof of the corresponding result in \cite{Boussaada2025OverOrder}, the proof of Theorem~\ref{thm:max-mult} in Chapter~\ref{chap:MID} provides additional details in order to properly justify the dimensions of the vector spaces introduced in the beginning of the proof. Section~\ref{sec:MID-complex} is based on \cite{Mazanti2020Spectral}, and its main result, Theorem~\ref{thm:complex-dominance}, is the main result of \cite{Mazanti2020Spectral}. However, that reference only provided a sketch of its proof, promising a full proof in a future paper that has not (yet) been completed. Section~\ref{sec:MID-complex} builds upon the ideas of the sketch of proof presented in \cite{Mazanti2020Spectral} in order to provide the promised full proof for the first time. Section~\ref{sec:MID-two-delays} reports on the results from \cite{Fueyo2021Insights, Fueyo2023Pole} without much novelty with respect to those references. Finally, Section~\ref{sec:MID-applis} collects and presents in a unified way a series of 10 examples and applications that appeared in the references \cite{Mazanti2020Qualitative, Mazanti2021Effects, Benarab2020MID, Mazanti2021Multiplicity, Boussaada2022Generic, Boussaada2025Decay, Fueyo2023Pole}.

In addition to the results from Chapters~\ref{chap:controllability}--\ref{chap:MID}, this monograph also provides, in Chapter~\ref{chap:delay-software}, a brief presentation of two software packages for time-delay systems, YALTAPy and P3$\delta$, to which the author has contributed. The former is the Python version of a Matlab package designed for the stability analysis of time-delay systems, including computation of stability windows and root loci for unstable roots as the delay varies in a given interval, while the latter is focused on the use of the MID property, and a similar property named coexisting real roots-induced dominance (CRRID\label{CRRID}), for the stabilization of delay-differential equations with a single delay. The exposition of Chapter~\ref{chap:delay-software} is based on \cite{Cavalera2022YALTAPy, Boussaada2022New, Boussaada2020Partial, Boussaada2021New}.

Far from simply answering questions, the contents of this monograph also highlight several open problems and conjectures on the control, stability, and stabilization of the classes of linear time-delay systems considered here. Chapter~\ref{chap:delay-open-problems} provides a nonexhaustive list of open problems that the author deems particularly interesting, presenting each one together with a discussion on what is known, possible techniques for addressing them, the difficulties involved, and some perspectives.

Finally, the author hopes that readers will find as much pleasure in reading this monograph as the author did in writing it.

\medskip

\hspace*{\fill} \emph{Paris, 11 February 2026}

\hspace*{\fill} \foreignlanguage{french}{\emph{23 pluviôse CCXXXIV}}

}

\pagestyle{fancy}
\renewcommand{\chaptermark}[1]{\markboth{\thechapter.\ #1}{\thechapter.\ #1}}
\renewcommand{\sectionmark}[1]{\markright{\thesection.\ #1}}
\fancyhf{}
\fancyhf[CEH]{\it \nouppercase{\leftmark}}
\fancyhf[COH]{\it \nouppercase{\rightmark}}
\fancyhf[ROF, LEF]{\it \thepage}
\renewcommand{\headrulewidth}{0.5pt}
\renewcommand{\footrulewidth}{0pt}
\pagenumbering{arabic}

\chapter{Introduction}
\label{chap:introduction}

\epigraph{\foreignlanguage{brazil}{Ai, que burro, dá zero pra ele!}}{\foreignlanguage{brazil}{\emph{Chaves}}}

\section{Time-delay systems}

Time delays, sometimes also referred to as time lags, appear naturally in the modeling of several phenomena, in particular when propagation of mass, energy, or information with finite propagation speed is involved, and when the focus of the study is not the propagation phenomenon itself but its impact on other systems. Models involving time delays are useful in several applications, including in biology, where delays may represent incubation periods, maturation times, or age structure; in engineering, where delays typically model latency effects due to signal propagation and processing speed; in chemistry, where delays usually represent reaction times; in economics, where delays may represent the time lag between an action being implement and its economic effects or the reaction time of economic agents to changes in their environment; or also in physics, where delays are typically linked to propagation phenomena. One usually refer to systems where time delays appear as \emph{time-delay systems} (TDSs\label{TDS}), but other terms are also used for this class of systems or some of its subclasses, such as \emph{functional differential equations} (FDEs\label{FDE}), \emph{differential-difference equations}, or \emph{delay-differential equations} (DDEs\label{DDE}), the latter being typically used for systems described by differential equations involving discrete delays.

Due to their applications and the challenging mathematical problems arising in their analysis, the theory of time-delay systems has seen an important development since the 1950s and 1960s \cite{Pinney1958Ordinary, Diekmann1995Delay, Bellman1963Differential, Krasovskii1963Stability, Stepan1989Retarded, MacDonald1989Biological, Kolmanovskii1986Stability, Hale1977Theory, Elsgolts1973Introduction, Wang1999Finite, Hale1993Introduction, Salamon1984Control, Halanay1966Differential, Gopalsamy1992Stability} and up to recent years \cite{Arino2006Delay, Gu2003Stability, Breda2015Stability, Li2015Analytic, Breda2023Controlling, Michiels2014Stability, Sipahi2011Stability, Insperger2011Semi, Valmorbida2019Delays, Witrant2016Recent, Valmorbida2022Accounting, Insperger2021Delay, Corduneanu2016Functional, Fridman2014Introduction}, even though important early works before the 1950s also exist, such as the works of Vito Volterra \cite{Volterra1928Theorie, Volterra1909Equazioni} highlighting the presence of time delays in models in biology and elasticity, where delays, called \emph{hereditary actions}, are described as being\footnote{``\foreignlanguage{french}{Actions actuelles dues à des causes qui ont existé dans le passé}'', in the original French text in \cite{Volterra1928Theorie}.} ``current actions due to causes that have existed in the past''; the work \cite{Minorsky1942Self} by Nicholas Minorsky highlighting the importance of taking the delay into account in the study of ship stabilization and steering; the work \cite{Hazen1934Theory} by Harold L.~Hazen on servomechanisms, highlighting the potential destabilizing effects that delays may have in the dynamics of a system; or the works \cite{Callender1936Time, Hartree1937Time}, showing the difficulty of handling delays in control loops for second- and third-order linear time-invariant systems.

\section{Retarded and neutral time-delay systems}

Time-delay systems are typically classified\footnote{Some classifications sometimes also include systems of \emph{advanced} type, which, roughly speaking, correspond to systems whose highest order of derivation appears only in delayed terms, such as a linear scalar equation of the form $x^\prime(t - 1) = a x(t) + b x(t - 2)$. For this equation, up to a time translation, the equation can be rewritten as $x^\prime(t) = a x(t + 1) + b x(t - 1)$, and hence the time derivative of $x$ at time $t$ depends on \emph{future} values of $t$ when $a \neq 0$. Advanced time-delay systems are ill-posed in general: for instance, the system $x^\prime(t) = x(t + 1)$ only admits globally defined solutions if its initial condition $x_0\colon [-1, 0] \to \mathbb R$ belongs to $\mathcal C^\infty([-1, 0], \mathbb R)$ and $x_0^{(k)}(0) = x_0^{(k)}(-1)$ for every $k \in \mathbb N$. In this monograph, we will not consider such systems, restricting our attention to the retarded and neutral cases.} into \emph{retarded} and \emph{neutral} types. A time-delay system of retarded type is usually informally defined as a system in which the highest order of derivation only appears in terms without delay. More precisely, a retarded time-delay system is a system that can be written under the form
\begin{equation}
\label{eq:intro-general-retarded}
x^\prime(t) = f(t, x_t), \qquad t \in \mathbb R_+,
\end{equation}
where $x(t) \in \mathbb R^d$ is the instantaneous state of the system and $x_t$ is the \emph{history function} at time $t$, defined for $s \in [-r, 0]$ by $x_t(s) = x(t + s)$, where $r > 0$ is a given number representing the maximal delay appearing in the system. Usually, $x_t$ is assumed to belong to some functional space, such as $\mathcal C([-r, 0], \mathbb R^d)$, and the dynamics of the equation are described by the function $f\colon \mathbb R_+ \times \mathcal C([-r, 0], \mathbb R^d) \to \mathbb R^d$. Important particular cases of \eqref{eq:intro-general-retarded} include systems with finitely many discrete delays, of the form
\[
x^\prime(t) = f(t, x(t), x(t - \tau_1(t)), \dotsc, x(t - \tau_N(t))),
\]
where $\tau_n\colon \mathbb R_+ \to [0, r]$ for $n \in \llbracket 1, N\rrbracket$ and $f \colon \mathbb R_+ \times (\mathbb R^d)^{N+1} \to \mathbb R^d$, and systems with distributed delays, of the form
\[
x^\prime(t) = \int_{-r}^0 f(t, s, x(t + s)) \diff s,
\]
where $f \colon \mathbb R_+ \times [-r, 0] \times \mathbb R^d \to \mathbb R^d$, among other classes.

In neutral systems, the highest order of derivation on the system appears both in a term without delay and a term with delay. More precisely, a neutral time-delay system is a system that can be written under the form
\begin{equation}
\label{eq:intro-general-neutral}
\frac{\df}{\df t} D(t, x_t) = f(t, x_t), \qquad t \in \mathbb R_+,
\end{equation}
where $x_t$ and $f$ are as in the retarded case and $D\colon \mathbb R_+ \times \mathcal C([-r, 0], \mathbb R^d) \to \mathbb R^d$ is assumed to be \emph{atomic at $0$}, a notion that, roughly speaking, means that $D(t, x_t)$ ``contains'' a term of the form $x(t)$. More precisely, we say that $D$ is atomic at $0$ if $D$ is twice Fréchet differentiable in its second variable with continuous derivatives and, for every $(t, \phi) \in \mathbb R_+ \times \mathcal C([-r, 0], \mathbb R^d)$, the matrix $\eta_{t, \phi}(0) - \eta_{t, \phi}(0^-)$ is invertible, where $\eta_{t, \phi}\colon [-r, 0] \to \mathcal M_d(\mathbb R)$ is a function of bounded variation representing the first Fréchet derivative of $D$ with respect to its second variable at $(t, \phi)$, $\partial_2 D(t, \phi)$, in the sense of Riesz representation theorem, i.e.,
\[
\partial_2 D(t, \phi) y = \int_{-r}^0 \diff[\eta_{t, \phi}(s)] y(s) \qquad \text{ for all } y \in \mathcal C([-r, 0], \mathbb R^d).
\]
Note that, with this definition, one can take $D(t, \phi) = \phi(0)$, in which case \eqref{eq:intro-general-neutral} reduces to \eqref{eq:intro-general-retarded}, i.e., retarded systems are a particular class of neutral systems. However, we will usually reserve the word ``neutral'' to neutral systems that are not retarded.

Similarly to retarded systems, important particular cases of neutral systems include systems with finitely many discrete delays, of the form
\[
x^\prime(t) + g(t, x^\prime(t - r_1(t)), \dotsc, x^\prime(t - r_M(t))) = f(t, x(t), x(t - \tau_1(t)), \dotsc, x(t - \tau_N(t))),
\]
where $\tau_n\colon \mathbb R_+ \to [0, r]$ for $n \in \llbracket 1, N\rrbracket$, $r_m\colon \mathbb R_+ \to (0, r]$ for $m \in \llbracket 1, M$, $f \colon \mathbb R_+ \times (\mathbb R^d)^{N+1} \to \mathbb R^d$, and $g \colon \mathbb R_+ \times (\mathbb R^d)^{M} \to \mathbb R^d$, and systems with distributed delays, of the form
\[
x^\prime(t) + \int_{-r}^0 g(t, s, x^\prime(t + s)) \diff s = \int_{-r}^0 f(t, s, x(t + s)) \diff s,
\]
where $f, g \colon \mathbb R_+ \times [-r, 0] \times \mathbb R^d \to \mathbb R^d$, among other classes.

Note that both retarded and neutral time-delay systems can be seen as systems in the infinite-dimensional Banach space $\mathcal C([-r, 0], \mathbb R^d)$. In particular, initial conditions for a time-delay system should be prescribed as continuous functions defined on the interval $[-r, 0]$. We refer the interested reader to \cite{Hale1993Introduction, Corduneanu2016Functional} for more details on retarded and neutral time-delay systems, including well-posedness results. We also remark that, although we presented here the theory only for delays taking values in a bounded interval, there also exists a longstanding theory of time-delay systems with unbounded delays, with classical results being presented in the survey \cite{Corduneanu1980Equations} and in \cite[Section~12.9]{Hale1993Introduction}. One may also replace the space of continuous functions $\mathcal C([-r, 0], \mathbb R^d)$ by other functional spaces, such as $L^p$ spaces, and we refer in particular to \cite{Kaashoek1992Characteristic, Delfour1972Controllability, Delfour1972Hereditary} for some interesting results on that topic.

\section{Continuous-time difference equations}
\label{sec:intro-cont-time-diff-eqn}

Another class of time-delay systems, which is sometimes seen as a particular or degenerate case of neutral systems, is that of \emph{continuous-time difference equations} (to which we will refer simply as \emph{difference equations} in the sequel, and which are sometimes also called \emph{difference delay equations} or \emph{delay difference equations}), which are equations of the form
\begin{equation}
\label{eq:intro-general-difference}
D(t, x_t) = 0, \qquad t \in \mathbb R_+,
\end{equation}
where $D\colon \mathbb R_+ \times \mathcal C([-r, 0], \mathbb R^d) \to \mathbb R^d$ is atomic at $0$. These equations do not involve time derivatives of the unknown function $x$, expressing instead just a relation between $x(t)$ and its past values in the interval $[t - r, t)$. Note that \eqref{eq:intro-general-difference} is very similar to \eqref{eq:intro-general-neutral} in the case $f = 0$: indeed, in this case, \eqref{eq:intro-general-neutral} reduces to $\frac{\df}{\df t} D(t, x_t) = 0$, i.e., $t \mapsto D(t, x_t)$ is constant, and \eqref{eq:intro-general-difference} considers the particular case where this constant is $0$.

In the linear time-invariant case, i.e., when $D$ is independent of its first variable and a bounded linear operator with respect to its second variable, \eqref{eq:intro-general-difference} can be rewritten, after isolating the atom of $D$ at $0$ and applying the Riesz--Markov--Kakutani representation theorem to the remaining part, as
\begin{equation}
\label{eq:intro-diff-eqn-integral-measure}
x(t) = \int_{-r}^0 \diff \mu(\theta) x(t + \theta), \qquad t \in \mathbb R_+,
\end{equation}
where $\mu$ is a matrix-valued Borel measure on $[-r, 0]$. Alternatively, \eqref{eq:intro-diff-eqn-integral-measure} can be rewritten in terms of a Lebesgue--Stieltjes integral as
\begin{equation}
\label{eq:intro-diff-eqn-integral-M}
x(t) = \int_{-r}^0 \diff M(\theta) x(t + \theta), \qquad t \in \mathbb R_+,
\end{equation}
where $M\colon [-r, 0] \to \mathcal M_d(\mathbb R)$ is the function of bounded variation obtained from the measure $\mu$ by setting $M(-r) = 0$ and $M(t) = \mu([-r, t])$ for $t \in (-r, 0]$. The general form of \eqref{eq:intro-diff-eqn-integral-measure} and \eqref{eq:intro-diff-eqn-integral-M}, with a matrix-valued function $M$, allows one to represent several phenomena, including pointwise delays, corresponding to discontinuity points of $M$, as well as a distributed delay, corresponding to the absolutely continuous part of $M$.

A particular case of \eqref{eq:intro-diff-eqn-integral-measure} that has been extensively considered in the literature (see, e.g., \cite{Hale2002Strong, Hale2003Stability, Hale1993Introduction, Henry1974Linear, Avellar1980Zeros, Avellar1990Difference, Chitour2023Hautus, Chitour2016Stability, Chitour2020Approximate, Cruz1970Stability, Mazanti2017Relative, Melvin1974Stability, Michiels2009Strong, Ngoc2015Exponential, OConnor1983Stabilization, Silkowski1976Star}) is when $\mu$ is a finite sum of Dirac measures on $[-r, 0)$ with constant matrix coefficients (or, equivalently, when $M$ from \eqref{eq:intro-diff-eqn-integral-M} is piecewise constant). When $\mu = \sum_{k = 1}^N A_k \delta_{-\tau_k}$ for some $N \in \mathbb N^\ast$, $\tau_1, \dotsc, \tau_N$ in $(0, r]$, and $A_1, \dotsc, A_N$ in $\mathcal M_d(\mathbb R)$, \eqref{eq:intro-diff-eqn-integral-measure} becomes the difference equation
\begin{equation}
\label{eq:intro-diff-eqn-sum}
x(t) = \sum_{k = 1}^N A_k x(t - \tau_k), \qquad t \in \mathbb R_+.
\end{equation}

The study of \eqref{eq:intro-general-difference} and its particular cases \eqref{eq:intro-diff-eqn-integral-measure}, \eqref{eq:intro-diff-eqn-integral-M}, and \eqref{eq:intro-diff-eqn-sum} has been motivated in the literature both from the fact that their properties, such as stability, allow one to obtain information on corresponding properties for neutral equations of the form \eqref{eq:intro-general-neutral} (see, e.g., \cite{Henry1974Linear, Hale2002Strong}, \cite[Chapter~9]{Hale1993Introduction}, and \cite[Section~1.4]{Mazanti2016Stability}, as well as Section~\ref{sec:stability-LTI-TDS} below), and also from the fact that some hyperbolic partial differential equations (PDEs\label{PDE}) on space dimension one can be put into this form by using the method of characteristics, a fact widely used in the literature (see, e.g., \cite{Auriol2019Explicit, Baratchart2021Sufficient, Bastin2016Stability, Cooke1968Differential, Slemrod1971Nonexistence, Miranker1961Periodic, Brayton1966Bifurcation, Brayton1967Nonlinear, Greenberg1984Effect, Kloss2012Flow, Fridman2010Bounds, Coron2015Dissipative, Gugat2010Stars, Coron2008Dissipative, Chitour2016Stability, Chitour2021One, Chitour2017Persistently, Balogoun2025Novel} and \cite[Section~1.4.3]{Mazanti2016Stability}). This connection allows properties of the PDE to be investigated via the corresponding delay system, and it has motivated many works on the difference equation \eqref{eq:intro-general-difference} and its particular cases along the years, such as \cite{Melvin1974Stability, Avellar1980Zeros, Cruz1970Stability, Henry1974Linear, Datko1977Linear, Hale1985Stability, Hale2002Strong, Michiels2009Strong, Chitour2016Stability, Mazanti2017Relative, Chitour2020Approximate, Chitour2023Hautus, Moreno1973Zeros, Hale2003Stability, Hale1975Parametric, Silkowski1976Star, Ngoc2015Exponential, Avellar1990Difference, Hale1967Class, Pandolfi1976Stabilization, Salamon1984Control, OConnor1983Function, OConnor1983Stabilization, Chyung1970Controllability, Pospisil2015Relative, Diblik2008Controllability}. While in many cases such a transformation yields difference equations with finitely many discrete delays only, it has been established in \cite{Auriol2019Explicit} that some hyperbolic PDEs with in-domain coupling terms can be transformed into difference equations with both discrete and distributed delays.

\section{Delay differential-algebraic equations}

\emph{Differential-algebraic equations} (DAEs\label{DAE}), also known as \emph{differential-algebraic systems}, are finite-dimensional dynamical systems in which part of the state is governed by an ordinary differential equation (ODE\label{ODE}), while the other part of the state follows an algebraic equation. In its most general form, a differential-algebraic system is a system that can be written under the form
\[
F(x^\prime(t), x(t), t) = 0
\]
for some function $F$ defined in an open subset of $\mathbb R^d \times \mathbb R^d \times \mathbb R$ and taking values in $\mathbb R^d$. When the parts of the state that are governed by the ODE and the algebraic equation do not change in time, a DAE is usually written under the form
\begin{equation}
\label{eq:general-DAE}
\left\{
\begin{aligned}
x^\prime(t) & = f(t, x(t), y(t)), \\
0 & = g(t, x(t), y(t)),
\end{aligned}
\right.
\end{equation}
where the state is now made of the pair of variables $(x, y)$, taking values in $\mathbb R^{d_x} \times \mathbb R^{d_y}$ for some nonnegative integers $d_x, d_y$, and $f$ and $g$ are defined in open subsets of $\mathbb R \times \mathbb R^{d_x} \times \mathbb R^{d_y}$ and take values, respectively, in $\mathbb R^{d_x}$ and $\mathbb R^{d_y}$. DAEs have been extensively studied in the delay-free setting (see, e.g., \cite{Brenan1996Numerical, Coleman1998Differential, Griepentrog1986Differential, Kumar1999Control, Kunkel2006Differential}) and they appear as natural models in several applications, such as in some electronic circuit models, in some control problems with constraints, or in the limiting behavior of singularly perturbed systems (see, e.g., \cite{Khalil1992Nonlinear, Kokotovic1999Singular, Avramovic1982Timescale, Kokotovic1968Singular, Kokotovic1975Singular, Kokotovic1976Singular, Kokotovic1972Singular, Naidu2002Singular, Saksena1984Singular, Friedrichs1946Singular, Prandtl1904Flussigkeitsbewegung, Tang2017Stability, Arias2023Frequency}).

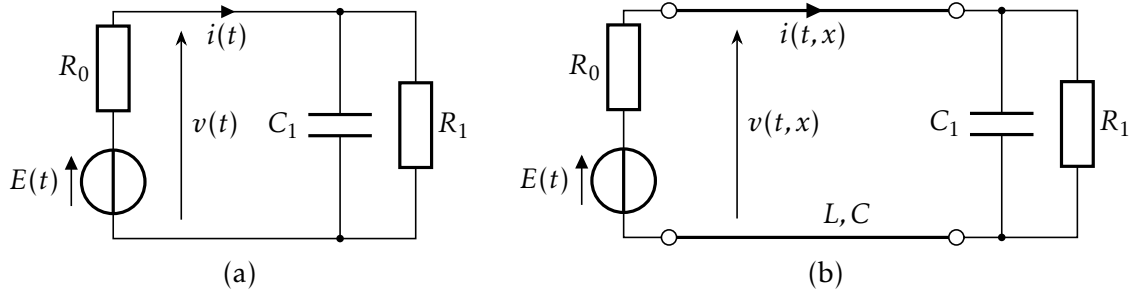
\begin{figure}[ht]
\centering
\begin{tabular}{@{} >{\centering} m{0.45\textwidth} @{} >{\centering} m{0.55\textwidth} @{}}
\begin{circuitikz}
\draw[semithick] (3, 0) -- (0, 0) to[V=$E(t)$] (0, 1.5) to[R=$R_0$] (0, 3) to[short, i_=$i(t)$] (3, 3);
\draw[semithick, -Stealth] (0.9, 0.25) -- node[midway, right] {$v(t)$} (0.9, 2.75);
\draw[semithick] (3, 0) to[C=$C_1$] (3, 3);
\draw[semithick] (3, 3) -- (4, 3) to[R=$R_1$] (4, 0) -- (3, 0);
\fill (3, 0) circle[radius=0.05] (3, 3) circle[radius=0.05];
\end{circuitikz} & \begin{circuitikz}
\draw[semithick] (0.6, 0) circle[radius=0.1] (0.5, 0) -- (0, 0) to[V=$E(t)$] (0, 1.5) to[R=$R_0$] (0, 3) -- (0.5, 3) (0.6, 3) circle[radius=0.1];
\draw[very thick] (0.7, 3) to[short, i_={$i(t, x)$}] (4.3, 3);
\draw[very thick] (0.7, 0) -- (4.3, 0);
\draw[semithick, -Stealth] (1.5, 0.25) -- node[midway, right] {$v(t, x)$} (1.5, 2.75);
\draw[semithick] (4.4, 0) circle[radius=0.1] (4.5, 0) -- (5, 0) to[C=$C_1$] (5, 3) -- (4.5, 3) (4.4, 3) circle[radius=0.1];
\draw[semithick] (5, 3) -- (6, 3) to[R=$R_1$] (6, 0) -- (5, 0);
\fill (5, 0) circle[radius=0.05] (5, 3) circle[radius=0.05];
\draw (2.5, 0) node[above right] {$L, C$};
\end{circuitikz} \tabularnewline
(a) & (b) \tabularnewline
\end{tabular}
\caption[Electric circuits modeled by differential-algebraic equations and delay dif\-fe\-ren\-tial-al\-ge\-braic equations]{(a) Electric circuit made of a voltage source $E(t)$ with internal resistance $R_0$ and a parallel connection between a resistor of resistance $R_1$ and a capacitor of capacitance $C_1$, and (b) the same electric circuit, but now with a transmission line of characteristic impedance and capacitance given by $L$ and $C$, respectively, between the voltage source and the parallel circuit made of the resistor and the capacitor}
\label{fig:circuits}
\end{figure}

A simple example of a differential-algebraic equation can be obtained when modeling the electric circuit from Figure~\ref{fig:circuits}(a), which is made of a voltage source $E(t)$ with internal resistance $R_0$ and a parallel connection between a resistor of resistance $R_1$ and a capacitor of capacitance $C_1$. Indeed, Kirchhoff's current law, combined with the equations of ideal resistors and capacitors, gives $i(t) = C_1 v^\prime(t) + \frac{v(t)}{R_1}$, while the analysis of the voltage source shows that $v(t) = E(t) - R_0 i(t)$. Hence, the pair $(v, i)$ satisfies the DAE
\[
\left\{
\begin{aligned}
v^\prime(t) & = -\frac{1}{R_1 C_1} v(t) + \frac{1}{C_1} i(t), \\
0 & = v(t) + R_0 i(t) - E(t).
\end{aligned}
\right.
\]

When propagation phenomena are introduced in a system modeled by a DAE, the propagation can sometimes be modeled simply by a time delay, and the system of equations one obtains when introducing delays in \eqref{eq:general-DAE} is called a \emph{delay differential-algebraic equation} (DDAE\label{DDAE}) or \emph{delay differential-algebraic system}. This is the case, for instance of the electric circuit represented in Figure~\ref{fig:circuits}(b), which represents the same circuit from Figure~\ref{fig:circuits}(a) but the direct connection between the voltage source and the parallel circuit made of $R_1$ and $C_1$ was replaced by a transmission line of characteristic impedance and capacitance given by $L$ and $C$, respectively. The voltage and the current at time $t$ across the transmission line are represented by the variables $v(t, x)$ and $i(t, x)$, respectively, where $x \in [0, 1]$ is the (normalized) position along the line. In this example, taken from \cite{Mazanti2021Effects}, the dynamics of the electric circuit can be modeled by a system of telegrapher's equations describing the evolution of $v$ and $i$ along the line, together with boundary conditions at $x = 0$ and $x = 1$ describing the relations between voltage and current on the voltage source ($x = 0$) and in the parallel circuit made of $R_1$ and $C_1$ ($x = 1$). More precisely, the circuit from Figure~\ref{fig:circuits}(b) can be modeled by
\begin{equation}
\label{ElectricCircuit}
\left\{
\begin{aligned}
& \partial_x v(t, x) + L \partial_t i(t, x) = 0, & \qquad & t \geq 0,\, x \in (0, 1), \\
& \partial_x i(t, x) + C \partial_t v(t, x) = 0, & & t \geq 0,\, x \in (0, 1), \\
& v(t, 0) = E(t) - R_0 i(t, 0), & & t \geq 0, \\
& i(t, 1) = C_1 \partial_t v(t, 1) + \tfrac{1}{R_1} v(t, 1), & & t \geq 0.
\end{aligned}
\right.
\end{equation}
Performing the classical change of variables into Riemann invariants
\[
u_1(t, x) = \frac{1}{2}\left[v(t, x) + \sqrt{\frac{L}{C}} i(t, x)\right], \qquad u_2(t, x) = \frac{1}{2}\left[v(t, x) - \sqrt{\frac{L}{C}} i(t, x)\right],
\]
and setting $y_1(t) = v(t, 1)$ and $y_2(t) = u_2(t, 1)$, one verifies that \eqref{ElectricCircuit} is rewritten as the delay differential-algebraic system
\begin{equation*}
\left\{
\begin{aligned}
y_1^\prime(t) & = -\frac{1}{C_1} \left(\frac{1}{R_1} + \sqrt{\frac{C}{L}}\right) y_1(t) - \frac{2}{C_1} \sqrt{\frac{C}{L}} \rho y_2(t - \tau) + \frac{2}{C_1} \sqrt{\frac{C}{L}} \frac{E(t - \tau/2)}{1 + R_0 \sqrt{\frac{L}{C}}}, \\
y_2(t) & = y_1(t) + \rho y_2(t - \tau) - \frac{E(t - \tau/2)}{1 + R_0 \sqrt{\frac{L}{C}}},
\end{aligned}
\right.
\end{equation*}
where $\rho = \frac{1 - R_0 \sqrt{\frac{L}{C}}}{1 + R_0 \sqrt{\frac{L}{C}}}$ is the reflection coefficient at the voltage source and $\tau = 2 \sqrt{L C}$ represents the round-trip travel time along the transmission line.

Another example of situation which is naturally modeled by a delay differential-algebraic system is the delayed output feedback control of linear finite-dimensional systems. Indeed, consider the linear control system
\begin{equation}
\label{LinControlSystem}
\left\{
\begin{aligned}
x^\prime(t) & = A x(t) + B u(t), \\
y(t) & = C x(t) + D u(t),
\end{aligned}
\right.
\end{equation}
where $x(t) \in \mathbb R^{d_x}$ is the state, $u(t) \in \mathbb R^{d_u}$ is the control, $y(t) \in \mathbb R^{d_y}$ is the output, $d_x$, $d_u$, and $d_y$ are positive integers, and $A$, $B$, $C$, and $D$ are matrices of appropriate dimensions. We consider the problem of stabilizing \eqref{LinControlSystem} by a delayed output feedback of the form $u(t) = K y(t - \tau)$, where $\tau > 0$ is the delay and $K$ is a matrix of appropriate dimension. The closed-loop system is then
\begin{equation*}
\left\{
\begin{aligned}
x^\prime(t) & = A x(t) + B K y(t - \tau), \\
y(t) & = C x(t) + D K y(t - \tau),
\end{aligned}
\right.
\end{equation*}
which is a delay differential-algebraic system with a single delay.

Just as in the above examples, many applications considered in the literature involving delay differential-algebraic systems involve in general either control systems with delays in the feedback loop or systems of partial differential equations representing propagation phenomena in one space dimension endowed with static and dynamic boundary conditions. We refer the interested reader to \cite{Brayton1968Small, Halanay1997Stability, Michiels2014Stability, Hale1993Introduction, Michiels2011Spectrum, Gumussoy2011Fixed, Rasvan1998Dynamical, Du2013Stability, Campbell2009Stability, Luzyanina2006Periodic, Shampine2006Delay, Zhu1997Asymptotic, Kosovalic2017Algebraic, Ha2012Analysis, Fridman2002Stability, Mattenet2022Improved, Zhang2023Coupled} for more details on the theory and applications of delay differential-algebraic systems, the subtleties in their analysis, and further examples of engineering systems which can be written as delay differential-algebraic systems.

\section{Controllability}

When systems of the form \eqref{eq:intro-general-retarded}, \eqref{eq:intro-general-neutral}, or \eqref{eq:intro-general-difference} also include a control input $u$, a natural structural question in control theory is to establish whether the system is \emph{controllable}: given an initial and a final condition for the system, is it possible to find a control $u$ steering the system from its initial condition to the desired final condition? There are many variants on this question, according to the functional state on which the system is considered, on whether one wishes to reach the desired target state in a given prescribed time $T$ or not, or on whether one wishes to reach exactly the final condition, which corresponds to the notion of exact controllability, or only an arbitrary neighborhood of it, corresponding to the notion of approximate controllability.

Many works have addressed controllability questions for retarded or neutral time-delay systems, most often in the linear case, in the approximate sense, in Hilbert spaces such as $L^2$ or $H^1$, and without a prescribed controllability time. Usual tools in this context are semigroups properties and Laplace transforms, which lead to Hautus-type conditions for controllability\footnote{Recall that, for a linear time-invariant finite-dimensional control system of the form $x^\prime(t) = A x(t) + B u(t)$, the Hautus criterion states that the system is controllable if and only if $\rank\begin{pmatrix} s \Id - A & B\end{pmatrix} = d$ for every $s \in \mathbb C$, or, equivalently, for every eigenvalue $s$ of $A$ (see, e.g., \cite[Lemma~3.3.7]{Sontag1998Mathematical}).}, i.e., rank conditions on certain matrix-valued holomorphic functions. Necessary and sufficient conditions for retarded systems can be found in \cite{Manitius1978Function}, with extensions to linear neutral systems with a single delay in \cite{OConnor1983Function, Jacobs1976Criteria}. A dual formulation of exact controllability, which can be seen as an observability inequality, was provided in \cite{Salamon1984Control}, which also presents some insights for approximate controllability results for linear neutral systems with several delays. Necessary and sufficient Hautus-type conditions for the approximate controllability of linear neutral systems were finally provided in \cite{Yamamoto1989Reachability}, making use of the realization theory developed in \cite{Yamamoto1981RealizationI, Yamamoto1981RealizationII}.

Controllability questions for difference equations are less present in the literature. In the linear time-invariant case with constant delays and no delays in the control input, a controlled difference equation can be written under the form
\begin{equation}
\label{eq:intro-linear-difference-equation}
 x(t)=\sum_{j=1}^N A_j x(t - \Lambda_j) + B u(t), \qquad t \in \mathbb R_+,
\end{equation}
where $N$ is the number of delays, $m$ is the number of control inputs, $x(t) \in \mathbb R^d$ is the instantaneous state at time $t$, $u(t) \in \mathbb R^m$ is the control at time $t$,  $A_1, \dotsc, A_N$ belong to $\mathcal M_d(\mathbb R)$, $B \in \mathcal M_{d, m}(\mathbb R)$, and $\Lambda_1, \dotsc, \Lambda_N$ are distinct positive delays, assumed, without loss of generality, to be ordered as $\Lambda_1 < \dotsb < \Lambda_N$.

One of the controllability notions that has been considered for \eqref{eq:intro-linear-difference-equation} in some references is that of relative controllability: instead of requiring that the final condition to be reached is a history function $x_T \in \mathcal C([-r, 0], \mathbb R^d)$, relative controllability studies the simpler question of determining whether, given an initial condition $x_0 \in \mathcal C([-r, 0], \mathbb R^d)$, a target instantaneous state $x_\ast$, and a positive time $T$, there exists a control $u$ such that the corresponding solution of \eqref{eq:intro-linear-difference-equation} satisfies $x(T) = x_\ast$. Relative controllability was originally introduced in the study of control systems with delays in the input \cite{Klamka1976Relative, Olbrot1972Controllability, Chyung1970Controllability, Balachandran1985Relative}, having been later extended to systems with delays in the state \cite{Diblik2008Controllability, Pospisil2015Relative, Mazanti2017Relative} and also to more general frameworks, such as stochastic control systems \cite{Karthikeyan2015Controllability} or fractional integro-differential systems \cite{Balachandran2016Relative}. The most general necessary and sufficient conditions for the relative controllability of \eqref{eq:intro-linear-difference-equation} can be found in \cite{Mazanti2017Relative}.

As for the exact and approximate controllability of \eqref{eq:intro-linear-difference-equation}, some partial results can be found in \cite{Chitour2020Approximate}, which considers controllability in the Hilbert space $L^2([-\Lambda_N,0],\mathbb{R}^d)$, providing a bound on the minimal time of controllability in the case where the delays are rationally dependent, a situation which is addressed by a standard state augmentation technique, reducing the matters at hand to the case of a single delay. That reference also addressed the first nontrivial case of two rationally independent delays when the system has dimension two and has a single scalar control, i.e., when $N = d = 2$ and $m = 1$, providing necessary and sufficient conditions for both approximate and exact controllability. However, the techniques of proof used in \cite{Chitour2020Approximate} become intractable when working in higher dimension or with more than two delays.

One of the main results of this monograph, presented in Chapter~\ref{chap:controllability} below and describing the results from \cite{Chitour2023Hautus}, provides necessary and sufficient conditions for the approximate controllability in time $T > 0$ of \eqref{eq:intro-linear-difference-equation} in $L^p([-\Lambda_N, 0], \mathbb R^d)$ for $p \in [1, +\infty)$ in terms of a Hautus-type criterion, showing in particular that approximate controllability in $L^p$ is independent of $p$. We will also provide in Chapter~\ref{chap:controllability} the other main results from \cite{Chitour2023Hautus}, showing in particular that approximate controllability in some time $T > 0$ is equivalent to approximate controllability in time $d \Lambda_N$, and providing necessary conditions for exact controllability in $L^1([-\Lambda_N, 0], \mathbb R^d)$.

\section{Realization theory}
\label{sec:Realization}

The approach used to address controllability of \eqref{eq:intro-linear-difference-equation} in Chapter~\ref{chap:controllability} is based on works on realization theory by Yutaka Yamamoto in \cite{Yamamoto1978Realization, Yamamoto1981RealizationI, Yamamoto1981RealizationII, Yamamoto1988Pseudo, Yamamoto1989Reachability, Yamamoto2007Coprimeness, Yamamoto2011Bezout, Yamamoto2008Behavioral, Yamamoto2016Behavioral, Yamamoto1999Coprimeness}, extending older ideas given in \cite{Kalman1972Realization, Kamen1976Module}.

Realization theory relies on the external description of control systems, in which a system is described as an input/output map. In this point of view of control systems, inputs are applied during a compact time interval of nonpositive times and outputs are observed at nonnegative times. The input space is a space of compactly supported functions defined in $\mathbb R_-$ and taking values in $\mathbb R^m$, such as $\Omega_p = \{u \in L^p(\mathbb R_-, \mathbb R^m) \suchthat \supp(u) \text{ is compact}\}$, the output space is a space of functions defined in $\mathbb R_+$ and taking values in $\mathbb R^\ell$, such as $L^p_{\loc}(\mathbb R_+, \mathbb R^\ell)$, and a control system is simply a map $f\colon \Omega_p \to L^p_{\loc}(\mathbb R_+, \mathbb R^\ell)$.

The original question in realization theory was, given a control system in its external description as a map $f\colon \Omega_p \to L^p_{\loc}(\mathbb R_+, \mathbb R^\ell)$ and assuming that $f$ is linear and time-invariant (in the sense that time shifts in the input yield outputs shifted in time by the same amount), determine whether the system can also be represented as a linear time-invariant finite-dimensional control system under state-space form, i.e., under the form
\[
\left\{
\begin{aligned}
x^\prime(t) & = A x(t) + B u(t), \\
y(t) & = C x(t) + D u(t),
\end{aligned}
\right.
\]
with input $u$, state $x$, and output $y$ taking values in $\mathbb R^m$, $\mathbb R^d$, and $\mathbb R^\ell$, respectively, and $A$, $B$, $C$, $D$ being constant matrices of appropriate dimensions; and, when this is possible, construct $(A, B, C, D)$ with a minimal number of states $d$ from information from the external representation $f$, such as its impulse response, when the input space allows for the representation of a Dirac mass at time $0$. This question was solved by Elmer G.~Gilbert and Rudolf E.~Kálmán in the works \cite{Gilbert1963Controllability, Kalman1963Mathematical}, and we refer the interested reader to \cite{Antoulas1991Linear} for more historical background on realization theory.

Y.~Yamamoto's approach to realization theory in \cite{Yamamoto1978Realization, Yamamoto1981RealizationI, Yamamoto1981RealizationII, Yamamoto1988Pseudo, Yamamoto1989Reachability, Yamamoto2007Coprimeness, Yamamoto2011Bezout, Yamamoto2008Behavioral, Yamamoto2016Behavioral, Yamamoto1999Coprimeness} is focused on input/output maps that may represent infinite-dimensional linear systems, such as time-delay systems. Those works focus mainly on linear time-invariant input/output systems represented as a convolution, i.e., the output $y \in L^p(\mathbb R_+, \mathbb R^\ell)$ is computed from the input $u \in \Omega_p$ as
\begin{equation}
\label{eq:intro-input-output}
y(t) = \int_0^{+\infty} \diff A(s) u(t - s), \qquad t \in \mathbb R_+,
\end{equation}
where $u$ is extended by $0$ to $\mathbb R$ and the impulse response $A$ is a $\ell \times m$ matrix-valued Radon measure in $\mathbb R_+$.

Most of the previously cited Yamamoto's works on realization theory consider the important particular case of pseudo-rational systems. The impulse response $A$ and the corresponding input/output system described by \eqref{eq:intro-input-output} are said to be \emph{pseudo-rational} (see \cite[Definition~2.5]{Yamamoto1988Pseudo}) if there exist $\ell \times \ell$ and $\ell \times m$ matrix-valued distributions $Q$ and $P$ defined in $\mathbb R$, with compact support in $\mathbb R_-$, and such that
\begin{enumerate}
\item $Q$ is invertible in the convolution algebra of $\ell \times \ell$ matrix-valued distributions with support bounded on the left, i.e., there exists a $\ell \times \ell$ matrix-valued distribution $Q^{-1}$ defined in $\mathbb R$ with support bounded on the left and such that $Q^{-1} * Q = \delta_0 \Id_\ell$;
\item $Q^{-1}$ is supported in $\mathbb R_+$;
\item $A$ is the truncation to $\mathbb R_+$ of $Q^{-1} * P$;
\item the distribution $\det Q$ is invertible in the convolution algebra of real-valued distributions with support bounded on the left; and
\item the order of the distribution $(\det Q)^{-1}$ is the opposite of the order of $\det Q$.
\end{enumerate}
A major motivation for considering pseudo-rational systems is that, as remarked in \cite[Example~2.7]{Yamamoto1988Pseudo}, retarded and neutral time-delay systems with finitely many discrete delays define pseudo-rational systems, and this fact readily extends to difference equations of the form \eqref{eq:intro-linear-difference-equation}.

The main results obtained in the series of works \cite{Yamamoto1978Realization, Yamamoto1981RealizationI, Yamamoto1981RealizationII, Yamamoto1988Pseudo, Yamamoto1989Reachability, Yamamoto2007Coprimeness, Yamamoto2011Bezout, Yamamoto2008Behavioral, Yamamoto2016Behavioral, Yamamoto1999Coprimeness} concerning controllability of pseudo-rational input/output systems are characterizations of some forms of controllability in terms of the solvability of a Bézout's identity. More precisely, given a pseudo-rational input/output system of the form \eqref{eq:intro-input-output} and $Q$ and $P$ as in the definition of pseudo-rational, \eqref{eq:intro-input-output} is exactly controllable from the origin in the space of distributions if and only if there exist $\ell \times \ell$ and $m \times \ell$ matrix-valued distributions $R$ and $S$ such that the Bézout's identity $Q * R + P * S = \delta_0 \Id_\ell$ is satisfied \cite[Theorem~3.4]{Yamamoto2008Behavioral}, while $L^2$ approximate controllability from the origin is equivalent to proving an approximate Bézout's identity, i.e., the existence of sequences of $\ell \times \ell$ and $m \times \ell$ matrix-valued distributions $(R_n)_{n \in \mathbb N}$ and $(S_n)_{n \in \mathbb N}$ such that $Q * R_n + P * S_n$ converges to $\delta_0 \Id_\ell$ in the distributional sense as $n$ tends to infinity \cite[Theorem~4.4]{Yamamoto1988Pseudo}. In realization theory, the solvability of a Bézout's identity (respectively, an approximate Bézout's identity) is usually called \emph{left coprimeness} (respectively, \emph{approximate left coprimeness}) (see, e.g., \cite{Polderman1998Introduction, Sontag1998Mathematical}). Note that one of the advantages of relating controllability issues and Bézout's identities is that, if the latter is solved, one can find, for a given target output, an explicit input which gives the desired output, as we describe below in Proposition~\ref{prop_fond_exact_controllability}.

\section{Stability analysis of linear time-invariant time-delay systems}
\label{sec:stability-LTI-TDS}

In addition to controllability, two other important questions on time-delay systems that have been extensively considered in the literature, and that we will address in this monograph, are their stability analysis and stabilization, which still attract much research effort and are active research fields \cite{Gu2003Stability, Michiels2014Stability, Sipahi2011Stability, Cruz1970Stability, Avellar1980Zeros, Datko1977Linear, Hale2003Stability, Henry1974Linear, Rodrigues1980Growth}. The main goal of most works on the topic is to provide conditions on a time delay-system under which one can guarantee that a given equilibrium point is asymptotically or exponentially stable. While, in stability analysis, one usually works with systems under the form \eqref{eq:intro-general-retarded}, \eqref{eq:intro-general-neutral}, or \eqref{eq:intro-general-difference}, in stabilization problems these equations also present a control input $u$, and the goal is typically to choose the control $u(t)$ at time $t$ in feedback form, i.e., as a function of the instantaneous state $x(t)$ and its history function $x_t$, in order to render asymptotically or exponentially stable some equilibrium point. As this monograph focuses mostly on linear time-delay systems with finitely many discrete delays, we provide a brief overview of some classical results on stability for this class of systems, based mainly on \cite[Chapters~7 and 9]{Hale1993Introduction} and \cite[Chapter~1]{Michiels2014Stability}.

Let us first consider the case of the neutral time-delay system
\begin{equation}
\label{eq:neutral-linear}
x^\prime(t) + \sum_{j=1}^N H_j x^\prime(t - r_j) = A_0 x(t) + \sum_{k=1}^M A_k x(t - \tau_k),
\end{equation}
where $r_1, \dotsc, r_N, \tau_1, \dotsc, \tau_M$ are positive real numbers and $H_1, \dotsc, H_N, A_0, \dotsc, A_M$ are matrices in $\mathcal M_d(\mathbb R)$. Similarly to the case of linear time-invariant finite-dimensional systems, the stability of the trivial solution of \eqref{eq:neutral-linear} can be studied through the zeros of its characteristic function. Indeed, let $\Delta\colon \mathbb C \to \mathbb C$ be the holomorphic function defined by
\begin{equation}
\label{eq:intro-Delta}
\Delta(s) = \det\left(s\left(\Id_d + \sum_{j=1}^N H_j e^{-r_j s}\right) - A_0 - \sum_{k=1}^M A_k e^{-\tau_k s}\right).
\end{equation}
Then the trivial solution of \eqref{eq:neutral-linear} is exponentially stable if and only if there exists $\alpha > 0$ such that all roots $s \in \mathbb C$ of $\Delta$ satisfy $\Real s \leq -\alpha$. The function $\Delta$ from \eqref{eq:intro-Delta} is a particular case of a \emph{quasipolynomial}, a class of holomorphic functions on which we will provide further properties in Section~\ref{sec:quasipolynomials} below. Roots of $\Delta$ are usually called \emph{characteristic roots} or \emph{spectral values} of \eqref{eq:neutral-linear}, and the set of all these roots is the \emph{spectrum} of \eqref{eq:neutral-linear}.

The \emph{spectral abscissa} of \eqref{eq:neutral-linear} is the number $\gamma$ defined by
\begin{equation}
\label{eq:spectral-abscissa}
\gamma = \sup\{\Real s \suchthat \Delta (s) = 0\},
\end{equation}
and $\gamma$ is necessarily finite, i.e., $\Delta$ admits at least one root and the set $\{\Real s \suchthat \Delta (s) = 0\}$ is upper bounded. Exponential stability of \eqref{eq:neutral-linear} can be rephrased in terms of its spectral abscissa by saying that the trivial solution of \eqref{eq:neutral-linear} is exponentially stable if and only if $\gamma < 0$. More precisely, it is shown in \cite[Chapter~9, Corollary~3.1]{Hale1993Introduction} that, for every $\varepsilon > 0$, there exists $C > 0$ such that, for every solution $x$ of \eqref{eq:neutral-linear}, one has $\abs{x(t)} \leq C e^{(\gamma + \varepsilon) t} \max_{\theta \in [-r, 0]} \abs{x(\theta)}$, where $r = \max\{r_1, \dotsc, r_N, \tau_1, \dotsc, \tau_M\}$.

We say that a root $s_0$ of $\Delta$ is \emph{dominant} if it attains the spectral abscissa, i.e., $\Real s_0 = \gamma$, and we say that it is \emph{strictly dominant} if in addition $\Real s < \Real s_0$ for every root $s$ of $\Delta$ with $s \neq s_0$ and $s \neq \overline s_0$. In particular, a sufficient condition for exponential stability of \eqref{eq:neutral-linear} is the existence of a dominant root with negative real part.

In the retarded case, i.e., when $H_1 = \dotsb = H_N = 0$ in \eqref{eq:neutral-linear}, for every $\alpha \in \mathbb R$, the set $\{s \in \mathbb C \suchthat \Delta(s) = 0 \text{ and } \Real s \geq \alpha\}$ is finite. In other words, if $(s_n)_{n \in \mathbb N}$ is a sequence of distinct spectral values of \eqref{eq:neutral-linear}, then necessarily $\Real s_n \to -\infty$ as $n \to +\infty$. In addition, the spectral abscissa $\gamma$ is a continuous function of the matrices $A_0, \dotsc, A_M$ and the delays $\tau_1, \dotsc, \tau_M$. These properties do not carry over to systems of neutral type: if the delays $r_1, \dotsc, r_N$ are all distinct and there exists $j \in \llbracket 1, N\rrbracket$ such that $H_j \neq 0$, then there exist $\zeta \in \mathbb R$ and a sequence $(s_n)_{n \in \mathbb N}$ of spectral values of \eqref{eq:neutral-linear} such that $\Real s_n \to \zeta$ and $\Imag s_n \to +\infty$ as $n \to +\infty$. In addition, $\gamma$ may fail to be continuous (see, e.g., \cite[Example~1.41]{Michiels2014Stability}).

Similar properties also hold for difference equations of the form
\begin{equation}
\label{eq:difference-linear}
x(t) = - \sum_{j=1}^N H_j x(t - r_j),
\end{equation}
in which case the corresponding characteristic function is
\begin{equation}
\label{eq:intro-Delta-difference}
\Delta_0(s) = \det\left(\Id_d + \sum_{j=1}^N H_j e^{-r_j s}\right).
\end{equation}
Once again, the trivial solution of \eqref{eq:difference-linear} is exponentially stable if and only if there exists $\alpha > 0$ such that all roots $s \in \mathbb C$ of $\Delta_0$ satisfy $\Real s \leq -\alpha$, and the spectral abscissa of \eqref{eq:difference-linear} is defined as in \eqref{eq:spectral-abscissa}, with $\Delta$ replaced by $\Delta_0$. One can easily construct examples illustrating that the spectral abscissa of \eqref{eq:difference-linear} can be discontinuous in the delays $r_1, \dotsc, r_N$ (see, e.g., \cite[Section~9.6]{Hale1993Introduction}, \cite[Example~1.37]{Mazanti2016Stability}, \cite[Example~1.36]{Michiels2014Stability}). In addition, the stability of \eqref{eq:neutral-linear} depends on that of \eqref{eq:difference-linear}: a necessary condition for the exponential stability of the trivial solution of \eqref{eq:neutral-linear} is the exponential stability of the trivial solution of \eqref{eq:difference-linear}.

Characterizing the stability of \eqref{eq:neutral-linear} or \eqref{eq:difference-linear} amounts to finding the location of the roots of the complex functions in \eqref{eq:intro-Delta} and \eqref{eq:intro-Delta-difference} and, for that purpose, one may use classical techniques in complex analysis, such as Cauchy's argument principle. This approach was adopted in \cite{Hassard1997Counting, Stepan1989Retarded}, providing results allowing one to count the number of roots of \eqref{eq:intro-Delta} with positive real part in the retarded case, i.e., $H_1 = \dotsb = H_N = 0$.

Let us finally consider the case of linear delay differential-algebraic equations with fi\-ni\-tely many discrete delays of the form\footnote{Despite its simple appearance and its resemblance to retarded time-delay systems, the DDAE \eqref{eq:linear-DDAE} contains a large class of equations as particular cases, including neutral time-delay systems. For instance, the equation
\[
E x^\prime(t) + \sum_{j=1}^N H_j x^\prime(t - r_j) = A_0 x(t) + \sum_{k=1}^M A_k x(t - \tau_k)
\]
(which reduces to the neutral time-delay system \eqref{eq:neutral-linear} when $E = \Id$) can be rewritten under the form \eqref{eq:linear-DDAE} by setting $y(t) = E x(t) + \sum_{j=1}^N H_j x(t - r_j)$ and rewriting the system as
\[
\left\{
\begin{aligned}
y^\prime(t) & = A_0 x(t) + \sum_{k=1}^M A_k x(t - \tau_k),\\
0 & = E x(t) + \sum_{j=1}^N H_j x(t - r_j) - y(t).
\end{aligned}
\right.
\]
}
\begin{equation}
\label{eq:linear-DDAE}
E x^\prime(t) = A_0 x(t) + \sum_{k=1}^M A_k x(t - \tau_k),
\end{equation}
where $\tau_1, \dotsc, \tau_M$ are positive real numbers and $E, A_0, \dotsc, A_M$ are matrices in $\mathcal M_d(\mathbb R)$. When $E$ is invertible, the above system can be rewritten, after multiplication by $E^{-1}$, as a retarded time-delay system, but the interesting case for delay differential-algebraic equations is when $E$ is singular. In many applications, $E$ will be a diagonal matrix, with diagonal entries in $\{0, 1\}$: a diagonal entry equal to $1$ represents a delay-differential equation, while a diagonal entry equal to $0$ corresponds to the ``algebraic'' part of the DDAE.

The characteristic function associated with \eqref{eq:linear-DDAE} is the holomorphic function $\Delta\colon \mathbb C \to \mathbb C$ defined by
\begin{equation}
\label{eq:intro-Delta-DDAE}
\Delta(s) = \det\left(s E - A_0 - \sum_{k=1}^M A_k e^{-\tau_k s}\right).
\end{equation}
Let $\nu \in \llbracket 0, d\rrbracket$ denote the dimension of the kernel of $E$ and $U, V \in \mathcal M_{d, \nu}(\mathbb R)$ be such that $E V = 0$ and $U^\transp E = 0$. Under the assumption that the matrix $U^\transp A_0 V \in \mathcal M_{\nu}(\mathbb R)$ is invertible, the roots of $\Delta$ characterize the asymptotic behavior of solutions of \eqref{eq:linear-DDAE} (see, e.g., \cite[Proposition~1]{Michiels2011Spectrum}): letting $\gamma$ be spectral abscissa of \eqref{eq:intro-Delta-DDAE}, defined as before by $\gamma = \sup\{\Real s \suchthat \Delta (s) = 0\}$, we have that the trivial solution of \eqref{eq:neutral-linear} is exponentially stable if and only if $\gamma < 0$, or, equivalently, if and only if there exists $\alpha > 0$ such that all roots $s \in \mathbb C$ of $\Delta$ satisfy $\Real s \leq -\alpha$.

The above discussion shows that, when one wishes to study the exponential stability of \eqref{eq:neutral-linear}, \eqref{eq:difference-linear}, or \eqref{eq:linear-DDAE}, one may proceed by studying their respective characteristic functions \eqref{eq:intro-Delta}, \eqref{eq:intro-Delta-difference}, or \eqref{eq:intro-Delta-DDAE}, which are all holomorphic functions. Results from complex analysis thus play an important role in the stability analysis of linear time-invariant time-delay systems. A particular result from complex analysis that shall be of much use in this monograph is Hurwitz's theorem, whose statement we now recall (see, e.g., \cite[Chapter~VII, Theorem~2.5]{Conway1978Functions} or \cite[p.~231]{Gamelin2001Complex}).

\begin{theorem}[Hurwitz's theorem]
\label{thm:Hurwitz}
Let $G \subset \mathbb C$ be open and connected and $(f_n)_{n \in \mathbb N}$ be a sequence of analytic functions defined in $G$ and converging to a nontrival analytic function $f\colon G \to \mathbb C$ uniformly on every compact subset of $G$. Then, for every open ball $B \subset G$ so that $\overline B \subset G$, if $f$ has no zeros at $\partial B$, then $f_n$ and $f$ have the same number of zeros in $B$ (counted with their multiplicities) for $n$ large enough.
\end{theorem}

\section{Stabilization of linear time-invariant time-delay systems}
\label{sec:stabilization-LTI-TDS}

As for the stabilizability problem, let us consider \eqref{eq:neutral-linear} with an additional control input $u$, i.e.,
\begin{equation}
\label{eq:neutral-control}
x^\prime(t) + \sum_{j=1}^N H_j x^\prime(t - r_j) = A_0 x(t) + \sum_{k=1}^M A_k x(t - \tau_k) + B u(t),
\end{equation}
where $u$ takes values in $\mathbb R^m$ and $B \in \mathcal M_{d, m}(\mathbb R)$.

One of the difficulties of the stabilization of time-delay systems is that, contrarily to the case of finite-dimensional systems in which the characteristic function is a polynomial, the quasipolynomial \eqref{eq:intro-Delta} has infinitely many roots. While, for a finite-dimensional linear system of the form $x^\prime(t) = A x(t) + B u(t)$, the linear feedback $u(t) = K x(t)$ allows one to assign all roots of the characteristic function of the resulting closed-loop system $x^\prime(t) = (A + B K) x(t)$ under the assumption that the pair $(A, B)$ is controllable, a similar result is not possible in general for \eqref{eq:neutral-control}, since a linear time-invariant feedback law of the form
\begin{equation}
\label{eq:feedback}
u(t) = K_0 x(t) + \sum_{j=1}^N K_{1, j} x^\prime(t - r_j) + \sum_{k=1}^M K_{2, k} x(t - \tau_k)
\end{equation}
contains only finitely many parameters, but the resulting characteristic function has infinitely many roots. In addition, there exist criteria to obtain information on the location of the roots of a polynomial in terms of its coefficients, such as Routh--Hurwitz criterion \cite{Barnett1983Polynomials}, but similar simple criteria based only on coefficients and delays are not available for \eqref{eq:intro-Delta} in general, except for some particular configurations (see, e.g., \cite{Hayes1950Roots, Cebotarev1949Routh}), and relating the coefficients and delays of \eqref{eq:neutral-linear} with its stability properties still a question of ongoing interest \cite{Gu2003Stability, Ma2019Delay}.

Some works have tried to implement stabilization strategies based on assigning roots of the characteristic function of \eqref{eq:neutral-control}--\eqref{eq:feedback} by selecting the coefficients of the matrices $K_0, K_{1, 1}, \dotsc, K_{1, N}, K_{2, 1}, \dotsc, K_{2, M}$, but assigning only finitely many spectral values \cite{Ackermann1972Entwurf, Vyhlidal2009Quasi, Ram2011Partial}, a technique sometimes called \emph{partial pole placement}. However, it is difficult in general to control the behavior of the other spectral values, and, even if one assigns spectral values with negative real part only, there might be other spectral values with nonnegative real part, preventing exponential stability. This motivates research on methods for assigning finitely many spectral values through a suitable choice of the parameters of the control law \eqref{eq:feedback} while also guaranteeing that the spectral abscissa of \eqref{eq:neutral-control}--\eqref{eq:feedback} is attained at one of the assigned values, guaranteeing thus that the real parts of all other spectral values are bounded by the largest real part of the assigned spectral values. Chapter~\ref{chap:MID} below presents one possible method for a subclass of systems of the form \eqref{eq:neutral-control}--\eqref{eq:feedback}, based on the assignment of spectral values with the maximal possible multiplicity.

Another approach for the stabilization of \eqref{eq:neutral-control} with a control law of the form \eqref{eq:feedback} was presented in \cite{Hale2002Strong}, which considers not only stabilizability of the system for a given set of delays $r_1, \dotsc, r_N, \tau_1, \dotsc, \tau_M$, but also robustness of stability with respect to these delays, in the sense of \ref{item:HS-local} from Proposition~\ref{prop:Hale-Silkowski}. The main results of \cite{Hale2002Strong} provide necessary and sufficient conditions, under the form of a Hautus-type test, for the existence of a stabilizing feedback law of a form slightly more general than \eqref{eq:feedback} which is also robust with respect to perturbations of the delays, treating also the case of stabilizability of difference equations.

Other works have considered the stabilization of \eqref{eq:neutral-control} with control laws more general than \eqref{eq:feedback}, depending on the whole history function $x_t$. An interesting idea in this sense is the finite spectrum assignment, also called finite pole placement, introduced in the late 1970s in \cite{Olbrot1978Stabilizability, Manitius1979Finite}. In this strategy, $u(t)$ is chosen as a continuous linear function of the history function $x_t$ in such a way as to guarantee that delays are compensated and the resulting system, after inserting $u$ in \eqref{eq:neutral-control}, is finite-dimensional. In this case, the feedback law can typically be interpreted as containing a prediction of future states of the system, which is able to compensate the delay and hence render the full system finite-dimensional. This strategy was carried out \cite{Manitius1979Finite} in the case of retarded systems with commensurable delays (i.e., the delays are all integer multiples of a same constant), as well as for the case of systems without delays in the state $x$ but with a delay in the input $u$, and other works have further explored this idea, such as \cite{Brethe1998Effective}, which presents a necessary and sufficient condition for finite spectrum assignment to be possible and provides a constructive procedure for constructing a suitable feedback control law. We refer to \cite{Wang1999Finite} for a presentation of the main ideas of this method as well as some comparisons with other methods devoted to the control of systems with input delays. Note, however, that the practical implementation of finite spectrum assignment may present many issues, since the system may remain unstable if the controlled feedback law is not implemented exactly, and even small numerical errors coming from the discretization of integrals appearing in the feedback control law may destabilize the system, as illustrated in \cite{Engelborghs2001Limitations}. This comes from the fact that finite spectrum assignment guarantees a finite spectrum for a certain control law, but approximations of this law will still yield an infinite spectrum, which may contain spectral values with positive real part even for arbitrarily small discretization parameters, a phenomenon known as spillover.

A numerical strategy for selecting parameters of a control law of the form \eqref{eq:feedback} in order to stabilize \eqref{eq:neutral-control} in the retarded case, known as the continuous pole placement, was proposed in \cite{Michiels2002Continuous}. This idea is based on the continuous dependence of spectral values on the parameters of \eqref{eq:feedback}, and consists, roughly speaking, in ``shifting'' to the left spectral values with positive real part in a ``quasi-continuous'' way, subject to the strong constraint that, during this shifting action, spectral values with negative real part remain in the left half-plane. We refer to \cite{Michiels2014Stability} and references therein for further insights on the number of controlled spectral values, as well as the interpretation of continuous pole placement as a local strategy to solve an appropriate optimization problem where the objective function, the spectral abscissa, is not differentiable. These ideas were extended to neutral systems in \cite{Michiels2005Eigenvalue}.

In addition to the methods described above, many other stabilizability methods have been proposed in the literature, such as the $D$-partition method from \cite{Neimark1948Structure, Neimark1998D, Neimark1949D}, the am\-pli\-tude-phase method from \cite{Cypkin1946Stability}, the $\tau$-partition method from \cite{Lee1969Tau}, the frequency-sweeping method \cite{Chen1995Frequency, Chen1995Computing}, as well as many other important ideas \cite{Walton1987Direct, Cooke1986Zeroes, Thowsen1981Analytic, Olgac2002Exact}. We refer the interested reader to \cite{Elsgolts1973Introduction, Pinney1958Ordinary, Stepan1989Retarded, Michiels2014Stability, Gu2003Stability, Sipahi2011Stability, Niculescu2021Stability, Benarab2022Over} for more details.

\section{Sensitivity with respect to perturbations}
\label{sec:sensitivity}

In many practical applications, systems are subject to uncertainties arising from modeling inaccuracies or physical limitations in measurement and implementation. In this context, small variations in the parameters can significantly affect the system dynamics. Therefore, it is essential to investigate how the behavior of the system changes under such variations, in order to ensure robustness, stability, and satisfactory performance even under nonideal conditions.

In the context of the present monograph, we will be interested on the sensitivity of the exponential stability property of a time-delay system. As discussed in Section~\ref{sec:stability-LTI-TDS}, exponential stability of linear time-invariant time-delay systems with finitely many discrete delays of the forms \eqref{eq:neutral-linear}, \eqref{eq:difference-linear} (or, equivalently, \eqref{eq:intro-diff-eqn-sum}), and \eqref{eq:linear-DDAE} can be studied through the roots of their respective characteristic functions, and more precisely through their respective spectral abscissas.

For linear difference equations of the form \eqref{eq:intro-diff-eqn-sum}, these sensitivity issues have been systematically addressed in the literature \cite{Hale1993Introduction, Avellar1980Zeros, Henry1974Linear, Hale2002Strong, Hale2003Stability, Michiels2009Strong, Chitour2016Stability, Melvin1974Stability, Cruz1970Stability, Silkowski1976Star, Avellar1990Difference, Moreno1973Zeros}. Variations of the parameters of \eqref{eq:intro-diff-eqn-sum} correspond to variations in the matrices $A_1, \dotsc, A_N$ and in the delays $\tau_1, \dotsc, \tau_N$. It can be shown that, if \eqref{eq:intro-diff-eqn-sum} is exponentially stable, then, for small changes in the matrices, the system remains exponentially stable \cite[Section~9.6]{Hale1993Introduction}. However, as briefly remarked in Section~\ref{sec:stability-LTI-TDS}, the situation is different with respect to perturbation in delays, as illustrated by the following example, adapted from \cite[Section~9.6]{Hale1993Introduction} (see also \cite{Silkowski1976Star, Avellar1980Zeros}).

\begin{example}
Let $a$ and $b$ be real numbers with $b + 1 < a < 2$ and $-1 < b < -\frac{a^2}{4}$ and consider the scalar difference equation
\begin{equation}
\label{eq:diff-eqn-perturbation}
x(t) = a x(t - 1) + b x(t - 2), \qquad t \in \mathbb R_+.
\end{equation}
According to \eqref{eq:intro-Delta-difference}, its characteristic function is the function $\Delta_0 \colon \mathbb C \to \mathbb C$ given, for $s \in \mathbb C$, by
\begin{equation}
\label{eq:perturbation-Delta0}
\Delta_0(s) = 1 - a e^{-s} - b e^{-2s}.
\end{equation}
In particular, $s \in \mathbb C$ is a root of $\Delta_0$ if and only if the quantity $\rho = e^s$ is a root of the polynomial $\rho^2 - a \rho - b = 0$. One can easily compute that, under our standing assumptions on $a$ and $b$, the roots of this polynomial are the pair of complex conjugate numbers $\frac{a \pm \sqrt{4 \abs{b} - a^2} i}{2}$, and, in particular, any root $\rho$ of the polynomial satisfies $\abs{\rho} = \sqrt{\abs{b}}$, which shows that any root $s$ of $\Delta_0$ satisfies $\Real s = \frac{1}{2} \ln \abs{b} < 0$. Thus, the spectral abscissa of \eqref{eq:diff-eqn-perturbation} is the negative number $\frac{1}{2} \ln \abs{b}$, and hence \eqref{eq:diff-eqn-perturbation} is exponentially stable.

Let us now consider the perturbation of \eqref{eq:diff-eqn-perturbation} given by
\begin{equation}
\label{eq:diff-eqn-perturbation-perturbed}
x(t) = a x\left(t - \frac{4k}{4k + 1}\right) + b x(t - 2), \qquad t \in \mathbb R_+,
\end{equation}
where $k \in \mathbb N^\ast$. Its characteristic function is the function $\Delta \colon \mathbb C \to \mathbb C$ given, for $s \in \mathbb C$, by
\begin{equation}
\label{eq:perturbation-Delta}
\Delta(s) = 1 - a e^{-\frac{4k}{4 k + 1} s} - b e^{-2s}.
\end{equation}
For every $\alpha \in \mathbb R$, we have
\[
\Delta\left(\alpha + \left(2 k + \frac{1}{2}\right) \pi i\right) = 1 - a e^{-\frac{4 k}{4 k + 1} \alpha} + b e^{-2 \alpha},
\]
and thus, in particular, $\Delta\left(\alpha + \left(2 k + \frac{1}{2}\right) \pi i\right) \in \mathbb R$ for every $\alpha \in \mathbb R$ and
\[
\Delta\left(\left(2 k + \frac{1}{2}\right) \pi i\right) = 1 - a + b < 0, \qquad \lim_{\alpha \to +\infty}\Delta\left(\alpha + \left(2 k + \frac{1}{2}\right) \pi i\right) = 1.
\]
By the intermediate value theorem, there exists $\alpha_k \in \mathbb R_+^\ast$ such that $\Delta\left(\alpha_k + \left(2 k + \frac{1}{2}\right) \pi i\right) = 0$. In particular, $\Delta$ admits a root with positive real part, showing that \eqref{eq:diff-eqn-perturbation-perturbed} is not exponentially stable. More precisely, \eqref{eq:diff-eqn-perturbation-perturbed} admits the unbounded solution $x$ given, for $t \in [-2, +\infty)$, by $x(t) = e^{\alpha_k t} \sin\left(\left(2 k + \frac{1}{2}\right) \pi t\right)$.

\begin{figure}[ht]
\centering
\resizebox{0.6\textwidth}{!}{\input{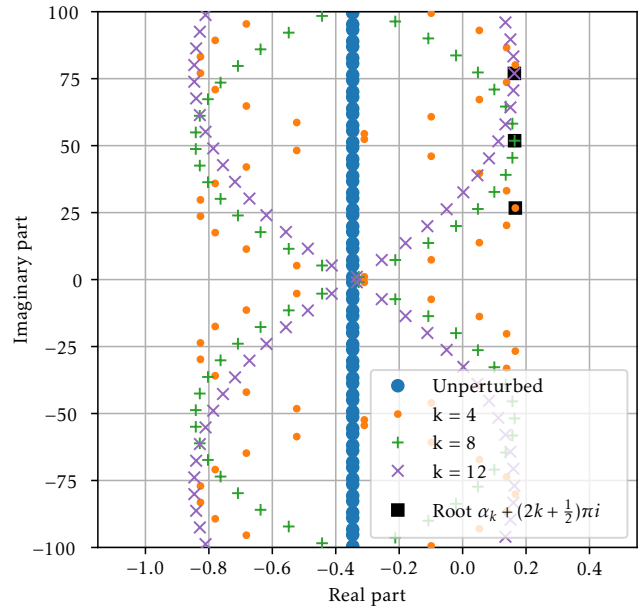}}
\caption{Roots of $\Delta_0$ from \eqref{eq:perturbation-Delta0} and of $\Delta$ from \eqref{eq:perturbation-Delta} for $k \in \{4, 8, 12\}$ and with $a = \frac{3}{4}$, $b = -\frac{1}{2}$}
\label{fig:spectrum-perturbation}
\end{figure}

Figure~\ref{fig:spectrum-perturbation} shows\footnote{Roots of $\Delta_0$ were computed analytically using their explicit expression, while those from $\Delta$ for $k \in \{4, 8, 12\}$ were computed in the rectangle $\{s \in \mathbb C \suchthat \abs{\Real s} \leq 5,\, \abs{\Imag s} \leq 100\}$ using the QPmR algorithm from \cite{Vyhlidal2014QPmR} and, more precisely, its Python implementation available at \url{https://github.com/DSevenT/QPmR}.} the roots of the characteristic equation $\Delta_0$ of the unperturbed system \eqref{eq:diff-eqn-perturbation} (blue circles) in the case $a = \frac{3}{4}$ and $b = -\frac{1}{2}$, as well as the roots of $\Delta$ for $k = 4$, $k = 8$, and $k = 12$ (orange points, green pluses, and purple crosses, respectively), the root $\alpha_k + \left(2 k + \frac{1}{2}\right) \pi i$ being marked by a black square for each such $k$. One observes that, as expected, the roots of $\Delta_0$ all lie in the same vertical axis, but $\Delta$ always admits roots with positive real part.
\end{example}

The discontinuity of the spectral abscissa for neutral systems and for difference equations means that a given system can be exponentially stable, but lose its exponential stability under arbitrarily small perturbations on the delays. For this reason, some works have considered more robust stability properties, and one of the main results in this sense for difference equations is the following one (see, e.g., \cite[Theorem~5.2]{Avellar1980Zeros} and \cite[Chapter~9, Theorem~6.1]{Hale1993Introduction}, \cite{Silkowski1976Star}), known as the Hale--Silkowski criterion, after \cite{Hale1975Parametric, Silkowski1976Star}.

\begin{proposition}
\label{prop:Hale-Silkowski}
Let $A = (A_1, \dotsc, A_N) \in \mathcal M_d(\mathbb R)$, $h = (h_1, \dotsc, h_N) \in \left(\mathbb R_+^\ast\right)^N$, $r = (r_1, \dotsc, r_N) \in \left(\mathbb R_+^\ast\right)^N$, and assume that $\{h_1, \dotsc, h_N\}$ is rationally independent. Define
\begin{equation}
\label{eq:rho-HS}
\rho_{\mathrm{HS}} = \max_{\theta = (\theta_1, \dotsc, \theta_N) \in [0, 2\pi]^N} \spr\left(\sum_{k=1}^N A_k e^{i \theta_k}\right).
\end{equation}
Then the following assertions are equivalent.
\begin{enumerate}
\item\label{item:HS-criterion} We have $\rho_{\mathrm{HS}} < 1$.
\item\label{item:HS-rat-indep} System~\eqref{eq:intro-diff-eqn-sum} is exponentially stable with delays $\tau_1 = h_1, \dotsc, \tau_N = h_N$.
\item\label{item:HS-local} There exists $\varepsilon > 0$ such that \eqref{eq:intro-diff-eqn-sum} is exponentially stable with delays $\tau_1 = r_1^\prime, \dotsc, \tau_N = r_N^\prime$, for all $r^\prime = (r_1^\prime, \dotsc, r_N^\prime) \in \left(\mathbb R_+^\ast\right)^N$ with $\abs{r_k^\prime - r_k} < \varepsilon$ for every $k \in \llbracket 1, N\rrbracket$.
\item\label{item:HS-global} System~\eqref{eq:intro-diff-eqn-sum} is exponentially stable with delays $\tau_1 = r_1^\prime, \dotsc, \tau_N = r_N^\prime$, for all $r^\prime = (r_1^\prime, \dotsc, r_N^\prime) \in \left(\mathbb R_+^\ast\right)^N$.
\end{enumerate}
\end{proposition}

The most striking fact about the above result is that exponential stability of \eqref{eq:intro-diff-eqn-sum} for \emph{some} rationally independent family of positive delays $h_1, \dotsc, h_N$ is equivalent to exponential stability of \eqref{eq:intro-diff-eqn-sum} for \emph{any} choice of positive delays $r_1^\prime, \dotsc, r_N^\prime$ (with no constraints on their rational independence), and that this can be characterized by $\rho_{\mathrm{HS}} < 1$, where the quantity $\rho_{\mathrm{HS}}$ is independent of the delays. As a consequence, robustness of exponential stability of \eqref{eq:intro-diff-eqn-sum} with respect to delay perturbations is a sole property of the matrices $A_1, \dotsc, A_N$.

Property~\ref{item:HS-local} from Proposition~\ref{prop:Hale-Silkowski} is known as the \emph{local strong stability} of \eqref{eq:intro-diff-eqn-sum}, while \ref{item:HS-global} is known as the \emph{global strong stability} (or simply \emph{strong stability}) of \eqref{eq:intro-diff-eqn-sum}. These definitions of strong stability extend to neutral time-delay systems of the form \eqref{eq:neutral-linear} by requiring the same kind of stability for the associated difference equation \eqref{eq:difference-linear} (see, e.g., \cite{Hale2002Strong}).

Some of the equivalences of Proposition~\ref{prop:Hale-Silkowski} had been established in the earlier work \cite{Melvin1974Stability}, but only for scalar difference equations, this simple case being known as the Melvin Criterion. Note that, in this case, $\rho_{\mathrm{HS}} = \sum_{k = 1}^N \abs{A_k}$. A generalization of Proposition~\ref{prop:Hale-Silkowski} to difference equations involving also distributed delays (i.e., of the form \eqref{eq:intro-diff-eqn-integral-M} with $M$ being the sum of an element of $W^{1, 1}([-r, 0], \mathcal M_d(\mathbb{R}))$ and a piecewise constant function) was provided in \cite{Hu2011Stability}, but only perturbations in the pointwise delays (i.e., the discontinuity points of $M$) are considered. In particular, the notion of strong stability from \cite{Hu2011Stability} does not coincide with the one provided later in Chapter~\ref{chapter:strong-stability}. We also mention that Proposition~\ref{prop:Hale-Silkowski} can be extended to cases in which the delays satisfy some rational dependence structure, as detailed in \cite{Michiels2009Strong, Chitour2016Stability}, and a similar result also holds for systems with time-varying matrices \cite{Chitour2016Stability}.

\section{Destabilizing and stabilizing effects of delays}

It has to be noted that, in many of the early works in the literature, in particular in control theory, time delays are regarded from the perspective of their destabilizing effect on the systems where they appear, which is seen as ``generally undesirable in practice'', quoting \cite{Minorsky1942Self}. Indeed, the linear system $x^\prime(t) = - x(t)$ is exponentially stable, but $x^\prime(t) = - x(t - \tau)$ loses this stability property for $\tau \geq \frac{\pi}{2}$: its characteristic equation is $\Delta(s) = s + e^{- s \tau}$ and, using \cite[Theorem~1]{Hayes1950Roots}, it follows that all roots of $\Delta$ have negative real part if and only if $\tau \in [0, \frac{\pi}{2})$.

On the other hand, introducing a delay in a system may help in its stabilization. Indeed, consider a controlled harmonic oscillator described by the equation
\begin{equation}
\label{eq:expl-good-delay}
x^\prime(t) + x(t) = u(t),
\end{equation}
where $x(t) \in \mathbb R$ is the position of the oscillator and $u(t) \in \mathbb R$ is the control. A classical way to exponentially stabilize \eqref{eq:expl-good-delay} is to choose a proportional-derivative controller,
\begin{equation}
\label{eq:expl-PD-feedback}
u(t) = - k_p x(t) - k_d x^\prime(t),
\end{equation}
where $k_p$ and $k_d$ are constants. It is immediate to verify that \eqref{eq:expl-PD-feedback} renders \eqref{eq:expl-good-delay} exponentially stable if and only if $k_p > -1$ and $k_d > 0$. However, in some applications, instantaneous measures of the velocity $x^\prime(t)$, or even the position $x(t)$, may not be available, and hence it may be impossible to implement the control law \eqref{eq:expl-PD-feedback}. An alternative approach in this case is to use the control law
\begin{equation}
\label{eq:expl-delayed-feedback}
u(t) = - k x(t - \tau),
\end{equation}
where $k \in \mathbb R$ and the delay $\tau > 0$ are parameters of the control law. It was shown in \cite{Sipahi2011Stability} that there are choices of $k$ and $\tau$ such that \eqref{eq:expl-good-delay} with the control law \eqref{eq:expl-delayed-feedback} is exponentially stable.

The beneficial effects of time delays in control are known at least since \cite{Tallman1958Analog}, where the authors use a delay-based controller to improve the stability of systems with oscillatory behavior and small damping, and \cite{Suh1979Proportional, Suh1980Use, Atay1999Balancing}, which consider the approximation of the derivative action of a proportional-integral-derivative controller by a finite difference involving a delay. Several other works have explored ideas from the $\tau$-decomposition method introduced in \cite{Lee1969Tau} to characterize regions in the space of parameters of exponential stability of time-delay systems by considering the time delay as a free parameter. For instance, it was shown in \cite{Niculescu2004Stabilizing, Kharitonov2005Static} that a chain of $n$ integrators can be stabilized by a controller including $n$ delays, and stabilizing effects of the delays on oscillators was further explored in \cite{Niculescu2010Delay}. One of the main advantages of exploiting delays in control design is the simplicity of their implementation, and we refer to \cite{Sipahi2011Stability, Michiels2014Stability, Breda2015Stability} for an overview of methods and techniques.

\section{Quasipolynomials}
\label{sec:quasipolynomials}

The characteristic functions from \eqref{eq:intro-Delta} and \eqref{eq:intro-Delta-difference} are particular cases of \emph{quasipolynomials} (sometimes also called \emph{exponential polynomials}), a class of holomorphic complex functions extensively studied in the literature \cite{Avellar1980Zeros}, \cite[Chapter~3]{Berenstein1995Complex}, \cite{Gromova1968Asymptotic, Pontrjagin1942Zeros, Neimark1948Structure, Cebotarev1949Routh, Stavskii1968Division, Ronkin1978Quasipolynomials, Kharitonov1991Stability, Michiels2014Stability, Wielonsky2001Rolle}. More precisely, a \emph{quasipolynomial} $Q$ is an entire function $Q\colon \mathbb C \to \mathbb C$ which can be written under the form
\begin{equation}
\label{GenericQuasipolynomial}
Q(s) = \sum_{k = 0}^\ell P_k(s) e^{\lambda_k s},
\end{equation}
where $\ell$ is a nonnegative integer, $\lambda_0, \dotsc, \lambda_\ell$ are pairwise distinct real numbers, and, for $k \in \llbracket 0, \ell\rrbracket$, $P_k$ is a nonzero polynomial with complex coefficients of degree $d_k \geq 0$. The integer $D = \ell + \sum_{k=0}^\ell d_k$ is called the \emph{degree} of $Q$.

Recall that, given a polynomial $P$ of a single variable, one can always multiply $P$ by a nonzero constant in order to normalize the coefficient of its term of highest order to $1$ without changing its roots, and the degree of $P$ coincides with the number of remaining coefficients in $P$ (taking into account coefficients equal to $0$ in terms of order lower to the maximal one). A similar interpretation can be given to the degree $D$ of a quasipolynomial $Q$: up to multiplying $Q$ by a nonzero constant, one can normalize a nonzero coefficient of $Q$ to $1$, for instance, the coefficient of the term of highest degree in $P_0$. Since each polynomial $P_k$ of degree $d_k$ has $d_k+1$ coefficients, the number of remaining coefficients is $\sum_{k=0}^\ell (d_k + 1) - 1 = \ell + \sum_{k = 0}^\ell d_k$, which is precisely the number $D$ defined as the degree of a quasipolynomial.

Contrarily to the case of polynomials, the degree of a quasipolynomial does not determine the number of its roots, which is infinite except in trivial cases. However, there does exist a link between the degree of a quasipolynomial and the number of its roots in horizontal strips of the complex plane, thanks to a classical result known as \emph{P\'{o}lya--Szeg\H{o} bound} (see \cite[Part~Three, Problem~206.2]{Polya1998ProblemsI}), which we state in the next proposition.

\begin{proposition}
\label{prop:Polya-Szego}
Let $Q$ be a quasipolynomial of degree $D$ given under the form \eqref{GenericQuasipolynomial}, $\alpha, \beta \in \mathbb R$ be such that $\alpha \leq \beta$, and $\lambda_\delta = \max_{j, k \in \llbracket 0, \ell\rrbracket} \lambda_j - \lambda_k$. Let $m_{\alpha, \beta}$ denote the number of roots of $Q$ contained in the set $\{s \in \mathbb C \suchthat \alpha \leq \Imag s \leq \beta\}$ counting multiplicities. Then
\[
\frac{\lambda_\delta (\beta - \alpha)}{2 \pi} - D \leq m_{\alpha, \beta} \leq \frac{\lambda_\delta (\beta - \alpha)}{2 \pi} + D.
\]
\end{proposition}

Note that, as a consequence of Proposition~\ref{prop:Polya-Szego}, one deduces that a quasipolynomial $Q$ always has infinitely many roots, except for the trivial case where $\lambda_\delta = 0$. Indeed, if $\lambda_\delta > 0$, Proposition~\ref{prop:Polya-Szego} ensures that $Q$ has at least one root in the horizontal strip $\{s \in \mathbb C \suchthat \alpha \leq \Imag s \leq \beta\}$ as soon as $\beta - \alpha > \frac{2 \pi D}{\lambda_\delta}$, and the conclusion follows since one can find infinitely many such strips that are also pairwise disjoint.

Given a root $s_0 \in \mathbb C$ of a quasipolynomial $Q$, one immediately obtains, by letting $\beta = \alpha = \Imag s_0$ in the statement of Proposition~\ref{prop:Polya-Szego}, the following link between the multiplicity of $s_0$ and the degree of $Q$.

\begin{corollary}
\label{coro:Polya-Szego}
Let $Q$ be a quasipolynomial of degree $D$. Then any root $s_0 \in \mathbb C$ of $Q$ has multiplicity at most $D$.
\end{corollary}

\section{Multiplicity-induced dominance}
\label{sec:intro-MID}

On the one hand, we have introduced in Section~\ref{sec:stability-LTI-TDS} the notion of dominant root of the characteristic function of a time-delay system. On the other hand, we have presented, in Corollary~\ref{coro:Polya-Szego}, a bound on the multiplicity of roots of quasipolynomials. Although a priori unrelated, the notions of root of large (or maximal) multiplicity and of dominant root turn out to have some interesting links: it was observed at least since the 1950s that, for some simple, low-order time-delay systems, minimizing the spectral abscissa with respect to some free parameters of the system usually yields a spectral abscissa that is attained at a multiple root, which is thus dominant \cite[Chapter~III, \S~10; Chapter~IV, \S~6; Chapter~V, \S~7]{Pinney1958Ordinary}. Analogous properties also exist for linear finite-dimensional systems without delays, whose characteristic functions are polynomials, as it was shown, for instance, in \cite{Blondel2012Explicit, Chen1979Output}, where it is proved that, if the spectral abscissa admits a minimizer among a class of polynomials with an affine constraint on their coefficients, then one such minimizer is a polynomial with a single root of maximal multiplicity.

Despite the early results of \cite{Pinney1958Ordinary}, a deeper and more general analysis of the links between multiple and dominant roots, to the author's knowledge, has been left unexplored until the 2010s, when more systematic studies were carried out motivated by the links between minimizing the spectral abscissa and multiple roots on the one hand \cite{Michiels2002Continuous, Ramirez2015Design, Ramirez2016Design, Ramirez2017Analytical, Villafuerte2009Proportional, Villafuerte2013Tuning, Vanbiervliet2008Nonsmooth}, and by the analysis of the behavior of multiple roots on the other hand \cite{Gu2015Migration, Boussaada2015Inverted, Boussaada2016Tracking, Irofti2016Codimension, Boussaada2016Characterizing, Boussaada2014ComputingI, Boussaada2014ComputingII, Boussaada2016Further, Irofti2016Migration, Irofti2017Geometric, Boussaada2020Review, Michiels2017Explicit, Irofti2018Some}, and the property that real roots of large (or maximal) multiplicity tend to be dominant came to be known as \emph{multiplicity-induced dominance} (MID). 

One of the main interests of having a property such as MID lies in the stabilization of time-delay systems: when a control system such as \eqref{eq:neutral-control} has the MID property, a possible strategy for its stabilization is to choose the control $u$ under the form \eqref{eq:feedback} by selecting the coefficients of the matrices $K_0, K_{1, 1}, \dotsc, K_{1, N}, K_{2, 1}, \dotsc, K_{2, M}$ in order to guarantee the existence of a negative real root of large enough multiplicity, and the MID property will ensure that this multiple root is dominant, yielding exponential stability of the system. Hence, the interest in investigating multiple roots does not rely on the multiplicity itself, but rather on its connection with dominance of this root and the corresponding consequences for stability.

Most works on the MID property consider time-delay systems of the form
\begin{equation}
\label{eq:DDE}
y^{(n)}(t) + \sum_{k=0}^{n-1} a_k y^{(k)}(t) + \sum_{k=0}^m \alpha_k y^{(k)}(t - \tau) = 0,
\end{equation}
where $n$ and $m$ are nonnegative integers with $n \geq m$, $y$ is the real-valued unknown function, $\tau > 0$ is the delay, and $a_0, \dotsc, a_{n-1}, \alpha_0, \dotsc, \alpha_m$ are real coefficients. The characteristic function of \eqref{eq:DDE} is
\begin{equation}
\label{eq:Delta-DDE}
\Delta(s) = s^n + \sum_{k=0}^{n-1} a_k s^k + e^{-s \tau} \sum_{k=0}^m \alpha_k s^k,
\end{equation}
which is a quasipolynomial of degree $D = n + m + 1$. The first works on the MID property considered only particular cases of \eqref{eq:DDE}, such as \cite{Boussaada2016Multiplicity}, which considers the case $(n, m) = (1, 0)$ and proves that any double root of \eqref{eq:Delta-DDE} is dominant thanks to a suitable factorization of $\Delta$; \cite{Boussaada2018Further}, which addresses the case $(n, m) = (2, 0)$ using a similar factorization technique; \cite{Boussaada2018Dominancy, Boussaada2018Towards, Boussaada2020Multiplicity}, which consider the case $(n, m) = (2, 1)$ and prove that any root of multiplicity $4$ of \eqref{eq:Delta-DDE} is dominant by using Cauchy's argument principle; \cite{Mazanti2020Qualitative}, which considers the same case $(n, m) = (2, 1)$ but makes use instead of a factorization of $\Delta$ to prove dominance of roots of multiplicity $4$; or \cite{Villafuerte2013Tuning}, which also considers the case $(n, m) = (2, 1)$ and makes use of the $D$ partition method from \cite{Neimark1948Structure, Neimark1998D, Neimark1949D} to describe stability regions for \eqref{eq:DDE} in terms of $\alpha_0$, $\alpha_1$, and $\tau$ assuming that $a_1$ and $a_0$ are given and positive, noticing that, in this situation, the choice of $\alpha_0$, $\alpha_1$, and $\tau$ minimizing the spectral abscissa yields a configuration in which the spectral abscissa is attained at a dominant triple root. We refer to \cite{Benarab2020MID, Ma2020Multiplicity, Mazanti2021Effects, Boussaada2021Insights, Boussaada2018Further, Ma2022PID, Benarab2023MultiplicityI, Benarab2023MultiplicityII, OaxacaAdams2022Study, Villafuerte2009Proportional, Ramirez2015Design, Ramirez2016Design, Ramirez2017Analytical} for other works on the MID property considering particular configurations of \eqref{eq:DDE} and to \cite{Michiels2023Complete, OaxacaAdams2023Controllers} for further results on the links between minimization of the spectral abscissa and multiple roots, which illustrate in particular that, even though multiple spectral values often correspond to minimizers of the spectral abscissa, this is not necessarily always the case.

When one looks into the details of the previously mentioned results on multiplicity-induced dominance for \eqref{eq:DDE} with small values of $n$ and $m$, a first observation is that, whenever a root reaches the maximal multiplicity $n + m + 1$ from Corollary~\ref{coro:Polya-Szego}, it is necessarily dominant, but additional conditions are often needed for multiple roots with smaller multiplicity. This has motivated the search for a general proof of dominance of roots of maximal multiplicity of \eqref{eq:DDE}, which was first obtained in \cite{Mazanti2021Multiplicity} in the case $m = n-1$ using the same factorization of \cite{Mazanti2020Qualitative} and noticing that, when $\Delta$ from \eqref{eq:Delta-DDE} admits a root of maximal multiplicity $n + m + 1 = 2 n$, it can be factored in terms of a particular Kummer confluent hypergeometric function, whose roots had already been previously studied in \cite{Wynn1973Zeros}. The case of any pair of nonnegative integers $(n, m)$ with $m \leq n$ was addressed later in \cite{Boussaada2022Generic} using similar techniques, but the factorization of $\Delta$ now involves a more general Kummer confluent hypergeometric function, to which the results of \cite{Wynn1973Zeros} do not apply, and the proof of the main result of \cite{Boussaada2022Generic} relies instead on \cite{Boussaada2022Some}, which corrects results on the location of zeros of Kummer functions previously stated in \cite{Tsvetkov1941Complex, Tsvetkov1941Roots} by using the techniques described in \cite{Hille1922Oscillation}. A slightly simplified proof of the main result of \cite{Boussaada2022Generic} was later provided in \cite{Boussaada2022Pade}, which we present below in Section~\ref{sec:MID-maximal}.

The main issue in applications of results on the MID property with maximal multiplicity to stabilizability is that, in many applications, the free parameters of the system are not enough to ensure the existence of a root of maximal multiplicity $D = n + m + 1$. A possible strategy in this case is to consider multiple roots with multiplicity smaller than $D$, but such roots are not always necessarily dominant. A general theory providing necessary and sufficient conditions for the validity of the MID property for roots with multiplicity smaller than $D$ is still missing, but several important cases were considered in the literature. The main results in this sense are the ones provided in \cite{Balogh2020Towards, Balogh2022Conditions, Sassi2022Observer, RojasRicca2022Multiplicity} for the case $m = n - 1$ for roots of multiplicity $n + 1$, providing sufficient conditions for the dominance of such multiple roots. Some partial results and insights in more general situations are provided in \cite{Boussaada2022MID, Boussaada2025Decay} and are the subject of Section~\ref{sec:MID-intermediate}.

Another possible strategy to leverage additional degrees of freedom in the stabilization of time-delay systems is to consider, instead of a single root of maximal multiplicity $D$, the coexistence of $D$ distinct real roots. The study of this situation for small values of $n$ and $m$ was done in \cite{Bedouhene2019Pole, Amrane2020Some, Schmoderer2023Quasipolynomials, Schmoderer2024Insights, Amrane2018Qualitative, Bedouhene2020Real}, where it is shown that, in many situations, when the quasipolynomial $\Delta$ from \eqref{eq:Delta-DDE} admits $D$ distinct real roots, the largest of those roots is also the spectral abscissa of \eqref{eq:DDE}. The most general among those results is that of \cite{Bedouhene2020Real}, which proves such a dominance property in the case of any nonnegative integer $n$, but with $m = 0$, for $D = n + 1$ coexisting real roots. The property that the rightmost of such coexisting real roots is dominant came to be known as \emph{coexisting real roots-induced dominance} (CRRID), and some of the above works have also considered situations in which the number of coexisting real roots is strictly less than $D$.

Extensions of the MID property were also studied in the literature, including to pairs of multiple complex conjugate roots in the case $(n, m) = (2, 1)$ \cite{Mazanti2020Spectral} and to equations of order $n = 1$ and with two delays of order $0$ \cite{Fueyo2021Insights, Fueyo2023Pole}. Both analyses rely on perturbation arguments: in \cite{Mazanti2020Spectral}, the authors start from a configuration with a dominant real root of multiplicity $4$ and prove that, if this root splits into two double complex conjugate roots, dominance is preserved; while, in \cite{Fueyo2021Insights, Fueyo2023Pole}, the authors prove that, if one of the delays tends to $0$, one retrieves the classical MID for a single-delay system, and the dominance of the multiple root is preserved as the delay increases from $0$. Both results are provided below, in Sections~\ref{sec:MID-complex} and \ref{sec:MID-two-delays}.

Since the MID property was developed keeping in mind applications to the stabilization of time-delay systems, several works have explicitly addressed its applications (as well as the applications of the CRRID property), such as to active vibration control \cite{Boussaada2017Delayed, Tliba2019Active, Tliba2020Low, Tliba2020Delayed, FalconPrado2021Active, Tliba2022Design, Boussaada2017Coalescence}, to delayed resonators \cite{Kure2018Spectral}, to drilling systems \cite{Auriol2021Torsional, Auriol2022Comparing, Boussaada2019Low}, to quadrotor drones \cite{CastilloZamora2023Time}, to a proton exchange membrane fuel cell system \cite{Ramirez2021Proportional}, to wave and telegrapher's equations \cite{Ammari2026Multiplicity, Mazanti2021Effects}, or to biomechanical models, such as a muscle-tendon model \cite{ElAti2020Analysis} or human balance models \cite{Benarab2022Rolling, Molnar2021Calculation, Mei2022Multiplicity}. A software implementing stabilization methods through the MID and CRRID properties was also developed \cite{Boussaada2020Partial, Boussaada2021New, Boussaada2022New} and will be presented below in Chapter~\ref{chap:delay-software}. Finally, we refer the interested reader to the surveys \cite{Niculescu2016Critical, Niculescu2021Stability, Niculescu2023Counting, Benarab2022Over} containing, among others, more details on the MID and CRRID properties and other stabilization strategies for time-delay systems.

\section{Stability of systems with time-dependent coefficients and delays}
\label{sec:intro-varying}

Even though we have mentioned, in Sections~\ref{sec:stability-LTI-TDS} and \ref{sec:stabilization-LTI-TDS}, several results on the stability and stabilization of linear time-delay systems with constant coefficients, there are also many classical stability and stabilization results available in the case of time-dependent coefficients (see, e.g., \cite[Part~III]{Briat2015Linear}, \cite[Chapter~6]{Gu2003Stability}, \cite[Chapter~4]{Halanay1966Differential}, \cite[Chapters~6 and 8 and Section~9.1]{Hale1993Introduction}). One of the main difficulties in the analysis of time-delay systems with time-varying parameters is that a spectral characterization of stability is not available, except for some particular situations, such as the case of periodic parameter variation, in which a Floquet-type theory can be developed \cite[Chapter~8]{Hale1993Introduction}.

One of the approaches to deal with time-varying matrices in \eqref{eq:neutral-linear}, \eqref{eq:difference-linear}, or \eqref{eq:neutral-control}, particularly useful when the time dependence of those matrices is not precisely known, is the robustness approach, which is the point of view adopted, for instance, in \cite[Chapter~8]{Briat2015Linear} and \cite[Chapter~6]{Gu2003Stability}. This approach consists in assuming that, although the time dependence of the matrices in \eqref{eq:neutral-linear}, \eqref{eq:difference-linear}, or \eqref{eq:neutral-control} is not precisely known, one knows that these matrices will always take values in some compact sets of matrices, and the goal is to study whether the system is stable or stabilizable for \emph{any} choice of time-dependent matrices taking values in these sets. Typical tools in this kind of analysis are Lyapunov and Lyapunov--Krasovski\u{\i} functions, which usually provide sufficient conditions for stability or stabilizability, often under the form of linear matrix inequalities (LMIs\label{LMI}). A notable exception is \cite{Chitour2016Stability}, which considers \eqref{eq:difference-linear} with time-varying matrices $H_1, \dotsc, H_N$, providing necessary and sufficient stability conditions based on a suitable representation formula for solutions, generalizing also the Hale--Silkowski criterion from Proposition~\ref{prop:Hale-Silkowski} to time-dependent matrices.

The case of time-varying delays has also been extensively studied in the literature \cite{Fridman2008Robust, Fridman2007Stability, Kao2007Stability, Shustin2007Delay, Fridman2006Input, Louisell2001Delay, Michiels2005Stabilization, Verriest2012State, Bonnet2020L2, Mazanti2014Stabilization}. Similar techniques based on Lyapunov and Lya\-pu\-nov--Krasovski\u{\i} functions can also be used to treat this case, providing sufficient conditions for stability or stabilizability which are often under the form of LMIs. In addition, many works also require, for technical reasons, conditions on the rate of variation of the delay, such as assuming that a given time-dependent delay $\tau(t)$ satisfies $\tau^\prime(t) < 1$ for every $t$, or $\tau^\prime(t) \leq c$ for some constant $c < 1$ and every $t$, the main interest of these assumptions being to ensure that $t \mapsto t - \tau(t)$ is increasing. Other strategies for dealing with time-dependent delays include perturbation approaches, used, for instance, in \cite{Bonnet2020L2, Mazanti2014Stabilization}, which consist in establishing first a stability or stabilization result for constant nominal values of the delays, and then proving that such stability properties are preserved if the delays are time-varying and take values in small intervals containing the nominal values. A sufficient condition for the stability of \eqref{eq:difference-linear} with time-varying matrices and delays was also provided in \cite{Avellar1990Difference}, where, inspired by the previous work \cite{Banas1982Existence} for (possibly nonlinear) systems with a single delay, the authors make use of a fixed-point technique to prove that solutions necessarily belong to a suitably chosen functional space containing only functions converging to the origin as time tends to infinity.

It has to be noted that, even if a time-delay system with constant parameters and delays is stable for all constant choices of parameters and delays in a given set, the system may become unstable if the parameters are allowed to vary in time, even if their values remain in the given set at all times. The following example, adapted from \cite[Section~4]{Avellar1990Difference}, illustrates this fact in a simple difference equation in dimension $2$.

\begin{example}
\label{expl:AvellarMarconato}
Let $\alpha \in (0, 1)$ and consider the matrices
\[
A_1 = \begin{pmatrix}
\alpha & 0 \\
0 & -\alpha \\
\end{pmatrix}, \qquad A_2 = \begin{pmatrix}
0 & \alpha \\
-\alpha & 0 \\
\end{pmatrix}.
\]
Let us compute the quantity $\rho_{\mathrm{HS}}$ from \eqref{eq:rho-HS} for this pair of matrices. For $\theta_1$ and $\theta_2$ in $[0, 2\pi]$, we have
\[
A_1 e^{i \theta_1} + A_2 e^{i \theta_2} = \begin{pmatrix}
\alpha e^{i \theta_1} & \alpha e^{i \theta_2} \\
-\alpha e^{i \theta_2} & -\alpha e^{i \theta_1} \\
\end{pmatrix},
\]
and an immediate computation shows that
\[
\spr\left(A_1 e^{i \theta_1} + A_2 e^{i \theta_2}\right) = \alpha \sqrt{\abs{e^{2 i \theta_1} - e^{2 i \theta_2}}},
\]
yielding that $\rho_{\mathrm{HS}} = \sqrt{2} \alpha$. In particular, by Proposition~\ref{prop:Hale-Silkowski}, for every $\alpha \in \left(0, \frac{\sqrt{2}}{2}\right)$, the difference equation
\begin{equation}
\label{eq:expl-Avellar}
x(t) = A_1 x(t - r_1) + A_2 x(t - r_2)
\end{equation}
is exponentially stable for every choice of positive delays $r_1, r_2$.

Let us now consider the difference equation \eqref{eq:expl-Avellar} with a time-varying delay $r_1$ and with $r_2 = 1$, i.e.,
\begin{equation}
\label{eq:expl-Avellar-varying}
x(t) = A_1 x(t - r_1(t)) + A_2 x(t - 1).
\end{equation}
We assume further that $r_1$ is continuous, $2$-periodic, takes values in $[1, 2]$, and $r_1(0) = 1$, $r_1(1) = 2$. Hence, for every $k \in \mathbb N$, we have
\begin{align*}
x(2k) & = (A_1 + A_2) x(2k - 1), \\
x(2k + 1) & = A_1 x(2k - 1) + A_2 x(2k) = A_0 x(2k - 1),
\end{align*}
where
\[
A_0 = A_1 + A_2(A_1 + A_2) = \begin{pmatrix}
\alpha - \alpha^2 & - \alpha^2 \\
-\alpha^2 & -\alpha - \alpha^2
\end{pmatrix}.
\]
An immediate inductive argument shows that
\[x(2k + 1) = A_0^{k+1} x(-1).\]
The eigenvalues $\lambda_1, \lambda_2$ and respective eigenvectors $v_1, v_2$ of $A_0$ are
\begin{align*}
\lambda_1 & = - \alpha^{2} - \alpha \sqrt{\alpha^{2} + 1}, & v_1 & = \begin{pmatrix} \alpha \left(\sqrt{\alpha^{2} + 1} - 1\right) \\ \alpha^2 \end{pmatrix},\\
\lambda_2 & = - \alpha^{2} + \alpha \sqrt{\alpha^{2} + 1}, & v_2 & = \begin{pmatrix} -\alpha \left(\sqrt{\alpha^{2} + 1} + 1\right) \\ \alpha^2\end{pmatrix}.
\end{align*}
In particular, we have $\lambda_1 < -1$ if $\alpha > \frac{\sqrt{3}}{3}$. Thus, if we consider a solution of \eqref{eq:expl-Avellar-varying} with initial condition such that $x(-1) = v_1$, we have
\begin{equation}
\label{eq:sol-Avellar-odd}
x(2k + 1) = \lambda_1^{k+1} v_1,
\end{equation}
and $\abs{x(2k + 1)} \to +\infty$ if $\alpha > \frac{\sqrt{3}}{3}$. Hence, for $\alpha \in \left(\frac{\sqrt{3}}{3}, \frac{\sqrt{2}}{2}\right)$, \eqref{eq:expl-Avellar} is exponentially stable for every choice of positive delays $r_1, r_2$, but \eqref{eq:expl-Avellar-varying} admits unbounded solutions.

We illustrate, in Figure~\ref{fig:AvellarMarconato}, the solution of \eqref{eq:expl-Avellar-varying} with $\alpha = \frac{2}{3}$, $r_1(t) = \frac{1}{2} (3 - \cos(\pi t))$, and initial condition $x_0(t) = (1 + t) (A_1 + A_2) v_1 - t v_1$ for $t < 0$, in which case
\[
\lambda_1 = -\frac{4 + 2 \sqrt{13}}{9} \approx -1.2457, \qquad v_1 = \begin{pmatrix} \frac{2}{9}(\sqrt{13} - 3) \\ \frac{4}{9} \end{pmatrix} \approx \begin{pmatrix} 0.13457 \\ 0.44444\end{pmatrix}.
\]
The initial condition in the interval $[-2, 0)$ is represented by a dashed black line, the solution $x$ in the interval $[0, 10]$ is represented by a solid thin purple line, the values of the solution at odd integer times are represented by blue circles, and the values of the solution at even integer times are represented by red squares. Note that the values of the solution at odd integer times satisfy \eqref{eq:sol-Avellar-odd}, while the values of the solution at even integer times satisfy $x(2k) = \lambda_1^k (A_1 + A_2) v_1$.

\begin{figure}[ht]
\centering
\resizebox{0.667\textwidth}{!}{\input{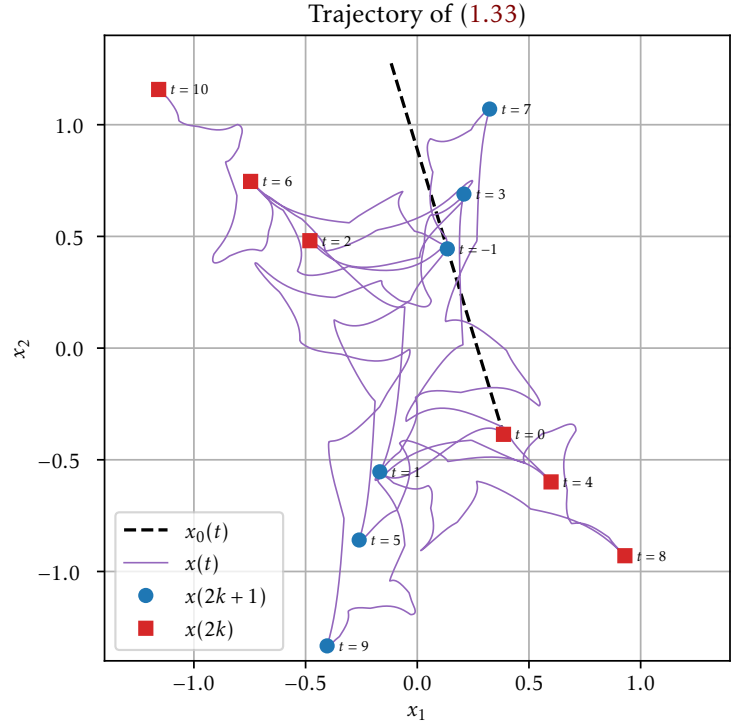}}
\caption[Unbounded trajectory of \eqref{eq:expl-Avellar-varying}]{Unbounded trajectory of \eqref{eq:expl-Avellar-varying} with initial condition $x_0(t) = (1 + t) (A_1 + A_2) v_1 - t v_1$ for $t < 0$}
\label{fig:AvellarMarconato}
\end{figure}
\end{example}

\section{Description of the main results}

Let us now describe the main results proved in the next chapters. Our first main results, presented in Chapter~\ref{chap:controllability}, concern the controllability of \eqref{eq:intro-linear-difference-equation} in $L^p$ spaces and are extracted from \cite{Chitour2023Hautus}: Theorem~\ref{main_result1} provides necessary and sufficient conditions for the approximate controllability of \eqref{eq:intro-linear-difference-equation} in $L^p$ for $p \in [1, +\infty)$, under the form of a Hautus-type condition, while Theorem~\ref{main_result2} provides necessary conditions for the exact controllability of \eqref{eq:intro-linear-difference-equation} in $L^1$, which turn out to be the ``closure'' of the Hautus-type conditions from Theorem~\ref{main_result1}. To obtain these two results, an important intermediate step, which has an interest on its own, is Theorem~\ref{lem:3}, which relates approximate (respectively, exact) controllability from the origin in free final time with approximate (respectively, exact) controllability in a specific prescribed time, which depends only on the dimension $d$ and on the maximal delay appearing in \eqref{eq:intro-linear-difference-equation}.

We then move, in Chapter~\ref{chapter:strong-stability}, to the sensitivity analysis of exponential stability of \eqref{eq:intro-diff-eqn-integral-M} with respect to perturbations in $M$, one of our main aims being to generalizing Proposition~\ref{prop:Hale-Silkowski} to this class of systems. After a short general discussion, Chapter~\ref{chapter:strong-stability}, which presents the contents of \cite{ChitourStrong, GoncalvesNetto2025Strong}, start by a well-posedness result for \eqref{eq:intro-diff-eqn-integral-M}. In order to consider the sensitivity of the exponential stability of \eqref{eq:intro-diff-eqn-integral-M} with respect to perturbations in $M$, the first step is to properly define what should be considered as a perturbation of $M$. For the particular case of \eqref{eq:intro-diff-eqn-sum} in which $M$ is given by
\begin{equation}
\label{eq_M-finite-sum}
M = \sum_{k = 1}^N A_k \1_{[-\tau_k, 0]},
\end{equation}
the matrices $A_1, \dotsc, A_N$ and the delays $\tau_1, \dotsc, \tau_N$ play different roles. In particular, a perturbation of $A_1, \dotsc, A_N$ yields a perturbation of $M$ in the total variation norm, but this is not the case for a variation of the delays $\tau_1, \dotsc, \tau_N$ (see Example~\ref{expl_var-delay}).

The first main result we prove in Chapter~\ref{chapter:strong-stability} concerns perturbations of $M$ in total variation norm, and can be seen as the counterpart of the fact that exponential stability of \eqref{eq:intro-diff-eqn-sum} is preserved under perturbations of the matrices $A_1, \dotsc, A_N$. In particular, we prove that, if $M$ is such that \eqref{eq:intro-diff-eqn-integral-M} is exponentially stable, then \eqref{eq:intro-diff-eqn-integral-M} will remain exponentially stable for all other matrix-valued functions close enough to $M$ in total variation norm (see Corollary~\ref{coro_tvstab}), at least for an important class of functions $M$, those without a singular part. As for the generalization of Proposition~\ref{prop:Hale-Silkowski} to \eqref{eq:intro-diff-eqn-integral-M}, we first provide a suitable notion of ``perturbation of the delays'' of \eqref{eq:intro-diff-eqn-integral-M} and generalize thus the definition of strong stability and of $\rho_{\mathrm{HS}}$ from \eqref{eq:rho-HS} to \eqref{eq:intro-diff-eqn-integral-M} (see Definition~\ref{def_strong-stab} and \eqref{GHSC}). We conjecture that the counterpart of Proposition~\ref{prop:Hale-Silkowski} holds true in this generalized setting (see Conjectures~\ref{GHS-Conjecture} and \ref{GHS-Conjecture-SDG}), and we prove this conjecture to be true in two cases: for scalar systems (Theorem~\ref{HS:scalar}) and for systems satisfying the spectrum-determined growth condition with a modified notion of strong stability (Theorem~\ref{GHS-Theorem}).

The next set of main results concern a difference equation with a single time-varying delay (namely, \eqref{eq:diff-eqn-N1-time} in Chapter~\ref{chap:diff-eqn-varying-delay}), and is extracted from \cite{MazantiWellPosedness}. This equation is studied in essentially three kinds of functional spaces: those of continuous functions, of regulated functions, and in $L^p$ for $p \in [1, +\infty]$, and the main results concern existence and uniqueness of solutions and their asymptotic behavior in each of these function spaces. More precisely, Theorems~\ref{thm:exist-continuous}, \ref{thm:exist-regulated}, \ref{thm:exist-L-infty}, and \ref{thm:exist-L-p} provide sufficient conditions for existence and uniqueness of continuous, regulated, $L^\infty$, and $L^p$, $p \in [1, +\infty)$, solutions, respectively. The analysis of whether those solutions converge to $0$ as time tends to infinity is done in Corollaries~\ref{Corol-1-C}, \ref{Prop-1-C}, and \ref{coro:cv-to-0-Linfty} and Theorem~\ref{Thm-1-L-P} for continuous, regulated, $L^\infty$, and $L^p$, $p \in [1, +\infty)$, solutions, respectively. Finally, the study of whether the previous convergence to $0$ holds with an exponential rate is the subject of Corollaries~\ref{coro:exp-cont}, \ref{coro:exp-G}, and \ref{coro:exp-Linfty} and Theorem~\ref{thm:exp-Lp} for continuous, regulated, $L^\infty$, and $L^p$, $p \in [1, +\infty)$, solutions, respectively.

All other main results of the monograph consider the multiplicity-induced dominance (MID) property described in Section~\ref{sec:intro-MID} and are extracted from \cite{Boussaada2022Some, Boussaada2022Generic, Mazanti2021Multiplicity, Boussaada2022Pade, Mazanti2021Effects, Benarab2020MID, Mazanti2020Qualitative, Boussaada2025OverOrder, Boussaada2025Decay, Boussaada2022MID, Mazanti2020Spectral, Fueyo2021Insights, Fueyo2023Pole, Niculescu2021Stability}. The main result characterizing the MID property for delay-differential equations of the form \eqref{eq:DDE} is Theorem~\ref{thm:main-MID}, which provides a link between a root of maximal multiplicity of the characteristic function \eqref{eq:Delta-DDE} of \eqref{eq:DDE}, Padé approximations of the exponential function, and Kummer confluent hypergeometric functions. As for roots with ``large enough'', but not necessarily maximal, multiplicity, we provide characterizations of the existence of such roots in terms of factorizations of the characteristic function \eqref{eq:Delta-DDE} in Theorem~\ref{thm:max-mult} and Corollary~\ref{coro:max-mult-Kummer}, and in terms of the unknown coefficients of the system in Theorem~\ref{thm:charact-EPF}. We also point out that a more detailed discussion on what is ``large enough'' multiplicity is provided in the beginning of Section~\ref{sec:MID-intermediate}.

All the previously discussed results on the MID property concern systems of the form \eqref{eq:DDE} and a \emph{real} root of ``large enough'' or maximal multiplicity. Theorem~\ref{thm:complex-dominance} considers instead the MID property for a complex conjugate pair of roots of maximal multiplicity for \eqref{eq:DDE} with $n = 2$ and $m = 1$, while Theorem~\ref{thm:main-two-delays} considers the MID property for a simple delay-differential equation with two delays, with maximal order of derivation equal to $1$ and no derivatives in the delayed terms. Both of these theorems are proved through a similar technique, which is based on considering a family of one-parameter systems, obtaining the MID property for a special value of the parameter by applying Theorem~\ref{thm:main-MID}, and then showing that the MID property is preserved as the parameter varies in its domain.

\clearpage

\section{Structure of the monograph}

In this monograph, we will study first the controllability of difference equations of the form \eqref{eq:intro-linear-difference-equation} in Chapter~\ref{chap:controllability}, proceed to study the sensitivity of exponential stability of \eqref{eq:intro-diff-eqn-integral-M} with respect to perturbations in $M$ in Chapter~\ref{chapter:strong-stability}, then study a difference equation with a single time-varying delay in Chapter~\ref{chap:diff-eqn-varying-delay}, and finally consider the MID property for \eqref{eq:DDE} and related equations in Chapter~\ref{chap:MID}. A brief comment on software for time-delay systems is presented in Chapter~\ref{chap:delay-software}, and the monograph is concluded by a discussion of open problems in Chapter~\ref{chap:delay-open-problems}.

\chapter{Controllability of difference equations}
\label{chap:controllability}

\epigraph{\foreignlanguage{brazil}{Cultivar a ciência pela utilidade prática, imediata, é desvirtuar a alma da própria ciência!}}{\foreignlanguage{brazil}{\emph{O Homem que Calculava}, Malba Tahan}}

The aim of this chapter is to present the results on approximate and exact controllability in $L^p$, $p \in [1, +\infty)$, of the linear difference equation \eqref{eq:intro-linear-difference-equation} obtained in \cite{Chitour2023Hautus}. We will always assume here that the delays $\Lambda_1, \dotsc, \Lambda_N$ from \eqref{eq:intro-linear-difference-equation} are ordered as $\Lambda_1 < \dotsb < \Lambda_N$.

The chapter is organized as follows. In order to state and prove our main results, we shall need some notation and elementary results from the theory of distributions, which are presented in Section~\ref{sec:distributions}. We then proceed to providing the definitions of controllability used in this chapter in Section~\ref{sec:statements}, where we also provide the statement of our main results, Theorems~\ref{main_result1} and \ref{main_result2}, as well as other important results proved in this chapter, Theorems~\ref{lem:3} and \ref{thm:Ck} and Propositions~\ref{Prop2:main_result0} and \ref{Prop2:main_result}. Section~\ref{sec:6} provides an explicit formula for solutions of \eqref{eq:intro-linear-difference-equation} in terms of their initial conditions and controls, and uses it to characterize controllability in terms of the range of an endpoint map, exploiting the properties of this map to provide a proof of Theorem~\ref{lem:3}. The main goal of Section~\ref{section_Hautus_criteria} is to provide the proofs of Theorems~\ref{main_result1}, \ref{main_result2}, and \ref{thm:Ck} and Propositions~\ref{Prop2:main_result0} and \ref{Prop2:main_result}, and this is done by first rewriting \eqref{eq:intro-linear-difference-equation} as an input/output system, as described in Section~\ref{sec:Realization} above, in order to make use of the results by Y.~Yamamoto on realization theory \cite{Yamamoto1978Realization, Yamamoto1981RealizationI, Yamamoto1981RealizationII, Yamamoto1988Pseudo, Yamamoto1989Reachability, Yamamoto2007Coprimeness, Yamamoto2011Bezout, Yamamoto2008Behavioral, Yamamoto2016Behavioral, Yamamoto1999Coprimeness}. We finally conclude this chapter by applications of our main results to the cases $(N, d, m) = (2, 2, 1)$ and $(N, d, m) = (2, 3, 1)$ in Section~\ref{sec:applications}, showing, in the case $(N, d, m) = (2, 2, 1)$, that one can retrieve the approximate controllability results previously proved in the $L^2$ framework in \cite{Chitour2020Approximate} and that the necessary conditions for $L^1$ exact controllability of \eqref{eq:intro-linear-difference-equation} from Theorem~\ref{main_result2} turn out to coincide with the necessary and sufficient conditions for $L^2$ exact controllability from \cite{Chitour2020Approximate}.

\section{Preliminaries on the theory of distributions}
\label{sec:distributions}

We now recall several notions from the theory of distributions that will be used in the sequel of the chapter. A detailed presentation of these notions with precise definitions can be found, e.g., in \cite{Schwartz1966Theorie, Yamamoto1988Pseudo, Yamamoto1989Reachability}.

Throughout this chapter, we use $\mathcal{D}(\mathbb{R})$ to denote the space of $\mathcal C^{\infty}$ functions defined on $\mathbb{R}$ with compact support, endowed with its canonical LF topology, and $\mathcal{D}^\prime(\mathbb{R})$ to denote the space of continuous linear forms acting on $\mathcal{D}(\mathbb{R})$, i.e., the space of all distributions on $\mathbb{R}$, endowed with the strong dual topology of uniform convergence on bounded subsets of $\mathcal{D}(\mathbb{R})$. The support of a distribution $\alpha \in \mathcal D^\prime(\mathbb R)$, denoted by $\supp\alpha$, is the complement of the largest open set $U$ such that $\scalprod{\alpha, f} = 0$ for every $f \in \mathcal D(\mathbb R)$ with $\supp f \subset U$. The order of a distribution $\alpha \in\mathcal{D}^\prime(\mathbb{R})$ is the smallest integer $\ell$ such that, for every compact set $K \subset \mathbb R$, there exists $C_K > 0$ such that
\begin{equation*}
\abs{\scalprod{\alpha,\psi}} \le C_K \sup_{x \in K} \abs{\psi^{(\ell)}(x)},
\end{equation*}
for all $\psi \in \mathcal D(\mathbb R)$ with compact support in $K$. For $x \in \mathbb R$, we denote the Dirac distribution at $x$ by $\delta_x \in \mathcal{D}^\prime(\mathbb{R})$, and we recall that $\delta_x$ has order zero.

We shall also use in this chapter the subspaces $\mathcal{E}^\prime(\mathbb{R}_{-})$, $\mathcal{E}^\prime(\mathbb{R}_{+})$, $\mathcal{D}^\prime_{+}(\mathbb{R})$, $M(\mathbb{R}_{-})$, $M(\mathbb{R}_{+})$, and $M_+(\mathbb{R})$ of $\mathcal D^\prime(\mathbb R)$, as well as the space of test functions $\mathcal D(\mathbb R_+)$ and its topological dual $\mathcal D^\prime(\mathbb R_+)$ (not to be confused with $\mathcal D^\prime_+(\mathbb R)$), whose definitions we provide in Table~\ref{tab:subspaces-distributions}.

\begin{table}[ht]
\caption{Definitions of some spaces of test functions and distributions}
\label{tab:subspaces-distributions}
\begin{tabular}{@{} >{\raggedright} p{0.18\textwidth} @{\hspace*{0.02\textwidth}} p{0.8\textwidth} @{}}
\toprule
\textbf{Space} & \textbf{Definition} \tabularnewline
\midrule
$\mathcal D(\mathbb R_+)$ & Subspace of $\mathcal D(\mathbb R)$ of functions whose support is included in $\mathbb R_+$. \tabularnewline
$\mathcal{E}^\prime(\mathbb{R}_{-})$, $\mathcal{E}^\prime(\mathbb{R}_{+})$ & Spaces of distributions having compact support in $\mathbb{R}_{-}$ and $\mathbb{R}_{+}$, respectively. \tabularnewline
$\mathcal{D}^\prime_{+}(\mathbb{R})$ & Space of distributions with support bounded on the left. \tabularnewline
$\mathcal D^\prime(\mathbb R_+)$ & Topological dual of $\mathcal D(\mathbb R_+)$. \tabularnewline
$M(\mathbb{R}_{-})$, $M(\mathbb{R}_{+})$ & Subspaces of $\mathcal{E}^\prime(\mathbb{R}_{-})$ and $\mathcal{E}^\prime(\mathbb{R}_{+})$, respectively, consisting of distributions of order zero. \tabularnewline
$M_+(\mathbb{R})$ & Subspace of $\mathcal{D}^\prime_{+}(\mathbb{R})$ of distributions of order zero. \tabularnewline
\bottomrule
\end{tabular}
\end{table}

Note that, according to the definitions in Table~\ref{tab:subspaces-distributions}, the support of a distribution in $\mathcal{E}^\prime(\mathbb{R}_{-})$ or $\mathcal{E}^\prime(\mathbb{R}_{+})$ may contain $0$, and we have the inclusions $\mathcal E^\prime(\mathbb R_+) \subset \mathcal D^\prime_+(\mathbb R)$ and $M(\mathbb R_+) \subset M_+(\mathbb R)$. Moreover, thanks to the Riesz representation theorem  (see, e.g., \cite[Theorem~6.19]{Rudin1987Real}), the spaces $M(\mathbb{R}_{-})$ and $M(\mathbb{R}_{+})$ are the sets of Radon measures with compact support contained in $\mathbb{R}_{-}$ and $\mathbb{R}_{+}$, respectively, while $M_+(\mathbb{R})$ is the space of Radon measures with support bounded on the left. We shall always consider that the space $\mathcal{D}^\prime_{+}(\mathbb{R})$ is endowed with the convolution product~$*$, which makes it a commutative algebra. With a slight abuse of notation, we also write $\mathcal D(\mathbb R)$, $\mathcal D^\prime(\mathbb R)$, $\mathcal D(\mathbb R_+)$, $\mathcal{E}^\prime(\mathbb{R}_{-})$, $\mathcal{E}^\prime(\mathbb{R}_{+})$, $\mathcal D^\prime_+(\mathbb R)$, $\mathcal D^\prime(\mathbb R_+)$, $M(\mathbb{R}_{-})$, $M(\mathbb{R}_{+})$, and $M_+(\mathbb R)$ to refer to sets of matrices whose entries belong to those respective spaces.

Given a Radon measure $\mu \in M_+(\mathbb{R})$, we use $\widehat{\mu}(s)$ to denote the two-sided Laplace transform of $\mu$ evaluated at $s\in \mathbb{C}$, that is,
\begin{equation}
\label{laplace_transform_radon_measure}
\widehat{\mu}(s)=\int_{-\infty}^{+\infty} \diff\mu(t) e^{-st},
\end{equation}
provided that the integral exists. In particular, the Laplace transform of a Dirac distribution $\delta_x$, $x \in \mathbb{R}$, is the holomorphic map
\begin{equation}
\label{laplace_transform_dirac}
\widehat{\delta}_x\colon s\mapsto e^{-s x},\qquad s\in \mathbb{C}.
\end{equation}
In the more general setting of a distribution $\alpha \in\mathcal{D}^\prime(\mathbb{R})$, the Laplace transform $\widehat{\alpha}(s)$ of $\alpha$ at $s\in\mathbb{C}$ is defined by
$\scalprod{\alpha, e^{- \cdot s}}$ when this quantity makes sense, i.e., when the linear form $\alpha$ on $\mathcal D(\mathbb R)$ can be extended by continuity (in the standard topology of $\mathcal C^\infty(\mathbb R)$) to the function $t \mapsto e^{-t s}$.

We now define the truncation to positive times of a distribution in $\mathcal D^\prime_+(\mathbb R)$, following the presentation of \cite{Yamamoto1988Pseudo, Yamamoto1989Reachability}. For $\alpha \in \mathcal{D}^\prime_{+}(\mathbb{R})$, we define its truncation $\pi\alpha \in \mathcal D^\prime(\mathbb R_+)$ by
\begin{equation*}
\scalprod{\pi \alpha, \psi} = \scalprod{\alpha, \psi},\qquad \psi \in \mathcal{D}(\mathbb{R}_+).
\end{equation*}
We have that $\pi \alpha$ is a well-defined element of $\mathcal{D}^\prime(\mathbb{R}_+)$. Furthermore, $\pi$ is continuous with respect to the strong dual topology of $\mathcal{D}^\prime_{+}(\mathbb{R})$ and $\mathcal{D}^\prime(\mathbb{R}_+)$. Note that, since any function $\psi \in \mathcal D(\mathbb R_+)$ satisfies $\psi(0) = 0$, the operator $\pi$ truncates distributions to positive times only and, in particular, $\pi\delta_0 = 0$. In addition, we have the following properties of $\pi$, whose proofs can be found in \cite[Lemmas~A2 and A3]{Yamamoto1988Pseudo}.

\begin{lemma}
\label{lem:sumlemmasYamamoto}
The following assertions hold true.
\begin{enumerate}
\item \label{lem:A_2_Yamamoto} For $ \alpha \in \mathcal{D}^\prime(\mathbb{R})$, we have $\pi(\alpha)=0$ if and only if $\supp(\alpha) \subset (-\infty,0]$.
\item \label{lem:A_3_Yamamoto} $\pi(\alpha * \pi \beta)=\pi(\alpha*\beta)$ for every $\alpha \in \mathcal{E}^\prime(\mathbb{R}_{-})$ and $\beta \in \mathcal{D}^\prime_{+}(\mathbb{R})$.
\end{enumerate}
\end{lemma}

\section{Basic definitions and statement of the main results}
\label{sec:statements}

We now turn to the main topic of this chapter, the study of the controllability of the linear difference equation \eqref{eq:intro-linear-difference-equation}. We start by recalling (see, e.g., \cite[Definition~2.1]{Mazanti2017Relative}) that, given $T > 0$, an initial condition $x_0\colon [-\Lambda_N, 0) \to \mathbb R^d$, and a control $u\colon [0, T] \to \mathbb R^m$, a function $x\colon [-\Lambda_N, T] \to \mathbb R^d$ is said to be a \emph{solution} of \eqref{eq:intro-linear-difference-equation} if $x(t) = x_0(t)$ for every $t \in [-\Lambda_N, 0)$ and the equality in \eqref{eq:intro-linear-difference-equation} is satisfied for every $t \in [0, T]$. Existence and uniqueness of solutions can be obtained by a direct step-by-step construction of the solution, as given in \cite[Proposition~3.2]{Chitour2016Stability}, and, as noted in \cite[Remark~2.3]{Mazanti2017Relative}, solutions will belong to $L^p$ if the initial condition and the control are in $L^p$. More precisely, we have the following result.

\begin{proposition}
Let $p \in [1, +\infty]$, $T > 0$, $u \in L^p([0,T],\mathbb{R}^m)$, and $x_0 \in L^p([-\Lambda_N,0],\mathbb{R}^d)$. Then there exists a unique solution $x \in L^p([-\Lambda_N, T],\mathbb{R}^d)$ of \eqref{eq:intro-linear-difference-equation} with initial condition $x_0$ and control $u$.
\end{proposition}

Given a function $x\colon [-\Lambda_N, T] \to \mathbb R^d$ and $t \in [0, T]$, we will consider the history function of $x$ at time $t$ as the function $x_t\colon [-\Lambda_N, 0] \to \mathbb R^d$ defined by $x_t(\theta) = x(t + \theta)$ for every $\theta \in [-\Lambda_N, 0]$. With this notation, we are now in position to provide the definitions of the standard notions of approximate and exact controllability in finite time $T > 0$ that we are interested to study in this chapter.

\begin{definition}\label{def:Lq-cont}
Let $p \in [1, +\infty]$ and $T > 0$. System~\eqref{eq:intro-linear-difference-equation} is said to be
\begin{enumerate}
\item \emph{$L^p$ approximately controllable in time $T>0$} if, for every $x_0$ and $\phi$ in $L^p([-\Lambda_N,0], \allowbreak\mathbb{R}^d)$ and $\epsilon > 0$, there exists $u \in L^p([0,T],\mathbb{R}^m)$ such that the solution $x$ of \eqref{eq:intro-linear-difference-equation} with initial condition $x_0$ and control $u$ satisfies
\[
\norm{x_T - \phi}_{L^p([-\Lambda_N, 0], \mathbb R^d)} < \epsilon.
\]
\item \emph{$L^p$ exactly controllable in time $T>0$} if, for every $x_0$ and $\phi$ in $L^p([-\Lambda_N,0],\mathbb{R}^d)$, there exists $u \in L^p([0,T],\mathbb{R}^m)$ such that the solution $x$ of \eqref{eq:intro-linear-difference-equation} with initial condition $x_0$ and control $u$ satisfies
\[
x_T = \phi.
\]
\end{enumerate}
\end{definition}

\begin{remark}
Using standard arguments in control theory, one immediately observes that, if \eqref{eq:intro-linear-difference-equation} is $L^p$ approximately (respectively, exactly) controllable in time $T>0$, then it is $L^p$ approximately (respectively, exactly) controllable in any time $T^\prime \in [T, +\infty)$.
\end{remark}

We present below weaker notions of controllability, corresponding to those used by Yamamoto in \cite{Yamamoto1989Reachability}, in which the controllability time is not fixed in advance, but may depend on the target to be reached, and the initial state is always assumed to be the origin.

\begin{definition}\label{def:contr-from-origin} Let $p \in [1, +\infty]$. System~\eqref{eq:intro-linear-difference-equation} is said to be
\begin{enumerate}
\item \emph{$L^p$ approximately controllable (from the origin)} if, for every $\phi \in L^p([-\Lambda_N,0],\mathbb{R}^d)$ and $\epsilon>0$, there exist $T > 0$ and $u \in L^p([0,T], \mathbb{R}^m)$ such that the solution $x$ of \eqref{eq:intro-linear-difference-equation} with initial condition $0$ and control $u$ satisfies
\[
\norm{x_{T} - \phi}_{L^p([-\Lambda_N, 0], \mathbb R^d)} < \epsilon.
\]
\item \emph{$L^p$ exactly controllable (from the origin)} if, for every $\phi \in L^p([-\Lambda_N,0],\mathbb{R}^d)$, there exist $T > 0$ and $u \in L^p([0, T], \mathbb{R}^m)$ such that the solution $x$ of \eqref{eq:intro-linear-difference-equation} with initial condition $0$ and control $u$ satisfies
\[
x_{T} = \phi.
\]
\end{enumerate}
\end{definition}

Note that obvious implications hold true between the different notions provided in Definitions~\ref{def:Lq-cont} and \ref{def:contr-from-origin}: approximate (respectively, exact) controllability in time $T$ implies approximate (respectively, exact) controllability from the origin. One of the main results of the present chapter is that one actually has equivalence between such notions, as stated in the next theorem, whose proof is the subject of Section~\ref{sec:6}.

\begin{theorem}
\label{lem:3}
Let $p \in [1,+\infty]$. System~\eqref{eq:intro-linear-difference-equation} is $L^p$ approximately (respectively, exactly) controllable from the origin if and only if it is $L^p$ approximately (respectively, exactly) controllable in time $T = d \Lambda_N$.
\end{theorem}

Before stating our two main results on approximate and exact controllability of \eqref{eq:intro-linear-difference-equation}, let us introduce the matrix-valued holomorphic map $H\colon \mathbb C \to \mathcal M_d(\mathbb R)$ defined, for $s \in \mathbb C$, by
\begin{equation}\label{eq:defH}
H(s) = \Id_d- \sum_{j=1}^N e^{- s \Lambda_j} A_j.
\end{equation}
The holomorphic function $H$ arises when one considers solutions of \eqref{eq:intro-linear-difference-equation} with $u \in L^p(\mathbb R_+,\allowbreak \mathbb R^m)$ and an initial condition equal to $0$. Indeed, extending $x$ and $u$ by zero for all negative times, \eqref{eq:intro-linear-difference-equation} is satisfied for every $t \in \mathbb R$ and then, taking the two-sided Laplace transform in \eqref{eq:intro-linear-difference-equation}, we obtain that
\begin{equation*}
\widehat{x}(s)= \mathcal{H}(s) \widehat{u}(s),
\end{equation*}
for $s$ in some right-half plane of $\mathbb{C}$, where $ \mathcal{H}(s)=H(s)^{-1} B$ is the transfer function of \eqref{eq:intro-linear-difference-equation}. We notice that the Laplace transform of $x$ exists for all $s$ in a suitable right-half plane of $\mathbb{C}$ because the solutions of \eqref{eq:intro-linear-difference-equation} grow at most exponentially.

Let us now state our main results on approximate and exact controllability of \eqref{eq:intro-linear-difference-equation}, whose proofs are the main subject of Section~\ref{section_Hautus_criteria}, and are given in Sections~\ref{subsec:Approximate_contr_HY} and \ref{subsec:Exact_contr_HY}, respectively.

\begin{theorem}
\label{main_result1}
Let $p \in [1,+\infty)$. System~\eqref{eq:intro-linear-difference-equation} is $L^p$ approximately controllable in time $d \Lambda_N$ if and only if the two following conditions hold true.
\begin{enumerate}
\item \label{assumption1bis} $\rank\begin{pmatrix}H(s) & B\end{pmatrix} = d$ for every $s \in \mathbb C$.
\item \label{assumption2bis} $\rank\begin{pmatrix}A_N & B\end{pmatrix} = d$.
\end{enumerate}
\end{theorem}

\begin{theorem}
\label{main_result2} If \eqref{eq:intro-linear-difference-equation} is $L^1$ exactly controllable in time $d \Lambda_N$, then the two following conditions hold true.
\begin{enumerate}
\item \label{assumption1-L1} $\rank\begin{pmatrix}M & B\end{pmatrix} = d$ for every $M \in \overline{H(\mathbb{C})}$.
\item \label{assumption2-L1} $\rank\begin{pmatrix}A_N & B\end{pmatrix} = d$.
\end{enumerate}
\end{theorem}

\begin{remark}
Theorem~\ref{main_result2} provides only necessary conditions for $L^1$ exact controllability. These conditions turn out to be also sufficient for $L^1$ exact controllability, as proved in \cite{Fueyo2024Corona}. In this document, we will provide only the proof of Theorem~\ref{main_result2}, and we invite the interested reader to check directly \cite{Fueyo2024Corona} for the proof of the converse statement.
\end{remark}

Theorem~\ref{main_result1}\ref{assumption1bis} can be seen as a Hautus-type condition, and it can be equivalently reformulated as follows.

\begin{proposition}
\label{Prop2:main_result0}
The following three statements are equivalent.
\begin{enumerate}
\item $\rank\begin{pmatrix}H(s) & B\end{pmatrix} = d$ for every $s \in \mathbb{C}$.

\item\label{app-b} For every $s \in \mathbb{C}$, we have
\[
 \inf \left\{ \abs*{g^\transp H(s)} + \abs*{g^\transp B} \suchthat  g\in \mathbb{C}^d,\, \abs{g} = 1 \right\} > 0.
\]

\item\label{app-b'} For every $s \in \mathbb{C}$, we have
\[
\det \left(H(s) H(s)^\ast + B B^\ast\right) > 0.
\]
\end{enumerate}
\end{proposition}

Similarly, we can reformulate \ref{assumption1-L1} from Theorem~\ref{main_result2} as follows.

\begin{proposition}
\label{Prop2:main_result}
Assume that $\rank\begin{pmatrix}A_N & B\end{pmatrix} = d$. Then the following three statements are equivalent.
\begin{enumerate}
\item\label{exa-a} $\rank\begin{pmatrix}M & B\end{pmatrix} = d$ for every $M\in \overline{H(\mathbb{C})}$.

\item\label{exa-b} There exists $\alpha>0$ such that, for every $s\in \mathbb{C}$,
\[
 \inf \left\{\abs*{g^\transp H(s)} + \abs*{g^\transp B} \suchthat g\in \mathbb{C}^d,\, \abs{g} = 1\right\} \ge \alpha.
\]

\item\label{exa-b'} There exists $\alpha>0$ such that, for every $s \in \mathbb{C}$,
\[
\det \left(H(s) H(s)^\ast + B B^\ast\right) \ge \alpha.
\]
\end{enumerate}
\end{proposition}

The proofs of Propositions~\ref{Prop2:main_result0} and \ref{Prop2:main_result} are given in Section~\ref{subsec:useful_properties}. These propositions allow one to verify \ref{assumption1bis} from Theorem~\ref{main_result1} or \ref{assumption1-L1} from Theorem~\ref{main_result2} through the study of the regions of the complex plane where the real-valued function $s \mapsto \det \left(H(s) H(s)^\ast + B B^\ast\right)$ vanishes or is small. It turns out that one can easily provide a lower bound for this function when $\Real s$ is large enough or small enough, and hence, if one wishes to study controllability of \eqref{eq:intro-linear-difference-equation} through Proposition~\ref{Prop2:main_result0}\ref{app-b'} or Proposition~\ref{Prop2:main_result}\ref{exa-b'}, one is only left to study the function $s \mapsto \det \left(H(s) H(s)^\ast + B B^\ast\right)$ in a vertical strip of the complex plane. More precisely, we have the following result.

\begin{proposition}
\label{prop_zero_Q0}
There exist $\beta_1, \beta_2, \rho\in\mathbb{R}$ such that
\begin{equation}
\label{eq:ineg-H-H*}
\det \left(H(s) H(s)^\ast + B B^\ast\right) \ge \abs{\det \left(H(s)\right)}^2 \ge \rho > 0
\end{equation}
for all $s \in \mathbb C$ such that $\Real s < \beta_1$ or $\Real s > \beta_2$.
\end{proposition}

\begin{proof}
The first inequality in \eqref{eq:ineg-H-H*} is trivial. As for the second one, notice first that one can assume, with no loss of generality, that $H(\cdot)\not\equiv \Id_d$. Moreover, for $\Real s$ large enough, $H(s)$ is arbitrarily close to the identity. This implies that there exists $\beta_2 > 0$ such that $\abs{\det \left(H(s)\right)} \geq 1/2$ for $\Real s> \beta_2$.

On the other hand, the holomorphic function $h$ defined by $h(s) = \det \left(H(s)\right)$ can be written as
\begin{equation*}
h(s) = \widetilde{h}(s) + h_n e^{-s(\Lambda \cdot n)},\qquad s\in \mathbb{C},
\end{equation*}
where $h_n$ is a nonzero real number, $n$ is in $\mathbb{N}^N$ and satisfies $0 \le \abs{n}_1 \le N$, and $\widetilde{h}$ is a holomorphic function such that $\displaystyle\lim_{\Real s \to -\infty} \widetilde{h}(s) e^{s(\Lambda \cdot n)}=0$. We then deduce that $h_n e^{-s(\Lambda \cdot n)}$ is the dominant term of $h(s)$ as $\Real s \to -\infty$, and it goes to $\infty$ in modulus. This implies that, for every $\rho > 0$, there exists $\beta_1 \in \mathbb{R}$ such that $\abs{h(s)}^2 \ge \rho$ for $\Real s < \beta_1$.
\end{proof}

Theorem~\ref{main_result1} is a complete characterization of $L^p$ approximate controllability, providing a necessary and sufficient condition which can be seen as the counterpart for difference equations of Hautus controllability criterion. In particular, the property of $L^p$ approximate controllability in time $d \Lambda_N$ does not depend on $p \in [1, +\infty)$, in the sense that, if it holds for some $p \in [1, +\infty)$ then it holds for all of them. On the contrary, Theorem~\ref{main_result2} provides only a necessary condition for $L^1$ exact controllability, which turns out to be also sufficient, as proved in \cite{Fueyo2024Corona}. We expect in fact that the conditions from Theorem~\ref{main_result2} are also necessary and sufficient for $L^p$ exact controllability for any $p \in [1, +\infty)$ and, in that direction, we propose the following conjecture.

\begin{conjecture}
\label{conj:HY-exact}
Let $p \in [1,+\infty)$. System~\eqref{eq:intro-linear-difference-equation} is $L^p$ exactly controllable in time $d \Lambda_N$ if and only if the two following conditions hold true.
\begin{enumerate}
\item\label{assumption1-conject} $\rank\begin{pmatrix}M & B\end{pmatrix} = d$ for every $M \in \overline{H(\mathbb{C})}$.
\item\label{assumption2-conject} $\rank\begin{pmatrix}A_N & B\end{pmatrix} = d$.
\end{enumerate}
\end{conjecture}

This is motivated firstly by the fact that this conjecture actually holds true for $L^1$ exact controllability, as well as for $L^2$ exact controllability in the case $N = d = 2$ and $m = 1$, as it follows from the results of \cite{Chitour2020Approximate} (see Section~\ref{sec:applications} for details). In addition, we provide below a result giving ground to the above conjecture. For that purpose, we introduce the following compatibility condition for the existence of regular enough solutions of \eqref{eq:intro-linear-difference-equation}.

\begin{definition}\label{def:compatibility}
Let $k \in \mathbb N$ and $x_0 \in \mathcal C^k([-\Lambda_N, 0], \mathbb R^d)$. We say that $x_0$ is \emph{$\mathcal C^k$-admissible} for \eqref{eq:intro-linear-difference-equation} if, for every integer $\ell$ with $0 \leq \ell \leq k$, we have
\[
x_0^{(\ell)}(0) = \sum_{j=1}^N A_j x_0^{(\ell)}(-\Lambda_j).
\]
\end{definition}

\begin{theorem}
\label{thm:Ck}
Assume that \ref{assumption1-conject} and \ref{assumption2-conject} from Conjecture~\ref{conj:HY-exact} hold true. Then there exists a nonnegative integer $k$ such that, for every $p \in [1, +\infty)$ and $x_0, \phi \in \mathcal C^k([-\Lambda_N, 0], \mathbb R^d)$ $\mathcal C^k$-admissible for \eqref{eq:intro-linear-difference-equation}, there exists $u \in L^p([-\Lambda_N, 0], \mathbb R^m)$ such that the solution $x$ of \eqref{eq:intro-linear-difference-equation} with initial condition $x_0$ and control $u$ satisfies $x_{d \Lambda_N} = \phi$.
\end{theorem}

The proof of Theorem~\ref{thm:Ck} is provided in Section~\ref{subsec:Ck}.

\section{Properties of the endpoint map and consequences for controllability}
\label{sec:6}

The goal of this section is threefold. Firstly, we give a variation-of-constants formula for \eqref{eq:intro-linear-difference-equation} allowing us to write the solution operator as the sum of a flow operator and an endpoint map. As a consequence, we can give characterizations of the controllability of \eqref{eq:intro-linear-difference-equation} in terms of the range of the endpoint map. Secondly, through a Cayley--Hamilton theorem for multivariate polynomials, we show that the range of the endpoint map is constant from the time $d\Lambda_N$ and, finally, we deduce that the approximate (exact respectively) controllability of \eqref{eq:intro-linear-difference-equation} is equivalent to the approximate controllability (exact respectively) in finite time $d\Lambda_N$, thus providing a proof for Theorem~\ref{lem:3}.

\subsection{Definitions and preliminary remarks}

We recall in Proposition~\ref{prop_var_const_formula} the explicit representation formula for solutions of \eqref{eq:intro-linear-difference-equation}, often called variation-of-constants formula (and sometimes flow formula), which can be found in \cite[Proposition~2.4]{Chitour2020Approximate} and \cite[Proposition~2.7]{Mazanti2017Relative}. Before the statement of the proposition, let us give the following definitions.

\begin{definition}
\label{def_Xi_n}
Let $N$ and $d$ be positive integers and $A_1, \dotsc, A_N$ be matrices in $\mathcal M_d(\mathbb R)$. The family of matrices $\Xi_n \in \mathcal{M}_{d}(\mathbb{R})$, $n \in \mathbb{Z}^N$, is defined inductively as
\begin{equation*}
\Xi_n = \begin{dcases*}
0, & if $n \in \mathbb{Z}^N \setminus \mathbb{N}^N$,\\
\Id_d, & if  $n=0$,\\
\sum_{k=1}^N A_k \Xi_{n-e_k}, & if $n \in \mathbb{N}^N$ and $\abs{n}_1>0$,
\end{dcases*}
\end{equation*}
where we recall that $e_k$ denotes the $k$-th canonical vector of $\mathbb Z^N$ for $k \in \llbracket 1,N\rrbracket$, i.e., all the coordinates of $e_k$ are zero except for the $k$-th one, which is equal to one.
\end{definition}

\begin{definition}
\label{def_Upsilon_E}
Let $p \in [1, +\infty]$ and $T \in [0,+\infty)$.
\begin{enumerate}
\item The \emph{flow operator $\Upsilon_{p}(T)\colon L^p([-\Lambda_N,0],\mathbb{R}^d) \to L^p([-\Lambda_N,0],\mathbb{R}^d)$ of \eqref{eq:intro-linear-difference-equation}} is the operator defined by
\[
\left(\Upsilon_{p}(T)x_0 \right)(t)=\sum_{\substack{(n,j) \in \mathbb{N}^N \times \llbracket1,N\rrbracket \\ -\Lambda_j\le  T+t-\Lambda \cdot n <0}}  \Xi_{n-e_j} A_j x_0(T+t-\Lambda \cdot n),
\]
for $x_0 \in L^p([-\Lambda_N,0],\mathbb{R}^d)$ and $t \in [-\Lambda_N,0]$.

\item\label{item_def_endpoint} The \emph{endpoint operator $E_{p}(T)\colon L^p([0,T],\mathbb{R}^m) \to L^p([-\Lambda_N,0],\mathbb{R}^d)$ of \eqref{eq:intro-linear-difference-equation}} is the operator defined by
\[
\left(E_{p}(T)u \right)(t)= \sum_{\substack{n \in \mathbb{N}^N \\ \Lambda \cdot n \le T+t}}  \Xi_n B u(T+t-\Lambda \cdot n),
\]
for $u \in L^p([0,T],\mathbb{R}^m)$ and $t\in [-\Lambda_N,0]$.
\end{enumerate}
\end{definition}

The main interest of the above definitions is the following representation formula for the solutions of \eqref{eq:intro-linear-difference-equation}.

\begin{proposition}[Variation-of-constants formula]
\label{prop_var_const_formula}
Let $p \in [1, +\infty]$, $T \geq 0$, $x_0 \in L^p([-\Lambda_N,0],\allowbreak \mathbb{R}^d)$, and $u \in L^p([0,T],\mathbb{R}^m)$. Then the unique solution $x \in L^p([-\Lambda_N, T], \mathbb R^d)$ of \eqref{eq:intro-linear-difference-equation} with initial condition $x_0$ and control $u$ satisfies
\begin{equation}
\label{eq:variation_constant}
x_t = \Upsilon_{p}(t) x_0 + E_{p}(t) u, \qquad \text{ for every } t \in [0, T].
\end{equation}
\end{proposition}

The operator $\Upsilon_{p}(T)$ can be interpreted as the solution operator of \eqref{eq:intro-linear-difference-equation} with $B = 0$, mapping an initial condition $x_0$ to the history function of the corresponding solution at time $T$, $x_T$. The operator $E_p(T)$ represents instead the endpoint map of \eqref{eq:intro-linear-difference-equation}, associating, with each control $u$, the history function $x_T$ at time $T$ of the solution $x$ of \eqref{eq:intro-linear-difference-equation} with initial condition $0$ and control $u$. In this sense, \eqref{eq:variation_constant} can be seen as a simple consequence of the linearity of \eqref{eq:intro-linear-difference-equation}.

Using classical techniques in control theory, the notions of approximate and exact controllability from Definitions~\ref{def:Lq-cont} and \ref{def:contr-from-origin} can be recast in terms of the endpoint operator as follows.

\begin{proposition}
\label{prop:first_charact_control}
Let $p \in [1, +\infty]$ and $T > 0$.
\begin{enumerate}
\item System~\eqref{eq:intro-linear-difference-equation} is $L^p$ approximately controllable in time $T>0$ if and only if $\Ran E_{p}(T)$ is dense in $L^p([-\Lambda_N, 0], \mathbb R^d)$.

\item System~\eqref{eq:intro-linear-difference-equation} is $L^p$ approximately controllable from the origin if and only if
\begin{equation*}
\overline{\bigcup_{T \ge 0} \Ran E_{p}(T)}=L^p([-\Lambda_N,0],\mathbb{R}^d).
\end{equation*}

\item System~\eqref{eq:intro-linear-difference-equation} is $L^p$ exactly controllable in time $T>0$ if and only if
\[\Ran E_{p}(T) = L^p([-\Lambda_N, 0], \mathbb R^d).\]

\item System~\eqref{eq:intro-linear-difference-equation} is $L^p$ exactly controllable from the origin if and only if
\begin{equation*}
\bigcup_{T \ge 0} \Ran E_{p}(T)=L^p([-\Lambda_N,0],\mathbb{R}^d).
\end{equation*}
\end{enumerate}
\end{proposition}

\subsection{Saturation of the range of the endpoint map in time \texorpdfstring{$d\Lambda_N$}{d Lambda N}}

In the case $N = 1$, i.e., when \eqref{eq:intro-linear-difference-equation} contains a single delay, an instrumental tool in the study of the controllability properties of \eqref{eq:intro-linear-difference-equation} is the Cayley--Hamilton theorem (see \cite[Remark~3.5]{Chitour2020Approximate}). It turns out that the following generalization of Cayley--Hamilton theorem in the case of multivariate polynomials plays a similar role in our subsequent arguments.

\begin{lemma}
\label{lem:1}
Let $\Xi_n$, $n\in \mathbb{Z}^N$, be the matrices introduced in Definition~\ref{def_Xi_n}. There exist real coefficients $\alpha_k$, for $k\in \mathbb{N}^N$ with $0<\abs{k}_1 \le d$, such that, for every $n\in \mathbb{N}^N$ with $\abs{n}_1\ge d$,
\begin{equation*}
\Xi_n=-\sum_{\substack{k \in \mathbb{N}^N \\ 0<\abs{k}_1 \le d}}\alpha_k \Xi_{n-k}.
\end{equation*}
\end{lemma}

\begin{proof}
For $t=(t_1,\dotsc,t_N) \in \mathbb{R}^N$, set
\[
A(t) = t_1 A_1 + t_2 A_2 + \dotsb + t_N A_N.
\]
One deduces, by Definition~\ref{def_Xi_n} and an immediate induction argument that, for every $j\in \mathbb{N}$ and $t \in \mathbb{R}^N$, we have
\begin{equation}
\label{eq:lem1:1.1}
A(t)^j = \sum_{\substack{n \in \mathbb{N}^N \\ \abs{n}_1 = j}} \Xi_n t^n,
\end{equation}
where $t^n$ is defined by $t^n=t_1^{n_1}  t_2^{n_2} \dotsm t_N^{n_N}$. Using Neumann series, we deduce from \eqref{eq:lem1:1.1} that, for $t$ small enough,
\begin{equation}
\label{eq:lem1:2}
\bigl(\Id_d-A(t)\bigr)^{-1}=\sum_{j\in \mathbb{N}} A(t)^j = \sum_{n\in \mathbb{N}^N} \Xi_{n}t^n.
\end{equation}
Notice that $P(t)=\det\bigl(\Id_d-A(t)\bigr)$ is a multivariate polynomial of degree $d$, that is,
\begin{equation}
\label{eq:lem1:3}
P(t) = \sum_{\substack{k \in \mathbb{N}^N \\ 0 \le \abs{k}_1 \le d}} \alpha_k t^k,
\end{equation}
for some real numbers $\alpha_k$ defined for $k \in \mathbb{N}^N$ such that $\abs{k}_1 \le d$ and with $\alpha_0=1$. We have that, on the one hand, there exist $M_k \in \mathcal{M}_{d, d}(\mathbb{R})$ for $k \in \mathbb{N}^N$ with $0 \le\abs{k}_1\le d-1$ such that
\begin{equation}
\label{eq:lem1:4}
\Adj\bigl(\Id_d-A(t)\bigr) = \sum_{\substack{k \in \mathbb{N}^N \\ 0 \le \abs{k}_1 \le d-1}} M_{k} t^k,
\end{equation}
where we recall that $\Adj(B)$ denotes the adjugate of the square matrix $B$; and, on the other hand, \eqref{eq:lem1:2} and \eqref{eq:lem1:3} lead to
\begin{equation}
\begin{split}
\Adj\bigl(\Id_d-A(t)\bigr) & = P(t)\bigl(\Id_d-A(t)\bigr)^{-1} = \sum_{\substack{k \in \mathbb{N}^N \\ 0 \le \abs{k}_1 \le d}} \alpha_{k} t^k\sum_{n \in \mathbb{N}^N} \Xi_{n}t^n
\\
& = \sum_{n \in \mathbb{N}^N} \sum_{\substack{k \in \mathbb{N}^N \\ 0 \le \abs{k}_1 \le d}} \alpha_k \Xi_{n} t^{n+k}.
\end{split}
\label{eq:lem1:5}
\end{equation}
The substitution $\ell = n+k$ in \eqref{eq:lem1:5} allows us to write
\begin{equation}
\label{eq:lem1:6}
\Adj\bigl(\Id_d-A(t)\bigr)= \sum_{n \in \mathbb{N}^N} \sum_{\substack{\ell \in \mathbb{N}^N,\,\ell - n \in \mathbb{N}^N \\ 0 \le \abs{\ell - n}_1 \le d}} \alpha_{\ell - n} \Xi_{n} t^{\ell}.
\end{equation}
We deduce from \eqref{eq:lem1:4}--\eqref{eq:lem1:6} that, for $\ell \in \mathbb{N}^N$ and $\abs{\ell}_1 \ge d$,
\begin{equation}
\label{eq:lem1:7}
\sum_{\substack{n \in \mathbb{N}^N,\,\ell-n \in \mathbb{N}^N \\ 0 \le \abs{\ell-n}_1 \le d}} \alpha_{\ell-n} \Xi_{n} =0.
\end{equation}
Setting $n^\prime=\ell-n$ in \eqref{eq:lem1:7}, one obtains
\begin{equation}
\label{eq:lem1:8}
\sum_{\substack{n^\prime \in \mathbb{N}^N,\,\ell-n^\prime \in \mathbb{N}^N \\ 0 \le \abs{n^\prime}_1 \le d}} \alpha_{n^\prime} \Xi_{\ell-n^\prime} =0.
\end{equation}
Since $\Xi_{\ell-n^\prime}=0$ for $\ell-n^\prime \in \mathbb{Z}^N \setminus \mathbb{N}^N$ (see Definition~\ref{def_Xi_n}), we deduce from \eqref{eq:lem1:8} that
\begin{equation*}
\Xi_{\ell}=-\sum_{\substack{n^\prime \in \mathbb{N}^N \\ 0 < \abs{n^\prime}_1 \le d}} \alpha_{n^\prime} \Xi_{\ell-n^\prime},
\end{equation*}
yielding the conclusion.
\end{proof}

We now use Lemma~\ref{lem:1} to prove that the range of the operator $E_{p}(T)$ is constant with respect to $T$ for $T \ge d\Lambda_N$.

\begin{theorem}
\label{lem:RanE_con}
Let $p \in [1, +\infty]$ and consider the endpoint operator from Definition~\ref{def_Upsilon_E}\ref{item_def_endpoint}. For every $T\in [d \Lambda_N,+\infty)$, we have
\begin{equation}
\label{eq:RanE_con1}
\Ran E_{p}(T) = \Ran E_{p}(d\Lambda_N).
\end{equation}
\end{theorem}

\begin{proof}
The proof is divided in two steps. We first prove that \eqref{eq:RanE_con1} is satisfied for $T \in [d \Lambda_N,d\Lambda_N+\delta]$ and some $\delta > 0$. We then deduce from the flow formula \eqref{eq:variation_constant} that \eqref{eq:RanE_con1} is actually satisfied for all $T \in [d \Lambda_N,+\infty)$.

Let $T > d \Lambda_N$ and $u \in L^{p}([0,T], \mathbb{R}^m)$. We define the functions $u_1$ and $u_2$ in $L^p([0, d \Lambda_N], \mathbb R^m)$ by setting, for $t \in [0, d \Lambda_N]$,
\begin{align}
u_1(t) & = u(t+T-d\Lambda_N), \\
u_2(t) & = - \sum_{\substack{0 < \abs{k}_1 \le d,\, k \in \mathbb{N}^N \\ t<\Lambda\cdot k \le t+T-d \Lambda_N }} \alpha_k u(t - \Lambda \cdot k+T-d \Lambda_N), \label{eq:RanE_con3}
\end{align}
where the coefficients  $\alpha_k$ are those from Lemma~\ref{lem:1} and the sum in \eqref{eq:RanE_con3} is understood to be zero when the indices are taken in an empty set. For $t \in [-\Lambda_N,0]$, we have
\begin{equation}\label{eq:RanE_con4}
\begin{split}
\left(E_{p}(T)u \right)(t)&= \sum_{\substack{n \in \mathbb{N}^N \\ \Lambda \cdot n \le {T+t}}}  \Xi_n B u(T+t-\Lambda \cdot n) \\
&=\sum_{\substack{n \in \mathbb{N}^N \\ \Lambda \cdot n \le d \Lambda_N+t}}  \Xi_n B u(T+t-\Lambda \cdot n)
+\sum_{\substack{n \in \mathbb{N}^N \\ d \Lambda_N+t<\Lambda \cdot n \le T+t}}  \Xi_n B u(T+t-\Lambda \cdot n).
\end{split}
\end{equation}
Notice that
\begin{equation}
\label{eq:RanE_con5}
\sum_{\substack{n \in \mathbb{N}^N \\ \Lambda \cdot n \le d \Lambda_N+t}}  \Xi_n B u(T+t-\Lambda \cdot n) = \left(E_{p}(d\Lambda_N)u_1\right)(t).
\end{equation}

Since $d \Lambda_N+t<\Lambda \cdot n$ for $n \in \mathbb{N}^N$ implies that $\abs{n}_1\ge d$,  we deduce from
 Lemma~\ref{lem:1} that
\begin{equation}
\label{eq:RanE_con6}
\sum_{\substack{n \in \mathbb{N}^N \\ d \Lambda_N+t<\Lambda \cdot n \le T+t}}  \Xi_n B u(T+t-\Lambda \cdot n) = -\sum_{\substack{n \in \mathbb{N}^N \\ d \Lambda_N+t<\Lambda \cdot n \le T+t}} \sum_{\substack{k \in \mathbb{N}^N \\ 0<\abs{k}_1 \le d }}\alpha_k \Xi_{n-k} B u(T+t-\Lambda \cdot n).
\end{equation}
The substitution $n^\prime=n-k$ in \eqref{eq:RanE_con6} yields
\begin{equation}
\label{eq:RanE_con7}
\sum_{\substack{n \in \mathbb{N}^N \\ d \Lambda_N+t<\Lambda \cdot n \le T+t}}  \Xi_n B u(T+t-\Lambda \cdot n)
=-\sum_{\substack{k \in \mathbb{N}^N \\ 0<\abs{k}_1 \le d }}
\sum_{\substack{n^\prime \in \mathbb{N}^N \\ d \Lambda_N+t<\Lambda \cdot (n^\prime+k) \le T+t}}  \alpha_k \Xi_{n^\prime} B u(T+t-\Lambda \cdot (n^\prime+k))
.
\end{equation}
Let $\delta>0$ be such that, for all $T\in [d \Lambda_N,d \Lambda_N+\delta]$, $t \in [-\Lambda_N,0]$, $k \in \mathbb{N}^N$ with $0<\abs{k}_1\le d$,
and $n^\prime \in \mathbb{N}^{\mathbb{N}}$, if $d \Lambda_N+t<\Lambda \cdot (n^\prime+k) \le T+t$, then $\Lambda\cdot n^\prime \le d \Lambda_N+t$. Letting $T\in [d \Lambda_N,d \Lambda_N+\delta]$, we can thus rewrite \eqref{eq:RanE_con7} as
\begin{equation}
\label{eq:RanE_con8}
\begin{split}
\MoveEqLeft[6] \sum_{\substack{n \in \mathbb{N}^N \\ d \Lambda_N+t<\Lambda \cdot n \le T+t}}  \Xi_n B u(T+t-\Lambda \cdot n)  \\
& = -\sum_{\substack{n^\prime \in \mathbb{N}^N \\ \Lambda \cdot n^\prime \le d \Lambda_N+t}}\Xi_{n^\prime}B
\sum_{\substack{k \in \mathbb{N}^N,\, 0<\abs{k}_1 \le d \\ d \Lambda_N+t<\Lambda \cdot (n^\prime+k) \le T+t}} \alpha_k  u(T+t-\Lambda \cdot (n^\prime+k))
\\
& = \left(E_{p}(d \Lambda_N)u_2\right)(t).
\end{split}
\end{equation}
Equations~\eqref{eq:RanE_con4}, \eqref{eq:RanE_con5}, and \eqref{eq:RanE_con8} prove that, for $T \in [d\Lambda_N,d\Lambda_N+\delta]$ and
$t \in [-\Lambda_N,0]$,
\begin{equation}
\label{eq:RanE_con9}
\left(E_{p}(T)u\right)(t) = \left(E_{p}(d\Lambda_N)u_1\right)(t) + \left(E_{p}(d \Lambda_N)u_2\right)(t).
\end{equation}
From \eqref{eq:RanE_con9}, we deduce that
\begin{equation}
\label{eq:RanE_con10}
\Ran E_{p}(T)=\Ran E_{p}(d\Lambda_N) \qquad \text{ for every } T \in [d \Lambda_N,d\Lambda_N+\delta].
\end{equation}

Let us now extend \eqref{eq:RanE_con10} to all  $T \in [d \Lambda_N,+\infty)$. Let $V=\Ran E_{p}(d \Lambda_N)$ and $x \in V$. Fix  $u\in L^p([0,d\Lambda_N],\mathbb{R}^m)$ such that $x=E_{p}(d\Lambda_N)u$. For $t\in [0,\delta]$, define $\widetilde{u}\in L^p([0,d\Lambda_N+t],\mathbb{R}^m)$ by setting $\widetilde{u}\rvert_{[0,d\Lambda_N]}=u$ and ${\widetilde{u}}(t^\prime)=0$ for $t^\prime \in [d \Lambda_N,d\Lambda_N+t]$. From the variation-of-constants formula \eqref{eq:variation_constant}, we have
\begin{equation*}
\Upsilon_{p}(t) x= E_{p}(d\Lambda_N+t)\widetilde{u} \in \Ran E_{p}(d\Lambda_N+t).
\end{equation*}
Thanks to \eqref{eq:RanE_con10}, we have proved that
\begin{equation}\label{eq:invariance}
\Upsilon_{p}(t)x \in V\qquad \text{ for every } t\in[0,\delta] \text{ and } x\in V.
\end{equation}

Let $T \in [d\Lambda_N+\delta,d\Lambda_N+2 \delta]$, $y \in \Ran E_{p}(T)$, and $u \in L^p([0,T],\mathbb{R}^m)$ such that $y=E_{p}(T)u$. Define $z=E_{p}(d\Lambda_N+\delta)u\rvert_{[0,d\Lambda_N+\delta]} \in V$. The variation-of-constants formula \eqref{eq:variation_constant} gives
\begin{equation*}
y=\Upsilon_{p}(T-d\Lambda_N-\delta)z+E_{p}(T-d\Lambda_N-\delta)\check{u},
\end{equation*}
where $\check u$ is defined by $\check{u}(\alpha)=u(\alpha+d\Lambda_N+\delta)$ for $\alpha \in [0,T-d\Lambda_N-\delta]$. We deduce from \eqref{eq:RanE_con10} and \eqref{eq:invariance} that $y \in \Ran E_{p}(d\Lambda_N)$, proving that
\begin{equation}
\label{eq_lem3:5}
\Ran E_{p}(d \Lambda_N)=\Ran E_{p}(T), \qquad T \in [d\Lambda_N+\delta,d \Lambda_N+2 \delta].
\end{equation}
The iteration of the same process proves that \eqref{eq_lem3:5} actually holds for all $T\ge d \Lambda_N$.
\end{proof}

\subsection{Upper bound on the minimal control time}

Thanks to Theorem~\ref{lem:RanE_con}, we are now in position to provide a simple proof of Theorem~\ref{lem:3}.

\begin{proof}[Proof of Theorem~\ref{lem:3}]
By Theorem~\ref{lem:RanE_con} and since $T\mapsto \Ran E_{p}(T)$ is monotone nondecreasing with respect to set inclusion, we have that
\[\bigcup_{T \ge 0} \Ran E_{p}(T)=\Ran E_{p}(d\Lambda_N).\]
The conclusion follows directly by Proposition~\ref{prop:first_charact_control}.
\end{proof}

One of the major questions concerning the approximate and exact controllability in finite time $T$ is to determine the \emph{minimal time} of controllability, in the sense of the following definition.

\begin{definition}
\label{def:Tmin-ap-ex-p}
Let $p \in [1, +\infty]$. We define the \emph{minimal times} of approximate and exact controllability, $T_{\min,\,\mathrm{ap},\,p}$ and $T_{\min,\,\mathrm{ex},\,p}$, respectively, by
\begin{align*}
T_{\min,\,\mathrm{ap},\,p}&=\inf \{T \in \mathbb{R}_+ \suchthat \text{System~\eqref{eq:intro-linear-difference-equation} is $L^p$ approximately controllable in time $T$}\},\\
T_{\min,\,\mathrm{ex},\,p}&=\inf \{T \in \mathbb{R}_+ \suchthat \text{System~\eqref{eq:intro-linear-difference-equation} is $L^p$ exactly controllable in time $T$} \},
\end{align*}
with the convention that $\inf\emptyset = +\infty$.
\end{definition}

Since an immediate inspection of \eqref{eq:intro-linear-difference-equation} shows that it is never approximately or exactly controllable before time $\Lambda_N$, we can recast Theorem~\ref{lem:3} in terms of minimal time of controllability as follows.

\begin{corollary}
\label{Corollary_min_con}
Let $p \in [1, +\infty]$. Both times $T_{\min,\,\mathrm{ap},\,p}$ and $T_{\min,\,\mathrm{ex},\,p}$ belong to the set  $[\Lambda_N,d \Lambda_N] \cup \{+\infty\}$.
\end{corollary}

In the case $N=d=2$, $m=1$, and $p = 2$, it follows from the results of \cite{Chitour2020Approximate} that the pair $(T_{\min,\,\mathrm{ap},\,2}, T_{\min,\,\mathrm{ex},\,2})$ is equal to $(+\infty, +\infty)$, $(2 \Lambda_N, +\infty)$, or $(2 \Lambda_N, 2 \Lambda_N)$. In the general case, however, $d \Lambda_N$ is not always minimal, as proved in \cite{Chitour2020Approximate} by looking at the case of commensurable delays. It was noticed in the same reference that $d \Lambda_N$ is minimal when controllability holds for systems with a single input and commensurable delays. We conjecture that this result holds true also for the case of not necessarily commensurable delays. It would be also interesting to investigate the question of equality between $T_{\min,\,\mathrm{ap},\,p}$ and $T_{\min,\,\mathrm{ex},\,p}$ when both of them are finite.

\section{Hautus--Yamamoto criteria for approximate and exact controllability}
\label{section_Hautus_criteria}

Our goal is now to analyze approximate and exact controllability of \eqref{eq:intro-linear-difference-equation} by using Y.~Yamamoto's controllability results based on realization theory from \cite{Yamamoto1978Realization, Yamamoto1981RealizationI, Yamamoto1981RealizationII, Yamamoto1988Pseudo, Yamamoto1989Reachability, Yamamoto2007Coprimeness, Yamamoto2011Bezout, Yamamoto2008Behavioral, Yamamoto2016Behavioral, Yamamoto1999Coprimeness}. We will consider from now on only the case $p \in [1, +\infty)$, since the arguments we will use do not readily extend to $L^\infty$.

\subsection{Input/output map for (\ref{eq:intro-linear-difference-equation})}
\label{subsec:real_theory}

When the initial condition is fixed at $0$ (which is the case in Definition~\ref{def:contr-from-origin}), \eqref{eq:intro-linear-difference-equation} describes a linear relation between the control $u$ on $[0, T]$ and the history function $x_T$ of the state $x$ at time $T$, described by the operator $E_p(T)$ introduced in Definition~\ref{def_Upsilon_E}. We now aim at rewriting this relation as an input/output map of the form \eqref{eq:intro-input-output}.

Notice first that \eqref{eq:intro-input-output} considers that inputs are applied during negative times, and the output is measured for positive times. Hence, we first rewrite \eqref{eq:intro-linear-difference-equation}, after suitable time translations, as
\begin{equation}
\label{syst_lin_avec_sortie}
\begin{dcases*}
x(t) = \sum_{j=1}^N A_j x(t-\Lambda_j)+Bu(t), & for $t\ge \inf \supp u$, \\
x(t) = 0, & for $t < \inf \supp u$, \\
y(t) = x(t-\Lambda_N), & for $t \in [0,+\infty)$,
\end{dcases*}
\end{equation}
where $u$ and $x$ are translated versions of the corresponding variables from \eqref{eq:intro-linear-difference-equation} and $y$ is the output. The choice of the output $y$ is justified by the fact that, if $u$ is a control in \eqref{eq:intro-linear-difference-equation} applied during a time interval $[0, T]$, then the corresponding control $u$ in \eqref{syst_lin_avec_sortie} is applied during the time interval $[-T, 0]$, and the restriction of $y$ to the interval $[0, \Lambda_N]$ coincides, up to a time translation, with the history function $x_T$ of the state $x$ of \eqref{eq:intro-linear-difference-equation}.

We fix $p \in [1, +\infty)$ and we consider that the input $u$ belongs to the space $\Omega_p$ defined by
\begin{equation*}
\Omega_{p}=\{u\in L^p(\mathbb{R}_-,\mathbb{R}^m) \suchthat \supp u \text{ is compact}\}.
\end{equation*}
As for the output $y$, we are interested in its values in the interval $[0, \Lambda_N]$, but an in\-put/out\-put map of the form \eqref{eq:intro-input-output} defines an output $y$ for all positive times. On the other hand, since the support of $u$ is included in $\mathbb R_-$, \eqref{syst_lin_avec_sortie} implies that
\[
x(t) = \sum_{j = 1}^N A_j x(t - \Lambda_j) \qquad \text{ for } t > 0,
\]
and thus, using the definition of the output $y$, we have
\begin{equation}
\label{eq:output-after-LambdaN}
y(t) = \sum_{j = 1}^N A_j y(t - \Lambda_j) \qquad \text{ for } t > \Lambda_N.
\end{equation}
In other words, $y$ is uniquely determined by its values on $[0, \Lambda_N]$ since, if it is known in this interval, it can be extended to $\mathbb R_+$ by using \eqref{eq:output-after-LambdaN}. This motivates the definition of the output space
\[
X^{Q,\,p} = \left\{y \in L^p_{\loc}\left(\mathbb{R}_+,\mathbb{R}^d\right)\suchthat y(t) = \sum_{j = 1}^N A_j y(t - \Lambda_j) \text{ for } t > \Lambda_N  \right\}.
\]
Note that, by construction, $X^{Q,\,p}$ can be identified with the space $L^p\left(\left[0,\Lambda_N\right],\mathbb{R}^d\right)$: we have $y \in  X^{Q,\,p}$ if and only if the restriction $y\rvert_{[0,\Lambda_N]}$ is in $L^p([0,\Lambda_N],\mathbb{R}^d)$ and $y$ is the unique extension of $y\rvert_{[0,\Lambda_N]}$ to $\mathbb R_+$ satisfying \eqref{eq:output-after-LambdaN}.

Let us define the Radon measures $Q$ and $P$ in $M(\mathbb R_-)$ by 
\begin{align}
Q & =\delta_{-\Lambda_N} \Id_d - \sum_{j=1}^N \delta_{-\Lambda_N+\Lambda_j} A_j, \label{eq:defQ} \\
P & = B \delta_0. \label{eq:defP}
\end{align}
Note that \eqref{eq:output-after-LambdaN} can be rewritten in a simpler manner using $Q$ as $\pi(Q*y)=0$ and thus the space $X^{Q,\,p}$ can be rewritten as
\begin{equation*}
X^{Q,\,p}=\left\{y \in L^p_{\loc}\left(\mathbb{R}_+,\mathbb{R}^d\right)\suchthat \pi(Q*y)=0  \right\}.
\end{equation*}

The distributions $Q$ and $P$ allow us to obtain the Radon measure $A$ representing \eqref{syst_lin_avec_sortie} as a convolution operator of the form \eqref{eq:intro-input-output}. Denoting by $\widetilde{y}$ the natural extension of the output $y$ on $\mathbb{R}$, i.e., $\widetilde{y}(t)=x(t-\Lambda_N)$ for $t \in \mathbb{R}$, the first equation of \eqref{syst_lin_avec_sortie} implies that
\begin{equation}
\label{obtention_convolution1}
\left(Q*\widetilde{y}\right)(t)=\left(P*u\right)(t),\qquad t \in \mathbb{R}.
\end{equation}
The distribution $Q$ is invertible over $\mathcal{D}^\prime_+(\mathbb{R})$ (with respect to the convolution product in the algebra $\mathcal D^\prime_+(\mathbb R)$) and its inverse distribution $Q^{-1}$ belongs to $M_+(\mathbb{R})$. More precisely, the distribution
\begin{equation}
\label{eq:explicit-Q-1}
Q^{-1} = \delta_{\Lambda_N}*  \left( \sum_{n=0}^{+\infty} \left(   \sum_{j=1}^N \delta_{\Lambda_j} A_j \right)^n \right)
\end{equation}
is easily seen to be the inverse of $Q$ by a Neumann series argument. We take the convolution product of \eqref{obtention_convolution1} on the left by $Q^{-1}$ and we obtain
\begin{equation}
\label{obtention_convolution2}
\widetilde{y}(t)=\left(Q^{-1}*P*u\right)(t),\qquad t \in \mathbb{R}.
\end{equation}
Applying the operator $\pi$ in \eqref{obtention_convolution2}, we obtain
\begin{equation}
\label{representation_radon_measure3}
y = \pi \left(Q^{-1}*P*u\right),
\end{equation}
which is precisely under the form \eqref{eq:intro-input-output} with $A = Q^{-1} * P$. Notice also that the measure $A$ belongs to $M_+(\mathbb{R})$ and, as a consequence, \eqref{syst_lin_avec_sortie} is pseudo-rational, as defined in Section~\ref{sec:Realization}.

\begin{remark}
Equation~\eqref{eq:intro-input-output} with $A = Q^{-1} * P$ provides an expression for the relation between the input $u$ and the output $y$ of \eqref{syst_lin_avec_sortie}. Other expressions for this same relation can be obtained by other means, such as by relying on the operator $E_{p}(T)$, as done in Section~\ref{sec:6}, or by expressing $y$ as a Stieltjes integral of $u$ with respect to the fundamental solution of the system, as done in the variation-of-constants formula in \cite{Hale1993Introduction}. While clearly equivalent, some representations can be more suitable than others for a given purpose. The representation through the operator $E_{p}(T)$ was useful in Section~\ref{sec:6} to prove Theorem~\ref{lem:3} while, in the remaining part of the section, we will use \eqref{eq:intro-input-output} to obtain the controllability results from Theorems~\ref{main_result1} and \ref{main_result2}.
\end{remark}

We now characterize the controllability notions from Definition~\ref{def:contr-from-origin} in terms of the above realization theory formalism.

\begin{proposition}
\label{def_reachability}
Let $p \in [1, +\infty)$.
\begin{enumerate}
\item System~\eqref{eq:intro-linear-difference-equation} is $L^p$ approximately controllable if and only if, for every $\phi \in  X^{Q,\,p}$, there exists a sequence of inputs $(u_n)_{n \in \mathbb{N}} \in (\Omega_p)^{\mathbb{N}}$ such that its associated sequence of outputs $(y_n)_{n \in \mathbb{N}} \in \left(X^{Q,\,p}\right)^{\mathbb{N}}$ through \eqref{syst_lin_avec_sortie} satisfies
\[
y_n \xrightarrow[n \to +\infty]{} \phi \qquad \text{ in } L^p_{\loc}\left(\mathbb{R}_+,\mathbb{R}^d\right).
\]
\item System~\eqref{eq:intro-linear-difference-equation} is $L^p$ exactly controllable if and only if, for every $\phi \in  X^{Q,\,p}$, there exists $u \in \Omega_p$ such that
its associated output $y$ through \eqref{syst_lin_avec_sortie} satisfies
\[
y = \phi.
\]
\end{enumerate}
\end{proposition}

\begin{remark}
Following Yamamoto (see, e.g., \cite{Yamamoto1989Reachability}), we will refer in the sequel to the characterization of approximate and exact controllability of \eqref{eq:intro-linear-difference-equation} given in Proposition~\ref{def_reachability} as \emph{approximate and exact controllability of \eqref{syst_lin_avec_sortie}}. Note also that, in \cite{Yamamoto1989Reachability}, approximate (respectively, exact) controllability of \eqref{eq:intro-linear-difference-equation} is also called \emph{quasi-rea\-cha\-bi\-lity} (respectively, \emph{reachability}).
\end{remark}

\subsection{Preliminary properties and proofs of Propositions~\ref{Prop2:main_result0} and \ref{Prop2:main_result}}
\label{subsec:useful_properties}

Before making use of the realization theory to provide criteria for approximate and exact controllabilities in Sections~\ref{subsec:Approximate_contr_HY} and \ref{subsec:Exact_contr_HY}, respectively, we provide in this section preliminary technical results relating properties of the (Laplace transform of the) distributions $Q$ and $P$ introduced in \eqref{eq:defQ} and \eqref{eq:defP} with $H$ associated with \eqref{eq:intro-linear-difference-equation} from \eqref{eq:defH}. In particular, we shall also obtain from those technical results the proofs of Propositions~\ref{Prop2:main_result0} and \ref{Prop2:main_result}.

We start by considering the Laplace transforms of the distributions $Q$ and $P$, which, using \eqref{laplace_transform_radon_measure} and \eqref{laplace_transform_dirac}, are given for $s \in \mathbb C$ by
\begin{equation*}
\widehat{Q}(s)=e^{s \Lambda_{N}}\Id_d- \sum_{j=1}^N e^{ s(\Lambda_N- \Lambda_j)} A_j \quad \text{ and }\quad \widehat{P}(s) = B.
\end{equation*}
Since $\widehat P$ is constant, we omit its argument $s$ in the sequel. Note that the function $H$ from \eqref{eq:defH} and the Laplace transform $\widehat{Q}$ of $Q$ satisfy the relation
\begin{equation}
\label{lien_entre_les_deux_fonctions}
H(s)=e^{-s \Lambda_N} \widehat{Q}(s),\qquad \text{ for every } s \in \mathbb{C}.
\end{equation}

Our next aim is to provide a link between the rank condition in Theorems~\ref{main_result1}\ref{assumption1bis} and \ref{main_result2}\ref{assumption1-L1} with a similar condition on $\widehat Q$. As a preliminary step, let us prove the following technical lemma.

\begin{lemma}
\label{lem:sous_suite1}
Let $H$ be the function defined in \eqref{eq:defH} and $\widehat Q$ be the Laplace transform of the distribution $Q$ defined in \eqref{eq:defQ}.
\begin{enumerate}
\item\label{item:suites-1} There exists $M\in \overline{H(\mathbb{C})}$ such that $\rank\begin{pmatrix}M & B\end{pmatrix} < d$ if and only if there exist $g \in \mathbb{C}^d$ and $(s_n)_{n \in \mathbb{N}} \in \mathbb{C}^{\mathbb{N}}$ with bounded real part such that $\abs{g}=1$, $g^\transp B=0$, and $\lim_{n\to\infty}g^\transp H(s_n)=0$.

\item\label{item:suites-2} Under the assumption that $\rank\begin{pmatrix}A_N & B\end{pmatrix} = d$, there exists $\widetilde{Q}\in \overline{\widehat{Q}(\mathbb{C})}$ such that
$\rank\begin{pmatrix}\widetilde{Q} & B\end{pmatrix} \allowbreak <d$ if and only if there exist $g \in \mathbb{C}^d$ and $(s_n)_{n \in \mathbb{N}} \in \mathbb{C}^{\mathbb{N}}$ with bounded real part such that $\abs{g}=1$, $g^\transp B=0$, and $\lim_{n\to\infty} g^\transp\widehat{Q}(s_n) = 0$.

\item\label{item:suites-3} Under the assumption that $\rank\begin{pmatrix}A_N & B\end{pmatrix} = d$, Proposition~\ref{Prop2:main_result}\ref{exa-b} is not satisfied if and only if there exist $g\in \mathbb{C}^d$ and $(s_n)_{n \in \mathbb N} \in \mathbb{C}^{\mathbb{N}}$ with bounded real part such that $\abs{g}=1$, $g^\transp B=0$, and $\lim_{n\to\infty} g^\transp H(s_n) = 0$.
\end{enumerate}
\end{lemma}

\begin{proof}
We start by proving \ref{item:suites-1}.  Let $g\in \mathbb{C}^d$ and $(s_n)_{n\in \mathbb{N}}$ with bounded real part be such that $\abs{g}=1$, $g^\transp B=0$ and $\lim_{n\to\infty}g^\transp H(s_n) =0$. Since $H(\cdot)$ is uniformly bounded on any bounded vertical strip, it follows that, up to a subsequence, $H(s_n)$ converges to some matrix $M \in \mathcal{M}_{d,d}(\mathbb{C})$ as $n\to\infty$. We deduce that $g^\transp M=0$ and $g^\transp B=0$, proving that $\rank\begin{pmatrix}M & B\end{pmatrix} < d$. Conversely, assume that $g \in \mathbb{C}^d$ and $M$ are such that $\abs{g}=1$, $g^\transp B=0$, $g^\transp M=0$, and $M=\lim_{n\to\infty}H(s_n)$ for some sequence  $(s_n)_{n \in \mathbb{N}} \in \mathbb{C}^{\mathbb{N}}$. The sequence $(s_n)_{n\in \mathbb{N}}$ has bounded real part since \begin{enumerate*}[label={(\roman*)}] \item $H(s)$ is nonsingular for $s$ out of a bounded vertical strip (see Proposition~\ref{prop_zero_Q0}), \item $H(s)$ converges to $\Id_d$ when the real part of $s$ tends to $+\infty$, and \item $H(s)$ diverges in norm when the real part of $s$ tends to $-\infty$  (cf.~the proof of Proposition~\ref{prop_zero_Q0}).\end{enumerate*} This concludes the proof of the converse implication.

The proof of \ref{item:suites-2} is similar to that of \ref{item:suites-1}, with $\widehat Q$ playing the role of $H$. The only difference is that, in the second part of the argument, one needs to use the assumption that $\rank\begin{pmatrix}A_N & B\end{pmatrix} = d$ in order to ensure that the sequence $(s_n)_{n \in \mathbb{N}} \in \mathbb{C}^{\mathbb{N}}$ has bounded real part, since $\widehat Q(s)$ converges to $-A_N$ (instead of $\Id_d$) when the real part of $s$ tends to $-\infty$ (instead of $+\infty$).

Let us prove \ref{item:suites-3}. If Proposition~\ref{Prop2:main_result}\ref{exa-b} is not satisfied, there exist a sequence $(s_n)_{n \in \mathbb{N}} \in \mathbb{C}^{\mathbb{N}}$ and a sequence of vectors $(g_n)_{n \in \mathbb{N}} \in \left(\mathbb{C}^d\right)^{\mathbb{N}}$ such that $\abs{g_n^\transp}=1$  for all $n \in \mathbb{N}$,  $\lim_{n\to\infty} g_n^\transp B = 0$, and $\lim_{n\to\infty} g_n^\transp H(s_n) = 0$. Proposition~\ref{prop_zero_Q0} and the fact that $\rank\begin{pmatrix}A_N & B \end{pmatrix} = d$ imply
that the real part of $(s_n)_{n \in \mathbb{N}}$ is bounded, and thus $(H(s_n))_{n \in \mathbb N}$ is also bounded.
 Without loss of generality, $\lim_{n\to\infty} g_n = g$ for some vector $g \in \mathbb{C}^d$ and we deduce that $\abs{g}=1$, $g^\transp B=0$, and  $\lim_{n\to\infty}g^\transp H(s_n) = 0$,
which
proves one of the two implications. The proof of the converse is obvious.
\end{proof}

Thanks to Lemma~\ref{lem:sous_suite1}, we are able to prove the following proposition establishing a link between the rank properties of $H$ and $\widehat{Q}$.

\begin{proposition}
\label{prop:lienQH}
Let $H$ be the function defined in \eqref{eq:defH} and $\widehat Q$ be the Laplace transform of the distribution $Q$ defined in \eqref{eq:defQ}.
\begin{enumerate}
\item \label{prop:lienHQ1} $\rank\begin{pmatrix}H(s) & B\end{pmatrix} = d$ for every $s \in \mathbb{C}$ if and only if $\rank\begin{pmatrix}\widehat{Q}(s) & B\end{pmatrix} = d$ for every $s \in \mathbb{C}$.
\item \label{prop:lienHQ2} Under the assumption that $\rank\begin{pmatrix}A_N & B\end{pmatrix} = d$, we have that $\rank\begin{pmatrix}M & B\end{pmatrix} = d$  for every $M\in \overline{H(\mathbb{C})}$ if and only if $\rank\begin{pmatrix}\widetilde{Q} & B\end{pmatrix} = d$ for every $\widetilde{Q}\in \overline{\widehat{Q}(\mathbb{C})}$.
\end{enumerate}
\end{proposition}

\begin{proof}
Item~\ref{prop:lienHQ1} is a straightforward consequence of \eqref{lien_entre_les_deux_fonctions}. Let us now prove \ref{prop:lienHQ2}. Note that, thanks to Lemma~\ref{lem:sous_suite1}, there exists $M\in \overline{H(\mathbb{C})}$ such that  $\rank\begin{pmatrix}M & B\end{pmatrix} < d$ if and only if there exists $g \in \mathbb{C}^d$ and $(s_n)_{n \in \mathbb{N}} \in \mathbb{C}^{\mathbb{N}}$ with bounded real part such that $\abs{g}=1$, $g^\transp B=0$, and $\lim_{n\to\infty}g^\transp H(s_n) = 0$. By \eqref{lien_entre_les_deux_fonctions}, we can rewrite the last relation as $\lim_{n\to\infty}g^\transp e^{-s_n \Lambda_N}\widehat{Q}(s_n) = 0$, which is equivalent to $\lim_{n\to\infty}g^\transp\widehat{Q}(s_n) = 0$, since the real part of $(s_n)_{n \in \mathbb{N}} \in \mathbb{C}^{\mathbb{N}}$ is bounded.
\end{proof}

A simple adaptation of Lemma~\ref{lem:sous_suite1}, giving an analogue of Lemma~\ref{lem:sous_suite1}\ref{item:suites-3} in terms of $\widehat Q$, allows us to prove the following result.

\begin{proposition}
\label{prop:eq_Q_corona_thm}
Let $H$ be the function defined in \eqref{eq:defH} and $\widehat Q$ be the Laplace transform of the distribution $Q$ defined in \eqref{eq:defQ}. Assume that $\rank\begin{pmatrix}A_N & B\end{pmatrix} = d$. Then the following statements are equivalent.
\begin{enumerate}
\item $\rank\begin{pmatrix}\widetilde{Q} & B\end{pmatrix} = d$ for every $\widetilde{Q}\in \overline{\widehat{Q}(\mathbb{C})}$.

\item There exists $\alpha>0$ such that, for every $s \in \mathbb{C}$,
\[
\inf \left\{ \abs*{g^\transp \widehat{Q}(s)} + \abs*{g^\transp B} \suchthat g\in \mathbb{C}^d,\, \abs{g^\transp} = 1 \right\} \ge \alpha.
\]
\end{enumerate}
\end{proposition}

We are now in position to provide the proofs of Propositions~\ref{Prop2:main_result0} and \ref{Prop2:main_result}.

\begin{proof}[Proof of Propositions~\ref{Prop2:main_result0} and~\ref{Prop2:main_result}]
We only provide an argument for the proof of Proposition~\ref{Prop2:main_result}, since the proof of Proposition~\ref{Prop2:main_result0} can be obtained following similar arguments, based on a simpler version of Lemma~\ref{lem:sous_suite1} adapted to the case of approximate controllability.

The equivalence between \ref{exa-a} and \ref{exa-b} is a consequence of Lemma~\ref{lem:sous_suite1}\ref{item:suites-1} and \ref{item:suites-3}.
Let us now prove that \ref{exa-b'} implies \ref{exa-b}.
If \ref{exa-b} is not satisfied, it follows by Lemma~\ref{lem:sous_suite1} that there exist  $g\in \mathbb{C}^d$ and a sequence $(s_n)_{n \in \mathbb{N}}\in \mathbb{C}^{\mathbb{N}}$ with a bounded real part such that $\abs{g}=1$, $g^\transp B=0$, and $\lim_{n\to\infty} g^\transp H(s_n) = 0$. Since $H(\cdot)$ is uniformly bounded on every bounded vertical strip of $\mathbb{C}$, we deduce that condition \ref{exa-b'} is not satisfied. Thus \ref{exa-b'} implies \ref{exa-b}. The converse implication is obvious.
\end{proof}

To conclude this section, we recall the following result from \cite[Lemma~4.3]{Yamamoto1988Pseudo}.

\begin{lemma}
\label{lem:4.3_Yamamoto}
Let $\Psi$ be an element in $ \mathcal{D}^\prime_{+}(\mathbb{R})$ such that $\pi(Q*\pi \Psi)=0$. Then there exists a sequence $(\Psi_n)_{n \in \mathbb N}$ in $\mathcal C^\infty(\mathbb R_+, \mathbb R^d)$ such that $\pi(Q * \Psi_n) = 0$ for every $n \in \mathbb N$ and $\lim_{n\to\infty} \Psi_n = \pi \Psi$ in $\mathcal D^\prime(\mathbb R_+)$.
\end{lemma}

\subsection{Approximate controllability}
\label{subsec:Approximate_contr_HY}

Our next aim is to provide a proof of Theorem~\ref{main_result1}. We start with the following approximate controllability criterion for \eqref{syst_lin_avec_sortie}, which is obtained by using the results of \cite{Yamamoto1988Pseudo, Yamamoto1989Reachability}.

\begin{theorem}
\label{th_Yamamoto1}
Let $p \in [1, +\infty)$. The following assertions are equivalent.
\begin{enumerate}
\item\label{item:Yamamoto-approx-2} System~\eqref{syst_lin_avec_sortie} is $L^2$ approximate controllable.
\item\label{item:Yamamoto-approx-p} System~\eqref{syst_lin_avec_sortie} is $L^p$ approximate controllable.
\item\label{item:Yamamoto-approx-Bezout} There exist sequences of matrices $(R_n)_{n \in \mathbb N}$ and $(S_n)_{n \in \mathbb N}$ in $\mathcal E^\prime(\mathbb R_-)$ of appropriate sizes such that
\begin{equation}
\label{eq:approx-Bezout}
\lim_{n \to +\infty} Q * R_n + P * S_n = \delta_0 \Id_d \qquad \text{ in } \mathcal D^\prime(\mathbb R).
\end{equation}
\item\label{item:Yamamoto-approx-Hautus} The two following conditions hold true.
\begin{enumerate}
\item $\rank\begin{pmatrix}\widehat{Q}(s) & B\end{pmatrix} = d$ for every $s \in \mathbb{C}$.

\item $\rank\begin{pmatrix}A_N & B\end{pmatrix} = d$.
\end{enumerate}
\end{enumerate}
\end{theorem}

\begin{proof}
Equivalence between \ref{item:Yamamoto-approx-2} and \ref{item:Yamamoto-approx-Hautus} follows from \cite[Corollary~4.10]{Yamamoto1989Reachability} after recalling that $P = B \delta_0$ and noting that $Q$ can be written as $Q = Q_0 + Q_1$, with $Q_0 = \delta_0 A_N$ and $\supp Q_1 \subset [-\Lambda_N, -\Lambda_N + \Lambda_{N-1}]$, while the equivalence between \ref{item:Yamamoto-approx-2} and \ref{item:Yamamoto-approx-Bezout} was shown in \cite[Theorem~4.4]{Yamamoto1988Pseudo}. The proof of the equivalence between \ref{item:Yamamoto-approx-2} and \ref{item:Yamamoto-approx-Bezout} in \cite[Theorem~4.4]{Yamamoto1988Pseudo} can be immediately generalized in order to also show the equivalence between \ref{item:Yamamoto-approx-p} and \ref{item:Yamamoto-approx-Bezout}. For sake of completeness, we provide this generalization here.

Assume that \ref{item:Yamamoto-approx-p} holds true. Since $Q \in \mathcal E^\prime(\mathbb R_-)$ and $Q^{-1} \in \mathcal D^\prime_+(\mathbb R)$, Lemma~\ref{lem:4.3_Yamamoto} applied to each column of $Q^{-1}$ ensures that there exists a sequence $(T_n)_{n \in \mathbb N}$ in $\mathcal C^\infty([0, +\infty), \mathcal M_d(\mathbb R))$ such that $\pi(Q * T_n) = 0$ for all $n \in \mathbb N$ and $T_n \to \pi Q^{-1}$ as $n \to +\infty$ in $\mathcal D^\prime(\mathbb R_+)$. Since $\pi(Q * T_n) = 0$, we have that each column of $T_n$ belongs to $X^{Q, p}$ and thus, by \ref{item:Yamamoto-approx-p}, we deduce that, for every $n \in \mathbb N$, there exists a sequence $(S_{n, k})_{k \in \mathbb N}$ of $m \times d$ matrices whose columns are all in $\Omega_p$ and such that $\pi(Q^{-1} * P * S_{n, k}) \to T_n$ as $k \to +\infty$, in the topology of $L^p_{\loc}\left(\mathbb{R}_+, \mathcal M_{d}(\mathbb{R})\right)$, and hence also in that of $\mathcal D^\prime(\mathbb R_+)$. Thus, there exists a sequence $(k_n)_{n \in \mathbb N}$ in $\mathbb N$ such that $k_n \to +\infty$ as $n \to +\infty$ and $\pi(Q^{-1} * P * S_{n, k_n}) \to \pi Q^{-1}$ as $n \to +\infty$ in $\mathcal D^\prime(\mathbb R_+)$. For $n \in \mathbb N$, we set $S_n = S_{n, k_n}$, so that
\[
\lim_{n \to +\infty} \pi(Q^{-1} * P * S_n) = \pi Q^{-1} \qquad \text{ in } \mathcal D^\prime(\mathbb R_+).
\]

For $n \in \mathbb N$, let $\Psi_n$ be the extension of $\pi(Q^{-1} * P * S_n) \in L^p_{\loc}\left(\mathbb{R}_+, \mathcal M_{d}(\mathbb{R})\right)$ to $L^p_{\loc}\left(\mathbb{R}, \mathcal M_{d}(\mathbb{R})\right)$ obtained by setting $\Psi_n(x) = 0$ for $x < 0$. By \eqref{eq:explicit-Q-1}, we have that $\inf \supp Q^{-1} = \Lambda_N > 0$, and we thus see that $\lim_{n \to +\infty} \Psi_n = Q^{-1}$ in $\mathcal D^\prime(\mathbb R)$. Let $R_n = \Psi_n - Q^{-1} * P * S_n$. Since $\Psi_n = Q^{-1} * P * S_n$ on $\mathbb R_+$, we deduce that $\supp R_n \subset \mathbb R_-$. In addition, $\supp(Q^{-1} * P * S_n)$ is bounded on the left, and thus $\supp R_n$ is compact. This shows that $R_n \in \mathcal E^\prime(\mathbb R_-)$. We have, in $\mathcal D^\prime(\mathbb R)$,
\[
Q^{-1} = \lim_{n \to +\infty} \Psi_n = \lim_{n \to +\infty} \left[Q^{-1} * P * S_n + R_n\right],
\]
and \ref{item:Yamamoto-approx-Bezout} follows by taking the convolution with $Q$ on the left in the above identity.

Assume now that \ref{item:Yamamoto-approx-Bezout} holds true and let $(R_n)_{n \in \mathbb N}$ and $(S_n)_{n \in \mathbb N}$ be as in \ref{item:Yamamoto-approx-Bezout}. Using the identification between $X^{Q, p}$ and $L^p([0, \Lambda_N], \mathbb R^d)$ described in Section~\ref{subsec:real_theory}, one deduces that $X^{Q, p} \cap \mathcal C^\infty(\mathbb R_+, \mathbb R^d)$ is dense in $X^{Q, p}$, since this subset of $X^{Q, p}$ can be identified with a subset of $L^p([0, \Lambda_N], \mathbb R^d)$ that contains $\mathcal C^\infty_{\mathrm{c}}((0, \Lambda_N), \mathbb R^d)$, which is dense in $L^p([0, \Lambda_N], \mathbb R^d)$. Hence, to prove \ref{item:Yamamoto-approx-p}, it suffices to show that, for every $\phi \in X^{Q, p} \cap \mathcal C^\infty(\mathbb R_+, \mathbb R^d)$, there exists a sequence of inputs $(u_n)_{n \in \mathbb N}$ in $\Omega_p$ such that $\pi(Q^{-1} * P * u_n)$ converges in $L^p_{\loc}(\mathbb R_+, \mathbb R^d)$ to $\phi$ as $n \to +\infty$.

Take $\phi \in X^{Q, p} \cap \mathcal C^\infty(\mathbb R_+, \mathbb R^d)$ and extend $\phi$ to a function $\widetilde\phi \in \mathcal C^\infty(\mathbb R, \mathbb R^d)$ such that $\supp{\widetilde\phi}$ is bounded on the left. Note that $\pi\widetilde\phi = \phi$. Let $u_n = S_n * Q * \widetilde\phi$ and note that, by definition, $u_n \in \mathcal C^\infty(\mathbb R, \mathbb R^m)$ and $\supp u_n$ is bounded on the left. Using Lemma~\ref{lem:sumlemmasYamamoto}\ref{lem:A_3_Yamamoto} and the fact that $\phi \in X^{Q, p}$, we have that
\[
\pi u_n = \pi(S_n * Q * \widetilde\phi) = \pi(S_n * \pi(Q * \widetilde\phi)) = \pi(S_n * \pi(Q * \pi\widetilde\phi)) = 0,
\]
showing that $\supp u_n \subset \mathbb R_-$, and thus $u_n \in \Omega_p$. By a similar reasoning, we also deduce that
\begin{equation}
\label{eq:pi-Rn-zero}
\pi(R_n * Q * \widetilde\phi) = \pi(R_n * \pi(Q * \pi\widetilde\phi)) = 0.
\end{equation}
Notice now that
\begin{align*}
Q^{-1} * P * u_n + R_n * Q * \widetilde\phi & = Q^{-1} * P * S_n * Q * \widetilde\phi + R_n * Q * \widetilde\phi \\
& = Q^{-1} * (P * S_n + Q * R_n) * Q * \widetilde\phi \to \widetilde\phi
\end{align*}
as $n \to +\infty$ and, since $\widetilde\phi \in \mathcal C^{\infty}(\mathbb R, \mathbb R^d)$ and its support is bounded on the left, one deduces, by standard properties of the convolution product between a distribution and a test function, that the above convergence is in particular uniform on all compact subsets of $\mathbb R$. Combining \eqref{eq:pi-Rn-zero} with the linearity of the truncation operator $\pi$ and its continuity with respect to uniform convergence on compact subsets, we deduce that
\[
\phi = \pi\widetilde\phi = \lim_{n \to +\infty} \pi(Q^{-1} * P * u_n + R_n * Q * \widetilde\phi) = \lim_{n \to +\infty} \pi(Q^{-1} * P * u_n),
\]
where the limit is uniform in compact subsets of $\mathbb R_+$. In particular, the above limit takes place also in the topology of $L^p_{\loc}(\mathbb R_+, \mathbb R^d)$, showing that \ref{item:Yamamoto-approx-p} holds true, as required.
\end{proof}

We are finally in position to provide a proof of Theorem~\ref{main_result1}.

\begin{proof}[Proof of Theorem~\ref{main_result1}]
Theorem~\ref{lem:3} proves that $L^p$ approximate controllability is equivalent to the $L^p$ approximate controllability in time $T=d \Lambda_N$. The conclusion follows by combining Theorem~\ref{th_Yamamoto1} and Proposition~\ref{prop:lienQH}\ref{prop:lienHQ1}.
\end{proof}

\begin{remark}
Definition~\ref{def:Lq-cont} introduces the notion of $L^p$ approximate controllability in time $T$ by using the strong topology of $L^p$. One can also introduce a similar notion using the weak topology of $L^p$: For $p \in [1, +\infty)$, we say that \eqref{eq:intro-linear-difference-equation} is \emph{$L^p$ weakly approximately controllable in time $T>0$} if, for every $x_0$ and $\phi$ in $L^p([-\Lambda_N,0], \allowbreak\mathbb{R}^d)$, there exists a sequence $(u_n)_{n \in \mathbb N}$ in $L^p([0,T],\mathbb{R}^m)$ such that the solution $x_n$ of \eqref{eq:intro-linear-difference-equation} with initial condition $x_0$ and control $u_n$ satisfies
\[
x_{n, T} \xrightharpoonup[n \to +\infty]{} \phi,
\]
where the arrow $\xrightharpoonup{}$ denotes weak convergence in $L^p([-\Lambda_N,0], \allowbreak\mathbb{R}^d)$. A similar definition can also be provided for $L^\infty$ weak approximate controllability in time $T$, replacing the above weak convergence by weak-$\ast$ convergence.

It turns out that Theorem~\ref{main_result1} remains true when approximate controllability is replaced by weak approximate controllability, for any $p \in [1, +\infty]$. Indeed, Definition~\ref{def:contr-from-origin} can be easily adapted to weak approximate controllability from the origin, and Theorem~\ref{lem:3} remains true for weak approximate controllability from the origin and weak approximate controllability in time $d \Lambda_N$, with the same proof, since it relies on Theorem~\ref{lem:RanE_con}, which is independent of the choice of topology in $L^p$. As for Theorem~\ref{th_Yamamoto1}, the proof of the equivalence between \ref{item:Yamamoto-approx-p} and \ref{item:Yamamoto-approx-Bezout} provided above can be straightforwardly adapted to the case where \ref{item:Yamamoto-approx-p} is replaced by $L^p$ weak approximate controllability, after noticing that weak convergence in $L^p$ for $p \in [1, +\infty)$ and weak-$\ast$ convergence in $L^\infty$ both imply convergence in the sense of distributions.
\end{remark}

\subsection{Exact controllability}
\label{subsec:Exact_contr_HY}

Let us now turn to the proof of the necessary condition for exact controllability from Theorem~\ref{main_result2}, which turns out to be much more subtle than that of Theorem~\ref{main_result1}. In the sequel, we start by expounding and improving some of the exact controllability results from \cite{Yamamoto1988Pseudo, Yamamoto2008Behavioral, Yamamoto2011Bezout}.

\subsubsection{Radon exact controllability and Bézout's identity}

The main idea in \cite{Yamamoto1988Pseudo, Yamamoto2008Behavioral, Yamamoto2011Bezout} to address exact controllability is to address this question in a distributional sense. Thanks to the specific features of \eqref{eq:intro-linear-difference-equation}, we are able to replace the general distributional framework of those works by the more structured setting of distributions of order zero (or, equivalently, Radon measures). More precisely, we first define the Radon measure space $\overline{X}^Q$ by
\begin{equation*}
\overline{X}^Q=\left\{ \pi \Psi\suchthat \Psi \in  \left(M(\mathbb{R}_+)\right)^d \text{ and } \pi(Q*\pi \Psi)=0  \right\}.
\end{equation*}
It is then straightforward to see that outputs of the input-output system defined in \eqref{representation_radon_measure3} corresponding to inputs in
$\left(M(\mathbb{R}_-)\right)^{m}$ belong to $\overline{X}^Q$.
We next extend the definition of exact controllability to Radon measures as follows.

\begin{definition}
\label{def_reachability_distribution}
System~\eqref{syst_lin_avec_sortie} is \emph{Radon exactly controllable} if, for every $\Psi \in (M(\mathbb{R}_+))^d$ such that $\pi \Psi \in \overline{X}^Q$, there exists $u \in \left(M(\mathbb{R}_-)\right)^m$ such that $\pi(Q^{-1}*P*u) = \pi \Psi$.
\end{definition}

We have the following characterization of Radon exact controllability in terms of a Bézout's identity over the ring of Radon measures.

\begin{proposition}
\label{prop_fond_exact_controllability}
System~\eqref{syst_lin_avec_sortie} is Radon exactly controllable if and only if there exist  two matrices $R$ and $S$ with entries in $M(\mathbb{R}_{-})$ such that
\begin{equation}
\label{eq:exact_cont0}
Q*R+P*S=\delta_0 \Id_d.
\end{equation}
Moreover, if \eqref{eq:exact_cont0} holds true, then,
for every target output $\pi \Psi \in \overline{X}^Q$ with $\Psi \in (M(\mathbb{R}_+))^d$, the input $\omega \in M_+(\mathbb R)$ given by
\begin{equation}\label{eq:motion-planning}
\omega = S * Q * \Psi
\end{equation}
steers the origin to the state $\pi \Psi$ along the system \eqref{syst_lin_avec_sortie}, i.e.,
\begin{equation}\label{eq:sol-con-B}
\pi(Q^{-1} * P * \omega) = \pi \Psi.
\end{equation}
\end{proposition}

\begin{proof}
We first prove that \eqref{eq:exact_cont0} is a necessary condition for Radon exact controllability of \eqref{syst_lin_avec_sortie}. The entries of $Q^{-1}$ are in the space $M_{+}(\mathbb{R})$ and we have $\pi (Q^{-1})=Q^{-1}$ so that the columns of $Q^{-1}$ are in $\overline{X}^Q$. Since \eqref{syst_lin_avec_sortie} is Radon exactly controllable, we can take $\Psi$ equal to each column of $Q^{-1}$ in Definition~\ref{def_reachability_distribution}
and thus deduce the existence of a matrix $S$ of size $m \times d$ and with entries in $M(\mathbb{R}_-)$ such that
\begin{equation}
\label{eq:bezout20}
\pi(Q^{-1}*P*S)
=\pi \left(Q^{-1}\right).
\end{equation}
We define $R$ by
\begin{equation}
\label{eq:bezout30}
R = Q^{-1}-Q^{-1}*P*S.
\end{equation}
From \eqref{eq:bezout20}, we deduce that $\pi(R)=0$ and Lemma~\ref{lem:sumlemmasYamamoto}\ref{lem:A_2_Yamamoto} proves that the support of $R$ is included in $\mathbb{R}_-$. Since the supports of $Q^{-1}$, $P$, and $S$ are bounded on the left, the same is true for the support of $R$, which is, therefore, compact. Moreover the entries of $R$ are distributions of order zero because the same  is true for $Q^{-1}$, $P$, and $S$. Hence the entries of $R$ belong to $M(\mathbb R_-)$, and we finally obtain \eqref{eq:exact_cont0} by taking the convolution product of \eqref{eq:bezout30} on the left by $Q$.

We now prove that \eqref{eq:exact_cont0} is sufficient for Radon exact controllability of \eqref{syst_lin_avec_sortie}. Let $R$ and $S$ be two  matrices with entries in the space $M(\mathbb{R}_{-})$ satisfying \eqref{eq:exact_cont0}.
Let $\Psi \in \left(M(\mathbb{R}_+)\right)^d$ be such that
\begin{equation}\label{eq:psiXQ}
\pi(Q*\pi\Psi)=0,
\end{equation}
so that $\pi\Psi \in \overline{X}^Q$. Set $\omega = S*Q*\Psi$. Since $Q$ and $S$ have entries in $M(\mathbb{R_-})$, we have that $\omega$ belongs to the space $M_+(\mathbb{R})$. Lemma~\ref{lem:sumlemmasYamamoto}\ref{lem:A_3_Yamamoto} implies that
\begin{equation*}
\pi\left( \omega\right)= \pi\left(S*\pi\left(Q*  \pi\Psi\right)\right)=0,
\end{equation*}
and we deduce that $\omega$ is in $\left(M(\mathbb{R}_-)\right)^{m}$ by Lemma~\ref{lem:sumlemmasYamamoto}\ref{lem:A_2_Yamamoto}. By definition of $\omega$ and \eqref{eq:exact_cont0}, we also have
\begin{equation}
\label{eq:bezout5}
Q^{-1} * P * \omega + R * Q * \Psi = Q^{-1} * (P * S + Q * R) * Q * \Psi = \Psi.
\end{equation}
We deduce from \eqref{eq:bezout5}, Lemma~\ref{lem:sumlemmasYamamoto}\ref{lem:A_3_Yamamoto}, and \eqref{representation_radon_measure3} that
\begin{align*}
\pi \Psi & = \pi\left(Q^{-1}*P*\omega\right)+\pi\left(R*Q*\Psi\right)\\
& =\pi\left(Q^{-1}*P*\omega\right)+\pi\left(R*\pi\left(Q*\pi \Psi \right)\right)\\
& = \pi\left(Q^{-1}*P*\omega\right),
\end{align*}
where the last equality follows from \eqref{eq:psiXQ}.
We have shown that $\pi \Psi$ is the output corresponding to $\omega$ for the input-output map~\eqref{representation_radon_measure3}, proving the Radon exact controllability of \eqref{syst_lin_avec_sortie}.
\end{proof}

\begin{remark}
Proposition~\ref{prop_fond_exact_controllability} can be seen as a Radon counterpart of the distributional Bézout's identity characterization of exact distributional controllability provided in \cite{Yamamoto2008Behavioral}.
\end{remark}

As a consequence of Proposition~\ref{prop_fond_exact_controllability}, we can now deduce that Radon exact controllability implies $L^p$ exact controllability.

\begin{corollary}
\label{cor:distru_exac_equiv}
If \eqref{syst_lin_avec_sortie} is Radon exactly controllable then it is $L^p$ exactly controllable for every $p \in [1,+\infty)$.
\end{corollary}

\begin{proof}
Let $R$ and $S$ be two  matrices with entries belonging to the space $M(\mathbb{R}_{-})$ satisfying \eqref{eq:exact_cont0}. Let $p \in [1,+\infty)$. Let $y \in  X^{Q,\,p}$ and denote by $\widetilde{y}$ the extension of $y$ on $(-\infty,+\infty)$ by setting $\widetilde{y}$ equal to zero on $(-\infty,0)$. Set $\omega=S*Q*\widetilde{y}$ as in \eqref{eq:motion-planning}. Then $\omega \in \Omega_p$ and it follows from \eqref{eq:sol-con-B} that $y = \pi\left(Q^{-1}*P*\omega\right)$, which proves $L^p$ exact controllability thanks to Proposition~\ref{def_reachability}.
\end{proof}

A challenging question is to investigate a converse to the previous corollary, i.e., whether $L^p$ exact controllability for some $p \in [1,+\infty)$ implies Radon exact controllability or not. We bring a positive answer to that question for $p=1$ in the following section.

\subsubsection{\texorpdfstring{$L^1$}{L1} exact controllability and Bézout's identity}

We next provide a necessary and sufficient condition for \eqref{syst_lin_avec_sortie} to be $L^1$ exactly controllable.

\begin{theorem}
\label{prop_fond_exact_controllability_L1}
System~\eqref{syst_lin_avec_sortie} is $L^1$ exactly controllable if and only if there exist  two matrices $R$ and $S$ with entries in $M(\mathbb{R}_{-})$ such that
\eqref{eq:exact_cont0} holds true.
\end{theorem}

\begin{proof}
We first notice that, by Proposition~\ref{prop_fond_exact_controllability} and Corollary~\ref{cor:distru_exac_equiv}, \eqref{eq:exact_cont0} implies that \eqref{syst_lin_avec_sortie} is $L^1$ exactly controllable.

Let us now prove the converse implication. We proceed in four steps.

\step{1} Let us define the map $\widetilde{G}\colon \widetilde{\Omega}_1 \longrightarrow L^1([0,\Lambda_N],\mathbb{R}^d)$ by
\[\widetilde{G}(u)(t)=\pi(Q^{-1}*P*u)(t),\qquad t\in [0,\Lambda_N],\, u \in \widetilde{\Omega}_1,\]
where $\widetilde{\Omega}_1$ denotes the subspace of $\Omega_1$ made of inputs with compact support in $[-d \Lambda_N,0]$, endowed with the norm $\norm{\cdot}_{L^1([-d\Lambda_N, 0],\mathbb R^m)}$. Firstly, we can see that the map $\widetilde{G}$ is a bounded operator because $Q^{-1}*P$ is a distribution with a finite number of Dirac distributions on each compact interval of $\mathbb{R}$. We deduce that $\widetilde{G}$ is a continuous linear map. Secondly, Theorem~\ref{lem:RanE_con} implies that \eqref{syst_lin_avec_sortie} is $L^1$ exactly controllable if and only if the map $\widetilde{G}$ is surjective. We can now apply the open mapping theorem (see, e.g, \cite[Theorem~4.13]{Rudin1991Functional}) and deduce that there exists $\delta>0$ such that
\begin{equation}
\label{eq:open_map_thereom}
\widetilde{G}(U) \supset \delta V,
\end{equation}
where $U$ and $V$ are the open unit balls of $ \widetilde{\Omega}_1$ and $L^1([0,\Lambda_N],\mathbb{R}^d)$ respectively.

\step{2} Since $\pi(Q*\pi (Q^{-1}))=\pi(\delta_0)=0$, Lemma~\ref{lem:4.3_Yamamoto} implies that there exists a sequence $\psi_n=(\psi_{n,1},\dotsc,\psi_{n,d}) \in (X^{Q,1})^d$,  $n \in \mathbb{N}$, such that $ \psi_n \rightarrow \pi (Q^{-1})$ in the distributional sense as $n\to\infty$. In other words, for $i,j \in \llbracket 1, d\rrbracket$, if we define the Radon measures $\left(Q^{-1}_{n}\right)_{i,j}(t)=\int_0^t \left(\psi_{n}\right)_{i,j}(x) \diff x$ for $t \in \mathbb{R}_+$ and $n \in \mathbb{N}$, we get that $\left(Q^{-1}_n\right)_{i,j}$ weak-* converges to $\left(\pi(Q^{-1})\right)_{i,j}$ in the sense of \cite[Section~4.3]{Maggi2012Sets}. By \cite[Remark~4.35]{Maggi2012Sets} and the Banach--Steinhaus theorem, we obtain that
\begin{equation*}
\sup_{n \in \mathbb{N}}\abs*{\left(Q^{-1}_n\right)_{i,j}}([0,\Lambda_N])<+\infty,
\end{equation*}
where
\begin{equation*}
\abs*{\left(Q^{-1}_n\right)_{i,j}}([0,\Lambda_N]) = \sup_{\substack{\varphi \in \mathcal C([0, \Lambda_N]) \\ \norm{\varphi}_{\mathcal C([0, \Lambda_N], \mathbb R)} \leq 1}} \abs*{\int_0^{\Lambda_N}  \varphi(t) \diff \left(Q^{-1}_n\right)_{i,j}(t)}.
\end{equation*}
By the Riesz representation theorem (see, e.g., \cite[Theorem~6.19]{Rudin1987Real}), $\abs*{\left(Q^{-1}_n\right)_{i,j}}([0, \Lambda_N])$ is the total variation of the measure on the interval $[0, \Lambda_N]$ of $\left(Q^{-1}_n\right)_{i,j}$ in the sense of \cite[Chapter~6, (3)]{Rudin1987Real}. In particular, we have that
\begin{equation*}
\abs*{\left(Q^{-1}_n\right)_{i,j}}([0, \Lambda_N])= \int_0^{\Lambda_N} \abs*{\left(\psi_n\right)_{i,j}(t)} \diff t,
\end{equation*}
so that each column of $\left(\psi_{n}\right)_{n \in \mathbb{N}} $  is uniformly bounded in $L^1([0, \Lambda_N],\mathbb{R}^d)$, that is, there exists $C>0$ such that, for all $j \in \llbracket 1, d\rrbracket$ and $n \in \mathbb{N}$, we have
\begin{equation*}
\norm{\psi_{n,j}}_{L^1([0, \Lambda_N],\mathbb{R}^d)} \le C.
\end{equation*}
Let $M^\prime>0$ be such that $\delta M^\prime >C$. We get from \eqref{eq:open_map_thereom} that
\begin{equation*}
\widetilde{G}(M^\prime U) \supset \delta M^\prime V
\end{equation*}
so that, for all $\psi_{n,j}$ with  $j \in \llbracket 1,d\rrbracket$ and $n \in \mathbb{N}$, there exists $S_{n,j} \in \widetilde{\Omega}_1$ such that
\begin{equation*}
\widetilde{G}(S_{n,j})= \psi_{n,j} \quad \text{ and } \quad  \norm{S_{n,j}}_{L^1([-d\Lambda_N,0], \mathbb R^m)} \le M^\prime.
\end{equation*}

\step{3} We define $S_n=(S_{n,1},\dotsc,S_{n,d})$.
 By construction,
 $S_n \in {\widetilde{\Omega}_1}^d$ and
\begin{equation}
\label{eq:bezout1}
\pi(Q^{-1}*P*S_n) \rightarrow \pi(Q^{-1}), \qquad \text{ as }n\to\infty,
\end{equation}
in a distributional sense. Since the columns of $S_n$, for $n \in\mathbb{N}$, are uniformly bounded for the norm in $\widetilde{\Omega}_1$, by the weak compactness of Radon measures (see for instance \cite[Theorem~4.33]{Maggi2012Sets}), there exist a matrix $S$ with entries in $M(\mathbb{R}_-)$ and a subsequence of $\left( S_n \right)_{n \in \mathbb{N}}$ (still denoted $\left(S_n\right)_{n \in \mathbb{N}}$ by abuse of notation) such that $\lim_{n \rightarrow + \infty}S_n = S$ in distributional sense. Since the convolution is continuous in distributional sense (see, for instance, \cite[\foreignlanguage{french}{Théorème}~7.4.9]{Bony2001Cours}), and $\pi$ is continuous with respect to the strong dual topology of $\mathcal D^\prime(\mathbb R)$, we deduce from \eqref{eq:bezout1} that
\begin{equation}
\label{eq:bezout2}
\pi(Q^{-1}*P*S)
=\pi \left(Q^{-1}\right).
\end{equation}

\step{4} We define $R$ by
\begin{equation}
\label{eq:bezout3}
R = Q^{-1}-Q^{-1}*P*S
\end{equation}
and we conclude the proof as in Proposition~\ref{prop_fond_exact_controllability}: by \eqref{eq:bezout2} we have $\pi(R)=0$, which, together with Lemma~\ref{lem:sumlemmasYamamoto}\ref{lem:A_2_Yamamoto}, implies that the support of $R$ is compact and contained in $(-\infty,0]$. Moreover, the entries of $R$ are distributions of order zero because the same  is true for $Q^{-1}$, $P$, and $S$. Bézout's identity is then obtained by taking the convolution product of \eqref{eq:bezout3} on the left by $Q$.
\end{proof}

\begin{remark}
\label{remark1}
It is an open question whether, for each $p \in [1,+\infty)$, Bézout's identity is equivalent to the $L^p$ controllability. We can stress that the proof given for the case $p=1$ does not straightforwardly extend to $p > 1$. Indeed, in Step~2, the convergence $\psi_n \to \pi(Q^{-1})$ in the distributional sense implies the boundedness of $(\psi_n)_{n \in \mathbb N}$ in $L^1([0, \Lambda_N], \mathbb R^d)$, but such a sequence may fail to be bounded in $L^p([0, \Lambda_N], \mathbb R^d)$ for $p > 1$.
\end{remark}

As an immediate consequence of Proposition~\ref{prop_fond_exact_controllability} and Theorem~\ref{prop_fond_exact_controllability_L1}, we obtain the following equivalence between Radon and $L^1$ exact controllability.

\begin{corollary}
\label{corollaire_L^1_equiv_radon}
System~\eqref{syst_lin_avec_sortie} is Radon exactly controllable if and only if it is $L^1$ exactly controllable.
\end{corollary}

The open question raised in Remark~\ref{remark1} is tantamount to know if, for difference delay systems of the type~\eqref{eq:intro-linear-difference-equation}, the $L^p$ exact controllability for some $p \in [1,+\infty)$ is equivalent to the $L^p$ exact controllability  for every $p \in [1,+\infty)$. It has to be pointed out that such a property holds when dealing with the exponential stability of such systems, as noticed in  \cite{Baratchart2021Sufficient, Chitour2016Stability}.

\subsubsection{Solvability of Bézout's identity and proof of Theorem~\ref{main_result2}}

Proposition~\ref{prop_fond_exact_controllability} and Theorem~\ref{prop_fond_exact_controllability_L1} reduced the exact controllability problem to the
problem of existence of solutions of a suitable Bézout's identity.

\begin{definition}
\label{def:solvable-Bezout}
We say that Bézout's identity \eqref{eq:exact_cont0} is \emph{solvable} if there exist two matrices $R$ and $S$ with entries in $M(\mathbb{R}_{-})$ that satisfy it.
\end{definition}

We next show that the solvability of Bézout's identity is equivalent to a corona problem, the latter having led a tremendous literature on the subject for some alike problems. We can cite for example the celebrated paper \cite{Carleson1962Interpolations} resolving the corona problem in one dimension for holomorphic bounded functions in the unit disk. However the present corona problem arising from Bézout's identity over a Radon measure algebra has not received a great attention and it is still an open question.
To the best of our knowledge, the closest result dealing with the solvability of a Bézout's identity is the one provided in \cite[Theorem~5.1]{Yamamoto2011Bezout} over the algebra of distributions, which gives a sufficient condition on $(Q,P)$ for such a solvability. In turn, solving a Bézout's identity over Radon measures is fundamentally different than over distributions because of the different topologies endowing these two spaces.

Our first result for Radon measures is the following necessary condition for the solvability of Bézout's identity.

\begin{proposition}
\label{prop_nece_bezout}
A necessary condition for the Bézout's identity \eqref{eq:exact_cont0} to be solvable is
the following.
\begin{enumerate}
\item\label{cdt_corona} There exists $\alpha>0$ such that, for every $s \in \mathbb{C}$,
\begin{equation}
    \label{eq_corona}
\inf \left\{ \abs*{g^\transp \widehat{Q}(s)} + \abs*{g^\transp B} \suchthat g\in \mathbb{C}^d,\, \abs{g^\transp} = 1 \right\} \ge \alpha.
\end{equation}
\end{enumerate}
\end{proposition}

\begin{proof}
Let $R$ and $S$ be the two matrices with entries in $M(\mathbb R_-)$ satisfying Bezout's identity from \eqref{eq:exact_cont0}. By contradiction, we assume that \ref{cdt_corona} is not satisfied, so that there exist a sequence $(s_n)_{n \in \mathbb{N}} \in \mathbb{C}^{\mathbb{N}}$ and a sequence
 $(g_n)_{n \in \mathbb{N}} \in \left(\mathbb{C}^d\right)^{\mathbb{N}}$ such that $\abs{g_n^\transp}=1$  for all $n \in \mathbb{N}$,  $\lim_{n\to\infty} g_n^\transp B = 0$, and $\lim_{n\to\infty} g_n^\transp \widehat{Q}(s_n) = 0$.

We claim that there exists $\widetilde{\alpha}>0$ such that $\Real s_n \le \widetilde{\alpha}$ for every $n\in \mathbb{N}$. Indeed, $\widehat{Q}(s)$ is equivalent to $e^{s \Lambda_N} \Id_d$ as $\Real s \to +\infty$, and unboundedness from above of $\Real s_n$ would contradict the relation $\lim_{n\to\infty} g_n^\transp \widehat{Q}(s_n) = 0$.
We next get, by a classical estimate of the Laplace transform of an element of $M(\mathbb{R}_-)$, that
there exists $C>0$ such that, for all $s \in \mathbb{C}$ with $\Real s \le \widetilde{\alpha}$,
\begin{equation}
\label{eq:PW1}
\max \left\{\abs*{\widehat{R}(s)}, \abs*{\widehat{S}(s)} \right\} \le C.
\end{equation}
We then deduce, using \eqref{eq:exact_cont0} and applying the Laplace transform, that
\begin{equation}
\label{eq:PW2}
g_n^\transp \widehat{Q}(s) \widehat{R}(s)+g_n^\transp B \widehat S(s) = g_n^\transp,\qquad n \in \mathbb{N},\, s\in\mathbb{C}.
\end{equation}

Equations~\eqref{eq:PW1} and \eqref{eq:PW2} imply that it is impossible to have both  $\lim_{n\to\infty} g_n^\transp B = 0$ and $\lim_{n\to\infty} g_n^\transp \widehat{Q}(s_n) = 0$.
We reached a contradiction, so that \ref{cdt_corona} is a necessary condition for the solvability of Bézout's identity.
\end{proof}

We can now undertake the proof of Theorem~\ref{main_result2}.

\begin{proof}[Proof of Theorem~\ref{main_result2}]
Theorem~\ref{lem:3} states that $L^1$ exact controllability is equivalent to the $L^1$ exact controllability in time $T=d \Lambda_N$. By Theorem~\ref{main_result1}, Theorem~\ref{main_result2}\ref{assumption2-L1} is necessary for $L^1$ approximate controllability, so that it is also necessary for $L^1$ exact controllability. Combining Theorem~\ref{prop_fond_exact_controllability_L1} and Proposition~\ref{prop_nece_bezout}, we deduce that \eqref{eq_corona} is a necessary condition for $L^1$ exact controllability, and the equivalence between \eqref{eq_corona} and Theorem~\ref{main_result2}\ref{assumption1-L1} is a consequence of Propositions~\ref{prop:lienQH}\ref{prop:lienHQ1} and \ref{prop:eq_Q_corona_thm}.
\end{proof}

The necessity of Proposition~\ref{prop_nece_bezout}\ref{cdt_corona} is the easy part for the solvability of Bézout's identity over the Radon measure algebra. The harder part would be to prove that it is also a sufficient condition to have the solvability of Bézout's identity. A first step of such a proof would consist in using exactly the trick from \cite{Fuhrmann1968Corona}, which reduces a corona matrix problem to the corresponding 1D version, which can be written as follows.

\begin{conjecture}
\label{conjecture1}
Let $k \in \mathbb{N}$ and $q_j \in  M(\mathbb{R}_-)$ for $j \in \llbracket 1, k\rrbracket$. Assume the existence of a $c>0$ such that
\begin{equation*}
\sum_{j=1}^k \abs{\widehat{q}_j(s)} \ge c>0, \qquad s \in \mathbb{C}.
\end{equation*}
Then there exists $p_j \in  M(\mathbb{R}_-)$ for $j \in \llbracket 1, k\rrbracket$ satisfying the equation
\begin{equation*}
\sum_{j=1}^k \widehat{q}_j(s)\widehat{p}_j(s) =1, \qquad s \in \mathbb{C}.
\end{equation*}
\end{conjecture}

\begin{remark}
In fact, due to the structure of the difference delay systems that we consider, i.e., with a finite number of delays, it would be sufficient to prove Conjecture~\ref{conjecture1} for the measures $q_j$, for all $j \in \llbracket 1, k\rrbracket$, belonging to the sub-algebra of Radon measures $M(\mathbb{R}_-)$ finitely generated by the Dirac measures $(\delta_{\Lambda_i-\Lambda_N})_{i \in \llbracket 1, N\rrbracket}$.
\end{remark}

Notice that a proof of Conjecture~\ref{conjecture1} would imply a positive answer to Conjecture~\ref{conj:HY-exact} in the case $p = 1$, yielding a sufficient and necessary criterion for the $L^1$ exact controllability.

\subsection{Sufficient exact controllability criterion for \texorpdfstring{$\mathcal C^k$}{Ck} functions}
\label{subsec:Ck}

We now turn to our last exact controllability result from Theorem~\ref{thm:Ck}, which relies on the sufficient condition for the resolution of Bézout's identity over the distributional algebra in dimension $d$ given in \cite{Yamamoto2011Bezout}. In order to translate the controllability notion of Theorem~\ref{thm:Ck} in terms of \eqref{syst_lin_avec_sortie}, we provide the following definition.

\begin{definition}
Let $k$ be an integer and set
\begin{equation*}
X_{k}^Q = \left\{y \in  \mathcal C^k(\mathbb{R}_+,\mathbb{R}^d) \suchthat \pi(Q*y)=0  \right\}.
\end{equation*}
System~\eqref{syst_lin_avec_sortie} is said to be \emph{$\mathcal C^k$ exactly controllable} if, for every $y \in X_{k}^Q$, there exists $u \in \mathcal C^0(\mathbb{R},\mathbb{R}^m)$ compactly supported in $\mathbb{R}_-$ such that $\pi(Q^{-1}*P*u)=y$.
\end{definition}

\begin{remark}
Note that, contrarily to the notions of controllability that we introduced previously, $\mathcal C^k$ exact controllability requires the control $u$ to belong only to the less regular space of continuous functions, instead of requiring it to belong to a space similar to that of the state $y$.
\end{remark}

The proposition below states a sufficient controllability criteria for the $\mathcal C^k$ functions.

\begin{proposition}\label{prop-Ckexact}
If the condition in Proposition~\ref{prop_nece_bezout}\ref{cdt_corona} holds true, then there exists $k \in \mathbb N$ such that \eqref{syst_lin_avec_sortie} is $\mathcal C^k$ exactly controllable.
\end{proposition}

\begin{proof}
In this proof, we fix $\abs{\cdot}$ as the Euclidean norm. We proceed in two steps.

\step{1} We aim to apply \cite[Theorem~5.1]{Yamamoto2011Bezout}. Note first that $\det \left(\widehat{Q}(\cdot)\right)$ is a quasipolynomial, as defined in Section~\ref{sec:quasipolynomials}, and thus, by Corollary~\ref{coro:Polya-Szego}, there exists a positive integer $D$ such that the zeros of $\det \left(\widehat{Q}(\cdot)\right)$ have multiplicity at most $D$.

Assume that the condition in Proposition~\ref{prop_nece_bezout}\ref{cdt_corona} holds true. Let $s \in \mathbb{C}$ be such that $\det \widehat{Q}(s)=0$ and denote by $\ell$ the dimension of the left kernel of $\widehat{Q}(s)$. Let $G_{s} \in \mathcal{M}_{\ell,d}(\mathbb{C})$ be such that its rows form an orthonormal basis of the left kernel of $ \widehat{Q}(s)$. Hence, $G_{s} \widehat{Q}(s)=0$ and $\abs{z^\transp G_s}= \abs{z}$ for every $z \in \mathbb C^\ell$.

We claim that the rows of the matrix $G_s \widehat{P}\in \mathcal{M}_{\ell,m}(\mathbb{C})$ are linearly independent. Indeed, if this were not the case, there would exist $z \in \mathbb{C}^{\ell}$ such that $\abs{z} = 1$ and $z^\transp G_s\widehat{P}=0$.
Thus we would have $z^\transp G_s \widehat{Q}(s)=0$, $z^\transp G_s\widehat{P}=0$, and $\abs{z^\transp G_s}=1$, in contradiction with Proposition~\ref{prop_nece_bezout}\ref{cdt_corona}.

We deduce that $G_{s}\widehat{P}$ has right inverses and we use $(G_{s}\widehat{P})^{-1}$ to denote its Moore--Penrose right inverse, i.e., $(G_{s}\widehat{P})^{-1} = (G_{s}\widehat{P})^\transp (G_{s}\widehat{P} \widehat{P}^\transp G_s^\transp)^{-1}$. Define
\begin{equation*}
\widehat{\Phi}(s)=\left(G_{s}\widehat{P}\right)^{-1} G_{s}
\end{equation*}
and notice that $\widehat{\Phi}(s)$ satisfies
\begin{equation}
\label{eq:centrale2}
G_s\widehat{P}\widehat{\Phi}(s)=G_s.
\end{equation}

We now claim that there exists $c>0$ (independent of the zero $s$ of $\det \widehat{Q}$) such that
\begin{equation}
\label{eq:centrale3}
\abs*{\widehat{\Phi}(s)} \le c.
\end{equation}
Indeed, arguing by contradiction yields a sequence $(s_n)_{n \in \mathbb{N}}$ of zeros of $\det \widehat Q$ such that the norm $\abs*{\left(G_{s_n} \widehat P\right)^{-1}}$ tends to infinity as $n \to +\infty$. We have
\[\abs*{\left(G_{s_n} \widehat P\right)^{-1}} = \frac{1}{\lambda_{\min}(G_{s_n}\widehat{P} \widehat{P}^\transp G_{s_n}^\transp)},\]
 where $\lambda_{\min}(\cdot)$ denotes the smallest eigenvalue of its argument. Let $(z_n)_{n \in \mathbb N}$ denote a sequence of vectors in $\mathbb C^\ell$ such that, for every $n \in \mathbb N$, $\abs{z_n^\transp} = 1$ and $z_n$ is an eigenvector of $G_{s_n}\widehat{P} \widehat{P}^\transp G_{s_n}^\transp$ associated with its smallest eigenvalue. Then
\[
\abs*{z_n^\transp G_{s_n} \widehat P}^2 = \lambda_{\min}(G_{s_n}\widehat{P} \widehat{P}^\transp G_{s_n}^\transp) \xrightarrow[n \to +\infty]{} 0,
\]
and, since $\abs{z_n^\transp G_{s_n}}= 1$ for every $n \in \mathbb N$, and $G_{s_n} \widehat Q(s_n) = 0$, we obtain a contradiction with Proposition~\ref{prop_nece_bezout}\ref{cdt_corona}, yielding \eqref{eq:centrale3}.

Assumptions (18)\footnote{Equation~(18) in \cite{Yamamoto2011Bezout} is the Laplace transform of (8) of the same article, which is an equation in $\left(\mathcal{E}^\prime(\mathbb{R}) \right)^{d \times d} / \left( Q \right)$, where $(Q)$ is the ideal generated over $\mathcal{E}^\prime(\mathbb{R}_-)$ by $Q$. Thus the equality in (18) has to be understood as an equality modulo the ideal generated by $\widehat Q(s)$, which is equivalent to requiring that $z^\transp \left(\widehat P(s) \widehat \Phi(s) - \Id_d\right) = 0$ for every $z^\transp$ in the left null space of $\widehat{Q}(s)$. Therefore, \cite[(18)]{Yamamoto2011Bezout} is equivalent to \eqref{eq:centrale2} above.} and (17) in the statement of \cite[Theorem~5.1]{Yamamoto2011Bezout} are satisfied thanks to \eqref{eq:centrale2} and \eqref{eq:centrale3}, respectively. We can then apply  \cite[Theorem~5.1]{Yamamoto2011Bezout},
which ensures the existence of two distribution $R$ and $S$ with entries in $\mathcal{E}^\prime(\mathbb{R}_{-})$ such that
\begin{equation*}
Q*R+P*S=\delta_0 \Id_d.
\end{equation*}

\step{2} Since $S$ is a distribution with compact support, it has a finite order $k$. Let $y \in  X^{Q}_{k}$ and denote by $\widetilde{y} \in \mathcal C^k(\mathbb{R},\mathbb{R}^d)$ an extension on $\mathbb{R}$ of $y$ having a support bounded on the left. Set $\omega=S*Q*\widetilde{y}$. Then $\omega$ belongs to $\mathcal C^0(\mathbb{R},\mathbb{R}^m)$ and it has a compact support included in $\mathbb{R}_-$. It follows from a similar argument given to obtain the equation \eqref{eq:sol-con-B} that $y=\pi\left(Q^{-1}*P*\omega\right)$, which proves the $\mathcal C^k$ exact controllability.
\end{proof}

Finally, the proof of Theorem~\ref{thm:Ck} can be obtained as a consequence of Proposition~\ref{prop-Ckexact}, as we now detail.

\begin{proof}[Proof of Theorem~\ref{thm:Ck}]
Notice that $y \in X_{k}^Q$ if and only if $[-\Lambda_N,0]\ni t\mapsto y(t+\Lambda_N)$ is $\mathcal C^k$-admissible for \eqref{eq:intro-linear-difference-equation} in the sense of Definition~\ref{def:compatibility}. Recall that, by Proposition~\ref{Prop2:main_result}, \ref{assumption1-conject} and \ref{assumption2-conject} from Conjecture~\ref{conj:HY-exact} imply that the condition from Proposition~\ref{prop_nece_bezout}\ref{cdt_corona} holds true. Hence, by Proposition~\ref{prop-Ckexact}, every $\mathcal C^k$-admissible function for \eqref{eq:intro-linear-difference-equation} is in the range of $E_p(T)$ for some $T>0$ possibly depending on $y$. By Theorem~\ref{lem:RanE_con}, moreover, the space of $\mathcal C^k$-admissible functions for \eqref{eq:intro-linear-difference-equation} is contained in $\Ran E_p(d\Lambda_N)$. The conclusion then follows by Proposition~\ref{prop_var_const_formula} and the fact that, if $x_0$ is $\mathcal C^k$-admissible for \eqref{eq:intro-linear-difference-equation}, then $\Upsilon_p(t)x_0$ is $\mathcal C^k$-admissible for \eqref{eq:intro-linear-difference-equation} for every $t\ge 0$.
\end{proof}

\section{Applications}
\label{sec:applications}

In this section, we show in simple cases how to derive from Theorems~\ref{main_result1} and \ref{main_result2} Kalman-type conditions for controllability in the single input case $m=1$. By a Kalman-type condition, we refer to a  frequency-free criterion for controllability, i.e., independent of the complex variable $s$. For instance, in the single-delay case $N=1$, approximate and exact controllability (in any fixed time $T \geq d \Lambda_1$) coincide and are equivalent to the classical \emph{Kalman rank condition} stating that the rank of the controllability matrix $\begin{pmatrix}B & A_1B & \cdots & A_1^{d-1}B\end{pmatrix}$ is equal to the state space dimension $d$. We first partially recover criteria given in \cite{Chitour2020Approximate} in the case of two delays, two space dimensions, and a single input, and then extend such a study to the case of two delays, three space dimensions, and a single input.

Up to a time rescaling, we can assume from now on that $(\Lambda_1,\Lambda_2)=(L,1)$ with $L \in (0,1)$. The case $L \in \mathbb Q$ was completely addressed in \cite{Chitour2020Approximate} and we assume for the rest of this section that $L$ is irrational.

\subsection{Two delays, two space dimensions, and a single input}
\label{sec:N=D=2,m=1}

In this section, we recover some results given in \cite[Theorem~4.1]{Chitour2020Approximate}, which  concern the particular case where $N=d=2$ and $m=1$.
The results in \cite{Chitour2020Approximate} are stated for matrices with complex coefficients, hence, what we actually recover here are some results of \cite[Theorem~4.1]{Chitour2020Approximate} restricted to the case of matrices with real coefficients. On the other hand, the results in \cite{Chitour2020Approximate} are stated only for $L^2$ controllability, but here we consider $L^p$ controllability for any $p \in [1, +\infty)$.

Fix then $A_1,A_2\in \mathcal{M}_{2,2}(\mathbb{R})$ and  $B\in \mathcal{M}_{2,1}(\mathbb{R})$. Given a matrix $A\in \mathcal{M}_{2,2}(\mathbb{R})$, we say that the pair $(A,B)$ is controllable if it satisfies the Kalman rank condition.

We split our analysis into three cases.

\case{1}{$\Ran A_2 \subset \Ran B$ or both pairs $(A_1,B)$, $(A_2,B)$ are not controllable} In both subcases, \eqref{eq:intro-linear-difference-equation} is not even $L^p$ approximately controllable. Indeed, in the first subcase, Theorems~\ref{main_result1}\ref{assumption2bis} and~\ref{main_result2}\ref{assumption2-L1} are not satisfied. In the second subcase, up to a linear change of variables, one can assume with no loss of generality that
\[
A_1=\begin{pmatrix}a_1&0\\a_2&a_3
\end{pmatrix}, \qquad
A_2=\begin{pmatrix}b_1&0\\b_2&b_3
\end{pmatrix}, \qquad
B=\begin{pmatrix}0\\1
\end{pmatrix},
\]
with $b_1\neq 0$ (otherwise we are back to the first subcase). Then the first coordinate $x_1$ of the state $x$ is not controllable since it satisfies the autonomous dynamics $x_1(t) = a_1 x_1(t - L) + b_1 x_1(t - 1)$.

\case{2}{$\Ran A_2  \not\subset \Ran B$ and exactly one of the pairs $(A_1,B)$, $(A_2,B)$ is controllable} In this case, \eqref{eq:intro-linear-difference-equation} is $L^p$ approximately controllable in time $2 \Lambda_2$ for every $p \in [1, +\infty)$. Indeed, notice first that $\rank\begin{pmatrix}A_2 & B\end{pmatrix} = 2$, and we are thus left to show that the condition in Theorem~\ref{main_result1}\ref{assumption1bis} holds true, which is equivalent to proving that the condition in Proposition~\ref{Prop2:main_result0}\ref{app-b} holds true. Assume, for instance, that $(A_1,B)$ is controllable and $(A_2,B)$ is not. Hence, up to a linear change of variables, we can assume with no loss of generality that
\[
A_1=\begin{pmatrix}a_1&a_2\\a_3&a_4
\end{pmatrix}, \qquad
A_2=\begin{pmatrix}b_1&0\\b_2&b_3
\end{pmatrix}, \qquad
B=\begin{pmatrix}0\\1
\end{pmatrix},
\]
with $a_2 \neq 0$ and $b_1 \neq 0$. For $g \in \mathbb C^2$ with $\abs{g} = 1$, one has either $g^\transp B \neq 0$ or $g^\transp = (\alpha, 0)$ for some $\alpha \in \mathbb C$ with $\abs{\alpha} = 1$. Hence, in order to show Proposition~\ref{Prop2:main_result0}\ref{app-b}, it suffices to show that $\abs{g^\transp H(s)}$ is nonzero for every $s \in \mathbb C$ and with $g^\transp = (1, 0)$. For every $s \in \mathbb C$, one checks that the second coordinate of $g^\transp H(s)$ is equal to $-a_2 e^{-s L}$, which never vanishes,
yielding the conclusion. The case where $(A_1,B)$ is not controllable and $(A_2,B)$ is can be handled similarly.

\case{3}{$(A_1,B)$ and $(A_2,B)$ are both controllable} Let $B^{\bot} \in \mathbb{R}^2$ be the unique vector such that $\det(B,B^{\bot})=1$ and $B^\transp B^{\bot}=0$. Set
\[
\beta= \frac{\det \begin{pmatrix}B & A_2B \end{pmatrix}}{\det \begin{pmatrix}B & A_1B \end{pmatrix}}, \qquad \alpha= \det  \begin{pmatrix}B & (A_2-\beta A_1)B^{\bot}\end{pmatrix}.
\]
Up to a linear change of coordinates, we can assume that
\begin{equation}
\label{eq:new_coordinates-221-case3}
A_1=\begin{pmatrix}
0 & 1  \\
a_1 & a_2
\end{pmatrix}, \qquad A_2=\begin{pmatrix}
\alpha & \beta  \\
b_1 & b_2
\end{pmatrix},  \qquad B=\begin{pmatrix}
0   \\
1
\end{pmatrix}.
\end{equation}

The holomorphic map $H$ is now given by
\begin{equation}
\label{eq_221}
H(s)=\Id_2- e^{-sL} A_1 - e^{- s}A_2=\begin{pmatrix}
1-\alpha e^{-s} & -e^{-sL}-\beta e^{-s}  \\
\ast & \ast
\end{pmatrix},\qquad s \in \mathbb{C}.
\end{equation}
Let $\alpha = \abs{\alpha} e^{i \theta}$ with $\theta \in \{0,\pi\}$. Since $\beta \neq 0$, we have that $\rank\begin{pmatrix}A_2 & B\end{pmatrix} = 2$. Then, by Theorem~\ref{main_result1} and \eqref{eq_221}, we have that \eqref{eq:intro-linear-difference-equation} is not $L^p$ approximately controllable in time $2\Lambda_2$ if and only if
\begin{align*}
\exists s\in\mathbb{C}\text{ s.t.\ } \rank\begin{pmatrix}{H}(s) & B\end{pmatrix} < 2 &\iff \exists s\in\mathbb{C}\text{ s.t.\ }1-\alpha e^{-s}=0 \text{ and  } e^{-sL}+\beta e^{-s}=0\\
&\iff 0 \in S,
\end{align*}
where $S = \left\{ \beta+ \abs{\alpha}^{1-L}e^{i (\theta+2k \pi)(1-L)} \suchthat  k \in \mathbb{Z} \right\}$. Since $L$ is irrational, notice that $\overline S$ is the circle in $\mathbb C$ of center $\beta$ and radius $\abs{\alpha}^{1-L}$, denoted hereafter by $C$. By Theorem~\ref{main_result2} and \eqref{eq_221}, one can prove in the same way that \eqref{eq:intro-linear-difference-equation} is not $L^1$ exactly controllable in time $2 \Lambda_2$ if $0 \in C$.

\medskip

In the third case above, the fact that there exists $s\in\mathbb{C}$ so that $1-\alpha e^{-s}=0$ and  $e^{-sL}+\beta e^{-s}=0$ can be equivalently written as
\begin{equation*}
M Y(s)=\begin{pmatrix}1\\0\end{pmatrix},\qquad
\hbox{ where }
M=\begin{pmatrix}0&\alpha\\1&\beta\end{pmatrix}
\hbox{ and } Y(s)=\begin{pmatrix}e^{-sL}\\e^{-s}\end{pmatrix}.
\end{equation*}

The above computations are a particular case of a more general situation, described next.

\begin{lemma}
\label{lem:M-Y-Z}
Let $M$ be a $2 \times 2$ matrix with complex coefficients, $Z \in \mathbb C^2 \setminus \{0\}$, $L \in (0, 1)\setminus \mathbb{Q}$, and define the map $Y\colon s \mapsto (e^{-s L}, e^{-s})^\transp$.
\begin{enumerate}
\item If $M$ is not invertible, then the following assertions are equivalent.
\begin{enumerate}
\item $Z \in \Ran M$.
\item $Z \in M Y(\mathbb C)$.
\item $Z \in M \overline{Y(\mathbb C)}$.
\end{enumerate}

\item If $M$ is invertible, set $M^{-1} Z = (a_1, a_2)^\transp$ and define
\[
S = \left\{a_1 - \abs{a_2}^L e^{i L (\theta_2 + 2 k \pi)} \suchthat k \in \mathbb Z\right\},
\]
where $a_2 = \abs{a_2} e^{i \theta_2}$ and $\theta_2 \in \mathbb R$. Then the following assertions hold true.
\begin{enumerate}
\item $Z \in M Y(\mathbb C)$ if and only if $0 \in S$.
\item $Z \in M \overline{Y(\mathbb C)}$ if and only if $0 \in C$, where $C$ is the circle in $\mathbb C$ of center $a_1$ and radius $\abs{a_2}^L$.
\end{enumerate}
\end{enumerate}
\end{lemma}

Our two main results allowed us to recover partially, in the case of matrices with real coefficients, the algebraic controllability result stated in \cite[Theorem~4.1]{Chitour2020Approximate} when $N=d=2$ and $m=1$. The latter considers $L^2$ (approximate and exact) controllability concepts while, in the present section, we are instead dealing with $L^p$ approximate controllability and $L^1$ exact controllability. We have shown that three possibilities occur:
\begin{enumerate*}
\item neither approximate nor exact controllability hold true;
\item $L^p$ approximate controllability holds true but we cannot say anything about  exact controllability;
\item there exists a countable subset $S$ of the complex plane, completely characterized by $A_1,A_2,B$, and $L$, which is dense in a circle $C$, such that $L^p$ approximate controllability is equivalent to the fact that $0\notin S$ and a necessary condition for $L^1$ exact controllability is $0\notin C$.
\end{enumerate*}

\begin{remark}
For Case~3, we can assume with no loss of generality (see \cite[Lemma~4.5]{Chitour2020Approximate}) that $a_1=a_2=b_1=b_2=0$ in \eqref{eq:new_coordinates-221-case3}. Using Theorem~\ref{prop_fond_exact_controllability_L1}, we have that \eqref{eq:intro-linear-difference-equation} is $L^1$ exactly controllable if and only if Bézout's identity \eqref{eq:exact_cont0} is solvable in $M(\mathbb R_-)$, that is, if and only if there exists $(s_1,s_2,r_1,r_2,r_3,r_4) \in M(\mathbb{R}_-)^6$ such that
\begin{equation}
\label{cas_particulier_bezout_main}
\left\{
\begin{aligned}
 & (\delta_{-1} - \alpha \delta_0) * r_1 - (\delta_{L - 1} + \beta \delta_0) * r_2 = \delta_0, \\
 & (\delta_{-1} - \alpha \delta_0) * r_3 - (\delta_{L - 1} + \beta \delta_0) * r_4 = 0, \\
 & \delta_{-1} * r_2 + s_1 = 0, \\
 & \delta_{-1} * r_4 + s_2 = \delta_0.
\end{aligned}
\right.
\end{equation}
The second and fourth equations of \eqref{cas_particulier_bezout_main} are satisfied if one chooses
\[
r_3 = \delta_{L - 1} + \beta \delta_0, \qquad
r_4 = \delta_{-1} - \alpha \delta_0, \qquad
s_2 = \delta_0 + \alpha \delta_{-1} - \delta_{-2}.
\]
In addition, as soon as the first equation of \eqref{cas_particulier_bezout_main} is satisfied, the third one can be satisfied by setting $s_1 = - \delta_{-1} * r_2$, and we hence focus in the first equation of \eqref{cas_particulier_bezout_main} in the sequel, i.e.,
\begin{equation}
\label{cas_particulier_bezout_main2}
(\delta_{-1} - \alpha \delta_0) * r_1 - (\delta_{L - 1} + \beta \delta_0) * r_2 = \delta_0.
\end{equation}

Let $q_1 = \delta_{-1} - \alpha \delta_0$ and $q_2 = -\delta_{L - 1} - \beta \delta_0$, and notice that their Laplace transforms are given by $\widehat q_1(s) = e^s - \alpha$ and $\widehat q_2(s) = -e^{(1-L)s} - \beta$ for $s \in \mathbb C$. By definition of $S$, we have that $0 \in S$ if and only if there exists $s \in \mathbb C$ such that $\widehat q_1(s) = \widehat q_2(s) = 0$, and one can show that $0 \in \overline S$ if and only if $\inf_{s \in \mathbb C} \abs{\widehat q_1(s)} + \abs{\widehat q_2(s)} = 0$. Hence, the condition $0 \notin \overline S$ is equivalent to the existence of a constant $c > 0$ such that $\abs{\widehat q_1(s)} + \abs{\widehat q_2(s)} \geq c$ for every $s \in \mathbb C$. If Conjecture~\ref{conjecture1} held true, then it would imply that there exist $r_1$ and $r_2$ in $M(\mathbb R_-)$ such that their Laplace transforms satisfy $\widehat q_1(s) \widehat r_1(s) + \widehat q_2(s) \widehat r_2(s) = 1$, which is exactly the Laplace transform of \eqref{cas_particulier_bezout_main2}.
In particular, this shows that, if $0 \notin \overline S$ and Conjecture~\ref{conjecture1} holds true, then \eqref{eq:intro-linear-difference-equation} in the case of the present example is $L^1$ exact controllable.

Note that the condition $0 \notin \overline S$ is shown in \cite{Chitour2020Approximate} to be a necessary and sufficient condition for $L^2$ exact controllability of \eqref{eq:intro-linear-difference-equation} in Case~3. Unfortunately, even for this simple example, it is not clear if $L^2$ exact controllability is equivalent to $L^1$ exact controllability.
\end{remark}

\subsection{Two delays, three space dimensions, and a single input: a geometric locus controllability result}

In this section, we assume that $N=2$, $d=3$, and $m=1$. In this case $A_1,A_2\in \mathcal{M}_{3,3}(\mathbb{R})$ and $B\in  \mathcal{M}_{3,1}(\mathbb{R})$. Up to a linear change of variables, we assume in the sequel that $B=(0,0,1)^\transp$. We use $r_0$ to denote the dimension of the vector space spanned by the three vectors $B, A_1 B, A_2 B$.

As in Section~\ref{sec:N=D=2,m=1}, we subdivide the discussion in three cases.

\case{1}{$\rank\begin{pmatrix}A_2 & B\end{pmatrix}\leq 2$ or $r_0=1$} In both subcases, \eqref{eq:intro-linear-difference-equation} is not even approximately controllable. Indeed, in the first subcase, the rank condition of Theorem~\ref{main_result1}\ref{assumption2bis} is violated. In the second subcase, the matrices $A_1$ and $A_2$ are of the form
\[
A_1=\begin{pmatrix}
\widetilde{A}_1 &0_{2,1}\\
\ast &a_1
\end{pmatrix},
\qquad A_2=\begin{pmatrix}
\widetilde{A}_2&0_{2,1}\\
\ast &a_2
\end{pmatrix}.
\]
Then the subsystem of \eqref{eq:intro-linear-difference-equation}  made of the first two coordinates is uncoupled to the third coordinate and is uncontrolled. The overall system is then
not controllable.

\case{2}{$\rank\begin{pmatrix}A_2 & B\end{pmatrix} = 3$ and $r_0=2$}
We only treat the case where $A_2B$ is not colinear to $B$, the other case being entirely similar. Up to a linear change of variables (in a basis having $A_2B$ as second vector and $B$ as third one), one can transform the matrices $A_1$ and $A_2$ into the form
\[
A_1 = \left(\begin{array}{cc}
\multirow{2}{*}{$\widetilde A_1$} & 0 \\
& \alpha \\
* & * \\
\end{array}\right),
\qquad A_2 = \left(\begin{array}{cc}
\multirow{2}{*}{$\widetilde{A}_2$} & 0\\
& 1 \\
\ast & 0 \\
\end{array}\right),
\]
where the matrices $\widetilde{A}_1$ and $\widetilde{A}_2$ are $2\times 2$.

The holomorphic map $H$ associated with \eqref{eq:intro-linear-difference-equation} is equal to
\begin{equation}
\label{eq_231}
H(s) = \left(\begin{array}{cc}
\multirow{2}{*}{$\Id_2- e^{-sL} \widetilde{A}_1 - e^{- s}\widetilde{A}_2$} & 0\\
& -(e^{-sL}\alpha+e^{-s}) \\
\ast & \ast \\
\end{array}\right).
\end{equation}
Then, by Theorem~\ref{main_result1} and \eqref{eq_231}, \eqref{eq:intro-linear-difference-equation} is not $L^p$ approximately controllable in time $3\Lambda_2$ if and only if there exists $s\in\mathbb{C}$ such that
$\rank\begin{pmatrix}{H}(s) & B\end{pmatrix} < 3$, i.e., either
\begin{enumerate}
\item\label{dim3-case2a} $e^{-sL}\alpha+e^{-s}=0$ and $\det(\Id_2- e^{-sL} \widetilde{A}_1 - e^{- s}\widetilde{A}_2)=0$
\end{enumerate}
or
\begin{enumerate}[resume]
\item\label{dim3-case2b} $e^{-sL}\alpha+e^{-s}\neq 0$ and the $(1,1)$ and $(1,2)$ coefficients of the $3\times 3$ matrix $H(s)$ are both equal to zero.
\end{enumerate}

In the subcase \ref{dim3-case2a}, the second condition is equivalent to the fact that
$\det\bigl(e^{sL}\Id_2-(\widetilde{A}_1 - \alpha\widetilde{A}_2)\bigr)=0$, i.e., $e^{sL}$ is
one of the eigenvalues of $\widetilde{A}_1 - \alpha\widetilde{A}_2$. Hence we are in the situation of Lemma~\ref{lem:M-Y-Z} with $M = \begin{pmatrix}\alpha & 1 \\ 1 & 0\end{pmatrix}$ and at most two vectors $Z$ corresponding to the nonzero eigenvalues of $\widetilde{A}_1 - \alpha\widetilde{A}_2$.
It follows that there exist at most two countable sets of complex numbers $S_1, S_2$, completely characterized in terms of $L$ and the coefficients of $A_1, A_2$, which are dense in two circles $C_1, C_2$, respectively,
and such that \eqref{eq:intro-linear-difference-equation} is $L^p$ approximately controllable if and only if $0$ does not belongs to the union of $S_1$  and $S_2$, while \eqref{eq:intro-linear-difference-equation} is not $L^1$ exactly controllable if $0$ belongs to the union of $C_1$  and $C_2$. (Note that each $S_i$ can be equal to $C_i$ if the latter reduces to a point.) It is not difficult to see that the subcase \ref{dim3-case2b} also boils down to a similar situation, but with at most one circle.

\case{3}{$\rank\begin{pmatrix}A_2 & B\end{pmatrix} = 3$ and $r_0=3$}
Up to a linear change of variables (in a basis having $A_1B$ as first vector, $A_2B$ as second vector, and $B$ as third one), one can transform the matrices $A_1$ and $A_2$ into the form
\[
A_1 = \left(\begin{array}{cc}
\multirow{2}{*}{$\widetilde{A}_1$} & 1 \\
& 0 \\
\ast & 0 \\
\end{array}\right),
\qquad A_2 = \left(\begin{array}{cc}
\multirow{2}{*}{$\widetilde{A}_2$} & 0 \\
& 1 \\
\ast & 0 \\
\end{array}\right),
\]
where the matrices $\widetilde{A}_1$ and  $\widetilde{A}_2$ are $2\times 2$.

The holomorphic map $H$ associated with \eqref{eq:intro-linear-difference-equation} is equal to
\begin{equation*}
H(s) = \left(\begin{array}{cc}
\multirow{2}{*}{$\Id_2- e^{-sL} \widetilde{A}_1 - e^{- s}\widetilde{A}_2$} & -e^{-sL} \\
& -e^{-s} \\
\ast & 0 \\
\end{array}\right).
\end{equation*}

We denote by $v_1$ and $v_2$ the columns of the matrix $\widetilde{H} (s)=\Id_2-e^{-sL} \widetilde{A}_{1}-e^{-s} \widetilde{A}_{2}$ and we set $v_3=\begin{pmatrix}e^{-sL} & e^{-s}\end{pmatrix}^\transp$. Note that the rank assumption on $\begin{pmatrix}A_2 & B\end{pmatrix}$ implies that the first row of $\widetilde A_2$ is not equal to zero.

It is immediate to see that there exists $s\in\mathbb{C}$ so that
$\rank\begin{pmatrix}{H}(s) & B\end{pmatrix} < 3$ if and only if both $v_1$ and $v_2$ are colinear to the nonzero vector $v_3$, i.e.,
$\det(v_1,v_3)=0$ and $\det(v_2,v_3)=0$. These conditions can be rewritten as two scalar equations
$v_3^\transp Q_1 v_3+e^{-s} = 0$ and $v_3^\transp Q_2v_3+e^{-s L}=0$, where $Q_1,Q_2$ are symmetric matrices with real coefficients and at least one between $Q_1$ and $Q_2$ is not equal to zero because of the rank assumption on $\begin{pmatrix}A_2 & B\end{pmatrix}$.

The two previous equations in the unknowns $e^{-sL}$ and $e^{-s}$ define two distinct
conic sections, at least one of them being nontrivial. Therefore, they have $k$ distinct intersection points, with $k \in \llbracket 0, 4\rrbracket$. For each of them, one can completely characterize, in terms of $L$ and the coefficient of $A_1,A_2$, a countable set of complex numbers $S_j$, $j \in \llbracket 1, k\rrbracket$, which is dense in a circle $C_j$, such that \eqref{eq:intro-linear-difference-equation} is $L^p$ approximately controllable if and only if $0$ does not belong to the union of the $S_j$, while \eqref{eq:intro-linear-difference-equation} is not $L^1$ exactly controllable if $0$ does belong to the union of the $C_j$.

\medskip

In conclusion, if $N=2$, $d=3$, and $m=1$, Theorem~\ref{main_result1} (respectively, Theorem~\ref{main_result2}) allows one to derive frequency-free necessary and sufficient (respectively, necessary) criteria for $L^p$ approximate (respectively, $L^1$ exact) controllability, and the results
are qualitatively similar to those obtained if $N=d=2$ and $m=1$, with the difference in Case~3 where now we may have up to four distinct countable sets $S_j$, $k \in \llbracket 1, k\rrbracket$, each of them dense in a circle $C_j$.

\chapter{Strong stability of linear difference equations}
\label{chapter:strong-stability}

\epigraph{\foreignlanguage{brazil}{Ao vencedor, as batatas!}}{\foreignlanguage{brazil}{\emph{Quincas Borba}, Machado de Assis}}

This chapter, which presents the main results from \cite{GoncalvesNetto2025Strong, ChitourStrong}, considers the sensitivity analysis of exponential stability of difference equations of the form \eqref{eq:intro-diff-eqn-integral-M} with respect to perturbations in the matrix-valued function $M$ defining the system. We prove that exponential stability of \eqref{eq:intro-diff-eqn-integral-M} is preserved under small enough perturbations of $M$ in total variation norm, at least for an important class of functions $M$ (namely, those without a singular part), generalizing thus the corresponding result for the particular case of \eqref{eq:intro-diff-eqn-sum} for perturbations in the matrices $A_1, \dotsc, A_N$ only, mentioned in Section~\ref{sec:sensitivity}. We then proceed to consider a generalization of the Hale--Silkowski criterion from Proposition~\ref{prop:Hale-Silkowski}. For that purpose, we first put forward an appropriate definition of ``perturbation on the delays'' for \eqref{eq:intro-diff-eqn-integral-M} and then propose a conjecture in the spirit of Proposition~\ref{prop:Hale-Silkowski}, which we prove in dimension $1$ or under a slightly modified form and with additional restrictions on $M$.

More precisely, we shall consider in this chapter the linear difference equation
\begin{equation}
\label{sys_general}
\Sigma(M)\colon\qquad x(t) = \int_{-1}^0 \diff M(\theta) x(t+\theta), \qquad t \in \mathbb R_+,
\end{equation}
where $x(t)\in \mathbb{R}^d$ and $M\colon [-1,0]\to \mathcal{M}_d(\mathbb{R})$ is a matrix-valued function with entries of bounded variation, assumed to be continuous from the right on $(-1,0)$, and the integral in the right-hand side of \eqref{sys_general} is taken in the Riemann--Stieltjes sense. System~$\Sigma(M)$ corresponds to \eqref{eq:intro-diff-eqn-integral-M} with $r = 1$, and this normalization is not a restriction since \eqref{eq:intro-diff-eqn-integral-M} can be reduced to \eqref{sys_general} by the change of time scale $y(t) = x(r t)$.

The chapter is organized as follows. Section~\ref{sec:well-posed-strong-stability} provides the definition of solution of \eqref{sys_general} used in this chapter, discusses an alternative equivalent formulation of \eqref{sys_general} in terms of a matrix-valued Radon measure, and provides a result on the existence and uniqueness of solutions of \eqref{sys_general}, Proposition~\ref{prop_wellposed}, under rather general assumptions. We then discuss, in Section~\ref{sec:definition-spaces}, the exponential behavior of \eqref{sys_general} and its relation with spectral properties of the system, providing the definitions of several functional spaces to be used in the sequel. Our first main result, Theorem~\ref{thm:spectral-abscissa-continuous}, dealing with the continuity of the spectral abscissa of \eqref{sys_general} with respect to $M$ in the topology of the total variation norm, is stated and proved in Section~\ref{sec_tv-norm}, and, as a consequence, we obtain, in Corollary~\ref{coro_tvstab}, that the exponential stability of System~\eqref{sys_general} is robust with respect to perturbations of $M$ in total variation norm. Section~\ref{sec_strong} considers a different kind of perturbation of $M$, generalizing the notion of perturbation of the delays of \eqref{eq:intro-diff-eqn-sum} to \eqref{sys_general}. After extending the definition of strong stability of \eqref{eq:intro-diff-eqn-sum} to \eqref{sys_general}, we state and prove our main results on this notion: Theorem~\ref{HS:scalar}, which is a generalization of the Hale--Silkowski criterion to \eqref{sys_general} in the scalar case, and Theorem~\ref{GHS-Theorem}, which concerns the vector-valued case. Numerical illustrations of our results are provided in Section~\ref{sec_numerical}, and the chapter is concluded by Appendix~\ref{sec:app}, containing some useful measure-theoretical or technical results.

\paragraph*{Notation.} In this chapter, given a bounded interval $[a, b] \subset \mathbb R$ with $a < b$ and a bounded function $\varphi\colon [a, b] \to \mathbb R^d$, we use $\norm{\varphi}_\infty$ to denote the quantity $\sup_{t \in [a, b]} \abs{\varphi(t)}$, where $\abs{\cdot}$ is an arbitrary norm in $\mathbb R^d$. The $\sigma$-algebra of all Borel subsets of the interval $[a, b]$ is denoted by $\mathfrak{B}_{[a,b]}$.

Although we consider in this chapter real solutions of \eqref{sys_general} only, due to our use of spectral methods, we also need to deal with complex vectors and matrices, and we shall thus need norms of both real and complex vectors and matrices. Given a norm $\abs{\cdot}_{\mathbb C^d}$ in $\mathbb C^d$, one can define a norm $\abs{\cdot}_{\mathbb R^d}$ in $\mathbb R^d$ as the restriction of $\abs{\cdot}_{\mathbb C^d}$ to $\mathbb R^d$, and then consider the respective induced matrix norms on $\mathcal{M}_d(\mathbb{C})$ and $\mathcal M_d (\mathbb{R})$, denoted by $\abs{\cdot}_{\mathcal M_d(\mathbb C)}$ and $\abs{\cdot}_{\mathcal M_d(\mathbb R)}$. We always assume in this chapter that the norm $\abs{\cdot}_{\mathbb C^d}$ is chosen in such a way that the following compatibility condition is satisfied: we have $\abs{A}_{\mathcal M_d(\mathbb C)} = \abs{A}_{\mathcal M_d(\mathbb R)}$ for every $A \in \mathcal M_d(\mathbb R)$. It is easy to find examples of norms in $\mathbb C^d$ such that this property is not satisfied. However, given any norm $\abs{\cdot}_{\mathbb R^d}$ in $\mathbb R^d$, one can extend it to a norm in $\mathbb C^d$ satisfying the compatibility condition by setting, for $x \in \mathbb{C}^d$,
\[\abs{x}_{\mathbb{C}^d}=\max_{z\in \mathbb{C}, \abs{z}\leq 1} \abs{\Real(zx)}_{\mathbb{R}^d},\]
where the real part of a complex vector should be understood as the componentwise real part. This is indeed an extension because $\abs{\Real(zx)}_{\mathbb{R}^d}=\abs{\Real(z)}\abs{x}_{\mathbb{R}^d}$ for $x\in \mathbb{R}^d$, and one can easily check that $\abs{\cdot}_{\mathbb C^d}$ satisfies the required compatibility condition. From now on, we use $\abs{\cdot}$ to denote any of the norms $\abs{\cdot}_{\mathbb{R}^d}$, $\abs{\cdot}_{\mathbb{C}^d}$, $\abs{\cdot}_{\mathcal M_d(\mathbb{R})}$, and $\abs{\cdot}_{\mathcal M_d(\mathbb{C})}$.

\section{Main system and its solutions}
\label{sec:well-posed-strong-stability}

Let us start by providing the definition of solution of System~\eqref{sys_general} that we will use in this chapter.

\begin{definition}
Given $\delta \in \mathbb R_+^\ast \cup \{+\infty\}$ (respectively, $\delta \in \mathbb R_+^\ast$), a function $x \in \mathcal C([-1,\delta), \mathbb{R}^d)$ (respectively, $x \in \mathcal C([-1,\delta], \mathbb{R}^d)$) is said to be a \emph{local solution} of \eqref{sys_general} in $[-1, \delta)$ (respectively, in $[-1, \delta]$) if $x$ satisfies \eqref{sys_general} for all $t \in [0, \delta)$ (respectively, $t \in [0, \delta]$). We say that $x$ is a \emph{global solution} of \eqref{sys_general} if it is a local solution in $[-1, +\infty)$. The restriction of $x$ to the interval $[-1, 0]$ is called the \emph{initial condition} of $x$.
\end{definition}

As usual with time-delay systems, given a function $x \in \mathcal C([-1, \delta), \mathbb{R}^d)$ for some $\delta \in \mathbb{R}_+^\ast \cup \{+\infty\}$ and $t \in [0, \delta)$ (or $x \in \mathcal C([-1, \delta], \mathbb{R}^d)$ for some $\delta \in \mathbb{R}_+^\ast$ and $t \in [0, \delta]$), the \emph{history function} of $x$ at time $t$ is the function $x_t \in \mathcal C([-1, 0], \mathbb{R}^d)$ defined by
\[
x_t(\theta) = x(t+\theta), \qquad \theta \in [-1, 0].
\]
With this notation, the initial condition of a solution $x$ of \eqref{sys_general} is simply $x_0$.

Note that, if $x$ is a solution of \eqref{sys_general}, then the initial condition $x_0$ necessarily satisfies
\begin{equation*}
x_0(0) = \int_{-1}^0 \diff M(\theta) x_0(\theta).
\end{equation*}
We will thus consider in the sequel only initial conditions in the set $\mathcal C_0$ defined by
\[
\mathcal C_0 = \left\{\eta \in \mathcal C([-1, 0], \mathbb{R}^d) \suchthat \eta(0) = \int_{-1}^0 \diff M(\theta) \eta(\theta)\right\}.
\]
Moreover, if $x$ is a solution of \eqref{sys_general}, then $x_t \in \mathcal C_0$ for every $t \in \mathbb R_+$ in the interval of definition of $x$.

System~\eqref{sys_general} involves the integral of a continuous function $x$ with respect to a matrix-valued function $M$, which can be understood as in the usual Riemann--Stieltjes sense, but also as the integral with respect to a matrix-valued measure associated with $M$, in the usual Lebesgue--Stieltjes sense. Later, we will deal with functions that may not be integrable in the Riemann--Stieltjes sense, making it necessary to consider integrals in the Lebesgue--Stieltjes sense. The use of these two notions of integration is consistent, since if a function is integrable in the Riemann--Stieltjes sense, then it is also integrable in the Lebesgue--Stieltjes sense, and the two integrals coincide. As both such interpretations are useful for our analysis, we next detail these points of view, starting by providing the definitions of total variation of a matrix-valued function and of a matrix-valued Borel measure.

\begin{definition}
Given a matrix-valued function $M\colon [a,b] \to \mathcal M_d (\mathbb{R})$ and a matrix-valued Borel measure $\mu\colon \mathfrak{B}_{[a, b]} \to \mathcal M_d (\mathbb{R})$ (i.e., each entry of $\mu$ is a Borel signed measure), we define the \emph{total variation} of $M$ and $\mu$ in $[a, b]$ as
\[\V_a^b(M) = \sup\left\{\sum_{i=1}^k \abs{M(t_{i+1})-M(t_{i})}\suchthat k\in \mathbb N^\ast \text{ and }a=t_1<t_2<\dotsb<t_{k+1}=b \right\} \]
and
\[\V_a^b(\mu) = \sup\left\{\sum_{i=1}^k \abs{\mu(E_i)}\suchthat k\in \mathbb N, E_1, \dotsc,E_{k}\in \mathfrak{B}_{[a,b]} \text{ are pairwise disjoint and} \bigcup_{j=1}^{k}E_j =[a,b]\right\}, \]
respectively.
We say that $M$ (respectively, $\mu$) is of \emph{bounded variation} in $[a, b]$ if $\V_{a}^{b}(M) < +\infty$ (respectively, $\V_{a}^{b}(\mu) < +\infty$).
\end{definition}

Note that, even though $\V_a^b(M)$ and $\V_a^b(\mu)$ depend on the choice of the induced matrix norm $\abs{ \cdot }$, the notion of bounded variation is independent of the choice of the norm. In addition, $M\colon [a,b]\to \mathcal{M}_d(\mathbb{R})$ (respectively, $\mu\colon \mathfrak{B}_{[a, b]} \to \mathcal M_d (\mathbb{R})$) is of bounded variation if and only if all its entries are of bounded variation.

\begin{definition}
A matrix-valued function of bounded variation $M\colon [a,b]\to \mathcal M_d (\mathbb{R})$ is called \emph{normalized} if $M(a)=0$ and $M$ is continuous from the right on $(a,b)$. The space of normalized bounded variation functions is denoted by $\NBV([a,b],\mathcal M_d (\mathbb{R}))$.
\end{definition}

\begin{remark}\label{remk:M-mu}
There is a one-to-one correspondence (see, e.g., \cite[Appendix~I]{Diekmann1995Delay}) between normalized functions of bounded variation $M\colon [a,b]\to \mathcal M_d(\mathbb{R})$ and matrix-valued Borel measures $\mu_M$ of bounded variation, expressed by
\[M(t)=\mu_M([a,t]), \qquad t \in (a, b].\]
In addition, as shown in Proposition~\ref{proposition:tvcoincides} in Appendix~\ref{app:measure}, we have 
$\V_{a}^{b}(M) = \V_{a}^{b}(\mu_M)$. Note also that, from the results in \cite[Appendix~I]{Diekmann1995Delay}, thanks to the normalization $M(a)=0$ for $M\in \NBV([a,b],\mathcal M_d (\mathbb{R}))$, we have that $\V_{a}^{b}(\cdot)$ is a norm in the space $\NBV([a,b],\mathcal M_d (\mathbb{R}))$, which renders it a Banach space.

We also remark that, in the context of \eqref{sys_general}, the normalization of $M$ is not a restriction, since the integral in \eqref{sys_general} does not change if a constant is added to $M$ or if $M$ is modified in all its discontinuity points in $(-1, 0)$ in order to become right-continuous.
\end{remark}

In the sequel, given $M \in \NBV([a, b], \mathcal M_d(\mathbb R))$, we denote by $\mu_M$ the unique Borel measure associated with $M$ in the sense of Remark~\ref{remk:M-mu}, or simply by $\mu$ when $M$ is clear from the context. In addition, when the interval $[a, b]$ is clear from the context, we denote $\V_{a}^{b}(M)$ and $\V_a^b(\mu)$ simply by $\V(M)$ and $\V(\mu)$, respectively.

Using Remark~\ref{remk:M-mu}, one can rewrite \eqref{sys_general} as
\begin{equation*}
x(t) = \int_{-1}^0 \diff \mu_M(\theta) x(t + \theta), \qquad t \in \mathbb R_+,
\end{equation*}
where the integral in the right-hand side is a Lebesgue integral with respect to the matrix-valued Borel measure $\mu_M$. This corresponds to the Lebesgue--Stieltjes sense of the integral in the right-hand side of \eqref{sys_general}, and it makes sense even for functions $x$ that are not necessarily continuous.

\begin{remark}
\label{remk:decomposition}
Given a matrix-valued Borel measure of finite total variation $\mu \colon \mathfrak{B}_{[a,b]}\to \mathcal M_d(\mathbb{R})$, applying \cite[Theorem~19.61]{Hewitt1965Real} componentwise, we can decompose  $\mu$ as 
\[\mu=\mu_{\mathrm{ac}}+\mu_{\mathrm{sc}}+\mu_{\mathrm d},\]
where $\mu_{\mathrm{ac}}$ is absolute continuous with respect to the Lebesgue measure, $\mu_{\mathrm{sc}}$ is continuous and singular with respect to the Lebesgue measure and $\mu_{\mathrm d}$ is a discrete measure. In view of Remark~\ref{remk:M-mu}, the function $M$ associated with $\mu$ is decomposed as
\[M(t)=N(t)+S(t)+D(t),\qquad t\in [a,b],\]
where $N$ is absolutely continuous, $S$ is continuous and satisfies $S^\prime(t) = 0$ for a.e.\ $t \in [-1, 0]$, and $D$ is a countable sum of indicator functions multiplied by matrices. The functions $N$, $S$, and $D$ are called the \emph{absolutely continuous}, \emph{singular}, and \emph{purely discontinuous} parts of $M$, respectively.
\end{remark}

We can now turn to our result on existence and uniqueness of solutions of \eqref{sys_general}. The proof we present below is based on the existence proof from \cite[Chapter~2, Theorem~8.1]{Hale1993Introduction} for a more general neutral equation. Since our equation is linear, we can adapt that strategy to prove not only local existence but also uniqueness and global existence.

\begin{proposition}[{\cite[Proposition~1]{GoncalvesNetto2025Strong}}]\label{prop_wellposed}
Let $M\colon [-1, 0] \to \mathcal M_d(\mathbb R)$ be of bounded variation, define $A_M = M(0) - M(0^-)$, and suppose that $\det (\Id-A_M) \ne 0$. Then, for every $\phi\in \mathcal C_0$, there exists a unique global solution $x\colon [-1,+\infty)\to \mathbb{R}^d$ of \eqref{sys_general} with initial condition $x_0 = \phi$.
\end{proposition}

\begin{proof}
With no loss of generality, we assume that $A_M = 0$, for otherwise we can define $N\colon [-1, 0] \to \mathcal M_d(\mathbb R)$ by
\begin{equation}
\label{eq:no-atom-at-0}
N(\theta)=\begin{dcases*}
(\Id - A_M)^{-1} M(\theta), & if $-1\leq \theta <0$,\\
(\Id - A_M)^{-1} M(0^-), & if $\theta=0$,
\end{dcases*}
\end{equation}
and remark that $N$ is continuous at $0$ and \eqref{sys_general} is equivalent to
\[
x(t)=\int_{-1}^0 \diff N(\theta)x(t+\theta),
\]
in the sense that $x$ is a solution of the above system if and only if it is a solution of \eqref{sys_general}.

Given $\alpha \in (0, 1]$, we define the set $\mathcal A(\alpha)$ by
\[
\mathcal A(\alpha) = \left\{\eta \in \mathcal C([-1, \alpha], \mathbb{R}^d) \suchthat \eta_0 = \phi\right\}
\]
and we endow $\mathcal A(\alpha)$ with the metric induced by the supremum norm in $\mathcal C([-1, \alpha], \mathbb{R}^d)$, which we denote by $\norm{\cdot}_{\alpha}$ in this proof. We define the operator $S\colon \mathcal{A}(\alpha)\to\mathcal{A}(\alpha)$ by setting
\begin{equation*}
S(x)(t) = \begin{dcases*}
\phi(t), & if $-1 \leq t \leq 0$, \\
\int_{-1}^0 \diff M(\theta) x(t + \theta), & if $0 < t \leq \alpha$,
\end{dcases*}
\end{equation*}
and we remark that $S(x) \in \mathcal A(\alpha)$ for every $x \in \mathcal A(\alpha)$ since $\phi \in \mathcal C_0$. Notice also that $x$ is a local solution of \eqref{sys_general} in the interval $[-1, \alpha]$ with initial condition $\phi$ if and only if $S(x) = x$.

For $x, y \in \mathcal{A}(\alpha)$, we have $S(x)(t) - S(y)(t) = 0$ for $-1 \leq t \leq 0$ and, for $0 < t \leq \alpha$, we have from Proposition~\ref{proposition:tvcoincides} in Appendix~\ref{app:measure} that
\begin{align*}
\abs{S(x)(t) - S(y)(t)} & = \abs*{\int_{-1}^0 \diff M(\theta) \bigl[x(t + \theta) - y(t + \theta)\bigr] } = \abs*{ \int_{-t}^0 \diff M(\theta) \bigl[x(t + \theta) - y(t + \theta)\bigr]} \\
& \leq \V_{-t}^{0}(M) \norm{x - y}_{\alpha} \leq \V_{-\alpha}^{0}(M) \norm{x - y}_{\alpha}.
\end{align*}
Therefore, for every $x, y$ in $\mathcal{A}(\alpha)$, one has
\[\norm{S(x)-S(y)}_{\alpha} \leq \V_{-\alpha}^{0}(M) \norm{x - y}_{\alpha}. \]
Since $\lim_{s\to 0^-} \V_{s}^{0}(M) = \abs{A_M} = 0$, one deduces that, for some $\alpha > 0$ small enough that depends only on $M$, $S$ is a contraction. It therefore has a unique fixed point $x \in \mathcal{A}(\alpha)$, which is the unique local solution of \eqref{sys_general} in $[-1, \alpha]$. As $\alpha > 0$ depends only on $M$ and not on $\phi$, one can iterate the above argument to extend the solution $x$ to $[-1, +\infty)$, obtaining a unique global solution.
\end{proof}

Motivated by Proposition~\ref{prop_wellposed}, we introduce the class
\[
\WellPosed = \{M \in \NBV([-1, 0], \mathcal M_d(\mathbb R)) \suchthat \det (\Id-A_M) \ne 0 \text{ where } A_M = M(0) - M(0^-)\}.
\]
In the sequel, with a slight abuse of notation, we will also say that a measure $\mu$ belongs to $\WellPosed$ if the function of normalized bounded variation $M$ associated with $\mu$ in the sense of Remark~\ref{remk:M-mu} belongs to $\WellPosed$. We shall only work from now on with System~\eqref{sys_general} under the assumption that $M \in \WellPosed$ in order to ensure existence and uniqueness of solutions.

\section{Exponential behavior of solutions}
\label{sec:definition-spaces}

Given $M \in \WellPosed$, one can establish, using the classical exponential bound of a strongly continuous semigroup (see, e.g., \cite[Section~12.3]{Hale1977Theory}), that solutions of \eqref{sys_general} satisfy an exponential bound, i.e., there exist positive constants $K,\alpha>0$ such that, for every $\phi\in \mathcal{C}_0$,
\begin{equation}
\label{eq_exp-estimate}
\norm{x_t}_\infty \leq K e^{\alpha t}\norm{\phi}_\infty,\qquad  t\in \mathbb R_+,
\end{equation}
where $x$ is the unique solution of \eqref{sys_general} with initial condition $\phi$. For sake of completeness, we provide a proof of \eqref{eq_exp-estimate} in Appendix~\ref{app:other}, relying on direct estimates from \eqref{sys_general}.

In order to study the stability properties of System~\eqref{sys_general}, we recall the definitions of exponential stability, exponential rate, characteristic function, and spectral abscissa.

\begin{definition}
\label{def_big-def}
Let $M \in \NBV([-1, 0], \mathcal M_d(\mathbb R))$.
\begin{enumerate}
\item For $M \in \WellPosed$, we define the \emph{exponential rate} $\omega_M$ of $M$ by
\begin{align*}
\omega_M = \inf \{&\alpha \in \mathbb{R}\suchthat \exists K>0 \text{ s.t.\ }\norm{x_t}_\infty\leq K e^{\alpha t}\norm{\phi}_\infty \text{ for all } t \in \mathbb R_+ \text{ and }\phi \in \mathcal{C}_0,\\
&\text{where $x$ is the unique solution of \eqref{sys_general} with } x_0=\phi\}.
\end{align*}
We say that $M$ (or $\Sigma(M)$) \emph{exponentially stable} if $\omega_M < 0$.

\item We define the \emph{characteristic function} of $M$ as the function $\Delta_M \colon \mathbb C \to \mathbb C$ given, for $s \in \mathbb C$, by
\[\Delta_M(s)=\det \left(\Id - \int_{-1}^0 e^{s\theta} \diff M(\theta)\right).\]

\item We define the \emph{spectral abscissa} $S_M$ of $M$ by
\[S_M=\sup \{\Real s \suchthat \Delta_M(s) = 0\}.\]

\item\label{item:def-sdg} For $M \in \WellPosed$, we say that $M$ satisfies the \emph{spectrum-determined growth condition} if we have the equality
\[\omega_M = S_M.\]
\end{enumerate}
\end{definition}

A series of remarks are now in order.

\begin{remark}
\label{remk:big-remk}
\begin{enumerate}
\item Although we have defined exponential rate, characteristic function, spectral abscissa, and the spectrum-determined growth condition as properties relating to some $M \in \NBV([-1, 0], \mathcal M_d(\mathbb R))$, we shall often talk about these properties also for a matrix-valued measure $\mu$ on $[-1, 0]$, and in this case they should be understood as the corresponding property for the function $M$ associated with $\mu$ in the sense of Remark~\ref{remk:M-mu}.

\item As a consequence of \eqref{eq_exp-estimate}, one deduces that the set over which the infimum is computed in the definition of $\omega_M$ is nonempty, and in particular we have $\omega_M < +\infty$.

\item As the notion of exponential rate is related to solutions of \eqref{sys_general}, we require in its definition that $M \in \WellPosed$ in order to ensure the well-posedness of \eqref{sys_general} from Proposition~\ref{prop_wellposed}. However, the assumption that $M \in \WellPosed$ is not necessary in order to define the characteristic function and the spectral abscissa.

\item\label{item:big-remk-ineq-SM-omegaM} The inequality $S_M \leq \omega_M$ is satisfied for every $M \in \WellPosed$ since, for every $s \in \mathbb C$ such that $\Delta_M(s) = 0$, denoting by $v \in \mathbb C^d$ a nonzero element in the kernel of $\Id - \int_{-1}^0 e^{s\theta} \diff M(\theta)$, we obtain that $t \mapsto e^{s t} v$ is a ($\mathbb C^d$-valued) solution of \eqref{sys_general}, so its real and imaginary parts are $\mathbb R^d$-valued solutions of \eqref{sys_general}, and the definition of exponential rate yields that $\Real s \leq \omega_M$.

\item\label{item_remk-sdg-condition-H} From \cite[Chapter~9, Theorem~4.4]{Hale1993Introduction}, it follows that a sufficient condition for a function $M \in \NBV([-1, 0], \mathcal M_d(\mathbb R))$ to satisfy the spectrum-determined growth condition is that the entries $M_{jk}$ of $M$ have an atom before they become constant, i.e., for all $(j, k) \in \llbracket 1, d\rrbracket^2$ there exists $t_{jk} \in [-1, 0)$ with $M_{jk}(t) = M_{jk}(t_{jk}^+)$ for $t \geq t_{jk}$ and $M_{jk}(t_{jk}^-) \neq M_{jk}(t_{jk}^+)$. Note that, in this case, we necessarily have $M(0) = M(0^-)$, and in particular $M \in \WellPosed$.

\item\label{item_N-subset-G} Another sufficient condition for $M \in \NBV([-1, 0], \mathcal M_d(\mathbb R))$ to satisfy the spectrum-de\-ter\-mined growth condition, provided in \cite[Chapter~12, Theorem~3.4]{Hale1977Theory}, is when $M$ is given by
\begin{equation}\label{eq_NS}
M(\theta) = \sum_{k=1}^{+\infty} A_k \1_{[-\tau_k,0]}(\theta) + N(\theta), \qquad \theta \in [-1, 0],
\end{equation}
where $\tau_k \in (0, 1)$ and $A_k\in \mathcal M_d (\mathbb{R})$ for $k \in \mathbb N^\ast$ and $N$ is absolutely continuous. In this case, we also have $M \in \WellPosed$, and the singular part of $M$ (in the sense of the decomposition of $M$ from Remark~\ref{remk:decomposition}) vanishes.

The discontinuity points of $M$ given by \eqref{eq_NS} are the values $-\tau_k$, $k \in \mathbb N$, which all belong to the open interval $(-1, 0)$. One can also add a discontinuity point at $-1$ in $M$ and still obtain the exact same results that we will prove in the sequel in this more general framework (namely Remark~\ref{remk:explicit-variation} and Propositions~\ref{prop_tvnosingular} and \ref{prop_positiveinf}). However, note that, since $M$ is defined only over $[-1, 0]$, such a discontinuity can only be achieved by adding a term of the form $A_0 \1_{(-1, 0]}$, i.e., with the indicator function of the interval $(-1, 0]$ open at $-1$.

Notice also that, if $M \in \WellPosed$ and it has no singular part in the sense of the decomposition from Remark~\ref{remk:decomposition}, then, using the transformation from \eqref{eq:no-atom-at-0}, System~\eqref{sys_general} can be rewritten in an equivalent way as a system for which $M$ is given by \eqref{eq_NS} (with possibly an additional term of the form $A_0 \1_{(-1, 0]}$).
\end{enumerate}
\end{remark}

We now define the three subsets of $\NBV([-1,0],\mathcal M_d (\mathbb{R}))$ given by
\begin{equation}
\label{eq_def-N-L-G}
\begin{aligned}
\NoSingular & = \{M \in \NBV([-1,0],\mathcal M_d (\mathbb{R})) \suchthat M \text{ is given by } \eqref{eq_NS}\}, \\
\LowerBound & = \left\{M\in\NBV([-1,0], \mathcal M_d (\mathbb{R})) \suchthat S_M < +\infty \text{ and } \inf_{\Real s \ge \alpha} \abs{\Delta_M(s)} > 0 \text{ for every } \alpha>S_M \right\}, \\
\SDG & = \{M \in \WellPosed \suchthat S_M = \omega_M \}.
\end{aligned}
\end{equation}
The notation of these sets is justified by the facts that all elements of $\NoSingular$ have \emph{no singular} part in the sense of Remark~\ref{remk:decomposition}, the set $\LowerBound$ is characterized by a \emph{lower bound} on its characteristic function on right half-planes, and $\SDG$ is the set of those $M$ satisfying the \emph{spectrum-de\-ter\-mined growth} condition. By Remark~\ref{remk:big-remk}.\ref{item_N-subset-G}, we have $\NoSingular \subset \SDG \subset \WellPosed$. In addition, we have $\NoSingular \subset \LowerBound \subset \WellPosed$. Indeed, Proposition~\ref{prop_positiveinf} in Appendix~\ref{app:other} shows the left inclusion, while, for the other one, it suffices to note that, as $s \to +\infty$ with $s \in \mathbb R$, we have $\Delta_M(s) \to \det(\Id - A_M)$, with $A_M = M(0) - M(0^-)$.

With a slight abuse of notation, we will sometimes say that a measure $\mu$ belongs to one of the classes $\NoSingular$, $\LowerBound$, or $\SDG$ if the function of normalized bounded variation $M$ associated with $\mu$ in the sense of Remark~\ref{remk:M-mu} belongs to the same class.

\section{Continuity of the spectral abscissa in total variation norm}
\label{sec_tv-norm}

In this section, we will prove that the spectral abscissa of $M$ is continuous in the total variation norm when $M$ belongs to $\LowerBound$. As a consequence, we will obtain that exponential stability of \eqref{sys_general} for some $M$ with no singular part implies exponential stability of \eqref{sys_general} for small enough perturbations of $M$ in total variation norm, still with no singular part.

We start with the following technical result relating the characteristic functions of two elements $M$ and $N$ from $\NBV([-1, 0], \mathcal M_d(\mathbb R))$ on a right half-plane.

\begin{lemma}\label{lemma_tvineq}
Given $\alpha\in \mathbb{R}$ and $R>0$, there exists $C = C(\alpha, R) > 0$ such that
\begin{equation}
\sup_{\Real s \geq \alpha} \abs{\Delta_M(s)-\Delta_N(s)} \leq C \V(M-N),
\end{equation}
for every $M,N\in \NBV([-1,0],\mathcal M_d (\mathbb{R}))$ such that $\V(M) \leq R$ and $\V(N) \leq R$.
\end{lemma}

\begin{proof}
Let $\alpha_-=\max\{0,-\alpha\}$ and $r = 1 + e^{\alpha_-} R$. We set $K_r=\{A\in \mathcal M_d(\mathbb{R})\suchthat \abs{A}\leq r\}$. Since $\det(\cdot)$ is continuously differentiable, by the mean value theorem there exists $\widetilde{C}=\widetilde{C}(r)>0$ such that 
\begin{equation}\label{ineq_detlipschitz}
\abs{\det (A)-\det (B)} \leq \widetilde{C} \abs{A - B}, \qquad \text{ for all } A, B\in K_r. 
\end{equation}
Note now that $\Id - \int_{-1}^0 e^{s\theta} \diff M(\theta)$ and $\Id - \int_{-1}^0 e^{s\theta} \diff N(\theta)$ belong to $K_r$ for all $s \in \mathbb C$ with $\Real s \ge \alpha$, since, using Proposition~\ref{proposition:tvcoincides} from Appendix~\ref{app:measure}, we have
\[\abs*{\Id - \int_{-1}^0 e^{s\theta} \diff M(\theta)}\leq \abs{\Id} + \abs*{ \int_{-1}^0 e^{s\theta} \diff M(\theta)}\leq 1 + \max_{\theta \in [-1,0]}\abs{e^{s\theta}} \V(M) \leq 1 + e^{\alpha_-} \V(M) \leq r,\]
and similarly for $N$. Using \eqref{ineq_detlipschitz} we have that, for all $s \in \mathbb C$ with $\Real s \ge \alpha$,
\[
\abs{\Delta_M(s)-\Delta_N(s)}\leq \widetilde{C} \abs*{ \int_{-1}^0 e^{s\theta} \diff(N-M)(\theta) } \leq \widetilde{C} e^{\alpha_-} \V(M-N),
\]
and the conclusion follows with $C = \widetilde{C} e^{\alpha_-}$.
\end{proof}

We are now ready to prove the main result of this section.

\begin{theorem}
\label{thm:spectral-abscissa-continuous}
The spectral abscissa map $\LowerBound \ni M \mapsto S_M\in \mathbb R \cup \{- \infty\}$ is continuous on $\LowerBound $ with respect to the total variation norm. 
\end{theorem}

\begin{proof}
We only treat here the case $S_M \in \mathbb R$, since the case $S_M = -\infty$ can be treated exactly as in the first part of the proof below, up to replacing $S_M + \varepsilon$ by $-K$ for an arbitrary $K > 0$.

Let $M\in \LowerBound$ and $(N_n)_{n\in \mathbb N}$ be a sequence in $\LowerBound$ such that $\lim_{n\to +\infty} \V(M-N_n)=0$. Given $\varepsilon > 0$, let
\[\delta = \inf_{\Real s \geq S_M+\varepsilon}\abs{\Delta_M(s)}\]
and note that $\delta > 0$ by definition of $\LowerBound$.
Since $\lim_{n\to +\infty} \V(M-N_n)=0$, we may assume without loss of generality that $\V(N_n) \leq 1 + \V(M)$ for every $n\in \mathbb N$. Let $\alpha = S_M + \varepsilon$ and $R=1 + \V(M)$. By Lemma~\ref{lemma_tvineq}, there exists $C(\alpha,R) > 0$ such that, for all $s \in \mathbb C$ with $\Real s \ge \alpha$,
\[\delta-\abs{\Delta_{N_n}(s)}\leq \abs{\Delta_M(s)}-\abs{\Delta_{N_n}(s)}\leq \abs{\Delta_M(s)-\Delta_{N_n}(s)} \leq C(\alpha,R) \V(M-N_n), \qquad \text{ for all } n\in \mathbb N.\]
Hence,
\[\inf_{\Real s \ge S_M + \varepsilon} \abs{\Delta_{N_n}(s)} \ge \delta - C(\alpha, R) \V(M-N_n), \qquad \text{ for all } n\in \mathbb N.\]
Since $\V(M-N_n)\to 0$ as $n\to +\infty$, there exists $n_1$ such that, for all $n\ge n_1$,
\[\inf_{\Real s \ge S_M + \varepsilon} \abs{\Delta_{N_n}(s)} > 0,\]
and consequently $S_{N_n} \leq S_M+ \varepsilon$ for every $n\ge n_1$.

On the other hand, by the same argument as above, for $\beta=S_M-\varepsilon$ and $R= 1+\V(M)$ there exists $C(\beta,R)>0$ such that
\[\sup_{\Real s \ge S_M-\varepsilon } \abs{\Delta_M(s) - \Delta_{N_n}(s)} \leq C(\beta,R)\V(M-N_n).\]
That is, $\Delta_{N_n}(s)\to \Delta_M (s)$ as $n\to +\infty$, uniformly on $\Real s \ge S_M-\varepsilon$. By the definition of $S_M$, there exists $s^*\in \mathbb{C}$ with $\Real s^* > S_M-\varepsilon$ and $\Delta_M(s^*)=0$. Take $\eta>0$ small enough so that the ball $B_{\mathbb C}(s^\ast, \eta)$ is contained in the half-plane $\{s \in \mathbb C \suchthat \Real s \ge S_M-\varepsilon\}$. By Hurwitz's theorem, there exists $n_2$ such that, for all $n\geq n_2$, $\Delta_{N_n}(\cdot)$ has at least one zero inside $B(s^\ast, \eta)$, and hence $S_{N_n} \geq S_M-\varepsilon$ for all $n\geq n_2$. Therefore, for every $n\geq \max \{n_1,n_2\}$ we have, 
\[S_M-\varepsilon \leq S_{N_n} \leq S_M+\varepsilon,\]
which implies that $\lim_{n\to +\infty}S_{N_n}=S_M$, completing the proof.
\end{proof}

We next deduce a consequence of the previous theorem for matrix-valued functions in $\LowerBound \cap \SDG$. 

\begin{corollary}\label{coro_tvstab}
Assume that $M\in \LowerBound \cap \SDG$ is exponentially stable. Then there exists $\varepsilon>0$ such that $\widetilde M$ is exponentially stable for every $\widetilde M \in \LowerBound \cap \SDG$ with $\V(M - \widetilde M) < \varepsilon$.
\end{corollary}

\begin{proof}
Since $M$ is exponentially stable, we have $\omega_M < 0$. From the fact that $M \in \SDG$, it follows that $S_M = \omega_M < 0$. Hence, by Theorem~\ref{thm:spectral-abscissa-continuous}, there exists $\varepsilon > 0$ such that, if $\widetilde M \in \LowerBound$ and $\V(\widetilde{M} - M) < \varepsilon$, then $S_{\widetilde{M}} < 0$. Consequently, $\widetilde{M}$ is exponentially stable for every $\widetilde{M} \in \LowerBound \cap \SDG$ satisfying $\V(M - \widetilde M) < \varepsilon$.
\end{proof}

\begin{remark}
\label{remk:explicit-variation}
Recalling that $\NoSingular \subset \LowerBound \cap \SDG$, we deduce that the previous corollary applies in particular to matrix-valued functions $M$ in $\NoSingular$, i.e., without singular part. In this case, we can explicitly compute the quantity $\V(M - \widetilde M)$ from the statement of the corollary. Indeed, take $M, \widetilde M \in \NoSingular$ and write
\[ M=\sum_{k=1}^{+\infty} A_k \1_{[-\tau_k,0]}+N, \qquad \widetilde M= \sum_{k=1}^{+\infty}  B_k \1_{[-\widetilde \tau_k,0]}+\widetilde N,\]
where $A_k, \tau_k, N, B_k, \widetilde \tau_k, \widetilde N$ are as in \eqref{eq_NS}, assuming further, with no loss of generality, that $\tau_k \neq \tau_j$ and $\widetilde \tau_k \neq \widetilde \tau_j$ if $k \neq j$. We rewrite $M$ and $\widetilde M$ as 
\[M=\sum_{k=1}^{+\infty} \widehat A_k \1_{[-h_k,0]}+N, \qquad \widetilde M=\sum_{k=1}^{+\infty} \widehat B_k \1_{[-h_k,0]}+\widetilde N,\]
where $\{h_k\suchthat k \in \mathbb N^\ast\}=\{\tau_k\suchthat k \in \mathbb N^\ast\} \cup \{\widetilde\tau_k\suchthat k \in \mathbb N^\ast\}$ and
\[\widehat A_k = \begin{dcases*}
A_j, & if $h_k = \tau_j$ for some $j\ge 1$, \\
0, & otherwise,
\end{dcases*} \qquad \widehat B_k = \begin{dcases*} B_j, & if $h_k = \widetilde\tau_j$ for some $j\ge 1$, \\
0, & otherwise.
\end{dcases*}\]
So, from Proposition~\ref{prop_tvnosingular} in Appendix~\ref{app:measure}, we have
\[\V(M - \widetilde M) =  \sum_{k=1}^{+\infty} \abs{\widehat A_k - \widehat B_k } + \int_{-1}^0 \abs{N^\prime(\theta)-\widetilde N^\prime(\theta)} \diff \theta. \]
\end{remark}

Despite the interest of Corollary~\ref{coro_tvstab}, it cannot be applied to the case where $M$ is purely discontinuous and one perturbs the discontinuity points of $M$, as the total variation norm may increase discontinuously under such perturbations. This is illustrated by the following example.

\begin{example}
\label{expl_var-delay}
Let $0<\tau_1<\tau_2<1$, $A$ be a nonzero element of $\mathcal M_d (\mathbb{R})$, and $M,N\colon [-1,0] \to \mathcal M_d (\mathbb{R})$ be given by $M=A \1_{[-\tau_1,0]}$ and $N=A \1_{[-\tau_2,0]}$. Then
\[\V(M-N)=2\abs{A}.\]
In this way, no matter how close $\tau_1$ is to $\tau_2$, $N$ is at a constant distance of $2 \abs{A}$ from $M$ in the sense of total variation, and in particular Corollary~\ref{coro_tvstab} does not guarantee the exponential stability of $N$ given the exponential stability of $M$.
\end{example}

The above example shows that the total variation norm is not adapted to understand if exponential stability is preserved or not with respect to perturbations in delays. This fact motivates the next section.

\section{Strong stability}
\label{sec_strong}

Our next aim is to generalize the Hale--Silkowski criterion from Proposition~\ref{prop:Hale-Silkowski} for system of the form \eqref{eq:intro-diff-eqn-sum} to more general systems of the form \eqref{sys_general}. The Hale--Silkowski criterion, and the corresponding notions of (local or global) strong stability from Proposition~\ref{prop:Hale-Silkowski}\ref{item:HS-local} and \ref{item:HS-global}, can be seen as a result on the robustness of exponential stability with respect to the delays of \eqref{eq:intro-diff-eqn-sum}. One should then start by providing a suitable notion of ``perturbation of the delays'' for systems under the more general form \eqref{sys_general}.

Note that, in \eqref{sys_general}, the delay comes from the fact that the unknown function $x$ is evaluated at $t + \theta$ in the right-hand side of the equation, for $\theta \in [-1, 0]$. The variable $\theta$ can be seen as a ``distributed delay'' in \eqref{sys_general}, and hence a natural candidate for being a ``perturbation of the delays'' of System~\eqref{sys_general} consists in replacing $\theta$ by $\varphi(\theta)$, where $\varphi\colon [-1, 0] \to [-1, 0]$ is a perturbation of the identity map. System~\eqref{sys_general} is then replaced by the perturbed system given by
\begin{equation}
\label{sys_perturbed}
x(t) = \int_{-1}^0 \diff M(\theta) x(t + \varphi(\theta)), \qquad t \in \mathbb R_+.
\end{equation}
Note that, in the case of System~\eqref{eq:intro-diff-eqn-sum} where $M$ is given by \eqref{eq_M-finite-sum} with delays $\tau_1, \dotsc, \tau_N$ in $(0, 1)$ and $A_1, \dotsc, A_N$ in $\mathcal M_d(\mathbb R)$, System~\eqref{sys_perturbed} corresponds exactly to perturbing the delays $\tau_1, \dotsc, \tau_N$ to $-\varphi(-\tau_1), \dotsc, -\varphi(-\tau_N)$. As $\varphi$ is not necessarily continuous, the integral in the right-hand side of \eqref{sys_perturbed} may not exist in the Riemann--Stieltjes sense, but it is well-defined in the Lebesgue--Stieltjes sense under the assumption that $\varphi$ is Borel-measurable.

More precisely, we define the set 
\[\mathbf{B}=\{\varphi\colon[-1,0]\to [-1,0]\suchthat \varphi \text{ is Borel-measurable}\}.\]
Using the correspondence between matrix-valued functions $M$ of normalized bounded variation and matrix-valued Borel measures with finite total variation from Remark~\ref{remk:M-mu}, one can rewrite the perturbed system \eqref{sys_perturbed} in terms of a pushforward of the measure $\mu_M$ corresponding to $M$ as
\begin{equation}
\label{sys_perturbed-measure}
\Sigma(\varphi_{\#}\mu_M)\colon \qquad x(t)=\int_{-1}^0 \diff(\varphi_{\#}\mu_M) (\theta)x(t+\theta), \qquad t \in \mathbb R_+,
\end{equation}
where $\varphi_{\#}\mu_M$ is the pushforward of $\mu_M$ by $\varphi\in \textbf{B}$, i.e., the measure defined by $\varphi_{\#}\mu_M(A)=\mu_M(\varphi^{-1}(A))$ for $A \in \mathfrak{B}_{[-1,0]}$. Indeed, for every $x \in \mathcal C([-1, +\infty), \mathbb R^d)$ and $t\in \mathbb R_+$, $x(t+\varphi(\cdot))$ is Borel-measurable and bounded, hence it is integrable with respect to $\mu_M$, and we have
\begin{equation*}
\int_{-1}^0 \diff(\varphi_{\#}\mu_M) (\theta)x(t+\theta) = \int_{-1}^0 \diff \mu_M(\theta)x(t+\varphi(\theta)) = \int_{-1}^0 \diff M(\theta)x(t+\varphi(\theta)),\qquad t\ge 0.
\end{equation*}
Note that System~\eqref{sys_perturbed-measure} is of the same form of System~\eqref{sys_general}, with $M$ replaced by the pushforward measure $\varphi_{\#} \mu_M$.

We are now in position to provide our generalization of the definition of strong stability to System~\eqref{sys_general}.

\begin{definition}[$\mathsf{K}$-strong stability]
\label{def_strong-stab}
Let $\mathsf{K} \subset \WellPosed$ be nonempty, $M \in \mathsf K$, and $\mu$ be the matrix-valued Borel measure associated with $M$ in the sense of Remark~\ref{remk:M-mu}.
\begin{enumerate}[ref={\arabic*)}]
\item We say that $\Sigma(M)$ is \emph{$\mathsf{K}$-locally strongly stable} if there exists $\varepsilon > 0$ such that $\Sigma(\varphi_{\#} \mu)$ is exponentially stable for all $\varphi \in \mathbf B$ with $\norm{ \varphi - \id }_\infty < \varepsilon$  and $\varphi_{\#} \mu \in \mathsf{K}$.
\item We say that $\Sigma(M)$ is \emph{$\mathsf{K}$-globally strongly stable} (or simply \emph{$\mathsf{K}$-strongly stable}) if $\Sigma(\varphi_{\#} \mu)$ is exponentially stable for all $\varphi\in\mathbf{B}$ such that $\varphi_{\#} \mu \in \mathsf{K}$.
\end{enumerate}
\end{definition}

\begin{remark}
In Definition~\ref{def_strong-stab}, we fix a class $\mathsf K$ of functions included in $\WellPosed$ and we restrict our attention to only those perturbations $\varphi \in \mathbf B$ that preserve $\mathsf K$, in the sense that $\varphi_{\#} \mu \in \mathsf{K}$ (where, once again, we do the slight abuse of notation of saying that a measure belongs to a subset of $\NBV([-1, 0], \mathcal M_d(\mathbb R))$ if the corresponding function of normalized bounded variation in the sense of Remark~\ref{remk:M-mu} belongs to the said subset). The main reason for this is that, although one would ideally like to take $\mathsf K = \WellPosed$, we are unable to provide a full generalization of the Hale--Silkowski criterion in this setting, except in the scalar case, and, since we make use of spectral methods, we are obliged to consider the subclass $\mathsf K = \SDG$.
\end{remark}

\begin{remark}
\label{remk:pushfwd-in-WP}
When $\mathsf K = \WellPosed$, one can easily verify that the condition $\varphi_{\#}\mu \in \WellPosed$ is satisfied, in particular, if $M \in \WellPosed$ and $\varphi(\theta) \neq 0$ for all $\theta \in [-1, 0)$. Moreover, if $M \in \NBV([-1, 0], \mathcal M_d(\mathbb R))$ and $\varphi\colon [-1, 0] \to (-1, 0)$ is Borel-measurable and its image is at most countable, one can easily check that $\varphi_{\#} \mu_M \in \NoSingular$.
\end{remark}

Now that we have provided a suitable generalization of the definition of strong stability to System~\eqref{sys_general}, in order to generalize the Hale--Silkowski criterion from Proposition~\ref{prop:Hale-Silkowski}, we also need to provide a suitable generalization of the quantity $\rho_{\mathrm{HS}}$ from \eqref{eq:rho-HS}. At the light of the form of \eqref{eq:rho-HS}, we propose, for a given $M \in \NBV([-1, 0], \mathcal M_d(\mathbb R))$, to use the quantity
\begin{equation}
\label{GHSC}
\rho_{\mathrm{HS}}(M) = \sup_{\xi\in \mathfrak C}  \spr \left(\int_{-1}^0 e^{i\xi(\theta)} \diff M(\theta)\right),
\end{equation}
where $\mathfrak C=\{\xi\colon [-1,0]\to \mathbb{R} \suchthat \xi\text{ is Borel-measurable}\}$. Note that the integral in the above expression reduces to the sum in \eqref{eq:rho-HS} when $M$ is given by \eqref{eq_M-finite-sum}, with $\theta_k = \xi(-\tau_k)$.

\begin{remark}
\label{remk:rho-HS-TV}
By Proposition~\ref{proposition:tvcoincides} in Appendix~\ref{app:measure}, we immediately have from \eqref{GHSC} that $\rho_{\mathrm{HS}}(M) \leq \V(M)$. In addition, one has equality in dimension $d = 1$ since, in that case, denoting by $P$ and $N$ the positive and negative sets of $\mu_M$ in the sense of the Hahn decomposition with $P \cup N = [-1, 0]$ and $P \cap N = \emptyset$ (see, e.g., \cite[Theorem~3.3]{Folland1999Real}), we have the equality $\V(M) = \int_{-1}^0 e^{i\xi(\theta)} \diff M(\theta)$ for the function $\xi\colon [-1,0]\to \mathbb{R}$ given by $\xi(\theta) = 0$ if $\xi \in P$ and $\xi(\theta) = \pi$ if $\xi \in N$.
\end{remark}

We then conjecture the following general result.
\begin{conjecture}
\label{GHS-Conjecture}
Let $M \in \WellPosed$. Then the following statements are equivalent.
\begin{enumerate}
\item We have $\rho_{\mathrm{HS}}(M) <1$.
\item System~$\Sigma(M)$ is $\WellPosed$-strongly stable.
\item System~$\Sigma(M)$ is $\WellPosed$-locally strongly stable.
\end{enumerate}
\end{conjecture}

Conjecture~\ref{GHS-Conjecture} was shown to be true in dimension $d = 1$ in \cite{GoncalvesNetto2025Strong}, and we present this result in Section~\ref{sec_main-result} below. The proof relies on the fact that, by Remark~\ref{remk:rho-HS-TV}, we have $\rho_{\mathrm{HS}}(M) = \V(M)$ in the scalar setting. Whether this conjecture holds true in higher dimension remains an open problem.

Note that, by comparing the expression of $\rho_{\mathrm{HS}}(M)$ with the characteristic function $\Delta_M$, it is tempting to try to study Conjecture~\ref{GHS-Conjecture} through spectral methods, i.e., the analysis of $\Delta_M$ and its roots. However, in order to relate properties of $\Delta_M$ with the asymptotic behavior of solutions of $\Sigma(M)$, one must require $M$ to satisfy the spectrum-determined growth condition. This leads to the idea that one might use spectral methods to prove the following modification of Conjecture~\ref{GHS-Conjecture}.

\begin{conjecture}
\label{GHS-Conjecture-SDG}
Let $M \in \SDG$. Then the following statements are equivalent.
\begin{enumerate}
\item We have $\rho_{\mathrm{HS}}(M) <1$.
\item System~$\Sigma(M)$ is $\SDG$-strongly stable.
\item System~$\Sigma(M)$ is $\SDG$-locally strongly stable.
\end{enumerate}
\end{conjecture}

We are (still) not able to prove Conjecture~\ref{GHS-Conjecture-SDG} in full generality. However, it can be proved under a modification of the definition of strong stability in order to include a uniform exponential rate. This modification, together with other partial results towards the proof of Conjecture~\ref{GHS-Conjecture-SDG}, are the subject of Section~\ref{sec_main-strong-stab} below.

\subsection{Strong stability criterion in the scalar case}
\label{sec_main-result}

To prove Conjecture~\ref{GHS-Conjecture} in dimension $d = 1$, we first provide the following partial result, which holds true independently of the dimension $d$.

\begin{proposition}
\label{prop:TV-less-1-implies-global-strong-stab}
Let $M \in \WellPosed$ and assume that $\V(M) < 1$. Then System~$\Sigma(M)$ is $\WellPosed$-strongly stable.
\end{proposition}

\begin{proof}
Let $\varphi \in \mathbf B$ be such that $\varphi_{\#} \mu \in \WellPosed$, where $\mu$ is the measure associated with $M$ in the sense of Remark~\ref{remk:M-mu}. Let $x \in \mathcal C([-1, +\infty), \mathbb R^d)$ be a solution of $\Sigma(\varphi_{\#} \mu)$. Then, by Proposition~\ref{proposition:tvcoincides} in Appendix~\ref{app:measure}, we have
\begin{equation*}
\abs{x(t)} = \abs*{\int_{-1}^0 \diff\mu(\theta) x(t + \varphi(\theta))} \leq \V(M) \norm{x_t}_\infty.
\end{equation*}
Hence, from \cite[Lemma~1]{Mazenc2015Trajectory}, $\Sigma(\varphi_{\#} \mu)$ is exponentially stable, yielding that $\Sigma(M)$ is $\WellPosed$-strongly stable.
\end{proof}

Let us now state and prove our main result in dimension $d = 1$.

\begin{theorem}\label{HS:scalar}
Let $M \in \WellPosed$ and assume that $d = 1$. Then the following statements are equivalent.
\begin{enumerate}
\item\label{item_main-HS} We have $\rho_{\mathrm{HS}}(M) <1$.
\item\label{item_main-glob} System~$\Sigma(M)$ is $\WellPosed$-strongly stable.
\item\label{item_main-loc} System~$\Sigma(M)$ is $\WellPosed$-locally strongly stable.
\end{enumerate}
\end{theorem}

\begin{proof}
The fact that \ref{item_main-HS} implies \ref{item_main-glob} follows from Proposition~\ref{prop:TV-less-1-implies-global-strong-stab}, and it is clear from the definition of strong stability that \ref{item_main-glob} implies \ref{item_main-loc}.

We are thus left to show that \ref{item_main-loc} implies \ref{item_main-HS}, and we proceed by contraposition. Recall that, by Remark~\ref{remk:rho-HS-TV}, we have $\rho_{\mathrm{HS}}(M) = \V(M)$. Assume that $\V(M) \geq 1$ and denote by $\mu$ the measure associated with $M$ in the sense of Remark~\ref{remk:M-mu}. Recall that, by the Hahn and Jordan decomposition theorems (see, e.g., \cite[Theorems~3.3 and 3.4]{Folland1999Real}), there exist Borel-measurable subsets $P$ and $N$ of $[-1, 0]$ and unique positive measures $\mu^+$ and $\mu^-$ such that $\mu = \mu^+ - \mu^-$, $\mu^-(P) = 0$, $\mu^+(N) = 0$, $P \cup N = [-1, 0]$, and $P \cap N = \emptyset$. In addition, $\V(\mu) = \mu^+([-1, 0]) + \mu^-([-1, 0]) = \mu^+(P) + \mu^-(N)$.

We need to show that, for every $\varepsilon>0$, there exists $\varphi \in \mathbf B$ such that $\norm{\varphi - \id}_{\infty}<\varepsilon$, $\varphi_{\#}\mu \in \WellPosed$, and $\Sigma(\varphi_{\#} \mu)$ is not exponentially stable. So, given $\varepsilon>0$, we take $m \in \mathbb N$ and a partition $-1 = t_0 < t_1 < \dotsb < t_{m+1} = 0$ of $[-1,0]$ into $m+1$ intervals such that
\[\max_{k \in \llbracket 0, m\rrbracket} \abs{t_{k+1}-t_k} < \varepsilon.\]
Let us define the intervals $I_k=[t_k,t_{k+1})$ for $k\in \llbracket 0, m-1\rrbracket$ and $I_m=[t_m,t_{m+1}]$, and the sets
\[P_k=P\cap I_k\quad\text{ and }\quad N_k=N\cap I_k, \qquad k \in \llbracket 0, m\rrbracket.\]
We define the function $\varphi\colon [-1,0]\to [-1,0]$ by setting
\[
\varphi(\theta)=\begin{dcases*}
-\tau_k^P, & if $\theta \in P_k$ for some $k \in \llbracket 0, m\rrbracket$,\\
-\tau_k^N, & if $\theta \in N_k$ for some $k \in \llbracket 0, m\rrbracket$,
\end{dcases*}
\]
where, for $k \in \llbracket 0, m\rrbracket$, $-\tau_k^P$ and $-\tau_k^N$ are nonzero and belong to $I_k$ and the family $\{\tau_0^P,\dotsc,\tau_m^P,\allowbreak\tau_0^N,\dotsc,\tau_m^N\}$ is rationally independent. By construction, we have $\varphi \in \mathbf B$,  $\varphi_{\#}\mu \in \WellPosed$ thanks to Remark~\ref{remk:pushfwd-in-WP}, $\norm{ \varphi - \id }_\infty < \varepsilon$, and $\Sigma(\varphi_{\#} \mu)$ becomes
\begin{equation}
\label{eq:system-not-strongly-stable}
x(t) = \sum_{k=0}^m \left(a_k x(t - \tau_k^P) + b_k x(t - \tau_k^N)\right),
\end{equation}
where $a_k = \mu(P_k) = \mu^+(P_k) \geq 0$ and $b_k = \mu(N_k) = -\mu^-(N_k) \leq 0$.

Note that \eqref{eq:system-not-strongly-stable} is of the form \eqref{eq:intro-diff-eqn-sum}. Since $\sum_{k=0}^m \left(\abs{ a_k} + \abs{ b_k}\right) = \sum_{k=0}^m \left(\mu^+(P_k) + \mu^-(N_k)\right) = \mu^+(P) + \mu^-(N) = \V(\mu) \geq 1$, we deduce that \eqref{eq:system-not-strongly-stable} does not satisfy Proposition~\ref{prop:Hale-Silkowski}\ref{item:HS-criterion}. Hence, by Proposition~\ref{prop:Hale-Silkowski}\ref{item:HS-rat-indep}, \eqref{eq:system-not-strongly-stable} is not exponentially stable, yielding the conclusion.
\end{proof}

\subsection{Strong stability results in higher dimension}
\label{sec_main-strong-stab}

We shall prove in the sequel a partial result on Conjecture~\ref{GHS-Conjecture-SDG}, in which we replace the strong stability of $\Sigma(M)$ by slightly stronger notions which take into account a uniform decay rate.

\begin{definition}[$\mathsf K$-strong stability with uniform rate]
\label{def_stron-stab-UR}
Let $\mathsf K \subset \WellPosed$ be nonempty, $M \in \mathsf K$, and $\mu$ be the matrix-valued Borel measure associated with $M$ in the sense of Remark~\ref{remk:M-mu}.
\begin{enumerate}
\item We say that $\Sigma(M)$ is \emph{$\mathsf K$-locally strongly stable with uniform rate} if there exist $\alpha > 0$ and $\varepsilon > 0$ such that $\omega_{\varphi_{\#}\mu}\leq -\alpha$ for all $\varphi \in \mathbf B$ with $\norm{ \varphi - \id }_\infty < \varepsilon$  and $\varphi_{\#} \mu \in \mathsf K$. 
\item We say that $\Sigma(M)$ is \emph{$\mathsf K$-globally strongly stable with uniform rate} (or simply \emph{$\mathsf K$-strongly stable with uniform rate}) if there exists $\alpha > 0$ such that $\omega_{\varphi_{\#}\mu} \leq -\alpha$ for all $\varphi\in\mathbf{B}$ with $\varphi_{\#} \mu \in \mathsf K$.
\end{enumerate}
\end{definition}

\begin{remark}
\label{remk:TV-less-1-implies-global-strong-stab-unif-rate}
If $M \in \WellPosed$ and $V(M) < 1$, then $\Sigma(M)$ is $\WellPosed$-strongly stable with uniform rate. Indeed, the proof of Proposition~\ref{prop:TV-less-1-implies-global-strong-stab} also shows that one can provide an upper bound on the exponential rate of trajectories of $\Sigma(\varphi_{\#} \mu_M)$ independently of $\varphi \in \mathbf B$ with $\varphi_{\#} \mu_M \in \WellPosed$.
\end{remark}

With respect to Definition~\ref{def_strong-stab}, Definition~\ref{def_stron-stab-UR} requires $\Sigma(\varphi_{\#} \mu)$ not only to be exponentially stable, but also for its exponential rate $\omega_{\varphi_{\#} \mu}$ to be upper bounded by a negative constant $-\alpha$, uniformly with respect to $\mu$. This uniformity requirement is justified in view of the following theorem, which states that Conjecture~\ref{GHS-Conjecture-SDG} is true when strong stability is replaced by strong stability with uniform rate.

\begin{theorem}\label{GHS-Theorem}
Let $M \in \SDG$. Then the following statements are equivalent.
\begin{enumerate}
\item\label{item:GHSthm-rhoHS-less-1} We have $\rho_{\mathrm{HS}}(M) < 1$.
\item\label{item:GHSthm-GSSUR} System~$\Sigma(M)$ is $\SDG$-strongly stable with uniform rate.
\item\label{item:GHSthm-LSSUR} System~$\Sigma(M)$ is $\SDG$-locally strongly stable with uniform rate.
\end{enumerate}
\end{theorem}

The remainder of Section~\ref{sec_main-strong-stab} is dedicated to the proof of Theorem~\ref{GHS-Theorem}. We first prove some properties of $\rho_{\mathrm{HS}}$ in Section~\ref{sec_rhoHS}, before studying the characteristic function of a perturbed system $\Sigma(\varphi_{\#}\mu)$ in Section~\ref{sec_Delta}, and concluding the proof in Section~\ref{sec_proof-GHS-Theorem}, where further remarks are also provided.

\subsubsection{Properties of \texorpdfstring{$\rho_{\mathrm{HS}}$}{rhoHS}}
\label{sec_rhoHS}

We start by proving some properties of the quantity $\rho_{\mathrm{HS}}$ defined in \eqref{GHSC}. Our first result concerns the particular case in which $M$ is given by \eqref{eq_M-finite-sum} and $\rho_{\mathrm{HS}}$ reduces to \eqref{eq:rho-HS}. In this case, as $\theta_1, \dotsc, \theta_N$ are only used in $e^{i \theta_1}, \dotsc, e^{i \theta_N}$ in \eqref{eq:rho-HS}, the maximum in this expression can be seen as a maximum over $N$ complex numbers in the unit circle. Using a suitable version of the maximum principle for the spectral radius matrix-valued functions, one can replace the maximum in \eqref{eq:rho-HS} for a maximum over $N$ complex numbers in the closed unit disk $\overline{B}_{\mathbb C}(0, 1)$. More precisely, we have the following result.

\begin{proposition}
\label{prop_max_principle}
Let $A = (A_1, \dotsc, A_N) \in \mathcal M_d (\mathbb{R})^N$. Then
\[\max_{t \in \overline{B}_{\mathbb C}(0, 1)^N} \spr \left( \sum_{k=1}^N t_k A_k  \right) = \max_{\theta\in[0,2\pi]^N} \spr \left( \sum_{k=1}^N A_k e^{i\theta_k} \right).\]
\end{proposition}

\begin{proof}
We take $f(z_1, \dotsc,z_N)=\sum_{k=1}^N z_k A_k$ and apply Lemma~\ref{lemma_max_principle} from Appendix~\ref{app:other}.
\end{proof}

Our next aim is to generalize Proposition~\ref{prop_max_principle} by proving that the function $e^{i \xi(\cdot)}$ in \eqref{GHSC} can be replaced by a function taking values in the closed unit disk. For that purpose, we define the set
\begin{equation*}
\mathfrak D = \left\{\eta\colon [-1,0]\to \overline{B}_{\mathbb C}(0, 1) \suchthat \eta \text{ is Borel-measurable}\right\}.
\end{equation*}
Our first result in that direction is the following.

\begin{lemma}
\label{lemm:partial_max_principle}
Let $M\in \NBV([-1, 0],\mathcal M_d (\mathbb{R}))$ and suppose that 
\[
\sup_{\xi\in \mathfrak C}  \spr \left(\int_{-1}^0 e^{i\xi(\theta)} \diff M(\theta)\right) <1.
\]
Then 
\[
\sup_{\eta\in \mathfrak D} \spr \left(\int_{-1}^0 \eta(\theta) \diff M(\theta)\right) \leq 1.
\]
\end{lemma}

\begin{proof}
By contradiction, suppose that 
\[
\sup_{\eta\in \mathfrak D} \spr \left(\int_{-1}^0 \eta(\theta) \diff M(\theta)\right) > 1.
\]
Then there exists $\eta \in \mathfrak D$ such that $\spr \left(\int_{-1}^0 \eta(\theta) \diff M(\theta)\right) > 1$. From \cite[Theorem~2.10]{Folland1999Real} there exists a sequence $(\eta_n)_{n \in \mathbb N}$ of simple functions such that $\abs{\eta_n(\theta)} \leq \abs{\eta(\theta)}$ for every $\theta\in [-1,0]$ and $n \in \mathbb N$ and $\eta_n \to \eta$ pointwise as $n \to +\infty$. Then, from Lebesgue's dominated convergence theorem, we obtain
\[
\lim_{n \to +\infty} \int_{-1}^0 \eta_n(\theta) \diff \mu_M (\theta) = \int_{-1}^0 \eta(\theta) \diff \mu_M(\theta).
\]
By the continuity of the spectral radius, we obtain that there exists a simple function $\widetilde\eta \in \mathfrak D$ satisfying $\spr \left(\int_{-1}^0 \widetilde\eta(\theta) \diff M(\theta)\right) > 1$. Let us write $\widetilde\eta = \sum_{k=1}^m\alpha_k \1_{E_k}$ with $m \in \mathbb N$, $\alpha_k \in \overline{B}_{\mathbb C}(0, 1)$ and $E_k \in \mathfrak B_{[-1,0]}$ for $k \in \llbracket 1, m\rrbracket$ such that $\bigcup_{k=1}^m E_k=[-1,0]$ and $E_j\cap E_k=\emptyset$ for $j\ne k$. Using Proposition~\ref{prop_max_principle}, we have
\begin{align*}
1<\spr \left(\int_{-1}^0 \widetilde\eta(\theta) \diff M(\theta)\right) &=\spr \left(\sum_{k=1}^m \alpha_k \mu_M(E_k) \right)\leq \max_{t\in \overline{B}_{\mathbb C}(0, 1)^m} \spr \left( \sum_{k=1}^m t_k \mu_M(E_k)  \right)\\
&=  \max_{\theta\in[0,2\pi]^m} \spr \left( \sum_{k=1}^m \mu_M(E_k) e^{i\theta_k} \right)\leq \sup_{\xi\in \mathfrak C}  \spr \left(\int_{-1}^0 e^{i\xi(\theta)} \diff M(\theta)\right)<1,
\end{align*}
which gives the desired contradiction. 
\end{proof}

We are now in position to generalize Proposition~\ref{prop_max_principle} to $\rho_{\mathrm{HS}}$ from \eqref{GHSC}.

\begin{proposition}\label{prop_max_principle1}
Let $M \in \NBV([-1, 0],\mathcal M_d (\mathbb{R}))$. Then
\begin{equation}
\label{eq_supC-equal-supD}
\rho_{\mathrm{HS}}(M) = \sup_{\xi\in \mathfrak C}  \spr \left(\int_{-1}^0 e^{i\xi(\theta)} \diff M(\theta)\right)=\sup_{\eta\in \mathfrak D} \spr \left(\int_{-1}^0 \eta(\theta)\diff M(\theta)\right).
\end{equation}
\end{proposition}

\begin{proof}
Set $\rho_0= \rho_{\mathrm{HS}}(M)$. Note that, for every $\varepsilon > 0$,
\[
\sup_{\xi\in \mathfrak C}  \spr \left(\int_{-1}^0 e^{i\xi(\theta)} \diff\left( \frac{M(\theta)}{\rho_0+\varepsilon}\right)\right) = \frac{\rho_0}{\rho_0+\varepsilon} < 1.
\]
Therefore, by Lemma~\ref{lemm:partial_max_principle}, we have
\[
\sup_{\eta\in \mathfrak D}  \spr \left(\int_{-1}^0 \eta(\theta) \diff\left( \frac{M(\theta)}{\rho_0+\varepsilon}\right)\right) \leq 1,
\]
hence
\[
\sup_{\eta\in \mathfrak D} \spr \left(\int_{-1}^0 \eta(\theta) \diff M(\theta)\right) \leq \rho_0 + \varepsilon,\qquad \text{ for all } \varepsilon > 0.
\]
Therefore, as $\varepsilon > 0$ is arbitrary, we deduce that
\[
\sup_{\xi\in \mathfrak C}  \spr \left(\int_{-1}^0 e^{i\xi(\theta)} \diff M(\theta)\right)\leq \sup_{\eta\in \mathfrak D} \spr \left(\int_{-1}^0 \eta(\theta)\diff M(\theta)\right) \leq \rho_0 = \sup_{\xi\in \mathfrak C}  \spr \left(\int_{-1}^0 e^{i\xi(\theta)} \diff M(\theta)\right),
\]
as required.
\end{proof}

Our next result shows that the supremum over $\mathfrak D$ in the right-hand side of \eqref{eq_supC-equal-supD} is attained. Although this result is not directly used in the proof of Theorem~\ref{GHS-Theorem}, it can be obtained as a byproduct of our analysis and we provide it here since we deem that it has an interest on itself.

\begin{proposition}
Let $M\in \NBV([-1, 0],\mathcal M_d (\mathbb{R}))$. Then there exists $\eta_\ast\in \mathfrak D$ such that
\[\sup_{\eta\in \mathfrak D} \spr \left(\int_{-1}^0 \eta(\theta) \diff M(\theta)\right)=\spr \left(\int_{-1}^0 \eta_\ast(\theta) \diff M(\theta)\right).\]
\end{proposition}
\begin{proof}
Consider a sequence $(\eta_n)_{n \in \mathbb N}$ in $\mathfrak D$ such that 
\[
\lim_{n\to +\infty}\spr \left(\int_{-1}^0 \eta_n(\theta) \diff M(\theta)\right)=\sup_{\eta\in \mathfrak D} \spr \left(\int_{-1}^0 \eta(\theta) \diff M(\theta)\right).
\]
Let $\mu$ be the matrix-valued measure associated with $M$ in the sense of Remark~\ref{remk:M-mu} and denote its entries by $\mu_{jk}$, $(j, k) \in \llbracket 1, d\rrbracket^2$. By the Jordan decomposition theorems, for each $(j, k) \in \llbracket 1, d\rrbracket^2$, the signed measure $\mu_{jk}$ can be written as $\mu_{jk} = \mu_{jk}^+ - \mu_{jk}^-$, where $\mu_{jk}^+$ and $\mu_{jk}^-$ are mutually singular positive measures. To simplify the notation in the sequel, we fix a reindexing $\mu_1, \dotsc, \mu_K$ of the family $\{\mu_{jk}^+, \mu_{jk}^- \suchthat (j, k) \in \llbracket 1, d\rrbracket^2\}$, with $K = 2 d^2$. Denoting the measurable space $([-1,0],\mathfrak{B}_{[-1,0]})$ by $\mathcal X$ and writing simply $L^\infty(\mu_k)$ for $L^\infty(\mathcal X,\mathbb{C},\mu_k)$, we have that $\eta_n\in L^\infty(\mu_k)$ and $\norm{\eta_n}_{L^\infty(\mu_k)}\leq 1$ for every $k\in \llbracket 1, K \rrbracket$. Using the fact that the topological dual of $L^1(\mu_k)$ can be identified with $L^\infty(\mu_k)$ and applying the Banach--Alaoglu--Bourbaki theorem, we deduce that, up to extracting a subsequence, there exist $n_1, \dotsc,n_K$ in $\mathfrak D$ such that
\[
\eta_n\xrightharpoonup[n \to +\infty]{*} n_k \text{ in } L^\infty(\mu_k),\qquad \text{ for all } k\in \llbracket 1, K\rrbracket.
\]
By Proposition~\ref{prop_measure2} in Appendix~\ref{app:measure}, there exists $\eta_\ast\in \mathfrak D$ such that, for all $ k\in \llbracket 1, K\rrbracket$, we have $\eta_n\xrightharpoonup{*} \eta_\ast$ in $L^\infty(\mu_k)$ as $n \to +\infty$, that is,
\[
\lim_{n\to +\infty}\int_{-1}^0\eta_n\phi \diff \mu_k =\int_{-1}^0\eta_\ast\phi \diff \mu_k, \qquad \text{ for all } \phi \in L^1(\mu_k).
\]
Taking $\phi\equiv1$ and using the continuity of the spectral radius, we obtain that
\[
\spr\left(\int_{-1}^0 \eta_\ast(\theta)\diff \mu(\theta)\right)=\lim_{n \to +\infty}\spr \left(\int_{-1}^0 \eta_n(\theta)\diff \mu(\theta)\right)=\sup_{\eta\in \mathfrak D} \spr \left(\int_{-1}^0 \eta(\theta) \diff M(\theta)\right),
\]
as required.
\end{proof}

\subsubsection{Characteristic function of perturbed systems}
\label{sec_Delta}

Out next preparatory results for the proof of Theorem~\ref{GHS-Theorem} concern properties of the characteristic functions of perturbed systems $\Sigma(\varphi_{\#} \mu_M)$ from \eqref{sys_perturbed-measure}. We first prove that the condition $\rho_{\mathrm{HS}}(M) < 1$ ensures that such characteristic functions admit no roots in a neighborhood of the imaginary axis, whose size is independent of the perturbation $\varphi$.

\begin{lemma}\label{lemma_noroots}
Let $M\in \NBV([-1, 0],\mathcal M_d (\mathbb{R}))$, denote by $\mu$ the matrix-valued Borel measure associated with $M$, and suppose that $\rho_{\mathrm{HS}}(M)<1$. Then there exists $\nu>0$ such that, for all $\varphi \in \mathbf B$, there are no roots of the characteristic function of $\varphi_{\#}\mu$ on the vertical strip $\{s \in \mathbb{C} \suchthat \abs{\Real s} \leq \nu\}$.
\end{lemma}

\begin{proof}
Note that, for a given $\varphi \in \mathbf B$, using the properties of a pushforward measure, we have that the characteristic function $\Delta_\varphi$ of $\varphi_{\#} \mu$ is given by
\[\Delta_\varphi(s) = \det\left(\Id-\int_{-1}^0 e^{s\varphi(\theta)} \diff M(\theta) \right), \qquad s \in \mathbb C.\]

Consider the set 
\[\mathcal{T} = \left\{T_{\varphi,b}(a)\in \mathcal{M}_d(\mathbb{R})\suchthat -1 \leq a \leq 1,\, b\in \mathbb{R},\text{ and } \varphi \in \mathbf B\right\},\]
where
\[
T_{\varphi,b}(a)=\int_{-1}^0 e^{(a + i b)\varphi(\theta)}\diff M(\theta).
\]
Note that $\mathcal{T}$ is a bounded set since, by Proposition~\ref{proposition:tvcoincides} in Appendix~\ref{app:measure}, for any $T_{\varphi,b}(a)\in \mathcal{T}$, we have
\[
\abs*{T_{\varphi,b}(a)}\leq \V_{-1}^{0}(M) \sup_{\theta\in[-1,0]}\abs*{e^{(a+ib)\varphi(\theta)}} \leq \V_{-1}^{0}(M) e.
\]
Therefore, the closure $\overline{\mathcal{T}}$ of $\mathcal{T}$ is compact and the restriction $\spr_{\rvert \overline{\mathcal{T}}}$ of the spectral radius function to $\overline{\mathcal T}$ is uniformly continuous. Let $\varepsilon=\frac{1-\rho_{\mathrm{HS}}(M)}{2}$. Then there exists $\delta>0$ such that, for any $A,B\in \overline{\mathcal{T}}$ such that $\abs{A-B}<\delta$, we have
\[\spr(A)<\spr(B)+\varepsilon.\]
On the other hand, for $T_{\varphi,b}(a),T_{\varphi,b}(0)\in \mathcal{T}$,
\begin{align*}
\abs*{ T_{\varphi,b}(a)-T_{\varphi,b}(0)}&=\abs*{\int_{-1}^0(e^{a\varphi(\theta)}-1)e^{ib\varphi(\theta)}\diff M(\theta)}\leq \V_{-1}^{0}(M)\sup_{\theta\in [-1,0]}\abs*{e^{a \varphi(\theta)} - 1}\\
&\leq \V_{-1}^{0}(M) \abs*{e^{-a} - 1}.
\end{align*}  
Let $\nu \in (0, 1]$ be such that, for every $a\in [-\nu,\nu]$, we have $\V_{-1}^{0}(M) \abs{e^{-a} - 1} < \delta$. Thus 
\[
\spr\left(T_{\varphi,b}(a)\right) < \spr\left(T_{\varphi,b}(0)\right) + \varepsilon \leq \frac{1 + \rho_{\mathrm{HS}}(M)}{2} < 1,
\]
for every $a\in [-\nu,\nu]$, $b\in \mathbb{R}$, and $\varphi \in \mathbf B$. Therefore, $\Id - T_{\varphi,b}(a)$ is invertible for all $a\in [-\nu,\nu]$, $b\in \mathbb{R}$, and $\varphi \in \mathbf B$, resulting in
\[\det \left(\Id - T_{\varphi,b}(a)\right) = \det\left(\Id- \int_{-1}^0 e^{(a + i b) \varphi(\theta)} \diff M(\theta)\right)\ne 0. \qedhere\]
\end{proof}

The last technical result on characteristic functions of perturbed systems $\Sigma(\varphi_{\#} \mu_M)$ from \eqref{sys_perturbed-measure} that we shall need is the following.

\begin{lemma}\label{lemma_roots}
Let $M \in \NBV([-1, 0],\mathcal M_d (\mathbb{R}))$, denote by $\mu$ the matrix-valued Borel measure associated with $M$, and suppose that $\rho_{\mathrm{HS}}(M) \geq 1$. Then, for every $\varepsilon > 0$ and $\delta > 0$, there exists $\varphi \in \mathbf B$ such that $\norm{\varphi-\id}_\infty < \varepsilon$, $\varphi_{\#} \mu \in \NoSingular$, and such that the characteristic function of $\varphi_{\#} \mu$ admits a root $s$ with $\Real s \geq \ln \rho_{\mathrm{HS}}(M) - \delta$.
\end{lemma}

\begin{proof}
To simplify the notation in the proof, set $\rho_0 = \rho_{\mathrm{HS}}(M)$. Let $\varepsilon > 0$ and $\delta > 0$ be given. 
Note that it suffices to prove the result for $\varepsilon>0$ and $\delta>0$ small enough, and thus we assume, without loss of generality, that $\varepsilon < 1$, and also that $\delta < 2 \ln \rho_0$ in the case $\rho_0 > 1$. We set
\begin{equation}
\label{eq_delta-1}
\delta_1 = \begin{dcases*}
\frac{\delta}{2}, & if $\rho_0 > 1$, \\
\frac{\varepsilon \delta}{16}, & if $\rho_0 = 1$.
\end{dcases*}
\end{equation}

By Proposition~\ref{prop_max_principle1}, there exists a simple function $\eta_0\in \mathfrak D$ such that
\begin{equation}\label{eq_lemma_expinstab1}
\spr\left(\int_{-1}^0 \eta_0 (\theta) \diff M(\theta)\right)\ge \rho_0 e^{-\delta_1}.
\end{equation}
Without loss of generality, we write $\eta_0$ as
\[\eta_0=\sum_{k=1}^N \xi_k \1_{E_k},\]
with $N \in \mathbb N$, $E_k \in \mathfrak B_{[-1, 0]}$ and $\sup E_k - \inf  E_k<\frac{\varepsilon}{2}$ for all $k \in \llbracket 1, N\rrbracket$, and $\max_k \abs{\xi_k}\leq 1$. Now we take $\tau_1, \dotsc,\tau_N\in (0, 1)$ rationally independent so that, for all $k \in \llbracket 1, N\rrbracket$, we have
\[ -\tau_k \in \left(\max\left\{\inf E_k-\frac{\varepsilon}{4}, -1\right\}, \min\left\{\sup E_k+\frac{\varepsilon}{4}, -\frac{\varepsilon}{8}\right\}\right).\]
Note that all the above intervals have length at least $\frac{\varepsilon}{8}$ and that
\begin{equation}
\label{eq_bound-tau}
\tau_k > \frac{\varepsilon}{8} \qquad \text{ for all } k \in \llbracket 1, N\rrbracket.
\end{equation}

We define $\varphi \colon [-1,0] \to (-1, 0)$ so that $\varphi_{\rvert E_k}=-\tau_k$ for all $k \in \llbracket 1, N\rrbracket$. It clearly holds that $\varphi \in \mathbf B$, $\norm{\varphi - \id}_\infty<\varepsilon$, and $\varphi_{\#} \mu \in \NoSingular$ thanks to Remark~\ref{remk:pushfwd-in-WP}. Setting $A_k=\int_{E_k} \diff M(\theta)$ for every $k\in \llbracket 1, N\rrbracket$, the perturbed system $\Sigma(\varphi_{\#} \mu_M)$ from \eqref{sys_perturbed-measure} can be simply written as
\begin{equation}
\label{eq_lemma_expinstab-sum}
x(t) = \sum_{k=1}^N A_k x(t-\tau_k).
\end{equation}
On the other hand, from \eqref{eq_lemma_expinstab1} we obtain
\[
\spr\left(\sum_{k=1}^N \xi_k A_k\right)\ge \rho_0 e^{-\delta_1}.
\]
Set $\beta = \rho_0 e^{-\delta_1}$ and $\alpha = \ln \rho_0 - \frac{\delta}{2}$. It follows from \eqref{eq_delta-1} and \eqref{eq_bound-tau} that $0 < \frac{e^{\alpha \tau_k}}{\beta} \leq 1$ for every $k \in \llbracket 1, N\rrbracket$ and 
\[
\spr\left(\sum_{k=1}^N \frac{\xi_ke^{\alpha \tau_k}}{\beta} e^{-\alpha\tau_k}A_k\right)\ge \frac{\rho_0 e^{-\delta_1}}{\beta} = 1.
\]
From Proposition~\ref{prop_max_principle}, we have
\[
\max_{\theta\in[0,2\pi]^N} \spr \left( \sum_{k=1}^N e^{-\alpha\tau_k}A_k e^{i\theta_k} \right) \geq 1.
\]
Thus, from Proposition~\ref{prop:Hale-Silkowski}, the system 
\begin{equation}\label{eq_lemma_expinstab2}
y(t)=\sum_{k=1}^N e^{-\alpha \tau_k}A_k y(t-\tau_k)
\end{equation}
is not exponentially stable, implying that its characteristic function admits a root $z$ with $\Real z > -\frac{\delta}{2}$. Note that Systems~\eqref{eq_lemma_expinstab-sum} and \eqref{eq_lemma_expinstab2} are related by the change of variables $x(t) = e^{\alpha t} y(t)$, and hence the spectrum of \eqref{eq_lemma_expinstab-sum} is simply a translation of that of \eqref{eq_lemma_expinstab2} by $\alpha$. In particular, the characteristic function of  $\varphi_{\#} \mu$ admits a root $s$ with $\Real s > \alpha - \frac{\delta}{2} = \ln \rho_0 - \delta$.
\end{proof}

As an easy consequence of Lemma~\ref{lemma_roots}, we obtain the following result on local strong stability. Although this result is not directly used in the proof of Theorem~\ref{GHS-Theorem}, we present it here since it allows one to highlight other links between strong stability and the quantity $\rho_{\mathrm{HS}}(M)$ from \eqref{GHSC}.

\begin{corollary}
\label{coro:lss-implies-rho-leq-1}
Let $\mathsf K \subset \WellPosed$ be such that $\mathsf K \supset \NoSingular$ and take $M \in \mathsf K$. Assume that $\Sigma(M)$ is $\mathsf K$-locally strongly stable. Then $\rho_{\mathrm{HS}}(M) \leq 1$.
\end{corollary}

\begin{proof}
We proceed by contraposition and assume that $\rho_{\mathrm{HS}}(M) > 1$. Take $\varepsilon > 0$ and set $\delta = \frac{1}{2} \ln \rho_{\mathrm{HS}}(M) > 0$. Then, by Lemma~\ref{lemma_roots}, there exists $\varphi \in \mathbf B$ such that $\norm{\varphi - \id}_\infty < \varepsilon$, $\varphi_{\#} \mu_M \in \NoSingular \subset \mathsf K$, and such that the characteristic function of $\varphi_{\#} \mu_M$ admits a root with positive real part. Hence $S_{\varphi_{\#} \mu_M} > 0$ and, by Remark~\ref{remk:big-remk}\ref{item:big-remk-ineq-SM-omegaM}, we deduce that $\omega_{\varphi_{\#} \mu_M} > 0$, implying that $\Sigma(\varphi_{\#} \mu_M)$ is not exponentially stable. Hence $\Sigma(M)$ is not $\mathsf K$-locally strongly stable.
\end{proof}

\subsubsection{Proof of Theorem~\ref{GHS-Theorem} and further comments}
\label{sec_proof-GHS-Theorem}

Combining the results from Section~\ref{sec_rhoHS} and \ref{sec_Delta}, we are in position to prove Theorem~\ref{GHS-Theorem}

\begin{proof}[Proof of Theorem~\ref{GHS-Theorem}]
It is clear from the definition that \ref{item:GHSthm-GSSUR} implies \ref{item:GHSthm-LSSUR}.

To show that \ref{item:GHSthm-rhoHS-less-1} implies \ref{item:GHSthm-GSSUR}, let $\varphi \in \mathbf B$ be such that $\varphi_{\#} \mu_M \in \SDG$. From Lemma~\ref{lemma_noroots}, there exists $\alpha>0$ such that there are no roots of the characteristic function $\Delta_\varphi$ of $\varphi_* \mu_M$, namely
\[
\Delta_\varphi(s)=\det\left(\Id-\int_{-1}^0 e^{s\varphi(\theta)} \diff M(\theta)\right) \qquad s \in \mathbb C,
\]
on the vertical strip $\{s \in \mathbb{C} \suchthat \abs{\Real s} \leq \alpha\}$. We claim that $s \in \mathbb C$ is not a root of $\Delta_\varphi$ if $\Real s \ge 0$. Indeed, given such an $s$, for $\eta\colon[-1,0]\to \overline{B}_{\mathbb C}(0, 1)$ given by $\eta(\theta)=e^{s\varphi(\theta)}$, it is clear that $\eta \in \mathfrak D$, so from the hypothesis and Proposition~\ref{prop_max_principle1} we obtain
\[
\spr \left(\int_{-1}^0 e^{s\varphi(\theta)} \diff M(\theta)\right)\leq\sup_{\eta\in \mathfrak D} \spr \left(\int_{-1}^0 \eta(\theta) \diff M(\theta)\right)<1.
\]
This implies that $1$ is not an eigenvalue of $\int_{-1}^0 e^{s\varphi(\theta)} \diff M(\theta)$, i.e., $s$ is not a root of $\Delta_\varphi$. Hence $\Delta_\varphi$ admits no roots in the half-plane $\{s \in \mathbb C \suchthat \Real s \geq -\alpha\}$ and, since $\varphi_{\#} \mu_M \in \SDG$ by assumption, we deduce that $\omega_{\varphi_{\#}\mu}\leq -\alpha$, concluding the proof of \ref{item:GHSthm-GSSUR}.

To show that \ref{item:GHSthm-LSSUR} implies \ref{item:GHSthm-rhoHS-less-1}, we argue by contradiction and hence assume that \ref{item:GHSthm-rhoHS-less-1} does not hold true, i.e., $\rho_{\mathrm{HS}}(M) \geq 1$. Since $M$ is locally strongly stable with uniform rate, we can find $\alpha>0$ and $\varepsilon>0$ such that, for all $\varphi\in \mathbf{B}$ with $\norm{\varphi - \id}_{\infty}<\varepsilon$ and $\varphi_{\#}\mu_M\in \SDG$, all the roots $s$ of 
\[
\Delta_\varphi(s)=\det\left(\Id-\int_{-1}^0 e^{s\varphi(\theta)}\diff M(\theta)\right)
\]
satisfy $\Real s \leq -\alpha$. This is in contradiction with Lemma~\ref{lemma_roots}.
\end{proof}

\begin{figure}[ht]
\centering
\begin{tikzpicture}
\node (rhoHS) at (0, 0) {\small $\rho_{\mathrm{HS}}(M) < 1$};
\node[text width=3cm, align = center] (SSUR) at (4.5, 0) {\small $\Sigma(M)$ is $\SDG$-strongly stable with uniform rate};
\node[text width=3cm, align = center] (LSSUR) at (9.5, 0) {\small $\Sigma(M)$ is $\SDG$-locally strongly stable with uniform rate};
\node[text width=3cm, align = center] (SS) at (4.5, -3) {\small $\Sigma(M)$ is $\SDG$-strongly stable};
\node[text width=3cm, align = center] (LSS) at (9.5, -3) {\small $\Sigma(M)$ is $\SDG$-locally strongly stable};
\node (rhoHSleq) at (13.5, -3) {\small $\rho_{\mathrm{HS}}(M) \leq 1$};

\draw[double equal sign distance, double, Implies-Implies] (rhoHS) --node[midway, above] {\small Thm.~\ref{GHS-Theorem}} (SSUR);
\draw[double equal sign distance, double, Implies-Implies] (SSUR) --node[midway, above] {\small Thm.~\ref{GHS-Theorem}} (LSSUR);
\draw[double equal sign distance, double, -Implies] (SSUR) -- (SS);
\draw[double equal sign distance, double, -Implies] (LSSUR) -- (LSS);
\draw[double equal sign distance, double, -Implies] (SS) -- (LSS);
\draw[double equal sign distance, double, -Implies] (LSS) --node[midway, above] {\small Cor.~\ref{coro:lss-implies-rho-leq-1}} (rhoHSleq);
\end{tikzpicture}
\caption{Diagram of the results of the chapter on Conjecture~\ref{GHS-Conjecture-SDG}}
\label{fig_results}
\end{figure}
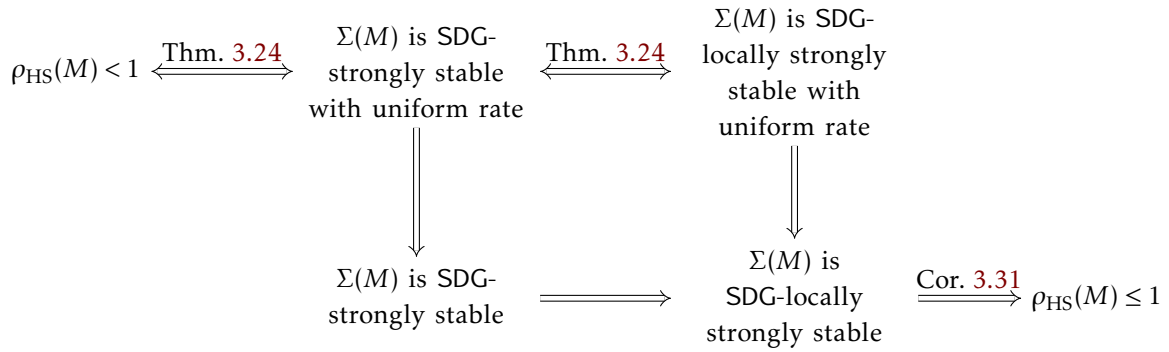

Regarding Conjecture~\ref{GHS-Conjecture-SDG}, the results presented in this chapter allow one to obtain the diagram given in Figure~\ref{fig_results} for a given $M \in \SDG$. Indeed, by Theorem~\ref{GHS-Theorem}, we have that $\rho_{\mathrm{HS}}(M) < 1$ implies $\SDG$-strong stability with uniform rate, which trivially implies $\SDG$-strong stability, the latter trivially implying $\SDG$-local strong stability. By Corollary~\ref{coro:lss-implies-rho-leq-1}, $\SDG$-local strong stability implies that $\rho_{\mathrm{HS}}(M) \leq 1$. Hence, Conjecture~\ref{GHS-Conjecture-SDG} would be proved if one could reinforce the conclusion of Corollary~\ref{coro:lss-implies-rho-leq-1} to $\rho_{\mathrm{HS}}(M) < 1$, at least in the case $\mathsf K = \SDG$. This would also imply that the notions of strong stability from Definition~\ref{def_strong-stab} coincide with the corresponding notions of strong stability with uniform rate from Definition~\ref{def_stron-stab-UR} in the case $\mathsf K = \SDG$. We refer to Section~\ref{sec:conjectures-strong-stab} for a deeper discussion on open problems on the strong stability of $\Sigma(M)$.

\section{Numerical illustrations}
\label{sec_numerical}

To illustrate the results of this chapter, we consider in this section System~\eqref{sys_general} in the scalare case $d = 1$ and under the assumption that the measure $\mu$ associated with $M$ through Remark~\ref{remk:M-mu} is absolutely continuous with respect to the Lebesgue measure with an affine density, given by $\theta \mapsto a + b \theta$ for some real constants $a$ and $b$. In this case, we have $M(\theta) = a(\theta + 1) + \frac{b}{2}(\theta^2 - 1)$ and \eqref{sys_general} becomes
\begin{equation}
\label{sys_a_b}
x(t) = \int_{-1}^0 x(t+\theta) (a+b\theta) \diff \theta.
\end{equation}
Our aim is to explicitly characterize the region of $\WellPosed$-strong stability described by Theorem~\ref{HS:scalar} in the $(a, b)$-plane, and also to numerically compute the region of exponential stability, in order to compare both.

\subsection{Region of strong stability}
\label{sec_example-strong}

Let $\theta_v = -\frac{a}{b}$ be the vertex of the parabola described by the function $M$ defined above (with the convention $\theta_v = 0$ if $b = 0$). We compute that $\V_{-1}^{0}(M) = \abs{M(-1) - M(\theta_v)} + \abs{M(\theta_v) - M(0)}$ if $\theta_v \in (-1,0)$ and $\V_{-1}^{0}(M) = \abs{M(-1) - M(0)}$ otherwise. Noticing that $\theta_v \in (-1, 0)$ if and only if $a (b - a) > 0$, we thus obtain, by Theorem~\ref{HS:scalar}, the following result.

\begin{proposition}
System~\eqref{sys_a_b} is $\WellPosed$-strongly stable if and only if $(a, b)$ belongs to the union of the sets
\begin{gather*}
\left\{(a, b) \in \mathbb{R}^2 \suchthat (a - b)^2 + a^2 < 2 \abs{b} \text{ and }  a (b - a) > 0\right\}, \\
\left\{(a, b) \in \mathbb{R}^2 \suchthat \abs*{a - \frac{b}{2}} < 1 \text{ and }  a (b - a) \leq 0\right\}.
\end{gather*}
\end{proposition}

\subsection{Region of exponential stability}
\label{sec_example-expo}

Let us now describe how we numerically approximate the region in the $(a, b)$-plane where \eqref{sys_a_b} is exponentially stable. From \cite[Chapter~12, Corollary~3.1]{Hale1977Theory}, we have that \eqref{sys_a_b} is exponentially stable if and only if there exists $\delta > 0$ such that all the zeros $s$ of the characteristic equation
\[\Delta_{a,b}(s) = 1 - \int_{-1}^0 (a + b \theta) e^{s \theta} \diff \theta\]
satisfy $s \lambda \leq -\delta$. After a straightforward computation, we have that $\Delta_{a,b}(s) = \frac{Q(s)}{s^2}$, where
\[Q(s) = s^2 - a s + b + \big((a - b)s - b \big)e^{-s}.\]
Since $\Delta_{a, b}$ is holomorphic over $\mathbb C$, $s=0$ is a double root of $Q$, and the other roots of $Q$ coincide with the roots of $\Delta_{a,b}$.

We first notice that one can exclude a half-space from the exponential stability region of \eqref{sys_a_b}.

\begin{proposition}
System~\eqref{sys_a_b} is not exponentially stable on the half-space $\{(a, b) \in \mathbb{R}^2 \suchthat 2 - 2 a + b \leq 0\}$.
\end{proposition}

\begin{proof} 
We have $Q(0) = Q^\prime(0) = 0$ and $Q^{\prime\prime}(0) = 2 - 2 a + b$, and thus $s=0$ is a triple root of $Q$ if $2-2a+b=0$, that is, $s=0$ is a root of $\Delta_{a,b}$ and \eqref{sys_a_b} is not exponentially stable. Now, if $Q^{\prime\prime}(0) = 2 - 2 a + b < 0$, we have that the restriction of $Q$ to the real line has a local maximum at $s=0$. So there is $s_1>0$ such that $Q(s_1) < 0$. On the other hand, 
$\lim_{s\to +\infty} Q(s)=+\infty$, so by continuity there is $s_2>s_1>0$ such that $Q(s_2)=0$. Hence $s_2$ is a positive real root of $\Delta_{a,b}$ and \eqref{sys_a_b} is not exponentially stable.
\end{proof}

We also remark that $Q$ can be seen as the characteristic function of a second-order time-delay system of retarded type and, as such, it admits only finitely many roots in any vertical strip of the complex plane (see, e.g., \cite[Corollary~1.9]{Michiels2014Stability}). In particular, \eqref{sys_a_b} is exponentially stable if and only if $\Real s < 0$ for all roots $s$ of $\Delta_{a, b}$. The function $Q$ admits the representation
\[
Q(s) = \det\left[s \Id - A_0 - A_1 e^{-s}\right]
\]
where
\[A_0=\begin{pmatrix}
0 & 1\\ b & -a
\end{pmatrix} \text{ and } A_1=\begin{pmatrix}
0 & 0\\ -b & a-b
\end{pmatrix},\]
and thus, by \cite[Proposition~1.10]{Michiels2014Stability}, we have that, if $s$ is a root of $Q$ with $\Real s \geq 0$, then
\[
\abs{s} \leq \abs{A_0} + \abs{A_1},
\]
where $\abs{\cdot}$ is an arbitrary induced matrix norm or the Frobenius norm. Taking into account that roots of $Q$ appear in complex conjugate pairs, we deduce that, given any root $s$ of $Q$ with nonnegative real part, either $s$ or $\overline s$ belong necessarily to the square $S_\alpha = \{s \in \mathbb C \suchthat 0 \leq \Real s \leq \alpha,\, 0 \leq \Imag s \leq \alpha\}$, where $\alpha = \abs{A_0} + \abs{A_1}$.

In order to numerically approximate the region of exponential stability of \eqref{sys_a_b} in the $(a, b)$-plane, we thus discretize the region of interest into a grid $(a_j, b_k)$, for $j$ and $k$ in suitable finite sets. For each pair $(a_j, b_k)$, the system is not exponentially stable if $2 - 2 a_j + b_k \leq 0$. Otherwise, we compute numerically the roots of $Q$ in the square $S_\alpha$ using the QPmR algorithm from \cite{Vyhlidal2014QPmR}\footnote{More precisely, we use the Python implementation available at \url{https://github.com/DSevenT/QPmR}.}. If $Q$ admits roots in this square with nonnegative real part other than the double root at $0$, then the system is not exponentially stable, otherwise it is exponentially stable.

\subsection{Comparisons}

\begin{figure}[htp]
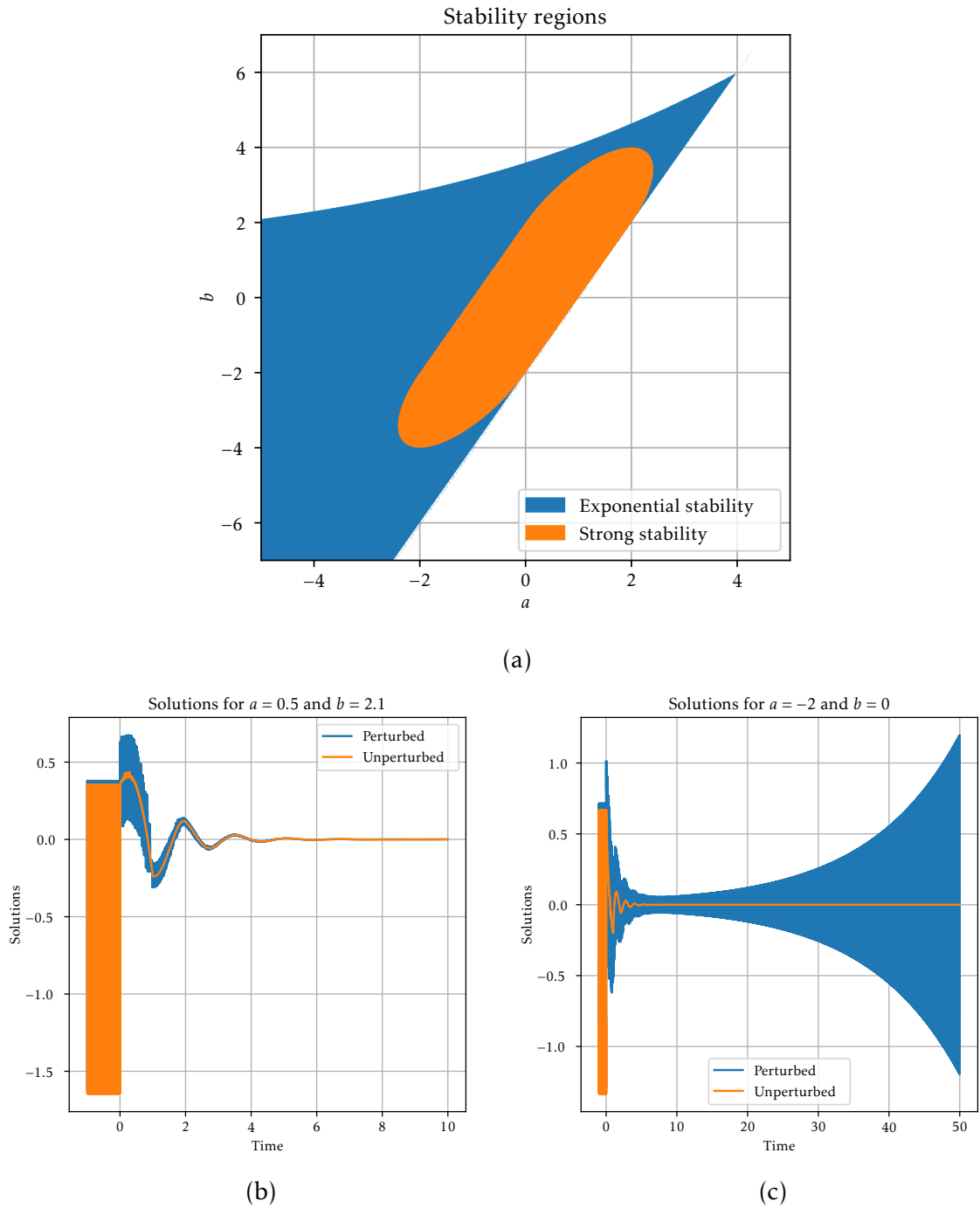

\centering
\begin{tabular}{@{} >{\centering} m{0.5\textwidth} @{} >{\centering} m{0.5\textwidth} @{}}
\multicolumn{2}{@{} >{\centering} m{\textwidth} @{}}{\resizebox{0.6666667\textwidth}{!}{\input{Figures/stability.pgf}}} \tabularnewline
\multicolumn{2}{@{} >{\centering} m{\textwidth} @{}}{(a)} \tabularnewline
\resizebox{0.5\textwidth}{!}{\input{Figures/trajectories_strongly_stable.pgf}} & \resizebox{0.5\textwidth}{!}{\input{Figures/trajectories_exponentially_stable.pgf}} \tabularnewline
(b) & (c) \tabularnewline
\end{tabular}
\caption[Stability regions and solutions of \eqref{sys_a_b} and \eqref{sys_a_b_perturbed} for two pairs $(a, b)$]{(a) Stability regions and (b, c) solutions of \eqref{sys_a_b} (orange curves) and of \eqref{sys_a_b_perturbed} (blue curves) (b) for the pair $(a, b) = (0.5, 2.1)$ in the $\WellPosed$-strong stability region and (c) for the pair $(a, b) = (-2, 0)$ in the exponential stability region, but outside of the $\WellPosed$-strong stability region}
\label{stab_regions}
\end{figure}

Figure~\ref{stab_regions}(a) represents the regions of exponential stability (in blue) and of $\WellPosed$-strong stability (in orange) of \eqref{sys_a_b}, computed according to the discussions in Sections~\ref{sec_example-strong} and \ref{sec_example-expo}.

In order to illustrate the effects of the lack of $\WellPosed$-strong stability, Figure~\ref{stab_regions}(b) and (c) provide numerical simulations of solutions of \eqref{sys_a_b} and of the perturbed system
\begin{equation}
\label{sys_a_b_perturbed}
x(t) = \int_{-1}^0 x(t + \varphi(\theta)) (a + b \theta) \diff\theta.
\end{equation}
We selected a perturbation $\varphi\colon [-1, 0] \to [-1, 0]$ which is constant in each interval of the form $(-\frac{j}{N}, -\frac{j-1}{N})$ for $j \in \llbracket 1, N\rrbracket$ and with $N = 25$, with randomly chosen values in each of these intervals in such a way that $\norm{ \varphi - \id}_\infty \leq \frac{1}{N} = 0.04$. As the lack of strong stability usually results from the apparition of unstable roots with large imaginary part in the perturbed system, we selected an initial condition of the form $x_0(t) = \cos(2 k \pi t) + c_0$, with $k = 127$ and $c_0 \in \mathbb R$ chosen so that $x_0 \in \mathcal C_0$. We observe, in Figure~\ref{stab_regions}(b) and (c), that the corresponding solutions of \eqref{sys_a_b} converge exponentially to the origin in both cases $(a, b) = (0.5, 2.1)$ and $(a, b) = (-2, 0)$, as expected, but, for the perturbed system \eqref{sys_a_b_perturbed}, the corresponding solution converges to the origin in the $\WellPosed$-strongly stable case $(a, b) = (0.5, 2.1)$, but it diverges in the case $(a, b) = (-2, 0)$ where $\WellPosed$-strong stability does not hold true.

{
\renewcommand{\thesection}{\thechapter.A}
\section{Appendix}
\label{sec:app}

\subsection{Results in measure theory}
\label{app:measure}

We gather here some measure-theoretical results needed in the chapter. The first one asserts that some classical properties for scalar-valued functions of bounded variation and their associated Borel measures remain true in the matrix-valued context of this chapter. In order to state that result, for a given matrix-valued Borel measure $\mu$ over the interval $[a, b]$, we set 
\[
\widetilde{\V}_a^b(\mu) = \sup\left\{\sum_{j=1}^{+\infty} \abs{\mu(E_j)} \suchthat  E_1, E_2, \dotsc \in \mathfrak{B}_{[a,b]} \text{ are pairwise disjoint and } \bigcup_{j=1}^{+\infty} E_j = [a,b]\right\}.
\]
We then have the following proposition (see also \cite[Chapter~3, Theorem~5.6]{Gripenberg1990Volterra} for the second part of the statement).

\begin{proposition}\label{proposition:tvcoincides}
Let $M \in \NBV([a, b], \mathcal M_d (\mathbb{R}))$ and $\mu_M$ be the measure associated with $M$ in the sense of Remark~\ref{remk:M-mu}. Then
\begin{equation}
\label{eq_variations-coincide}
\V_{a}^{b}(M) = \V_{a}^{b}(\mu_M) = \widetilde{\V}_a^b(\mu_M).
\end{equation}
Moreover, for $\phi\colon [a,b]\to \mathbb{R}^d$ and $\psi\colon [a,b]\to \mathbb{C}$ bounded Borel-measurable functions, we have
\[\abs*{\int_a^b \diff M(\theta)\phi(\theta)} \leq \V_{a}^{b}(M) \norm{\phi}_{\infty} \quad \text{ and } \quad \abs*{\int_a^b \psi(\theta) \diff M(\theta)} \leq \V_{a}^{b}(M) \norm{\psi}_{\infty}.\]
\end{proposition}

\begin{proof}
Given a partition $a=t_1<t_2<\dotsb<t_{k+1}=b$ of $[a, b]$, we have
\[M(t_{j+1})-M(t_{j})=\begin{dcases*}
M(t_{2}) = \mu_M([t_1,t_2]), & if $j = 1$, \\
\mu_M([a,t_{j+1}])-\mu_M([a,t_j])=\mu_M((t_j,t_{j+1}]), & if $j \in \llbracket 2, k\rrbracket$.
\end{dcases*}\]
So the sets defined by $E_1=[t_1,t_2]$ and $E_j=(t_j, t_{j + 1}]$, $j \in \llbracket 2, k\rrbracket$, are pairwise disjoint, belong to $\mathfrak{B}_{[a,b]}$, and $\bigcup_{i=1}^k E_i=[a,b]$. Hence $\V_{a}^{b}(M)\leq \V_{a}^{b}(\mu_M)$.

To prove the converse inequality, we define $T_M(x)=\V_{a}^{x}(M)$ for $x \in [a, b]$. We claim that $\abs{M(s)-M(t)}\leq T_M(s)-T_M(t)$ for every $s,t\in [a,b]$ with $s\geq t$. This is clearly true for $t = a$ or $t = s$ and, for $t > a$ and $s > t$, it is not difficult to see that
\begin{align*}
T_M(s) & =\sup\left\{\sum_{j=1}^k \abs{M(t_{j+1})-M(t_{j})}\suchthat k\in \mathbb N^\ast \text{ and }a=t_1<t_2<\dotsb<t_{k+1}=s \right\}\\
& = \sup\left\{\sum_{j=1}^k \abs{M(t_{j+1})-M(t_{j})}\suchthat k\in \mathbb N^\ast \text{ and }a=t_1<\dotsb<t_{\ell+1}=t<\dotsb<t_{k+1}=s \right\}.
\end{align*}
From this, it follows that
\begin{align*}
T_M(s) & \ge \sup\left\{\sum_{j=1}^{\ell} \abs{M(t_{j+1})-M(t_{j})}\suchthat \ell\in \mathbb N^\ast \text{ and }a=t_1<\dotsb<t_{\ell+1}=t \right\} + \abs{M(s)-M(t)}\\  & = T_M(t)+\abs{M(s)-M(t)},
\end{align*}
as desired. Using the continuity of $\abs{\cdot}$, we obtain, for $s\geq t$, that
\[\abs {M(\widehat{s})-M(\widehat{t})} \leq T_M(\widehat{s})-T_M(\widehat{t})\]
where $\widehat{s}\in \{s^-,s^+,s\}$ and $\widehat{t}\in \{t^-,t^+,t\}$. Let $I$ be an interval contained in $[a,b]$ with endpoints $s\geq t$. Using the continuity from below and above of $\mu_M$ and $\mu_{T_M}$, we obtain, after some calculations, that
\[\mu_M(I)=M(\widehat{s}(I))-M(\widehat{t}(I)) \text{ and } \mu_{T_M}(I)=T_M(\widehat{s}(I))-T_M(\widehat{t}(I)),\]
where
\[\widehat t(I) = 
\begin{dcases*}
t^-, & if $t\ne a$ and $t \in I$,\\
t^+, & if $t=a$ and $t\notin I$,\\
t, & otherwise,
\end{dcases*}
\qquad \widehat s(I) = 
\begin{dcases*}
s, & if $s \in I$ and $I \ne [a, a]$, \\
s^-, & if $s\notin I$ and $I \ne [a, a]$,\\
s^+, & if $I = [a, a]$.
\end{dcases*}
\]
Therefore
\begin{equation}
\label{eq_mu_M-leq-mu_TM}
\abs{\mu_M(I)} \leq \mu_{T_M}(I),
\end{equation}
for every interval $I$ contained in $[a,b]$. The set $\mathcal A$ of all pairwise disjoint finite unions of intervals of $[a,b]$ is an algebra over $[a,b]$ and, thanks to \eqref{eq_mu_M-leq-mu_TM}, any element $E$ of $\mathcal A$ satisfies $\abs{\mu_M(E)} \leq \mu_{T_M}(E)$. If we set $\mathcal{D} = \left\{E\in \mathfrak{B}_{[a,b]}\suchthat \abs{\mu_M(E)}\leq \mu_{T_M}(E)\right\}$, then $\mathcal{A}\subset \mathcal{D}$. Given $\{A_j\}_{j=1}^{+\infty} \subset \mathcal{D}$ with $A_1\subset A_2\subset \dotsb$, using the continuity from below of $\mu_M$ and $\mu_{T_M}$, we have
\begin{gather*}
\abs*{\mu_M\Bigg(\bigcup_{j=1}^{+\infty} A_j\Bigg)}=\abs*{\lim_{j\to +\infty}\mu_M (A_j)}=\lim_{j\to +\infty}\abs*{\mu_M (A_j)}\leq \lim_{j\to +\infty}\mu_{T_M} (A_j)=\mu_{T_M}\Bigg(\bigcup_{j=1}^{+\infty} A_j\Bigg).
\end{gather*}
Similarly, given $\{B_j\}_{j=1}^{+\infty}\subset\mathcal D$ with $B_1\supset B_2\supset \dotsb$, we have
\[
\abs*{\mu_M\Bigg(\bigcap_{j=1}^{+\infty} B_j\Bigg)}\leq \mu_{T_M}\Bigg(\bigcap_{j=1}^{+\infty} B_j\Bigg).
\]
Therefore $\mathcal{D}$ is a monotone class containing $\mathcal{A}$ and, from the monotone class theorem for sets, the smallest monotone class containing $\mathcal{A}$ is precisely the $\sigma$-algebra generated by $\mathcal{A}$, which is exactly $\mathfrak{B}_{[a,b]}$. Hence $\mathcal{D}=\mathfrak{B}_{[a,b]}$. Finally,
\begin{align*}
\V_{a}^{b}(\mu_M) &=\sup\left\{\sum_{j=1}^k \abs{\mu_M(E_j)}\suchthat k\in \mathbb N, E_1, \dotsc,E_{k}\in \mathfrak{B}_{[a,b]} \text{ are pairwise disjoint, } \bigcup_{j=1}^{k}E_j =[a,b]\right\}\\
& \leq \sup\left\{\sum_{j=1}^k \mu_{T_M}(E_j)\suchthat k\in \mathbb N, E_1, \dotsc,E_{k}\in \mathfrak{B}_{[a,b]} \text{ are pairwise disjoint, } \bigcup_{j=1}^{k}E_j =[a,b]\right\}\\
& =\sup\left\{ \mu_{T_M}\left(\bigcup_{j=1}^k E_j\right)\suchthat k\in \mathbb N, E_1, \dotsc,E_{k}\in \mathfrak{B}_{[a,b]} \text{ are pairwise disjoint, } \bigcup_{j=1}^{k}E_j =[a,b]\right\}\\
& =\mu_{T_M}([a,b])=T_M(b)=\V_{a}^{b}(M),
\end{align*}
yielding the first equality in \eqref{eq_variations-coincide}. 

It is clear that $\V_{a}^{b}(\mu_M)\leq \widetilde{\V}_a^b(\mu_M)$. Given a family $\{E_j\}_{j \in \mathbb N}$ of pairwise disjoint sets in $\mathfrak B_{[a, b]}$ such that $\bigcup_{j=0}^\infty E_j = [a,b]$, we denote $\widetilde{E_k}=\bigcup_{j=k}^{+\infty} E_j$ for every $k \in \mathbb N$. Then, for every $k \in \mathbb N$, we have
\[\sum_{j=1}^{k} \abs{\mu_M(E_j)}+\abs{\mu_M(\widetilde E_{k+1})}\leq \V_{a}^{b}(\mu_M).\]
From the continuity from below of $\mu_M$, we have $\lim_{k\to +\infty} \abs{\mu_M(\widetilde E_{k+1})} = 0$, so, taking the limit as $k \to +\infty$ in the above expression, we get
\[\sum_{i=1}^{+\infty} \abs{\mu_M(E_i)}\leq \V_{a}^{b}(\mu_M).\]
Since the family $\{E_j\}_{j \in \mathbb N}$ of pairwise disjoint sets in $\mathfrak B_{[a, b]}$ such that $\bigcup_{j=0}^\infty E_j = [a,b]$ is arbitrary, we deduce that $\widetilde{\V}_a^b(\mu_M) \leq \V_{a}^{b}(\mu_M)$, which gives the second equality in \eqref{eq_variations-coincide}. 

Proceeding now as in \cite[Section~4.1]{Cohn2013Measure}, let $\phi\colon[a, b] \to \mathbb{R}^d$ be a Borel-measurable simple function. Writing $\phi=\sum_{j=1}^k \phi_j \1_{E_j}$ for some $E_1, \dotsc, E_k$ in $\mathfrak B_{[a, b]}$ and $\phi_1, \dotsc, \phi_k$ in $\mathbb R^d$, we have
\[
\abs*{\int_a^b \diff M(\theta)\phi(\theta)}=\abs*{\sum_{j=1}^k\mu_M(E_j)\phi_j}\leq \sum_{j=1}^k \abs{\mu_M(E_j)}\abs{\phi_j} \leq \V_{a}^{b}(M) \norm{\phi}_\infty.
\]
Since every bounded Borel-measurable function is a limit of Borel-measurable simple functions in the uniform norm, we obtain desired result. The second inequality is analogous.
\end{proof}

We now provide the explicit expression of the total variation of a matrix-valued function $M$ belonging to the class $\NoSingular$ from \eqref{eq_def-N-L-G}. 

\begin{proposition}\label{prop_tvnosingular}
Let $M \in \NoSingular$ and write
\[
M(\theta) = \sum_{k=1}^{+\infty} A_k \1_{[-\tau_k,0]}(\theta)+N(\theta),
\]
where $\tau_k \in (0, 1)$ and $A_k\in \mathcal M_d (\mathbb{R})$ for $k \in \mathbb N^\ast$, $\tau_k \neq \tau_j$ if $k \neq j$, and $N$ is absolutely continuous. Then $\sum_{k=1}^{+\infty} \abs{A_k} < +\infty$ and
\begin{equation}\label{eq_tvequality}
\V_{-1}^{0}(M) = \sum_{k=1}^{+\infty} \abs{A_k} + \int_{-1}^0 \abs{N^\prime(\theta)} \diff\theta.
\end{equation}
\end{proposition}

\begin{proof}
Let $\widetilde N (\theta)=M(\theta)-N(\theta)$. Then $\V(\widetilde N)\leq \V(M) + \V(N) < +\infty$. The measure $\mu_{\widetilde N}$ associated with $\widetilde N$ in the sense of Remark~\ref{remk:M-mu} is given by
\[
\mu_{\widetilde N}= \sum_{k=1}^{+\infty} A_k \delta_{-\tau_k}.
\]
If we set $E_k=\{-\tau_k\}$ for $k\ge 2$ and $E_1=[-1,0] \setminus \bigcup_{k=2}^{+\infty} E_k$, then, from Proposition~\ref{proposition:tvcoincides},
\[
\sum_{k=1}^{+\infty} \abs{A_k} = \sum_{k=1}^{+\infty} \abs{\mu_{\widetilde N} (E_k)} \leq \V(\widetilde N) < +\infty.
\]

Let us prove first that, for $m \in \mathbb N^\ast$, if $M_m$ is the function defined by
\[
M_m(\theta)=\sum_{k=1}^m A_k \1_{[-\tau_k,0]}(\theta)+N(\theta),\qquad \theta\in[-1,0],
\]
then
\begin{equation}\label{eq_tvequalitytruncated}
\V_{-1}^{0}(M_m)  = \sum_{k=1}^m \abs{A_k} + \int_{-1}^0 \abs{N^\prime(\theta)} \diff\theta.
\end{equation}
Indeed, up to reordering $\tau_1, \dotsc, \tau_m$, we assume that $-1< -\tau_1<\dotsb< -\tau_m<0$ and we set by convention $\tau_0 = 1$. We have
\begin{equation*}
\V_{-1}^{0}(M_m) =  \sum_{k=0}^{m-1} \V_{-\tau_k}^{-\tau_{k+1}}(M_m)+\V_{-\tau_m}^{0}(M_m).
\end{equation*}
We first show that, for $k \in \llbracket 0, m-1\rrbracket$, we have
\[
\V_{-\tau_k}^{-\tau_{k+1}}(M_m)=\abs{A_{k+1}}+\V_{-\tau_k}^{-\tau_{k+1}}(N).
\]
Given $\varepsilon>0$, by continuity of $N$, there exists $\delta>0$ such that, for every $t\in(-\tau_{k+1}-\delta,-\tau_{k+1})$ we have $\abs{N(-\tau_{k+1})-N(t)}<\frac{\varepsilon}{4}$. Take a partition $-\tau_{k}=t_1<\dotsb<t_{n+1}=-\tau_{k+1}$ of $[-\tau_k,-\tau_{k+1}]$ such that $\sum_{j=1}^n \abs{N(t_{j+1})-N(t_{j})} \geq \V_{-\tau_k}^{-\tau_{k+1}}(N) - \frac{\varepsilon}{2}$. Without loss of generality, we may assume $t_n\in(-\tau_{k+1}-\delta,-\tau_{k+1})$. Then
\begin{align*}
\sum_{j=1}^n \abs{M_m(t_{j+1})-M_m(t_j)} &= \sum_{j=1}^{n-1}\abs{N(t_{j+1})-N(t_j)}+\abs{A_{k+1}+N(t_{n+1})-N(t_n)} \\
&\geq\sum_{j=1}^{n-1} \abs{N(t_{j+1})-N(t_j)} + \abs{A_{k+1}} -  \abs{N(t_{n+1})-N(t_n)}\\
&=\sum_{j=1}^{n} \abs{N(t_{j+1})-N(t_j)} + \abs{A_{k+1}} -  2\abs{N(-\tau_{k+1})-N(t_n)}\\
&\geq \V_{-\tau_k}^{-\tau_{k+1}}(N) -\frac{\varepsilon}{2} + \abs{A_{k+1}} -2\frac{\varepsilon}{4} \\
&= \V_{-\tau_k}^{-\tau_{k+1}}(N)  + \abs{A_{k+1}} -\varepsilon.
\end{align*}
Thus
\[\V_{-\tau_k}^{-\tau_{k+1}}(M_m)  \geq  \V_{-\tau_k}^{-\tau_{k+1}}(N)  + \abs{A_{k+1}} -\varepsilon,\]
for every $\varepsilon>0$, which implies
\[\V_{-\tau_k}^{-\tau_{k+1}}(M_m) \geq  \V_{-\tau_k}^{-\tau_{k+1}}(N)  + \abs{A_{k+1}},\qquad k \in \llbracket 0, m-1\rrbracket.\]
The reverse inequality is clear, hence we have
\begin{equation}\label{eq_prop_tvnosingula1}
\V_{-\tau_k}^{-\tau_{k+1}}(M_m) =  \V_{-\tau_k}^{-\tau_{k+1}}(N) + \abs{A_{k+1}},\qquad k \in \llbracket 0, m-1\rrbracket.
\end{equation}
On the other hand, since $M_m = N$ on $[-\tau_m, 0]$, we have $\V_{-\tau_m}^{0}(M_m) = \V_{-\tau_m}^{0}(N)$. Summing \eqref{eq_prop_tvnosingula1} over $k$ in $\llbracket 0, m-1\rrbracket$ and using the previous fact, we obtain the equality \eqref{eq_tvequalitytruncated} after recalling that $\V_{-1}^{0}(N) = \int_{-1}^0 \abs{N^\prime(\theta)} \diff\theta$ since $N$ is absolutely continuous.

To conclude, note that
\[\V(M-M_m)=\V\left(\sum_{k=m+1}^{+\infty} A_k \1_{[-\tau_k,0]}\right) \leq \sum_{k=m+1}^\infty \abs{A_k}\xrightarrow[m \to +\infty]{} 0.\]
Since $\abs{\V_{-1}^{0}(M)-\V_{-1}^{0}(M_m)}\leq \V_{-1}^{0}(M-M_m)$, we have
\[\V_{-1}^{0}(M)=\lim_{m\to +\infty}\V_{-1}^{0}(M_m)=\lim_{m\to +\infty} \left(\sum_{k=1}^m \abs{A_k} + \int_{-1}^0 \abs{N^\prime(\theta)} \diff\theta\right),\]
yielding \eqref{eq_tvequality}.
\end{proof}

Our next two results concern the uniqueness of weak-$*$ $L^\infty$ limits with respect to different measures. The first result, Lemma~\ref{lemma_measure1} below, is an intermediate step to get Proposition~\ref{prop_measure2}. To state and prove them, we use the standard notation $\mu \ll \nu$ to say that the measure $\mu$ is absolutely continuous with respect to the measure $\nu$, and $\mu \perp \nu$ to denote that the measures $\mu$ and $\nu$ are mutually singular.

\begin{lemma}\label{lemma_measure1}
Let $\mathcal X = (X, \mathfrak S)$ be a measurable space, $\nu$ and $\mu$ be finite and positive measures on $\mathcal X$, and $(\eta_n)_{n\in\mathbb N}$ be a sequence of complex-valued measurable functions defined on $X$ such that
\[
\eta_n\xrightharpoonup[n \to +\infty]{*} f \text{ in } L^\infty(\mathcal X,\mathbb{C},\nu)\quad\text{ and }\quad\eta_n\xrightharpoonup[n \to +\infty]{*} g \text{ in } L^\infty(\mathcal X,\mathbb{C},\mu)
\]
for some measurable functions $f$ and $g$. Suppose in addition that there exist positive measures $\mu_{\mathrm{ac}}$ and $\overline\mu$ such that $\mu=\mu_{\mathrm{ac}}+\overline\mu$, $\mu_{\mathrm{ac}}\ll\nu$, and $\mu_{\mathrm{ac}} \perp \overline\mu$. Then $f=g$ $\mu_{\mathrm{ac}}$-a.e.
\end{lemma}
\begin{proof}
By hypothesis, we have
\[\int_X\eta_n\varphi \diff\mu \xrightarrow[n \to +\infty]{} \int_X g\varphi \diff\mu,\qquad \text{ for all } \varphi\in L^1(\mathcal X,\mathbb{C},\mu),\]
or, equivalently,
\begin{equation}
\label{eq_convergence}
\int_X \eta_n\varphi \diff\mu_{\mathrm{ac}}+\int_X \eta_n\varphi \diff\overline\mu  \xrightarrow[n \to +\infty]{} \int_X g\varphi \diff\mu_{\mathrm{ac}}+\int_X g\varphi \diff\overline\mu, \qquad \text{ for all } \varphi\in L^1(\mathcal X,\mathbb{C},\mu).
\end{equation}
Since $\mu_{\mathrm{ac}} \perp \overline\mu$, there exist $A_1$ and $A_2$ in $\mathfrak S$ such that $A_1 \cap A_2 = \emptyset$, $A_1 \cup A_2 = X$, and $\mu_{\mathrm{ac}}(A_1) = \overline\mu(A_2) = 0$. For any $\varphi\in L^1(\mathcal X,\mathbb{C},\mu_{\mathrm{ac}})$ it is clear that the function $\widetilde\varphi$ defined by $\widetilde\varphi = \varphi \1_{A_2}$ belongs to $L^1(\mathcal X,\mathbb{C},\mu)$. So, from \eqref{eq_convergence}, we obtain, for every $\varphi\in L^1(\mathcal X,\mathbb{C},\mu_{\mathrm{ac}})$,
\[
\int_X \eta_n\varphi \diff\mu_{\mathrm{ac}}=\int_X \eta_n\widetilde\varphi \diff\mu_{\mathrm{ac}}+ \underbrace{\int_X \eta_n\widetilde\varphi \diff\overline\mu}_{=0} \xrightarrow[n \to +\infty]{} \int_X g\widetilde\varphi \diff\mu_{\mathrm{ac}}+\underbrace{\int_X g\widetilde\varphi \diff\overline\mu}_{=0}=\int_X g\varphi \diff\mu_{\mathrm{ac}}.
\]
Thus one has $\eta_n\xrightharpoonup[n \to +\infty]{*}g$ in $L^\infty(\mathcal X,\mathbb{C},\mu_{\mathrm{ac}})$. 

On the other hand, since $\mu_{\mathrm{ac}}\ll \nu$, the Radon--Nikodym derivative $\frac{\df\mu_{\mathrm{ac}}}{\df\nu}$ is well-defined $\nu$-a.e.\ and belongs to $L^1(\mathcal X,\mathbb{C},\nu)$. We set $\psi = \frac{\df\mu_{\mathrm{ac}}}{\df\nu}$. From the convergence
\[ 
\int_X \eta_n\varphi \diff\mu_{\mathrm{ac}}\xrightarrow[n \to +\infty]{} \int_X g\varphi \diff\mu_{\mathrm{ac}},\qquad\text{ for all } \varphi\in L^1(\mathcal X,\mathbb{C},\mu_{\mathrm{ac}}),
\]
we have in particular
\[
\int_X \eta_n\varphi \psi \diff\nu \xrightarrow[n \to +\infty]{} \int_X g\varphi\psi \diff\nu,\qquad\text{ for all } \varphi\in \{\1_A\suchthat A\in \mathfrak S\}.
\]
Since $\varphi\psi \in L^1(\mathcal X,\mathbb{C},\nu)$ for every $\varphi\in \{\1_A\suchthat A\in \mathfrak S \}$ and $\eta_n\xrightharpoonup[n \to +\infty]{*} f \text{ in } L^\infty(\mathcal X,\mathbb{C},\nu)$, we obtain
\[
\int_X \eta_n\varphi \psi \diff\nu \xrightarrow[n \to +\infty]{} \int_X f\varphi\psi \diff\nu,
\]
so, from the uniqueness of the limit,
\[\int_X g\varphi\psi \diff\nu=\int_X f\varphi\psi \diff\nu, \qquad\text{ for all } \varphi\in \{\1_A\suchthat A\in \mathfrak S \}. \]
Therefore,
\begin{equation}\label{eq_lemma_measure1}
\int_X  (f-g) \varphi \psi \diff\nu=\int_X  (f - g) \varphi \diff\mu_{\mathrm{ac}} = 0, \qquad \text{ for all } \varphi\in \{\1_A\suchthat A\in \mathfrak S \}.
\end{equation}
We claim that $f=g$ $\mu_{\mathrm{ac}}$-a.e. Indeed, we set $D=\{\xi \in X\suchthat f(\xi)\ne g(\xi)\}$ and suppose by contradiction that $\mu_{\mathrm{ac}}(D)>0$. Then we have $\mu_{\mathrm{ac}}(D_{\Real})>0$ or $\mu_{\mathrm{ac}}(D_{\Imag})>0$, where $D_{\Real}=\{\xi \in X \suchthat \Real f(\xi)\ne \Real g(\xi)\}$ and $D_{\Imag}=\{\xi \in X \suchthat \Imag f(\xi)\ne \Imag g(\xi)\}$. Without loss of generality, suppose that $\mu_{\mathrm{ac}}(D_{\Real})>0$, so $\mu_{\mathrm{ac}}(D_{\Real}^+)>0$ or $\mu_{\mathrm{ac}}(D_{\Real}^-)>0$, where $D_{\Real}^+=\{\xi \in X \suchthat \Real f(\xi)>\Real g(\xi)\}$ and $D_{\Real}^-=\{\xi \in X \suchthat \Real f(\xi)<\Real g(\xi)\}$. If we assume $\mu_{\mathrm{ac}}(D_{\Real}^+)>0$, then
\[\int_X (\Real f - \Real g) \1_{D_{\Real}^+} \diff\mu_{\mathrm{ac}}>0\]
which contradicts \eqref{eq_lemma_measure1}. The same contradiction happens if we assume $\mu_{\mathrm{ac}}(D_{\Real}^-) > 0$, and therefore the assertion holds.
\end{proof}

We now generalize Lemma~\ref{lemma_measure1} to the case of an arbitrary number of measures.

\begin{proposition}\label{prop_measure2}
Let $\mathcal X = (X, \mathfrak S)$ be a measurable space, $\mu_1,\dotsc,\mu_K$  be finite and positive measures on $\mathcal X$, and $(\eta_n)_{n\in\mathbb N}$ be a sequence of complex-valued measurable functions defined on $X$. Suppose that, for every $k\in \llbracket 1, K\rrbracket$, there exists a measurable function $\overline{\eta}_k$ such that
\[\eta_n\xrightharpoonup[n \to +\infty]{*} \overline{\eta}_k\text{ in } L^\infty(\mathcal X,\mathbb{C},\mu_k).\]
Then there exists a measurable function $\overline{\eta}$ such that $\overline{\eta}=\overline{\eta}_k$ $\mu_k$-a.e., for every $k\in \llbracket 1, K\rrbracket$.
\end{proposition}

\begin{proof}
Let us prove this by induction in $K$. For $K=1$, we just take $\overline{\eta}=\overline{\eta}_1$. Suppose that the result holds for the family $\mu_1, \dotsc, \mu_N$ for some $N \in \llbracket 1, K - 1\rrbracket$. Thus, by the induction hypothesis, there exists a measurable function $\widetilde{\eta}$ such that 
\[\eta_n\xrightharpoonup[n \to +\infty]{*} \widetilde{\eta}\text{ in } L^\infty(\mathcal X,\mathbb{C},\mu_k), \qquad \text{ for all } k\in \llbracket 1, N\rrbracket \]
and 
\[\eta_n\xrightharpoonup[n \to +\infty]{*} \overline{\eta}_{N+1}\text{ in } L^\infty(\mathcal X,\mathbb{C},\mu_{N+1}).\]
Now, we decompose $\mu_{N+1}$ using the Radon--Nikodym--Lebesgue decomposition as follows:
\begin{enumerate}
\item We first write $\mu_{N+1} = \mu_{\mathrm{ac}}^N + \mu_s^N$, where $\mu_{\mathrm{ac}}^N\ll \mu_N$ and $\mu_{s}^N \perp \mu_N$;
\item Recursively, for $j$ equal successively to $N, N-1, \dotsc, 2$, we decompose $\mu_s^j=\mu_{\mathrm{ac}}^{j-1}+\mu_s^{j-1}$, where $\mu_{\mathrm{ac}}^{j-1}\ll \mu_{j-1}$ and $\mu_s^{j-1}\perp \mu_{j-1}$.
\end{enumerate}
Thus, we have
\[\mu_{N+1}=\mu_{\mathrm{ac}}^N+\mu_{\mathrm{ac}}^{N-1}+ \dotsb +\mu_{\mathrm{ac}}^1+\mu_s^1.\]
Now, for $k \in \llbracket 1, N\rrbracket$, we have $\mu_s^k\perp\mu_k$, so we take $\widetilde{A}_k, \widetilde{B}_k$ in $\mathfrak S$ such that $\widetilde{A}_k\cap\widetilde{B}_k=\emptyset$, $\widetilde{A}_k\cup \widetilde{B}_k=X$, and $\mu_{s}^k(\widetilde{A}_k) = \mu_k(\widetilde{B}_k) = 0$. We define the sets $A_j,B_j$ for $j \in \llbracket 1, N\rrbracket$ by $A_j = \widetilde{A}_j\cup\dotsb \cup \widetilde{A}_N $ and $B_j = \widetilde{B}_j\cap\dotsb \cap \widetilde{B}_N$, and by construction we have $A_j \cap B_j = \emptyset$ and $A_j \cup B_j = X$.
Consider $\overline{\eta}$ given by
\[\overline{\eta}(\xi)=\begin{dcases*}
\widetilde{\eta}(\xi),& if $\xi\in A_1$, \\ 
\overline{\eta}_{N+1}(\xi),& if $\xi\in B_1$.
\end{dcases*}\]
To conclude the proof, we will show that $\overline{\eta}=\overline{\eta}_k$ $\mu_k$-a.e.\ for every $k\in \llbracket 1, N+1\rrbracket$. For every $k \in \llbracket 1, N\rrbracket$, note that the set $B=\{\xi \in X\suchthat \overline{\eta}(\xi)\ne \widetilde{\eta}(\xi)\}$ is a subset of $B_1$, thus $\mu_k(B)\leq \mu_k(B_1)\leq \mu_k(\widetilde{B}_k)=0$. Consequently, $\overline{\eta} = \widetilde{\eta} = \overline\eta_k$ $\mu_k$-a.e.\ for every $k \in \llbracket 1, N\rrbracket$. We are thus left to show that $\overline\eta = \overline\eta_{N + 1}$ $\mu_{N+1}$-a.e.

Since $\mu_s^k=\mu_{\mathrm{ac}}^{k-1}+ \dotsb + \mu_{\mathrm{ac}}^1+\mu_s^1$ for $k \in \llbracket 1, N\rrbracket$ and all those measures are positive, then $\mu_s^k(\widetilde{A}_k)=0$ implies $\mu_s^1(\widetilde{A}_k)=0$ for every $k\in \llbracket 1, N\rrbracket$. So $\mu_s^1(A_1)=0$. 

For every $k\in \llbracket 1, N\rrbracket$, we write $\mu_{N+1}=\overline \mu^k+\mu_{\mathrm{ac}}^k$, where $\overline \mu^k=\mu_{\mathrm{ac}}^N+\dotsb+\mu_{\mathrm{ac}}^{k+1}+\mu_{s}^k$, and we set $A^k = \widetilde B_k \cup \widetilde A_{k+1}\cup \dotsb \cup \widetilde A_N$ and $B^k = X \setminus A^k = \widetilde A_k \cap \widetilde B_{k+1} \cap \dotsb \cap \widetilde B_N$. Since $\mu_s^j = \mu_{\mathrm{ac}}^{j-1}+ \dotsb+\mu_{\mathrm{ac}}^1+\mu_s^1$ for $j \in \llbracket 1, N\rrbracket$ and all those measures are positive, we deduce that, if $j \in \llbracket k + 1, N\rrbracket$, then $\mu_{\mathrm{ac}}^k(\widetilde A_j) \leq \mu_s^j(\widetilde A_j) = 0$. In addition, since $\mu_{\mathrm{ac}}^k \ll \mu_k$ and $\mu_k(\widetilde B_k) = 0$, we have $\mu_{\mathrm{ac}}^k(\widetilde B_k) = 0$. Those two facts show that $\mu_{\mathrm{ac}}^k(A^k) = 0$. In addition, for $j \in \llbracket k + 1, N\rrbracket$, we have from the definition of $B^k$ that $\mu_{\mathrm{ac}}^j(B^k) \leq \mu_{\mathrm{ac}}^j(\widetilde B_j) = 0$ and $\mu_{s}^k(B^k) \leq \mu_s^k(\widetilde A_k) = 0$, showing that $\overline \mu^k (B^k) = 0$. Hence $\mu_{\mathrm{ac}}^k \perp \overline\mu^k$ and, applying Lemma~\ref{lemma_measure1} with $f = \widetilde\eta$, $g = \overline\eta_{N + 1}$, $\nu = \mu_k$, and $\mu = \mu_{N+1} = \mu_{\mathrm{ac}}^k + \overline\mu^k$, we have that $\widetilde{\eta}=\overline{\eta}_{N+1}$ $\mu_{\mathrm{ac}}^k$-a.e.\ for every $k \in \llbracket 1, N\rrbracket$. Thus, for $A=\{\xi \in X\suchthat \widetilde{\eta}(\xi)\ne \overline{\eta}_{N+1}(\xi)\}$, one has $\mu_{\mathrm{ac}}^k(A)=0$, for every $k \in \llbracket 1, N\rrbracket$. Letting
\[A^* = \{\xi\in X\suchthat \overline{\eta}(\xi)\ne \overline{\eta}_{N+1}(\xi)\}\]
and noting that $A^* = A_1 \cap A$, we obtain
\[\mu_{N+1}(A^*)\leq \mu_{\mathrm{ac}}^N(A)+ \dotsb+\mu_{\mathrm{ac}}^1(A)+\mu_s^1(A_1)=0.\]
Therefore $\overline\eta=\overline\eta_{N+1}$ $\mu_{N+1}$-a.e., as required.
\end{proof}

\subsection{Other technical results}
\label{app:other}

We next provide an elementary argument for the exponential bound \eqref{eq_exp-estimate} for solutions of \eqref{sys_general}. Recall that such a bound can also be obtained from standard functional-analytic tools, as detailed in \cite[Section~12.3]{Hale1977Theory}.

\begin{proof}[Proof of \eqref{eq_exp-estimate}]
As in the proof of Proposition~\ref{prop_wellposed}, we assume, with no loss of generality, that $\lim_{s \to 0^-} \V_{s}^{0}(M) = \abs{M(0) - M(0^-)} = 0$. Let $x$ be the solution of \eqref{sys_general} with initial condition $\phi$. For $\eta>0$, define $y_\eta \colon [-1, +\infty) \to \mathbb R^d$ by $y_\eta(t) = e^{-\eta t} x(t)$. Then $y_\eta$ satisfies
\[
y_\eta(t) = \int_{-1}^0 \diff M_\eta(\theta) y_\eta(t+\theta),\qquad t \in \mathbb R_+,
\]
where
\[
M_\eta(\theta) = \int_{-1}^{\theta} e^{\eta \tau} \diff M(\tau), \qquad \theta\in [-1,0].
\]
From the definition of $M_\eta$ and Proposition~\ref{proposition:tvcoincides}, it follows that, for every $s \in (-1,0)$,
\[
\V_{-1}^{0}(M_\eta) 
= \V_{-1}^{s}(M_\eta) + \V_{s}^{0}(M_\eta) 
\le e^{\eta s}\V_{-1}^{s}(M) + \V_{s}^{0}(M).
\]
Since $\V_{s}^{0}(M) \to 0$ as $s \to 0^-$, fix $\overline{s} \in (-1, 0)$ so that $\V_{\overline{s}}^{0}(M)< \frac{1}{2}$. Then, for 
$\eta$ sufficiently large,
$e^{\eta \overline{s}}\V_{-1}^{\overline{s}}(M)< \frac{1}{2}$ and hence $\V_{-1}^{0}(M_\eta) < 1$. Applying Proposition~\ref{proposition:tvcoincides} once more, we have, for all $t \in \mathbb R_+$,
\[
\abs{y_\eta(t)} 
= \abs*{ \int_{-1}^0 \diff M_\eta(\theta) y_\eta(t+\theta) }
\le \V_{-1}^{0}(M_\eta) \norm{ (y_\eta)_t }_\infty.
\]
By \cite[Lemma~1]{Mazenc2015Trajectory}, one deduces that there exist positive constants $\widetilde{K}$ and $\beta$ such that
\begin{equation*}
\norm{ (y_\eta)_t }_\infty 
\le \widetilde{K} e^{-\beta t} \norm{ (y_\eta)_0 }_\infty,
\qquad t \ge 0.
\end{equation*}
Since $x(t) = e^{\eta t} y_\eta(t)$, one gets the conclusion.
\end{proof}

Our next result proves that the classes $\NoSingular$ and $\LowerBound$ from \eqref{eq_def-N-L-G} satisfy $\NoSingular \subset \LowerBound$.

\begin{proposition}\label{prop_positiveinf}
For every $M \in \NoSingular$, we have $M \in \LowerBound$.
\end{proposition}

\begin{proof}
We write $M$ as in \eqref{eq_NS}, with $\tau_k \in (0, 1)$ and $A_k\in \mathcal M_d (\mathbb{R})$ for $k \in \mathbb N^\ast$ and with $N$ absolutely continuous. Recall that the characteristic function $\Delta_M$ of $M$ is the function given by
\[
\Delta_M(s) = \det\left(\Id - \sum_{k = 1}^{+\infty} A_k e^{-s \tau_k} - \int_{-1}^0 N^\prime(\theta) e^{s\theta} \diff\theta\right), \qquad s \in \mathbb C.
\]
Since $\NoSingular \subset \WellPosed$, we have $S_M \leq \omega_M < +\infty$. Suppose by contradiction that there exists $\alpha > S_M$ such that $\inf_{\Real s \ge \alpha} \abs{\Delta_M(s)} = 0$. Then there exists a sequence $(s_n)_{n\in \mathbb N}$ such that $\Delta_M(s_n)\to 0$ as $n\to +\infty$ with $\Real s_n\ge \alpha$ for all $n \in \mathbb N$. Since $\Real s_n \ge \alpha > S_M$, it follows that $\abs{s_n} \to +\infty$, otherwise we would obtain a subsequence converging to a zero of $\Delta_M$ located in the half-plane $\{z \in \mathbb{C} \suchthat \Real z \geq \alpha\}$, which yields a contradiction. If $(\Real s_n)_{n \in \mathbb N}$ were unbounded, then, up to a subsequence, we would have $\Real s_n\to +\infty$ as $n \to +\infty$, which would imply, by the dominated convergence theorem, that $\Delta_M(s_n)\to \det(\Id)=1$ as $n \to +\infty$. Thus $(\Real s_n)_{n \in \mathbb N}$ is bounded and converges, up to a subsequence, to some $\beta\ge \alpha$. Since $\abs{s_n} \to +\infty$, we have $\abs{\Imag s_n} \to +\infty$ and, up to replacing by the complex conjugate, we may assume $\Imag s_n\to +\infty$. Setting $x_n = \Real s_n$ and $y_n = \Imag s_n$, we then have $x_n\to \beta$ and $y_n\to +\infty$ as $n \to +\infty$. Define $\Delta_n \colon \mathbb C \to \mathbb C$ by setting, for $s \in \mathbb C$,
\[
\Delta_n(s)=\Delta_M(s+iy_n)=\det \left(\Id - \sum_{k=1}^{+\infty} A_k e^{-s\tau_k}e^{-iy_n\tau_k}-\int_{-1}^0 N^\prime(\theta) e^{s\theta}e^{iy_n\theta} \diff \theta\right).
\]
Up to extracting subsequences through a diagonal procedure, we get that, for every $k \in \mathbb N^\ast$, the sequence $(e^{-iy_n\tau_k})_{n \in \mathbb N}$ converges to $e^{i\theta_k}$ for some angle $\theta_k \in [0, 2\pi)$. By the Riemann--Lebesgue lemma, we have
\[
\int_{-1}^0 N^\prime(\theta)e^{s\theta}e^{iy_n\theta} \diff \theta \xrightarrow[n \to +\infty]{} 0,
\]
uniformly with respect $s$ belonging to any compact subset of the complex plane. As a consequence,
\[
\Delta_n(s) \xrightarrow[n \to +\infty]{} \Delta_{\ast}(s)=\det\left(\Id - \sum_{k=1}^{+\infty} A_k e^{-s\tau_k}e^{i\theta_k}\right),
\]
uniformly with respect to $s$ belonging to a compact neighborhood of $\beta$. Observe that
\[
\abs{\Delta_M(x_n+iy_n)-\Delta_\ast(\beta)}\leq \abs{\Delta_M(x_n+iy_n)-\Delta_M(\beta+iy_n)}+\abs{\Delta_n(\beta)-\Delta_\ast(\beta)}.
\]
Since $\det(\cdot)$ is uniformly continuous on compact sets, we have that $\lim_{n\to +\infty} \abs{\Delta_M(x_n+iy_n)-\Delta_M(\beta+iy_n)} = 0$, so that $\Delta_\ast (\beta)=\lim_{n\to +\infty} \Delta_M(s_n)=0$. By Hurwitz's theorem, for $\delta=\alpha-S_M > 0$, there exists $n \in \mathbb N$ and $s^\prime \in \mathbb{C}$ such that $\abs{s^\prime - \beta} < \delta$ and $\Delta_n(s^\prime) = 0$. Therefore $s^* = s^\prime + i y_n$ is a root of $\Delta_M$ satisfying $\Real s^*= \Real s^\prime > \beta - \delta \ge \alpha - \delta = S_M$, yielding the desired contradiction. 
\end{proof}

The last technical result we present is a generalization of the maximum modulus principle to matrix-valued functions of several variables, replacing the modulus both by a norm and by the spectral radius.

\begin{lemma}\label{lemma_max_principle}
Let $C_1, \dotsc, C_N\subset\mathbb{C}$ be nonempty open sets and $f\colon C_1\times\dotsb\times C_N \to \mathcal{M}_d(\mathbb{C})$ be analytic. For every $j\in \llbracket 1, N\rrbracket$, suppose that $D_j$ is a nonempty, open, bounded, and connected set such that $\overline D_j \subset C_j$. Then
\begin{equation}\label{eq_lemma_max_modulus_principle}
\max_{z\in \overline D_1\times \dotsb\times \overline D_N} \abs{f(z)} = \max_{z\in \partial D_1\times \dotsb\times \partial D_N} \abs{f(z)}.
\end{equation}
Moreover, 
\begin{equation}\label{eq_lemma_max_modulus_principle1}
\max_{z\in \overline D_1\times \dotsb\times \overline D_N}\spr(f(z)) = \max_{z\in \partial D_1\times \dotsb\times \partial D_N}  \spr(f(z)).
\end{equation}
\end{lemma}

\begin{proof}
Let $z^*=(z^*_1, \dotsc,z^*_N)\in \overline D_1\times \dotsb\times \overline D_N$ be such that 
\[
\abs{f(z^*)}=\max_{z\in \overline D_1\times \dotsb\times \overline D_N} \abs{f(z)}.
\]
Let us show that there exists $\overline z= (\overline z_1, \dotsc,\overline z_N)\in \partial D_1\times \dotsb\times \partial D_N$ such that $\abs{f(z^*)}=\abs{f(\overline z)}$. For $j\in \llbracket 1, N\rrbracket$, if $z^*_j\in \partial D_j$ we set $\overline z_j = z^*_j$, otherwise we define $g\colon C_j \to \mathcal{M}_d(\mathbb{C})$ by
\[g(\xi)=f(z^*_1, \dotsc,z^*_{j-1},\xi,z^*_{j+1}, \dotsc,z^*_N).\]
So $g$ is analytic and, from the maximum modulus principle theorem for vector-valued functions (see \cite[Section~III.14, p.~203]{Dunford1988Linear}), we have that $\abs{g(\cdot)}$ is constant on $D_j$, hence for any $\overline z_j\in \partial D_j$, we have $\abs{g(\overline z_j)} = \abs{g(z^*_j)}$. If we apply this argument for any coordinate, we get the desired $\overline z= (\overline z_1, \dotsc,\overline z_N)\in \partial D_1\times \dotsb\times \partial D_N $ and this shows \eqref{eq_lemma_max_modulus_principle}.

Suppose by contradiction that \eqref{eq_lemma_max_modulus_principle1} does not hold. We set
\[\varepsilon=\max_{z\in \overline D_1\times \dotsb\times \overline D_N}\spr(f(z)) -\max_{z\in \partial D_1\times \dotsb\times \partial D_N}  \spr(f(z))>0.\]
From Gelfand's formula, we have  
\[\lim_{n\to +\infty} \abs{A^n}^{\frac{1}{n}}=\spr(A),\]
and the limit is uniform on compact sets of matrices (see, e.g., \cite[Lemma~4.1]{Colonius2019Decay}). Since $\overline D = \overline D_1\times \dotsb \times \overline D_N$ is compact and $f$ is continuous, then $f(\overline D)$ is compact, so there exists $M\in \mathbb N^\ast$ such that
\begin{equation}\label{eq_lemma_max_modulus_principle2}
\abs*{\spr(f(z))-\abs{f(z)^n}^{\frac{1}{n}}} < \frac{\varepsilon}{2}
\end{equation}
for every $n\geq M$ and $z\in \overline D$. Let $z^*\in \overline D$ be such that $\spr(f(z^*)) = \max_{z\in \overline D}\spr (f(z))$ and $z_*\in \partial D = \partial D_1\times \dotsb\times \partial D_N$ such that $\spr(f(z_*)) = \max_{z\in \partial D}\spr (f(z))$. From \eqref{eq_lemma_max_modulus_principle2} and recalling that $\spr(f(z^*))=\spr(f(z_*))+\varepsilon$ we have
\begin{align*}
\abs{f(z^*)^M}^\frac{1}{M}&>\spr(f(z^*))-\frac{\varepsilon}{2}\\
&=\spr(f(z_*))+\frac{\varepsilon}{2}=\max_{z\in \partial D} \spr(f(z))+\frac{\varepsilon}{2}\\
&\geq \max_{z\in \partial D}\abs{f(z)^M}^\frac{1}{M}.
\end{align*}
Therefore $\abs{f(z^*)^M}>\max_{z\in \partial D}\abs{f(z)^M}$, which contradicts \eqref{eq_lemma_max_modulus_principle} applied to $F(z)=f(z)^M$.
\end{proof}
}

\chapter{Stability of difference equations with a time-varying delay}
\label{chap:diff-eqn-varying-delay}

\epigraph{\foreignlanguage{brazil}{Todos os meus movimentos são friamente calculados}}{\foreignlanguage{brazil}{\emph{Chapolin Colorado}}}

The goal of this chapter, which reports the results presented in \cite{MazantiWellPosedness}, is to study in details the effects of a time-varying parameter in the stability of a difference equation of the form \eqref{eq:difference-linear}. As recalled in Section~\ref{sec:intro-varying}, time-varying parameters may have an important effect in the stability of a system, and most common techniques used in the literature provide only sufficient conditions for exponential stability. This chapter aims at providing more precise results in a simpler case, which should be seen as a toy example before treating more general situations. Since the case where the matrices of \eqref{eq:difference-linear} vary in time has already been extensively studied previously in \cite{Chitour2016Stability}, we focus here in the case of time-varying delays and, for simplicity, we consider the case of a single time-varying delay and of a constant matrix.

In other words, we consider here a difference equation of the form
\begin{equation}
\label{eq:diff-eqn-N1-time}
x(t) = A x(t - \tau(t)),
\end{equation}
where the unknown function $x$ takes values in $\mathbb R^d$ for some positive integer $d$, $A \in \mathcal M_d(\mathbb R)$, and $\tau\colon \mathbb R_+ \to \mathbb R_+^\ast$ is a time-dependent delay. Despite its simple-looking form, the analysis of \eqref{eq:diff-eqn-N1-time} is quite subtle, as we detail in this chapter, since even basic questions, such as existence and uniqueness of solutions, have nontrivial answers depending on the functional space being considered. We will focus in this chapter on continuous, regulated, and $L^p$ solutions of \eqref{eq:diff-eqn-N1-time}, which will be prolerly introduced later in Definition~\ref{def:solutions}.

The sequel of this chapter is organized as follows. Section~\ref{sec:prelim} provides some preliminary results on preserving regulatedness or $L^p$ integrability of a function by right composition, which have an independent interest. We then present, in Section~\ref{sec:Hypo-def}, the main hypotheses to be used in this chapter, together with our main definitions. Results on the existence and uniqueness of solutions of \eqref{eq:diff-eqn-N1-time} are provided in Section~\ref{sec:well-posed-varying-delay}, while the analysis of the asymptotic behavior of solutions is given in Section~\ref{sec:stability}. Finally, Section~\ref{sec:appli-varying-delay} provides applications of our results to equations with state-dependent delays and to a transport equation with a time-varying velocity.

\clearpage

\section{Preliminary results}
\label{sec:prelim}

Note that \eqref{eq:diff-eqn-N1-time} can be rewritten as $x(t) = A x(\sigma_1(t))$, where $\sigma_1(t) = t - \tau(t)$. If one wishes to consider solutions $t \mapsto x(t)$ of \eqref{eq:diff-eqn-N1-time} with a certain regularity in $t$, it is important to ensure that such a regularity is preserved by the composition $x \circ \sigma_1$. For instance, a necessary and sufficient condition for $x \circ \sigma_1$ to be continuous for every continuous $x$ is the continuity of $\sigma_1$. This section studies such a question for other types of regularity. More precisely, in Section~\ref{sec:regulated}, we provide the definition of regulated functions and identify the class of functions preserving regulatedness by right composition, while, in Section~\ref{sec:Lp}, we address a similar question for functions in $L^p$.

\subsection{Regulated functions}
\label{sec:regulated}

\begin{definition}
\label{def:regulated}
Let $a, b \in \mathbb R$ with $a < b$ and $f \colon [a, b] \to \mathbb R$.
\begin{enumerate}
\item We say that $f$ is \emph{regulated} if, for every $t_0 \in [a, b)$, the lateral limit $f(t_0^+)$ exists and, for every $t_0 \in (a, b]$, the lateral limit $f(t_0^-)$ exists. The set of regulated functions defined in $[a, b]$ is denoted by $\mathcal G([a, b], \mathbb R)$.

\item\label{item:well-regulated} We say that $f$ is \emph{well-regulated} if $f$ is regulated and, in addition, the following two assertions hold true:
\begin{enumerate}[label={(\roman*)}, ref={\roman*}]
\item\label{item:well-regulated-+} For every $t_0 \in [a, b)$, there exists $\varepsilon \in (0, b - t_0)$ such that one of the following three assertions holds true:
\begin{enumerate}[label={(\theenumii-\arabic*)}, itemindent=2\parindent]
\item\label{item:well-regulated-++} $f(t) > f(t_0^+)$ for every $t \in (t_0, t_0 + \varepsilon)$;
\item\label{item:well-regulated-+0} $f(t) = f(t_0^+)$ for every $t \in (t_0, t_0 + \varepsilon)$;
\item\label{item:well-regulated-+-} $f(t) < f(t_0^+)$ for every $t \in (t_0, t_0 + \varepsilon)$.
\end{enumerate}

\item\label{item:well-regulated--} For every $t_0 \in (a, b]$, there exists $\varepsilon \in (0, t_0 - a)$ such that one of the following three assertions holds true:
\begin{enumerate}[label={(\theenumii-\arabic*)}, itemindent=2\parindent]
\item $f(t) > f(t_0^-)$ for every $t \in (t_0 - \varepsilon, t_0)$;
\item $f(t) = f(t_0^-)$ for every $t \in (t_0 - \varepsilon, t_0)$;
\item $f(t) < f(t_0^-)$ for every $t \in (t_0 - \varepsilon, t_0)$.
\end{enumerate}
\end{enumerate}
The set of well-regulated functions defined in $[a, b]$ is denoted by $\mathcal{WG}([a, b], \mathbb R)$.
\end{enumerate}
\end{definition}

We extend the above definitions to functions defined in an interval of the form $[a, +\infty)$ by saying that $f$ is regulated (respectively, well-regulated) in $[a, +\infty)$ if it is regulated (respectively, well-regulated) in $[a, b]$ for every $b > 0$, and we denote the corresponding functional space by\footnote{In the literature, $\mathcal G([a, +\infty), \mathbb R)$ is usually defined as the set of \emph{bounded} regulated functions in order for it to be a Banach space, but we do not make this requirement here.} $\mathcal G([a, +\infty), \mathbb R)$ (respectively, $\mathcal{WG}([a, +\infty), \mathbb R)$). We also extend the definitions similarly to functions defined in an interval of the form $(-\infty, a]$ or in $\mathbb R$. A vector-valued function is said to be regulated if all of its components are regulated, and the space of regulated functions defined in a nonempty closed interval $I \subset \mathbb R$ and taking values in $\mathbb R^d$ is denoted by $\mathcal G(I, \mathbb R^d)$. Finally, we denote by $\mathcal G(I, A)$ (respectively, $\mathcal{WG}(I, A)$) the set of all regulated (respectively, well-regulated) functions defined in the nonempty closed interval $I$ and taking values in $A \subset \mathbb R^d$ (respectively, $A \subset \mathbb R$).

\begin{remark}
If $I \subset \mathbb R$ is a nonempty compact interval, then $\mathcal G(I, \mathbb R^d)$ is a Banach space when endowed with the norm $\norm{f}_{\mathcal G(I, \mathbb R^d)} = \sup_{t \in I} \abs{f(t)}$. If $I$ is unbounded, the latter function is no longer a norm in $\mathcal G(I, \mathbb R^d)$ since the latter set may contain unbounded functions, but $\mathcal G(I, \mathbb R^d)$ is a Fréchet space when endowed with the family of seminorms $\norm{f}_{\mathcal G(I, \mathbb R^d), n} = \sup_{t \in I \cap [-n, n]} \abs{f(t)}$ for $n \in \mathbb N^\ast$.
\end{remark}

\begin{remark}
\label{remk:cinfty-not-well-regulated}
Well-regulated functions include several classes of functions, such as monotone functions, piecewise monotone functions (with a finite number of changes in monotonicity in every bounded interval), or analytic functions. However, continuous functions, and even $\mathcal C^{\infty}$ functions, may fail to be well-regulated. Indeed, the function $f\colon [-1, 1] \to \mathbb R$ defined by $f(t) = e^{-1/t^2} \sin(1/t)$ for $t \neq 0$ and $f(0) = 0$ belongs to $\mathcal C^{\infty}([-1, 1], \mathbb R)$ but is not well-regulated.
\end{remark}

\begin{remark}
\label{remk:sum-not-well-regulated}
The sum of well-regulated functions may fail to be well-regulated. Indeed, the function $f\colon \mathbb R \to \mathbb R$ defined by $f(0) = 0$ and $f(t) = -t \left(1 - \frac{1}{2} \sin\left(\frac{1}{t}\right)\right)$ for $t \neq 0$ is well-regulated, as well as the identity function $\id$, but $f + \id$ is not well-regulated.
\end{remark}

The main interest of the class of well-regulated functions is the next result, which identifies well-regulated functions as the class of functions preserving regulatedness by right composition.

\begin{theorem}
\label{thm:well-regulated}
Let $I, J$ be two nonempty closed intervals in $\mathbb R$ and $d$ be a positive integer. Then
\begin{equation}
\label{eq:charact-well-regulated}
\mathcal{WG}(J, I) = \{\varphi\colon J \to I \suchthat f \circ \varphi \in \mathcal G(J, \mathbb R^d) \text{ for every } f \in \mathcal G(I, \mathbb R^d)\}.
\end{equation}
\end{theorem}

\begin{proof}
With no loss of generality, we assume in this proof that $d = 1$.

Let $\varphi \in \mathcal{WG}(J, I)$ and $f \in \mathcal G(I, \mathbb R)$. Let $t_0$ be either an interior point of $J$ or its left extremity, if the latter is finite. Since $\varphi$ is well-regulated, $\varphi(t_0^+)$ exists and there exists $\eta > 0$ such that $\varphi$ satisfies \ref{item:well-regulated-++}, \ref{item:well-regulated-+0}, or \ref{item:well-regulated-+-} from Definition~\ref{def:regulated}\ref{item:well-regulated} in the interval $(t_0, t_0 + \eta)$.

Let us consider first the case \ref{item:well-regulated-++}, and assume thus that $\varphi(t) > \varphi(t_0^+)$ for every $t \in (t_0, t_0 + \eta)$. Let $L = \lim_{s \to \varphi(t_0^+)^+} f(s)$, which exists since $f$ is regulated. Take $\varepsilon > 0$. Then there exists $\delta > 0$ such that $\abs{f(s) - L} < \varepsilon$ for every $s \in (\varphi(t_0^+), \varphi(t_0^+) + \delta)$. On the other hand, by definition of $\varphi(t_0^+)$, there exists $\rho > 0$ such that $\abs{\varphi(t) - \varphi(t_0^+)} < \delta$ for every $t \in (t_0, t_0 + \rho)$. Letting $\widetilde\rho = \min\{\rho, \eta\}$, we deduce that, for every $t \in (t_0, t_0 + \widetilde\rho)$, we have $\varphi(t) \in (\varphi(t_0^+), \varphi({t_0^+}) + \delta)$, and hence $\abs{f(\varphi(t)) - L} < \varepsilon$, proving that $\lim_{t \to t_0^+} f(\varphi(t)) = L$. Cases~\ref{item:well-regulated-+0} and \ref{item:well-regulated-+-} are treated by similar arguments, setting $L = f(\varphi(t_0^+))$ and $L = f(\varphi(t_0^+)^-)$, respectively. Hence, $\lim_{t \to t_0^+} f(\varphi(t))$ exists.

One shows by a completely similar argument that $\lim_{t \to t_0^-} f(\varphi(t))$ exists for every $t_0$ in the interior of $J$ or equal to its right extremity if it is finite. Hence, $f \circ \varphi$ is regulated, as required.

Assume now that $\varphi\colon J \to I$ is such that $f \circ \varphi \in \mathcal G(J, \mathbb R)$ for every $f \in \mathcal G(I, \mathbb R)$. Notice that, since the identity function $\id$ is regulated, we have $\varphi = \id \circ \varphi \in \mathcal G(J, \mathbb R)$, i.e., $\varphi$ is regulated.

Assume, to obtain a contradiction, that \eqref{item:well-regulated-+} or \eqref{item:well-regulated--} is false. We only treat the first case, the other case being treated similarly. Then there exists $t_0$ in the interior of $J$ or equal to its left extremity if the latter is finite and a sequence $(t_n)_{n \in \mathbb N^\ast}$ in $J$ with $t_n \to t_0^+$ as $n \to +\infty$ such that at least two of the sets
\begin{align*}
I_+ & = \{n \in \mathbb N^\ast \suchthat \varphi(t_n) > \varphi(t_0^+)\}, \\
I_0 & = \{n \in \mathbb N^\ast \suchthat \varphi(t_n) = \varphi(t_0^+)\}, \\
I_- & = \{n \in \mathbb N^\ast \suchthat \varphi(t_n) < \varphi(t_0^+)\},
\end{align*}
are infinite. Let $f\colon I \to \mathbb R$ be defined for $s \in I$ by $f(s) = \sign(s - \varphi(t_0^+))$, and notice that $f \in \mathcal G(I, \mathbb R)$. We have
\begin{align*}
I_+ & = \{n \in \mathbb N^\ast \suchthat f(\varphi(t_n)) = 1\}, \\
I_0 & = \{n \in \mathbb N^\ast \suchthat f(\varphi(t_n)) = 0\}, \\
I_- & = \{n \in \mathbb N^\ast \suchthat f(\varphi(t_n)) = -1\},
\end{align*}
and since at least two of the above sets are infinite, the sequence $(f(\varphi(t_n)))_{n \in \mathbb N}$ does not converge. Since $t_n \to t_0^+$, this implies that $\lim_{t \to t_0^+} f(\varphi(t))$ does not exist, and hence $f \circ \varphi$ is not regulated, which is a contradiction since $\varphi$ belongs to the set in the right-hand side of \eqref{eq:charact-well-regulated}.
\end{proof}

\begin{remark}
The above proof can be easily adapted to regulated functions taking values in a Banach space $\mathsf X$. Indeed, the proof that any well-regulated function preserves regulatedness by right composition is identical, while, for the proof of the converse statement, one replaces in the argument provided above the identity function and the sign function by the functions $t \mapsto t v$ and $t \mapsto \sign(t) v$, respectively, where $v$ is an arbitrary nonzero element of $\mathsf X$.
\end{remark}

\begin{remark}
This property of well-regulatedness plays important role, since the composition of regulated functions in general is not a regulated function. For instance, consider $f,g \colon [0,1] \to \mathbb R$ given by $f(x)= x\sin(1/x)$ and $g(x) = \sign (x)$. Both of these functions are regulated, however the composition $g \circ f$ has no one-sided limits at $0$. Also, even a
composition of a regulated and continuous functions need not to
be regulated.
\end{remark}

\subsection{Right composition of a function in \texorpdfstring{$L^p$}{Lp}}
\label{sec:Lp}

We are now interested in identifying the class of functions $g$ such that $f \circ g$ belongs to $L^p$ for every $f \in L^p$, i.e., the operator of composition with $g$ on the right preserves the $L^p$ space. As a preliminary step, we provide the following result on the $(\mathfrak L, \mathfrak B)$-measurability of the composition $f \circ g$.

\begin{proposition}
\label{prop:preserves-measurability}
For $i \in \{1, 2\}$, let $X_i \subset \mathbb R^{d_i}$, $m \in \mathbb N^\ast$, and $g\colon X_1 \to X_2$. The following assertions are equivalent.
\begin{enumerate}
\item\label{item:preserves-measurability} For every $f\colon X_2 \to \mathbb R^m$ $(\mathfrak L, \mathfrak B)$-measurable, $f \circ g$ is $(\mathfrak L, \mathfrak B)$-measurable.
\item\label{item:g-is-LL-measurable} The function $g$ is $(\mathfrak L, \mathfrak L)$-measurable.
\end{enumerate}
\end{proposition}

\begin{proof}
The fact that \ref{item:g-is-LL-measurable} implies \ref{item:preserves-measurability} follows immediately from the definition of measurable functions. Conversely, assume that \ref{item:preserves-measurability} holds true and let $A \subset X_2$ be Lebesgue measurable. Let $v \in \mathbb R^m \setminus \{0\}$ and consider the function $f\colon X_2 \to \mathbb R^m$ be given for $x \in X_2$ by $f(x) = \1_A(x) v$. By assumption, $f \circ g$ is $(\mathfrak L, \mathfrak B)$-measurable and thus, in particular, $(f \circ g)^{-1}(\{v\}) = g^{-1}(A)$ is Lebesgue measurable, yielding \ref{item:g-is-LL-measurable}.
\end{proof}

We are now in position to characterize the set of functions preserving $L^p$ by right composition.

\begin{theorem}
\label{thm:composition-Lp}
For $i \in \{1, 2\}$, let $X_i \subset \mathbb R^{d_i}$ be an open set, $m \in \mathbb N^\ast$, and $g\colon X_1 \to X_2$.
\begin{enumerate}[ref=\alph*]
\item\label{item:composition-Linfty-equivs} The following assertions are equivalent.
\begin{enumerate}[label=\textup{(\theenumi-\roman*)}, leftmargin=*]
\item\label{item:composition-Linfty} For every $f \in L^\infty(X_2, \mathbb R^m)$, we have $f \circ g \in L^\infty(X_1, \mathbb R^m)$.
\item\label{item:g-absolutely-continuous} The function $g$ is $(\mathfrak L, \mathfrak L)$-measurable and the measure $g_\# \mathcal L$ is absolutely continuous with respect to $\mathcal L$.
\end{enumerate}
\item\label{item:composition-Lp-equivs} Assume further that $\mathcal L(X_1) < +\infty$ and take $p \in [1, +\infty)$. Then the following assertions are equivalent.
\begin{enumerate}[label=\textup{(\theenumi-\roman*)}, leftmargin=*]
\item\label{item:composition-Lp} For every $f \in L^p(X_2, \mathbb R^m)$, we have $f \circ g \in L^p(X_1, \mathbb R^m)$.
\item\label{item:g-radon-nikodym} The function $g$ is $(\mathfrak L, \mathfrak L)$-measurable, the measure $g_\# \mathcal L$ is absolutely continuous with respect to $\mathcal L$, and the Radon--Nikodym derivative of $g_\# \mathcal L$ with respect to $\mathcal L$ is bounded.
\end{enumerate}
\end{enumerate}
\end{theorem}

\begin{proof}
Let us first prove \eqref{item:composition-Linfty-equivs}. Assume that \ref{item:g-absolutely-continuous} holds true and let $f \in L^\infty(X_2, \mathbb R^m)$. Recall that, by Proposition~\ref{prop:preserves-measurability}, $f \circ g$ is $(\mathfrak L, \mathfrak B)$-measurable. In addition, since $f$ is bounded, so is $f \circ g$. One is only left to show that $f \circ g$ defines an element in $L^\infty(X_1, \mathbb R^m)$, i.e., that, if $\widetilde f = f$ a.e.\ in $X_2$, then $\widetilde f \circ g = f \circ g$ a.e.\ in $X_1$. Let $N = \{x \in X_2 \suchthat \widetilde f(x) \neq f(x)\}$ and notice that $\mathcal L(N) = 0$. In addition, $\widetilde f \circ g(x) \neq f \circ g(x)$ if and only if $x \in g^{-1}(N)$. Since $g_\# \mathcal L$ is absolutely continuous with respect to $\mathcal L$, we have $\mathcal L(g^{-1}(N)) = 0$, showing that $\widetilde f \circ g = f \circ g$ a.e.\ in $X_1$ and yielding the conclusion.

Assume now that \ref{item:composition-Linfty} holds true. In particular, by Proposition~\ref{prop:preserves-measurability}, $g$ is $(\mathfrak L, \mathfrak L)$-mea\-su\-ra\-ble. Let $N \subset X_2$ be a Lebesgue-measurable set satisfying $\mathcal L(N) = 0$, $v \in \mathbb R^m \setminus \{0\}$, and $f\colon X_2 \to \mathbb R^m$ be given for $x \in X_2$ by $f(x) = \1_N(x) v$. Then $f \in L^\infty(X_2, \mathbb R^m)$, $f = 0$ a.e.\ in $X_2$, and thus $f \circ g = 0$ a.e.\ in $X_1$, showing that $\mathcal L(g^{-1}(N)) = \mathcal L(\{x \in X_1 \suchthat f \circ g(x) = v\}) = 0$. This proves that $g_\# \mathcal L$ is absolutely continuous with respect to $\mathcal L$.

Let us now prove \eqref{item:composition-Lp-equivs}. Assume that \ref{item:g-radon-nikodym} holds true and let $f \in L^p(X_2, \mathbb R^m)$. As before, the $(\mathfrak L, \mathfrak B)$-measurability of $f \circ g$ follows from Proposition~\ref{prop:preserves-measurability} and, arguing as in the proof of \eqref{item:composition-Linfty-equivs}, we obtain that, if $\widetilde f = f$ a.e.\ in $X_2$, then $\widetilde f \circ g = f \circ g$ a.e.\ in $X_1$. In addition, denoting by $\varphi\colon X_2 \to \mathbb R^m$ the Radon--Nikodym derivative of $g_\# \mathcal L$ with respect to $\mathcal L$, we have
\[
\int_{X_1} \abs{f \circ g(x)}^p \diff x = \int_{X_2} \abs{f(x)}^p \diff g_\# \mathcal L(x) = \int_{X_2} \abs{f(x)}^p \varphi(x) \diff x < +\infty,
\]
since $\varphi$ is bounded. Hence $f \circ g \in L^p(X_1, \mathbb R^m)$.

Assume now that \ref{item:composition-Lp} holds true. By Proposition~\ref{prop:preserves-measurability}, $g$ is $(\mathfrak L, \mathfrak L)$-measurable and, arguing as in the proof of \eqref{item:composition-Linfty-equivs}, $g_\# \mathcal L$ is absolutely continuous with respect to $\mathcal L$. Since $\mathcal L(X_1) < +\infty$, the measure $g_{\#} \mathcal L$ is finite, and hence its Radon--Nikodym derivative $\varphi\colon X_2 \to \mathbb R_+ \cup \{+\infty\}$ with respect to $\mathcal L$ exists.

Let $L\colon L^1(X_2, \mathbb R) \to \mathbb R$ be the linear map defined for $\psi \in L^1(X_2, \mathbb R)$ by
\[
L \psi = \int_{X_2} \psi(x) \varphi(x) \diff x.
\]
Let us prove first that $L$ is well-defined. Indeed, write $\psi = \psi_+ - \psi_-$, where $\psi_+ = \max(\psi, 0)$ and $\psi_- = \max(-\psi, 0)$, and note that $\psi_+$ and $\psi_-$ are nonnegative functions in $L^1(X_2, \mathbb R)$. Let $f_+$ and $f_-$ be functions in $L^p(X_2, \mathbb R^m)$ such that $\psi_+ = \abs{f_+}^p$ and $\psi_- = \abs{f_-}^p$. By \ref{item:composition-Lp}, both $f_+ \circ g$ and $f_- \circ g$ belong to $L^p(X_1, \mathbb R^m)$, and we compute
\[
\int_{X_2} \psi_{\pm}(x) \varphi(x) \diff x = \int_{X_2} \abs{f_{\pm}(x)}^p \diff g_\# \mathcal L(x) = \int_{X_1} \abs{f_{\pm} \circ g(x)}^p \diff x < +\infty.
\]
Hence, $\psi \varphi \in L^1(X_2, \mathbb R)$, and thus $L \psi$ is well-defined.

For $n \in \mathbb N^\ast$, let $\varphi_n = \min(\varphi, n)$. Let $L_n\colon L^1(X_2, \mathbb R) \to \mathbb R$ be the continuous linear operator defined for $\psi \in L^1(X_2, \mathbb R)$ by $L_n \psi = \int_{X_2} \psi(x) \varphi_n(x) \diff x$. For a given $\psi \in L^1(X_2, \mathbb R)$, we have
\begin{align*}
L_n \psi & = \int_{X_2} \psi(x) \varphi_n(x) \diff x \\
& = \int_{X_2} \psi_+(x) \varphi_n(x) \diff x - \int_{X_2} \psi_-(x) \varphi_n(x) \diff x \\
& \xrightarrow[n \to +\infty]{} \int_{X_2} \psi_+(x) \varphi(x) \diff x - \int_{X_2} \psi_-(x) \varphi(x) \diff x \\
& = \int_{X_2} \psi(x) \varphi(x) \diff x = L \psi,
\end{align*}
where the convergence holds true thanks to the Lebesgue monotone convergence theorem. Hence, by the Banach--Steinhaus theorem, $L$ is a continuous linear operator, which implies that $L$ belongs to the topological dual of $L^1(X_2, \mathbb R)$, and thus $\varphi \in L^\infty(X_2, \mathbb R)$, as required.
\end{proof}

\section{Hypotheses and definitions}
\label{sec:Hypo-def}

Let us now present and comment the main definitions and hypotheses used in the sequel of the chapter. We start with the definition of delay function.

\begin{definition}
A function $\tau\colon \mathbb R_+ \to \mathbb R_+^\ast$ is said to be a \emph{delay function}.
\end{definition}

The study of existence and uniqueness of solutions of \eqref{eq:diff-eqn-N1-time} and their asymptotic behavior requires several assumptions. For simplicity, we state all such assumptions here, and refer to them whenever needed in the sequel.

\begin{hypotheses}
\label{hypo:hypos}
Let $A \in \mathcal M_d(\mathbb R^d)$, $p \in [1, +\infty]$, and $\tau$ be a delay function and define the function $\sigma_1\colon \mathbb R_+ \to \mathbb R$ by $\sigma_1(t) = t - \tau(t)$.
\begin{hypolist}
\item\label{hypo:tau-infimum} For every $t \geq 0$, we have $\inf_{s \in [0, t]} \tau(s) > 0$.
\item\label{hypo:tau-cont} The function $\tau$ is continuous.
\item\label{hypo:regulated} The function $\sigma_1$ is well-regulated.
\item\label{hypo:tau-measurable} The function $\tau$ is $(\mathfrak L, \mathfrak L)$-measurable.
\item\label{hypo:null} We have $\mathcal L(\sigma_1^{-1}(N)) = 0$ for every $N \in \mathfrak L$ such that $\mathcal L(N) = 0$.
\item\label{hypo:radon-nik-bdd} Hypotheses~\ref{hypo:tau-measurable} and \ref{hypo:null} hold true and, for every $T > 0$, the Radon--Nikodym derivative of the measure $\left(\sigma_1\rvert_{[0, T]}\right)_{\#}\mathcal L$ with respect to $\mathcal L$ is bounded.
\item\label{hypo:spr-A-less-1} We have $\spr(A) < 1$.
\item\label{hypo:delay-infinity} We have $\lim_{t \to +\infty} \sigma_1(t) = +\infty$.
\item\label{hypo:radon-nik-and-rho} Hypotheses~\ref{hypo:tau-measurable}, \ref{hypo:null}, and \ref{hypo:delay-infinity} hold true, $p < +\infty$, and the Radon--Nikodym derivative $\varphi$ of ${\sigma_1}_{\#}\mathcal L$ with respect to $\mathcal L$ is bounded and satisfies $\norm{\varphi}_{L^\infty(\mathbb R, \mathbb R)} \spr(A)^p < 1$.
\item\label{hypo:delay-bounded} The delay function $\tau$ is bounded.
\end{hypolist}
\end{hypotheses}

Assumption~\ref{hypo:tau-infimum} is used for obtaining the representation formula \eqref{eq:explicit-formula} given in Theorem~\ref{theorem-1}, and it is important for the uniqueness of solutions of \eqref{eq:diff-eqn-N1-time} (see Remark~\ref{remk:H2}). The continuity of $\tau$ assumed in \ref{hypo:tau-cont} is important if one wants to ensure the existence of continuous solutions of \eqref{eq:diff-eqn-N1-time} (see Theorem~\ref{thm:exist-continuous}), while, in view of Theorem~\ref{thm:well-regulated}, \ref{hypo:regulated} is used to ensure existence of solutions of \eqref{eq:diff-eqn-N1-time} in the class of regulated functions (see Theorem~\ref{thm:exist-regulated}).

In view of Theorem~\ref{prop:preserves-measurability}, the measurability of $\tau$ assumed in \ref{hypo:tau-measurable} is important to obtain existence of measurable solutions of \eqref{eq:diff-eqn-N1-time} (see Theorem~\ref{thm:exist-measurable}). Assumption~\ref{hypo:null} states that $\sigma_1$ preserves the class of sets of Lebesgue measure zero, and it turns out to be a necessary assumption to provide a notion of solution verifying \eqref{eq:diff-eqn-N1-time} only for almost every $t \in \mathbb R_+$ (see Theorem~\ref{thm:exist-ae} and Proposition~\ref{prop:h2-necessary}). This is fundamental in order to consider solutions of \eqref{eq:diff-eqn-N1-time} in $L^p$ spaces. Note that \ref{hypo:tau-measurable} ensures that ${\sigma_1}_{\#} \mathcal L$ is a measure in $(\mathbb R, \mathfrak L)$ and, under this assumption, \ref{hypo:null} is equivalent to requiring ${\sigma_1}_{\#} \mathcal L$ to be absolutely continuous with respect to the Lebesgue measure, which, by Theorem~\ref{thm:composition-Lp}\eqref{item:composition-Linfty-equivs}, is fundamental for ensuring that any $L^\infty$ function remains in $L^\infty$ after a right composition with $\sigma_1$. In addition, \ref{hypo:tau-measurable} and \ref{hypo:null} ensure the existence of the Radon--Nikodym derivative in \ref{hypo:radon-nik-bdd}, since $\left(\sigma_1\rvert_{[0, T]}\right)_{\#}\mathcal L$ is a finite measure. Finally, taking into account Theorem~\ref{thm:composition-Lp}, \ref{hypo:radon-nik-bdd} allows one to prove existence of solutions of \eqref{eq:diff-eqn-N1-time} in $L^p$ for $p \in [1, +\infty)$.

Note that \ref{hypo:tau-cont} is equivalent to assuming continuity of $\sigma_1$, and \ref{hypo:tau-measurable} is equivalent to assuming $(\mathfrak L, \mathfrak L)$-measurability of $\sigma_1$, but, due to Remark~\ref{remk:sum-not-well-regulated}, \ref{hypo:regulated} is not equivalent to assuming well-regulatedness of $\tau$. 

Assumptions~\ref{hypo:spr-A-less-1}--\ref{hypo:delay-bounded} are useful for the analysis of the asymptotic behavior of solutions of \eqref{eq:diff-eqn-N1-time} carried out in Section~\ref{sec:stability}. Indeed, in order to obtain that solutions of \eqref{eq:diff-eqn-N1-time} converge to $0$ as time tends to infinity, it seems natural from the form of \eqref{eq:diff-eqn-N1-time} to require \ref{hypo:spr-A-less-1}, since one may expect the existence of solutions diverging to infinity in norm in the case $\spr(A) > 1$, and not converging to $0$ in the case $\spr(A) = 1$. In addition, if $\sigma_1(t) = t - \tau(t)$ were to remain bounded as $t \to +\infty$, then, from \eqref{eq:diff-eqn-N1-time}, $x(t)$ would depend, as $t \to +\infty$, only on the past values of $x(\cdot)$ computed in a bounded time interval, and one might expect solutions of \eqref{eq:diff-eqn-N1-time} to not necessarily converge to $0$ in this case. These intuitions will be made precise later on in Theorem~\ref{thm:cv-to-0} and Proposition~\ref{prop:necessary-asympt}, which show that, under \ref{hypo:tau-infimum}, assumptions~\ref{hypo:spr-A-less-1} and \ref{hypo:delay-infinity} are necessary\footnote{The necessity of \ref{hypo:delay-infinity} requires also the assumption that $A$ is not nilpotent.} and sufficient to have $\abs{x(t)} \to 0$ as $t \to +\infty$ for any solution $x$ of \eqref{eq:diff-eqn-N1-time} with bounded initial condition.

In addition to the convergence of $\abs{x(t)}$ to $0$ as time tends to infinity, we are also interested in Section~\ref{sec:stability} in the convergence to $0$ of functional norms of $x(\cdot)$ computed at time intervals tending to infinity. While \ref{hypo:spr-A-less-1} and \ref{hypo:delay-infinity} remain sufficient for the $L^\infty$ norm, we will show (Example~\ref{example-3.1}) that this is not the case for $L^p$ norms with $p \in [1, +\infty)$, and it turns out that requiring also \ref{hypo:radon-nik-and-rho} to hold true is sufficient to retrieve convergence to $0$ also in the $L^p$. Note that, in \ref{hypo:radon-nik-and-rho}, the measure ${\sigma_1}_{\#}\mathcal L$ is necessarily $\sigma$-finite thanks to \ref{hypo:delay-infinity}, and thus the Radon--Nikodym derivative $\varphi$ exists.

Finally, assumption~\ref{hypo:delay-bounded} is useful in Section~\ref{sec:exponential} to obtain exponential decay rates of solutions of \eqref{eq:diff-eqn-N1-time}. Note that, trivially, \ref{hypo:delay-bounded} implies \ref{hypo:delay-infinity}.

\begin{remark}
A sufficient condition for \ref{hypo:radon-nik-bdd} to hold is that, in addition to \ref{hypo:tau-measurable} and \ref{hypo:null}, ${\sigma_1}_{\#} \mathcal L$ is $\sigma$-finite and its Radon--Nikodym derivative $\varphi$ with respect to $\mathcal L$ is bounded on $(-\infty, T]$ for every $T > 0$. Indeed, using the fact that $\sigma_1(t) < t$ for every $t \in \mathbb R_+$, we deduce that, for every $T > 0$, the Radon--Nikodym derivative of $\left(\sigma_1\rvert_{[0, T]}\right)_{\#} \mathcal L$ is supported on $(-\infty, T]$ and upper bounded by $\varphi$ on this set.

However, the previous condition is not necessary for \ref{hypo:radon-nik-bdd} to hold true: for instance, if $\tau(t) = \floor{t + 1}$ for $t \in \mathbb R_+$, we have $\sigma_1(t) = t - \floor{t + 1}$ for $t \in \mathbb R_+$, and we can verify that \ref{hypo:radon-nik-bdd} is satisfied, but ${\sigma_1}_{\#}\mathcal L$ is not $\sigma$-finite, and hence it does not admit a Radon--Nikodym derivative with respect to $\mathcal L$.
\end{remark}

\begin{remark}
\label{remk:radon-nikodym-derivative-tau}
If $\sigma_1$ is defined from a delay function $\tau$ as in Hypotheses~\ref{hypo:hypos} and it is a diffeomorphism between $\mathbb R_+$ and an interval $I \subset \mathbb R$, then ${\sigma_1}_{\#} \mathcal L$ is absolutely continuous with respect to $\mathcal L$, it is $\sigma$-finite, and its Radon--Nikodym derivative with respect to $\mathcal L$ is the function $\varphi\colon \mathbb R \to \mathbb R_+$ defined by $\varphi(t) = \abs{(\sigma_1^{-1})^\prime(t)} = \frac{1}{\abs{\sigma_1^\prime(\sigma_1^{-1}(t))}}$ if $t \in I$ and $\varphi(t) = 0$ otherwise. Given $T > 0$, $\varphi$ is bounded on $(-\infty, T]$ if and only if $\abs{\sigma_1^\prime}$ is bounded from below on $(-\infty, T]$ by a positive constant.

In terms of $\tau$, a sufficient condition for \ref{hypo:radon-nik-bdd} to hold true is that $\tau$ is differentiable and, for every $T > 0$ there exists $\alpha \in [0, 1)$ such that $\tau^\prime(t) \leq \alpha$ for every $t \in [0, T]$. In that case, $\varphi$ is bounded on $(-\infty, T]$ by $\frac{1}{1 - \alpha}$.
\end{remark}

We now introduce several definitions used in the sequel of the chapter.

\begin{definition}
\label{def:h}
Let $\tau$ be a delay function. The \emph{largest delay function} associated with $\tau$ is the function $h\colon \mathbb R_+ \to \mathbb R_+^\ast \cup \{+\infty\}$ defined for $t \in \mathbb R_+$ by
\[
h(t) = t - \inf_{s \in [t, +\infty)} (s - \tau(s)).
\]
\end{definition}

The largest delay function $h$ introduced in Definition~\ref{def:h} can be interpreted as follows: if a function $x$ satisfies \eqref{eq:diff-eqn-N1-time} for every $t \geq t_0$, for some $t_0 \in \mathbb R$, then the right-hand side of \eqref{eq:diff-eqn-N1-time} depends on evaluations of $x$ at the interval $\lvert t_0 - h(t_0), +\infty)$. In particular, if one is to define a solution of \eqref{eq:diff-eqn-N1-time} in an interval of the form $[t_0, +\infty)$, then the corresponding initial condition should be a function defined in an interval of the form $\lvert t_0 - h(t_0), t_0)$. Notice also that, from the definition of $h$, the function $t_0 \mapsto t_0 - h(t_0)$ is nondecreasing.

\begin{definition}
Let $\tau$ be a delay function and $h$ be its associated largest delay function.
\begin{enumerate}
\item Given $x\colon \lvert -h(0), +\infty) \to \mathbb R^d$ and $t \geq 0$, we denote by $x_t\colon \lvert -h(t), 0) \to \mathbb R^d$ the history function of $x$ at $t$, defined by $x_t(s) = x(t + s)$ for every $s \in \lvert -h(0), 0)$. Whenever needed, we also make the slight abuse of notation of considering $x_t$ to be defined on $\lvert -h(t), 0]$ by setting $x_t(0) = x(t)$.

\item For every $t \in \mathbb R_+$, we set $\mathsf C_t = \{x \in \mathcal C(\lvert -h(t), 0], \mathbb R^d) \suchthat x(0) = A x(-\tau(t))\}$ and $\mathsf C_t^{\mathrm b} = \{x \in \mathsf C_t \suchthat x \text{ is bounded}\}$.

\item For every $t \in \mathbb R_+$, we set $\mathsf G_t = \{x \in \mathcal G(\lvert -h(t), 0], \mathbb R^d) \suchthat x(0) = A x(-\tau(t))\}$ and $\mathsf G_t^{\mathrm b} = \{x \in \mathsf G_t \suchthat x \text{ is bounded}\}$.
\end{enumerate}
\end{definition}

Note that the sets $\mathsf C_t$ and $\mathsf G_t$ depend on $\tau$ and $A$, but we do not make this dependence explicit in their notation. We also remark that the interval $\lvert -h(t), 0]$ is not compact when $h(t) = +\infty$, motivating the boundedness requirements in the definitions of $\mathsf C_t^{\mathrm b}$ and $\mathsf G_t^{\mathrm b}$. In particular, $\mathsf C_t^{\mathrm b}$ and $\mathsf G_t^{\mathrm b}$ are Banach spaces with norms
\[
\norm{x}_{\mathsf C_t^{\mathrm b}} = \sup_{t \in \lvert -h(t), 0]} \abs{x(t)} \quad \text{ and } \quad \norm{x}_{\mathsf G_t^{\mathrm b}} = \sup_{t \in \lvert -h(t), 0]} \abs{x(t)},
\]
respectively.

\section{Existence and uniqueness of solutions}
\label{sec:well-posed-varying-delay}

We next address questions on the existence and uniqueness of solutions of \eqref{eq:diff-eqn-N1-time}, and we start by introducing the three main notions of solutions used along this chapter, namely continuous, regulated, and $L^p$ solutions, as well as some other auxiliary notions.

\begin{definition}
\label{def:solutions}
Let $\tau$ be a delay function and $h$ be its associated maximal delay function.
\begin{enumerate}
\item We say that $x\colon \lvert -h(0), +\infty) \to \mathbb R^d$ is a \emph{solution} of \eqref{eq:diff-eqn-N1-time} if the equality $x(t) = A x(t - \tau(t))$ is satisfied for every $t \in \mathbb R_+$.

\item A solution $x$ of \eqref{eq:diff-eqn-N1-time} is said to be a \emph{continuous} (respectively, \emph{regulated}) \emph{solution} of \eqref{eq:diff-eqn-N1-time} if it is continuous (respectively, regulated).

\item A solution $x$ of \eqref{eq:diff-eqn-N1-time} is said to be a \emph{measurable solution} of \eqref{eq:diff-eqn-N1-time} if it is $(\mathfrak L, \mathfrak B)$-mea\-su\-rable.

\item We say that $x\colon\lvert -h(0), +\infty) \to \mathbb R^d$ is an \emph{almost everywhere (a.e.)\ solution} of \eqref{eq:diff-eqn-N1-time} if, for every function $\widetilde x\colon \lvert -h(0), +\infty) \to \mathbb R^d$ equal to $x$ almost everywhere on $\lvert -h(0), +\infty)$, we have $\widetilde x(t) = A \widetilde x(t - \tau(t))$ for almost every $t \in \mathbb R_+$.

\item For $p \in [1, +\infty]$, a measurable a.e.\ solution $x$ of \eqref{eq:diff-eqn-N1-time} is said to be an \emph{$L^p$ solution} of \eqref{eq:diff-eqn-N1-time} if $x \in L^p_{\loc}(\lvert -h(0), +\infty), \mathbb R^d)$.
\end{enumerate}

If $x$ is a solution, a measurable solution, an a.e.\ solution, or an $L^p$ solution of \eqref{eq:diff-eqn-N1-time}, we say that its \emph{initial condition} is its restriction to the interval $\lvert -h(0), 0)$. If $x$ is a continuous or regulated solution of \eqref{eq:diff-eqn-N1-time}, we say that its \emph{initial condition} is its restriction to the interval $\lvert -h(0), 0]$.
\end{definition}

\subsection{Representation formula for solutions}
\label{sec:explicit-formula}

The goal of this section is to provide a representation formula for solutions of \eqref{eq:diff-eqn-N1-time}, expressing a solution $x$ at a given time $t$ in terms of its initial condition $x_0$ evaluated at well-identified times and suitable powers of the matrix $A$. The formula we aim at obtaining is similar in spirit to that provided in Proposition~\ref{prop_var_const_formula} for \eqref{eq:intro-linear-difference-equation}, obtained by iterating the difference equation until one reaches the initial condition. The major subtleties with respect to Proposition~\ref{prop_var_const_formula} come from the fact that the delay is no longer constant and, to deal with this setting, we start by a suitable definition of iterates of the function $t \mapsto t - \tau(t)$.

\begin{definition}
\label{def:Dn-sigman}
Let $\tau$ be a delay function. For $n \in \mathbb N$, we define the \emph{iterated delay functions} $\sigma_n\colon D_n \to \mathbb R$ inductively by setting $D_0 = \mathbb R$ and $\sigma_0(t) = t$ for $t \in \mathbb R$ and defining, for $n \in \mathbb N^\ast$,
\begin{equation}
\label{eq:defi-Dn}
D_n = \{t \in \mathbb R_+ \suchthat t - \tau(t) \in D_{n-1}\}
\end{equation}
and
\begin{equation*}
\sigma_n(t) = \sigma_{n-1}(t - \tau(t)) \qquad \text{ for } t \in D_n.
\end{equation*}
\end{definition}

Note that $D_1 = \mathbb R_+$ and $\sigma_1(t) = t - \tau(t)$. Moreover, we may have $D_n = \emptyset$ for some $n \in \mathbb N$: this is the case, for instance, if the delay function $\tau$ is defined for $t \in \mathbb R_+$ by $\tau(t) = t + 1$, in which case $\sigma_1(t) = -1$ for every $t \in \mathbb R_+$ and hence $D_2 = \emptyset$. The empty set, however, is not necessarily attained by some $D_n$: for instance, if $\tau(t) = 1$ for every $t \in \mathbb R_+$, one immediately computes that $D_n = [n-1, +\infty)$ for every $n \in \mathbb N^\ast$.

We now provide immediate properties of iterated delay functions. The first one is an immediate consequence of \eqref{eq:defi-Dn}.

\begin{lemma}
\label{lemma:D_n-sigma1}
Let $\tau$ be a delay function, $t \in \mathbb R_+$, and $n \in \mathbb N$. Then $t \in D_{n+1}$ if and only if $\sigma_1(t) \in D_n$.
\end{lemma}

The next result shows that the sequence of sets $(D_n)_{n \in \mathbb N}$ is nonincreasing.

\begin{lemma}
\label{lemma:Dn-decreasing}
Let $\tau$ be a delay function. Then, for every $n \in \mathbb N^\ast$, we have $D_n \subset D_{n-1}$.
\end{lemma}

\begin{proof}
We prove the result by induction on $n$. The case $n = 1$ is trivially true. Now, let $n \in \mathbb N^\ast$ be such that $D_n \subset D_{n-1}$. Let $t \in D_{n+1}$ and note that, by definition, $t \in \mathbb R_+$ and $t - \tau(t) \in D_n$, and hence $t - \tau(t) \in D_{n-1}$, implying that $t \in D_n$. Thus $D_{n+1} \subset D_n$, yielding the conclusion.
\end{proof}

\begin{lemma}
\label{lemma:caract-Dn-sigman}
Let $\tau$ be a delay function. Then, for every $n \in \mathbb N$ and $k \in \llbracket 0, n\rrbracket$, we have
\begin{equation}
\label{eq:rec-Dn}
D_n = \{t \in D_k \suchthat \sigma_k(t) \in D_{n-k}\}
\end{equation}
and
\begin{equation}
\label{eq:rec-sigman}
\sigma_n(t) = \sigma_{n-k}(\sigma_k(t)) \qquad \text{ for every } t \in D_n.
\end{equation}
\end{lemma}

\begin{proof}
We prove the statement by strong induction on $n$. Notice first that the result is immediate for $n = 0$ and $n = 1$. Let $N \in \mathbb N^\ast$ be such that \eqref{eq:rec-Dn} and \eqref{eq:rec-sigman} hold for every $n \in \llbracket 0, N-1\rrbracket$ and $k \in \llbracket 0, n\rrbracket$. We notice that \eqref{eq:rec-Dn} and \eqref{eq:rec-sigman} trivially hold with $n$ replaced by $N$ and $k$ replaced by either $0$ or $N$. We then fix $k \in \llbracket 1, N-1\rrbracket$.

By \eqref{eq:defi-Dn}, we have
\[
D_{N} = \{t \in D_1 \suchthat \sigma_1(t) \in D_{N-1}\}.
\]
Using the induction assumption \eqref{eq:rec-Dn} with $n$ replaced by $N-1$ and $k$ replaced by $k-1$, we deduce that $\sigma_1(t) \in D_{N-1}$ if and only if $\sigma_1(t) \in D_{k-1}$ and $\sigma_{k-1}(\sigma_1(t)) \in D_{N - k}$, i.e.,
\[
D_{N} = \{t \in D_1 \suchthat \sigma_1(t) \in D_{k-1} \text{ and } \sigma_{k-1}(\sigma_1(t)) \in D_{N - k}\}.
\]
By \eqref{eq:defi-Dn}, we have $t \in D_k$ if and only if $t \in D_1$ and $\sigma_1(t) \in D_{k-1}$, and thus
\[
D_{N} = \{t \in D_k \suchthat \sigma_{k-1}(\sigma_1(t)) \in D_{N - k}\}.
\]
Since $k \in \llbracket 1, N-1\rrbracket$, the induction assumption \eqref{eq:rec-sigman} also shows that $\sigma_{k-1}(\sigma_1(t))$ for every $t \in D_k$, yielding
\[
D_{N} = \{t \in D_k \suchthat \sigma_{k}(t) \in D_{N - k}\},
\]
as required.

Finally, to obtain \eqref{eq:rec-sigman} with $n$ replaced by $N$, notice that, by definition, for every $t \in D_N$, we have $\sigma_N(t) = \sigma_{N-1}(\sigma_1(t))$. On the other hand, by the induction assumption \eqref{eq:rec-sigman} with $n$ replaced by $N-1$ and $k$ replaced by $k-1$, we deduce that $\sigma_{N-1}(s) = \sigma_{N-k}(\sigma_{k-1}(s))$ for every $s \in D_{N-1}$. Since $\sigma_1(t) \in D_{N-1}$ for every $t \in D_N$, we finally deduce that, for every $t \in D_N$, we have
\[
\sigma_{N}(t) = \sigma_{N-1}(\sigma_1(t)) = \sigma_{N-k}(\sigma_{k-1}(\sigma_1(t))) = \sigma_{N-k}(\sigma_k(t)),
\]
as required.
\end{proof}

\begin{lemma}
\label{lemma:intersect-Dn-empty}
Let $\tau$ be a delay function satisfying \ref{hypo:tau-infimum}. Then
\[\bigcap_{n \in \mathbb N} D_n = \emptyset.\]
\end{lemma}

\begin{proof}
Assume, to obtain a contradiction, that there exists $t \in \bigcap_{n \in \mathbb N} D_n$, and consider the sequence $(t_k)_{k \in \mathbb N}$ defined for $k \in \mathbb N$ by $t_k = \sigma_k(t)$. By Lemma~\ref{lemma:caract-Dn-sigman}, we deduce that $t_k \in \bigcap_{n \in \mathbb N} D_n$ for every $k \in \mathbb N$. In particular, $t_k \geq 0$ for every $k \in \mathbb N$. Since $t_k = \sigma_1(\sigma_{k-1}(t)) = t_{k-1} - \tau(t_{k-1})$ and $\tau$ takes values in $\mathbb R_+^\ast$, we also have that $(t_k)_{k \in \mathbb N}$ is a decreasing sequence.

Let $\delta = \inf_{s \in [0, t]} \tau(s) > 0$. Since $t_k \in [0, t]$ for every $k \in \mathbb N$, we deduce that $\tau(t_{k-1}) \geq \delta$ for every $k \in \mathbb N^\ast$, and thus $t_k \leq t_{k-1} - \delta$, which yields $t_k \leq t - k \delta$ for every $k \in \mathbb N$. This, however, contradicts the fact that $t_k \geq 0$ for every $k \in \mathbb N$, establishing the result.
\end{proof}

As an immediate consequence of Lemmas~\ref{lemma:Dn-decreasing} and \ref{lemma:intersect-Dn-empty}, we obtain the following result.

\begin{corollary}
\label{coro:integer-iterations}
Let $\tau$ be a delay function satisfying \ref{hypo:tau-infimum}. Then there exists a unique function $\mathbf n\colon \mathbb R \to \mathbb N$ such that, for every $t \in \mathbb R$, we have $t \in D_{\mathbf n(t)} \setminus D_{\mathbf n(t) + 1}$.
\end{corollary}

Our next result provides several properties of the function $\mathbf n$.

\begin{lemma}
\label{lemma:properties_n}
Let $\tau$ be a delay function satisfying \ref{hypo:tau-infimum} and $\mathbf n\colon \mathbb R \to \mathbb N$ be the function from Corollary~\ref{coro:integer-iterations}. Then, for every $t \in \mathbb R_+$, we have
\begin{enumerate}
\item\label{item:prop-n-nonzero} $\mathbf n(t) \geq 1$,
\item\label{item:n_sigma_1} $\mathbf n(\sigma_1(t)) = \mathbf n(t) - 1$, 
\item\label{item:sigma-k-nonnegative} $\sigma_k(t) \in \mathbb R_+$ for every $k \in \llbracket 0, \mathbf n(t) - 1\rrbracket$, and
\item\label{item:sigma-nt-negative} $\sigma_{\mathbf n(t)}(t) \in \lvert -h(0), 0)$.
\end{enumerate}
\end{lemma}

\begin{proof}
Since $t \in \mathbb R_+$ and $D_0 \setminus D_1 = \mathbb R_-^\ast$, we have $t \notin D_0 \setminus D_1$, and hence $\mathbf n(t) \neq 0$, implying \ref{item:prop-n-nonzero}.

To prove \ref{item:n_sigma_1}, note that, by Corollary~\ref{coro:integer-iterations}, we have
\begin{equation}
\label{eq:n_sigma_1}
\sigma_1(t) \in D_{\mathbf n(\sigma_1(t))} \setminus D_{\mathbf n(\sigma_1(t)) + 1}.
\end{equation}
Recall that, by Lemma~\ref{lemma:D_n-sigma1}, for every $n \in \mathbb N$, we have $t \in D_{n+1}$ if and only if $\sigma_1(t) \in D_n$. Hence, \eqref{eq:n_sigma_1} is equivalent to
\[
t \in D_{\mathbf n(\sigma_1(t)) + 1} \setminus D_{\mathbf n(\sigma_1(t)) + 2}.
\]
Since the function $\mathbf n$ from Corollary~\ref{coro:integer-iterations} is unique, we deduce that $\mathbf n(t) = \mathbf n(\sigma_1(t)) + 1$, yielding the conclusion.

For $k \in \llbracket 0, \mathbf n(t) - 1\rrbracket$, it follows from the definition of $\mathbf n$ and Lemma~\ref{lemma:Dn-decreasing} that $t \in D_{k+1}$. By \eqref{eq:rec-Dn}, $D_{k+1} = \{ s \in D_k \suchthat \sigma_k (s) \in D_1\}$, which implies that $\sigma_k (t) \in D_1 = \mathbb R_+$, proving \ref{item:sigma-k-nonnegative}. 

Let us finally prove \ref{item:sigma-nt-negative}. Note that $t \in D_{\mathbf n(t)}$, but $t \notin D_{\mathbf n(t) +1}$. By \eqref{eq:rec-Dn}, we have $D_{\mathbf n(t) + 1} = \{ s \in D_{\mathbf n(t)} \suchthat \sigma_{\mathbf n(t)} (s) \in D_1 \}.$ Therefore, $t \notin D_{\mathbf n(t) + 1}$ implies that $t \notin D_{\mathbf n(t)}$ or $\sigma_{\mathbf n(t)} (t) \notin D_1 = \mathbb R_+$. Since $t \in D_{\mathbf n(t)}$, then the only possibility is $\sigma_{\mathbf n(t)} (t) \notin D_1 = \mathbb R_+$. This implies that $\sigma_{\mathbf n(t)} (t) < 0$.

On the other hand, $h(0) = - \inf_{s \geq 0} \sigma_1 (s) = - \inf_{s \geq 0} (s - \tau(s)) \in (0, +\infty]$. Therefore, $\sigma_1 (s) \geq -h(0)$ for every $s \in \mathbb R_+$. By \ref{item:prop-n-nonzero}, we have $\mathbf n(t) \geq 1$, then, by \eqref{eq:rec-sigman},
\[\sigma_{\mathbf n(t)} (t) = \sigma_1 (\sigma_{\mathbf n(t) - 1} (t)) \geq - h(0),\]
proving the desired result.
\end{proof}

Another useful result regarding the function $\mathbf n$ is the following.

\begin{lemma}
\label{lemma:n-loc-bounded}
Let $\tau$ be a delay function satisfying \ref{hypo:tau-infimum} and $\mathbf n\colon \mathbb R \to \mathbb N$ be the function from Corollary~\ref{coro:integer-iterations}. Then, for every $T > 0$, there exists $N_T > 0$ such that $\mathbf n(t) \leq N_T$ for every $t \in (-\infty, T]$.
\end{lemma}

\begin{proof}
Let $T > 0$ and set $\delta = \inf_{t \in [0, T]} \tau(t)$. Note that $\delta > 0$ thanks to \ref{hypo:tau-infimum}.

We claim that, for every $k \in \mathbb N^\ast$, we have $D_k \cap [0, T] \subset [(k-1)\delta, T]$. Indeed, this is trivially verified for $k = 1$. Now, let $k \in \mathbb N^\ast$ bet such that $D_k \cap [0, T] \subset [(k-1)\delta, T]$ and take $t \in D_{k+1} \cap [0, T]$. By \eqref{eq:defi-Dn}, we have $\sigma_1(t) \in D_k$. Since $k \neq 0$, we have $D_k \subset \mathbb R_+$, hence $\sigma_1(t) \geq 0$. Moreover, $\sigma_1(t) = t - \tau(t) < t \leq T$, hence $\sigma_1(t) \in D_k \cap [0, T] \subset [(k-1)\delta, T]$. Thus $t = \sigma_1(t) + \tau(t) \geq (k-1)\delta + \delta = k \delta$, proving that $t \in [k\delta, T]$. Hence, by induction, $D_k \cap [0, T] \subset [(k-1)\delta, T]$ for every $k \in \mathbb N^\ast$.

Let $t \in (-\infty, T]$. If $t < 0$, we have $\mathbf n(t) = 0$. Otherwise, for $t \in [0, T]$, we have $t \in D_{\mathbf n(t)} \cap [0, T] \subset [(\mathbf n(t) - 1)\delta, T]$, and thus $t \geq (\mathbf n(t) - 1)\delta$, implying that $\mathbf n(t) \leq \frac{t}{\delta} + 1 \leq \frac{T}{\delta} + 1$. This implies the conclusion with $N_T = \frac{T}{\delta} + 1$.
\end{proof}

Our next result provides an expression for a solution $x$ of \eqref{eq:diff-eqn-N1-time} in terms of the $k$-th power of $A$ and $\sigma_k$.

\begin{lemma}\label{lemma-solution} Let $\tau$ be a delay function satisfying \ref{hypo:tau-infimum} and $\mathbf n\colon \mathbb R \to \mathbb N$ be the function from Corollary~\ref{coro:integer-iterations}. If $x$ is a solution of \eqref{eq:diff-eqn-N1-time}, then 
\begin{equation}\label{equation-lemma}
x(t) = A^k x(\sigma_k (t)), \qquad \text{ for every } t \in \mathbb R_+ \text{ and } k \in \llbracket 0, \mathbf n(t)\rrbracket.
\end{equation}
\end{lemma}

\begin{proof}
Fix $t \in \mathbb R_+$. We prove \eqref{equation-lemma} by induction on $k$. Consider first $k =0$, then $A^0 = \Id$ and $\sigma_0 = \id$, thus clearly we get that \eqref{equation-lemma} is satisfied.

Now, let $k \in \llbracket 0, \mathbf n(t) - 1\rrbracket$ be such that the equation \eqref{equation-lemma} is fulfilled, that is,
\begin{equation}\label{equation-induction-sigma-k}
x(t) = A^k x(\sigma_k(t)).
\end{equation}
By Lemma~\ref{lemma:properties_n}\ref{item:sigma-k-nonnegative}, we have $\sigma_k(t) \geq 0$. On the other hand, since $x$ is a solution of \eqref{eq:diff-eqn-N1-time}, we have
\begin{equation*}
x(\sigma_k (t)) = A x(\sigma_k(t) - \tau(\sigma_k (t)))) = A x(\sigma_1 \circ \sigma_k(t)) = A x(\sigma_{k+1} (t)).
\end{equation*}
Combining with \eqref{equation-induction-sigma-k}, we get
\[x(t) = A^{k+1} x(\sigma_{k+1} (t)),\]
proving our desired result.
\end{proof}

Taking $k = \mathbf n(t)$ in \eqref{equation-lemma} and noticing that $\sigma_{\mathbf n(t)} \in \lvert -h(0), 0)$ thanks to Lem\-ma~\ref{lemma:properties_n}\ref{item:sigma-nt-negative}, we obtain our main result of this section, containing a representation formula for solutions of \eqref{eq:diff-eqn-N1-time} in terms of the initial condition $x_0$, the matrix $A$, and the delay $\tau$.

\begin{theorem}\label{theorem-1} Let $\tau$ be a delay function satisfying \ref{hypo:tau-infimum} and $\mathbf n\colon \mathbb R \to \mathbb N$ be the function from Corollary~\ref{coro:integer-iterations}. If $x$ is a solution of \eqref{eq:diff-eqn-N1-time} with initial condition $x_0$, then
\begin{equation}
\label{eq:explicit-formula}
x(t) = A^{\mathbf n(t)} x_0(\sigma_{\mathbf n(t)}(t)), \qquad \text{ for every } t \in \mathbb R_+.
\end{equation}
\end{theorem}

Note that the right-hand side of \eqref{eq:explicit-formula} is a function of $t$ uniquely determined from $x_0$, $A$, and the delay $\tau$ satisfying \ref{hypo:tau-infimum}. In particular, we deduce at once the following uniqueness result.

\begin{corollary}
\label{coro:uniqueness}
Let $\tau$ be a delay function satisfying \ref{hypo:tau-infimum}, $h$ be its associated largest delay function, and $x_0\colon \lvert -h(0), 0) \to \mathbb R^d$. Then \eqref{eq:diff-eqn-N1-time} admits at most one solution with initial condition $x_0$.
\end{corollary}

\begin{remark}
\label{remk:H2}
Assumption~\ref{hypo:tau-infimum} is important to ensure uniqueness of solutions in Corollary~\ref{coro:uniqueness}. Consider, for instance, the delay function $\tau$ defined by $\tau(0) = 1$ and $\tau(t) = \left(1 - \frac{1}{e}\right)t$ for $t > 0$. In particular, \eqref{eq:diff-eqn-N1-time} reads $x(t) = A x(\frac{t}{e})$ for $t > 0$. Assume that $A$ is not nilpotent. Let $x_0\colon [-1, 0) \to \mathbb R^d$ be a given initial condition. Let $\lambda \in \mathbb C \setminus \{0\}$ be a nonzero eigenvector of $A$, $v \in \mathbb C^d \setminus \{0\}$ be a corresponding eigenvector, $\alpha \in \mathbb C$ be such that $e^\alpha = \lambda$, and $\rho \in \mathbb C$. Let $x\colon [-1, +\infty) \to \mathbb R^d$ be the function defined by
\[
x(t) = \begin{dcases*}
x_0(t), & if $t \in [-1, 0)$, \\
A x_0(-1), & if $t = 0$, \\
\Real\left(\rho e^{\alpha \ln t} v\right), & if $t > 0$.
\end{dcases*}
\]
It is immediate to verify that $x$ is a solution of \eqref{eq:diff-eqn-N1-time} with initial condition $x_0$ and, since $\rho \in \mathbb C$ is arbitrary, \eqref{eq:diff-eqn-N1-time} admits infinitely many solutions with the same initial condition $x_0$.

Assumption~\ref{hypo:tau-infimum}, however, is not necessary for the uniqueness of solutions of \eqref{eq:diff-eqn-N1-time}. Indeed, consider the delay function $\tau$ given by $\tau(0) = 1$ and $\tau(t) = t$ for $t > 0$. In this case, \eqref{eq:diff-eqn-N1-time} reads $x(0) = A x(-1)$ and $x(t) = A x(0)$ for $t > 0$. It is then immediate to verify that, given an initial condition $x_0\colon [-1, 0) \to \mathbb R^d$, \eqref{eq:diff-eqn-N1-time} admits a unique solution with initial condition $x_0$, which is given by $x(t) = x_0(t)$ for $t \in [-1, 0)$, $x(0) = A x_0(-1)$, and $x(t) = A^2 x_0(-1)$ for $t > 0$.
\end{remark}

Since continuous, measurable, and regulated solutions of \eqref{eq:diff-eqn-N1-time} are particular cases of solutions of \eqref{eq:diff-eqn-N1-time}, the representation formula \eqref{eq:explicit-formula} from Theorem~\ref{theorem-1} also applies to such notions of solutions, and Corollary~\ref{coro:uniqueness} also implies uniqueness in such classes of solutions. It is also interesting to consider the counterparts of Theorem~\ref{theorem-1} and Corollary~\ref{coro:uniqueness} for a.e.\ solutions of \eqref{eq:diff-eqn-N1-time}. This is the topic of our next result.

\begin{theorem}
\label{thm:explicit-ae}
Let $\tau$ be a delay function satisfying \ref{hypo:tau-infimum} and \ref{hypo:null} and $\mathbf n\colon \mathbb R \to \mathbb N$ be the function from Corollary~\ref{coro:integer-iterations}. If $x$ is an a.e.\ solution of \eqref{eq:diff-eqn-N1-time} with initial condition equal to $x_0$ a.e., then
\begin{equation}
\label{eq:explicit-formula-ae}
x(t) = A^{\mathbf n(t)} x_0(\sigma_{\mathbf n(t)}(t)), \qquad \text{ for almost every } t \in \mathbb R_+.
\end{equation}
In particular, one has uniqueness of solutions of \eqref{eq:diff-eqn-N1-time} with respect to a.e.\ equality: given two a.e.\ solutions $x_1, x_2\colon \lvert -h(0), +\infty) \to \mathbb R^d$ of \eqref{eq:diff-eqn-N1-time} with $x_1(t) = x_2(t)$ for a.e.\ $t \in \lvert -h(0), 0)$, we have $x_1(t) = x_2(t)$ for a.e.\ $t \in \lvert -h(0), +\infty)$.
\end{theorem}

\begin{proof}
Let
\begin{equation}
\label{eq:def-N0-N1}
\begin{aligned}
N_0 & = \{t \in \lvert -h(0), 0) \suchthat x(t) \neq x_0(t)\}, \\
N_1 & = \{t \in \mathbb R_+ \suchthat x(t) \neq A x(t - \tau(t))\}, \\
N & = N_0 \cup N_1, \\
M & = \{t \in \lvert -h(0), +\infty) \suchthat x(t) \neq A^{\mathbf n(t)} x_0(\sigma_{\mathbf n(t)}(t))\},
\end{aligned}
\end{equation}
and note that, by assumption, $\mathcal L(N_0) = \mathcal L(N_1) = \mathcal L(N) = 0$. Showing \eqref{eq:explicit-formula-ae} amounts to showing that $\mathcal L(M) = 0$. We claim that
\begin{equation}
\label{eq:M-inclusion}
M \subset \bigcup_{k = 0}^\infty \sigma_k^{-1}(N).
\end{equation}
Indeed, assume that $t \in \mathbb R$ is such that $t \notin \bigcup_{k = 0}^\infty \sigma_k^{-1}(N)$. Then, for every $k \in \mathbb N$, we have $t \notin D_k$ or $\sigma_k(t) \notin N$. In particular, for every $k \in \llbracket 0, \mathbf n(t)\rrbracket$, we have $\sigma_k(t) \notin N$. Using Lemma~\ref{lemma:properties_n}\ref{item:sigma-k-nonnegative} and \ref{item:sigma-nt-negative}, \eqref{eq:rec-sigman}, and \eqref{eq:def-N0-N1}, we deduce that $x(\sigma_k(t)) = A x(\sigma_{k+1}(t))$ for every $k \in \llbracket 0, \mathbf n(t) - 1\rrbracket$ and $x(\sigma_{\mathbf n(t)}(t)) = x_0(\sigma_{\mathbf n(t)}(t))$. An immediate inductive argument shows that $x(t) = x(\sigma_0(t)) = A^{\mathbf n(t)} x(\sigma_{\mathbf n(t)}(t)) = A^{\mathbf n(t)} x_0(\sigma_{\mathbf n(t)}(t))$, showing that $t \notin M$, and hence concluding the proof of \eqref{eq:M-inclusion}.

Notice now that, by \eqref{eq:rec-sigman}, we have $\sigma_k^{-1}(N) = \sigma_{k-1}^{-1}(\sigma_1^{-1}(N))$ for every $k \in \mathbb N^\ast$, and thus, using \ref{hypo:null} and an immediate inductive argument, we deduce that $\mathcal L(\sigma_k^{-1}(N)) = 0$ for every $k \in \mathbb N^\ast$. Hence, \eqref{eq:M-inclusion} implies that $\mathcal L(N) = 0$, yielding the conclusion.

To prove uniqueness, note that, since \eqref{eq:diff-eqn-N1-time} is linear, if $x_1, x_2\colon \lvert -h(0), +\infty) \to \mathbb R^d$ are a.e.\ solutions of \eqref{eq:diff-eqn-N1-time} with $x_1(t) = x_2(t)$ for a.e.\ $t \in \lvert -h(0), 0)$, then $x_1 - x_2$ is an a.e.\ solution of \eqref{eq:diff-eqn-N1-time} with initial condition equal a.e.\ to $0$. Hence, by \eqref{eq:explicit-formula-ae}, $x_1(t) - x_2(t) = 0$ for a.e.\ $t \in \lvert -h(0), +\infty)$, yielding the conclusion.
\end{proof}

\subsection{Existence of solutions}

We now turn to the question of existence of solutions of \eqref{eq:diff-eqn-N1-time}. Our first result is the following consequence of the representation formula for the solutions from Theorem~\ref{theorem-1}.

\begin{theorem}
\label{thm:exist-general}
Let $\tau$ be a delay function satisfying \ref{hypo:tau-infimum}, $h$ be its associated largest delay function, and $x_0\colon \lvert -h(0), 0) \to \mathbb R^d$. Then there exists a unique solution $x\colon \lvert -h(0), +\infty) \to \mathbb R^d$ of \eqref{eq:diff-eqn-N1-time} with initial condition $x_0$.
\end{theorem}

\begin{proof}
Uniqueness of the solution was already established in Corollary~\ref{coro:uniqueness}, and so we are left to prove existence. Let $x\colon \lvert -h(0), +\infty) \to \mathbb R^d$ be defined by
\begin{equation}
\label{eq:def-x-sol}
x(t) = \begin{dcases*}
x_0(t), & if $t \in \lvert -h(0), 0)$, \\
A^{\mathbf n(t)} x_0(\sigma_{\mathbf n(t)}(t)), & if $t \geq 0$,
\end{dcases*}
\end{equation}
where $\mathbf n\colon \mathbb R \to \mathbb N$ is the function from Corollary~\ref{coro:integer-iterations}. Note that, thanks to Lem\-ma~\ref{lemma:properties_n}\ref{item:sigma-nt-negative}, we have $\sigma_{\mathbf n(t)}(t) \in \lvert -h(0), 0)$ for every $t \geq 0$, and thus $x$ is well-defined. We want to prove that $x$ satisfies $x(t) = A x(t - \tau(t))$ for every $t \geq 0$.

Let $t \geq 0$ and assume first that $\sigma_1(t) < 0$. Then $\mathbf n(t) = 1$, and we have $x(\sigma_1(t)) = x_0(\sigma_{\mathbf n(t)}(t))$. Thus, using \eqref{eq:def-x-sol}, we have
\[
x(t) = A^{\mathbf n(t)} x_0(\sigma_{\mathbf n(t)}(t)) = A x(\sigma_1(t)) = A x(t - \tau(t)).
\]
Consider now the case $\sigma_1(t) \geq 0$. Using \eqref{eq:def-x-sol} and Lemmas~\ref{lemma:caract-Dn-sigman} and \ref{lemma:properties_n}\ref{item:n_sigma_1}, we have
\begin{align*}
x(\sigma_1(t)) & = A^{\mathbf n(\sigma_1(t))} x_0(\sigma_{\mathbf n(\sigma_1(t))} (\sigma_1(t))) \\
& = A^{\mathbf n(t) - 1} x_0(\sigma_{\mathbf n(t) - 1} (\sigma_1(t))) \\
& = A^{\mathbf n(t) - 1} x_0(\sigma_{\mathbf n(t)}(t)).
\end{align*}
Hence
\[
A x(t - \tau(t)) = A x(\sigma_1(t)) = A^{\mathbf n(t)} x_0(\sigma_{\mathbf n(t)}(t)) = x(t),
\]
where the last equality follows from \eqref{eq:def-x-sol}.
\end{proof}

\begin{remark}
As highlighted in Remark~\ref{remk:H2}, \ref{hypo:tau-infimum} is important in Theorem~\ref{thm:exist-general} to deduce uniqueness of solutions, but it is not needed for getting existence of solutions. Without \ref{hypo:tau-infimum}, one cannot guarantee that the result of Lemma~\ref{lemma:intersect-Dn-empty} is true. However, setting $D_\infty = \bigcap_{n \in \mathbb N} D_n \subset \mathbb R_+$, we still get the conclusion of Corollary~\ref{coro:integer-iterations} with a function $\mathbf n$ defined on $\mathbb R \setminus D_\infty$ instead of $\mathbb R$. The results of Lemmas~\ref{lemma:properties_n} and \ref{lemma-solution} and Theorem~\ref{theorem-1} remain true without assuming \ref{hypo:tau-infimum} after replacing $t \in \mathbb R_+$ by $t \in \mathbb R_+ \setminus D_\infty$ in their statements and proofs. Then, given an initial condition $x_0\colon \lvert -h(0), 0) \to \mathbb R^d$, we define $x\colon \lvert -h(0), +\infty) \to \mathbb R^d$ by
\[
x(t) = \begin{dcases*}
x_0(t), & if $t \in \lvert -h(0), 0)$, \\
A^{\mathbf n(t)} x_0(\sigma_{\mathbf n(t)}(t)), & if $t \in \mathbb R_+ \setminus D_\infty$, \\
0, & if $t \in D_\infty$.
\end{dcases*}
\]
By the same argument used in the proof of Theorem~\ref{thm:exist-general}, we obtain that $x(t) = A x(t - \tau(t))$ for every $t \in \mathbb R_+ \setminus D_\infty$. In addition, if $t \in D_\infty$, then $t - \tau(t) = \sigma_1(t) \in D_\infty$, and thus the identity $x(t) = A x(t - \tau(t))$ is verified for every $t \in D_\infty$. This shows that $x$ is a solution of \eqref{eq:diff-eqn-N1-time} with initial condition $x_0$, and hence proves the existence assertion of Theorem~\ref{thm:exist-general} without the need of \ref{hypo:tau-infimum}.

In this chapter, we are only interested in cases where the solution of \eqref{eq:diff-eqn-N1-time} is unique, and hence we shall always assume \ref{hypo:tau-infimum} in the sequel. However, by arguments similar to those in this remark, some of our next results on existence of solutions can also be generalized to settings where \ref{hypo:tau-infimum} is not assumed to hold true.
\end{remark}

Our next result considers the existence of continuous solutions of \eqref{eq:diff-eqn-N1-time}.

\begin{theorem}
\label{thm:exist-continuous}
Let $\tau$ be a delay function satisfying \ref{hypo:tau-infimum} and \ref{hypo:tau-cont}, $h$ be its associated largest delay function, and $x_0 \in \mathsf C_0$. Then there exists a unique continuous solution $x \in \mathcal C(\lvert -h(0), +\infty), \mathbb R^d)$ of \eqref{eq:diff-eqn-N1-time} with initial condition $x_0$. In addition, $x_t \in \mathsf C_t$ for every $t \in \mathbb R_+$.
\end{theorem}

\begin{proof}
By Theorem~\ref{thm:exist-general}, \eqref{eq:diff-eqn-N1-time} admits a unique solution $x\colon \lvert -h(0), +\infty) \to \mathbb R^d$ with initial condition $x_0$, and we are thus only left to show that such a solution is continuous. To do so, we will prove that $x$ is continuous in $\lvert -h(0), T]$ for every $T > 0$.

Fix $T > 0$ and let $\widetilde\delta = \frac{1}{2}\inf_{t \in [0, T]} \tau(t)$. Note that, by \ref{hypo:tau-infimum}, we have $\widetilde\delta > 0$. Set $\delta = \frac{T}{\ceil{T / \widetilde\delta}}$ and note that $0 < \delta \leq \widetilde\delta$. By definition, $\tau(t) > \delta$ for every $t \in [0, T]$ and $T/\delta$ is a positive integer, which we denote by $M$.

We prove by induction that, for every $k \in \llbracket 0, M\rrbracket$, $x$ is continuous on $\lvert -h(0), k \delta]$, which will conclude the proof of the result since $\lvert -h(0), M \delta] = \lvert -h(0), T]$. This is trivially true for $k = 0$ since $x$ coincides with $x_0$ in $\lvert -h(0), 0)$ by definition, their values also coincide at $t = 0$ since $x(0) = A x_0(-\tau(0)) = x_0(0)$, and $x_0$ is continuous by assumption.

Let $k \in \llbracket 0, M-1\rrbracket$ be such that $x$ is continuous on $\lvert -h(0), k \delta]$. Define $\widehat x = x \rvert_{\lvert -h(0), k \delta]}$ and $y = x \rvert_{[k \delta, (k+1) \delta]}$ and notice that, thanks to the definition of $\delta$, for every $t \in [k \delta, (k+1) \delta]$, we have $y(t) = A \widehat x(t - \tau(t))$. Since $\widehat x$ and $\tau$ are continuous, $y$ is thus continuous as well. Hence $x$ is continuous on $\lvert -h(0), (k+1) \delta]$, as required.

Finally, for every $t \in \mathbb R_+$, we have $x_t(0) = x(t) = A x(t - \tau(t)) = A x_t(-\tau(t))$, showing that $x_t \in \mathsf C_t$.
\end{proof}

It is immediate to verify the importance of \ref{hypo:tau-cont} in Theorem~\ref{thm:exist-continuous}, as one can easily come up with examples of discontinuous delay functions yielding a discontinuous solution of \eqref{eq:diff-eqn-N1-time} for some initial condition in $\mathsf C_0$. A natural question is whether \ref{hypo:tau-cont} is necessary, and our next example provides a negative answer to that question.

\begin{example}
\label{expl:exist-continuous}
Let $\tau\colon \mathbb R_+ \to \mathbb R_+^\ast$ be the delay function defined by
\[
\tau(t) = \begin{dcases*}
t + 1, & if $t \in [0, 1)$, \\
2, & if $t \in [1, 2)$, \\
t - 1, & if $t \in [2, +\infty)$.
\end{dcases*}
\]
Note that $\tau$ is discontinuous at $t = 2$. Let $x_0 \in \mathsf C_0$ and $x\colon [-1, +\infty) \to \mathbb R^d$ be the unique solution of \eqref{eq:diff-eqn-N1-time} with initial condition $x_0$ provided by Theorem~\ref{thm:exist-general}. Then $x$ is continuous in $[-1, +\infty)$. Indeed, a straightforward computation from \eqref{eq:diff-eqn-N1-time} shows that
\begin{equation}
\label{eq:x-continuous}
x(t) = \begin{dcases*}
x_0(t), & if $t \in [-1, 0)$, \\
A x_0(-1), & if $t \in [0, 1)$, \\
A x_0(t - 2), & if $t \in [1, 2)$, \\
A^2 x_0(-1), & if $t \in [2, +\infty)$,
\end{dcases*}
\end{equation}
and it is immediate to verify that $x$ is continuous using the compatibility condition $x_0(0) = A x_0(-1)$, since $x_0 \in \mathsf C_0$.
\end{example}

Our next result provides conditions for the existence and uniqueness of regulated solutions of \eqref{eq:diff-eqn-N1-time}.

\begin{theorem}
\label{thm:exist-regulated}
Let $\tau$ be a delay function satisfying \ref{hypo:tau-infimum} and \ref{hypo:regulated}, $h$ be its associated largest delay function, and $x_0 \in \mathsf G_0$. Then there exists a unique regulated solution $x \in \mathcal G(\lvert -h(0), +\infty), \mathbb R^d)$ of \eqref{eq:diff-eqn-N1-time} with initial condition $x_0$. In addition, $x_t \in \mathsf G_t$ for every $t \in \mathbb R_+$.
\end{theorem}

\begin{proof}
As in the proof of Theorem~\ref{thm:exist-continuous}, we are only left to show that the unique solution $x\colon \lvert -h(0), +\infty) \to \mathbb R^d$ of \eqref{eq:diff-eqn-N1-time} with initial condition $x_0$ provided by Theorem~\ref{thm:exist-general} is regulated. The argument is very similar to that in the proof of Theorem~\ref{thm:exist-continuous}: we set $T > 0$, $\delta > 0$, and $M \in \mathbb N^\ast$ as in that proof, and we remark that the restriction of $x$ to $\lvert -h(0), 0]$ is regulated by assumption. If $k \in \llbracket 0, M-1\rrbracket$ is such that the restriction of $x$ to $\lvert -h(0), k \delta]$ is regulated, we define $\widehat x$ and $y$ as in the proof of Theorem~\ref{thm:exist-continuous}. For every $t \in [k\delta, (k+1)\delta]$, we have $y(t) = A \widehat x(t - \tau(t))$. Since $\widehat x$ is regulated and $t \mapsto t - \tau(t)$ is well-regulated, we deduce from Theorem~\ref{thm:well-regulated} that $y$ is regulated, yielding that the restriction of $x$ to $\lvert -h(0), (k+1)\delta]$ is regulated, as required.

Finally, for every $t \in \mathbb R_+$, we have $x_t(0) = x(t) = A x(t - \tau(t)) = A x_t(-\tau(t))$, showing that $x_t \in \mathsf G_t$.
\end{proof}

Assumption~\ref{hypo:regulated} is important in Theorem~\ref{thm:exist-regulated} to guarantee regulatedness of the solution $x$ by using Theorem~\ref{thm:well-regulated}. Note that, due to Remark~\ref{remk:sum-not-well-regulated}, the well-regulatedness assumption cannot be done in $\tau$, since this is not sufficient to guarantee the well-regulatedness of $t \mapsto t - \tau(t)$. Even though \ref{hypo:regulated} is important, this assumption is not necessary, since one can construct an example similar in spirit to Example~\ref{expl:exist-continuous}, as shown below.

\begin{example}
Let $\varphi\colon [1, +\infty) \to [0, 1)$ be a continuous function that is not well-regulated (such a function can be constructed, for instance, from a suitable translation of the function from Remark~\ref{remk:cinfty-not-well-regulated}). We define the delay function $\tau\colon \mathbb R_+ \to \mathbb R_+^\ast$ by
\[
\tau(t) = \begin{dcases*}
t + 1, & if $t \in [0, 1)$, \\
t - \varphi(t), & if $t \in [1, \infty)$.
\end{dcases*}
\]
Note that $\tau$ does not satisfy \ref{hypo:regulated}. Let $x_0 \in \mathsf G_0$ and $x\colon [-1, +\infty) \to \mathbb R^d$ be the unique solution of \eqref{eq:diff-eqn-N1-time} with initial condition $x_0$ provided by Theorem~\ref{thm:exist-general}. A straightforward computation from \eqref{eq:diff-eqn-N1-time} shows that
\begin{equation*}
x(t) = \begin{dcases*}
x_0(t), & if $t \in [-1, 0)$, \\
A x_0(-1), & if $t \in [0, 1)$, \\
A x(\varphi(t)) = A^2 x_0(-1), & if $t \in [1, +\infty)$,
\end{dcases*}
\end{equation*}
where we use the fact that $\varphi$ takes values in $[0, 1)$. Hence $x$ is regulated in $[-1, +\infty)$.
\end{example}

We now turn to the existence of measurable solutions of \eqref{eq:diff-eqn-N1-time}.

\begin{theorem}
\label{thm:exist-measurable}
Let $\tau$ be a delay function satisfying \ref{hypo:tau-infimum} and \ref{hypo:tau-measurable}, $h$ be its associated largest delay function, and $x_0\colon \lvert -h(0), 0) \to \mathbb R^d$ be $(\mathfrak L, \mathfrak B)$-measurable. Then there exists a unique measurable solution $x\colon \lvert -h(0), +\infty) \to \mathbb R^d$ of \eqref{eq:diff-eqn-N1-time} with initial condition $x_0$.
\end{theorem}

\begin{proof}
As in the proof of Theorems~\ref{thm:exist-continuous} and \ref{thm:exist-regulated}, we only need to show that unique solution $x\colon \lvert -h(0), +\infty) \to \mathbb R^d$ of \eqref{eq:diff-eqn-N1-time} with initial condition $x_0$ provided by Theorem~\ref{thm:exist-general} is $(\mathfrak L, \mathfrak B)$-measurable. We proceed by an inductive argument similar to those of the proofs of Theorems~\ref{thm:exist-continuous} and \ref{thm:exist-regulated}.

Let $T > 0$, $\delta > 0$, and $M \in \mathbb N^\ast$ be as in those proofs. Clearly, the restriction of $x$ to $\lvert -h(0), 0]$ is $(\mathfrak L, \mathfrak B)$-measurable, since it coincides with $x_0$ in $\lvert -h(0), 0)$ and $x_0$ is $(\mathfrak L, \mathfrak B)$-measurable.

Now, let $k \in \llbracket 0, M-1\rrbracket$ be such that the restriction of $x$ to $\lvert -h(0), k \delta]$ is $(\mathfrak L, \mathfrak B)$-measurable and define $\widehat x$ and $y$ as in the proofs of Theorems~\ref{thm:exist-continuous} and \ref{thm:exist-regulated}. For every $t \in [k\delta, (k+1)\delta]$, we have $y(t) = A \widehat x(t - \tau(t))$. Since $\widehat x$ is $(\mathfrak L, \mathfrak B)$-measurable and $t \mapsto t - \tau(t)$ is $(\mathfrak L, \mathfrak L)$-measurable, we deduce that $y$ is $(\mathfrak L, \mathfrak B)$-measurable, yielding the $(\mathfrak L, \mathfrak B)$-measurability of the restriction of $x$ to $\lvert -h(0), (k+1)\delta]$, as required.
\end{proof}

The counterpart in Theorem~\ref{thm:exist-measurable} of assumption \ref{hypo:tau-cont} from Theorem~\ref{thm:exist-continuous} and of assumption \ref{hypo:regulated} from Theorem~\ref{thm:exist-regulated} is \ref{hypo:tau-measurable}. Similarly to what was shown in Example~\ref{expl:exist-continuous}, we provide now an example showing that \ref{hypo:tau-measurable} is not necessary for the existence of measurable solutions.

\begin{example}
Let $P \subset [2, +\infty)$ be a set that is not Lebesgue measurable. Let $\tau\colon \mathbb R_+ \to \mathbb R_+^\ast$ be the delay function defined by
\[
\tau(t) = \begin{dcases*}
t + 1, & if $t \in [0, 1)$, \\
2, & if $t \in [1, 2)$, \\
t - \frac{2}{3} + \frac{1}{3} \1_P(t), & if $t \in [2, +\infty)$,
\end{dcases*}
\]
where $\1_P$ is the indicator function of $P$. Since $P \notin \mathfrak L$, the function $\tau$ is not $(\mathfrak L, \mathfrak L)$-measurable. Let $x_0\colon [-1, 0) \to \mathbb R^d$ be a $(\mathfrak L, \mathfrak B)$-mea\-su\-ra\-ble function and $x\colon [-1, +\infty) \to \mathbb R^d$ be the unique solution of \eqref{eq:diff-eqn-N1-time} with initial condition $x_0$ provided by Theorem~\ref{thm:exist-general}. Then $x$ is $(\mathfrak L, \mathfrak B)$-mea\-su\-ra\-ble in $[-1, +\infty)$. Indeed, similarly to Example~\ref{expl:exist-continuous}, $x$ satisfies \eqref{eq:x-continuous} for $t \in [-1, 2)$. For $t \geq 2$, we have from \eqref{eq:diff-eqn-N1-time} that $x(t) = A x(\frac{2}{3} - \frac{1}{3} \1_P(t))$ and, since $\frac{2}{3} - \frac{1}{3} \1_P(t) \in [0, 1)$ for every $t \geq 2$, we have $x(t) = A^2 x_0(-1)$. Hence $x$ satisfies \eqref{eq:x-continuous} for every $t \in [-1, +\infty)$, and it is immediate to verify that the right-hand side of \eqref{eq:x-continuous} defines a $(\mathfrak L, \mathfrak B)$-measurable function of $t$.
\end{example}

Our next result concerns existence and uniqueness of a.e.\ solutions.

\begin{theorem}
\label{thm:exist-ae}
Let $\tau$ be a delay function satisfying \ref{hypo:tau-infimum} and \ref{hypo:null}, $h$ be its associated largest delay function, and $x_0\colon \lvert -h(0), 0) \to  \mathbb R^d$. Then the unique solution $x\colon \lvert -h(0), +\infty) \to \mathbb R^d$ of \eqref{eq:diff-eqn-N1-time} with initial condition $x_0$ from Theorem~\ref{thm:exist-general} is also an a.e.\ solution of \eqref{eq:diff-eqn-N1-time} with initial condition $x_0$. In addition, such an a.e.\ solution is unique with respect to a.e.\ equality (in the sense of Theorem~\ref{thm:explicit-ae}).
\end{theorem}

\begin{proof}
Notice first that the assertion on uniqueness follows from Theorem~\ref{thm:explicit-ae}.

Let $x\colon \lvert -h(0), +\infty) \allowbreak \to \mathbb R^d$ be the solution of \eqref{eq:diff-eqn-N1-time} with initial condition $x_0$, whose existence and uniqueness is asserted in Theorem~\ref{thm:exist-general}. Let us prove that $x$ is also an a.e.\ solution of \eqref{eq:diff-eqn-N1-time}. Let $\widetilde x\colon \lvert -h(0), +\infty) \to \mathbb R^d$ be such that $\widetilde x(t) = x(t)$ for a.e.\ $t \in \lvert -h(0), +\infty)$. Let $N_0 = \{t \in \lvert -h(0), +\infty) \suchthat x(t) \neq \widetilde x(t)\}$ and notice that $x(t - \tau(t)) \neq \widetilde x(t - \tau(t))$ if and only if $t \in \sigma_1^{-1}(N_0)$. Let $N = N_0 \cup \sigma_1^{-1}(N_0)$ and note that, thanks to \ref{hypo:null}, $\mathcal L(N) = 0$. For every $t \in \mathbb R_+ \setminus N$, we have
\[
\widetilde x(t) = x(t) = A x(t - \tau(t)) = A \widetilde x(t - \tau(t)).
\]
Hence $x$ is an a.e.\ solution of \eqref{eq:diff-eqn-N1-time} with initial condition $x_0$.
\end{proof}

It turns out that \ref{hypo:null} is also a necessary condition for the existence of a.e.\ solutions of \eqref{eq:diff-eqn-N1-time} in the nontrivial case where $A \neq 0$, as detailed in the next result.

\begin{proposition}
\label{prop:h2-necessary}
Let $\tau$ be a delay function and assume that $A \neq 0$ and there exists an a.e.\ solution of \eqref{eq:diff-eqn-N1-time}. Then $\tau$ satisfies \ref{hypo:null}.
\end{proposition}

\begin{proof}
Fix an a.e.\ solution $x\colon \lvert -h(0), +\infty) \to \mathbb R^d$ of \eqref{eq:diff-eqn-N1-time}. Let $v \in \mathbb R^d \setminus \{0\}$ be such that $A v \neq 0$. Since the image of $\sigma_1$ is included in $\lvert -h(0), +\infty)$, it suffices to show that $\mathcal L(\sigma_1^{-1}(N)) = 0$ for every Lebesgue-measurable set $N \subset \lvert -h(0), +\infty)$ such that $\mathcal L(N) = 0$.

Let $N \subset \lvert -h(0), +\infty)$ be Lebesgue measurable with $\mathcal L(N) = 0$. Let $\widetilde x\colon \lvert -h(0),\allowbreak +\infty) \to \mathbb R^d$ be given by $\widetilde x(t) = x(t) + v$ if $t \in N$ and $\widetilde x(t) = x(t)$ otherwise. Hence $\widetilde x$ is equal a.e.\ to $x$ and, since $x$ is an a.e.\ solution of \eqref{eq:diff-eqn-N1-time}, there exist Lebesgue-measurable sets $N_0$ and $\widetilde N_0$ included in $\mathbb R_+$ with $\mathcal L(N_0) = \mathcal L(\widetilde N_0) = 0$ and such that
\begin{equation*}
\begin{aligned}
\widetilde x(t) & = A \widetilde x(\sigma_1(t)), & \qquad & \text{ for every } t \in \mathbb R_+ \setminus \widetilde N_0, \\
x(t) & = A x(\sigma_1(t)), & & \text{ for every } t \in \mathbb R_+ \setminus N_0. \\
\end{aligned}
\end{equation*}
For every $t \in \sigma_1^{-1}(N) \setminus (N_0 \cup \widetilde N_0)$, we then have
\[
\widetilde x(t) = A \widetilde x(\sigma_1(t)) = A x(\sigma_1(t)) + A v = x(t) + A v \neq x(t).
\]
Since $\widetilde x$ is equal a.e.\ to $x$, we deduce that $\mathcal L(\sigma_1^{-1}(N) \setminus (N_0 \cup \widetilde N_0)) = 0$. In addition,
\[
\sigma_1^{-1}(N) = \left(\sigma_1^{-1}(N) \setminus (N_0 \cup \widetilde N_0)\right) \cup \left(\sigma_1^{-1}(N) \cap (N_0 \cup \widetilde N_0)\right).
\]
Since $\mathcal L(\sigma_1^{-1}(N) \cap (N_0 \cup \widetilde N_0)) = 0$, we finally deduce that $\mathcal L(\sigma_1^{-1}(N)) = 0$.
\end{proof}

We now turn to the question of existence of $L^p$ solutions of \eqref{eq:diff-eqn-N1-time}. The conditions for the existence of such solutions turn out to depend on whether $p$ is finite or infinite, and we start with the case $p = +\infty$.

{\sloppy
\begin{theorem}
\label{thm:exist-L-infty}
Let $\tau$ be a delay function satisfying \ref{hypo:tau-infimum}, \ref{hypo:tau-measurable}, and \ref{hypo:null}, $h$ be its associated largest delay function, and $x_0 \in L^\infty(\lvert -h(0), 0], \mathbb R^d)$. Then there exists a unique $L^\infty$ solution $x \in L^\infty_{\loc}(\lvert -h(0), +\infty),\allowbreak \mathbb R^d)$ of \eqref{eq:diff-eqn-N1-time} with initial condition $x_0$.
\end{theorem}
}

\begin{proof}
Since $L^\infty$ solutions are particular cases of a.e.\ solutions, the assertion on the uniqueness of an $L^\infty$ solution follows from the uniqueness assertion for a.e.\ solutions from Theorem~\ref{thm:exist-ae}. We are thus left to prove the existence of an $L^\infty$ solution with a given initial condition.

Let $x_0 \in L^\infty(\lvert -h(0), 0], \mathbb R^d)$. Up to replacing $x_0$ by a suitable element in its equivalence class with respect to a.e.\ equality, we can assume that $x_0$ is bounded, i.e., there exists $M > 0$ such that $\abs{x_0(t)} \leq M$ for every $t \in \lvert -h(0), 0]$.

Since $x_0 \in L^\infty(\lvert -h(0), 0], \mathbb R^d)$, we have in particular that $x_0$ is $(\mathfrak L, \mathfrak B)$-measurable and thus, by Theorem~\ref{thm:exist-measurable}, \eqref{eq:diff-eqn-N1-time} admits a unique measurable solution $x\colon \lvert -h(0),\allowbreak +\infty) \to \mathbb R^d$ with initial condition $x_0$. By Theorem~\ref{thm:exist-ae}, $x$ is also an a.e.\ solution of \eqref{eq:diff-eqn-N1-time} with initial condition $x_0$. By Theorem~\ref{theorem-1}, we have $x(t) = A^{\mathbf n(t)} x_0(\sigma_{\mathbf n(t)}(t))$ for every $t \in \mathbb R_+$, and thus $\abs{x(t)} \leq \abs{A}^{\mathbf n(t)} M$ for every $t \geq 0$. The conclusion now follows thanks to Lemma~\ref{lemma:n-loc-bounded}.
\end{proof}

We conclude this section by the following result on the existence and uniqueness of solutions in $L^p$ for $p \in [1, +\infty)$.

\begin{theorem}
\label{thm:exist-L-p}
Let $p \in [1, +\infty)$,  $\tau$ be a delay function satisfying \ref{hypo:tau-infimum} and \ref{hypo:radon-nik-bdd}, $h$ be its associated largest delay function, and $x_0 \in L^p(\lvert -h(0), 0], \mathbb R^d)$. Then there exists a unique $L^p$ solution $x \in L^p_{\loc}(\lvert -h(0), +\infty),\allowbreak \mathbb R^d)$ of \eqref{eq:diff-eqn-N1-time} with initial condition $x_0$.
\end{theorem}

\begin{proof}
As in the proof of Theorem~\ref{thm:exist-L-infty}, uniqueness of solutions is a consequence of Theorem~\ref{thm:exist-ae}.

Let $x_0 \in L^p(\lvert -h(0), 0], \mathbb R^d)$. Arguing as in the proof of Theorem~\ref{thm:exist-L-infty}, \eqref{eq:diff-eqn-N1-time} admits a unique solution $x\colon \lvert -h(0), +\infty) \to \mathbb R^d$ with initial condition $x_0$, which is both a measurable and an a.e.\ solution. We now take $T > 0$ and define $\delta > 0$ and $M \in \mathbb N^\ast$ as in the proof of Theorem~\ref{thm:exist-continuous}. Let $\varphi\colon \mathbb R \to \mathbb R_+$ be the Radon--Nikodym derivative of $\left(\sigma_1\rvert_{[0, T]}\right)_{\#}\mathcal L$ with respect to $\mathcal L$ and denote by $K > 0$ an upper bound of $\varphi$. We prove by induction on $k$ that $x \in L^p(\lvert -h(0), k \delta], \mathbb R^d)$ for every $k \in \llbracket 0, M\rrbracket$, which will conclude the proof since $M\delta = T$ and $T > 0$ is arbitrary. The case $k = 0$ is clear since $x$ coincides with $x_0$ in $\lvert -h(0), 0)$.

Let $k \in \llbracket 0, M-1\rrbracket$ be such that $x \in L^p(\lvert -h(0), k \delta], \mathbb R^d)$. Define $\widehat x$ and $y$ as in the proof of Theorem~\ref{thm:exist-continuous} and note that, as in that proof, $\sigma_1([k\delta, (k+1)\delta]) \subset \lvert -h(0), k\delta]$ and $y(t) = A \widehat x(\sigma_1(t))$ for every $t \in [k\delta, (k+1)\delta]$. By \ref{hypo:radon-nik-bdd} and Theorem~\ref{thm:composition-Lp}\eqref{item:composition-Lp-equivs}, we deduce that $y \in L^p([k\delta, (k+1)\delta], \mathbb R^d)$. Hence $x \in L^p(\lvert -h(0), (k+1) \delta], \mathbb R^d)$, yielding the conclusion.
\end{proof}

Contrarily to the case of $L^\infty$ from Theorem~\ref{thm:exist-L-infty}, assumptions \ref{hypo:tau-infimum}, \ref{hypo:tau-measurable}, and \ref{hypo:null} are not sufficient to ensure existence of solutions in $L^p$ with $p \in [1, +\infty)$, as we detail in the following example.

\begin{example}
\label{expl:H6-important-Lp}
Let $\tau$ be the delay function defined for $t \in \mathbb R_+$ by $\tau(t) = 2 - \frac{1}{t + 1}$. Then $\sigma_1(t) = t - 2 + \frac{1}{t + 1}$ for $t \in \mathbb R_+$. One easily verifies that $\sigma_1$ is increasing, maps $\mathbb R_+$ onto $[-1, +\infty)$, and is a diffeomorphism between $\mathbb R_+^\ast$ and $(-1, +\infty)$, with $\sigma_1^{-1}(s) = \frac{s + 1 + \sqrt{(s + 1)(s + 5)}}{2}$ and $(\sigma_1^{-1})^\prime(s) = \frac{1}{2} + \frac{\sqrt{s + 5}}{4 \sqrt{s + 1}} + \frac{\sqrt{s + 1}}{4 \sqrt{s + 5}}$ for every $s \in (-1, +\infty)$. These facts also allow one to show that \ref{hypo:tau-infimum}, \ref{hypo:tau-measurable}, and \ref{hypo:null} are verified, but \ref{hypo:radon-nik-bdd} is not. One computes in addition that $\mathbf n(t) = 1$ if and only if $t \in [0, \frac{1 + \sqrt{5}}{2})$.

Let $p \in [1, +\infty)$ and $x_0 \in L^p([-1, 0], \mathbb R^d)$. By combining Theorems~\ref{thm:exist-general}, \ref{thm:exist-measurable}, and \ref{thm:exist-ae}, \eqref{eq:diff-eqn-N1-time} admits a unique solution $x\colon [-1, +\infty) \to \mathbb R^d$ with initial condition $x_0$, which is also both a measurable and a.e.\ solution. An immediate computation from \eqref{eq:diff-eqn-N1-time} shows that, for every $t \in [0, \frac{1 + \sqrt{5}}{2})$, we have $x(t) = A x_0\left(t - 2 + \frac{1}{t + 1}\right)$. Hence,
\begin{align*}
\int_0^{\frac{1 + \sqrt{5}}{2}} \abs*{x(t)}^p \diff t & = \int_0^{\frac{1 + \sqrt{5}}{2}} \abs*{A x_0\left(t - 2 + \frac{1}{t + 1}\right)}^p \diff t \\
& = \int_{-1}^{0} \abs*{A x_0(s)}^p \left(\frac{1}{2} + \frac{\sqrt{s + 5}}{4 \sqrt{s + 1}} + \frac{\sqrt{s + 1}}{4 \sqrt{s + 5}}\right) \diff s,
\end{align*}
and the above integral may be equal to $+\infty$, in which case $x \notin L^p_{\loc}([-1, +\infty), \mathbb R^d)$.

For instance, assume that $A \neq 0$ and let $\alpha \in \left[\frac{1}{2p}, \frac{1}{p}\right)$ and $v \in \mathbb R^d \setminus \{0\}$ be such that $A v \neq 0$. Let $x_0$ be defined by $x_0(t) = \frac{1}{(t + 1)^\alpha} v$ for $t \in (-1, 0)$ and $x_0(-1) = 0$. Then $x_0 \in L^p([-1, 0], \mathbb R^d)$ but
\[
\int_{-1}^{0} \abs*{A x_0(s)}^p \frac{\sqrt{s + 5}}{4 \sqrt{s + 1}} \diff s = \abs{A v}^p \int_{-1}^0 \frac{\sqrt{s + 5}}{4 (s + 1)^{\alpha p + \frac{1}{2}}} \diff s = +\infty
\]
since $\alpha p + \frac{1}{2} \geq 1$.
\end{example}

Despite the importance of \ref{hypo:radon-nik-bdd} for Theorem~\ref{thm:exist-L-p}, illustrated by Theorem~\ref{thm:composition-Lp}\eqref{item:composition-Lp-equivs} and Example~\ref{expl:H6-important-Lp}, \ref{hypo:radon-nik-bdd} is not necessary for the existence of $L^p$ solutions of \eqref{eq:diff-eqn-N1-time}, as illustrated by the following example.

\begin{example}
Let $\tau$ be a delay function such that the function $\sigma_1 = \id - \tau$ is given, for $t \in \mathbb R_+$, by
\[
\sigma_1(t) = \begin{dcases*}
\sqrt[3]{t - 1} - 1, & if $t \in [0, 2)$, \\
(t - 3)^3 + 1, & if $t \in [2, 4)$, \\
t - 2, & if $t \in [4, +\infty)$.
\end{dcases*}
\]
Note that ${\sigma_1}_\# \mathcal L$ is $\sigma$-finite and absolutely continuous with respect to $\mathcal L$ and the function $\sigma_1\colon \mathbb R_+ \to [-2, +\infty)$ is increasing and invertible with
\[
\sigma_1^{-1}(s) = \begin{dcases*}
(s + 1)^3 + 1, & if $s \in [-2, 0)$, \\
\sqrt[3]{s - 1} + 3, & if $s \in [0, 2)$, \\
s + 2, & if $s \in [2, +\infty)$.
\end{dcases*}
\]
Hence, one can compute that the Radon--Nikodym derivative $\varphi$ of ${\sigma_1}_\# \mathcal L$ with respect to $\mathcal L$ is the unbounded function $\varphi$ defined for $s \in \mathbb R \setminus \{1\}$ by
\[
\varphi(s) = \begin{dcases*}
0, & if $s \in (-\infty, -2)$, \\
3 (s + 1)^2, & if $s \in [-2, 0)$, \\
\frac{1}{3 (s - 1)^{2/3}}, & if $s \in [0, 2) \setminus \{1\}$, \\
1, & if $s \in [2, +\infty)$.
\end{dcases*}
\]

Even though $\varphi$ is unbounded, \eqref{eq:diff-eqn-N1-time} admits a unique $L^p$ solution for a given initial condition in $L^p$. Indeed, given $x_0 \in L^p([-2, 0], \mathbb R^d)$, it follows from Theorems~\ref{thm:exist-measurable} and \ref{thm:exist-ae} that \eqref{eq:diff-eqn-N1-time} admits a unique solution with initial condition $x_0$, which is both a measurable solution and an a.e.\ solution, and, by Theorem~\ref{thm:explicit-ae}, this solution satisfies \eqref{eq:explicit-formula-ae}. Explicit computations show that, for every $n \in \mathbb N^\ast$, we have $D_n = [2(n-1), +\infty)$, so that $\mathbf n(t) = \floor{t / 2} + 1$ for every $t \in \mathbb R_+$. In addition, by an immediate inductive argument, we have, for $n \geq 2$ and $t \in D_n$,
\[
\sigma_n(t) = \begin{dcases*}
t - 2n, & if $t \in [2(n-1), 2n) \cup [2(n+1), +\infty)$, \\
(t - 2n - 1)^3 + 1, & if $t \in [2n, 2(n+1))$.
\end{dcases*}
\]
In particular, for $t \in \mathbb R_+$, we have
\[
\sigma_{\mathbf n(t)}(t) = \begin{dcases*}
\sqrt[3]{t - 1} - 1, & if $t \in [0, 2)$, \\
t - 2n, & if $t \in [2(n-1), 2n)$ for some $n \in \mathbb N$ with $n \geq 2$.
\end{dcases*}
\]
Hence, by \eqref{eq:explicit-formula-ae}, for a.e.\ $t \in \mathbb R_+$, we have
\[
x(t) = \begin{dcases*}
A x_0(\sqrt[3]{t - 1} - 1), & if $t \in [0, 2)$, \\
A^n x_0(t - 2n), & if $t \in [2(n-1), 2n)$ for some $n \in \mathbb N$ with $n \geq 2$,
\end{dcases*}
\]
and thus one immediately verifies that $x \in L^p_{\loc}([-2, +\infty), \mathbb R^d)$, yielding that $x$ is an $L^p$ solution of \eqref{eq:diff-eqn-N1-time} with initial condition $x_0$.
\end{example}

\section{Asymptotic behavior of solutions}
\label{sec:stability}

Now that we have established results on the existence of solutions of \eqref{eq:diff-eqn-N1-time} in Section~\ref{sec:well-posed-varying-delay}, we turn to the analysis of their asymptotic behavior, with the aim of providing sufficient (and, whenever possible, necessary) conditions for convergence of solutions to $0$.

\subsection{Convergence to the origin}

Our first result provides sufficient conditions for solutions of \eqref{eq:diff-eqn-N1-time} with bounded initial condition to converge to $0$.

\begin{theorem}
\label{thm:cv-to-0}
Let $A \in \mathcal M_d(\mathbb R)$ satisfy \ref{hypo:spr-A-less-1}, $\tau$ be a delay function satisfying \ref{hypo:tau-infimum} and \ref{hypo:delay-infinity}, $h$ be its associated largest delay function, and $x\colon \lvert -h(0), +\infty) \to \mathbb R^d$ be a solution of \eqref{eq:diff-eqn-N1-time} with bounded initial condition $x_0\colon \lvert -h(0), 0) \to \mathbb R^d$. Then $x(t) \to 0$ as $t \to +\infty$.
\end{theorem}

\begin{proof}
Fix a norm $\abs{\cdot}$ in $\mathbb R^d$ such that $\abs{A} < 1$, which exists since $\spr(A) < 1$ (see, e.g., \cite[Lem\-ma~5.6.10]{Horn2013Matrix}).

We first prove that $\limsup_{t \to +\infty} \abs{x(t)} < +\infty$. Let $M > 0$ be such that $\abs{x_0(t)} \leq M$ for every $t \in \lvert -h(0), 0)$. By Theorem~\ref{theorem-1}, we then have, for every $t \in \mathbb R_+$,
\[
\abs{x(t)} \leq \abs{A^{\mathbf n(t)}} \abs{x_0(\sigma_{\mathbf n(t)}(t))} \leq M,
\]
where $\mathbf n$ and $(\sigma_n)_{n \in \mathbb N}$ are defined as in Section~\ref{sec:explicit-formula}. Hence $\limsup_{t \to +\infty} \abs{x(t)} \leq M$.

For every $t \geq 0$, using \eqref{eq:diff-eqn-N1-time}, we have \(\abs{x(t)} \leq \abs{A} \abs{x(t - \tau(t))}\). In particular,
\[
\limsup_{t \to +\infty}\abs{x(t)} \leq \abs{A} \limsup_{t \to +\infty}\abs{x(t - \tau(t))},
\]
and, thanks to \ref{hypo:delay-infinity}, we deduce that
\[
\limsup_{t \to +\infty}\abs{x(t)} \leq \abs{A} \limsup_{t \to +\infty}\abs{x(t)}.
\]
Since $\abs{A} < 1$ and $\limsup_{t \to +\infty} \abs{x(t)} < +\infty$, we deduce that
\[\limsup_{t \to +\infty}\abs{x(t)} = 0,\]
yielding the conclusion.
\end{proof}

Note that convergence to $0$ in Theorem~\ref{thm:cv-to-0} holds in the sense of the pointwise value $x(t)$ converging to $0$ in the topology of $\mathbb R^d$. If we consider a continuous solution of \eqref{eq:diff-eqn-N1-time}, we also have convergence to $0$ in the norm of $\mathsf C_t^{\mathrm b}$.

\begin{corollary}
\label{Corol-1-C}
Let $A \in \mathcal M_d(\mathbb R)$ satisfy \ref{hypo:spr-A-less-1}, $\tau$ be a delay function satisfying \ref{hypo:tau-infimum} and \ref{hypo:delay-infinity}, $h$ be its associated largest delay function, and $x \in \mathcal C(\lvert -h(0), +\infty), \mathbb R^d)$ be a continuous solution of \eqref{eq:diff-eqn-N1-time} with bounded initial condition $x_0 \in \mathsf C_0^{\mathrm b}$. Then $\norm{x_t}_{\mathsf C_t^{\mathrm b}} \to 0$ as $t \to +\infty$.
\end{corollary}

\begin{proof}
By Definition~\ref{def:h} and assumption \ref{hypo:delay-infinity}, we have $t - h(t) = \inf_{s \in [t, +\infty)} \sigma_1(s) \to +\infty$ as $t \to +\infty$. Hence, by Theorem~\ref{thm:cv-to-0}, we have
\[
\norm{x_t}_{\mathsf C_t^{\mathrm b}} = \sup_{s \in \lvert -h(t), 0]} \abs{x_t(s)} = \sup_{s \in \lvert t - h(t), t]} \abs{x(s)} \to 0
\]
as $t \to +\infty$.
\end{proof}

The same idea used in the above proof also works for regulated solutions, yielding the following result.

\begin{corollary}
\label{Prop-1-C}
Let $A \in \mathcal M_d(\mathbb R)$ satisfy \ref{hypo:spr-A-less-1}, $\tau$ be a delay function satisfying \ref{hypo:tau-infimum} and \ref{hypo:delay-infinity}, $h$ be its associated largest delay function, and $x \in \mathcal G(\lvert -h(0), +\infty), \mathbb R^d)$ be a regulated solution of \eqref{eq:diff-eqn-N1-time} with bounded initial condition $x_0 \in \mathsf G_0^{\mathrm b}$. Then $\norm{x_t}_{\mathsf G_t^{\mathrm b}} \to 0$ as $t \to +\infty$.
\end{corollary}

A similar result also holds true for solutions in $L^\infty$, but one must take additional care due to the fact that solutions in $L^\infty$ satisfy \eqref{eq:diff-eqn-N1-time} only for a.e.\ $t \in \mathbb R_+$.

\begin{corollary}
\label{coro:cv-to-0-Linfty}
Let $A \in \mathcal M_d(\mathbb R)$ satisfy \ref{hypo:spr-A-less-1}, $\tau$ be a delay function satisfying \ref{hypo:tau-infimum}, \ref{hypo:null}, and \ref{hypo:delay-infinity}, $h$ be its associated largest delay function, and $x \in L^\infty_{\loc}(\lvert -h(0), +\infty),\allowbreak \mathbb R^d)$ be an $L^\infty$ solution of \eqref{eq:diff-eqn-N1-time} with bounded initial condition $x_0 \in L^\infty(\lvert -h(0), 0),\allowbreak \mathbb R^d)$. Then $\norm{x_t}_{L^\infty(\lvert -h(t), 0), \mathbb R^d)} \to 0$ as $t \to +\infty$.
\end{corollary}

\begin{proof}
Up to modifying $x_0$ on a set of Lebesgue measure zero, there exists $M > 0$ such that $\abs{x_0(t)} \leq M$ for every $t \in \lvert -h(0), 0)$. Denote by $\widetilde x\colon \lvert -h(0), +\infty) \to \mathbb R^d$ the unique solution of \eqref{eq:diff-eqn-N1-time} with initial condition $x_0$, whose existence is guaranteed by Theorem~\ref{thm:exist-general}. By Theorem~\ref{thm:cv-to-0}, we have $\widetilde x(t) \to 0$ as $t \to +\infty$. In addition, as in the proof of Corollary~\ref{Corol-1-C}, we have $t - h(t) \to +\infty$ as $t \to +\infty$, and thus $\sup_{s \in \lvert t - h(t), t)} \abs{\widetilde x(s)} \to 0$ as $t \to +\infty$. The conclusion follows since, by Theorem~\ref{thm:explicit-ae}, we have $\widetilde x = x$ a.e.\ in $\lvert -h(0), +\infty)$.
\end{proof}

To obtain convergence to $0$ of solutions, Theorem~\ref{thm:cv-to-0} requires assumptions \ref{hypo:spr-A-less-1} and \ref{hypo:delay-infinity}. We now show that \ref{hypo:spr-A-less-1} is also a necessary condition to have convergence of solutions to $0$, and the same is also true of \ref{hypo:delay-infinity} under the additional assumption that $A$ is not nilpotent.

\begin{proposition}
\label{prop:necessary-asympt}
Let $A \in \mathcal M_d(\mathbb R)$, $\tau$ be a delay function satisfying \ref{hypo:tau-infimum}, and $h$ be its associated largest delay function. Assume that, for every solution $x\colon \lvert -h(0), +\infty)\allowbreak \to \mathbb R^d$ of \eqref{eq:diff-eqn-N1-time} with bounded initial condition, we have $x(t) \to 0$ as $t \to +\infty$. Then \ref{hypo:spr-A-less-1} holds true. If, in addition, $A$ is not nilpotent, then \ref{hypo:delay-infinity} also holds true.
\end{proposition}

\begin{proof}
Let $\lambda \in \sigma(A)$ and $v \in \mathbb C^d \setminus \{0\}$ be an eigenvector of $A$ associated with $\lambda$. Consider the solution\footnote{Although we have so far only considered solutions taking values in $\mathbb R^d$, it is immediate to verify that the previous results extend with the same proofs to solutions taking values in $\mathbb C^d$.} $x\colon \lvert -h(0), +\infty) \to \mathbb C^d$ of \eqref{eq:diff-eqn-N1-time} with constant initial condition equal to $v$ and note that, due to the linearity of \eqref{eq:diff-eqn-N1-time}, the functions $t \mapsto \Real x(t)$ and $t \mapsto \Imag x(t)$ are solutions of \eqref{eq:diff-eqn-N1-time} taking values in $\mathbb R^d$ and with bounded initial conditions. Hence, by assumption, $\Real x(t) \to 0$ and $\Imag x(t) \to 0$ as $t \to +\infty$, showing that $x(t) \to 0$ as $t \to +\infty$. On the other hand, using Theorem~\ref{theorem-1}, we compute $x(t) = \lambda^{\mathbf{n}(t)} v$ for every $t \in \mathbb R_+$, hence $\lambda^{\mathbf n(t)} \to 0$ as $t \to +\infty$, implying that $\abs{\lambda} < 1$. As a consequence, $A$ satisfies \ref{hypo:spr-A-less-1}.

If, in addition, $A$ is not nilpotent, then there exists $\lambda \in \sigma(A)$ with $\lambda \neq 0$. Defining $v \in \mathbb C^d \setminus \{0\}$ and $x\colon \lvert -h(0), +\infty) \to \mathbb C^d$ from $\lambda$ as before and recalling that $x(t) \to 0$ as $t \to +\infty$ and $x(t) = \lambda^{\mathbf{n}(t)} v$ for every $t \in \mathbb R_+$, we obtain in particular that $x(\sigma_1(t)) = \lambda^{\mathbf{n}(\sigma_1(t))} v$ for every $t \in \mathbb R_+$. Indeed, this is a consequence of the previous fact if $\sigma_1(t) \geq 0$, and it is trivially verified if $\sigma_1(t) < 0$ since $\mathbf n(s) = 0$ for every $s \in \mathbb R_-^\ast$. Thus, for every $t \in \mathbb R_+$, we have $x(t) = A x(\sigma_1(t)) = \lambda^{\mathbf{n}(\sigma_1(t)) + 1} v$, showing that $\lambda^{\mathbf{n}(\sigma_1(t)) + 1} \to 0$ as $t \to +\infty$ since $v \neq 0$. As $\lambda \neq 0$, we deduce that $\mathbf{n}(\sigma_1(t)) \to +\infty$ as $t \to +\infty$, and thus, by Lemma~\ref{lemma:n-loc-bounded}, we conclude that $\sigma_1(t) \to +\infty$ as $t \to +\infty$, as required.
\end{proof}

A natural question is whether Corollary~\ref{coro:cv-to-0-Linfty} remains true if $L^\infty$ is replaced with $L^p$ for some $p \in [1, +\infty)$. Our next example shows that this is not the case: in dimension $d = 1$, for every $p \in [1, +\infty)$ and $A \in (-1, 1) \setminus \{0\}$, there exists a delay function $\tau$ satisfying \ref{hypo:tau-infimum}, \ref{hypo:radon-nik-bdd}, and \ref{hypo:delay-infinity}, such that, for every initial condition $x_0 \in L^p(\lvert -h(0), 0], \mathbb R)$, the unique solution $x \in L^p_{\loc}(\lvert -h(0), +\infty), \mathbb R)$ of \eqref{eq:diff-eqn-N1-time} does not converge to $0$ in $L^p$ norm as $t \to +\infty$.

\begin{example}\label{example-3.1}
Consider \eqref{eq:diff-eqn-N1-time} in dimension $d = 1$ and let $p \in [1, +\infty)$, $A \in (-1, 1) \setminus \{0\}$, $\alpha \in (1 - \abs{A}^{p}, 1)$, and $\tau$ be the delay function given for $t \in \mathbb R_+$ by $\tau(t) = \alpha t + 1$. Note that $\sigma_1(t) = (1 - \alpha) t - 1$ for every $t \in \mathbb R_+$ and that $\tau$ satisfies \ref{hypo:tau-infimum}, \ref{hypo:radon-nik-bdd}, and \ref{hypo:delay-infinity}. In addition, since $\sigma_1$ is nondecreasing, an immediate computation using Definition~\ref{def:h} shows that the largest delay function associated with $\tau$ is $h = \tau$.

Let $x_0 \in L^p([-1, 0), \mathbb R) \setminus \{0\}$ and $x$ be the unique $L^p$ solution of \eqref{eq:diff-eqn-N1-time} with initial condition $x_0$, which exists thanks to Theorem~\ref{thm:exist-L-p}. For every $t \geq \frac{1}{1 - \alpha}$, we compute
\begin{align*}
\norm{x_t}_{L^p([-h(t), 0), \mathbb R)}^p & = \int_{t - h(t)}^t \abs{x(s)}^p \diff s = \abs{A}^p \int_{\sigma_1(t)}^t \abs{x(s - \tau(s))}^p \diff s \\
& = \frac{\abs{A}^p}{1 - \alpha} \int_{\sigma_1(t) - \tau(\sigma_1(t))}^{\sigma_1(t)} \abs{x(s)}^p \diff s = \frac{\abs{A}^p}{1 - \alpha} \norm*{x_{\sigma_1(t)}}_{L^p([-h(\sigma_1(t)), 0), \mathbb R)}^p.
\end{align*}
Let $t_0 = 0$ and $t_n = \frac{1 + t_{n-1}}{1 - \alpha}$ for $n \in \mathbb N^\ast$. Noticing that $\sigma_1(t_n) = t_{n-1}$ and that $t_n = \sum_{k=1}^n \frac{1}{(1 - \alpha)^k}$ for every $n \in \mathbb N^\ast$, we deduce that $t_n \to +\infty$ as $n \to +\infty$ and thus
\[
\norm{x_{t_n}}_{L^p([-h(t_n), 0), \mathbb R)}^p = \left(\frac{\abs{A}^p}{1 - \alpha}\right)^n \norm{x_0}_{L^p([-1, 0), \mathbb R)}^p \to +\infty
\]
as $n \to +\infty$, since $\frac{\abs{A}^p}{1 - \alpha} > 1$ and $\norm{x_0}_{L^p([-1, 0), \mathbb R)} > 0$. Hence the conclusion of Corollary~\ref{coro:cv-to-0-Linfty} does not hold when $L^\infty$ is replaced by $L^p$.

One might wonder whether the previous divergence to $+\infty$ of the norm is related to the fact that $h(t_n) = \tau(t_n) \to +\infty$ as $n \to +\infty$, i.e., the norm $\norm{x_{t_n}}_{L^p([-h(t_n), 0), \mathbb R)}$ is being computed in intervals whose sizes tend to $+\infty$ as $n \to +\infty$. To show that this is not the case, let $\beta \in (0, 1/p)$ be such that $(1 - \alpha)^{\beta p} < \abs{A}^p$, which exists since $1 - \alpha < \abs{A}^p$, and consider the initial condition $x_0$ given by $x_0(s) = \frac{1}{\abs{s}^\beta}$, which belongs to $L^p([-1, 0), \mathbb R)$ since $\beta < 1/p$. For every $t \geq \frac{1}{1 - \alpha}$, $n \in \mathbb N$, and $a \in (0, h(t)]$, we have
\begin{align*}
\norm{x_t}_{L^p([-a, 0), \mathbb R)}^p & = \abs{A}^p \int_{t - a}^t \abs{x(s - \tau(s))}^p \diff s = \frac{\abs{A}^p}{1 - \alpha} \int_{\sigma_1(t-a)}^{\sigma_1(t)} \abs{x(s)}^p \diff s \\
& = \frac{\abs{A}^p}{1 - \alpha} \int_{\sigma_1(t) - a(1 - \alpha)}^{\sigma_1(t)} \abs{x(s)}^p \diff s = \frac{\abs{A}^p}{1 - \alpha} \norm*{x_{\sigma_1(t)}}_{L^p([-a(1 - \alpha), 0), \mathbb R)}^p,
\end{align*}
and thus, by an immediate inductive argument, we deduce that
\[
\norm{x_{t_n}}_{L^p([-1, 0), \mathbb R)}^p = \left(\frac{\abs{A}^p}{1 - \alpha}\right)^n \norm*{x_0}_{L^p([-(1 - \alpha)^n, 0), \mathbb R)}^p,
\]
where $(t_n)_{n \in \mathbb N}$ is the same sequence as before. In addition,
\[
\norm*{x_0}_{L^p([-(1 - \alpha)^n, 0), \mathbb R)}^p = \int_{-(1 - \alpha)^n}^{0} \frac{1}{\abs{s}^{\beta p}} \diff s = \frac{(1 - \alpha)^{n(1-\beta p)}}{1 - \beta p},
\]
and hence
\[
\norm{x_{t_n}}_{L^p([-1, 0), \mathbb R)}^p = \frac{1}{1 - \beta p} \left(\frac{\abs{A}}{(1 - \alpha)^{\beta}}\right)^{p n} \to +\infty
\]
as $n \to +\infty$, since $\frac{\abs{A}}{(1 - \alpha)^\beta} > 1$.
\end{example}

Our next result provides sufficient conditions for convergence of solutions to $0$ in $L^p$ when $p \in [1, +\infty)$.

\begin{theorem}
\label{Thm-1-L-P}
Let $p \in [1, +\infty)$, $A \in \mathcal M_d(\mathbb R)$, $\tau$ be a delay function satisfying \ref{hypo:tau-infimum} and \ref{hypo:radon-nik-and-rho}, and $h$ be its associated largest delay function. Let $x \in L^p_{\loc}(\lvert -h(0), +\infty),\allowbreak \mathbb R^d)$ be an $L^p$ solution of \eqref{eq:diff-eqn-N1-time} with initial condition $x_0 \in L^p(\lvert -h(0), 0), \mathbb R^d)$. Then $x \in L^p(\lvert -h(0), +\infty), \mathbb R^d)$ and $\norm{x_t}_{L^p(\lvert -h(t), 0), \mathbb R^d)} \allowbreak\to 0$ as $t \to +\infty$.
\end{theorem}

\begin{proof}
Let $\varphi$ denote the Radon--Nikodym derivative of ${\sigma_1}_{\#} \mathcal L$ with respect to $\mathcal L$. The inequality $\norm{\varphi}_{L^\infty(\mathbb R, \mathbb R)} \spr(A)^p < 1$ yields the existence of a norm $\abs{\cdot}$ in $\mathbb R^d$ for which we have $\norm{\varphi}_{L^\infty(\mathbb R, \mathbb R)} \abs{A}^p \allowbreak < 1$. We fix such a norm in the sequel.

Let $(D_n)_{n \in \mathbb N}$ be the nonincreasing sequence of sets from Definition~\ref{def:Dn-sigman} and define $D_\ast = \lvert -h(0), +\infty)$. Thanks to Lemma~\ref{lemma:D_n-sigma1}, we have $\sigma_1^{-1}(D_{n-1} \setminus D_n) = D_n \setminus D_{n+1}$ for every $n \in \mathbb N^\ast$. In addition, since $\sigma_1$ takes values in $D_\ast$ and $D_n \subset \mathbb R_+ \subset D_\ast$ for every $n \in \mathbb N^\ast$, we have $\sigma_1^{-1}(D_\ast \cap (D_{n-1} \setminus D_n)) = \sigma_1^{-1}(D_{n-1} \setminus D_n) = D_n \setminus D_{n+1} = D_\ast \cap (D_n \setminus D_{n+1})$ for every $n \in \mathbb N^\ast$.

Let $n \in \mathbb N^\ast$. Using \eqref{eq:diff-eqn-N1-time}, we have
\begin{align*}
\int_{D_\ast \cap (D_n \setminus D_{n+1})} \abs{x(t)}^p \diff t & \leq \abs{A}^p \int_{D_\ast \cap (D_n \setminus D_{n+1})} \abs{x(\sigma_1(t))}^p \diff t \displaybreak[0] \\
& = \abs{A}^p \int_{\sigma_1^{-1}(D_\ast \cap (D_{n-1} \setminus D_{n}))} \abs{x(\sigma_1(t))}^p \diff t \displaybreak[0] \\
& = \abs{A}^p \int_{D_\ast \cap (D_{n-1} \setminus D_{n})} \abs{x(t)}^p \diff {\sigma_1}_{\#} \mathcal L(t) \displaybreak[0] \\
& = \abs{A}^p \int_{D_\ast \cap (D_{n-1} \setminus D_{n})} \varphi(t) \abs{x(t)}^p \diff t \displaybreak[0] \\
& \leq \norm{\varphi}_{L^\infty(\mathbb R, \mathbb R)} \abs{A}^p \int_{D_\ast \cap (D_{n-1} \setminus D_{n})} \abs{x(t)}^p \diff t.
\end{align*}
Hence, for every $n \in \mathbb N$, we have
\begin{align*}
\int_{D_\ast \cap (D_n \setminus D_{n+1})} \abs{x(t)}^p \diff t & \leq \left(\norm{\varphi}_{L^\infty(\mathbb R, \mathbb R)} \abs{A}^p\right)^n \int_{D_\ast \cap (D_{0} \setminus D_{1})} \abs{x(t)}^p \diff t \displaybreak[0] \\
& = \left(\norm{\varphi}_{L^\infty(\mathbb R, \mathbb R)} \abs{A}^p\right)^n \norm{x_0}_{L^p(\lvert -h(0), 0), \mathbb R^d)}^p.
\end{align*}
Thus
\begin{equation*}
\int_{-h(0)}^{+\infty} \abs{x(t)}^p \diff t = \sum_{n = 0}^{\infty} \int_{D_\ast \cap (D_n \setminus D_{n+1})} \abs{x(t)}^p \diff t \leq \frac{\norm{x_0}_{L^p(\lvert -h(0), 0), \mathbb R^d)}^p}{1 - \norm{\varphi}_{L^\infty(\mathbb R, \mathbb R)} \abs{A}^p},
\end{equation*}
showing that $x \in L^p(\lvert -h(0), +\infty), \mathbb R^d)$.

Finally, since $\lim_{t \to +\infty} \sigma_1(t) = +\infty$, it follows from Definition~\ref{def:h} that $\lim_{t \to +\infty} t - h(t) = +\infty$, and thus, using the fact that $x \in L^p(\lvert -h(0), +\infty), \mathbb R^d)$, we deduce that
\[
\norm{x_t}_{L^p(\lvert -h(t), 0), \mathbb R^d)} = \norm{x(\cdot)}_{L^p(\lvert t-h(t), t), \mathbb R^d)} \leq \norm{x(\cdot)}_{L^p(\lvert t-h(t), +\infty), \mathbb R^d)} \xrightarrow[t \to +\infty]{} 0,
\]
as required.
\end{proof}

\begin{remark}
Under the assumptions of Theorem~\ref{Thm-1-L-P}, we have $\norm{\varphi}_{L^\infty(\mathbb R, \mathbb R)} \geq 1$, and in particular \ref{hypo:spr-A-less-1} is necessarily satisfied. Indeed, by \ref{hypo:delay-infinity}, there exists $T_0 \in \mathbb R_+^\ast$ such that $\sigma_1(t) \geq 0$ for every $t \in [T_0, +\infty)$. On the other hand, thanks to the definition of $\sigma_1$, we have $\sigma_1(t) < t$ for every $t \in \mathbb R_+$. In particular, for every $t \in [T_0, +\infty)$, we have $[T_0, t] \subset \sigma_1^{-1}([0, t])$. We also have
\[
\mathcal L(\sigma_1^{-1}([0, t])) = \int_0^t \diff {\sigma_1}_{\#} \mathcal L(s) = \int_0^t \varphi(s) \diff s \leq t \norm{\varphi}_{L^\infty(\mathbb R, \mathbb R)}.
\]
Hence, for every $t \in [T_0, +\infty)$, we have
\[
\norm{\varphi}_{L^\infty(\mathbb R, \mathbb R)} \geq \frac{1}{t} \mathcal L(\sigma_1^{-1}([0, t])) \geq \frac{1}{t} \mathcal L([T_0, t]) = \frac{t - T_0}{t} \xrightarrow[t \to +\infty]{} 1,
\]
showing that $\norm{\varphi}_{L^\infty(\mathbb R, \mathbb R)} \geq 1$. Together with \ref{hypo:radon-nik-and-rho}, one deduces that $\spr(A) < 1$.
\end{remark}

\begin{remark}
Recalling Remark~\ref{remk:radon-nikodym-derivative-tau}, one observes that, if the delay function $\tau$ is differentiable and there exists $\alpha \in (0, 1)$ such that $\tau^\prime(t) \leq \alpha$ for every $t \in \mathbb R_+$, then $\sigma_1$ is a diffeomorphism between $\mathbb R_+$ and its image, and the Radon--Nikodym derivative of ${\sigma_1}_{\#}\mathcal L$ with respect to $\mathcal L$ is bounded by $\frac{1}{1 - \alpha}$. In this case, \ref{hypo:radon-nik-and-rho} is satisfied as soon as $\frac{\spr(A)^p}{1 - \alpha} < 1$.
\end{remark}

\begin{remark}
Example~\ref{example-3.1} illustrates the importance of \ref{hypo:radon-nik-and-rho} in Theorem~\ref{Thm-1-L-P}. Indeed, the delay function $\tau$ from Example~\ref{example-3.1} satisfies \ref{hypo:tau-infimum}, \ref{hypo:tau-measurable}, \ref{hypo:null}, and \ref{hypo:delay-infinity}, and we also have that ${\sigma_1}_{\#} \mathcal L$ is $\sigma$-finite. Recalling Remark~\ref{remk:radon-nikodym-derivative-tau}, the Radon--Nikodym derivative $\varphi$ of ${\sigma_1}_{\#}\mathcal L$ with respect to $\mathcal L$ is the function $\varphi$ given by $\varphi(t) = \frac{1}{1 - \alpha}$ if $t \geq -1$ and $\varphi(t) = 0$ otherwise. In particular, it is bounded, but we have $\norm{\varphi}_{L^\infty(\mathbb R, \mathbb R)} \spr(A)^p = \frac{\abs{A}^p}{1 - \alpha} > 1$, and the conclusion of Theorem~\ref{Thm-1-L-P} does not hold in that example.
\end{remark}

\subsection{Exponential stability}
\label{sec:exponential}

We now turn to the analysis of the exponential stability of \eqref{eq:diff-eqn-N1-time}, i.e., whether its solutions converge exponentially to zero in a given norm. The following technical assumption will be crucial in the analysis carried out here.

\begin{hypotheses}
Let $\tau$ be a delay function.
\begin{hypolist}
\item\label{hypo:lower-bound-n} The function $\tau$ satisfies \ref{hypo:tau-infimum} and there exist $\alpha > 0$ and $\beta \in \mathbb R$ such that $\mathbf n(t) \geq \alpha t + \beta$ for every $t \in \mathbb R_+$, where $\mathbf n\colon \mathbb R \to \mathbb N$ is the function defined from $\tau$ in Corollary~\ref{coro:integer-iterations}.
\end{hypolist}
\end{hypotheses}

Hypothesis~\ref{hypo:lower-bound-n} is not easy to verify in practice, since $\mathbf n$ may be difficult to compute explicitly. However, it is satisfied in the important particular case of a bounded delay function $\tau$, as stated in the following result.

\begin{proposition}
\label{prop:tau-bounded-implies-n-lower-bounded}
Let $\tau$ be a delay function satisfying \ref{hypo:delay-bounded} and assume that $\tau$ is upper boun\-ded by some positive constant $\tau_{\max}$. Then $\tau$ satisfies \ref{hypo:lower-bound-n} with $\alpha = \frac{1}{\tau_{\max}}$ and $\beta = 0$.
\end{proposition}

\begin{proof}
Let $(D_n)_{n \in \mathbb N}$ be the sequence of sets from Definition~\ref{def:Dn-sigman}. We claim that, for every $n \in \mathbb N$, we have $D_n \supset [(n-1) \tau_{\max}, +\infty)$. Indeed, since $D_0 = \mathbb R$ and $D_1 = \mathbb R_+$, the previous property clearly holds for $n \in \{0, 1\}$. Assume now that $n \in \mathbb N$ is such that $D_n \supset [(n-1) \tau_{\max}, +\infty)$ and take $t \in [n \tau_{\max}, +\infty)$. Then $t - \tau(t) \in [(n-1) \tau_{\max}, +\infty) \subset D_n$ and, using \eqref{eq:defi-Dn}, we deduce that $t \in D_{n+1}$, yielding the desired result.

Now, for $t \in \mathbb R_+$, it follows from the definition of $\mathbf n$ that $t \notin D_{\mathbf n(t) + 1}$, implying that $t \notin [\mathbf n(t) \tau_{\max}, +\infty)$. Hence $t < \mathbf n(t) \tau_{\max}$, yielding the conclusion.
\end{proof}

Even though \ref{hypo:delay-bounded} is a sufficient condition for \ref{hypo:lower-bound-n}, it is not necessary, as shown in the following example.

\begin{example}
\label{expl:exponential-unbounded}
Let $\tau\colon \mathbb R_+ \to \mathbb R_+^\ast$ be the delay function defined by
\[
\tau(t) = \begin{dcases*}
k, & if there exists $k \in \mathbb N^\ast$ such that $t \in [2^k, 2^k + 1)$, \\
1, & otherwise.
\end{dcases*}
\]
An illustration of the graph of $\tau$ is provided in Figure~\ref{fig:tau-expl}. Note that $\tau$ satisfies \ref{hypo:tau-infimum} and \ref{hypo:delay-infinity}, but not \ref{hypo:delay-bounded}.

\begin{figure}[ht]
\centering
\begin{tikzpicture}[x=0.2cm, y=0.2cm]
\draw[black!15!white, step=0.2cm] (-0.4, -0.4) grid (67.4, 6.4);
\draw[->, thick] (-1, 0) -- (68, 0) node[right] {$t$};
\draw[->, thick] (0, -1) -- (0, 7) node[left] {$\tau$};
\foreach \i in {5, 10, ..., 65} {
  \node[below] at ({\i}, 0) {\scriptsize $\i$};
}
\node[below left] at (0, 0) {\scriptsize $0$};
\node[left] at (0, 5) {\scriptsize $5$};

\draw[blue, thick] (0, 1) -- (2, 1);
\foreach \k in {1, 2, ..., 6} {
	\draw[blue, thick] ({2^\k}, {\k}) -- ({2^\k + 1}, {\k});
	\draw[blue, thick] ({2^\k + 1}, 1) -- ({min(2^(\k + 1), 67.4)}, 1);
}
\end{tikzpicture}
\caption{Illustration of the graph of the delay function $\tau$ from Example~\ref{expl:exponential-unbounded}}
\label{fig:tau-expl}
\end{figure}
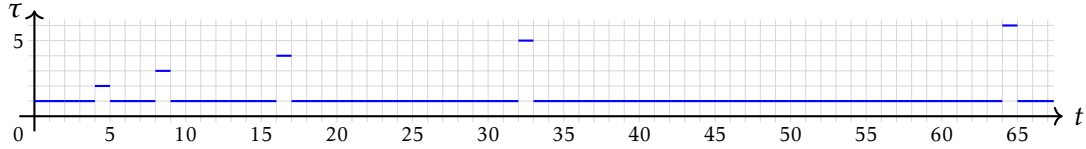

Let $\mathbf n$ be the function defined from $\tau$ as in Corollary~\ref{coro:integer-iterations} and note that, thanks to Lem\-ma~\ref{lemma:properties_n}\ref{item:n_sigma_1}, for every $t \in \mathbb R_+$, we have
\begin{equation}
\label{eq:n-recursive}
\mathbf n(t) = \begin{dcases*}
\mathbf n(t - k) + 1, & if there exists $k \in \mathbb N^\ast$ such that $t \in [2^k, 2^k + 1)$, \\
\mathbf n(t - 1) + 1, & otherwise.
\end{dcases*}
\end{equation}

We claim that
\begin{equation}
\label{eq:n-explicit}
\mathbf n(t) = \begin{dcases*}
0, & if $t < 0$, \\
1, & if $t \in [0, 1)$, \\
\floor{t} - \frac{1}{2}\left(\floor{\log_2 t} - 2\right)\left(\floor{\log_2 t} + 1\right), & if $t \geq 1$.
\end{dcases*}
\end{equation}
It follows from the definition of $\mathbf n$ that $\mathbf n(t) = 0$ for $t < 0$, and one easily deduces that $\mathbf n(t) = 1$ for $t \in [0, 1)$ using \eqref{eq:n-recursive}. To obtain the expression of $\mathbf n(t)$ for $t \geq 1$, we show by strong induction on $\ell$ that, for every $\ell \in \mathbb N^\ast$,
\begin{equation}
\label{eq:n-induction}
\mathbf n(t) = \floor{t} - \frac{1}{2}\left(\floor{\log_2 t} - 2\right)\left(\floor{\log_2 t} + 1\right) \qquad \text{ for every } t \in [\ell, \ell + 1).
\end{equation}
Indeed, for $t \in [1, 2)$, we have $\mathbf n(t) = \mathbf n(t - 1) + 1 = 2$, while, on the other hand, $\floor{\log_2 t} = 0$, and thus \eqref{eq:n-induction} is satisfied for $\ell = 1$. Assume now that $\ell_0 \in \mathbb N^\ast$ with $\ell_0 \geq 2$ is such that \eqref{eq:n-induction} is satisfied for every $\ell \in \llbracket 1, \ell_0-1\rrbracket$ and take $t \in [\ell_0, \ell_0 + 1)$. Let $k_0 = \floor{\log_2 \ell_0}$. If $\ell_0 \in \llbracket 2^{k_0} + 1, 2^{k_0 + 1} - 1\rrbracket$, then $t \in [2^{k_0} + 1, 2^{k_0 + 1})$ and thus $\floor{\log_2 t} = \floor{\log_2 (t-1)} = k_0$. Using \eqref{eq:n-recursive} and \eqref{eq:n-induction} with $\ell = \ell_0 - 1$, we have
\begin{align*}
\mathbf n(t) & = \mathbf n(t-1) + 1 = \floor{t - 1} - \frac{1}{2}\left(k_0 - 2\right)\left(k_0 + 1\right) + 1 \\
& = \floor{t} - \frac{1}{2}\left(k_0 - 2\right)\left(k_0 + 1\right),
\end{align*}
showing \eqref{eq:n-induction} for $\ell_0$. Otherwise, we have $\ell_0 = 2^{k_0}$. In this case, $t - k_0 \in [2^{k_0} - k_0, 2^{k_0} - k_0 + 1) \subset [2^{k_0 - 1}, 2^{k_0})$ since $k_0 \geq 1$, and thus $\floor{\log_2 t} = k_0$ and $\floor{\log_2 (t-k_0)} = k_0 - 1$. Hence, using \eqref{eq:n-recursive} and \eqref{eq:n-induction} with $\ell = \ell_0 - k_0 \geq 1$, we have
\begin{align*}
\mathbf n(t) & = \mathbf n(t - k_0) + 1 = \floor{t - k_0} - \frac{1}{2}\left(k_0 - 3\right)k_0 + 1 \\
 & = \floor{t - k_0} - \frac{1}{2}\left(k_0 - 2\right)\left(k_0 + 1\right) + k_0 \\
 & = \floor{t} - \frac{1}{2}\left(k_0 - 2\right)\left(k_0 + 1\right),
\end{align*}
showing \eqref{eq:n-induction} for $\ell_0$ and concluding the proof of the inductive step.

Using \eqref{eq:n-explicit}, one can now prove that
\begin{equation}
\label{eq:n-lower-bound}
\mathbf n(t) \geq \frac{t}{2} - 1 \qquad \text{ for every } t \in \mathbb R.
\end{equation}
Indeed, this is immediate for $t \in (-\infty, 1]$. For $t \geq 1$, since $\floor{t} \geq t - 1$, it suffices to prove that $\left(\floor{\log_2 t} - 2\right)\left(\floor{\log_2 t} + 1\right) \leq t$. We have $\left(\floor{\log_2 t} - 2\right)\left(\floor{\log_2 t} + 1\right) \leq \left(\floor{\log_2 t}\right)^2 - 2 \leq (\log_2 t)^2 - 2$ for every $t \geq 1$, and we can prove that $(\log_2 t)^2 - 2 \leq t$ for every $t \geq 1$ thanks to an analysis of the sign changes of the derivative of $t \mapsto t - (\log_2 t)^2 + 2$, yielding \eqref{eq:n-lower-bound}. Hence, $\tau$ satisfies \ref{hypo:lower-bound-n}.
\end{example}

Our first main result of this section is the following theorem, showing that \ref{hypo:tau-infimum}, \ref{hypo:spr-A-less-1}, and \ref{hypo:lower-bound-n} are sufficient conditions for the exponential convergence to $0$ of $\abs{x(t)}$ for solutions with bounded initial condition.

\begin{theorem}
\label{thm:exp-to-0}
Let $A \in \mathcal M_d(\mathbb R)$ satisfy \ref{hypo:spr-A-less-1}, $\tau$ be a delay function satisfying \ref{hypo:tau-infimum} and \ref{hypo:lower-bound-n}, and $h$ be its associated largest delay function. There exist positive constants $C$ and $\gamma$ such that, for every solution $x\colon \lvert -h(0), +\infty) \to \mathbb R^d$ of \eqref{eq:diff-eqn-N1-time} with a bounded initial condition $x_0\colon \lvert -h(0), 0) \to \mathbb R^d$, we have
\[\abs{x(t)} \leq C e^{-\gamma t} \sup_{s \in \lvert -h(0), 0)} \abs{x_0(s)} \qquad \text{ for every } t \in \mathbb R_+.\]
\end{theorem}

\begin{proof}
Let $\alpha$ and $\beta$ be the constants from \ref{hypo:lower-bound-n}, and let $\abs{\cdot}$ be a norm in $\mathbb R^d$ such that $\abs{A} < 1$, which exists thanks to \ref{hypo:spr-A-less-1}. Using Theorem~\ref{theorem-1}, \ref{hypo:lower-bound-n}, and the fact that $\abs{A} < 1$, we have, for every $t \in \mathbb R_+$,
\[
\abs{x(t)} \leq \abs{A}^{\mathbf n(t)} \abs{x_0(\sigma_{\mathbf n(t)} (t))} \leq \abs{A}^{\alpha t + \beta} \sup_{s \in \lvert -h(0), 0)} \abs{x_0(s)},
\]
yielding the conclusion with $C = \abs{A}^\beta > 0$ and $\gamma = -\alpha \ln\abs{A} > 0$. Note that, thanks to the equivalence of all norms in $\mathbb R^d$, the conclusion also holds true if $\abs{\cdot}$ is replaced by any other norm in $\mathbb R^d$, up to modifying the constant $C$.
\end{proof}

\begin{remark}
Hypothesis~\ref{hypo:lower-bound-n} is important in the proof of Theorem~\ref{thm:exp-to-0} in order to bound $\abs{A}^{\mathbf n(t)}$ by an exponential in $t$. However, it is not necessary: if, for instance, we replace \ref{hypo:lower-bound-n} by the more general lower bound $\mathbf n(t) \geq \alpha t + \beta - \ln p(t)$ for some $\alpha > 0$, $\beta \in \mathbb R$, and a polynomial $p$ satisfying $p(t) \geq 1$ for every $t \in \mathbb R_+$, then $\abs{A}^{\mathbf n(t)}$ can be upper bounded by $C_0 p(t)^{\nu} e^{-\gamma_0 t}$, with $C_0 = \abs{A}^\beta > 0$, $\nu = -\ln\abs{A} > 0$, and $\gamma_0 = \alpha \nu > 0$, and $C_0 p(t)^{\nu} e^{-\gamma_0 t}$ can be upper bounded by $C e^{-\gamma t}$ for any choice of $\gamma \in (0, \gamma_0)$, by choosing $C = \max_{t \in \mathbb R_+} C_0 p(t)^\nu e^{-(\gamma_0 - \gamma) t}$.
\end{remark}

As a consequence of Theorem~\ref{thm:exp-to-0}, we obtain the following criterion for the exponential convergence to $0$ of continuous solutions of \eqref{eq:diff-eqn-N1-time}.

\begin{corollary}
\label{coro:exp-cont}
Let $A \in \mathcal M_d(\mathbb R)$ satisfy \ref{hypo:spr-A-less-1}, $\tau$ be a delay function satisfying \ref{hypo:tau-infimum} and \ref{hypo:delay-bounded}, and $h$ be its associated largest delay function. Then there exist positive constants $C$ and $\gamma$ such that, for every continuous solution $x \in \mathcal C(\lvert -h(0), +\infty), \mathbb R^d)$ of \eqref{eq:diff-eqn-N1-time} with a bounded initial condition $x_0 \in \mathsf C_0^{\mathrm b}$, we have
\[\norm{x_t}_{\mathsf C_t^{\mathrm b}} \leq C e^{-\gamma t} \norm{x_0}_{\mathsf C_0^{\mathrm b}} \qquad \text{ for every } t \in \mathbb R_+.\]
\end{corollary}

\begin{proof}
Since $\tau$ satisfies \ref{hypo:tau-infimum} and \ref{hypo:delay-bounded}, it follows from Proposition~\ref{prop:tau-bounded-implies-n-lower-bounded} that \ref{hypo:lower-bound-n} also holds true. In addition, setting $\tau_{\max} = \sup_{t \in \mathbb R_+} \tau(t)$, we have from Definition~\ref{def:h} that $h(t) \leq \tau_{\max}$ for every $t \in \mathbb R_+$. Let $C_0$ and $\gamma$ denote the constants whose existence is asserted in Theorem~\ref{thm:exp-to-0} and assume, with no loss of generality, that $C_0 \geq 1$. Since $x$ is a continuous solution of \eqref{eq:diff-eqn-N1-time} with a bounded initial condition $x_0$, it follows from Theorem~\ref{thm:exp-to-0} that
\[\abs{x(t)} \leq C_0 e^{-\gamma t} \norm{x_0}_{\mathsf C_0^{\mathrm b}}\]
for every $t \in \mathbb R_+$ and, since $C_0 \geq 1$ and $\gamma > 0$, the above inequality trivially holds true also for $t \in \lvert -h(0), 0)$. Then, for every $t \in \mathbb R_+$, we have
\[
\norm{x_t}_{\mathsf C_t^{\mathrm b}} = \sup_{s \in \lvert t - h(t), t]} \abs{x(s)} \leq C_0 e^{-\gamma (t - h(t))} \norm{x_0}_{\mathsf C_0^{\mathrm b}} \leq C e^{-\gamma t} \norm{x_0}_{\mathsf C_0^{\mathrm b}},
\]
with $C = C_0 e^{\gamma \tau_{\max}}$.
\end{proof}

Assumption~\ref{hypo:delay-bounded} is useful in Corollary~\ref{coro:exp-cont} since it ensures that the function $h$ from Definition~\ref{def:h} is upper bounded by $\tau_{\max}$, and hence, in the above proof, one can bound $C_0 e^{-\gamma (t - h(t))}$ by $C_0 e^{\gamma \tau_{\max}} e^{-\gamma t}$. However, \ref{hypo:delay-bounded} is not necessary for the conclusions of Corollary~\ref{coro:exp-cont} to hold true, as illustrated in the following example.

\begin{example}
Let $\tau$ be the delay function from Example~\ref{expl:exponential-unbounded}, $h$ be its associated largest delay function, and $A \in \mathcal M_d(\mathbb R)$ be a matrix satisfying \ref{hypo:spr-A-less-1}. Note that $A$ and $\tau$ satisfy the assumptions of Theorem~\ref{thm:exp-to-0}, but not those of Corollary~\ref{coro:exp-cont}. We will prove, however, that the conclusion of Corollary~\ref{coro:exp-cont} remains true in this case.

Note that, by the definition of $\tau$, we have, for every $t \in \mathbb R_+$,
\begin{equation}
\label{eq:estim-inf}
\inf_{s \in [t, +\infty)} (s - \tau(s)) = \min\left(\inf_{s \in [t, +\infty) \setminus \bigcup_{k=1}^\infty [2^k, 2^k+1)} (s - 1), \inf_{\substack{(s, k) \in [t, +\infty) \times \mathbb N^\ast \\ \floor{s} = 2^k}} (s - k)\right).
\end{equation}
We clearly have
\begin{equation}
\label{eq:estim-inf-1}
\inf_{s \in [t, +\infty) \setminus \bigcup_{k=1}^\infty [2^k, 2^k+1)} (s - 1) \geq t - 1.
\end{equation}
On the other hand, for every $(s, k) \in [t, +\infty) \times \mathbb N^\ast$ with $\floor{s} = 2^k$, we have $s \geq 2$ and $k = \log_2 \floor{s} < \log_2 (s+1)$, and thus $s - k > s - \log_2 (s+1) \geq 2 - \log_2 3 > 0$ since $z \mapsto z - \log_2 (z+1)$ is increasing in $[1, +\infty)$. We then have
\begin{equation}
\label{eq:estim-inf-2}
\inf_{\substack{(s, k) \in [t, +\infty) \times \mathbb N^\ast \\ \floor{s} = 2^k}} (s - k) \geq t - \log_2 (t+1).
\end{equation}
Indeed, if $t \geq 1$, this follows from the fact that $z \mapsto z - \log_2 (z+1)$ is increasing in $[1, +\infty)$ while, if $t \in [0, 1)$, the inequality is trivial since $t - \log_2 (t+1) \leq 0$ but the left-hand side of \eqref{eq:estim-inf-2} is positive.

Inserting \eqref{eq:estim-inf-1} and \eqref{eq:estim-inf-2} into \eqref{eq:estim-inf} and using Definition~\ref{def:h}, we deduce that
\[
h(t) \leq 1 + \log_2(t+1) \qquad \text{ for every } t \in \mathbb R_+.
\]
Hence, applying Theorem~\ref{thm:exp-to-0} and proceeding as in the proof of Corollary~\ref{coro:exp-cont}, we deduce that, for every continuous solution $x \in \mathcal C(\lvert -h(0), +\infty), \mathbb R^d)$ of \eqref{eq:diff-eqn-N1-time} with a bounded initial condition $x_0 \in \mathsf C_0^{\mathrm b}$, we have
\[\norm{x_t}_{\mathsf C_t^{\mathrm b}} \leq C_0 e^{-\gamma_0 (t - h(t))} \norm{x_0}_{\mathsf C_0^{\mathrm b}} \qquad \text{ for every } t \in \mathbb R_+,\]
where $C_0$ and $\gamma_0$ are the constants from Theorem~\ref{thm:exp-to-0}. For every $t \in \mathbb R_+$, we have
\[C_0 e^{-\gamma_0 (t - h(t))} \leq C_0 e^{\gamma_0} (t+1)^{\frac{1}{\ln 2}} e^{-\gamma_0 t} \leq C e^{-\gamma t},\]
with an arbitrary $\gamma \in (0, \gamma_0)$ and $C > 0$ given by
\[C = C_0 e^{\gamma_0} \max_{t \in \mathbb R_+} (t+1)^{\frac{1}{\ln 2}} e^{-(\gamma_0 - \gamma) t}.\]
Hence, the conclusion of Corollary~\ref{coro:exp-cont} holds true.
\end{example}

The proof of Corollary~\ref{coro:exp-cont} can be straightforwardly adapted to regulated solutions, yielding the following result.

\begin{corollary}
\label{coro:exp-G}
Let $A \in \mathcal M_d(\mathbb R)$ satisfy \ref{hypo:spr-A-less-1}, $\tau$ be a delay function satisfying \ref{hypo:tau-infimum} and \ref{hypo:delay-bounded}, and $h$ be its associated largest delay function. Then there exist positive constants $C$ and $\gamma$ such that, for every regulated solution $x \in \mathcal G(\lvert -h(0), +\infty), \mathbb R^d)$ of \eqref{eq:diff-eqn-N1-time} with a bounded initial condition $x_0 \in \mathsf G_0^{\mathrm b}$, we have
\[\norm{x_t}_{\mathsf G_t^{\mathrm b}} \leq C e^{-\gamma t} \norm{x_0}_{\mathsf G_0^{\mathrm b}} \qquad \text{ for every } t \in \mathbb R_+.\]
\end{corollary}

We also have the following counterpart of Corollaries~\ref{coro:exp-cont} and \ref{coro:exp-G} for solutions in $L^\infty$.

\begin{corollary}
\label{coro:exp-Linfty}
Let $A \in \mathcal M_d(\mathbb R)$ satisfy \ref{hypo:spr-A-less-1}, $\tau$ be a delay function satisfying \ref{hypo:tau-infimum}, \ref{hypo:null}, and \ref{hypo:delay-bounded}, and $h$ be its associated largest delay function. Then there exist positive constants $C$ and $\gamma$ such that, for every $L^\infty$ solution $x \in L^\infty_{\loc}(\lvert -h(0), +\infty), \mathbb R^d)$ of \eqref{eq:diff-eqn-N1-time} with an initial condition $x_0 \in L^\infty(\lvert -h(0), 0), \mathbb R^d)$, we have
\[\norm{x_t}_{L^\infty(\lvert -h(t), 0), \mathbb R^d)} \leq C e^{-\gamma t} \norm{x_0}_{L^\infty(\lvert -h(0), 0), \mathbb R^d)} \qquad \text{ for every } t \in \mathbb R_+.\]
\end{corollary}

The proof of Corollary~\ref{coro:exp-Linfty} is slightly more subtle than that of Corollaries~\ref{coro:exp-cont} and \ref{coro:exp-G}, since $L^\infty$ solutions only satisfy \eqref{eq:diff-eqn-N1-time} for almost every $t \in \mathbb R_+$. It is omitted here since it can be done through an easy adaptation of the argument already provided in the proof of Corollary~\ref{coro:cv-to-0-Linfty}: one can select a measurable bounded function $\widetilde x_0$ in the equivalence class of $x_0$ with respect to almost everywhere equality and consider the solution $\widetilde x$ of \eqref{eq:diff-eqn-N1-time} with initial condition $\widetilde x_0$, to which one can apply Theorem~\ref{thm:exp-to-0}. One obtains the exponential convergence to $0$ of $\norm{\widetilde x_t}_{L^\infty(\lvert -h(t), 0), \mathbb R^d)}$ by proceeding as in the proof of Corollary~\ref{coro:exp-cont}, and the conclusion of Corollary~\ref{coro:exp-Linfty} follows since $x = \widetilde x$ almost everywhere thanks to Theorem~\ref{thm:explicit-ae}.

Even though \ref{hypo:tau-infimum}, \ref{hypo:spr-A-less-1}, and \ref{hypo:delay-bounded} are sufficient to obtain exponential stability of \eqref{eq:diff-eqn-N1-time} for continuous and regulated solutions, and also for solutions in $L^\infty$ if one also assumes \ref{hypo:null}, they are not sufficient to ensure the exponential stability of \eqref{eq:diff-eqn-N1-time} in $L^p$ for $p \in [1, +\infty)$, as illustrated by the following example.

\begin{example}
\label{expl:not-exponentially-stable-Lp}
Consider \eqref{eq:diff-eqn-N1-time} in dimension $d = 1$ and let $p \in [1, +\infty)$, $A \in (-1, 1) \setminus \{0\}$, and $\alpha \in (1 - \abs{A}^{p}, 1)$. Let $(t_n)_{n \in \mathbb N}$ be the sequence given by $t_n = \frac{n}{1 - \alpha}$ for $n \in \mathbb N$. Let $\tau$ be the delay function defined by $\tau(t) = \alpha (t - t_n) + 1$ if $t \in [t_n, t_{n+1})$ for some $n \in \mathbb N$. It is easy to verify that $\tau$ is well-defined, satisfies \ref{hypo:tau-infimum}, \ref{hypo:radon-nik-bdd}, and \ref{hypo:delay-bounded}, and in addition $\tau(t) \in [1, \frac{1}{1 - \alpha})$ for every $t \in \mathbb R_+$. We denote by $h$ the largest delay function associated with $\tau$.

For $n \in \mathbb N$ and $t \in [t_n, t_{n+1})$, we have $\sigma_1(t) = (1 - \alpha)(t - t_n) + t_n - 1$. Thus, $\sigma_1$ is increasing, since it is affine and increasing in any interval of the form $[t_n, t_{n+1})$, $n \in \mathbb N$, and its jump at $t_{n+1}$ is $\sigma_1(t_{n+1}) - \sigma_1(t_{n+1}^-) = t_{n+1} - 1 - t_n = \frac{1}{1 - \alpha} - 1 > 0$. In particular, it follows from Definition~\ref{def:h} that $h = \tau$. Using the above expression of $\sigma_1$, we also obtain easily that, for every $n \in \mathbb N$ and $t \in \mathbb R_+$, we have $t \in [t_n, t_{n+1})$ if and only if $\sigma_1(t) \in (t_{n-1}, t_n)$, with the convention $t_{-1} = -\infty$.

We claim that, for every $L^p$ solution $x$ of \eqref{eq:diff-eqn-N1-time}, $n \in \mathbb N$, and $k \in \llbracket 0, n\rrbracket$, we have
\begin{equation}
\label{eq:expl-Lp}
\int_{\sigma_1(t_n)}^{t_n} \abs{x(s)}^p \diff s = \left(\frac{\abs{A}^p}{1 - \alpha}\right)^k \int_{t_{n - k} - (1 - \alpha)^k}^{t_{n-k}} \abs{x(s)}^p \diff s.
\end{equation}
Indeed, this is true for $k = 0$ since $\sigma_1(t_n) = t_n - 1$. Assuming that $k \in \llbracket 0, n-1\rrbracket$ is such that \eqref{eq:expl-Lp} holds true, we remark that $t_{n-k} - (1 - \alpha)^k \geq t_{n-k} - \frac{1}{1 - \alpha} = t_{n - k - 1}$, and thus, for every $s \in (t_{n - k} - (1 - \alpha)^k, t_{n-k})$, we have in particular $s \in (t_{n-k-1}, t_{n-k})$, and thus $\sigma_1(s) = (1 - \alpha)(s - t_{n-k-1}) + t_{n-k-1} - 1$. Using \eqref{eq:diff-eqn-N1-time} and performing a change of variables, we obtain that
\begin{align*}
\int_{t_{n - k} - (1 - \alpha)^k}^{t_{n-k}} \abs{x(s)}^p \diff s & = \abs{A}^p \int_{t_{n - k} - (1 - \alpha)^k}^{t_{n-k}} \abs{x((1 - \alpha)(s - t_{n-k-1}) + t_{n-k-1} - 1)}^p \diff s \\
& = \frac{\abs{A}^p}{1 - \alpha} \int_{t_{n - k - 1} - (1 - \alpha)^{k+1}}^{t_{n-k-1}} \abs{x(s)}^p \diff s,
\end{align*}
and thus, inserting into \eqref{eq:expl-Lp}, we deduce that \eqref{eq:expl-Lp} holds with $k$ replaced by $k+1$, yielding the conclusion.

As in Example~\ref{example-3.1}, take $\beta \in (0, 1/p)$ such that $(1 - \alpha)^{\beta p} < \abs{A}^p$ and let $x_0\colon [-1, 0) \to \mathbb R$ be given for $s \in [-1, 0)$ by $x_0(s) = \frac{1}{\abs{s}^\beta}$. Thanks to the choice of $\beta$, we have $x_0 \in L^p([-1, 0), \mathbb R)$. Let $x$ be the unique $L^p$ solution of \eqref{eq:diff-eqn-N1-time} with initial condition $x_0$, which exists thanks to Theorem~\ref{thm:exist-L-p}. Given $n \in \mathbb N$, we apply \eqref{eq:expl-Lp} to $x$ with $k = n$, obtaining that
\begin{align*}
\norm{x_{t_n}}_{L^p([-h(t_n), 0), \mathbb R)}^p & = \int_{\sigma_1(t_n)}^{t_n} \abs{x(s)}^p \diff s = \left(\frac{\abs{A}^p}{1 - \alpha}\right)^n \int_{-(1 - \alpha)^n}^{0} \abs{x_0(s)}^p \diff s \\
& = \frac{1}{1 - \beta p} \left(\frac{\abs{A}^p}{(1 - \alpha)^{\beta p}}\right)^n \xrightarrow[n \to +\infty]{} +\infty.
\end{align*}
This shows that $\norm{x_t}_{L^p([-h(t), 0), \mathbb R)}^p$ cannot converge exponentially to $0$ as $t \to +\infty$.
\end{example}

Our last result in this section provides sufficient conditions for exponential stability in $L^p$.

\begin{theorem}
\label{thm:exp-Lp}
Let $p \in [1, +\infty)$, $A \in \mathcal M_d(\mathbb R)$, $\tau$ be a delay function satisfying \ref{hypo:tau-infimum}, \ref{hypo:radon-nik-and-rho}, and \ref{hypo:delay-bounded}, and $h$ be its associated largest delay function. Then there exist positive constants $C$ and $\gamma$ such that, for every $L^p$ solution $x \in L^p_{\loc}(\lvert -h(0), +\infty),\allowbreak \mathbb R^d)$ with initial condition $x_0 \in L^p(\lvert -h(0), 0), \mathbb R^d)$, we have $x \in L^p(\lvert -h(0), +\infty), \mathbb R^d)$ and
\begin{equation}
\label{eq:estim-exp-Lp}
\norm{x_t}_{L^p(\lvert -h(t), 0), \mathbb R^d)} \leq C e^{-\gamma t} \norm{x_0}_{L^p(\lvert -h(0), 0), \mathbb R^d)} \qquad \text{ for every } t \in \mathbb R_+.
\end{equation}
\end{theorem}

\begin{proof}
Note that the assumptions of Theorem~\ref{Thm-1-L-P} are satisfied and thus, in particular, the conclusion that $x \in L^p(\lvert -h(0), +\infty), \mathbb R^d)$ follows from that theorem. Note also that, thanks to \ref{hypo:delay-bounded}, it follows from Definition~\ref{def:h} that $h(t) \leq \tau_{\max}$ for every $t \in \mathbb R_+$, where $\tau_{\max} = \sup_{t \in \mathbb R_+} \tau(t)$.

Proceeding as in the proof of Theorem~\ref{Thm-1-L-P} and using its notations, we obtain the inequality
\[
\int_{D_\ast \cap (D_n \setminus D_{n+1})} \abs{x(t)}^p \diff t \leq \left(\norm{\varphi}_{L^\infty(\mathbb R, \mathbb R)} \abs{A}^p\right)^n \norm{x_0}_{L^p(\lvert -h(0), 0), \mathbb R^d)}^p
\]
for every $n \in \mathbb N$.

Note that the inequality in \eqref{eq:estim-exp-Lp} is trivially true for $t = 0$ if one chooses $C \geq 1$. Take $t > 0$ and set $I_t = (t - h(t), t)$. Using Proposition~\ref{prop:tau-bounded-implies-n-lower-bounded}, we have $\mathbf n(s) \geq \frac{s}{\tau_{\max}} \geq \frac{t - h(t)}{\tau_{\max}} \geq \frac{t}{\tau_{\max}} - 1$ for every $s \in I_t$ and thus, in particular, $I_t \subset D_{N_t}$, where $N_t = \ceil*{\frac{t}{\tau_{\max}}} - 1 \in \mathbb N$. Hence
\begin{align*}
\norm{x_t}_{L^p(\lvert -h(t), 0), \mathbb R^d)}^p & = \int_{I_t} \abs{x(s)}^p \diff s = \sum_{n = N_t}^{\infty} \int_{I_t \cap (D_n \setminus D_{n+1})} \abs{x(s)}^p \diff s \\
& \leq \sum_{n = N_t}^{\infty} \int_{D_\ast \cap (D_n \setminus D_{n+1})} \abs{x(s)}^p \diff s \\
& \leq \sum_{n = N_t}^{\infty} \left(\norm{\varphi}_{L^\infty(\mathbb R, \mathbb R)} \abs{A}^p\right)^n \norm{x_0}_{L^p(\lvert -h(0), 0), \mathbb R^d)}^p \\
& = \frac{\left(\norm{\varphi}_{L^\infty(\mathbb R, \mathbb R)} \abs{A}^p\right)^{N_t}}{1 - \norm{\varphi}_{L^\infty(\mathbb R, \mathbb R)} \abs{A}^p} \norm{x_0}_{L^p(\lvert -h(0), 0), \mathbb R^d)}^p,
\end{align*}
and, since $\norm{\varphi}_{L^\infty(\mathbb R, \mathbb R)} \abs{A}^p < 1$ and $N_t \geq \frac{t}{\tau_{\max}} - 1$, we obtain the conclusion with
\[C = \left[\norm{\varphi}_{L^\infty(\mathbb R, \mathbb R)} \abs{A}^p (1 - \norm{\varphi}_{L^\infty(\mathbb R, \mathbb R)} \abs{A}^p)\right]^{-1} \geq 1 \quad \text{ and } \quad \gamma = -\frac{\ln \left(\norm{\varphi}_{L^\infty(\mathbb R, \mathbb R)} \abs{A}^p\right)}{\tau_{\max}} > 0.\qedhere\]
\end{proof}

\section{Applications}\label{sec:appli-varying-delay}

In this section, we describe some applications of the previous results to other contexts, by exploring relations between the difference equation \eqref{eq:diff-eqn-N1-time} and other classes of systems.

\subsection{Difference equation with a single state-\texorpdfstring{de\-pen\-dent}{dependent} delay}
\label{sec:appli-varying-delay.1}

Consider the equation
\begin{equation}\label{equation-SD-1}
x(t) = A x(t - \tau(t, x_t)),
\end{equation}
where $A \in \mathcal M_d(\mathbb R)$ and $\tau(t, x_t)$ is a delay that depends both on the time $t$ and the state $x_t$. Equations with state-dependent delays have a wide range of applications, including in electrodynamics, mechanics, control theory, neural networks, and biology, and their analysis has attracted much research effort in recent years \cite{Hartung2006Functional, Driver1984Neutral, Hernandez2023Explicit, Hartung2021Differentiability, Henriquez2024Existence, Hernandez2020Well}. Even important elementary questions, such as existence and uniqueness of solutions, are highly nontrivial for systems with state-dependent delays, and much work remains to be done for other questions such as stability analysis, with some recent results provided, for instance, in \cite{Xu2021Exponential, BekiarisLiberis2013Compensation, Li2019Algebraic, Li2020Lyapunov, Fiter2013Stability}.

Here, we are interested just in illustrating that some stability results for \eqref{equation-SD-1} can be deduced from those for \eqref{eq:diff-eqn-N1-time} obtained in Section~\ref{sec:stability}. Questions on existence and uniqueness of solutions of \eqref{equation-SD-1}, even though interesting and challenging, and sharper stability results for \eqref{equation-SD-1}, are out of the scope of the present section.

The result we provide here on \eqref{equation-SD-1} is the following one, providing sufficient conditions for the exponential convergence to zero of solutions of \eqref{equation-SD-1}.

\begin{theorem}
\label{thm:state-dep}
Let $A \in \mathcal M_d(\mathbb R)$, $I \subset \mathbb R_+^\ast$ be a compact interval, $\tau_{\min} = \inf I > 0$, $\tau_{\max} = \sup I < +\infty$, and $\tau\colon \mathbb R_+ \times \mathsf L \to I$, where $\mathsf L = \mathcal C([-\tau_{\max}, 0], \mathbb R^d)$ (respectively, $\mathsf L = \mathcal G([-\tau_{\max}, 0], \mathbb R^d)$). If $A$ satisfies \ref{hypo:spr-A-less-1}, then there exist positive constants $C$ and $\gamma$ such that, for every continuous (respectively, regulated) solution $x$ of \eqref{equation-SD-1}, denoting by $x_0$ its initial condition, we have
\[
\norm{x_t}_{\mathsf L} \leq C e^{-\gamma t} \norm{x_0}_{\mathsf L} \qquad \text{ for every } t \in \mathbb R_+.
\]
\end{theorem}

\begin{proof}
Let $x$ be a continuous (respectively, regulated) solution of \eqref{equation-SD-1} and consider the function $\tau_\ast\colon \mathbb R_+ \to I$ given by $\tau_\ast(t) = \tau(t, x_t)$. Since $\tau_\ast$ takes values in the compact set $I \subset \mathbb R_+^\ast$, it clearly satisfies assumptions \ref{hypo:tau-infimum} and \ref{hypo:delay-bounded}. Let $h_\ast$ denote the largest delay function associated with $\tau_\ast$. By Theorem~\ref{thm:exp-to-0}, there exist positive constants $C_0$ and $\gamma$ such that
\begin{equation}
\label{eq:estim-exp-SD-proof}
\abs{x(t)} \leq C_0 e^{-\gamma t} \norm{x_0}_{\mathsf L} \qquad \text{ for every } t \in \mathbb R_+.
\end{equation}
In addition, the proof of Theorem~\ref{thm:exp-to-0}, together with the statement of Proposition~\ref{prop:tau-bounded-implies-n-lower-bounded}, shows that the constants $C_0$ and $\gamma$ can be chosen to depend on $\tau_\ast$ only through $\tau_{\max}$, and, in particular, they are independent of the solution $x$. We increase if needed the constant $C_0$ in order to have $C_0 \geq 1$, and hence the inequality in \eqref{eq:estim-exp-SD-proof} trivially holds true for all $t \in [-\tau_{\max}, 0)$ as well. We then have
\[
\norm{x_t}_{\mathsf L} = \sup_{s \in [t - \tau_{\max}, t]} \abs{x(s)} \leq \sup_{s \in [t - \tau_{\max}, t]} C_0 e^{-\gamma s} \norm{x_0}_{\mathsf L} \leq C e^{-\gamma t} \norm{x_0}_{\mathsf L},
\]
with $C = C_0 e^{\gamma \tau_{\max}}$.
\end{proof}

\begin{remark}
Theorem~\ref{thm:state-dep} addresses only the case of continuous or regulated solutions of \eqref{equation-SD-1}, and its strategy relies on fixing first a solution $x$ and then considering the time-dependent delay $\tau_\ast(t) = \tau(t, x_t)$, to which the results of Section~\ref{sec:stability} can be applied. One may follow a similar strategy for $L^p$ solutions of \eqref{equation-SD-1}, but the major difficulty is that the stability results from Section~\ref{sec:stability} concerning $L^p$ solutions all require the time-dependent delay to satisfy \ref{hypo:null}, and it is not straightforward to verify whether $\tau_\ast$ satisfies it or not for a given solution $x$ of \eqref{equation-SD-1}. Even though this could possibly be verified in some simple examples, providing general assumptions on the state-dependent delay $\tau$ ensuring that $\tau_\ast$ satisfies \ref{hypo:null} for every solution $x$ of \eqref{equation-SD-1} is a difficult and interesting question, which is out of the scope of the present work.
\end{remark}

\subsection{Transport equation}
\label{sec:appli-varying-delay.2}

Consider the transport equation
\begin{equation}
\label{eq:transport}
\left\{
\begin{aligned}
& \partial_t u(t, x) + \lambda(t, x) \partial_x u(t, x) = 0, & \qquad & t > 0,\, x \in (0, 1), \\
& u(t, 0) = A u(t, 1), & & t \geq 0, \\
& u(0, x) = u_0(x), & & x \in [0, 1],
\end{aligned}
\right.
\end{equation}
where $u\colon \mathbb R_+ \times [0, 1] \to \mathbb R^d$ is the unknown function, $u_0\colon [0, 1] \to \mathbb R^d$ is given, $A \in \mathcal M_d(\mathbb R)$, and $\lambda\colon \mathbb R_+ \times [0, 1] \to \mathbb R_+^\ast$. Let us first recall what is meant by a classical solution of \eqref{eq:transport}.

\begin{definition}
A function $u \in \mathcal C^0(\mathbb R_+ \times [0, 1], \mathbb R^d) \cap \mathcal C^1(\mathbb R_+^\ast \times (0, 1), \mathbb R^d)$ is said to be a \emph{classical solution} of \eqref{eq:transport} with \emph{initial condition} $u_0$ if all three equalities of \eqref{eq:transport} are satisfied in their respective domains.
\end{definition}

Note that a necessary condition for the existence of a classical solution is that $u_0$ is continuous with $u_0(0) = A u_0(1)$.

We now want to relate classical solutions of \eqref{eq:transport} with solutions of a certain difference equation of the form \eqref{eq:diff-eqn-N1-time}. For that purpose, we consider for simplicity the following assumption on $\lambda$.

\clearpage

\begin{hypotheses}\mbox{}
\begin{hypolist}
\item\label{hypo:lambda} There exist positive constants $\lambda_{\min}$, $\lambda_{\max}$, and $L$, with $\lambda_{\min} \leq \lambda_{\max}$, such that $\lambda \in \mathcal C^0(\mathbb R^2, [\lambda_{\min}, \lambda_{\max}])$, $\lambda$ is differentiable with respect to its second variable $x$, and $\partial_x \lambda \in \mathcal C^0(\mathbb R^2, [-L, L])$.
\end{hypolist}
\end{hypotheses}

Recall that, thanks to the classical theory of ordinary differential equations, under assumption \ref{hypo:lambda}, the flow of the differential equation $x^\prime(t) = \lambda(t, x(t))$ is a well-defined $\mathcal C^1$ function. More precisely, we have the following result.

\begin{lemma}
\label{lem:flow}
Let $\lambda$ be a function satisfying \ref{hypo:lambda}. Then there exists a unique function $\Phi \in \mathcal C^1(\mathbb R^3, \mathbb R)$ satisfying the following properties.
\begin{enumerate}
\item\label{item:diff-eqn-flow} $\partial_t \Phi(t, t_0, x_0) = \lambda(t, \Phi(t, t_0, x_0))$ for every $(t, t_0, x_0) \in \mathbb R^3$.
\item $\Phi(t_0, t_0, x_0) = x_0$ for every $(t_0, x_0) \in \mathbb R^2$.
\end{enumerate}
In addition, $\Phi$ also satisfies the following properties.
\begin{enumerate}[resume]
\item\label{item:semigroup} $\Phi(t_2, t_1, \Phi(t_1, t_0, x_0)) = \Phi(t_2, t_0, x_0)$ for every $(t_0, t_1, t_2, x_0) \in \mathbb R^4$.
\item\label{item:flow-reverse} For every $(t_0, t_1, x_0, x_1) \in \mathbb R^4$, we have $\Phi(t_1, t_0, x_0) = x_1$ if and only if $\Phi(t_0, t_1,\allowbreak x_1) = x_0$.
\item\label{item:deriv-x0} $\partial_{x_0} \Phi(t, t_0, x_0) = \exp\left(\int_{t_0}^t \partial_x \lambda(s, \Phi(s, t_0, x_0)) \diff s\right)$ for every $(t, t_0, x_0) \in \mathbb R^3$.
\item\label{item:deriv-t0} $\partial_{t_0} \Phi(t, t_0, x_0) = -\lambda(t_0, x_0) \partial_{x_0} \Phi(t, t_0, x_0)$ for every $(t, t_0, x_0) \in \mathbb R^3$.
\item\label{item:flow-increasing} For every $(t_0, x_0) \in \mathbb R^2$, the function $t \mapsto \Phi(t, t_0, x_0)$ is increasing and, for every $t \in \mathbb R$,
\[
\lambda_{\min} \abs{t - t_0} \leq \abs{\Phi(t, t_0, x_0) - x_0} \leq \lambda_{\max} \abs{t - t_0}.
\]
\end{enumerate}
\end{lemma}

The function $\Phi$ from the statement of Lemma~\ref{lem:flow} is referred to in the sequel as the \emph{flow} of $\lambda$. The proof of Lemma~\ref{lem:flow} can be obtained as an immediate consequence of classical results for ordinary differential equations, such as those from \cite[Chapter~1]{Hale1969Ordinary}, and is thus omitted here.

In order to relate solutions of \eqref{eq:transport} with solutions of a difference equation of the form \eqref{eq:diff-eqn-N1-time}, we need the following technical result.

\begin{lemma}
\label{lem:lambda}
Let $\lambda$ satisfy \ref{hypo:lambda} and $\Phi$ be the flow of $\lambda$. Then the following assertions hold true.
\begin{enumerate}
\item\label{item:T0} There exists a unique $T_0 > 0$ such that $\Phi(-T_0, 0, 1) = 0$.

\item\label{item:X0} The function $X_0 \in \mathcal C^1(\mathbb R^2, \mathbb R)$ defined for $(t, x) \in \mathbb R^2$ by $X_0(t, x) = \Phi(0, t, x)$ satisfies the following properties.
\begin{enumerate}[label={\textnormal{(\roman*)}}, ref={\theenumi-(\roman*)}, leftmargin=*]
\item\label{item:X0-decr} $\partial_t X_0(t, x) < 0$ for every $(t, x) \in \mathbb R^2$.
\item For every $t \in [-T_0, 0]$, we have $X_0(t, 0) \in [0, 1]$.
\end{enumerate}

\item\label{item:R} There exists a function $R \in \mathcal C^1(\mathbb R^2, \mathbb R)$ such that $\Phi(R(t, x), t, x) = 0$ for every $(t, x) \in \mathbb R^2$ and, in addition, $R$ satisfies the following properties.
\begin{enumerate}[label={\textnormal{(\roman*)}}, ref={\theenumi-(\roman*)}, leftmargin=*]
\item\label{item:pde-R} For every $(t, x) \in \mathbb R^2$, we have $\partial_t R(t, x) + \lambda(t, x) \partial_x R(t, x) = 0$.
\item\label{item:R-incr-t} For every $x \in \mathbb R$, $t \mapsto R(t, x)$ is increasing.
\item\label{item:R-decr-x} For every $t \in \mathbb R$, $x \mapsto R(t, x)$ is decreasing, with $R(t, 0) = t$. In particular, $R(t, x) < t$ for every $x > 0$.
\item\label{item:R-bounds} For every $(t, x) \in \mathbb R^2$, we have
\[
\frac{\abs{x}}{\lambda_{\max}} \leq \abs{t - R(t, x)} \leq \frac{\abs{x}}{\lambda_{\min}}.
\]
\item\label{item:range-R} For every $(t, x) \in \mathbb R_+ \times [0, 1]$, we have $R(t, x) \in [-T_0, +\infty)$.
\item\label{item:range-R0} For every $x \in [0, 1]$, we have $R(0, x) \in [-T_0, 0]$.
\item\label{item:R-X0} For every $t \in \mathbb R_+$, if $R(t, 1) < 0$, then $X_0(t, 1) \in (0, 1]$.
\item\label{item:X0_1_0} For every $t \in \mathbb R$, we have $X_0(t, 1) = X_0(R(t, 1), 0)$.
\end{enumerate}

\item\label{item:tau} The function $\tau\colon \mathbb R \to \mathbb R$ defined for $t \in \mathbb R$ by $\tau(t) = t - R(t, 1)$ satisfies the following properties.
\begin{enumerate}[label={\textnormal{(\roman*)}}, ref={\theenumi-(\roman*)}, leftmargin=*]
\item\label{item:tau-bounds} $\tau \in \mathcal C^1\left(\mathbb R, \left[\frac{1}{\lambda_{\max}}, \frac{1}{\lambda_{\min}}\right]\right)$.
\item\label{item:tau-prime} There exists a constant $\alpha \in (0, 1)$ such that $\tau^\prime(t) \leq \alpha$ for every $t \in \mathbb R$.
\item\label{item:delay-incr} The function $t \mapsto t - \tau(t)$ is increasing.
\item\label{item:delay-lower} We have $t - \tau(t) \geq -T_0$ for every $t \geq 0$.
\item\label{item:diffeo} For every $t \geq 0$, $R(t, \cdot)$ is a diffeomorphism between $(0, 1)$ and $(t - \tau(t), t)$, with inverse given by $s \mapsto \Phi(t, s, 0)$. In addition, there exist positive constants $\beta_0, \beta_1$ such that, for every $(t, x) \in \mathbb R_+ \times (0, 1)$, we have $\beta_0 \leq \abs{\partial_x R(t, x)} \leq \beta_1$.
\end{enumerate}
\end{enumerate}
\end{lemma}

\begin{proof}
Item \ref{item:T0} is an immediate consequence of Lemma~\ref{lem:flow}\ref{item:flow-increasing}, the intermediate value theorem, and the facts that $\Phi(0, 0, 1) = 1$ and $\lim_{t \to -\infty}\Phi(t, 0, 1) = -\infty$.

Note that, by Lemma~\ref{lem:flow}\ref{item:deriv-x0} and \ref{item:deriv-t0}, for every $(t, x) \in \mathbb R^2$, we have
\[\partial_t X_0(t, x) = -\lambda(t, x) \partial_{x_0} \Phi(0, t, x) < 0.\]
This also shows that $t \mapsto X_0(t, x)$ is decreasing for every $x \in \mathbb R$. Since $X_0(0, 0) = 0$ and $X_0(-T_0, 0) = \Phi(0, -T_0, 0) = 1$ by Lemma~\ref{lem:flow}\ref{item:flow-reverse}, we deduce that $X_0(t, 0) \in [0, 1]$ for every $t \in [-T_0, 0]$, concluding the proof of \ref{item:X0}.

By Lemma~\ref{lem:flow}\ref{item:flow-increasing}, for every $(t, x) \in \mathbb R^2$, the function $s \mapsto \Phi(s, t, x)$ is increasing with $\lim_{s \to -\infty} \Phi(s, t, x) = -\infty$ and $\lim_{s \to +\infty} \Phi(s, t, x) = +\infty$, and thus there exists a unique $R(t, x) \in \mathbb R$ such that $\Phi(R(t, x), t, x) = 0$. The fact that $(t, x) \mapsto R(t, x)$ defines a $\mathcal C^1$ function follows as an immediate consequence of the implicit function theorem and the fact that $\partial_t \Phi(t, t_0, x_0) > 0$ for every $(t, t_0, x_0) \in \mathbb R^3$. In addition, the implicit function theorem, combined with Lem\-ma~\ref{lem:flow}\ref{item:diff-eqn-flow}, \ref{item:deriv-x0}, and \ref{item:deriv-t0} also yields that, for every $(t, x) \in \mathbb R^2$,
\begin{align}
\partial_t R(t, x) & = \frac{\lambda(t, x)}{\lambda(R(t, x), 0)} \exp\left(\int_{t}^{R(t, x)} \partial_x \lambda(s, \Phi(s, t, x)) \diff s\right), \label{eq:diff-t-R}\\
\partial_x R(t, x) & = -\frac{1}{\lambda(R(t, x), 0)} \exp\left(\int_{t}^{R(t, x)} \partial_x \lambda(s, \Phi(s, t, x)) \diff s\right), \label{eq:diff-x-R}
\end{align}
which implies \ref{item:pde-R}, \ref{item:R-incr-t}, and that $x \mapsto R(t, x)$ is decreasing. Since $\Phi(R(t, 0),\allowbreak t, 0) = 0$, we deduce that $R(t, 0) = t$, and thus \ref{item:R-decr-x} also holds true, while \ref{item:R-bounds} follows easily from Lemma~\ref{lem:flow}\ref{item:flow-increasing}.

Combining \ref{item:R-incr-t} and \ref{item:R-decr-x}, we deduce that $\min_{(t, x) \in \mathbb R_+ \times [0, 1]} R(t, x) = R(0, 1) \allowbreak = - T_0$, yielding \ref{item:range-R}. Moreover, \ref{item:range-R0} follows immediately from the fact that $R(0, 0) = 0$ and $R(0, 1) = -T_0$.

To prove \ref{item:R-X0}, take $t \in \mathbb R_+$ such that $R(t, 1) < 0$. Since $t \geq 0$ and $s \mapsto X_0(s, 1)$ is decreasing by \ref{item:X0-decr}, we have that $X_0(t, 1) \leq X_0(0, 1) = 1$. On the other hand, using Lemma~\ref{lem:flow}\ref{item:semigroup}, we have $X_0(t, 1) = \Phi(0, R(t, 1), \Phi(R(t, 1), t, 1)) = \Phi(0, R(t, 1), 0)$ and, by Lemma~\ref{lem:flow}\ref{item:deriv-t0}, we have that $s \mapsto \Phi(0, s, 0)$ is decreasing. Thus, since $R(t, 1) < 0$, we have
\[X_0(t, 1) = \Phi(0, R(t, 1), 0) > \Phi(0, 0, 0) = 0,\]
yielding the conclusion of \ref{item:R-X0}.

To conclude the proof of \ref{item:R}, we use Lemma~\ref{lem:flow}\ref{item:semigroup} to deduce that, for every $t \in \mathbb R$, we have $X_0(t, 1) = \Phi(0, t, 1) = \Phi(0, R(t, 1), \Phi(R(t, 1), t, 1)) = \Phi(0, R(t, 1), 0) = X_0(R(t, 1), 0)$, yielding \ref{item:X0_1_0}.

Concerning the function $\tau$, note that \ref{item:tau-bounds} is an immediate consequence of \ref{item:R-decr-x} and \ref{item:R-bounds}. In addition, combining \ref{item:R-bounds} and \eqref{eq:diff-t-R}, we have
\[
\tau^\prime(t) = 1 - \partial_t R(t, 1) \leq 1 - \frac{\lambda_{\min}}{\lambda_{\max}} e^{-\frac{L}{\lambda_{\min}}} \in (0, 1),
\]
yielding the conclusion of \ref{item:tau-prime}. Item \ref{item:delay-incr} is an immediate consequence of \ref{item:R-incr-t}, while \ref{item:delay-lower} is a particular case of \ref{item:range-R}. Finally, \ref{item:diffeo} follows from \eqref{eq:diff-x-R} and the facts that $R(t, 0) = t$ (from \ref{item:R-decr-x}) and $R(t, 1) = t - \tau(t)$, with $\beta_0 = \frac{1}{\lambda_{\max}} e^{-\frac{L}{\lambda_{\min}}}$ and $\beta_1 = \frac{1}{\lambda_{\min}} e^{\frac{L}{\lambda_{\min}}}$.
\end{proof}

Using Lemma~\ref{lem:lambda}, we deduce the following equivalence between classical solutions of \eqref{eq:transport} and solutions of a difference equation of the form \eqref{eq:diff-eqn-N1-time}.

\begin{proposition}
\label{prop:transport}
Let $A \in \mathcal M_d(\mathbb R)$, $\lambda$ satisfy \ref{hypo:lambda}, and $T_0$, $X_0$, $R$, and $\tau$ be as in the statement of Lemma~\ref{lem:lambda}. If $u \in \mathcal C^0(\mathbb R_+ \times [0, 1], \mathbb R^d) \cap \mathcal C^1(\mathbb R_+^\ast \times (0, 1), \mathbb R^d)$ is a classical solution of \eqref{eq:transport} with initial condition $u_0$, then the function $v \in \mathcal C^0([-T_0, +\infty), \mathbb R^d)$ defined for $t \in [-T_0, +\infty)$ by
\begin{equation}
\label{eq:def-v}
v(t) = \begin{dcases*}
u(t, 0), & if $t \geq 0$, \\
u_0(X_0(t, 0)), & if $t \in [-T_0, 0)$
\end{dcases*}
\end{equation}
is a continuous solution of
\begin{equation}
\label{eq:diff-transport}
\left\{
\begin{aligned}
v(t) & = A v(t - \tau(t)), & \qquad & t \geq 0, \\
v(t) & = v_0(t), & & t \in [-T_0, 0),
\end{aligned}
\right.
\end{equation}
with $v_0(t) = u_0(X_0(t, 0))$ for $t \in [-T_0, 0)$, and in addition $v \in \mathcal C^1((-T_0, +\infty), \mathbb R^d)$. Conversely, if $v \in \mathcal C^0([-T_0, +\infty), \mathbb R^d) \cap \mathcal C^1((-T_0, +\infty), \mathbb R^d)$ is a continuous solution of \eqref{eq:diff-transport}, then the function $u\colon \mathbb R_+ \times [0, 1] \to \mathbb R^d$ defined by $u(t, x) = v(R(t, x))$ is a classical solution of \eqref{eq:transport} with initial condition $u_0(x) = v_0(R(0, x))$.

In addition, for every $t \geq 0$, we have
\begin{equation}
\label{eq:transp-norm-infty}
\norm{u(t, \cdot)}_{L^\infty([0, 1], \mathbb R^d)} = \norm{v_t}_{L^\infty([-\tau(t), 0], \mathbb R^d)}
\end{equation}
and, for every $t \geq 0$ and $p \in [1, +\infty)$, we have
\begin{equation}
\label{eq:transp-norm-p}
\beta_1^{-1/p} \norm{v_t}_{L^p([-\tau(t), 0], \mathbb R^d)} \leq \norm{u(t, \cdot)}_{L^p([0, 1], \mathbb R^d)} \leq \beta_0^{-1/p} \norm{v_t}_{L^p([-\tau(t), 0], \mathbb R^d)},
\end{equation}
where $\beta_0$ and $\beta_1$ are the constants from Lemma~\ref{lem:lambda}\ref{item:diffeo}.
\end{proposition}

\begin{proof}
Let $u \in \mathcal C^0(\mathbb R_+ \times [0, 1], \mathbb R^d) \cap \mathcal C^1(\mathbb R_+^\ast \times (0, 1), \mathbb R^d)$ be a classical solution of \eqref{eq:transport} with initial condition $u_0$ and $v$ be defined from $u$ and $u_0$ as in \eqref{eq:def-v}. Clearly, $v$ is continuous, since $v(0) = u(0, 0) = u_0(0)$ and $v(0^-) = u_0(X_0(0, 0)) = u_0(0)$. By definition, we have $v(t) = v_0(t)$ for every $t \in [-T_0, 0)$. A straightforward computation using Lemma~\ref{lem:flow}\ref{item:diff-eqn-flow} shows that, for every $(t, t_0, x_0) \in \mathbb R^3$ such that $(t, \Phi(t, t_0, x_0)) \in \mathbb R_+^\ast \times (0, 1)$, we have
\begin{equation}
\label{eq:characteristics}
\frac{\df}{\df t} u(t, \Phi(t, t_0, x_0)) = 0.
\end{equation}
In particular, for every $t \geq 0$, if $R(t, 1) \geq 0$, then $(s, \Phi(s, R(t, 1), 0)) \in \mathbb R_+^\ast \times (0, 1)$ for every $s \in (R(t, 1), t)$, and thus $u(t, 1) = u(t, \Phi(t, R(t, 1), 0)) = u(R(t, 1), 0) = u(t - \tau(t), 0)$. Hence, $v(t) = u(t, 0) = A u(t, 1) = A v(t - \tau(t))$, establishing the first case of \eqref{eq:diff-transport} for every $t \geq 0$ such that $R(t, 1) \geq 0$.

Now, let $t \geq 0$ be such that $R(t, 1) < 0$, which implies, by Lemma~\ref{lem:lambda}\ref{item:R-X0}, that $X_0(t, 1) \in (0, 1]$. In addition, for every $s \in (0, t)$, we have $(s, \Phi(s, 0, X_0(t, 1))) \in \mathbb R_+^\ast \times (0, 1)$ since $\Phi(t, 0, X_0(t, 1)) = 1$. Hence, using also Lemma~\ref{lem:lambda}\ref{item:X0_1_0}, we have
\[u(t, 1) = u(t, \Phi(t, 0, X_0(t, 1))) = u(0, X_0(t, 1)) = u_0(X_0(t, 1)) = u_0(X_0(R(t, 1), 0)) = v_0(t - \tau(t)),\]
showing that $u(t, 0) = A u(t, 1) = A v(t - \tau(t))$, and concluding thus the proof that $v$ satisfies \eqref{eq:diff-transport}.

To show that $v$ belongs to $\mathcal C^1((-T_0, +\infty), \mathbb R^d)$, let $\xi\colon\mathbb R_+ \to (0, 1]$ be the function given for $t \in \mathbb R_+^\ast$ by $\xi(t) = \frac{1}{t + 1}$. We claim that, for every $t \in (-T_0, +\infty)$, there exists a unique $S(t) \in \mathbb R_+^\ast$ such that $\Phi(S(t), t, 0) = \xi(S(t))$. Indeed, on the one hand, we have $\Phi(0, -T_0, 0) = 1$ and, since by Lemma~\ref{lem:flow}\ref{item:deriv-t0} the function $s \mapsto \Phi(0, s, 0)$ is decreasing, we deduce that $\Phi(0, t, 0) < 1$ for every $t \in (-T_0, +\infty)$. Thus, $\Phi(0, t, 0) - \xi(0) < 0$ for every $t \in (-T_0, +\infty)$. On the other hand, $\lim_{s \to +\infty} \Phi(s, t, 0) - \xi(s) = +\infty$, and thus the existence of $S(t)$ follows from the intermediate value theorem. In addition, for every $s \in \mathbb R_+^\ast$, we have
\begin{equation}
\label{eq:deriv-St}
\frac{\df}{\df s} \left[\Phi(s, t, 0) - \xi(s)\right] = \lambda(s, \Phi(s, t, 0)) + \frac{1}{(t+1)^2} > 0,
\end{equation}
thus $s \mapsto \Phi(s, t, 0) - \xi(s)$ is increasing in $\mathbb R_+^\ast$, yielding the uniqueness of $S(t)$. We also obtain from \eqref{eq:deriv-St} and the implicit function theorem that $S \in \mathcal C^1((-T_0, +\infty), \mathbb R_+^\ast)$.

For $t \geq 0$, using \eqref{eq:characteristics}, we have that $v(t) = u(t, 0) = u(S(t), \Phi(S(t), t, 0)) = u(S(t), \xi(S(t)))$, since $(s, \Phi(s, t, 0)) \in \mathbb R_+^\ast \times (0, 1)$ for every $s \in (t, S(t))$. For $t \in (-T_0, 0)$, still using \eqref{eq:characteristics}, we have that
\begin{align*}
v(t) & = u(0, X_0(t, 0)) = u(S(t), \Phi(S(t), 0, X_0(t, 0))) \\
& = u(S(t), \Phi(S(t), t, 0)) = u(S(t), \xi(S(t)))
\end{align*}
since $(s, \Phi(s, 0, X_0(t, 0))) \in \mathbb R_+^\ast \times (0, 1)$ for every $s \in (0, S(t))$. Thus, we have $v(t) = u(S(t), \xi(S(t)))$ for every $t \in (-T_0, +\infty)$ and, since $S \in \mathcal C^1((-T_0, +\infty), \mathbb R_+^\ast)$, $\xi \in \mathcal C^1(\mathbb R_+^\ast, (0, 1))$, and $u \in \mathcal C^1(\mathbb R_+^\ast \times (0, 1), \mathbb R^d)$, we deduce that $v \in \mathcal C^1((-T_0, +\infty), \mathbb R^d)$.

Conversely, consider a continuous solution $v \in \mathcal C^0([-T_0, +\infty), \mathbb R^d) \cap \mathcal C^1((-T_0, +\infty), \mathbb R^d)$ of \eqref{eq:diff-transport} with initial condition $v_0$ and define $u\colon \mathbb R_+ \times [0, 1] \to \mathbb R^d$ by $u(t, x) = v(R(t, x))$ and $u_0\colon [0, 1] \to \mathbb R^d$ by $u_0(x) = v_0(R(0, x))$. Since $R \in \mathcal C^1(\mathbb R^2, \mathbb R)$, we have $u \in \mathcal C^0(\mathbb R_+ \times [0, 1], \mathbb R^d) \cap \mathcal C^1(\mathbb R_+^\ast \times (0, 1), \mathbb R^d)$ and, using Lemma~\ref{lem:lambda}\ref{item:pde-R}, it is immediate to verify that
\[\partial_t u(t, x) + \lambda(t, x) \partial_x u(t, x) = 0 \qquad \text{ for every } (t, x) \in \mathbb R_+^\ast \times (0, 1).\]
For $x \in [0, 1]$, we have $u(0, x) = v(R(0, x)) = u_0(x)$ and, for $t \geq 0$, we have $u(t, 0) = v(R(t, 0)) = v(t) = A v(R(t, 1)) = u(t, 1)$. Hence, $u$ is a classical solution of \eqref{eq:transport} with initial condition $u_0$.

Let us finally show \eqref{eq:transp-norm-infty} and \eqref{eq:transp-norm-p}. Let $p \in [1, +\infty)$, $t \geq 0$, and $x \in (0, 1)$. We have $u(t, x) = v(R(t, x))$, and thus, using Lemma~\ref{lem:lambda}\ref{item:diffeo}, we deduce that
\[
\norm{u(t, \cdot)}_{L^\infty([0, 1], \mathbb R^d)} = \sup_{x \in (0, 1)} \abs{v(R(t, x))} = \sup_{s \in (t - \tau(t), t)} \abs{v(s)} = \norm{v_t}_{L^\infty([t - \tau(t), t], \mathbb R^d)},
\]
yielding \eqref{eq:transp-norm-infty}, and
\[
\norm{u(t, \cdot)}_{L^p([0, 1], \mathbb R^d)}^p = \int_0^1 \abs{v(R(t, x))}^p \diff x = \int_{t - \tau(t)}^{t} \frac{\abs{v(s)}^p}{\abs{\partial_x R(t, \raisebox{\depth}{\scalebox{1}[-1]{R}}(t, s))}} \diff s,
\]
where $\raisebox{\depth}{\scalebox{1}[-1]{R}}(t, \cdot)$ is the inverse of $R(t, \cdot)$. Hence \eqref{eq:transp-norm-p} follows from Lemma~\ref{lem:lambda}\ref{item:diffeo}.
\end{proof}

Using Proposition~\ref{prop:transport} and the results from Section~\ref{sec:exponential}, we obtain at once the following result.

\begin{corollary}
Let $A \in \mathcal M_d(\mathbb R)$ satisfy \ref{hypo:spr-A-less-1} and $\lambda$ satisfy \ref{hypo:lambda}. Then there exist positive constants $C$ and $\gamma$ such that, for every classical solution $u \in \mathcal C^0(\mathbb R_+ \times [0, 1], \mathbb R^d) \cap \mathcal C^1(\mathbb R_+^\ast \times (0, 1), \mathbb R^d)$ of \eqref{eq:transport} with initial condition $u_0$, we have
\begin{equation}
\label{eq:transport-Linfty}
\norm{u(t, \cdot)}_{L^\infty([0, 1], \mathbb R^d)} \leq C e^{-\gamma t} \norm{u_0}_{L^\infty([0, 1], \mathbb R^d)} \qquad \text{ for every } t \in \mathbb R_+.
\end{equation}

Moreover, let $\alpha$ be the constant from Lemma~\ref{lem:lambda}\ref{item:tau-prime}, $p \in [1, +\infty)$, and assume that $\spr(A)^p < 1-\alpha$. Then there exist positive constants $C$ and $\gamma$ such that, for every classical solution $u \in \mathcal C^0(\mathbb R_+ \times [0, 1], \mathbb R^d) \cap \mathcal C^1(\mathbb R_+^\ast \times (0, 1), \mathbb R^d)$ of \eqref{eq:transport} with initial condition $u_0$, we have
\begin{equation}
\label{eq:transport-Lp}
\norm{u(t, \cdot)}_{L^p([0, 1], \mathbb R^d)} \leq C e^{-\gamma t} \norm{u_0}_{L^p([0, 1], \mathbb R^d)} \qquad \text{ for every } t \in \mathbb R_+.
\end{equation}
\end{corollary}

\begin{proof}
Let $\tau$ be defined from $\lambda$ as in Lemma~\ref{lem:lambda}\ref{item:tau} and note that, by item \ref{item:tau-bounds} of that lemma, $\tau$ satisfies \ref{hypo:tau-infimum} and \ref{hypo:delay-bounded}. In addition, Lemma~\ref{lem:lambda}\ref{item:delay-incr} implies, from Definition~\ref{def:h}, that the largest delay function $h$ associated with $\tau$ is $h = \tau$. Moreover, combining Lemma~\ref{lem:lambda}\ref{item:tau-prime} and Remark~\ref{remk:radon-nikodym-derivative-tau}, we deduce that $\tau$ satisfies \ref{hypo:radon-nik-bdd}. Then \eqref{eq:transport-Linfty} is a consequence of Corollary~\ref{coro:exp-Linfty}.

We also note that, combining once again Lemma~\ref{lem:lambda}\ref{item:tau-prime} and Remark~\ref{remk:radon-nikodym-derivative-tau}, the fact that $\spr(A)^p < 1-\alpha$ implies that \ref{hypo:radon-nik-and-rho} is satisfied, and thus \eqref{eq:transport-Lp} is a consequence of Theorem~\ref{thm:exp-Lp}.
\end{proof}

\chapter{Stability and stabilization of time-delay systems through the MID property}
\label{chap:MID}

\epigraph{\foreignlanguage{brazil}{É preciso desconfiar sete vezes do cálculo e cem vezes do matemático}}{\foreignlanguage{brazil}{\emph{O Homem que Calculava}, Malba Tahan}}

This chapter deals with the stability analysis and the stabilization of time-delay systems through spectral methods, its aim being to establish the multiplicity-induced dominance (MID) property discussed in Section~\ref{sec:intro-MID} for some classes of time-delay systems. Its contents are extracted from the works \cite{Boussaada2022Some, Boussaada2022Generic, Mazanti2021Multiplicity, Boussaada2022Pade, Mazanti2021Effects, Benarab2020MID, Mazanti2020Qualitative, Boussaada2025OverOrder, Boussaada2025Decay, Boussaada2022MID, Mazanti2020Spectral, Fueyo2021Insights, Fueyo2023Pole, Niculescu2021Stability}.

In most of the chapter, we will work with time-delay systems described by the scalar delay-differential equation \eqref{eq:DDE} with a single, constant delay, or with particular cases of this equation, the exceptions being Section~\ref{sec:MID-DDAE}, in which we will consider a delay differential-algebraic equation, and Section~\ref{sec:MID-two-delays}, in which we will consider a first order delay-dif\-fe\-ren\-tial equation with two delays.

After presenting some required preliminary results in Section~\ref{sec:MID-prelim}, we start exposing the main results of this chapter in Section~\ref{sec:MID-maximal}, in which we prove the MID property for real spectral values of maximal multiplicity of the DDE \eqref{eq:DDE}. Some partial results for real spectral values with multiplicity smaller than the maximal one are provided in Section~\ref{sec:MID-intermediate}. We discuss dominance of pairs of complex conjugate roots for low-order DDEs in Section~\ref{sec:MID-complex}, and we address the case of a first-order DDE with two delays in Section~\ref{sec:MID-two-delays}. Finally, we present through several examples in Section~\ref{sec:MID-applis} how the MID property can be used in order to stabilize control systems with delays. Appendix~\ref{sec:MID-app} provides technical results used throughout the chapter.

\section{Preliminary results}
\label{sec:MID-prelim}

We now present some preliminary results on Padé approximants and Kummer and Whittaker functions useful for the sequel of this chapter.

\clearpage

\subsection{Padé approximants}
\label{sec:Pade}

Studied at least since the early works on rational approximations of functions by G.~Frobenius in \cite{Frobenius1881Relationen} and named after H.~Padé thanks to his major works on the topic \cite{Pade1892Representation}, Padé approximants can be seen as a generalization of Taylor approximations of a function in which, instead of approximating a given function $f$ around some point $z_0$ in its domain by polynomials, one wishes to approximate $f$ around $z_0$ by a rational function, i.e., by a function $R$ expressed as $R = \frac{Q}{P}$, where $Q$ and $P$ are polynomials. These approximants appear when studying many classical topics of Mathematics and Engineering and, for this reason, they have found several applications along the years, including in the analytic theory of continued fractions, the analysis of orthogonal families of polynomials, the approximation of the frequency response of control systems with constant time delays, as well as some topics in numerical analysis, signal processing, and theoretical physics. We refer the interested reader to \cite{Glader1991Approximation, Allemandou1966Low, Baker1996Pade, Baker1970Pade, Brezinski1991History, Brezinski1980Pade, Brezinski2002Computational, Lam1990Convergence, Lam1993Model, Perron1977Lehre, Saff1975Zeros, Saff1977Zeros, Saff1978Zeros, Stewart1960Generalized} for deeper presentations of Padé approximants and their applications.

We now present the main elements on Padé approximants that we will need in the sequel of this chapter, based on the presentations of \cite{Brezinski1980Pade, Saff1975Zeros, Baker1996Pade, Perron1977Lehre}. We start with the precise definition of Padé approximant used here. For simplicity, we only consider here Padé approximants around the origin.

\begin{definition}
\label{def:Pade}
Let $U \subset \mathbb C$ be a neighborhood of $0$, $f\colon U \to \mathbb C$, and $(m, n) \in \mathbb N^2$. We say that a pair of polynomials $(Q, P)$ with $\deg Q \leq m$ and $\deg P \leq n$ is a \emph{Padé approximant of $f$ of order $(m, n)$} if $P(0) \neq 0$ and
\begin{equation}
\label{eq:def-Pade}
f(z) - \frac{Q(z)}{P(z)} = O(z^{n + m + 1}) \qquad \text{ as } z \to 0.
\end{equation}
\end{definition}

A series of remarks are now in order.

\begin{remark}
Padé approximants of a given order are not necessarily unique since, if $(Q, P)$ is a Padé approximant of a given function, then $(\lambda Q, \lambda P)$ is a Padé approximant of the same function and with the same order for every $\lambda \in \mathbb C \setminus \{0\}$. A standard normalization for Padé approximants is to require that $P(0) = 1$.
\end{remark}

\begin{remark}
If $(Q, 1)$ is a Padé approximant of a given function $f$ of order $(m, 0)$, then $Q$ coincides with the Taylor--Maclaurin approximation of $f$ of order $m$.
\end{remark}

\begin{remark}
\label{remk:Pade-alternative}
Under the assumption that $P(0) \neq 0$, \eqref{eq:def-Pade} is equivalent to
\begin{equation}
\label{eq:error-Pade}
P(z) f(z) - Q(z) = O(z^{n + m + 1}) \qquad \text{ as } z \to 0.
\end{equation}

Equation~\eqref{eq:error-Pade}, without the requirement that $P(0) \neq 0$, was used as a definition of Padé approximants in early works, including \cite{Frobenius1881Relationen, Pade1892Representation}, and it has the advantage that, if $f$ is analytic on $U$, then, for every $(m, n) \in \mathbb N^2$, there exists a pair of polynomials $(Q, P)$ with $\deg Q \leq m$, $\deg P \leq n$, and $P$ not identically zero such that \eqref{eq:error-Pade} is satisfied. The issue with this alternative definition is that \eqref{eq:def-Pade} may fail to be satisfied: for instance, if $f(z) = 1 + z^2$, then one immediately verifies that the polynomials $Q$ and $P$ defined by $Q(z) = z$ and $P(z) = z$ satisfy \eqref{eq:error-Pade} with $n = m = 1$, but \eqref{eq:def-Pade} is not satisfied. On the other hand, using Definition~\ref{def:Pade}, Padé approximants of analytic functions may fail to exist at some orders: one can easily verify that the previous function $f$ has no Padé approximant of order $(1, 1)$.
\end{remark}

Several applications require Padé approximants of the exponential function, in particular in numerical analysis, signal processing, and control theory \cite{Saff1975Zeros, Saff1977Zeros, Saff1978Zeros, Allemandou1966Low, Lam1990Convergence, Lam1993Model, Stewart1960Generalized, Glader1991Approximation} and, for this reason, such approximants were extensively studied in the literature, including in early works such as \cite{Pade1892Representation}. We collect, in the next proposition, some of the results presented in \cite[\S~42]{Perron1977Lehre} (see also \cite[Section~1.2]{Baker1996Pade} and \cite[\S~66]{Pade1892Representation}).

\begin{proposition}
\label{prop:Pade-exp}
Let $m$ and $n$ be nonnegative integers and define
\begin{align*}
Q(z) & = \sum_{k=0}^m \frac{m! (n+m-k)!}{(m-k)! (n+m)!} \frac{z^k}{k!},\\
P(z) & = \sum_{k=0}^n \frac{n! (n+m-k)!}{(n-k)! (n+m)!} \frac{(-z)^k}{k!}.
\end{align*}
Then $(Q, P)$ is a Padé approximant of the exponential function $z \mapsto e^z$ of order $(m, n)$, and it is the unique such approximant with $P(0) = 1$. In addition, we have
\[
P(z) e^z - Q(z) = (-1)^n \frac{z^{n + m + 1}}{(n + m)!} \int_0^1 t^n (1-t)^m e^{z t} \diff t.
\]
\end{proposition}

In addition to providing explicit expressions for the Padé approximant of a given order of the exponential function, Proposition~\ref{prop:Pade-exp} also provides an explicit expression for the approximation error $P(z) e^z - Q(z)$, which not only is of order $z^{n + m + 1}$, as required by \eqref{eq:def-Pade}, but can also be written as $z^{n + m + 1}$ multiplied by an holomorphic function of $z$ expressed in terms of an integral. This holomorphic function turns out to be a Kummer confluent hypergeometric function, which we present in details next.

\subsection{Kummer and Whittaker functions}
\label{sec:Kummer-Whittaker}

We present in this section the definitions and main properties of Kummer and Whittaker functions used in the sequel, reporting here in particular the main results of \cite{Boussaada2022Some} on the location of their zeros.

Kummer and Whittaker functions are two kinds of special functions, named respectively after E.~E.~Kummer and E.~T.~Whittaker, who introduced these functions in \cite{Kummer1837Integralibus, Whittaker1903Expression}. Together with other special functions, such as Tricomi functions and Coulomb wave functions, they are part of the family of \emph{confluent hypergeometric functions}, a class of functions widely studied in the literature in particular due to its connections with problems in mathematical physics \cite{OldeDaalhuis2010Confluent, Slater1972Confluent, Buchholz1969Confluent, Tricomi1954Funzioni, Slater1960Confluent, Erdelyi1981Higher}.

\subsubsection{Elementary properties}

Let us start by providing the definitions and some classical results on Kummer and Whittaker functions, extracted from \cite{OldeDaalhuis2010Confluent}. Recall that, for $\alpha \in \mathbb C$ and $k \in \mathbb N$, $(\alpha)_k$ is the Pochhammer symbol for the ascending factorial, defined inductively as $(\alpha)_0 = 1$ and $(\alpha)_{k+1} = (\alpha+k) (\alpha)_k$ for $k \in \mathbb N$.

\begin{definition}[{\cite[(13.2.2) and (13.14.2)]{OldeDaalhuis2010Confluent}}]
\label{defi:Kummer-Whittaker}
Let $a, b, \kappa, \mu \in \mathbb C$ and assume that $-b \notin \mathbb N$ and $-2\mu \notin \mathbb N^\ast$.
\begin{enumerate}
\item The \emph{Kummer function} $M(a, b, \cdot)$ is the complex function defined for $z \in \mathbb C$ by the series
\begin{equation}
\label{DefiConfluent}
M(a, b, z) = \sum_{k=0}^{\infty} \frac{(a)_k}{(b)_k} \frac{z^k}{k!}.
\end{equation}

\item The \emph{Whittaker function} $M_{\kappa, \mu}$ is the complex function defined from the Kummer function $M\left(\frac{1}{2} + \mu - \kappa, 1 + 2 \mu, \cdot\right)$ by
\begin{equation}
\label{eq:def-Whittaker}
M_{\kappa, \mu}(z) = e^{-\frac{1}{2} z} z^{\mu + \frac{1}{2}} M\left(\frac{1}{2} + \mu - \kappa, 1 + 2 \mu, z\right).
\end{equation}
\end{enumerate}
\end{definition}

Note that the series in \eqref{DefiConfluent} converges for every $z \in \mathbb C$, and hence the Kummer function $M(a, b, \cdot)$ is defined on the whole complex plane $\mathbb C$ and satisfies
\begin{equation}
\label{eq:Kummer-in-zero}
M(a, b, z) = 1 + O(z) \qquad \text{ as } z \to 0.
\end{equation}
As for the Whittaker function, one observes that, if $\mu + \frac{1}{2}$ is not an integer, $M_{\kappa, \mu}$ is a multi-valued complex function with branch point at $z = 0$, while it is a single-valued function if $\mu + \frac{1}{2}$ is an integer. Despite being possibly multi-valued, it follows immediately from \eqref{eq:def-Whittaker} that the zeros of $M_{\kappa, \mu}$ do not depend on which branch of $M_{\kappa, \mu}$ is being taken into account: if $z_0$ is a zero of some branch of $M_{\kappa, \mu}$, then either $\mu + \frac{1}{2} > 0$ and $z_0 = 0$, or $M\left(\frac{1}{2} + \mu - \kappa, 1 + 2 \mu, z_0\right) = 0$, and in both cases $z_0$ will be a zero of all branches of $M_{\kappa, \mu}$. It also follows that the nontrival zeros of $M_{\kappa, \mu}$ coincide with those of $M\left(\frac{1}{2} + \mu - \kappa, 1 + 2 \mu, \cdot\right)$.

Both Kummer and Whittaker functions satisfy some classical second-order ordinary differential equations (ODEs), known as Kummer and Whittaker differential equations, which are presented in the next result.

\begin{proposition}[{\cite[(13.2.1) and (13.14.1)]{OldeDaalhuis2010Confluent}}]
Let $a, b, \kappa, \mu \in \mathbb C$ and assume that $-b \notin \mathbb N$ and $-2\mu \notin \mathbb N^\ast$.
\begin{enumerate}
\item The Kummer function $M(a, b, \cdot)$ satisfies the \emph{Kummer differential equation}
\begin{equation}
\label{KummerODE}
z \frac{\partial^2 M}{\partial z^2}(a, b, z) + (b - z) \frac{\partial M}{\partial z}(a, b, z) - a M(a, b, z) = 0.
\end{equation}

\item The Whittaker function $M_{\kappa, \mu}$ satisfies the \emph{Whittaker differential equation}
\begin{equation}
\label{WhittakerODE}
M_{\kappa, \mu}^{\prime\prime}(z) + \left(-\frac{1}{4} + \frac{\kappa}{z} + \frac{\frac{1}{4} - \mu^2}{z^2}\right)M_{\kappa, \mu}(z) = 0.
\end{equation}
\end{enumerate}
\end{proposition}

Note that both \eqref{KummerODE} and \eqref{WhittakerODE} may admit other solutions than $M(a, b, \cdot)$ and $M_{\kappa, \mu}$, which are also sometimes called Kummer or Whittaker functions (or also Tricomi functions in the case of \eqref{KummerODE}), but we use here the names \emph{Kummer} and \emph{Whittaker functions} to refer to the functions from Definition~\ref{defi:Kummer-Whittaker} only.

An important result for the sequel is the following one, providing an integral representation for Kummer functions for some values of their parameters.

\begin{proposition}[{\cite[(13.2.4) and (13.4.1)]{OldeDaalhuis2010Confluent}}]
\label{prop:Kummer-integral}
Let $a, b \in \mathbb C$ be such that $\Real b > \Real a > 0$. Then, for every $z \in \mathbb C$, we have
\begin{equation}
\label{eq:Kummer-integral}
M(a, b, z) = \frac{\Gamma(b)}{\Gamma(a) \Gamma(b - a)} \int_0^1 t^{a-1} (1-t)^{b-a-1} e^{zt} \diff t,
\end{equation}
where $\Gamma$ denotes the Gamma function.
\end{proposition}

One of the basic relations satisfied by Kummer functions, known as a \emph{Kummer transformation}, is that, for every $a, b \in \mathbb C$ such that $-b \notin \mathbb N$, we have
\begin{equation}
\label{eq:Kummer-transf-1}
M(a, b, z) = e^z M(b - a, b, -z), \qquad \text{ for every } z \in \mathbb C
\end{equation}
(see, e.g., \cite[(13.2.39)]{OldeDaalhuis2010Confluent}). As a consequence of \eqref{eq:Kummer-transf-1}, we obtain the following immediate symmetry property of roots of Whittaker functions.

\begin{proposition}
\label{PropWhittakerSymmetricRoots}
Let $\kappa, \mu \in \mathbb C$ and assume that $-2\mu \notin \mathbb N^\ast$. If $z \in \mathbb C \setminus \{0\}$ is a root of the Whittaker function $M_{\kappa, \mu}$, then $-z$ is a root of $M_{-\kappa, \mu}$.
\end{proposition}

\begin{proof}
Let $z \in \mathbb C \setminus \{0\}$ be a root of $M_{\kappa, \mu}$. Since $z \neq 0$, we deduce from \eqref{eq:def-Whittaker} that $M(\frac{1}{2} + \mu - \kappa, \allowbreak 1 + 2 \mu, z) = 0$. Using \eqref{eq:Kummer-transf-1}, we obtain that $M(\frac{1}{2} + \mu + \kappa, 1 + 2 \mu, -z) = 0$, and thus, by \eqref{eq:def-Whittaker}, $-z$ is a root of $M_{-\kappa, \mu}$.
\end{proof}

\subsubsection{Counterexamples to some previous results on the location of zeros of Whittaker functions}
\label{sec:counterexamples}

Some works in the literature have addressed the problem of getting information on the zeros of Kummer and Whittaker functions, and several results in that sense are provided in the classical references \cite{OldeDaalhuis2010Confluent, Slater1972Confluent, Buchholz1969Confluent, Tricomi1954Funzioni, Slater1960Confluent, Erdelyi1981Higher}, including results on the number of their real zeros and on the asymptotic behavior of their complex zeros at complex infinity.

We next state other results from \cite{Tsvetkov1941Roots, Saff1978Zeros} on the location of complex zeros of Whittaker functions in the case where the parameters $\kappa$ and $\mu$ are real numbers. We start with the following result from \cite{Tsvetkov1941Roots}.

\begin{proposition}[{\cite[Theorem~7]{Tsvetkov1941Roots}}]
\label{PropTsvetkoff1}
Let $\kappa, \mu \in \mathbb R$ be such that $\mu > -\frac{1}{2}$ and $z_0 \in \mathbb C \setminus \{0\}$ be such that $M_{\kappa, \mu}(z_0) = 0$.
\begin{enumerate}
\item\label{PropTsvetkoff-k-g-0} If $\kappa > 0$, then $\Real z_0 > 2 \kappa$.
\item\label{PropTsvetkoff-k-l-0} If $\kappa < 0$, then $\Real z_0 < 2 \kappa$.
\item If $\kappa = 0$, then $\Real z_0 = 0$.
\end{enumerate}
\end{proposition}

We also state the following two results from \cite{Saff1978Zeros}.

\begin{proposition}[{\cite[Proposition~3.2]{Saff1978Zeros}}]
\label{PropSaffVarga}
Let $\kappa, \mu \in \mathbb R$ be such that $\mu > 0$ and $z_0 \in \mathbb C \setminus \{0\}$ be such that $M_{\kappa, \mu}(z_0) = 0$.
\begin{enumerate}
\item\label{PropSaffVarga-k-g-0} If $\kappa > 0$, then $\Real z_0 > 2 \kappa$ and $(\Imag z_0)^2 > (\Real z_0)^2 \frac{4 \mu^2 - 4 \kappa^2 - 1}{4 \kappa^2}$.
\item\label{PropSaffVarga-k-l-0} If $\kappa < 0$, then $\Real z_0 < 2 \kappa$ and $(\Imag z_0)^2 > (\Real z_0)^2 \frac{4 \mu^2 - 4 \kappa^2 - 1}{4 \kappa^2}$.
\item If $\kappa = 0$, then $\Real z_0 = 0$ and $(\Imag z_0)^2 > 4 \mu^2 - 1$.
\end{enumerate}
\end{proposition}

\begin{proposition}[{\cite[Corollary~3.3]{Saff1978Zeros}}]
\label{CoroSaffVarga}
Let $\kappa, \mu \in \mathbb R$ be such that $\mu > \frac{1}{2}$ and $z_0 \in \mathbb C \setminus \{0\}$ be such that $M_{\kappa, \mu}(z_0) = 0$. Then $\abs{z_0} > \sqrt{4 \mu^2 - 1}$.
\end{proposition}

Unfortunately, no proof of Proposition~\ref{PropTsvetkoff1} is provided in the short note \cite{Tsvetkov1941Roots}. A more precise result, containing Proposition~\ref{PropTsvetkoff1} as a particular case, was also stated in \cite[Theorems~2 and 3]{Tsvetkov1941Complex}, once again without proof. Both \cite{Tsvetkov1941Roots, Tsvetkov1941Complex} only mention that the idea of the proof follows by using the strategy of \cite{Hille1922Oscillation}. Some of the results from \cite{Tsvetkov1941Complex, Tsvetkov1941Roots} were later described by F.~G.~Tricomi in \cite{Tricomi1950Abzahlung} as being\footnote{In the original German, \foreignlanguage{ngerman}{``die wichtigen aber (meines Wissens) bisher nicht geprüften Resultate von Herrn TSVETTKOFF [...]''.}} ``important'' but ``unchecked'', and it turns out that other theorems from \cite{Tsvetkov1941Complex, Tsvetkov1941Roots} were shown to have issues: \cite[pages 248--249]{Saff1978Zeros} provides a counterexample to \cite[Theorem~11]{Tsvetkov1941Roots}. On the other hand, the proofs of Propositions~\ref{PropSaffVarga} and \ref{CoroSaffVarga} provided in \cite[Proposition~3.2 and Corollary~3.3]{Saff1978Zeros} rely on Proposition~\ref{PropTsvetkoff1} above, combined with the method of \cite{Hille1922Oscillation} invoked in \cite{Tsvetkov1941Complex, Tsvetkov1941Roots}.

It turns out that Propositions~\ref{PropTsvetkoff1}--\ref{CoroSaffVarga} are not correct, and we now illustrate that by counterexamples. We start by the following example, which contradicts the assertions of Propositions~\ref{PropTsvetkoff1} and \ref{PropSaffVarga} for some particular values of $\kappa$ and $\mu$.

\begin{example}
\label{ExplTsvetkovFalse}
Let $\mu \in (-\frac{1}{2}, +\infty)$ and take $\kappa = \mu + \frac{3}{2}$. Let $z_0 \in \mathbb C \setminus \{0\}$ be such that $M_{\kappa, \mu}(z_0) = 0$. Then, by \eqref{eq:def-Whittaker}, we have $M(-1, 1 + 2 \mu, z_0) = 0$. From \eqref{DefiConfluent}, we have $M(-1, 1 + 2 \mu, z) = 1 - \frac{z}{1 + 2 \mu}$, and thus $z_0 = 1 + 2 \mu$. In particular, $\Real z_0 = 1 + 2 \mu < 2 \kappa = 3 + 2 \mu$, and thus Proposition~\ref{PropTsvetkoff1}\ref{PropTsvetkoff-k-g-0} and Proposition~\ref{PropSaffVarga}\ref{PropSaffVarga-k-g-0} are not verified.
\end{example}

\begin{remark}
The above example only contradicts the conclusion that $\Real z_0 > 2 \kappa$ from Proposition~\ref{PropSaffVarga}\ref{PropSaffVarga-k-g-0} and, as we will prove later in Proposition~\ref{PropTsvetkoffCorrected}, the second inequality in the conclusion of Proposition~\ref{PropSaffVarga}\ref{PropSaffVarga-k-g-0} is actually true. In addition, as a consequence of Proposition~\ref{PropWhittakerSymmetricRoots}, Example~\ref{ExplTsvetkovFalse} also provides a counterexample to Proposition~\ref{PropTsvetkoff1}\ref{PropTsvetkoff-k-l-0} and Proposition~\ref{PropSaffVarga}\ref{PropSaffVarga-k-l-0}.
\end{remark}

\begin{figure}[ht]
\centering
\resizebox{0.5\textwidth}{!}{\input{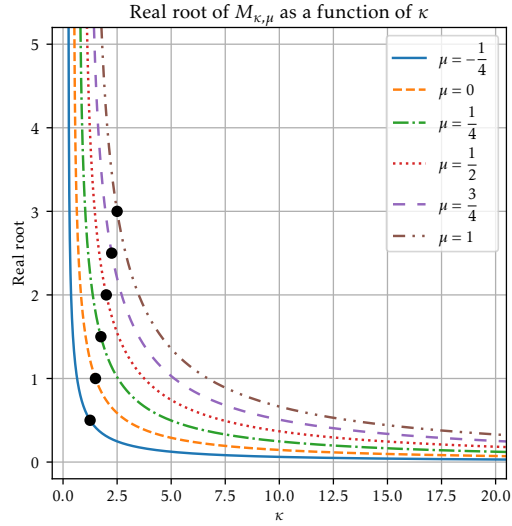}}
\caption[Real root $z_\mu(\kappa)$ of $M_{\kappa, \mu}$ satisfying $z_\mu(\mu + \frac{3}{2}) = 1 + 2 \mu$]{Real root $z_\mu(\kappa)$ of $M_{\kappa, \mu}$ satisfying $z_\mu(\mu + \frac{3}{2}) = 1 + 2 \mu$ for six different values of $\mu$}
\label{FigRootWhittaker}
\end{figure}

Let us consider the root $z_0 = 1 + 2 \mu$ of $M_{\kappa, \mu}$ for $\kappa = \mu + \frac{3}{2}$ from Example~\ref{ExplTsvetkovFalse}. Since this root is real and simple, for every $\mu > -\frac{1}{2}$, there exists an interval $I_\mu$ containing $1 + 2 \mu$ and a curve $\kappa \mapsto z_\mu(\kappa) \in \mathbb R$ defined on $I_\mu$ such that, for every $\kappa \in I_\mu$, $z_\mu(\kappa)$ is a real root of $M_{\kappa, \mu}$ with $z_\mu(\mu + \frac{3}{2}) = 1 + 2\mu$. These curves were computed numerically\footnote{Numerical computations were done in Python using the function \texttt{root\_scalar} from the module \texttt{sci\-py.\allowbreak{}op\-ti\-mize} from \cite{Virtanen2020SciPy}.} for $\mu \in \left\{-\frac{1}{4}, 0, \frac{1}{4}, \frac{1}{2}, \frac{3}{4}, 1\right\}$ and are represented in the $(\kappa, z_0)$ plane in Figure~\ref{FigRootWhittaker}. The black dots $(\kappa, z_0)$ correspond to $\kappa = \mu + \frac{3}{2}$ and the root $z_0 = 1 + 2 \mu$ from Example~\ref{ExplTsvetkovFalse}.

An inspection of Figure~\ref{FigRootWhittaker} leads to the conjecture that the maximal interval $I_\mu$ on which $z_\mu$ is defined is $I_\mu = (\mu + \frac{1}{2}, +\infty)$ and that $z_\mu(\kappa) \to 0$ as $\kappa \to +\infty$ and $z_\mu(\kappa) \to +\infty$ as $\kappa \to \mu + \frac{1}{2}$. In particular, if this conjecture is true, then one cannot expect to correct Proposition~\ref{PropTsvetkoff1}\ref{PropTsvetkoff-k-g-0} by replacing the term $2 \kappa$ by any function of $\kappa$ which remains lower bounded as $\kappa \to +\infty$.

Let us now present our counterexample to Proposition~\ref{CoroSaffVarga}.

\begin{example}
\label{ExplFalseSaffVarga}
Let $\mu \in (-\frac{1}{2}, +\infty)$ and take $\kappa = \mu + \frac{5}{2}$. Let $z_0 \in \mathbb C \setminus \{0\}$ be such that $M_{\kappa, \mu}(z_0) = 0$. Then, by \eqref{eq:def-Whittaker}, we have $M(-2, 1 + 2 \mu, z_0) = 0$. From \eqref{DefiConfluent}, we have $M(-2, 1 + 2 \mu, z) = 1 - \frac{2 z}{1 + 2 \mu} + \frac{z^2}{(1 + 2\mu)(2 + 2\mu)}$, whose roots are the positive real numbers $z_+ = 2 + 2\mu + \sqrt{2 + 2\mu}$ and $z_- = 2 + 2\mu - \sqrt{2 + 2\mu}$, hence $z_0 \in \{z_+, z_-\}$. Standard computations show that $z_- = \abs{z_-} \leq \sqrt{4 \mu^2 - 1}$ if and only if $\mu \geq \frac{5 + 3 \sqrt{33}}{16}$. In particular, the bound $\abs{z_0} > \sqrt{4 \mu^2 - 1}$ is not satisfied for a nonzero real root of $M_{\mu + \frac{5}{2}, \mu}$ as soon as $\mu \geq \frac{5 + 3 \sqrt{33}}{16}$.
\end{example}

\subsubsection{Location of zeros of Whittaker and Kummer functions}
\label{sec:correction}

Our statement below provides a correction of Propositions~\ref{PropTsvetkoff1}--\ref{CoroSaffVarga}.

\begin{proposition}
\label{PropTsvetkoffCorrected}
Let $\kappa \in \mathbb R$, $\mu \in \mathbb R_+^\ast$, and $z_0 \in \mathbb C \setminus \{0\}$ be such that $M_{\kappa, \mu}(z_0) = 0$.
\begin{enumerate}
\item\label{PropTsvetkoffCorrected-Sign} If $\mu \geq \frac{1}{2}$, then $\sign(\kappa) = \sign(\Real z_0)$.
\item\label{PropTsvetkoffCorrected-k-neq-0} If $\kappa \neq 0$, then $\Real z_0 \neq 0$ and $4 \kappa^2 (\Imag z_0)^2 > (\Real z_0)^2 \left(4\mu^2 - 4\kappa^2 - 1\right)$.
\item\label{PropTsvetkoffCorrected-z0-not-real} If $z_0 \notin \mathbb R$, then $\abs{z_0}^2 > 4 \mu^2 - 1$.
\end{enumerate}
\end{proposition}

As discussed in Section~\ref{sec:counterexamples}, even though G.~E.~Tsvetkov does not provide the proofs of the results in \cite{Tsvetkov1941Complex, Tsvetkov1941Roots}, the reader of those works is referred to techniques developed by E.~Hille in \cite{Hille1922Oscillation}, which have also been explored in \cite{Saff1978Zeros}. Following Tsvetkov's suggestion, we shall provide here a proof of Proposition~\ref{PropTsvetkoffCorrected} based on the techniques from \cite{Hille1922Oscillation}. Our proof will follow closely that provided in \cite[Proposition~3.2]{Saff1978Zeros} for the last inequalities in Proposition~\ref{PropSaffVarga}\ref{PropSaffVarga-k-g-0} and \ref{PropSaffVarga-k-l-0} above, but exploring some additional properties in order to avoid invoking Proposition~\ref{PropTsvetkoff1} as done in \cite[Proposition~3.2]{Saff1978Zeros}. Before turning to the proof itself, let us comment on Hille's approach.

In \cite{Hille1922Oscillation}, Hille studies the distribution of zeros of functions of a complex variable satisfying linear second-order homogeneous differential equations with variable coefficients, as is the case for the degenerate Whittaker function $M_{\kappa,\mu}$, which satisfies \eqref{WhittakerODE}. Using a suitable integration of the equation, called \emph{Green--Hille transform}, together with some further conditions on the behavior of the function, Hille shows how to discard regions in the complex plane from including complex roots.

More precisely, consider the general homogeneous second-order differential equation
\begin{equation}
\label{eq:Hille}
\frac{\df}{\df z}\left[K(z) \frac{\df\varphi}{\df z}(z)\right] + G(z)\varphi(z) = 0,
\end{equation}
where $z$ is the complex variable and the functions $G$ and $K$ are assumed analytic in some nonempty open set $\Omega \subset \mathbb C$ such that $K$ does not vanish in $\Omega$. Equation~\eqref{eq:Hille} can be written in $\Omega$ as a second-order system on the unknown functions $\varphi_1(z)=\varphi(z)$ and $\varphi_2(z) = K(z) \frac{\df\varphi}{\df z}(z)$, i.e.,
\begin{equation}
\label{eq:Hille-system}
\left\{
\begin{aligned}
\varphi_1^\prime(z) & = \frac{\varphi_2(z)}{K(z)}, \\
\varphi_2^\prime(z) & = - G(z) \varphi_1(z).
\end{aligned}
\right.
\end{equation}
Any solution $(\varphi_1, \varphi_2)$ of \eqref{eq:Hille-system} is necessarily analytic in its domain of definition $\widehat\Omega \subset \Omega$, where $\widehat\Omega$ is assumed to be open.

Let $\gamma \in \mathcal C^1([a, b], \widehat\Omega)$, where $a, b$ are real numbers with $a < b$, and denote the endpoints of $\gamma$ by $z_1 = \gamma(a)$ and $z_2 = \gamma(b)$. Multiplying the first equation in \eqref{eq:Hille-system} by $\overline{\varphi_2(z)}$, the second one by $\overline{\varphi_1(z)}$, and integrating along $\gamma$, we obtain that
\[
\begin{aligned}
\int_\gamma \varphi_1^\prime(z) \overline{\varphi_2(z)} \diff z & = \int_\gamma \frac{\abs{\varphi_2(z)}^2}{K(z)} \diff z, \\
\int_\gamma \overline{\varphi_1(z)} \varphi_2^\prime(z) \diff z & = - \int_\gamma G(z) \abs{\varphi_1(z)}^2 \diff z.
\end{aligned}
\]
On the other hand, we have
\[
\int_\gamma \overline{\varphi_1(z)} \varphi_2^\prime(z) \diff z + \overline{\int_\gamma \varphi_1^\prime(z) \overline{\varphi_2(z)} \diff z} = \overline{\varphi_1(z_2)} \varphi_2(z_2) - \overline{\varphi_1(z_1)} \varphi_2(z_1).
\]
Hence
\begin{equation}\label{Green}
\overline{\varphi_1(z_2)} \varphi_2(z_2) - \overline{\varphi_1(z_1)} \varphi_2(z_1) = \overline{\int_\gamma \frac{\abs{\varphi_2(z)}^2}{K(z)} \diff z} - \int_\gamma G(z) \abs{\varphi_1(z)}^2 \diff z.
\end{equation}
Equation~\eqref{Green} is known as the \emph{Green--Hille transform} of \eqref{eq:Hille}.

With this preliminary exposition of Hille's approach, we are in position to provide the proof of Proposition~\ref{PropTsvetkoffCorrected}.

\begin{proof}[Proof of Proposition~\ref{PropTsvetkoffCorrected}]
Notice first that the function $G_{\kappa, \mu}$ defined by
\[
G_{\kappa, \mu}(z) = -\frac{1}{4} + \frac{\kappa}{z} + \frac{\frac{1}{4} - \mu^2}{z^2},
\]
which appears in the Whittaker equation \eqref{WhittakerODE}, is analytic in the set $\Omega = \mathbb C \setminus \{0\}$. Recall that, even though the Whittaker function $M_{\kappa, \mu}$ is possibly multi-valued, its zeros do not depend on which branch of $M_{\kappa, \mu}$ is considered. Applying the Green--Hille transform \eqref{Green} to the Whittaker equation \eqref{WhittakerODE}, we obtain that, for every branch of $M_{\kappa, \mu}$, for every open subset $\widehat\Omega \subset \Omega$ on which the considered branch of $M_{\kappa, \mu}$ is analytic, for every real numbers $a, b$ with $a < b$, and for every $\gamma \in \mathcal C^1([a, b], \widehat\Omega)$, we have
\begin{equation}\label{GreenW}
\overline{M_{\kappa, \mu}(z_2)} M^\prime_{\kappa, \mu}(z_2) - \overline{M_{\kappa, \mu}(z_1)} M^\prime_{\kappa, \mu}(z_1) = \overline{\int_\gamma \abs{M^\prime_{\kappa, \mu}(z)}^2 \diff z} - \int_\gamma G_{\kappa, \mu}(z) \abs{M_{\kappa, \mu}(z)}^2 \diff z,
\end{equation}
where $z_1 = \gamma(a)$ and $z_2 = \gamma(b)$.

Fix a branch of $M_{\kappa, \mu}$ which is analytic in a set $\widehat\Omega$ equal to the complex plane minus a half-line starting at the origin and such that $z_0 \in \widehat\Omega$. For an arbitrary $\varepsilon \in (0, 1)$, consider the curve $\gamma \in \mathcal C^1([\varepsilon, 1], \widehat\Omega)$ corresponding to the straight line from $\varepsilon z_0$ to $z_0$ with constant speed, i.e, $\gamma(t) = t z_0$ for $t \in [\varepsilon, 1]$. Then \eqref{GreenW} reads
\[
-\overline{M_{\kappa, \mu}(\varepsilon z_0)} M^\prime_{\kappa, \mu}(\varepsilon z_0) = \overline z_0 \int_{\varepsilon}^1 \abs{M^\prime_{\kappa, \mu}(t z_0)}^2 \diff t - z_0 \int_\varepsilon^1 G_{\kappa, \mu}(t z_0) \abs{M_{\kappa, \mu}(t z_0)}^2 \diff t.
\]
Using \eqref{eq:def-Whittaker} and \eqref{eq:Kummer-in-zero}, we obtain that $M_{\kappa, \mu}(z) = O(\abs{z}^{\mu + \frac{1}{2}})$ and $M^\prime_{\kappa, \mu}(z) = O(\abs{z}^{\mu - \frac{1}{2}})$ as $z \to 0$. Since $\mu > 0$, letting $\varepsilon \to 0$ in the above equality, we deduce that
\[
\overline z_0 \int_0^1 \abs{M^\prime_{\kappa, \mu}(t z_0)}^2 \diff t = z_0 \int_0^1 G_{\kappa, \mu}(t z_0) \abs{M_{\kappa, \mu}(t z_0)}^2 \diff t,
\]
i.e.,
\begin{equation}
\label{eq:main-Tsvetkov}
\overline z_0 \int_0^1 \abs{M^\prime_{\kappa, \mu}(t z_0)}^2 \diff t = \int_0^1 \left(-\frac{z_0}{4} + \frac{\kappa}{t} + \frac{\frac{1}{4} - \mu^2}{t^2 z_0}\right) \abs{M_{\kappa, \mu}(t z_0)}^2 \diff t.
\end{equation}

In the sequel of the proof, we set $x_0 = \Real z_0$ and $y_0 = \Imag z_0$. Note that \eqref{eq:main-Tsvetkov} can be rewritten as
\begin{equation}
\label{eq:Tsvetkov-a-b-c}
A \overline z_0 + B z_0 + \frac{D \left(\mu^2 - \frac{1}{4}\right)}{z_0} = C \kappa,
\end{equation}
where the positive constants $A, B, C, D$ are given by
\begin{align*}
A & = \int_0^1 \abs{M^\prime_{\kappa, \mu}(t z_0)}^2 \diff t, & B & = \frac{1}{4} \int_0^1 \abs{M_{\kappa, \mu}(t z_0)}^2 \diff t, \\
C & = \int_0^1 \frac{\abs{M_{\kappa, \mu}(t z_0)}^2}{t} \diff t, & D & = \int_0^1 \frac{\abs{M_{\kappa, \mu}(t z_0)}^2}{t^2} \diff t.
\end{align*}
Taking the real part of \eqref{eq:Tsvetkov-a-b-c}, we deduce that
\[
\left(A + B + \frac{D \left(\mu^2 - \frac{1}{4}\right)}{x_0^2 + y_0^2}\right) x_0 = C \kappa.
\]
In particular, if $\mu \geq \frac{1}{2}$, we deduce that $x_0$ and $\kappa$ must have the same sign, showing \ref{PropTsvetkoffCorrected-Sign}. Moreover, if $\kappa \neq 0$, then $C \kappa \neq 0$, showing that $x_0 \neq 0$ and giving thus the first part of the conclusion of \ref{PropTsvetkoffCorrected-k-neq-0}.

Let us now prove the inequality in \ref{PropTsvetkoffCorrected-k-neq-0}. For $\kappa \neq 0$, taking the real part of \eqref{eq:main-Tsvetkov} and dividing by $x_0$, we obtain that
\[
\int_0^1 \abs{M^\prime_{\kappa, \mu}(t z_0)}^2 \diff t = \int_0^1 \frac{- \abs{z_0}^2 t^2 + 4 \kappa \frac{\abs{z_0}^2}{x_0} t - (4 \mu^2 - 1)}{4 t^2 \abs{z_0}^2} \abs{M_{\kappa, \mu}(t z_0)}^2 \diff t.
\]
Since the left-hand side of the above equality is positive, so must be its right-hand side, and in particular the integrand in the right-hand side is positive in a set of points $t$ of positive measure. Since $4 t^2 \abs{z_0}^2 > 0$ for every $t \in (0, 1)$, there exists $t_\ast \in (0, 1)$ such that
\[
- \abs{z_0}^2 t_\ast^2 + 4 \kappa \frac{\abs{z_0}^2}{x_0} t_\ast - (4 \mu^2 - 1) > 0.
\]
As the left-hand side of the above inequality is a polynomial of degree $2$ in $t_\ast$ with negative leading coefficient, its discriminant must necessarily be positive, which shows, using also that $z_0 \neq 0$, that
\[
4 \kappa^2 \frac{\abs{z_0}^2}{x_0^2} > (4 \mu^2 - 1),
\]
yielding the conclusion.

To prove \ref{PropTsvetkoffCorrected-z0-not-real}, we assume that $y_0 \neq 0$ and we take the imaginary part of \eqref{eq:main-Tsvetkov}, obtaining, after dividing by $-y_0$, that
\[
\int_0^1 \abs{M^\prime_{\kappa, \mu}(t z_0)}^2 \diff t = \int_0^1 \left(\frac{1}{4} - \frac{4 \mu^2 - 1}{4 t^2 \abs{z_0}^2}\right) \abs{M_{\kappa, \mu}(t z_0)}^2 \diff t.
\]
Since the left-hand side of the above equality is positive, so is its right-hand side, and thus there must exist $t^\ast \in (0, 1)$ such that
\[
\frac{1}{4} - \frac{4 \mu^2 - 1}{4 t_\ast^2 \abs{z_0}^2} > 0.
\]
Since $t \mapsto \frac{1}{4} - \frac{4 \mu^2 - 1}{4 t^2 \abs{z_0}^2}$ is increasing in $(0, 1]$, we also have
\[
\frac{1}{4} - \frac{4 \mu^2 - 1}{4 \abs{z_0}^2} > 0,
\]
yielding the conclusion.
\end{proof}

Now that we have proved Proposition~\ref{PropTsvetkoffCorrected}, a few remarks are in order.
\begin{remark}
If $\mu, \kappa \in \mathbb R$ are such that $4\mu^2 - 4\kappa^2 < 1$, the inequality in Proposition~\ref{PropTsvetkoffCorrected}\ref{PropTsvetkoffCorrected-k-neq-0} is satisfied for every $z_0 \in \mathbb C \setminus \{0\}$. That item only provides nontrivial information on roots $z_0$ of $M_{\kappa, \mu}$ when $4\mu^2 - 4\kappa^2 \geq 1$. Similarly, the inequality in Proposition~\ref{PropTsvetkoffCorrected}\ref{PropTsvetkoffCorrected-z0-not-real} is satisfied for every $z_0 \in \mathbb C \setminus \{0\}$ when $\mu \in (0, \frac{1}{2}]$.
\end{remark}

\begin{remark}
Even though we are reporting here the results from \cite{Boussaada2022Some}, the statement of Proposition~\ref{PropTsvetkoffCorrected} presented here is slightly more general than its counterpart \cite[Proposition~3.3]{Boussaada2022Some}. Indeed, \cite[Proposition~3.3]{Boussaada2022Some} always assumes that $\mu \geq \frac{1}{2}$. In Proposition~\ref{PropTsvetkoffCorrected}, we assume $\mu \geq \frac{1}{2}$ in \ref{PropTsvetkoffCorrected-Sign} and the inequalities in \ref{PropTsvetkoffCorrected-k-neq-0} and \ref{PropTsvetkoffCorrected-z0-not-real} are trivial unless $\mu \geq \frac{1}{2}$, but Proposition~\ref{PropTsvetkoffCorrected}\ref{PropTsvetkoffCorrected-k-neq-0} brings the additional information that $\kappa \neq 0$ implies $\Real z_0 \neq 0$ for all $\mu > 0$, and not only $\mu \geq \frac{1}{2}$.
\end{remark}

\begin{remark}
It turns out that the inequality $\abs{z} > \sqrt{4 \mu^2-1}$ on roots of $M_{\kappa, \mu}$ holds, as shown in Proposition~\ref{PropTsvetkoffCorrected}, for \emph{nonreal} roots of $M_{\kappa, \mu}$, but not for \emph{nontrivial} roots of $M_{\kappa, \mu}$, as stated in Proposition~\ref{CoroSaffVarga}, and, indeed, the counterexample provided in Counterexample~\ref{ExplFalseSaffVarga} presents a \emph{real} root which does not satisfy that inequality. The proof provided in \cite[Corollary~3.3]{Saff1978Zeros} for Proposition~\ref{CoroSaffVarga} above is based on the inequality $\abs{\Real z} > 2\abs{\kappa}$ which follows from Proposition~\ref{PropTsvetkoff1}, and hence the issue in the statement of Proposition~\ref{CoroSaffVarga} can be seen as a propagation of the issue in the statement of Proposition~\ref{PropTsvetkoff1}.
\end{remark}

\begin{remark}
Proposition~2.1 in \cite{Saff1978Zeros} makes use of Proposition~\ref{CoroSaffVarga} in its proof. However, its conclusion is still correct. Indeed, in the proof of \cite[Proposition~2.1]{Saff1978Zeros}, the authors apply Proposition~\ref{CoroSaffVarga} to the Whittaker function $M_{\kappa, \mu}$ with $\kappa = \frac{n - \nu}{2}$ and $\mu = \frac{n + \nu + 1}{2}$, where $n$ and $\nu$ are nonnegative integers with $n + \nu > 0$. In that case $M_{\kappa, \mu}$ does not admit nontrivial real roots. Indeed, if $n = \nu$, this is simply a consequence of Proposition~\ref{PropTsvetkoffCorrected}\ref{PropTsvetkoffCorrected-Sign} while, if $n \neq \nu$, we have $4 \mu^2 - 4 \kappa^2 - 1 = 2n + 2\nu + 4 n \nu > 0$, and thus, by Proposition~\ref{PropTsvetkoffCorrected}\ref{PropTsvetkoffCorrected-k-neq-0}, if $M_{\kappa, \mu}$ admitted a nontrivial real root $x$, we would have $x^2 \left(4{\mu}^{2}-4{\kappa}^{2} - 1\right) < 0$, which is a contradiction. In particular, all nontrivial roots of $M_{\kappa, \mu}$ are nonreal, and thus, by Proposition~\ref{PropTsvetkoffCorrected}, any such root $z$ satisfies the bound $\abs{z} > \sqrt{4 \mu^2 - 1}$, as required.
\end{remark}

Thanks to the connection between Whittaker and Kummer functions from \eqref{eq:def-Whittaker}, an immediate consequence of Proposition~\ref{PropTsvetkoffCorrected} on the location of zeros of Kummer functions is given in the following corollary.

\begin{corollary}
\label{CorZerosKummer}
Let $a, b\in \mathbb R$ be such that $b > 1$ and $z_0 \in \mathbb C \setminus \{0\}$ be such that $M(a, b, z_0) = 0$.
\begin{enumerate}
\item\label{CorZerosKummer-sign} If $b \geq 2$, then $\sign(b - 2a) = \sign(\Real z_0)$.
\item If $b \neq 2 a$, then $\Real z_0 > 0$ and $(b - 2a)^2 (\Imag z_0)^2 > (\Real z_0)^{2} \left(4a(b-a) - 2b\right)$.
\item If $z_0 \notin \mathbb R$, then $\abs{z_0}^2 > b(b - 2)$.
\end{enumerate}
\end{corollary}

\begin{remark}
\label{remk:Wynn}
In the case $a \in \{\alpha, \alpha + 1\}$ and $b = 2\alpha + 1$ for some $\alpha > -\frac{1}{2}$, the conclusion of Corollary~\ref{CorZerosKummer}\ref{CorZerosKummer-sign} was shown by P.~Wynn in \cite[Theorem~1]{Wynn1973Zeros}. The techniques used in that reference do not rely on Hille's approach, but use instead a continued fraction representation of a ratio of Kummer functions (see \cite[Theorem~B(ii)]{Wynn1973Zeros}). Note that Corollary~\ref{CorZerosKummer} does not cover all cases of Wynn's result, since the assumption $b \geq 2$ is equivalent to $\alpha \geq \frac{1}{2}$, whereas Wynn's result is obtained for all $\alpha > -\frac{1}{2}$.
\end{remark}

\begin{figure}[ht]
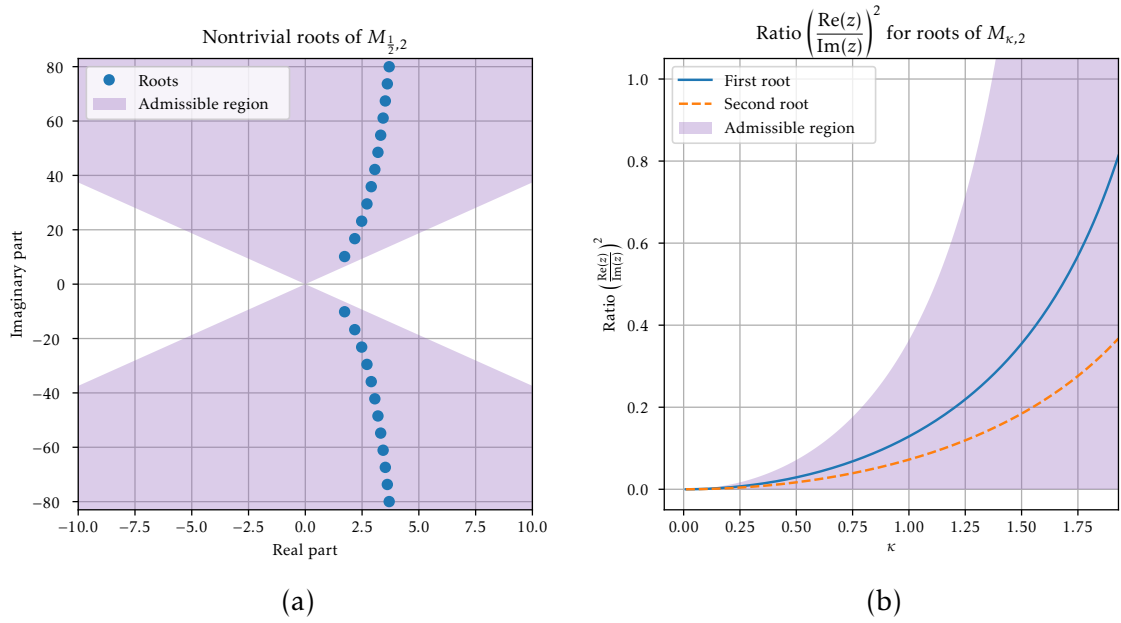

\centering
\begin{tabular}{@{} c @{} c @{}}
\resizebox{0.5\textwidth}{!}{\input{Figures/WhittakerBound.pgf}} & \resizebox{0.5\textwidth}{!}{\input{Figures/WhittakerBoundK.pgf}} \tabularnewline
(a) & (b) \tabularnewline
\end{tabular}
\caption[Nontrivial roots of the Whittaker function $M_{\frac{1}{2}, 2}$ and ratio $\left(\frac{\Real z}{\Imag z}\right)^2$ for the first two roots of $M_{\kappa, 2}$ as a function of $\kappa$]{(a) Nontrivial roots of the Whittaker function $M_{\frac{1}{2}, 2}$ (blue circles) and region of the complex plane (in light violet) where the bound from Proposition~\ref{PropTsvetkoffCorrected}\ref{PropTsvetkoffCorrected-k-neq-0} is satisfied. (b) Ratio $\left(\frac{\Real z}{\Imag z}\right)^2$ for the first two roots of $M_{\kappa, 2}$ as a function of $\kappa$ (solid blue and dashed orange lines) and region (in light violet) described by the bound from Proposition~\ref{PropTsvetkoffCorrected}\ref{PropTsvetkoffCorrected-k-neq-0}}
\label{FigWhittakerBound}
\end{figure}

To conclude this section, let us provide a few illustrations of the inequality from Proposition~\ref{PropTsvetkoffCorrected}\ref{PropTsvetkoffCorrected-k-neq-0}\footnote{For these illustrations, we have used the implementation of Kummer functions provided by the function \texttt{hyp1f1} from Python's \texttt{scipy.special} module \cite{Virtanen2020SciPy} and all roots were computed numerically using Python's \texttt{cxroots} module \cite{Parini2018cxroots}.}. When $4(\mu^2 - \kappa^2) \geq 1$, Proposition~\ref{PropTsvetkoffCorrected}\ref{PropTsvetkoffCorrected-k-neq-0} provides a nontrivial bound on the location of nonzero roots of $M_{\kappa, \mu}$ in the complex plane. Figure~\ref{FigWhittakerBound}(a) illustrates that bound for the function $M_{\frac{1}{2}, 2}$, in which case the bound from Proposition~\ref{PropTsvetkoffCorrected}\ref{PropTsvetkoffCorrected-k-neq-0} reads $(\Imag z)^2 > 14 (\Real z)^2$. The region of the complex plane where that bound is satisfied is represented in light violet in the figure, whereas the roots of $M_{\frac{1}{2}, 2}$ in the rectangle $\{z \in \mathbb C \suchthat \abs{\Real z} \leq 10 \text{ and } \abs{\Imag z} \leq 80\}$ are represented by blue circles.

We also represent, in Figure~\ref{FigWhittakerBound}(b), the ratio $\left(\frac{\Real z}{\Imag z}\right)^2$ for the first two roots (ordered according to their distance to the origin and considering only once each complex conjugate pair) of $M_{\kappa, \mu}$ when $\mu = 2$ and $\kappa$ varies in the interval $\left(0, \frac{\sqrt{15}}{2}\right)$, on which the bound from Proposition~\ref{PropTsvetkoffCorrected}\ref{PropTsvetkoffCorrected-k-neq-0} is nontrivial. Proposition~\ref{PropTsvetkoffCorrected}\ref{PropTsvetkoffCorrected-k-neq-0} states in this case that $\left(\frac{\Real z}{\Imag z}\right)^2 < \frac{4 \kappa^2}{4 (\mu^2 - \kappa^2) - 1}$ for any root $z$ of $M_{\kappa, \mu}$, and the region described by this inequality is represented in light violet color in the figure.

\section{MID for real roots of maximal multiplicity}
\label{sec:MID-maximal}

Using the preliminary results from Section~\ref{sec:MID-prelim}, we will now show the validity of the MID property for quasipolynomials of the form \eqref{eq:Delta-DDE}, i.e., that, if \eqref{eq:Delta-DDE} admits a root of maximal multiplicity, then necessarily it is dominant. We will also provide consequences of this result for the exponential stability of some time-delay systems.

\subsection{Main result}

As presented in Section~\ref{sec:intro-MID}, the MID property was established in some works in the literature for low-order cases, i.e., for small values of $n$ and $m$ in \eqref{eq:Delta-DDE}, using techniques that are not easy to generalize to arbitrary $n$ and $m$ with $m \leq n$. This latter case was only treated after it was shown in \cite{Mazanti2020Qualitative, Benarab2020MID, Mazanti2021Effects}, for low-order cases, that the quasipolynomial \eqref{eq:Delta-DDE} with a root of maximal multiplicity can be factorized in terms of a Kummer function. Such a factorization also works for arbitrary $n$ and $m$, and the first work to exploit it to prove the MID property in a more general context was \cite{Mazanti2021Multiplicity}, which only considered the case $m = n-1$, in which case the Kummer function appearing in the factorization is $M(n, 2n + 1, \cdot)$. The MID property in this case then follows from the result of \cite{Wynn1973Zeros} on the the location of the roots of this Kummer function (see Remark~\ref{remk:Wynn} above). After Proposition~\ref{PropTsvetkoffCorrected} and Corollary~\ref{CorZerosKummer} above were proved in \cite{Boussaada2022Some}, the MID property in full generality for \eqref{eq:Delta-DDE} was established as a consequence in \cite{Boussaada2022Generic}. Later, \cite{Boussaada2022Pade} improved some parts of the proof from \cite{Boussaada2022Generic} by exploiting the links between the existence of a root of maximal multiplicity of \eqref{eq:Delta-DDE} and the error term in the Padé approximation of the exponential function. It is this last, simpler version of the proof that we present here.

Our main result on the multiplicity-induced dominance for \eqref{eq:Delta-DDE} is the following, which highlights the links between roots of maximal multiplicities of the quasipolynomial $\Delta$ from \eqref{eq:Delta-DDE}, Kummer functions, and Padé approximations of the exponential function. More precisely, our main result characterizes the case where $\Delta$ admits a root of maximal multiplicity, equal to its degree $D$, through a suitable fatorization in terms of a Kummer function by using the results on the Padé approximants for the exponential function recalled in Section~\ref{sec:Pade}.

\begin{theorem}
\label{thm:main-MID}
Let $n, m$ be nonnegative integers, $\tau > 0$, $a_0, \dotsc, a_{n-1}, \alpha_0, \dotsc, \alpha_m$ be real numbers, and consider the quasipolynomial $\Delta$ given by \eqref{eq:Delta-DDE}. Let $D = n + m + 1$ denote the degree of $\Delta$ and take $s_0 \in \mathbb R$. Let $P$ and $Q$ be the polynomials defined by
\[
P(s) = s^n + \sum_{k=0}^{n-1} a_k s^k, \qquad Q(s) = \sum_{k=0}^m \alpha_k s^k
\]
and define the quasipolynomial $\widetilde\Delta$ and the polynomials $\widetilde P$ and $\widetilde Q$ by setting, for $z \in \mathbb C$,
\begin{subequations}
\begin{align}
\widetilde\Delta(z) & = \tau^n \Delta\left(s_0 + \tfrac{z}{\tau}\right), \label{eq:DeltaDeltaTilde}\\
\widetilde P(z) & = \tau^n P\left(s_0 + \tfrac{z}{\tau}\right), \label{eq:PPtilde} \\
\widetilde Q(z) & = e^{-s_0 \tau} \tau^n Q\left(s_0 + \tfrac{z}{\tau}\right). \label{eq:QQtilde}
\end{align}
\end{subequations}
Then the following assertions are equivalent.
\begin{enumerate}
\item\label{item:max-mult} The real number $s_0$ is a root of multiplicity $D$ of $\Delta$.
\item\label{item:Pade} The pair $(-\widetilde Q, \widetilde P)$ is a Padé approximant of the exponential function $z \mapsto e^z$ of order $(m, n)$.
\item\label{item:Coeffs} The coefficients $a_0, \dotsc, a_{n-1}, \alpha_0, \dotsc, \alpha_m$ of $\Delta$ are given by
\begin{equation}
\label{Coeffs}
\left\{
\begin{aligned}
a_k & = (-1)^{n-k} n! \sum_{j=k}^{n} \frac{\binom{j}{k} \binom{n + m - j}{m} s_0^{j-k}}{j! {\tau}^{n-j}}, & \qquad & \text{ for } k \in \llbracket 0, n-1\rrbracket, \\
\alpha_k & = (-1)^{n-1} \frac{e^{s_0 \tau}}{k!} \sum_{j=k}^{m} \frac{(n + m - j)! (-s_0)^{j-k}}{(m - j)! (j - k)! \tau^{n-j}}, & & \text{ for } k \in \llbracket 0, m\rrbracket.
\end{aligned}
\right.
\end{equation}
\item\label{item:DeltaTilde} The quasipolynomial $\widetilde\Delta$ is given by
\begin{equation}
\label{eq:DeltaTilde-as-Kummer}
\widetilde\Delta(z) = \frac{n! z^{D}}{D !} M(m + 1, D + 1, -z).
\end{equation}
\end{enumerate}
In addition, if the above equivalent assertions are satisfied, then
\begin{enumerate}[label={\textnormal{(\roman*)}}, ref={(\roman*)}, widest*=3]
\item\label{item:m-less-n} if $m < n$, we have $\Real s < s_0$ for every $s \in \mathbb C \setminus \{s_0\}$ such that $\Delta(s) = 0$;
\item if $m = n$, we have $\Real s = s_0$ for every $s \in \mathbb C$ such that $\Delta(s) = 0$;
\item\label{item:m-greater-n} if $m > n$, we have $\Real s > s_0$ for every $s \in \mathbb C \setminus \{s_0\}$ such that $\Delta(s) = 0$.
\end{enumerate}
In particular, if $m \leq n$, then $s_0$ is dominant.
\end{theorem}

\begin{proof}
Assume that \ref{item:max-mult} is satisfied. Thanks to \eqref{eq:DeltaDeltaTilde}, $0$ is a root of multiplicity $D$ of $\widetilde\Delta$. We claim that $\widetilde P(0) \neq 0$. Indeed, assume, to obtain a contradiction, that $\widetilde P(0) = 0$. Since $\widetilde\Delta(z) = \widetilde P(z) + e^{-z} \widetilde Q(z)$ and $\widetilde\Delta(0) = 0$, we deduce that $\widetilde Q(0) = 0$. Hence, the functions $\widehat P$ and $\widehat Q$ defined respectively by $\widehat P(z) = \frac{\widetilde P(z)}{z}$ and $\widehat Q(z) = \frac{\widetilde Q(z)}{z}$ are both polynomials, with $\deg \widehat P = n - 1$ and $\deg \widehat Q \leq m - 1$. In particular, the function $\widehat\Delta$ defined by $\widehat\Delta(z) = \frac{\widetilde\Delta(z)}{z} = \widehat P(z) + e^{-z} \widehat Q(z)$ is a quasipolynomial of degree at most $n + m - 1 = D - 2$. On the other hand, since $\widetilde\Delta(z) = z \widehat\Delta(z)$ and $0$ is a root of multiplicity $D$ of $\widetilde\Delta$, we deduce that $0$ must be a root of multiplicity $D - 1$ of $\widehat\Delta$, contradicting Corollary~\ref{coro:Polya-Szego}. This establishes that $\widetilde P(0) \neq 0$.

Now, since $0$ is a root of multiplicity $D = n + m + 1$ of $\widetilde\Delta$, we have $\widetilde\Delta(z) = O(z^{n + m + 1})$ as $z \to 0$, that is,
\begin{equation}
\label{eq:tildeP-tildeQ-Pade}
\widetilde P(z) e^{z} + \widetilde Q(z) = O(z^{n + m + 1}) \qquad \text{ as } z \to 0.
\end{equation}
Since $\widetilde P(0) \neq 0$, it follows from Remark~\ref{remk:Pade-alternative} that $(-\widetilde Q, \widetilde P)$ is a Padé approximant of $z \mapsto e^z$ of order $(m, n)$, implying \ref{item:Pade}.

Conversely, if \ref{item:Pade} is satisfied, then $\widetilde P(0) \neq 0$ and, by Remark~\ref{remk:Pade-alternative}, \eqref{eq:tildeP-tildeQ-Pade} is satisfied, showing that $0$ is a root of multiplicity at least $D = n + m + 1$ of $\widetilde\Delta$. By Corollary~\ref{coro:Polya-Szego}, we deduce that $0$ is a root of multiplicity $D$ of $\widetilde\Delta$, yielding \ref{item:max-mult} thanks to \eqref{eq:DeltaDeltaTilde}.

Let us now prove the equivalence between \ref{item:Pade} and \ref{item:Coeffs}. Recall that the polynomial $\widetilde P$ from \eqref{eq:PPtilde} is monic. By Proposition~\ref{prop:Pade-exp}, \ref{item:Pade} is equivalent to saying that $(-\widetilde Q, \widetilde P)$ is proportional to the pair of polynomials $(Q, P)$ from the statement of that proposition, with a proportionality constant chosen in such a way as for $\widetilde P$ to be monic. That is, \ref{item:Pade} holds true if and only if
\begin{align*}
\widetilde Q(z) & = -\lambda \sum_{k=0}^m \frac{m! (n+m-k)!}{(m-k)! (n+m)!} \frac{z^k}{k!},\\
\widetilde P(z) & = \lambda \sum_{k=0}^n \frac{n! (n+m-k)!}{(n-k)! (n+m)!} \frac{(-z)^k}{k!},
\end{align*}
where $\lambda \in \mathbb R^\ast$ is the unique constant such that $\widetilde P$ is monic. A straightforward computation shows that $\lambda = (-1)^n \frac{(n+m)!}{m!}$, and thus \ref{item:Pade} is equivalent to
\begin{subequations}
\begin{align}
\widetilde Q(z) & = (-1)^{n-1} \sum_{k=0}^m \frac{(n + m - k)!}{(m - k)!} \frac{z^k}{k!}, \label{eq:explicit-Qtilde}\\
\widetilde P(z) & = (-1)^{n} \sum_{k=0}^n \frac{n! (n + m - k)!}{(n - k)! m!} \frac{(-z)^k}{k!}. \label{eq:explicit-Ptilde}
\end{align}
\end{subequations}
Now, if $\widetilde Q$ is given by \eqref{eq:explicit-Qtilde}, we have from \eqref{eq:QQtilde} that
\begin{subequations}
\begin{align}
Q(s) & = \frac{e^{s_0 \tau}}{\tau^n} \widetilde Q(\tau(s - s_0)) = (-1)^{n-1} \frac{e^{s_0 \tau}}{\tau^n} \sum_{j=0}^m \frac{(n + m - j)!}{(m - j)!} \frac{\tau^j (s - s_0)^j}{j!} \notag \\
& = (-1)^{n-1} e^{s_0 \tau} \sum_{j=0}^m \frac{(n + m - j)!}{(m - j)! j! \tau^{n-j}} \sum_{k=0}^j \frac{j!}{k! (j-k)!} s^k (-s_0)^{j-k} \notag \\
& = (-1)^{n-1} e^{s_0 \tau} \sum_{k=0}^m \frac{1}{k!} \left[\sum_{j=k}^m \frac{(n + m - j)! (-s_0)^{j-k}}{(m - j)! (j-k)! \tau^{n-j}}\right] s^k \label{eq:explicit-Q}
\end{align}
and, conversely, if $Q$ is given by the last expression in \eqref{eq:explicit-Q}, then, by \eqref{eq:QQtilde}, $\widetilde Q$ is given by \eqref{eq:explicit-Qtilde}. Similarly, if $\widetilde P$ is given by \eqref{eq:explicit-Ptilde}, we have from \eqref{eq:PPtilde} that
\begin{align}
P(s) & = \frac{1}{\tau^n} \widetilde P(\tau(s - s_0)) = \frac{(-1)^{n}}{\tau^n} \sum_{j=0}^n \frac{(-1)^j n! (n + m - j)!}{(n - j)! m!} \frac{\tau^j (s - s_0)^j}{j!} \notag \\
& = \sum_{j=0}^n \frac{(-1)^{n-j} n! (n + m - j)!}{\tau^{n-j} (n - j)! m! j!} \sum_{k=0}^j \frac{j!}{k! (j-k)!} s^k (-s_0)^{j-k} \notag \\
& = \sum_{k=0}^n (-1)^{n-k} n! \left[\sum_{j=k}^n \frac{\binom{j}{k}\binom{n + m - j}{m} s_0^{j-k}}{j! \tau^{n-j}}\right] s^k \label{eq:explicit-P}
\end{align}
\end{subequations}
and, conversely, if $P$ is given by the last expression in \eqref{eq:explicit-P}, then, by \eqref{eq:PPtilde}, $\widetilde P$ is given by \eqref{eq:explicit-Ptilde}. Equations~\eqref{eq:explicit-Q} and \eqref{eq:explicit-P} are clearly equivalent to \eqref{Coeffs}, and hence \ref{item:Pade} is equivalent to \ref{item:Coeffs}.

If \ref{item:Pade} is satisfied, then, by Proposition~\ref{prop:Pade-exp}, we deduce that
\[
e^z \widetilde\Delta(z) = \widetilde P(z) e^z + \widetilde Q(z) = (-1)^n \lambda \frac{z^{D}}{(D-1)!} \int_0^1 t^n (1-t)^m e^{z t} \diff t,
\]
where $\lambda = (-1)^n \frac{(n+m)!}{m!}$ is as before. By Proposition~\ref{prop:Kummer-integral} and \eqref{eq:Kummer-transf-1}, we deduce that $\widetilde\Delta$ is given by \eqref{eq:DeltaTilde-as-Kummer}, and so \ref{item:DeltaTilde} is satisfied. Conversely, if \ref{item:DeltaTilde} is satisfied, then $0$ is a root of multiplicity at least $D$ of $\widetilde\Delta$ and, since $\widetilde\Delta$ is a quasipolynomial of degree $D$, it follows from Corollary~\ref{coro:Polya-Szego} that $0$ is a root of multiplicity $D$ of $\widetilde\Delta$, showing that \ref{item:max-mult} is satisfied.

Finally, if \ref{item:max-mult}--\ref{item:DeltaTilde} are satisfied, then \ref{item:m-less-n}--\ref{item:m-greater-n} follow as an immediate consequence of \eqref{eq:DeltaTilde-as-Kummer} and Corollary~\ref{CorZerosKummer}\ref{CorZerosKummer-sign} since $D + 1 - 2(m+1) = n - m$.
\end{proof}

A series of remarks are now in order.

\begin{remark}
The equivalence between \ref{item:max-mult} and \ref{item:Coeffs} was proved in \cite[Theorem~3.1]{Mazanti2021Multiplicity} in the case $m = n - 1$ by some long and tedious computations making use of several special relations satisfied by binomial coefficients, presented in \cite[Section~2.4]{Mazanti2021Multiplicity}. A simplification of this proof was proposed in \cite[Theorem~12]{Boussaada2022Generic}, now in the general case of any $n$ and $m$, making use of the factorization in \ref{item:DeltaTilde} in order to get the equivalence between \ref{item:max-mult} and \ref{item:Coeffs}. The proof we presented above exploits the Padé approximation of the exponential function from Proposition~\ref{prop:Pade-exp} to further simplify the proof, obtaining also as a byproduct the equivalence between \ref{item:max-mult}, \ref{item:Coeffs}, and \ref{item:DeltaTilde} with \ref{item:Pade}.
\end{remark}

\begin{remark}
The fact that $s_0 \in \mathbb R$ is only needed in the proof to ensure that the coefficients $a_0, \dotsc, a_{n-1}, \alpha_0, \dotsc, \alpha_m$ given by \eqref{Coeffs} are all real. Theorem~\ref{thm:main-MID} remains valid, with the exact same proof, if we consider $s_0 \in \mathbb C$ and allow for complex coefficients $a_0, \dotsc, a_{n-1}, \alpha_0, \dotsc, \alpha_m$.
\end{remark}

\begin{remark}
\label{remk:S0ANMinus1Tau}
By considering the first equation in \eqref{Coeffs} with $k = n - 1$, one obtains the simple and interesting relation between $s_0$, $\tau$, and $a_{n-1}$ given by
\begin{equation}
\label{RelationS0ANMinus1Tau}
s_0 = -\frac{a_{n-1}}{n} - \frac{m+1}{\tau}.
\end{equation}
\end{remark}

\begin{remark}
Let $\Delta$ be the quasipolynomial from \eqref{eq:Delta-DDE}. Assume that $\Delta$ has real coefficients and that it admits a root $s_0 \in \mathbb C$ of maximal multiplicity $D$. Then, as a consequence of \eqref{RelationS0ANMinus1Tau}, we necessarily have $s_0 \in \mathbb R$. In other words, nonreal roots of $\Delta$ cannot have maximal multiplicity $D$ if $\Delta$ has real coefficients.
\end{remark}

\begin{remark}
If the coefficients of the quasipolynomial $\Delta$ from \eqref{eq:Delta-DDE} are given by \eqref{Coeffs}, then $s_0$ is the unique real root of $\Delta$, and, more precisely, it is the unique root of $\Delta$ on the horizontal strip $\{s \in \mathbb C \suchthat \abs{\Imag(s)} < \frac{2\pi}{\tau}\}$ of the complex plane. Indeed, applying Proposition~\ref{prop:Polya-Szego} with $\alpha = 0$ and any $\beta \in (0, \frac{2\pi}{\tau})$, we deduce that the number of roots of $\Delta$ in the strip $S_\beta = \{s \in \mathbb C \suchthat 0 \leq \Imag(s) \leq \beta\}$, counted according to their multiplicity, is at most $D$, since $\frac{\tau(\beta - \alpha)}{2 \pi} < 1$ in this case. Since $s_0 \in S_\beta$ and its multiplicity as a root of $\Delta$ is exactly $D$, this shows that no other root of $\Delta$ can belong to $S_\beta$. The conclusion now follows since $\beta \in (0, \frac{2\pi}{\tau})$ is arbitrary, and by remarking that $\Delta$ has only real coefficients and thus its roots appear in complex-conjugate pairs.
\end{remark}

\subsection{Characterization of the spectrum in the neutral case}

In the neutral case, i.e., when $m = n$ in \eqref{eq:Delta-DDE}, Theorem~\ref{thm:main-MID} asserts that, under any of the equivalent assertions \ref{item:max-mult}--\ref{item:DeltaTilde}, all roots of the quasipolynomial $\Delta$ from \eqref{eq:Delta-DDE} have the same real part, equal to $s_0$. We can further characterize the imaginary parts of the roots of $\Delta$ in this case as solutions of a certain transcendental equation, as described in the next result.

\begin{proposition}
\label{prop:spectrum-neutral}
Let $n \in \mathbb N$, $\tau \in \mathbb R_+^\ast$, and $s_0 \in \mathbb R$. Consider the quasipolynomial $\Delta$ given by \eqref{eq:Delta-DDE} in the case $m = n$ and assume that its coefficients $a_0, \dotsc, a_{n-1}, \alpha_0, \dotsc, \alpha_n$ are given by \eqref{Coeffs}. Then the set of roots of $\Delta$ is $\{s_0 + i \frac{\zeta}{\tau} \suchthat \zeta \in \Xi_n\}$, where
\begin{equation*}
\Xi_n = \left\{\zeta \in \mathbb R \suchthat \tan\left(\frac{\zeta}{2}\right) = \frac{\displaystyle \zeta \sum_{\ell = 0}^{\floor*{\frac{n-1}{2}}} (-1)^\ell \frac{(2n - 2\ell - 1)!}{(2\ell + 1)! (n - 2\ell - 1)!} \zeta^{2\ell}}{\displaystyle \sum_{\ell = 0}^{\floor*{\frac{n}{2}}} (-1)^\ell \frac{(2n - 2\ell)!}{(2\ell)! (n - 2\ell)!} \zeta^{2\ell}}\right\}.
\end{equation*}
\end{proposition}

\begin{proof}
Let $\widetilde\Delta$ be defined from $\Delta$ as in \eqref{eq:DeltaDeltaTilde} and note that $\widetilde\Delta(z) = \widetilde P(z) + e^z \widetilde Q(z)$ for $z \in \mathbb C$, where the polynomials $\widetilde P$ and $\widetilde Q$ are defined from $\Delta$ as in the statement of Theorem~\ref{thm:main-MID}. By Theorem~\ref{thm:main-MID}, $(-\widetilde Q, \widetilde P)$ is a Padé approximation of $z \mapsto e^z$ of order $(n, n)$ and thus, by Proposition~\ref{prop:Pade-exp} and arguing as in the proof of Theorem~\ref{thm:main-MID}, we obtain that $\widetilde Q$ and $\widetilde P$ are given by \eqref{eq:explicit-Qtilde} and \eqref{eq:explicit-Ptilde}, respectively, with $m = n$, yielding that
\begin{equation}
\label{eq:DeltaTilde-neutral}
\widetilde\Delta(z) = \sum_{k=0}^n (-1)^{n - k} \frac{(2n - k)!}{k! (n - k)!} z^k + (-1)^{n - 1} e^{-z} \sum_{k = 0}^n \frac{(2n - k)!}{k! (n - k)!} z^k.
\end{equation}
It follows from Theorem~\ref{thm:main-MID} that all roots $z$ of $\widetilde\Delta$ lie in the imaginary axis. Using \eqref{eq:DeltaTilde-neutral}, we have that $z = i \zeta$ for $\zeta \in \mathbb R$ is a root of $\widetilde\Delta$ if and only if
\[
e^{i \frac{\zeta}{2}} \sum_{k=0}^n (-1)^{k} \frac{(2n - k)!}{k! (n - k)!} (i\zeta)^k = e^{-i\frac{\zeta}{2}} \sum_{k = 0}^n \frac{(2n - k)!}{k! (n - k)!} (i\zeta)^k.
\]
The right-hand side of the above equality is equal to the complex conjugate of the left-hand side. Hence, the above equality is equivalent to
\[
\Imag\left(e^{i \frac{\zeta}{2}}\sum_{k=0}^n (-1)^{k} \frac{(2n - k)!}{k! (n - k)!} (i\zeta)^k\right) = 0,
\]
which happens if and only if
\[
\sin\left(\frac{\zeta}{2}\right) \sum_{\ell = 0}^{\floor*{\frac{n}{2}}} (-1)^\ell \frac{(2n - 2\ell)!}{(2\ell)! (n - 2\ell)!} \zeta^{2\ell} = \zeta \cos\left(\frac{\zeta}{2}\right) \sum_{\ell = 0}^{\floor*{\frac{n-1}{2}}} (-1)^\ell \frac{(2n - 2\ell - 1)!}{(2\ell + 1)! (n - 2\ell - 1)!} \zeta^{2\ell},
\]
i.e., if and only if $\zeta \in \Xi_n$. The conclusion follows since, by \eqref{eq:DeltaDeltaTilde}, $z$ is a root of $\widetilde\Delta$ if and only if $s_0 + \frac{z}{\tau}$ is a root of $\Delta$.
\end{proof}

\begin{remark}
The proof of Proposition~\ref{prop:spectrum-neutral} was provided first in \cite[Theorem~4(c)]{Mazanti2021Effects} in the case $n = m = 1$, in which $\Xi_1 = \left\{\zeta \in \mathbb R \suchthat \tan\left(\frac{\zeta}{2}\right) = \frac{\zeta}{2}\right\}$. Its generalization to any $n \in \mathbb N$ was provided in \cite{Boussaada2022Generic}, from where we have extracted the above proof.
\end{remark}

\subsection{Application to the exponential stability of delay-differential equations}

Using the fact that \eqref{eq:Delta-DDE} is the characteristic equation of the delay-differential equation \eqref{eq:DDE}, we have the following immediate consequence of Theorem~\ref{thm:main-MID}, Remark~\ref{remk:S0ANMinus1Tau}, and the results on the stability of time-delay systems recalled in Section~\ref{sec:stability-LTI-TDS}.

\begin{corollary}
Let $n, m$ be nonnegative integers with $m \leq n$, $\tau > 0$, and $s_0 \in \mathbb R$. Consider the delay-differential equation \eqref{eq:DDE} with coefficients $a_0, \dotsc, a_{n-1}, \alpha_0, \dotsc, \alpha_m$ given by \eqref{Coeffs}. Then the trivial solution of \eqref{eq:DDE} is exponentially stable if and only if $s_0 < 0$ or, equivalently, $a_{n-1} > -\frac{n (m+1)}{\tau}$.
\end{corollary}

\subsection{Application to the exponential stability of a class of delay differential-algebraic equations}
\label{sec:MID-DDAE}

It turns out that equations of the form \eqref{eq:Delta-DDE} are also characteristic equations of a class of DDAEs, as it was already remarked in \cite{Mazanti2021Effects} for some low-order cases. The material we now present consists in a generalization of some arguments provided in \cite{Mazanti2021Effects}.

Consider a DDAE of the form
\begin{equation}
\label{eq:DDAE-MID}
\left\{
\begin{aligned}
& x^{(n)}(t) + \sum_{j=0}^{n-1} a_j x^{(j)}(t) = b \sum_{j = 1}^m k_j y_j(t - \tau), \\
& y_j(t) = \sum_{\ell = 0}^{n-1} c_{j,\ell} x^{(\ell)}(t) + d_j \sum_{j^\prime = 1}^m k_{j^{\prime}} y_{j^\prime}(t - \tau), \qquad j \in \llbracket 1, m\rrbracket,
\end{aligned}
\right.
\end{equation}
where $\tau > 0$ is a delay, $a_0, \dotsc, a_{n-1}, b, k_1, \dotsc, k_m, c_{1,0}, \dotsc, c_{m, n-1}, d_1, \dotsc, d_m$ are real coefficients, and the state of the system at time $t$ is made of the scalar $x(t) \in \mathbb R$ and the vector $y(t) = (y_1(t), \dotsc, y_m(t))^\transp \in \mathbb R^m$.

The motivation for considering \eqref{eq:DDAE-MID} comes from the output feedback control of the controlled DDAE
\begin{equation}
\label{eq:DDAE-open-loop}
\left\{
\begin{aligned}
& x^{(n)}(t) + \sum_{j=0}^{n-1} a_j x^{(j)}(t) = b u(t - \tau), \\
& y_j(t) = \sum_{\ell = 0}^{n-1} c_{j,\ell} x^{(\ell)}(t) + d_j u(t - \tau), \qquad j \in \llbracket 1, m\rrbracket,
\end{aligned}
\right.
\end{equation}
where $x(t) \in \mathbb R$ is the scalar state variable of the system at time $t$, $u(t) \in \mathbb R$ is the value at time $t$ of a scalar control, which can only be applied to the system after the delay $\tau$, and $y(t) = (y_1(t), \dotsc, y_m(t))^\transp \in \mathbb R^m$ is an $m$-dimensional output measured at time $t$. This system reduces to \eqref{eq:DDAE-MID} when one performs a linear output feedback of the form
\begin{equation}
\label{eq:control-DDAE}
u(t) = \sum_{j=1}^m k_j y_j(t).
\end{equation}

Note that \eqref{eq:DDAE-MID} can be rewritten as
\begin{equation}
\label{eq:DDAE-MID-matrix-1}
\left\{
\begin{aligned}
X^\prime(t) & = A X(t) + B K y(t - \tau), \\
y(t) & = C X(t) + D K y(t - \tau),
\end{aligned}
\right.
\end{equation}
where
\begin{gather}
X(t) = \begin{pmatrix}
x(t) \\
x^\prime(t) \\
\vdots \\
x^{(n-1)}(t) \\
\end{pmatrix},\qquad A = \begin{pmatrix}
0 & 1 & \cdots & 0 \\
\vdots & \ddots & \ddots & \vdots \\
0 & 0 & \cdots & 1 \\
-a_0 & -a_1 & \cdots & -a_{n-1} \\
\end{pmatrix}, \qquad B = \begin{pmatrix}
0 \\
\vdots \\
0 \\
b \\
\end{pmatrix}, \label{eq:matrices-DDAE-1} \displaybreak[0] \\
C = \begin{pmatrix}
c_{1, 0} & c_{1, 1} & \cdots & c_{1, n-1} \\
c_{2, 0} & c_{2, 1} & \cdots & c_{2, n-1} \\
\vdots & \vdots & \ddots & \vdots \\
c_{m, 0} & c_{m, 1} & \cdots & c_{m, n-1} \\
\end{pmatrix}, \qquad D = \begin{pmatrix}
d_1 \\
d_2 \\
\vdots \\
d_m
\end{pmatrix}, \qquad K = \begin{pmatrix}k_1 & k_2 & \cdots & k_m\end{pmatrix}. \label{eq:matrices-DDAE-1bis}
\end{gather}
System~\eqref{eq:DDAE-MID} can also be put under the form \eqref{eq:linear-DDAE}: setting
\begin{equation}
\label{eq:matrices-DDAE-2}
Z(t) = \begin{pmatrix}X(t) \\ y(t) \end{pmatrix} \in \mathbb R^{n + m}, \qquad E = \begin{pmatrix}\Id_n & 0 \\ 0 & 0 \\\end{pmatrix}, \qquad A_0 = \begin{pmatrix}A & 0 \\ C & -\Id_m\\\end{pmatrix}, \qquad A_1 = \begin{pmatrix}0 & B K \\ 0 & D K \\\end{pmatrix},
\end{equation}
we have that \eqref{eq:DDAE-MID-matrix-1} can be rewritten as
\[
E Z^\prime(t) = A_0 Z(t) + A_1 Z(t - \tau),
\]
which is under the form \eqref{eq:linear-DDAE} with $M = 1$. Defining $U \in \mathcal M_{n + m, m}(\mathbb R)$ by
\[
U = \begin{pmatrix}
0 \\
-\Id_m
\end{pmatrix},
\]
we remark that $U^\transp E = 0$, $E U = 0$, and the matrix $U^\transp A_0 U = - \Id_m$ is invertible. Hence, thanks to the results recalled in Section~\ref{sec:stability-LTI-TDS}, the exponential stability of \eqref{eq:DDAE-MID} can be studied through the roots of its characteristic function
\begin{equation}
\label{eq:Delta-DDAE}
\Delta(s) = \det\left(s E - A_0 - A_1 e^{-s \tau}\right).
\end{equation}
It turns out that $\Delta$ has the form \eqref{eq:Delta-DDE}, as shown in the next result.

\begin{proposition}
\label{prop:Delta-DDAE-simplified}
Let $\tau > 0$, $E$, $A_0$, and $A_1$ be defined as in \eqref{eq:matrices-DDAE-1}--\eqref{eq:matrices-DDAE-2}, and $\Delta\colon \mathbb C \to \mathbb C$ be the function defined in \eqref{eq:Delta-DDAE}. Then
\begin{equation}
\label{eq:Delta-DDAE-simplified}
\Delta(s) = P(s) - e^{-s\tau} \left[P(s) K D + b \sum_{j=0}^{n-1} s^j K C e_{j+1}\right],
\end{equation}
where $P(s) = \det(s \Id_n - A) = s^n + \sum_{j=1}^{n-1} a_j s^j$ and $e_1, \dotsc, e_n$ denotes the canonical basis of $\mathbb R^n$.
\end{proposition}

\begin{proof}
We only show \eqref{eq:Delta-DDAE-simplified} for $s \in \mathbb C$ such that $P(s) \neq 0$. Indeed, as both left- and right-hand sides of \eqref{eq:Delta-DDAE-simplified} define holomorphic functions on $\mathbb C$, it suffices to prove \eqref{eq:Delta-DDAE-simplified} for all $s \in \mathbb C$ except for those $s$ in the finite set $\{s \in \mathbb C \suchthat P(s) = 0\}$. We then fix in the sequel $s \in \mathbb C$ such that $P(s) \neq 0$.

Note that $\Delta$ can be written, in terms of the matrices $A$, $B$, $C$, $D$, and $K$ from \eqref{eq:matrices-DDAE-1}, as
\[
\Delta(s) = \det\begin{pmatrix}s \Id_n - A & - B K e^{-s\tau} \\ -C & \Id_m - D K e^{-s\tau}\end{pmatrix}.
\]
We now express the determinant in the right-hand side of the above expression using the Schur complement of the corresponding matrix. More precisely, note that
\begin{multline*}
\begin{pmatrix}s \Id_n - A & - B K e^{-s\tau} \\ -C & \Id_m - D K e^{-s\tau}\end{pmatrix} \\ = \begin{pmatrix}s \Id_n - A & 0 \\ -C & \Id_m - \left(D + C (s \Id_n - A)^{-1} B\right) K e^{-s\tau}\end{pmatrix} \begin{pmatrix}\Id_n & -(s \Id_n - A)^{-1} B K e^{-s\tau} \\ 0 & \Id_m \end{pmatrix},
\end{multline*}
and thus
\[
\Delta(s) = P(s) \det\left(\Id_m - \left(D + C (s \Id_n - A)^{-1} B\right) K e^{-s\tau}\right).
\]
Since $\left(D + C (s \Id_n - A)^{-1} B\right) K$ is a rank-one matrix, we have, from Cauchy's formula for the determinant of a rank-one perturbation (see, e.g., \cite[(0.8.5.11)]{Horn2013Matrix}), that
\[
\det\left(\Id_m - \left(D + C (s \Id_n - A)^{-1} B\right) K e^{-s\tau}\right) = 1 - e^{-s \tau} K \left(D + C (s \Id_n - A)^{-1} B\right).
\]
Hence
\[
\Delta(s) = P(s) - e^{-s\tau} P(s) K D - e^{-s\tau} K C \Adj (s\Id_n - A) B,
\]
where $\Adj (s\Id_n - A)$ denotes the adjugate of the matrix $s\Id_n - A$. A straightforward computation shows that $\Adj (s\Id_n - A) B = b \sum_{j=0}^{n-1} s^j e_{j+1}$, where $e_1, \dotsc, e_n$ denotes the canonical basis of $\mathbb R^n$. This finally gives \eqref{eq:Delta-DDAE-simplified}.
\end{proof}

As an immediate consequence of Theorem~\ref{thm:main-MID} and Proposition~\ref{prop:Delta-DDAE-simplified}, we obtain the following result on the exponential stability of \eqref{eq:DDAE-MID}.

\begin{corollary}
\label{coro:DDAE}
Let $n$ be a nonnegative integer, $\tau > 0$, and $s_0 \in \mathbb R$. Let $a_0, \dotsc,\allowbreak a_{n-1},\allowbreak b,\allowbreak k_1, \dotsc,\allowbreak k_m,\allowbreak c_{1,0}, \dotsc,\allowbreak c_{m, n-1},\allowbreak d_1, \dotsc,\allowbreak d_m$ be real coefficients satisfying
\[
\left\{
\begin{aligned}
a_k & = (-1)^{n-k} \frac{1}{k!} \sum_{j=k}^{n} \frac{(2n-j)! s_0^{j-k}}{(n-j)! (j-k)! {\tau}^{n-j}}, & \qquad & \text{ for } k \in \llbracket 0, n-1\rrbracket, \\
K D & = (-1)^n e^{s_0 \tau}, & & \\
K (a_k D + b C e_{k+1}) & = (-1)^n \frac{e^{s_0 \tau}}{k!} \sum_{j=k}^{n} \frac{(2n - j)! (-s_0)^{j-k}}{(n - j)! (j - k)! \tau^{n-j}}, & & \text{ for } k \in \llbracket 0, n-1\rrbracket,
\end{aligned}
\right.
\]
where the matrices $C$, $D$, and $K$ are defined as in \eqref{eq:matrices-DDAE-1bis} and $e_1, \dotsc, e_n$ denotes the canonical basis of $\mathbb R^n$. Then the trivial solution of \eqref{eq:DDAE-MID} is exponentially stable if and only if $s_0 < 0$ or, equivalently, $a_{n-1} > -\frac{n (m+1)}{\tau}$.
\end{corollary}

\section{Some insights on MID for roots of intermediate multiplicities}
\label{sec:MID-intermediate}

Theorem~\ref{thm:main-MID} implies that, given $n, m \in \mathbb N$, $\tau \in \mathbb R_+^\ast$, and $s_0 \in \mathbb R$, there exists a unique choice of coefficients $a_0, \dotsc, a_{n-1}, \alpha_0, \dotsc, \alpha_m$ of the quasipolynomial $\Delta$ from \eqref{eq:Delta-DDE} such that $s_0$ is a root of $\Delta$ of multiplicity $D = n + m + 1$, this choice being the one given in \eqref{Coeffs}. Even though Theorem~\ref{thm:main-MID} can be used in many situations to stabilize systems with delays, as we will describe later in Section~\ref{sec:MID-applis}, there are many applications in which the values of some of the coefficients $a_0, \dotsc, a_{n-1}, \alpha_0, \dotsc, \alpha_m$ is fixed from the modeling of the system and cannot be modified, and in these cases it may be impossible to choose the free parameters of the system in such a way that all equations in \eqref{Coeffs} are satisfied.

A natural idea to overcome this difficulty is to relax the constraint of \emph{maximal} multiplicity of $s_0$, and place instead a root $s_0$ with ``large enough'' multiplicity. This brings the natural question of whether the multiplicity-induced dominance property holds for roots of ``large enough'' multiplicity. In other words, given a quasipolynomial $\Delta$ of the form \eqref{eq:Delta-DDE}, if $\Delta$ admits a root $s_0$ of multiplicity $\Lambda$ ``large enough'', but strictly smaller than the maximal multiplicity $D = n + m + 1$, is $s_0$ necessarily dominant? A first remark is that such a property cannot hold in full generality, as illustrated by the following examples.

\begin{example}
\label{expl:not-dominant}
Consider the quasipolynomial $\Delta$ from \eqref{eq:Delta-DDE} with $n = 3$, $m = 2$, $\tau = 1$, and
\begin{align*}
a_0 & = \frac{-2 e^4 + 8 e^3 - 10 e^2 + 8 e}{e^4 - 6 e^3 + 12 e^2 - 8 e - 1}, & \alpha_0 & = \frac{2 e^4 - 8 e^3 + 10 e^2 - 8 e}{e^4 - 6 e^3 + 12 e^2 - 8 e - 1}, \\
a_1 & = \frac{5 e^4 - 22 e^3 + 31 e^2 - 22 e}{e^4 - 6 e^3 + 12 e^2 - 8 e - 1}, & \alpha_1 & = \frac{- 3 e^4 + 14 e^3 - 21 e^2 + 14 e}{e^4 - 6 e^3 + 12 e^2 - 8 e - 1}, \\
a_2 & = \frac{-4 e^4 + 21 e^3 - 36 e^2 + 23 e}{e^4 - 6 e^3 + 12 e^2 - 8 e - 1}, & \alpha_2 & = \frac{- 3 e^3 + 10 e^2 - 5 e}{e^4 - 6 e^3 + 12 e^2 - 8 e - 1}.
\end{align*}
This quasipolynomial has degree $D = 6$. One can verify, by straightforward computations, that $\Delta$ has a root of multiplicity $3$ at $s = 0$, a root of multiplicity $2$ at $s = 1$, and a root of multiplicity $1$ at $s = 2$. In particular, the root at $s = 0$ has multiplicity $3$ but it is not dominant, and, similarly, the root at $s = 1$ has multiplicity $2$ but it is not dominant either.

Using Proposition~\ref{prop:Polya-Szego}, one deduces that the roots at $0$, $1$ and $2$ are the only roots of $\Delta$ in the horizontal strip $\{s \in \mathbb C \suchthat \abs{\Imag s} < 2\pi\}$ and, in particular, they are the only real roots of $\Delta$. A numerical representation\footnote{Roots were computed numerically using Python's \texttt{cxroots} module \cite{Parini2018cxroots}, and the search for roots was restricted to the rectangle $\{s \in \mathbb C \suchthat \abs{\Real s} \leq 10,\, \abs{\Imag s} \leq 100\}$.} of the roots of $\Delta$ is provided in Figure~\ref{fig:expl-not-dominant}, in which simple roots are represented by blue circles, the double root at $s = 1$ is represented by an orange square, and the triple root at $s = 0$ is represented by a green ``x''.

\begin{figure}[ht]
\centering
\resizebox{0.5\textwidth}{!}{\input{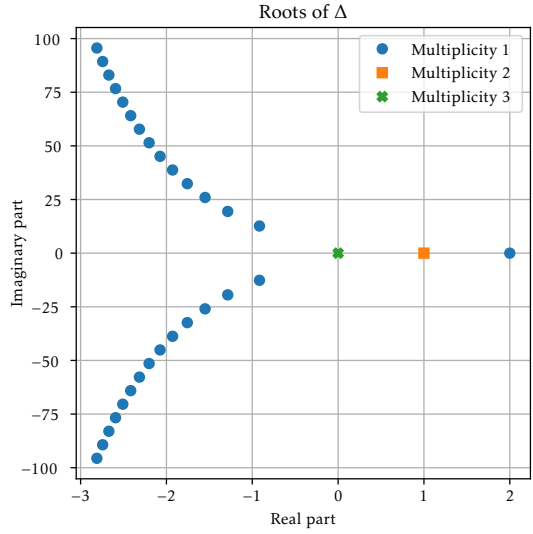}}
\caption{Roots of the quasipolynomial $\Delta$ from Example~\ref{expl:not-dominant}}
\label{fig:expl-not-dominant}
\end{figure}
\end{example}

One might argue that, in the previous example, the reason for the root of multiplicity $3$ at $s = 0$ not being dominant is that its multiplicity is not ``large enough'' with respect to the maximal multiplicity $D = 6$. However, even roots of multiplicity $D-1$ may fail to be dominant, as illustrated in the following example.

\begin{example}
\label{expl:not-dominant-parametric}
Let $\eta \in \mathbb R^\ast$ and consider the quasipolynomial $\Delta$ from \eqref{eq:Delta-DDE} with $n = 2$, $m = 1$, $\tau = 1$, and
\begin{equation}
\label{eq:coeffs-eta}
\begin{aligned}
a_0 & = \frac{2 \eta \left(\eta - 1 + e^{-\eta}\right)}{\eta - 2 + e^{-\eta} (\eta + 2)}, & \qquad \alpha_0 & = -\frac{2 \eta \left(\eta - 1 + e^{-\eta}\right)}{\eta - 2 + e^{-\eta} (\eta + 2)}, \\
a_1 & = \frac{- \eta^{2} + 2 - 2 (\eta + 1) e^{-\eta}}{\eta - 2 + e^{-\eta} (\eta + 2)}, & \alpha_1 & = \frac{- \eta^{2} + 2 \eta - 2 + 2 e^{-\eta}}{\eta - 2 + e^{-\eta} (\eta + 2)}.
\end{aligned}
\end{equation}
This quasipolynomial has degree $D = 4$ and one can verify that $\Delta$ has a root of multiplicity $3$ at $s = 0$ and a root of multiplicity $1$ at $s = \eta$. When $\eta > 0$, the root at $s = 0$ has multiplicity $3$, but it is not dominant. As in Example~\ref{expl:not-dominant}, one can verify, using Proposition~\ref{prop:Polya-Szego}, that the roots at $0$ and $\eta$ are the only roots of $\Delta$ in the horizontal strip $\{s \in \mathbb C \suchthat \abs{\Imag s} < 2 \pi\}$.

Note that, as $\eta \to 0$, the coefficients $a_0, a_1, \alpha_0, \alpha_1$ converge, respectively, to
\[a_0^\ast = 6, \qquad a_1^\ast = -4, \qquad \alpha_0^\ast = -6, \qquad \alpha_1^\ast = -2,\]
and hence the quasipolynomial $\Delta$ converges to the quasipolynomial
\begin{equation}
\label{eq:Delta-limit-eta-0}
\Delta_\ast(s) = s^2 - 4 s + 6 - e^{-s}(2 s + 6),
\end{equation}
uniformly on every compact subset of $\mathbb C$. The quasipolynomial $\Delta_\ast$ has degree $4$ and one immediately verifies that it admits a root of multiplicity $4$ at $s = 0$. Then, by Theorem~\ref{thm:main-MID}, the root $0$ of $\Delta_\ast$ is strictly dominant. In addition, as before, $0$ is the unique root of $\Delta_\ast$ in the horizontal strip $\{s \in \mathbb C \suchthat \abs{\Imag s} < 2 \pi\}$.

\begin{figure}[ht]
\centering
\resizebox{0.5\textwidth}{!}{\import{Figures}{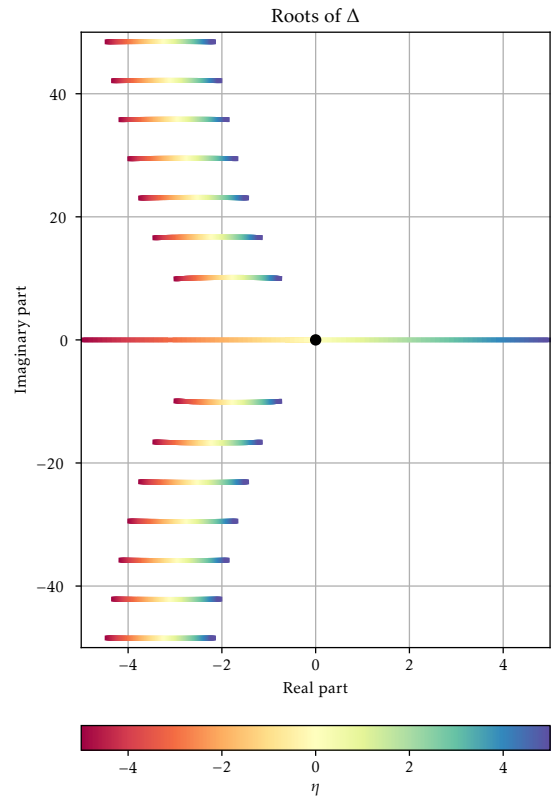}}
\caption{Roots of the quasipolynomial $\Delta$ from Example~\ref{expl:not-dominant-parametric} in terms of $\eta$}
\label{fig:expl-not-dominant-free-root}
\end{figure}

Figure~\ref{fig:expl-not-dominant-free-root} represents the roots\footnote{Roots were computed numerically using Python's \texttt{cxroots} module \cite{Parini2018cxroots}, and the search for roots was restricted to the rectangle $\{s \in \mathbb C \suchthat \abs{\Real s} \leq 10,\, \abs{\Imag s} \leq 50\}$.} of $\Delta$ in the rectangle $\{s \in \mathbb C \suchthat \abs{\Real s} \leq 5,\, \abs{\Imag s} \leq 50\}$ as $\eta$ varies in the interval $[-5, 5]$, with shades between red and yellow representing negative values of $\eta$ and shades between yellow and blue representing positive values. For $\eta = 0$, the roots represented in Figure~\ref{fig:expl-not-dominant-free-root} are those of the quasipolynomial $\Delta_\ast$ from \eqref{eq:Delta-limit-eta-0}. The black circle represents the root at $0$, which is of multiplicity $3$ for $\eta \neq 0$ and of multiplicity $4$ for $\eta = 0$. One observes that the root at $0$ is not dominant when $\eta > 0$.

Figure~\ref{fig:expl-not-dominant-free-root} suggests that $0$ may be dominant whenever $\eta \in (-\infty, 0]$, since all other roots represented in Figure~\ref{fig:expl-not-dominant-free-root} have negative real part in that case. This is indeed true, as one can verify by applying \cite[Theorem~4.2]{Boussaada2020Multiplicity}: using \eqref{eq:coeffs-eta}, one can compute that the quantity $s_+$ from \cite[(4.11)]{Boussaada2020Multiplicity} is equal to zero and deduce that \cite[($\star_+$)]{Boussaada2020Multiplicity} is satisfied. The dominance of the triple root at $0$ comes thus as a consequence of \cite[Theorem~4.2~iii)]{Boussaada2020Multiplicity}.
\end{example}

At the light of Examples~\ref{expl:not-dominant} and \ref{expl:not-dominant-parametric}, one can reformulate the question of the MID property for roots with intermediate multiplicity as: given a quasipolynomial $\Delta$ under the form \eqref{eq:Delta-DDE} admitting a root $s_0$ of ``large enough'' multiplicity, under which additional conditions is $s_0$ a dominant root of $\Delta$? As stated in Section~\ref{sec:intro-MID}, this question has been partially addressed in some previous works, most often in particular cases such as when $m = n - 1$ and for roots of multiplicity $n + 1$, providing only sufficient conditions for dominance \cite{Balogh2020Towards, Balogh2022Conditions, Sassi2022Observer, RojasRicca2022Multiplicity}, with a more detailed analysis of the case $n = 2$, $m = 1$ for a root of multiplicity $3$ provided in \cite[Theorem~4.2]{Boussaada2020Multiplicity}. A general answer to this question is still lacking in the literature, with existing works dealing only with particular cases or with sufficient conditions that are not expected to be sharp. In the sequel, we shall report some of the main results of the recent works \cite{Boussaada2025OverOrder, Boussaada2025Decay, Boussaada2022MID} on this topic.

So far, we have been quite vague in what we mean by a ``large enough'' multiplicity. Following \cite{Boussaada2025OverOrder}, we shall consider in the sequel roots of \emph{overorder} multiplicity, i.e., roots whose multiplicity is strictly larger than the integer $n$ from \eqref{eq:DDE} (which is the largest order of derivation in \eqref{eq:DDE} under the standard assumption that $m \leq n$, whence the name \emph{overorder}). We will first provide, in Section~\ref{sec:iMID-overorder}, a characterization of quasipolynomials with a real root of overorder multiplicity similar to the equivalence between \ref{item:max-mult} and \ref{item:DeltaTilde} in Theorem~\ref{thm:main-MID}. Further properties of quasipolynomials with roots of overorder multiplicity will then be considered in Sections~\ref{sec:control-oriented}. While these properties do not constitute a proof of dominance on themselves, they are useful for proving dominance in some examples, as we will illustrate later in Section~\ref{sec:expl-integrator}.

\subsection{Characterization of roots of overorder multiplicity}
\label{sec:iMID-overorder}

The goal of this section is to generalize the equivalence between \ref{item:max-mult} and \ref{item:DeltaTilde} in Theorem~\ref{thm:main-MID} to the case of roots of overorder multiplicity. The first main result is the following, which shows that, if a quasipolynomial of the form \eqref{eq:Delta-DDE} admits a root $s_0$ of overorder multiplicity $\Lambda$, then it can be factorized as the product of $(s - s_0)^\Lambda$ and a holomorphic function which can be expressed as an integral.

\begin{theorem}
\label{thm:max-mult}
Let $n$, $m$ be nonnegative integers with $n \geq m$, $\tau > 0$, $a_0, \dotsc, a_{n-1}, \alpha_0, \dotsc, \alpha_m$ be real numbers, and consider the quasipolynomial $\Delta$ from \eqref{eq:Delta-DDE}. Let $D = n + m + 1$ denote the degree of $\Delta$ and take $s_0 \in \mathbb R$ and $\Lambda \in \llbracket n + 1, D\rrbracket$. Then the real number $s_0$ is a root of multiplicity at least $\Lambda$ of $\Delta$ if and only if there exists a polynomial $p$ of degree at most $D - \Lambda$ with $p(0) = 1$ such that
\begin{equation}
\label{eq:Delta-max-mult}
\Delta(s) = \frac{\tau^{\Lambda - n} (s - s_0)^{\Lambda}}{(\Lambda - n - 1)!}  \int_0^1 t^{\Lambda - n - 1} (1 - t)^{\Lambda - m - 1} p(t) e^{-t \tau(s - s_0)} \diff t, \qquad \text{ for all } s \in \mathbb C.
\end{equation}
\end{theorem}

\begin{proof}
Let $\mathcal H$ denote the set of all holomorphic functions $\Delta\colon \mathbb C \to \mathbb C$, seen here as a real vector space. For fixed $n, m \in \mathbb N$ with $n \geq m$, $\tau > 0$, $s_0 \in \mathbb R$, and $\Lambda \in \llbracket n+1, D\rrbracket$, we introduce the sets
\begin{align*}
\mathcal V_\tau & = \left\{\Delta \in \mathcal H \suchthat \Delta \text{ is of the form \eqref{eq:Delta-DDE} for some real numbers } a_0, \dotsc, a_{n-1}, \alpha_0, \dotsc, \alpha_{m}\right\}, \\
\mathcal V_{s_0, \tau}^\Lambda & = \left\{\Delta \in \mathcal V_\tau \suchthat \Delta^{(k)}(s_0) = 0 \text{ for all }  k \in \llbracket 0, \Lambda - 1\rrbracket\right\}, \\
\mathcal W_{s_0, \tau}^\Lambda & = \{\Delta \in \mathcal H \suchthat \Delta \text{ is of the form \eqref{eq:Delta-max-mult} for some polynomial } p\colon \mathbb R \to \mathbb R \text{ of degree at most}\\
& \hphantom{{} = \{\Delta \in \mathcal H \suchthat {}} D - \Lambda \text{ with } p(0) = 1\}.
\end{align*}
Note that $\mathcal V_\tau$, $\mathcal V_{s_0, \tau}^\Lambda$, and $\mathcal W_{s_0, \tau}^\Lambda$ are all affine subspaces of $\mathcal H$. These spaces also depend on $n$ and $m$, but we omit these dependencies from their notation. We will prove the result by showing that $\mathcal V_{s_0, \tau}^\Lambda = \mathcal W_{s_0, \tau}^\Lambda$.

\step{1}{Dimension of $\mathcal V_\tau$}

Note that $\mathcal V_\tau = \{s \mapsto s^n\} + \widehat{\mathcal V}_\tau$, where $\widehat{\mathcal V}_\tau$ is the vector subspace of $\mathcal H$ spanned by the family of $D$ functions $\mathcal F_\tau$ defined by $\mathcal F_\tau = \{s \mapsto s^k \suchthat k \in \llbracket 0, n-1\rrbracket\} \cup \{s \mapsto e^{-s \tau} s^k \suchthat k \in \llbracket 0, m\rrbracket\}$. We claim that the family of functions $\mathcal F_\tau$ is linearly independent. Indeed, let $(a_0, \dotsc, a_{n-1}, \alpha_0, \dotsc, \alpha_m) \in \mathbb R^D$ be such that
\[
\sum_{k=0}^{n-1} a_k s^k + e^{-s \tau} \sum_{k=0}^m \alpha_k s^k = 0, \qquad \text{ for every } s \in \mathbb C.
\]
Let $P_0$ and $P_1$ be the polynomials defined for $s \in \mathbb C$ by
\[
P_0(s) = \sum_{k=0}^{n-1} a_k s^k, \qquad P_1(s) = \sum_{k=0}^m \alpha_k s^k.
\]
Hence
\begin{equation}
\label{eq:identically-zero}
P_0(s) + e^{-s \tau} P_1(s) = 0 \qquad \text{ for every } s \in \mathbb C.
\end{equation}
By taking derivatives of this identity and using the fact that $s \mapsto e^{-s \tau}$ never vanishes, one deduces that a complex number $s_\ast$ is a root of multiplicity $K \in \mathbb N$ of $P_0$ if and only if $s_\ast$ is a root of the same multiplicity $K$ of $P_1$. In particular, $P_0$ and $P_1$ are the same polynomial up to a multiplicative constant, i.e., there exists a constant $\lambda \in \mathbb C$ such that $P_0(s) = \lambda P_1(s)$ for every $s \in \mathbb C$. Hence, using \eqref{eq:identically-zero}, we deduce that
\[
(\lambda + e^{-s \tau}) P_1(s) = 0 \qquad \text{ for every } s \in \mathbb C.
\]
The function $s \mapsto \lambda + e^{-s \tau}$ has at most countably many zeros, and hence $P_1$ must vanish outside of these zeros. Since $P_1$ is a polynomial, this implies that $P_1$ is identically zero, and thus all its coefficients are zero. Since $P_0(s) = \lambda P_1(s)$ for every $s \in \mathbb C$, all the coefficients of $P_0$ are also zero. Thus, $a_0 = \dotsb = a_{n-1} = \alpha_0 = \dotsb = \alpha_m = 0$, yielding that $\mathcal F_\tau$ is linearly independent. Hence $\mathcal F_\tau$ is a basis of $\widehat{\mathcal V}_\tau$, showing that $\dim \mathcal V_\tau = \dim \widehat{\mathcal V}_\tau = D$.

\step{2}{Dimension of $\mathcal V_{s_0, \tau}^\Lambda$}

To compute the dimension of $\mathcal V_{s_0, \tau}^\Lambda$, let us consider the linear map $\mathcal N_{s_0, \tau}\colon \mathcal H \to \mathcal H$ defined by $(\mathcal N_{s_0, \tau} \Delta)(z) = \tau^n \Delta(s_0 + \frac{z}{\tau})$ for every $\Delta \in \mathcal H$ and $z \in \mathbb C$. One can easily verify that $\mathcal N_{s_0, \tau}$ is an isomorphism, with inverse given by $(\mathcal N_{s_0, \tau}^{-1} \Delta)(s) = \frac{1}{\tau^n} \Delta(\tau(s - s_0))$ for $\Delta \in \mathcal H$ and $s \in \mathbb C$. It is also immediate to verify that $\mathcal N_{s_0, \tau}$ is a bijection between $\mathcal V_{s_0, \tau}^\Lambda$ and $\mathcal V_{0, 1}^\Lambda$, and thus $\dim \mathcal V_{s_0, \tau}^\Lambda = \dim \mathcal V_{0, 1}^\Lambda$.

Fix $\Delta_0 \in \mathcal V_{0, 1}^\Lambda$ and note that $\mathcal V_{0, 1}^\Lambda = \Delta_0 + \widehat{\mathcal V}_{0, 1}^\Lambda$, where
\[
\widehat{\mathcal V}_{0, 1}^\Lambda = \left\{\Delta \in \widehat{\mathcal V}_1 \suchthat \Delta^{(k)}(0) = 0 \text{ for all } k \in \llbracket 0, \Lambda-1\rrbracket\right\}.
\]
Introducing the linear map $\Xi^\Lambda\colon \widehat{\mathcal V}_1 \to \mathbb R^\Lambda$ given by $\Xi^\Lambda \Delta = \left(\Delta(0), \Delta^\prime(0), \dotsc, \Delta^{(\Lambda-1)}(0)\right)$, one observes at once that $\widehat{\mathcal V}_{0, 1}^\Lambda = \Ker \Xi^\Lambda$. Considering the basis of $\widehat{\mathcal V}_1$ given by the functions in $\{s \mapsto \frac{s^k}{k!} \suchthat k \in \llbracket 0, n-1\rrbracket\} \cup \{s \mapsto e^{-s} \frac{s^k}{k!} \suchthat k \in \llbracket 0, m\rrbracket\}$ and the canonical basis of $\mathbb R^\Lambda$, the linear map $\Xi^\Lambda$ can be represented by the matrix $[\Xi^\Lambda] \in \mathcal M_{\Lambda, D}(\mathbb R)$ given by
\begin{equation}
\label{eq:Xi-Lambda-matrix}
[\Xi^\Lambda] = \begin{pmatrix}
\Id_n & A \\
0_{\Lambda - n, n} & B
\end{pmatrix}
\end{equation}
where the matrices\footnote{Note that we use here the convention that indices of matrices start at $0$, instead of $1$. This is done to simplify the expressions of $A_{j, k}$ and $B_{j, k}$.} $A = (A_{j, k})_{\substack{j \in \llbracket 0, n-1\rrbracket \\ k \in \llbracket 0, m\rrbracket}} \in \mathcal M_{n, m+1}(\mathbb R)$ and $B = (B_{j, k})_{\substack{j \in \llbracket 0, \Lambda-n-1\rrbracket \\ k \in \llbracket 0, m\rrbracket}} \in \mathcal M_{\Lambda-n, m+1}(\mathbb R)$ are defined by
\[
\begin{aligned}
A_{j, k} & = (-1)^{j-k}\binom{j}{k}, & \qquad j & \in \llbracket 0, n-1\rrbracket, & k \in \llbracket 0, m\rrbracket, \\
B_{j, k} & = (-1)^{n+j-k}\binom{n+j}{k}, & j & \in \llbracket 0, \Lambda-n-1\rrbracket, & k \in \llbracket 0, m\rrbracket,
\end{aligned}
\]
with the convention that $\binom{j}{k} = 0$ if $k > j$ (recall the notation in pages~\pageref{notation}--\pageref{end:notation}). That is,
\begin{align*}
A & = \begin{pmatrix}
\binom{0}{0} & 0 & 0 & \cdots & 0 \\
-\binom{1}{0} & \binom{1}{1} & 0 & \cdots & 0 \\
\binom{2}{0} & -\binom{2}{1} & \binom{2}{2} & \ddots & \vdots \\
\vdots & \vdots & \vdots & \ddots & 0 \\
(-1)^{n-1} \binom{n-1}{0} & (-1)^{n-2} \binom{n-1}{1} & (-1)^{n-3} \binom{n-1}{2} & \cdots & (-1)^{n-1-m} \binom{n-1}{m}
\end{pmatrix}, \displaybreak[0] \\
B & = \begin{pmatrix}
(-1)^{n} \binom{n}{0} & (-1)^{n - 1} \binom{n}{1} & \cdots & (-1)^{n - m} \binom{n}{m} \\
(-1)^{n+1} \binom{n+1}{0} & (-1)^{n} \binom{n+1}{1} & \cdots & (-1)^{n + 1 - m} \binom{n + 1}{m} \\
\vdots & \vdots & \ddots & \vdots \\
(-1)^{\Lambda - 1} \binom{\Lambda - 1}{0} & (-1)^{\Lambda - 2} \binom{\Lambda - 1}{1} & \cdots & (-1)^{\Lambda - 1 - m} \binom{\Lambda - 1}{m} \\
\end{pmatrix}.
\end{align*}

Let $C = (C_{j, k})_{\substack{j \in \llbracket 0, m\rrbracket \\ k \in \llbracket 0, \Lambda - n - 1}} \in \mathcal M_{m+1, \Lambda - n}(\mathbb R)$ be the matrix defined by
\[
C_{j, k} = (-1)^{n + k} \binom{n + m - j + k}{m - j}, \qquad j \in \llbracket 0, m\rrbracket, \, k \in \llbracket 0, \Lambda - n - 1\rrbracket,
\]
i.e.,
\[
C = \begin{pmatrix}
(-1)^{n} \binom{n + m}{m} & (-1)^{n + 1} \binom{n + m + 1}{m} & \cdots & (-1)^{\Lambda - 1} \binom{\Lambda + m - 1}{m} \\
(-1)^{n} \binom{n + m - 1}{m - 1} & (-1)^{n + 1} \binom{n + m}{m - 1} & \cdots & (-1)^{\Lambda - 1} \binom{\Lambda + m - 2}{m - 1} \\
\vdots & \vdots & \ddots & \vdots \\
(-1)^{n} \binom{n}{0} & (-1)^{n + 1} \binom{n + 1}{0} & \cdots & (-1)^{\Lambda - 1} \binom{\Lambda - 1}{0} \\
\end{pmatrix}
\]
Let us compute the product $B C \in \mathcal M_{\Lambda-n, \Lambda-n}(\mathbb R)$. For $(j, k) \in \llbracket 0, \Lambda - n - 1\rrbracket^2$, the $(j, k)$th element of $B C$, denoted by $(B C)_{j, k}$, is given by
\[
(B C)_{j, k} = (-1)^{j + k} \sum_{\ell = 0}^m (-1)^{\ell}\binom{n+j}{\ell} \binom{n + m - \ell + k}{m - \ell} = (-1)^{j + k} \sum_{\ell = 0}^m (-1)^{\ell} \binom{n+j}{\ell} \binom{n + k + m - \ell}{n + k}.
\]
Applying Proposition~\ref{PropHorribleSum} from Appendix~\ref{app:binomial} to the above sum, we deduce that
\[
(B C)_{j, k} = (-1)^{j + k} \binom{m + k - j}{m},
\]
that is,
\[
B C = \begin{pmatrix}
\binom{m}{m} & -\binom{m + 1}{m} & \binom{m + 2}{m} & \cdots & (-1)^{\Lambda - n - 1} \binom{m + \Lambda - n - 1}{m} \\
0 & \binom{m}{m} & -\binom{m + 1}{m} & \ddots & \vdots \\
0 & 0 & \binom{m}{m} & \ddots & \binom{m + 2}{m} \\
\vdots & \vdots & \ddots & \ddots & -\binom{m + 1}{m} \\
0 & 0 & \cdots & 0 & \binom{m}{m} \\
\end{pmatrix}.
\]
Hence $\rank B C = \Lambda - n$, showing that $\rank B = \Lambda - n$. Thus, by \eqref{eq:Xi-Lambda-matrix}, we have that $\rank \Xi^\Lambda = \Lambda$ and, by the rank--nullity theorem, we deduce that $\dim \mathcal V_{s_0, \tau}^\Lambda = \dim \widehat{\mathcal V}_{0, 1}^\Lambda = \dim \Ker \Xi^\Lambda = D - \Lambda$.

\step{3}{Dimension of $\mathcal W_{s_0, \tau}^\Lambda$}

As in Step~2, one can easily verify that the isomorphism $\mathcal N_{s_0, \tau}$ is a bijection between $\mathcal W_{s_0, \tau}^\Lambda$ and $\mathcal W_{0, 1}^\Lambda$, hence $\dim \mathcal W_{s_0, \tau}^\Lambda = \dim \mathcal W_{0, 1}^\Lambda$.

For $k \in \llbracket 0, D - \Lambda\rrbracket$, let $p_k$ be the polynomial defined for $t \in [0, 1]$ by $p_k(t) = t^k$ and $\Delta_k \in \mathcal H$ be the function defined by
\[
\Delta_k(s) = \frac{s^{\Lambda}}{(\Lambda - n - 1)!} \int_0^1 t^{\Lambda - n - 1} (1 - t)^{\Lambda - m - 1} p_k(t) e^{-t s} \diff t, \qquad \text{ for all } s \in \mathbb C.
\]
It follows from the definition of $\mathcal W_{0, 1}^\Lambda$ that $\mathcal W_{0, 1}^\Lambda = \Delta_0 + \Span \{\Delta_1, \dotsc, \Delta_{D - \Lambda}\}$.

The family $\{\Delta_1, \dotsc, \Delta_{D - \Lambda}\}$ is linearly independent. Indeed, if $\gamma_1, \dotsc, \gamma_{D - \Lambda}$ are such that $\gamma_1 \Delta_1 + \dotsb + \gamma_{D - \Lambda} \Delta_{D - \Lambda} = 0$ in $\mathcal H$, then
\begin{equation}
\label{eq:Laplace-q}
\int_0^1 q(t) e^{-t s} \diff t = 0 \text{ for all } s \in \mathbb C,
\end{equation}
where the polynomial $q\colon \mathbb R \to \mathbb R$ is defined by
\[
q(t) = t^{\Lambda - n - 1} (1 - t)^{\Lambda - m - 1} \sum_{k=1}^{D - \Lambda} \gamma_k t^k \qquad \text{ for } t \in \mathbb R.
\]
Equation~\eqref{eq:Laplace-q} states that the one-sided Laplace transform of the function $t \mapsto q(t) \1_{[0, 1]}(t)$ is identically zero, and thus, by the injectivity of the one-sided Laplace transform (see, e.g., \cite[Chapter~VI, Section~4, Theorem~2]{Yosida1980Functional}), we deduce that $q$ is identically zero. Since $t \mapsto t^{\Lambda - n - 1} (1 - t)^{\Lambda - m - 1}$ is not identically zero, we deduce that the polynomial $t \mapsto \sum_{k=1}^{D - \Lambda} \gamma_k t^k$ is identically zero, implying that $\gamma_1 = \dotsb = \gamma_{D - \Lambda} = 0$. Hence, the family $\{\Delta_1, \dotsc, \Delta_{D - \Lambda}\}$ is linearly independent.

As a conclusion, we obtain that $\dim \mathcal W_{s_0, \tau}^\Lambda = \dim \mathcal W_{0, 1}^\Lambda = D - \Lambda$.

\step{4}{$\mathcal W_{s_0, \tau}^\Lambda \subset \mathcal V_{s_0, \tau}^\Lambda$ and conclusion}

Let us now show that $\mathcal W_{s_0, \tau}^\Lambda \subset \mathcal V_{s_0, \tau}^\Lambda$. Let $\Delta \in \mathcal W_{s_0, \tau}^{\Lambda}$ and $p$ be a polynomial of degree at most $D - \Lambda$ with $p(0) = 1$ be such that $\Delta$ is given by \eqref{eq:Delta-max-mult}. Define $q(t) = t^{\Lambda - n - 1} (1 - t)^{\Lambda - m - 1} p(t)$ and notice that $q(0) = q^\prime(0) = \dotsb = q^{(\Lambda - n - 2)}(0) = 0$, $q^{(\Lambda - n - 1)}(0) = (\Lambda - n - 1)!$, and $q(1) = q^\prime(1) = \dotsb = q^{(\Lambda - m - 2)}(1) = 0$. By Proposition~\ref{PropIntegralPExp} in Appendix~\ref{app:integralPExp}, we deduce that
\begin{align}
\Delta(s) & = \frac{\tau^{\Lambda - n} (s - s_0)^{\Lambda}}{(\Lambda - n - 1)!} \sum_{k = 0}^{\Lambda - 1} \frac{q^{(k)}(0) - q^{(k)}(1) e^{-\tau(s - s_0)}}{\tau^{k+1} (s - s_0)^{k+1}} \notag \\
& = (s - s_0)^{n} + \sum_{k = 0}^{n-1} \frac{q^{(\Lambda - k - 1)}(0)}{\tau^{n - k} (\Lambda - n - 1)!} (s - s_0)^{k} - e^{-\tau(s - s_0)} \sum_{k = 0}^{m} \frac{q^{(\Lambda - k - 1)}(1)}{\tau^{n - k} (\Lambda - n - 1)!} (s - s_0)^{k}, \label{eq:Delta-expanded}
\end{align}
and thus $\Delta \in \mathcal V_\tau$. In addition, it follows immediately from \eqref{eq:Delta-max-mult} that $s_0$ is a root of multiplicity at least $\Lambda$ of $\Delta$, and thus $\mathcal W_{s_0, \tau}^{\Lambda} \subset \mathcal V_{s_0, \tau}^{\Lambda}$.

Finally, since $\mathcal W_{s_0}^{\Lambda}$ and $\mathcal V_{s_0}^{\Lambda}$ are both affine spaces with the same dimension and $\mathcal W_{s_0, \tau}^{\Lambda} \subset \mathcal V_{s_0, \tau}^{\Lambda}$, we conclude that $\mathcal W_{s_0}^{\Lambda} = \mathcal V_{s_0}^{\Lambda}$, yielding the conclusion.
\end{proof}

\begin{remark}
\label{remk:P0-P1-q}
Given a quasipolynomial $\Delta$ under the form \eqref{eq:Delta-DDE} admitting a root $s_0$ of multiplicity at least $\Lambda$, one can obtain a link between the polynomials $P_0\colon s \mapsto s^n + \sum_{k=0}^{n-1} a_k s^k$ and $P_1\colon s \mapsto \sum_{k=0}^m \alpha_k s^k$ appearing in \eqref{eq:Delta-DDE} and the polynomial $p$ appearing in \eqref{eq:Delta-max-mult} through \eqref{eq:Delta-expanded}. Indeed, defining the polynomial $q$ from $p$ as in Step~4 of the above proof, it follows from \eqref{eq:Delta-DDE} and \eqref{eq:Delta-expanded} that
\begin{align*}
P_0(s) & = (s - s_0)^n + \sum_{k = 0}^{n - 1} \frac{q^{(\Lambda - k - 1)}(0)}{\tau^{n - k} (\Lambda - n - 1)!} (s - s_0)^{k}, \\
P_1(s) & = - e^{\tau s_0} \sum_{k = 0}^{m} \frac{q^{(\Lambda - k - 1)}(1)}{\tau^{n - k} (\Lambda - n - 1)!} (s - s_0)^{k}.
\end{align*}
Recalling that $q(t) = t^{\Lambda-n-1} (1 - t)^{\Lambda - m - 1} p(t)$, one may further express $P_0$ and $P_1$ in terms of $p$. Indeed, we have
\begin{align*}
\frac{q^{(\Lambda-k-1)}(0)}{(\Lambda - n - 1)!} & = \sum_{j=0}^{n-k} \binom{\Lambda-k-1}{n-k} \binom{n-k}{j} \frac{(-1)^j  (\Lambda-m-1)!}{(\Lambda-m-j-1)!} p^{(n-k-j)}(0), \\
\frac{q^{(\Lambda-k-1)}(1)}{(\Lambda - n - 1)!} & = \sum_{j=0}^{m-k} \binom{\Lambda-k-1}{m-k} \binom{m-k}{j} \frac{(-1)^{\Lambda-m-1} (\Lambda-m-1)!}{(\Lambda-n-j-1)!} p^{(m-k-j)} (1),
\end{align*}
and thus
\begin{align*}
P_0(s) & = (s - s_0)^{n} + \sum_{k = 0}^{n-1} \sum_{j=0}^{n-k} \frac{(-1)^j}{\tau^{n - k}} \binom{\Lambda-k-1}{n-k} \binom{n-k}{j} \frac{(\Lambda-m-1)!}{(\Lambda-m-j-1)!} p^{(n-k-j)}(0) (s - s_0)^{k}, \\
P_1(s) & = - e^{\tau s_0} \sum_{k = 0}^{m} \sum_{j=0}^{m-k} \frac{ (-1)^{\Lambda-m-1}}{\tau^{n - k}} \binom{\Lambda-k-1}{m-k} \binom{m-k}{j}  \frac{(\Lambda-m-1)!}{(\Lambda-n-j-1)!} p^{(m-k-j)}(1) (s - s_0)^{k}.
\end{align*}
\end{remark}

For a given quasipolynomial under the form \eqref{eq:Delta-DDE}, Theorem~\ref{thm:max-mult} establishes an equivalence between this quasipolynomial admitting a root $s_0$ of multiplicity at least $\Lambda$ and it being expressed under the form \eqref{eq:Delta-max-mult}. Combining with Proposition~\ref{prop:Kummer-integral}, one can express the integral in \eqref{eq:Delta-max-mult} in terms of a sum of Kummer functions, obtaining thus a generalization of the equivalence between \ref{item:max-mult} and \ref{item:DeltaTilde} from Theorem~\ref{thm:main-MID}.

\begin{corollary}
\label{coro:max-mult-Kummer}
Let $n$, $m$ be nonnegative integers with $n \geq m$, $\tau > 0$, $a_0, \dotsc, a_{n-1}, \alpha_0, \dotsc, \alpha_m$ be real numbers, and consider the quasipolynomial $\Delta$ from \eqref{eq:Delta-DDE}. Let $D = n + m + 1$ denote the degree of $\Delta$ and take $s_0 \in \mathbb R$ and $\Lambda \in \llbracket n + 1, D\rrbracket$. Then the real number $s_0$ is a root of multiplicity at least $\Lambda$ of $\Delta$ if and only if there exist real numbers $\sigma_1, \dotsc, \sigma_{D - \Lambda}$ such that
\begin{equation}
\label{eq:Delta-max-mult-Kummer}
\Delta(s) = \tau^{\Lambda - n} (s - s_0)^{\Lambda} \sum_{k=0}^{D - \Lambda} \sigma_k \frac{(k + \Lambda - n - 1)!}{(k + 2\Lambda - D)!} M(k + \Lambda - n, k + 2\Lambda - n - m, -\tau(s - s_0)),
\end{equation}
with $\sigma_0 = 1$.
\end{corollary}

\begin{proof}
By Theorem~\ref{thm:max-mult}, $s_0$ is a root of multiplicity at least $\Lambda$ of $\Delta$ if and only if there exists a polynomial $p$ of degree $D - \Lambda$ with $p(0) = 1$ and such that \eqref{eq:Delta-max-mult} holds. Such a polynomial can be written in a unique manner as
\begin{equation}\label{p1stparam}
p(t) = 1 + \sum_{k = 1}^{D - \Lambda} \sigma_k t^k
\end{equation}
for some real numbers $\sigma_1, \dotsc, \sigma_{D - \Lambda}$. Hence, \eqref{eq:Delta-max-mult} can be rewritten as
\[
\Delta(s) = \frac{\tau^{\Lambda - n} (s - s_0)^{\Lambda}}{(\Lambda - n - 1)!} \sum_{k=0}^{D - \Lambda} \sigma_k \int_0^1 t^{k + \Lambda - n - 1} (1 - t)^{\Lambda - m - 1} e^{-t \tau(s - s_0)} \diff t,
\]
with $\sigma_0 = 1$. The equivalence between this formula and \eqref{eq:Delta-max-mult-Kummer} follows from Proposition~\ref{prop:Kummer-integral}, yielding the conclusion.
\end{proof}

\begin{remark}
In the above proof, instead of expressing $p$ as in \eqref{p1stparam}, one could have chosen a different basis of the space of polynomials, expressing $p$, for instance, as
\[
p(t) = \sum_{k = 0}^{D - \Lambda} \pi_k (1-t)^k
\]
for some real numbers $\pi_0, \dotsc, \pi_{D - \Lambda}$, the condition $p(0) = 1$ being satisfied if and only if $\sum_{k=0}^{D - \Lambda} \pi_k = 1$. Hence, proceeding as in the above proof, one deduces that the real numbers $s_0$ is a root of multiplicity at least $\Lambda$ of $\Delta$ if and only if there exist real numbers $\pi_0, \pi_1, \dotsc, \pi_{D-\Lambda}$ with $\sum_{k=0}^{D-\Lambda} \pi_k = 1$ such that
\begin{equation*}
\Delta(s) = \tau^{\Lambda - n} (s - s_0)^{\Lambda} \sum_{k=0}^{n+m+1-\Lambda} \pi_k \frac{(\Lambda-m+k-1)!}{(2\Lambda-m-n+k-1)!} M(\Lambda-n, 2\Lambda-m-n+k, -\tau(s - s_0)).
\end{equation*}
The coefficients $\sigma_0, \dotsc, \sigma_{D - \Lambda}$ and $\pi_0, \dotsc, \pi_{D - \Lambda}$ are related by the expressions
\[\sigma_k = (-1)^k \sum_{\ell = k}^{D - \Lambda} \binom{\ell}{k} \pi_\ell \quad \text{ and } \quad \pi_k = (-1)^k \sum_{\ell = k}^{D - \Lambda} \binom{\ell}{k} \sigma_\ell \qquad \text{ for } k \in \llbracket 0, D - \Lambda\rrbracket.\]
\end{remark}

\subsection{A control-oriented setting}
\label{sec:control-oriented}

In typical control problems involving time-delay systems of the form \eqref{eq:DDE}, some of the coefficients $a_0, \dotsc, a_{n-1}, \alpha_0, \dotsc, \alpha_m$ are fixed, corresponding to parameters coming from the physical modeling of the system, while other coefficients can be freely chosen, corresponding to parameters coming, for instance, from an implemented linear feedback controller (see \cite{Bhattacharyya1995Robust} for some illustrations of such applications). In some applications, even the delay $\tau$ can be seen as a control parameter that can be chosen, as described in \cite{Suh1979Proportional, Suh1980Use, Abdallah1993Delayed, Lee1969Tau, Michiels2014Stability}. In this section, we discuss how one can compute the ``free'' coefficients of \eqref{eq:DDE} in order to impose that a given real number $s_0$ is a root of multiplicity at least $\Lambda$ of the quasipolynomial $\Delta$ from \eqref{eq:Delta-DDE}.

Given a desired root $s_0$ and a its desired multiplicity $\Lambda$, we assume that $\Lambda-1$ coefficients of \eqref{eq:DDE} can be freely chosen, the other $D - \Lambda + 1$ coefficients being fixed (where $D = n + m + 1$ denotes the degree of \eqref{eq:Delta-DDE}). In this configuration, which appears in several practical applications (see, e.g., \cite{Seyranian2003Multiparameter}), the $\Lambda-1$ degrees of freedom in the coefficients of the system are typically not enough to ensure the existence of a root of multiplicity $\Lambda$ by themselves, and an additional constraint on the root $s_0$ and the delay $\tau$ must be imposed. This constraint can be formulated as a certain meromorphic function on $(s_0, \tau)$, depending on the known coefficients of \eqref{eq:DDE}, being equal to $0$. This meromorphic function will be called \emph{elimination-produced function}, and the main goal of this section is to describe it, first in a general setting, and then in the particular case where the known coefficients are the coefficients of the monomials of highest degree in the nondelayed part of \eqref{eq:Delta-DDE}, in which case we will prove that this function is a polynomial.

\subsubsection{Elimination-produced function}
\label{sec:iMID-EPF}

Let $n$ and $m$ be nonnegative integers with $n \geq m$, $D = n + m + 1$, $\Lambda \in \llbracket n + 1, D\rrbracket$, $I \subset \llbracket 1, D\rrbracket$ be a set of cardinality $\Lambda - 1$, and $\overline I = \llbracket 1, D\rrbracket \setminus I$. We define the holomorphic functions $\ell\colon \mathbb C^2 \to \mathcal M_{1, D}(\mathbb C)$, $\ell_I\colon \mathbb C^2 \to \mathcal M_{1, \Lambda-1}(\mathbb C)$, $\ell_{\overline I}\colon \mathbb C^2 \to \mathcal M_{1, D-\Lambda+1}(\mathbb C)$, $A\colon \mathbb C^2 \to \mathcal M_{\Lambda-1, D}(\mathbb C)$, $A_I\colon \mathbb C^2 \allowbreak \to \mathcal M_{\Lambda-1, \Lambda-1}(\mathbb C)$, $A_{\overline I}\colon \mathbb C^2 \allowbreak \to \mathcal M_{\Lambda-1, D - \Lambda + 1}(\mathbb C)$, and $b\colon \mathbb C \to \mathcal M_{\Lambda-1, 1}(\mathbb C)$, the real number $b_\Lambda$, and the matrices $X \in \mathcal M_{D, 1}(\mathbb C)$, $X_I \in \mathcal M_{\Lambda - 1, 1}(\mathbb C)$, and $X_{\overline I} \in \mathcal M_{D - \Lambda + 1, 1}(\mathbb C)$ as follows. 
 
For $(s_0, \tau) \in \mathbb C^2$, we set
\begin{equation}
\label{l}
\ell(s_0, \tau) = \begin{pmatrix}s_0^{n-1} & \cdots & s_0 & 1 & s_0^m e^{-s_0 \tau} & \cdots & s_0 e^{-s_0 \tau} & e^{-s_0 \tau}\end{pmatrix},
\end{equation}
$\ell_I(s_0, \tau)$ (respectively, $\ell_{\overline I}(s_0, \tau)$) is the row matrix whose entries are the columns of $\frac{\partial^{\Lambda-1} \ell}{\partial s_0^{\Lambda-1}}(s_0, \tau)$ whose indices belong to $I$ (respectively, to $\overline I$),
\begin{equation}
\label{MatA}
A(s_0, \tau) = \begin{pmatrix}
\ell(s_0, \tau) \\
\frac{\partial \ell}{\partial s_0}(s_0, \tau) \\
\vdots \\
\frac{\partial^{\Lambda-2} \ell}{\partial s_0^{\Lambda-2}}(s_0, \tau)
\end{pmatrix},
\end{equation}
$A_I(s_0, \tau)$ (respectively, $A_{\overline I}(s_0, \tau)$) is the matrix whose columns are the columns of $A(s_0, \tau)$ whose indices belong to $I$ (respectively, to $\overline I$), $b(s_0)$ is the column matrix with $\Lambda - 1$ rows whose first entry is $-s_0^n$ and such that its $k$th entry $b_k(s_0)$, for $k \in \llbracket 2, \Lambda-1\rrbracket$, satisfies $b_k(s_0) = b_{k-1}^\prime(s_0)$. We also set
\begin{equation*}
b_\Lambda = \begin{dcases*}
-n!, &  if $\Lambda = n+1$,\\
0, & otherwise.
\end{dcases*}
\end{equation*}
Finally, we let $X$ be the column matrix
\begin{equation*}
X = \begin{pmatrix}a_{n-1} & \cdots & a_1 & a_0 & \alpha_m & \cdots & \alpha_1 & \alpha_0\end{pmatrix}^\transp,
\end{equation*}
and we denote by $X_I$ (respectively, $X_{\overline I}$) the column matrix whose rows are the rows of $X$ whose indices belong to $I$ (respectively, to $\overline I$).

The set $I$ introduced above represents the indices of the coefficients of \eqref{eq:DDE} that are assumed to be free, when arranged in the column matrix of coefficients $X$, 
and $\overline I$ represents the indices of the known coefficients of \eqref{eq:DDE} in this same vector. Note that the order of the coefficients in the vector $X$ corresponds to the order of the entries of the row matrix $\ell(s_0, \tau)$, so that the quasipolynomial $\Delta$ from \eqref{eq:Delta-DDE} satisfies $\Delta(s_0) = s_0^n + \ell(s_0, \tau) X$.

With the above notation, we are now in position to introduce the elimination-produced function.

\begin{definition}
\label{def:EPF}
Assume that $\det A_I(\cdot, \cdot)$ is not identically zero. The \emph{elimination-produced function} is the function $\mathbb P$ defined for $(s_0, \tau) \in \mathbb C\times\mathbb R_+$ such that $\det A_I(s_0, \tau) \neq 0$ by
\begin{equation}
\label{eq:def-EPF}
\mathbb P(s_0, \tau) = \ell_{\overline I}(s_0, \tau) X_{\overline I} + \ell_{I}(s_0, \tau) A_I(s_0, \tau)^{-1} (b(s_0) - A_{\overline I}(s_0, \tau) X_{\overline I}) - b_\Lambda.
\end{equation}
\end{definition}

\begin{remark}
Note that $\det A_I(s_0, \tau)$ is a quasipolynomial in the variable $s_0$, depending only on $\tau$, $n$, $m$, $\Lambda$, and $I$. In particular, when $\det A_I$ is not identically zero, the elimination-produced function $\mathbb P$ is meromorphic and can be expressed as a fraction of quasipolynomials, its denominator being $\det A_I(s_0, \tau)$.
\end{remark}

\begin{example}
To illustrate the above notation and Definition~\ref{def:EPF}, consider the case $n = 2$ and $m = 1$, in which \eqref{eq:Delta-DDE} reads
\[
\Delta(s) = s^2 + a_1 s + a_0 + e^{-s \tau} (\alpha_1 s + \alpha_0),
\]
and the degree of $\Delta$ is $D = n + m + 1 = 4$. We are interested in this example in roots of multiplicity $\Lambda = 3$. We assume that $a_1$ and $\alpha_0$ are known, and $a_0$ and $\alpha_1$ are free. In this case, we have $I = \{2, 3\}$, which correspond to the indices of $a_0$ and $\alpha_1$ in the vector $(a_1, a_0, \alpha_1, \alpha_0)$, and thus $\overline I = \{1, 4\}$. We compute
\begin{align*}
\ell(s_0, \tau) & = \begin{pmatrix} s_0 & 1 & s_0 e^{-s_0 \tau} & e^{-s_0 \tau}\end{pmatrix}, \displaybreak[0] \\
\ell_I(s_0, \tau) & = \begin{pmatrix} 0 & (\tau^2 s_0 - 2 \tau) e^{-s_0 \tau}\end{pmatrix}, \displaybreak[0] \\
\ell_{\overline I}(s_0, \tau) & = \begin{pmatrix} 0 & \tau^2 e^{-s_0 \tau}\end{pmatrix}, \displaybreak[0] \\
A(s_0, \tau) & = \begin{pmatrix}
s_0 & 1 & s_0 e^{-s_0 \tau} & e^{-s_0 \tau} \\
1 & 0 & (1 - s_0 \tau) e^{-s_0 \tau} & -\tau e^{-s_0 \tau}
\end{pmatrix}, \displaybreak[0] \\
A_I(s_0, \tau) & = \begin{pmatrix}
1 & s_0 e^{-s_0 \tau} \\
0 & (1 - s_0 \tau) e^{-s_0 \tau}
\end{pmatrix}, \displaybreak[0] \\
A_{\overline I}(s_0, \tau) & = \begin{pmatrix}
s_0 & e^{-s_0 \tau} \\
1 & -\tau e^{-s_0 \tau}
\end{pmatrix}, \displaybreak[0] \\
b(s_0) & = \begin{pmatrix}
-s_0^2 \\
- 2 s_0
\end{pmatrix}, \displaybreak[0] \\
b_\Lambda & = -2, \\
X_{\overline I} & = \begin{pmatrix}
a_1 \\
\alpha_0
\end{pmatrix}. \\
\intertext{In particular, we have}
\det A_I(s_0, \tau) & = (1 - s_0 \tau) e^{-s_0 \tau},
\end{align*}
which is not identically zero. Thus, the elimination-produced function $\mathbb P$ is given by
\[
\mathbb P(s_0, \tau) = \frac{a_{1} s_{0} \tau^{2} - 2 a_{1} \tau + \alpha_{0} \tau^{2} e^{-s_{0} \tau} + 2 s_{0}^{2} \tau^{2} - 2 s_{0} \tau - 2}{s_{0} \tau - 1}.
\]
Note that, among the coefficients $a_1, a_0, \alpha_1, \alpha_0$, the only ones appearing in the expression of $\mathbb P$ are the known coefficients $a_1$ and $\alpha_0$.
\end{example}

The main interest of the notation we introduced above, and in particular of the eli\-mi\-na\-tion-produced function from Definition~\ref{def:EPF}, is that it allows one to express conditions for a real number $s_0$ to be a root of multiplicity $\Lambda$ of $\Delta$ in terms of the free coefficients of \eqref{eq:DDE}, as expressed in the following theorem.

\begin{theorem}
\label{thm:charact-EPF}
Let $n, m$ be nonnegative integers with $n \geq m$, $\tau > 0$, $a_0, \dotsc, a_{n-1}, \alpha_0, \dotsc, \alpha_m$ be real numbers, and consider the quasipolynomial $\Delta$ from \eqref{eq:Delta-DDE}. Let $D = n + m + 1$ denote the degree of $\Delta$ and take $s_0 \in \mathbb R$ and  $\Lambda \in \llbracket n+1, D\rrbracket$. Let $I \subset \llbracket 1, D\rrbracket$ be a set of cardinality $\Lambda-1$ and $\overline I = \llbracket 1, D\rrbracket \setminus I$. Consider the functions $\ell_I$, $\ell_{\overline I}$, $A_I$, $A_{\overline I}$, $b$, $\mathbb P$, the real number $b_\Lambda$, and the column matrices $X_I$ and $X_{\overline I}$ defined in \eqref{l}--\eqref{eq:def-EPF}.

Assume that $\det A_I(s_0, \tau) \neq 0$. Then $s_0$ is a root of multiplicity at least $\Lambda$ of $\Delta$ if and only if
\begin{equation}
\label{eq:multiplicity-M-EPF}
X_I = A_I(s_0, \tau)^{-1} \left(b(s_0) - A_{\overline I}(s_0, \tau)X_{\overline I}\right) \quad \text{ and } \quad \mathbb P(s_0, \tau) = 0.
\end{equation}
\end{theorem}

\begin{proof}
The number $s_0$ is a root of multiplicity at least $\Lambda$ of $\Delta$ if and only if $\Delta^{(k)}(s_0) = 0$ for every $k \in \llbracket 0, \Lambda-1\rrbracket$, and, recalling that $\Delta(s) = s^n + \ell(s, \tau) X$, these $\Lambda$ equations can be rewritten as
\[
\left\{
\begin{aligned}
A(s_0, \tau) X & = b(s_0), \\
\frac{\partial^{\Lambda-1}\ell}{\partial s_0^{\Lambda-1}}(s_0, \tau) X & = b_\Lambda,
\end{aligned}
\right.
\]
where $\ell$ and $A$ are defined in \eqref{l} and \eqref{MatA}. Splitting $X$ into $X_I$ and $X_{\overline I}$, the above system can be rewritten as
\[
\left\{
\begin{aligned}
A_I(s_0, \tau) X_I + A_{\overline I}(s_0, \tau) X_{\overline I} & = b(s_0), \\
\ell_I(s_0, \tau) X_I + \ell_{\overline I}(s_0, \tau) X_{\overline I} & = b_\Lambda,
\end{aligned}
\right.
\]
and, using the assumption $\det A_I(s_0, \tau) \neq 0$, we obtain that the above system is equivalent to \eqref{eq:multiplicity-M-EPF}, as required.
\end{proof}

Note that the first equation of \eqref{eq:multiplicity-M-EPF} provides explicit values of the unknown coefficients from $X_I$ in terms of $s_0$, $\tau$, and the known coefficients in $X_{\overline I}$, while the second one provides an additional relation that must be satisfied between $s_0$, $\tau$, and the known coefficients in $X_{\overline I}$ in order to ensure multiplicity at least $\Lambda$ for $s_0$.

\begin{remark}
In the above framework, we have considered that the values of some of the coefficients of \eqref{eq:DDE} are known. Our approach can also be generalized to the case where, instead of knowing values of the coefficients, we have other information on them, such as the fact that they satisfy some linear dependence relation, written under the form $C X = d$ for some matrix $C \in \mathcal M_{D - \Lambda + 1, D}(\mathbb C)$ of full rank and some $d \in \mathbb R^{D - \Lambda + 1}$. Indeed, in this case, we decompose the space $\mathbb R^D$ into the direct sum $\Ker C \oplus (\Ker C)^\perp$, and we replace the decomposition into indices in $I$ and $\overline I$ done in the construction above by projections into $\Ker C$ and $(\Ker C)^\perp$.
\end{remark}

\subsubsection{Elimination-produced polynomial}

In Section~\ref{sec:iMID-EPF}, $I$ is an arbitrary subset of $\llbracket 1, D\rrbracket$ with $\Lambda-1$ elements, meaning that we can consider any choice of free parameters of the quasipolynomial $\Delta$ from \eqref{eq:Delta-DDE}. With the aim of obtaining more precise results in a particular situation, we now consider the case where $I = \llbracket D - \Lambda + 2, D\rrbracket$, and thus $\overline I = \llbracket 1, D - \Lambda + 1\rrbracket$. In other words, the $D - \Lambda + 1$ coefficients of highest degree of the polynomial $s^n + \sum_{k = 0}^{n - 1} a_k s^k$ from \eqref{eq:Delta-DDE} are all known and fixed, and the $\Lambda - 1$ free parameters are the coefficients of the polynomial $\sum_{k=0}^m \alpha_k s^k$ from \eqref{eq:Delta-DDE} and those of lowest degree of $s^n + \sum_{k = 0}^{n - 1} a_k s^k$ (except in the particular case $n = m$ and $\Lambda = n + 1$, in which case the coefficient $\alpha_m$ is also assumed to be known).

Our first result in this setting is the following characterization of the eli\-mi\-na\-tion-pro\-duced function, which implies that it is a polynomial in $s_0$ and $\tau$.

\begin{proposition}
\label{prop:EPP}
Let $n, m$ be nonnegative integers with $n \geq m$, $\tau > 0$, $a_0, \dotsc, a_{n-1}, \alpha_0, \dotsc, \alpha_m$ be real numbers, and consider the quasipolynomial $\Delta$ from \eqref{eq:Delta-DDE}. Let $D = n + m + 1$ denote the degree of $\Delta$ and $\Lambda \in \llbracket n+1, D\rrbracket$ be such that $\Lambda \geq m + 2$. Let $\mathbb P$ be the elimination-produced function from Definition~\ref{def:EPF} in the case $I = \llbracket D - \Lambda + 2, D\rrbracket$. Then
\begin{equation}
\label{eq:EPP}
\mathbb P(s_0, \tau) = \sum_{k=0}^{m+1} \binom{m+1}{k} \tau^k P_{F, \Lambda}^{(\Lambda - 1 - k)} (s_0),
\end{equation}
where the polynomial $P_{F, \Lambda}$ is defined for $s \in \mathbb C$ by
\begin{equation}
\label{eq:P_FM}
P_{F, \Lambda}(s) = s^n + \sum_{k = 1}^{D - \Lambda + 1} a_{n-k} s^{n-k}.
\end{equation}
\end{proposition}

\begin{proof}
Note that, with the notation of Section~\ref{sec:iMID-EPF}, we have
\begin{equation}
\label{eq:proof-EPP-part1}
b(s_0) - A_{\overline I}(s_0, \tau) X_{\overline I} = \begin{pmatrix}
- P_{F, \Lambda}(s_0) \\
- P_{F, \Lambda}^\prime(s_0) \\
\vdots \\
- P_{F, \Lambda}^{(\Lambda-2)}(s_0)
\end{pmatrix}
\end{equation}
and
\begin{equation}
\label{eq:proof-EPP-part2}
\ell_{\overline I}(s_0, \tau) X_{\overline I} - b_\Lambda = P_{F, \Lambda}^{(\Lambda-1)}(s_0).
\end{equation}
Let
\[
\widetilde\ell(s_0, \tau) = \begin{pmatrix}
s_0^{\Lambda-m-3} & \cdots & s_0 & 1 & s_0^m e^{-s_0 \tau} & \cdots & s_0 e^{-s_0 \tau} & e^{-s_0 \tau}
\end{pmatrix}
\]
and note that
\[
A_I(s_0, \tau) = \begin{pmatrix}
\widetilde\ell(s_0, \tau) \\
\frac{\partial\widetilde\ell}{\partial s_0}(s_0, \tau) \\
\vdots \\
\frac{\partial^{\Lambda-2}\widetilde\ell}{\partial s_0^{\Lambda-2}}(s_0, \tau) \\
\end{pmatrix}
\]
The matrix $A_I(s_0, \tau)$ admits the block decomposition
\begin{equation}
\label{eq:AI-block}
A_I(s_0, \tau) = \begin{pmatrix}
A_{I, 0}(s_0) & A_{I, 1}(s_0, \tau) \\
0 & A_{I, 2}(s_0, \tau) \\
\end{pmatrix},
\end{equation}
where
\begin{align*}
A_{I, 0}(s_0) & = \begin{pmatrix}
s_0^{\Lambda - m - 3} & \cdots & s_0 & 1 \\
(\Lambda - m - 3) s_0^{\Lambda - m - 4} & \cdots & 1 & 0 \\
\vdots & \iddots & \vdots & \vdots \\
(\Lambda - m - 3)! & \cdots & 0 & 0 \\
\end{pmatrix} \in \mathcal M_{\Lambda - m - 2}(\mathbb R), \displaybreak[0] \\
A_{I, 1}(s_0, \tau) & = \begin{pmatrix}
\widehat\ell(s_0, \tau) \\
\frac{\partial\widehat\ell}{\partial s_0}(s_0, \tau) \\
\vdots \\
\frac{\partial^{\Lambda - m - 3}\widehat\ell}{\partial s_0^{\Lambda - m - 3}}(s_0, \tau) \\
\end{pmatrix} \in \mathcal M_{\Lambda - m - 2, m + 1}(\mathbb R), \displaybreak[0] \\
A_{I, 2}(s_0, \tau) & = \begin{pmatrix}
\frac{\partial^{\Lambda - m - 2}\widehat\ell}{\partial s_0^{\Lambda - m - 2}}(s_0, \tau) \\
\vdots \\
\frac{\partial^{\Lambda - 2}\widehat\ell}{\partial s_0^{\Lambda - 2}}(s_0, \tau) \\
\end{pmatrix} \in \mathcal M_{m + 1, m + 1}(\mathbb R), \displaybreak[0] \\
\intertext{and}
\widehat\ell(s_0, \tau) & = \begin{pmatrix} s_0^m e^{-s_0 \tau} & \cdots & s_0 e^{-s_0 \tau} & e^{-s_0 \tau}
\end{pmatrix} \in \mathcal M_{1, m + 1}(\mathbb R),
\end{align*}
Clearly, $A_{I, 0}(s_0)$ is invertible. Set $\rho(s_0, \tau) = e^{s_0 \tau} \frac{\partial^{\Lambda - m - 2}\widehat\ell}{\partial s_0^{\Lambda - m - 2}}(s_0, \tau)$ and note that
\begin{equation}
\label{eq:rho-row}
\rho(s_0, \tau) = \begin{pmatrix}
\mathfrak q_m(s_0, \tau) & \cdots & \mathfrak q_1(s_0, \tau) & \mathfrak q_0(s_0, \tau) \\
\end{pmatrix},
\end{equation}
where, for $j \in \llbracket 0, m\rrbracket$, the function $s_0 \mapsto \mathfrak q_j(s_0, \tau)$ is a polynomial of degree $j$, with coefficients depending on $\tau$, and its term of highest degree is $(-\tau)^{\Lambda - m - 2} s_0^j$. Thus
\[
\begin{pmatrix}
\rho(s_0, \tau) \\
\frac{\partial\rho}{\partial s_0}(s_0, \tau) \\
\vdots \\
\frac{\partial^m\rho}{\partial s_0^m}(s_0, \tau) \\
\end{pmatrix} = 
\begin{pmatrix}
\mathfrak q_m(s_0, \tau) & \cdots & \mathfrak q_2(s_0, \tau) & \mathfrak q_1(s_0, \tau) & (-\tau)^{\Lambda - m - 2} \\
\frac{\partial\mathfrak q_m}{\partial s_0}(s_0, \tau) & \cdots & \frac{\partial\mathfrak q_2}{\partial s_0}(s_0, \tau) & (-\tau)^{\Lambda - m - 2} & 0 \\
\frac{\partial^2\mathfrak q_m}{\partial s_0^2}(s_0, \tau) & \cdots & 2 (-\tau)^{\Lambda - m - 2} & 0 & 0 \\
\vdots & \iddots & \vdots & \vdots & \vdots \\
m! (-\tau)^{\Lambda - m - 2} & \cdots & 0 & 0 & 0
\end{pmatrix},
\]
and, in particular, since $\tau > 0$, the above matrix is invertible. On the other hand, taking derivatives directly in the definition of $\rho$, we have
\[
\begin{pmatrix}
\rho(s_0, \tau) \\
\frac{\partial\rho}{\partial s_0}(s_0, \tau) \\
\vdots \\
\frac{\partial^m\rho}{\partial s_0^m}(s_0, \tau) \\
\end{pmatrix} = e^{s_0 \tau} \begin{pmatrix}
1 & 0 & 0 & \cdots & 0 & 0 \\
\tau & 1 & 0 & \cdots & 0 & 0 \\
\vdots & \ddots & \ddots & \ddots & \ddots & \vdots \\
\binom{k}{k} \tau^k & \binom{k}{k-1} \tau^{k-1} & \ddots & 1 & \ddots & 0 \\
\vdots & \vdots & \ddots & \ddots & \ddots & \vdots \\
\binom{m}{m} \tau^m & \binom{m}{m-1} \tau^{m-1} & \cdots & \binom{m}{2} \tau^2 & m \tau & 1 \\
\end{pmatrix} A_{I, 2}(s_0, \tau).
\]
Hence, $A_{I, 2}(s_0, \tau)$ is invertible and, thanks to the block decomposition \eqref{eq:AI-block}, $A_I(s_0, \tau)$ is invertible.

Thanks to \eqref{eq:rho-row}, we have $\frac{\partial^{m+1}\rho}{\partial s_0^{m+1}}(s_0, \tau) = 0$ and, using the definition of $\rho$, this implies that
\[
\sum_{k=0}^{m+1} \binom{m+1}{k} \tau^k \frac{\partial^{\Lambda - k - 1}\widehat\ell}{\partial s_0^{\Lambda - k - 1}}(s_0, \tau) = 0,
\]
which can be rewritten in matrix form as
\[
\begin{pmatrix}
\binom{m+1}{m+1} \tau^{m+1} & \binom{m+1}{m} \tau^m & \cdots & \binom{m+1}{1} \tau
\end{pmatrix} A_{I, 2}(s_0, \tau) = - \frac{\partial^{\Lambda - 1}\widehat\ell}{\partial s_0^{\Lambda - 1}}(s_0, \tau).
\]
Hence
\begin{equation*}
\begin{pmatrix}
0 & \cdots & 0 & \binom{m+1}{m+1} \tau^{m+1} & \binom{m+1}{m} \tau^m & \cdots & \binom{m+1}{1} \tau
\end{pmatrix} A_{I}(s_0, \tau) = - \begin{pmatrix} 0 & \frac{\partial^{\Lambda - 1}\widehat\ell}{\partial s_0^{\Lambda - 1}}(s_0, \tau) \end{pmatrix} = - \ell_I(s_0, \tau),
\end{equation*}
so that
\begin{equation}
\label{eq:proof-EPP-part3}
\ell_I(s_0, \tau) A_I(s_0, \tau)^{-1} = - \begin{pmatrix}
0 & \cdots & 0 & \binom{m+1}{m+1} \tau^{m+1} & \cdots & \binom{m+1}{1} \tau
\end{pmatrix}.
\end{equation}
Inserting \eqref{eq:proof-EPP-part1}, \eqref{eq:proof-EPP-part2}, and \eqref{eq:proof-EPP-part3} into \eqref{eq:def-EPF}, we finally deduce \eqref{eq:EPP}.
\end{proof}

When the assumptions of Proposition~\ref{prop:EPP} are satisfied, the elimination-produced function from \eqref{eq:EPP} is a polynomial in $(s_0, \tau)$. In this case, we refer to $\mathbb P$ as the \emph{elimination-produced polynomial}.

\begin{remark}
Proposition~\ref{prop:EPP} requires the additional assumption that $\Lambda \geq m+2$. This is always satisfied if $n > m$, since $\Lambda \geq n+1$, and hence the only case not covered by Proposition~\ref{prop:EPP} is the case of neutral systems with the smallest overorder multiplicity, that is, $n = m$ and $\Lambda = n + 1$. In this case, $\mathbb P$ is not necessarily a polynomial. Indeed, consider the case $n = m = 1$, i.e.,
\[
\Delta(s) = s + a_0 + e^{-s \tau} (\alpha_1 s + \alpha_0),
\]
in which case $D = 3$. The lowest overorder multiplicity is $\Lambda = 2$ and, assuming $a_0$ and $\alpha_1$ to be known and $\alpha_0$ to be free (i.e., $I = \{3\}$ and $\overline I = \{1, 2\}$), straightforward computations from its definition show that the elimination-produced function is
\[
\mathbb P(s_0, \tau) = 1 + \tau s_0 + \tau a_0 + \alpha_1 e^{-s_0 \tau},
\]
which is not a polynomial in $(s_0, \tau)$.
\end{remark}

Note that the elimination-produced polynomial $\mathbb P$ from \eqref{eq:EPP} depends on the desired minimal multiplicity $\Lambda$. An interesting fact, proven in the next result, is that the elimination-produced polynomial for a root of multiplicity at least $\Lambda + 1$ can be obtained by taking the derivative of the corresponding elimination-produced polynomial for a root of multiplicity $\Lambda$.

\begin{proposition}
\label{prop:EPP-derivative}
Let $n, m$ be nonnegative integers with $n \geq m$, $\tau > 0$, $a_0, \dotsc, a_{n-1}, \alpha_0, \dotsc, \alpha_m$ be real numbers, and consider the quasipolynomial $\Delta$ from \eqref{eq:Delta-DDE}. Let $D = n + m + 1$ denote the degree of $\Delta$ and $\Lambda \in \llbracket n+1, D - 1\rrbracket$ be such that $\Lambda \geq m+2$. Let $\mathbb P_\Lambda$ (respectively, $\mathbb P_{\Lambda+1}$) be the elimination-produced polynomial from Proposition~\ref{prop:EPP} for a root of multiplicity at least $\Lambda$ (respectively, at least $\Lambda+1$). Then
\[
\mathbb P_{\Lambda+1}(s_0, \tau) = \frac{\partial \mathbb P_\Lambda}{\partial s_0}(s_0, \tau).
\]
\end{proposition}

\begin{proof}
Let $P_{F, \Lambda}$ and $P_{F, \Lambda+1}$ be defined as in \eqref{eq:P_FM}. Then
\[
P_{F, \Lambda}(s) = P_{F, \Lambda+1}(s) + a_{\Lambda - m - 2} s^{\Lambda - m - 2}.
\]
Hence, for every integer $k \in \llbracket 0, m+1\rrbracket$, we have $\Lambda - m - 2 < \Lambda - k$, and thus
\[
P_{F, \Lambda}^{(\Lambda - k)}(s) = P_{F, \Lambda+1}^{(\Lambda - k)}(s).
\]
It then follows from \eqref{eq:EPP} that
\begin{equation*}
\mathbb P_{\Lambda+1}(s_0, \tau) = \sum_{k=0}^{m+1} \binom{m+1}{k} \tau^k P_{F, \Lambda+1}^{(\Lambda - k)} (s_0) = \sum_{k=0}^{m+1} \binom{m+1}{k} \tau^k P_{F, \Lambda}^{(\Lambda - k)} (s_0) = \frac{\partial \mathbb P_\Lambda}{\partial s_0}(s_0, \tau). \qedhere
\end{equation*}
\end{proof}

Our next provides an explicit expression of the elimination-produced polynomial $\mathbb P$ from \eqref{eq:EPP} in terms of the known coefficients $a_{\Lambda - m - 2}, \dotsc, a_{n-1}$ of $\Delta$.

\begin{proposition}
\label{prop:EPP-explicit}
Let $n, m$ be nonnegative integers with $n \geq m$, $\tau > 0$, $a_0, \dotsc, a_{n-1},\allowbreak \alpha_0, \dotsc, \alpha_m$ be real numbers, and consider the quasipolynomial $\Delta$ from \eqref{eq:Delta-DDE}. Let $D = n + m + 1$ denote the degree of $\Delta$ and $\Lambda \in \llbracket n+1, D\rrbracket$ be such that $\Lambda \geq m+2$. Let $\mathbb P$ be the elimination-produced function from Definition~\ref{def:EPF} in the case $I = \llbracket D - \Lambda + 2, D\rrbracket$. Then
\begin{equation*}
\mathbb P(s_0, \tau) = \sum_{j = 0}^{D - \Lambda + 1} \sum_{\ell = j + \Lambda - m - 2}^{n} \binom{m + 1}{\Lambda - \ell - 1 + j} \frac{\ell!}{j!} \tau^{\Lambda - \ell - 1 + j} a_\ell s_0^j,
\end{equation*}
with the convention $a_n = 1$.
\end{proposition}

The proof of the above proposition follows by inserting \eqref{eq:P_FM} into \eqref{eq:EPP} and standard changes of summation variables in order to group the terms in powers of $s_0$.

As a final result concerning the elimination-produced polynomial, we prove the following proposition, which provides a link, at least in some special cases, between the polynomial $t \mapsto \mathbb P(s_0, \tau t)$ and the polynomial appearing in the integral in \eqref{eq:Delta-max-mult}.

\begin{proposition}
\label{prop:EPP-factorization}
Let $n$, $m$ be nonnegative integers with $n \geq 1$ and $m \in \{n - 1, n\}$, $\tau > 0$, $a_0, \dotsc, a_{n-1},\allowbreak \alpha_0, \dotsc, \alpha_m$ be real numbers, and consider the quasipolynomial $\Delta$ from \eqref{eq:Delta-DDE}. Let $D = n + m + 1$ denote the degree of $\Delta$ and take $s_0 \in \mathbb R$. Set $\Lambda = m + 2$ and denote by $\mathbb P$ be the elimination-produced function from Definition~\ref{def:EPF} in the case $I = \llbracket n + 1, D\rrbracket$. Then the real number $s_0$ is a root of multiplicity at least $\Lambda$ of $\Delta$ if and only $\mathbb P(s_0, \tau) = 0$ and
\begin{equation}
\label{eq:Delta-EPP}
\Delta(s) = \frac{\tau (s - s_0)^{m + 2}}{(m + 1)!} \int_0^1 \mathbb P(s_0, \tau t) e^{-t \tau(s - s_0)} \diff t, \qquad \text{ for all } s \in \mathbb C.
\end{equation}
\end{proposition}

\begin{proof}
Note that $\Lambda \in \llbracket n + 1, D\rrbracket$, and thus, by Proposition~\ref{prop:EPP}, $\mathbb P$ is a polynomial. In particular, by \eqref{eq:EPP}, $t \mapsto \mathbb P(s_0, \tau t)$ is a polynomial of degree at most $m + 1$. The polynomial $P_{F, \Lambda}$ from \eqref{eq:P_FM} is given, for $s \in \mathbb C$, by
\begin{equation}
\label{eq:P_FLmabda-is-P0}
P_{F, \Lambda}(s) = s^n + \sum_{k = 0}^{n-1} a_k s^k,
\end{equation}
and thus, using \eqref{eq:EPP}, we obtain at once that, if $m = n - 1$, we have $\mathbb P(s_0, 0) = n!$, while, for $m = n$, we have $\mathbb P(s_0, 0) = 0$ and $\left.\frac{\df}{\df t} \mathbb P(s_0, \tau t)\right\rvert_{t = 0} = (n + 1)! \tau$.

Assume first that $\mathbb P(s_0, \tau) = 0$ and that $\Delta$ is given by \eqref{eq:Delta-EPP}. Let $q$ be the polynomial defined by $q(t) = \mathbb P(s_0, \tau t)$ and note that $q(1) = 0$ and, in the case $m = n - 1$, we have $q(0) = n!$, while, in the case $m = n$, we have $q(0) = 0$ and $q^\prime(0) = (n + 1)! \tau$. Hence, there exists a polynomial $p$ of degree at most $n - 1$ and with $p(0) = 1$ such that
\[
\mathbb P(s_0, \tau t) = \begin{dcases*}
n! (1 - t) p(t), & if $m = n - 1$, \\
(n + 1)! \tau t (1 - t) p(t), & if $m = n$.
\end{dcases*}
\]
Inserting into \eqref{eq:Delta-EPP}, we deduce that, for $s \in \mathbb C$, we have
\begin{equation}
\label{eq:explicit-Delta-n-1-n}
\Delta(s) = \begin{dcases*}
\tau (s - s_0)^{m + 2} \int_0^1 (1 - t) p(t) e^{-t \tau(s - s_0)} \diff t, & if $m = n - 1$, \\
\tau^2 (s - s_0)^{m + 2} \int_0^1 t (1 - t) p(t) e^{-t \tau(s - s_0)} \diff t, & if $m = n$,
\end{dcases*}
\end{equation}
and thus $\Delta$ has the form \eqref{eq:Delta-max-mult}. Hence, by Theorem~\ref{thm:max-mult}, we deduce that $s_0$ is a root of multiplicity at least $\Lambda$ of $\Delta$.

Conversely, assume that $s_0$ is a root of multiplicity at least $\Lambda$ of $\Delta$. Then, by Theorem~\ref{thm:charact-EPF}, we have $\mathbb P(s_0, \tau) = 0$ and, by Theorem~\ref{thm:max-mult}, there exists a polynomial $p$ of degree at most $n - 1$ with $p(0) = 0$ such that $\Delta$ is given by \eqref{eq:explicit-Delta-n-1-n} for all $s \in \mathbb C$. Let $q$ be the polynomial defined by $q(t) = (1-t) p(t)$ if $m = n-1$ and $q(t) = t (1 - t) p(t)$ if $m = n$ and note that, in both cases, $q$ is of degree at most $m + 1$ and \eqref{eq:explicit-Delta-n-1-n} can be rewritten as
\begin{equation}
\label{eq:combined-n-1-n-Delta}
\Delta(s) = \tau^{m + 2 - n} (s - s_0)^{m + 2} \int_0^1 q(t) e^{-t \tau(s - s_0)} \diff t.
\end{equation}
By Remark~\ref{remk:P0-P1-q} and using \eqref{eq:P_FLmabda-is-P0}, we have that
\begin{equation}
\label{eq:P0-q}
P_{F, \Lambda}(s) = (s - s_0)^n + \sum_{k = 0}^{n - 1} \frac{q^{(m + 1 - k)}(0)}{\tau^{n - k}} (s - s_0)^{k}.
\end{equation}
On the other hand, since $P_{F, \Lambda}$ is a polynomial of degree $n$, it coincides with its Taylor expansion of order $n$ at $s_0$, i.e.,
\begin{equation}
\label{eq:P0-Taylor}
P_{F, \Lambda}(s) = \sum_{k=0}^n \frac{P_{F, \Lambda}^{(k)}(s_0)}{k!} (s - s_0)^k.
\end{equation}
Using the uniqueness of the Taylor expansion at a given point and combining \eqref{eq:P0-q} and \eqref{eq:P0-Taylor}, we deduce that
\[
P_{F, \Lambda}^{(k)}(s_0) = \frac{k!}{\tau^{n - k}} q^{(m + 1 - k)}(0), \qquad \text{ for } k \in \llbracket 0, n \rrbracket.
\]
Hence, by \eqref{eq:EPP}, we have, for $t \in [0, 1]$,
\begin{equation}
\label{eq:EPP-and-q}
\mathbb P(s_0, \tau t) = (m + 1)! \tau^{m + 1 - n} \sum_{k=0}^{m+1} \frac{t^k}{k!} q^{(k)}(0) = (m + 1)! \tau^{m + 1 - n} q(t),
\end{equation}
where the last equality is justified by the fact that $q$ is a polynomial of degree at most $m + 1$, and hence it is equal to its Taylor expansion at $0$ of order $m + 1$. As a consequence, we obtain that $q(t) = \frac{1}{(m + 1)! \tau^{m + 1 - n}} \mathbb P(s_0, \tau t)$ and, inserting into \eqref{eq:combined-n-1-n-Delta}, we obtain that $\Delta$ is given by \eqref{eq:Delta-EPP}, as required.
\end{proof}

\begin{remark}
An important step of the proof of Proposition~\ref{prop:EPP-factorization} is to establish \eqref{eq:EPP-and-q}, which states that the polynomials $t \mapsto \mathbb P(s_0, \tau t)$ and $t \mapsto q(t) = t^{\Lambda - n - 1} (1 - t)^{\Lambda - m - 1} p(t)$ are the same, up to a multiplicative constant. By \eqref{eq:EPP}, the polynomial $t \mapsto \mathbb P(s_0, \tau t)$ has degree at most $m + 1$, while, using the fact from Theorem~\ref{thm:max-mult} that $p$ has degree at most $D - \Lambda$, we deduce that $q$ has degree at most $\Lambda - 1$. For an equality such as \eqref{eq:EPP-and-q} to hold, it is natural to expect the maximal degrees of $q$ and $t \mapsto \mathbb P(s_0, \tau t)$ to coincide, and this is the case if an only if $\Lambda = m + 2$, justifying thus the corresponding assumption in the statement of Proposition~\ref{prop:EPP-factorization}. Theorems~\ref{thm:max-mult} and \ref{thm:charact-EPF} and Proposition~\ref{prop:EPP} all require $\Lambda \in \llbracket n+1, D\rrbracket$ and $n \geq m$, and, for $\Lambda = m + 2$, these conditions are equivalent to $n \geq 1$ and $m \in \{n - 1, n\}$, justifying thus the corresponding assumptions in Proposition~\ref{prop:EPP-factorization}. In other words, the conditions $m \in \{n - 1, n\}$, $n \geq 1$, and $\Lambda = m + 2$ assumed in Proposition~\ref{prop:EPP-factorization} are the only setting in which the maximal degrees of the polynomials $t \mapsto \mathbb P(s_0, \tau t)$ and $q$ coincide. Finding relations between $t \mapsto \mathbb P(s_0, \tau t)$ and $q$ in a more general setting is an interesting open question.
\end{remark}

Note that Theorems~\ref{thm:max-mult} and \ref{thm:charact-EPF} and Corollary~\ref{coro:max-mult-Kummer} provide characterizations of the existence of a root $s_0$ of multiplicity at least $\Lambda$, generalizing thus the equivalences between \ref{item:max-mult}--\ref{item:DeltaTilde} from Theorem~\ref{thm:main-MID}, while Propositions~\ref{prop:EPP}, \ref{prop:EPP-derivative}, \ref{prop:EPP-explicit}, and \ref{prop:EPP-factorization} only provide additional properties of the elimination-produced function from Definition~\ref{def:EPF} in the particular case where $I = \llbracket D - \Lambda + 2, D\rrbracket$. In particular, there are no results guaranteeing dominance of the root $s_0$ of multiplicity at least $\Lambda$. As recalled before, a general result on dominance of such a multiple root is still lacking in the literature, with existing works dealing only with particular cases or with sufficient conditions that are not expected to be sharp. However, the characterizations of the elimination-produced function from Propositions~\ref{prop:EPP}, \ref{prop:EPP-derivative}, \ref{prop:EPP-explicit}, and \ref{prop:EPP-factorization} can be quite useful when analyzing dominance in particular situations. We shall illustrate this fact later in Section~\ref{sec:expl-integrator}.

\section{The case of complex conjugate roots}
\label{sec:MID-complex}

The results in Sections~\ref{sec:MID-maximal} and \ref{sec:MID-intermediate} only deal with the MID property for \emph{real} roots of maximal or large enough multiplicity. A natural question is whether a dominance property can be proved for a pair of complex conjugate roots with large enough multiplicity. This section reports the contents of \cite{Mazanti2020Spectral} on this topic, which provides an affirmative answer to this question for systems of the form \eqref{eq:DDE} when $n = 2$ and $m = 1$ and a pair of complex conjugate roots of multiplicity $2$ each. Note that, in this configuration, the sum of the multiplicities of both roots is equal to the degree $D = 4$ of the corresponding quasipolynomial.

\subsection{Statement of the main result}

The main result we will prove in this section is the following.

\begin{theorem}
\label{thm:complex-dominance}
Let $\tau > 0$, $a_0, a_1, \alpha_0, \alpha_1$ be real numbers, and consider the quasipolynomial $\Delta$ from \eqref{eq:Delta-DDE} with $n = 2$ and $m = 1$. Take $s_0 \in \mathbb C \setminus \mathbb R$ and set $\sigma_0 = \Real s_0$ and $\theta_0 = \Imag s_0 \neq 0$. Then $s_0$ and $\overline s_0$ are roots of $\Delta$ of multiplicity $2$ each if and only if the coefficients $a_0, a_1, \alpha_0, \alpha_1$, the values $\sigma_0$ and $\theta_0$, and the delay $\tau$ satisfy the relations
\begin{subequations}
\label{eq:coeffs-complex-dominance}
\begin{align}
a_1 & = -2 \sigma_0 - 2 \theta_0 \frac{\tau \theta_0 - \sin(\tau \theta_0)\cos(\tau \theta_0)}{\tau^2 \theta_0^2 - \sin^2(\tau \theta_0)}, \displaybreak[0] \\
a_0 & = \sigma_0^2 + 2 \sigma_0 \theta_0\frac{\tau \theta_0 - \sin(\tau \theta_0)\cos(\tau \theta_0)}{\tau^2 \theta_0^2 - \sin^2(\tau \theta_0)} + \theta_0^2 \frac{\tau^2 \theta_0^2 + \sin^2(\tau \theta_0)}{\tau^2 \theta_0^2 - \sin^2(\tau \theta_0)}, \displaybreak[0] \\
\alpha_1 & = 2 \theta_0 e^{\sigma_0 \tau} \frac{\tau \theta_0 \cos(\tau \theta_0) - \sin(\tau \theta_0)}{\tau^2 \theta_0^2 - \sin^2(\tau \theta_0)}, \displaybreak[0] \\
\alpha_0 & = 2 \theta_0 e^{\sigma_0 \tau} \left(\sigma_0 \frac{\sin(\tau \theta_0) - \tau \theta_0 \cos(\tau \theta_0)}{\tau^2 \theta_0^2 - \sin^2(\tau \theta_0)} - \frac{\tau \theta_0^2 \sin(\tau \theta_0)}{\tau^2 \theta_0^2 - \sin^2(\tau \theta_0)}\right).
\end{align}
\end{subequations}
In addition, if \eqref{eq:coeffs-complex-dominance} is satisfied, then $\Real s < \sigma_0$ for every $s \in \mathbb C \setminus \{s_0, \overline s_0\}$ such that $\Delta(s) = 0$.
\end{theorem}

\begin{remark}
The expressions of $a_1, a_0, \alpha_1, \alpha_0$ in \eqref{eq:coeffs-complex-dominance} are well-defined for every $\theta_0 \in \mathbb R^\ast$ and $\tau > 0$, since $\sin^2(\tau\theta_0) = \tau^2\theta_0^2$ if and only if $\tau\theta_0 = 0$. These expressions are even functions of $\theta_0$, as one might expect by symmetry since one is placing both roots $s_0$ and $\overline s_0$. Moreover, as $\theta_0 \to 0$, we have
\begin{align*}
a_1 & \to -\frac{4}{\tau} - 2 \sigma_0, & \qquad a_0 & \to \frac{6}{\tau^2} + \frac{4}{\tau} \sigma_0 + \sigma_0^2, \displaybreak[0] \\
\alpha_1 & \to -\frac{2}{\tau} e^{\sigma_0 \tau}, & \alpha_0 & \to \frac{2}{\tau} e^{\sigma_0 \tau} \left(\sigma_0 - \frac{3}{\tau}\right),
\end{align*}
which are precisely the necessary and sufficient conditions for the existence of a root of multiplicity $4$ at $\sigma_0$ from \eqref{Coeffs}.
\end{remark}

\subsection{Remarks and auxiliary results}
\label{sec:MID-complex-auxiliary}

Let us comment on the strategy of proof of Theorem~\ref{thm:complex-dominance}. The first step is to proceed as in the proof of Theorem~\ref{thm:main-MID} and consider a normalized quasipolynomial $\widetilde\Delta$ defined similarly to \eqref{eq:DeltaDeltaTilde}. More precisely, we define $\widetilde\Delta\colon \mathbb C \to \mathbb C$ for $z \in \mathbb C$ by
\begin{equation}
\label{eq:complex-change-of-vars}
\widetilde\Delta(z) = \tau^2 \Delta\left(\sigma_0 + \tfrac{z}{\tau}\right),
\end{equation}
and we remark that, on the one hand, $\Delta$ admits $s_0$ and $\overline s_0$ as roots of multiplicity $2$ each if and only if $\widetilde\Delta$ admits $i \tau \theta_0$ and $-i \tau \theta_0$ as roots of multiplicity $2$ each and, on the other hand, the coefficients of $\Delta$ satisfy \eqref{eq:coeffs-complex-dominance} if and only if $\widetilde\Delta$ is given by
\begin{multline*}
\widetilde\Delta(z) = z^2  - 2 \tau \theta_0 \frac{\tau \theta_0 - \sin(\tau \theta_0)\cos(\tau \theta_0)}{\tau^2 \theta_0^2 - \sin^2(\tau \theta_0)} z + \tau^2 \theta_0^2 \frac{\tau^2 \theta_0^2 + \sin^2(\tau \theta_0)}{\tau^2 \theta_0^2 - \sin^2(\tau \theta_0)} \\ + 2 \tau \theta_0 e^{-z} \left(\frac{\tau \theta_0 \cos(\tau \theta_0) - \sin(\tau \theta_0)}{\tau^2 \theta_0^2 - \sin^2(\tau \theta_0)} z - \frac{\tau^2 \theta_0^2 \sin(\tau \theta_0)}{\tau^2 \theta_0^2 - \sin^2(\tau \theta_0)}\right).
\end{multline*}

At the light of these observations, we introduce the function $\widehat\Delta\colon \mathbb C \times \mathbb R^\ast \to \mathbb C$ defined for $(z, \vartheta) \in \mathbb C \times \mathbb R^\ast$ by
\begin{multline}
\label{eq:complex-widehat-Delta}
\widehat\Delta(z, \vartheta) = z^2  - 2 \vartheta \frac{\vartheta - \sin\vartheta \cos\vartheta}{\vartheta^2 - \sin^2 \vartheta} z + \vartheta^2 \frac{\vartheta^2 + \sin^2 \vartheta}{\vartheta^2 - \sin^2 \vartheta} \\ + 2 \vartheta e^{-z} \left(\frac{\vartheta \cos\vartheta - \sin\vartheta}{\vartheta^2 - \sin^2 \vartheta} z - \frac{\vartheta^2 \sin\vartheta}{\vartheta^2 - \sin^2 \vartheta}\right),
\end{multline}
and we remark that Theorem~\ref{thm:complex-dominance} is proved if we show that, for every $\vartheta \in \mathbb R^\ast$, a quasipolynomial of the form \eqref{eq:Delta-DDE} with $n = 2$, $m = 1$, and $\tau = 1$ admits $i \vartheta$ and $-i \vartheta$ as roots of multiplicity $2$ each if and only if this quasipolynomial is given by $\widehat\Delta(\cdot, \vartheta)$, and that $\Real z < 0$ for every $z \in \mathbb C \setminus \{i\vartheta, -i\vartheta\}$ such that $\widehat\Delta(z, \vartheta) = 0$.

Note that the function $\widehat\Delta$ can be extended by continuity to $\mathbb C \times \mathbb R$, and its extension, still denoted by $\widehat\Delta$, satisfies
\begin{equation}
\label{eq:complex-widehat-Delta-0}
\widehat\Delta(z, 0) = z^2 - 4 z + 6 - e^{-z} (2 z + 6)
\end{equation}
for every $z \in \mathbb C$. Note that, by Theorem~\ref{thm:main-MID}, $\widehat\Delta(\cdot, 0)$ is the unique quasipolynomial of the form \eqref{eq:Delta-DDE} with $n = 2$, $m = 1$, and $\tau = 1$ admitting $0$ as a root of multiplicity $4$, and that this root is strictly dominant.

Before turning to the proof of Theorem~\ref{thm:complex-dominance}, we provide two auxiliary results that will be useful for that proof.

\begin{lemma}
\label{lemm:F-positive}
Let $F\colon \mathbb R \to \mathbb R$ be given for $x \in \mathbb R$ by
\[
F(x) = x^2 + x \sin x \cos x - 2 \sin^2 x.
\]
Then $F(x) = 0$ if and only if $x = 0$.
\end{lemma}

\begin{proof}
One immediately verifies that $F(0) = 0$. Since $F$ is even, it suffices to show that $F(x) \neq 0$  for every $x > 0$. To do that, we first rewrite $F$ as
\begin{equation}
\label{FSin2x}
F(x) = x^2 + \frac{x \sin 2x}{2} + \cos 2 x - 1.
\end{equation}
Hence, for $x \geq 0$, one has the inequality
\[
F(x) \geq x^2 - \frac{x}{2} - 2.
\]
The right-hand side of the above inequality is strictly positive for every $x > \frac{1 + \sqrt{33}}{4}$, and thus $F(x) > 0$ for $x > \frac{1 + \sqrt{33}}{4}$. We are left to prove that $F(x) \neq 0$ for $x \in \left(0, \frac{1 + \sqrt{33}}{4}\right]$.

Recall that, for $x \geq 0$, one has the classical inequalities
\[
\sin x \geq x - \frac{x^3}{3!} + \frac{x^5}{5!} - \frac{x^7}{7!}, \qquad \cos x \geq 1 - \frac{x^2}{2!} + \frac{x^4}{4!} - \frac{x^6}{6!},
\]
which can be obtained by iterative integration on the interval $[0, x]$ of the inequality $\cos x \leq 1$. Using the above inequalities in \eqref{FSin2x}, one obtains that, for $x \geq 0$,
\[
F(x) \geq \frac{2 x^6}{45}\left(1 - \frac{2 x^2}{7}\right).
\]
In particular, $F(x) > 0$ for $x \in \left(0, \sqrt{7/2}\right)$. Since $\sqrt{7/2} > \frac{1 + \sqrt{33}}{4}$, this completes the proof.
\end{proof}

We also need the following results on the roots of the quasipolynomial $\widehat\Delta(\cdot, \vartheta)$ from \eqref{eq:complex-widehat-Delta} and \eqref{eq:complex-widehat-Delta-0}.

\begin{lemma}
\label{lemm:complex-preliminary-bound}
Let $\widehat\Delta$ be the function defined by \eqref{eq:complex-widehat-Delta} and \eqref{eq:complex-widehat-Delta-0} and take $\vartheta_\ast > 0$. There exists $R > 0$ such that, for every $\vartheta \in [-\vartheta_\ast, \vartheta_\ast]$, if $z_0$ is a root of $\widehat\Delta(\cdot, \vartheta)$ with $\Real z_0 \geq 0$, then $\abs{z_0} \leq R$.
\end{lemma}

\begin{proof}
Let $A_0, A_1\colon \mathbb R \to \mathcal M_2(\mathbb R)$ be the functions defined for $\vartheta \neq 0$ by
\begin{equation*}
A_0(\vartheta) = \begin{pmatrix}
0 & 1 \\
-\vartheta^2 \frac{\vartheta^2 + \sin^2 \vartheta}{\vartheta^2 - \sin^2 \vartheta} & 2 \vartheta \frac{\vartheta - \sin\vartheta \cos\vartheta}{\vartheta^2 - \sin^2 \vartheta} \\
\end{pmatrix}, \qquad
A_1(\vartheta) = \begin{pmatrix}
0 & 0 \\
\frac{2 \vartheta^3 \sin \vartheta}{\vartheta^2 - \sin^2 \vartheta} & -2 \vartheta \frac{\vartheta \cos \vartheta - \sin\vartheta}{\vartheta^2 - \sin^2 \vartheta} \\
\end{pmatrix},
\end{equation*}
and extended by continuity for $\vartheta = 0$. Notice that $\widehat\Delta(z, \vartheta) = \det(z \Id - A_0(\vartheta) - A_1(\vartheta) e^{-z})$ for every $z \in \mathbb C$ and $\vartheta \in \mathbb R$. Hence any root $z_0$ of $\widehat\Delta(\cdot, \vartheta)$ is an eigenvalue of the matrix $A_0(\vartheta) + A_1(\vartheta) e^{-z_0}$, and thus $\abs{z_0} \leq \spr(A_0(\vartheta) + A_1(\vartheta) e^{-z_0})$. Recall that, for any square matrix $M$, one has $\spr(M) \leq \abs{M}$, where $\abs{\cdot}$ is any matrix norm induced by a vector norm. Then
\begin{equation}
\label{BoundNorms}
\abs{z_0} \leq \abs{A_0(\vartheta)} + \abs{e^{-z_0}} \abs{A_1(\vartheta)}.
\end{equation}
Since $\vartheta \mapsto \abs{A_j(\vartheta)}$ is continuous in $\mathbb R$ for $j \in \{0, 1\}$, given $\vartheta_\ast > 0$, there exists $M > 0$ such that $\abs{A_j(\vartheta)} \leq M$ for every $\vartheta \in [-\vartheta_\ast, \vartheta_\ast]$ and $j \in \{0, 1\}$. Letting $R = 2 M$, we obtain from \eqref{BoundNorms} that, for every $\vartheta \in [-\vartheta_\ast, \vartheta_\ast]$, if $z_0$ is a root of $\widehat\Delta(\cdot, \vartheta)$ with $\Real z_0 \geq 0$, then $\abs{z_0} \leq M + \abs{e^{-z_0}} M \leq 2 M = R$, yielding the conclusion.
\end{proof}

\begin{lemma}
\label{lemm:unique-roots-on-imaginary-axis}
Let $\widehat\Delta$ be the function defined by \eqref{eq:complex-widehat-Delta} and \eqref{eq:complex-widehat-Delta-0} and $\vartheta \in \mathbb R$. Then $i\vartheta$ and $-i\vartheta$ are the unique roots of $\widehat\Delta(\cdot, \vartheta)$ on the imaginary axis.
\end{lemma}

\begin{proof}
For $\vartheta = 0$, this result is a consequence of Theorem~\ref{thm:main-MID}, so we are only left to treat the case $\vartheta \neq 0$.

Let $z$ be a root of $\widehat\Delta(\cdot, \vartheta)$ on the imaginary axis and define $\omega = \Imag z$. Then $\widehat\Delta(i\omega, \vartheta) = 0$, which means, using \eqref{eq:complex-widehat-Delta}, that
\begin{equation*}
-\omega^2 - 2 \vartheta \frac{\vartheta - \sin\vartheta \cos\vartheta}{\vartheta^2 - \sin^2 \vartheta} i \omega + \vartheta^2 \frac{\vartheta^2 + \sin^2 \vartheta}{\vartheta^2 - \sin^2 \vartheta} = - e^{-i\omega} \Biggl(2 \vartheta \frac{\vartheta \cos \vartheta - \sin\vartheta}{\vartheta^2 - \sin^2 \vartheta} i\omega - \frac{2 \vartheta^3 \sin \vartheta}{\vartheta^2 - \sin^2 \vartheta}\Biggr).
\end{equation*}
Multiplying both sides by $\vartheta^2 - \sin^2\vartheta$ and then taking the square of their absolute values, we obtain
\begin{multline*}
\bigl(\vartheta^2 (\vartheta^2 + \sin^2 \vartheta) - \omega^2 (\vartheta^2 - \sin^2 \vartheta)\bigr)^2 + 4 \vartheta^2 \omega^2 (\vartheta - \sin\vartheta \cos\vartheta)^2 \\ = 4 \vartheta^2 \omega^2 (\vartheta \cos \vartheta - \sin\vartheta)^2 + 4 \vartheta^6 \sin^2\vartheta.
\end{multline*}
Expanding the terms in the above equality, one arrives, after some tedious but straightforward computations, that it is equivalent to $(\omega^2 - \vartheta^2)^2 = 0$, and thus the only possible values for $\omega$ are $\vartheta$ and $-\vartheta$, as required.
\end{proof}

\subsection{Proof of Theorem~\ref{thm:complex-dominance}}
\label{sec:MID-complex-proof}

Thanks to the preliminary remarks and auxiliary results from Section~\ref{sec:MID-complex-auxiliary}, we are now in position to prove Theorem~\ref{thm:complex-dominance}.

\begin{proof}[Proof of Theorem~\ref{thm:complex-dominance}]
As stated in Section~\ref{sec:MID-complex-auxiliary}, to prove the first part of the theorem, it suffices to show that, given $\vartheta \in \mathbb R^\ast$, a quasipolynomial $\Delta$ of the form \eqref{eq:Delta-DDE} with $n = 2$, $m = 1$, and $\tau = 1$ admits $i \vartheta$ and $-i \vartheta$ as roots of multiplicity $2$ each if and only if $\Delta$ coincides with the quasipolynomial $\widehat\Delta(\cdot, \vartheta)$ from \eqref{eq:complex-widehat-Delta}.

Consider thus a quasipolynomial $\Delta$ given by
\[
\Delta(z) = z^2 + a_1 z + a_0 + e^{-z} (\alpha_1 z + \alpha_0), \qquad z \in \mathbb C,
\]
for some real numbers $a_0, a_1, \alpha_0, \alpha_1$. Note that, since $a_0, a_1, \alpha_0, \alpha_1$ are assumed to be real coefficients, it suffices to consider the root $i \vartheta$, since $i \vartheta$ is a root of multiplicity $2$ of $\Delta$ if and only if its complex conjugate $-i\vartheta$ is a root of multiplicity $2$ of $\Delta$.

The multiplicity of $i\vartheta$ as a root of $\Delta$ is at least $2$ if and only if $\Delta(i \vartheta) = \Delta^\prime(i \vartheta) = 0$, that is, if and only if the coefficients $a_0, a_1, \alpha_0, \alpha_1$ satisfy the linear system
\[
\left\{
\begin{aligned}
-\vartheta^2 + i a_1 \vartheta + a_0 + e^{-i\vartheta} (i\alpha_1 \vartheta + \alpha_0) & = 0, \\
2 i \vartheta + a_1 + e^{-i\vartheta} (- i\alpha_1 \vartheta - \alpha_0 + \alpha_1) & = 0,
\end{aligned}
\right.
\]
which can be rewritten, separating real and imaginary parts, as the system
\[
\left\{
\begin{aligned}
-\vartheta^2 + a_0 + \alpha_0 \cos\vartheta + \alpha_1 \vartheta \sin\vartheta & = 0, \\
a_1 \vartheta + \alpha_1 \vartheta\cos\vartheta - \alpha_0 \sin\vartheta & = 0, \\
a_1 + (\alpha_1 - \alpha_0)\cos\vartheta - \alpha_1 \vartheta \sin\vartheta & = 0, \\
2\vartheta + (\alpha_0 - \alpha_1) \sin\vartheta - \alpha_1 \vartheta \cos\vartheta & = 0.
\end{aligned}
\right.
\]
Solving the above system with respect to $a_1, a_0, \alpha_1, \alpha_0$, we obtain that its unique solution is given by
\begin{subequations}
\label{ExpressionsBBeta}
\begin{align}
a_1 & = - 2 \vartheta \frac{\vartheta - \sin\vartheta \cos\vartheta}{\vartheta^2 - \sin^2 \vartheta}, \displaybreak[0] \\
a_0 & = \vartheta^2 \frac{\vartheta^2 + \sin^2 \vartheta}{\vartheta^2 - \sin^2 \vartheta}, \displaybreak[0] \\
\alpha_1 & = 2 \vartheta \frac{\vartheta \cos\vartheta - \sin\vartheta}{\vartheta^2 - \sin^2 \vartheta}, \displaybreak[0] \\
\alpha_0 & = - \frac{2 \vartheta^3 \sin\vartheta}{\vartheta^2 - \sin^2 \vartheta},
\end{align}
\end{subequations}
which are precisely the coefficients of \eqref{eq:complex-widehat-Delta}.

We now verify that, under \eqref{ExpressionsBBeta}, $i\vartheta$ is a root of $\Delta$ of multiplicity exactly $2$, i.e., that $\Delta^{\prime\prime}(i\vartheta) \neq 0$. Indeed, since
\[
\Delta^{\prime\prime}(z) = 2 + e^{-z} (\alpha_1 z + \alpha_0 - 2 \alpha_1),
\]
one computes
\begin{align*}
\Real \Delta^{\prime\prime}(i \vartheta) & = 2 + (\alpha_0 - 2 \alpha_1) \cos\vartheta + \alpha_1 \vartheta \sin\vartheta, \displaybreak[0] \\
\Imag \Delta^{\prime\prime}(i \vartheta) & = \alpha_1 \vartheta \cos\vartheta - (\alpha_0 - 2 \alpha_1) \sin \vartheta.
\end{align*}
In particular, under \eqref{ExpressionsBBeta}, one has
\begin{align*}
\Imag \Delta^{\prime\prime}(i \vartheta) & = \alpha_1 (\vartheta\cos\vartheta + 2 \sin\vartheta) - \alpha_0 \sin\vartheta \displaybreak[0] \\
& = \frac{2 \vartheta(\vartheta^2 \cos^2\vartheta + \vartheta \sin\vartheta \cos\vartheta - 2 \sin^2\vartheta) + 2 \vartheta^3 \sin^2\vartheta}{\vartheta^2 - \sin^2\vartheta} \displaybreak[0] \\
& = 2 \vartheta\frac{\vartheta^2 + \vartheta \sin\vartheta \cos\vartheta - 2 \sin^2\vartheta}{\vartheta^2 - \sin^2\vartheta}.
\end{align*}
Since $\vartheta \neq 0$, it follows from Lemma~\ref{lemm:F-positive} that $\Imag \Delta^{\prime\prime}(i\vartheta) \neq 0$, which shows that, under \eqref{ExpressionsBBeta}, $\vartheta$ is a root of $\Delta$ of multiplicity exactly $2$, yielding the desired conclusion.

Let us now prove the last part of the statement of the theorem, which amounts to showing that, for every $\vartheta \in \mathbb R^\ast$ and $z \in \mathbb C \setminus \{i\vartheta, -i\vartheta\}$ such that $\widehat\Delta(z, \vartheta) = 0$, we have $\Real z < 0$, where $\widehat\Delta$ is given by \eqref{eq:complex-widehat-Delta}. Note that, since $\widehat\Delta(\cdot, \vartheta) = \widehat\Delta(\cdot, -\vartheta)$, it suffices to prove the previous statement for $\vartheta \in \mathbb R_+^\ast$. Assume, to obtain a contradiction, that the set
\begin{equation}
\label{eq:empty-set}
\left\{\vartheta \in \mathbb R_+^\ast \suchthat \text{there exists } z \in \mathbb C \setminus \{i\vartheta, -i\vartheta\} \text{ such that } \Real z \geq 0 \text{ and } \widehat\Delta(z, \vartheta) = 0\right\}
\end{equation}
is nonempty, and denote by $\vartheta_\ast \in \mathbb R_+$ its infimum. Then there exists a sequence $(\vartheta_n)_{n \in \mathbb N}$ in $\mathbb R_+^\ast$ and a sequence $(z_n)_{n \in \mathbb N}$ in $\mathbb C$ such that $\vartheta_n \to \vartheta_\ast$ as $n \to +\infty$ and, for every $n \in \mathbb N$, we have $z_n \notin \{i\vartheta_n, -i \vartheta_n\}$, $\Real z_n \geq 0$, and $\widehat\Delta(z_n, \vartheta_n) = 0$. Up to replacing $z_n$ by $\overline z_n$, we assume, with no loss of generality, that $\Imag z_n \geq 0$.

Since the sequence $(\vartheta_n)_{n \in \mathbb N}$ converges, it is bounded, and thus, by Lemma~\ref{lemm:complex-preliminary-bound}, the sequence $(z_n)_{n \in \mathbb N}$ is also bounded. Hence, up to extracting a subsequence (which we still denote by the same notation for simplicity), there exists $z_\ast \in \mathbb C$ such that $z_n \to z_\ast$ as $n \to +\infty$. We thus have $\Real z_\ast \geq 0$, $\Imag z_\ast \geq 0$, and, by continuity of $\widehat\Delta$, we also have $\widehat\Delta(z_\ast, \vartheta_\ast) = 0$. We split the sequel of the contradiction argument in two cases, according to whether $\Real z_\ast > 0$ or $\Real z_\ast = 0$.

Assume first that $\Real z_\ast > 0$. In this case, we have $\vartheta_\ast > 0$, since, as recalled in Section~\ref{sec:MID-complex-auxiliary}, $\widehat\Delta(\cdot, 0)$ admits $0$ as a strictly dominant root of multiplicity $4$, as a consequence of Theorem~\ref{thm:main-MID}, and hence it cannot admit a root $z_\ast$ with $\Real z_\ast > 0$. Now, for $n \in \mathbb N$, let $\widehat\vartheta_n = \left(1 - \frac{1}{n+2}\right)\vartheta_\ast$, so that $0 < \widehat\vartheta_n < \vartheta_\ast$ and $\widehat\vartheta_n \to \vartheta_\ast$ as $n \to +\infty$. Since $z_\ast$ is an isolated zero of $\widehat\Delta(\cdot, \vartheta_\ast)$ and $\widehat\Delta(\cdot, \widehat\vartheta_n)$ converges to $\widehat\Delta(\cdot, \vartheta_\ast)$ as $n \to +\infty$ uniformly on compact subsets of $\mathbb C$, we obtain from Theorem~\ref{thm:Hurwitz} that there exists $N \in \mathbb N$ such that, for every $n \in \mathbb N$ with $n \geq N$, there exists $\widehat z_n \in \mathbb C$ with $\Real \widehat z_n > 0$ and $\widehat\Delta(\widehat z_n, \widehat\vartheta_n) = 0$. Thus, $\widehat\vartheta_n$ belongs to the set in \eqref{eq:empty-set}, which contradicts the fact that $\vartheta_\ast$ is the infimum of that set.

Let us now consider the case $\Real z_\ast = 0$. By Lemma~\ref{lemm:unique-roots-on-imaginary-axis}, one deduces that $z_\ast = i \vartheta_\ast$. Let $\varepsilon > 0$ be such that $i \vartheta_\ast$ is the unique root of $\widehat\Delta(\cdot, \vartheta_\ast)$ in $\overline B_{\mathbb C}(i \vartheta_\ast, \varepsilon)$, which exists since all zeros of $\widehat\Delta(\cdot, \vartheta_\ast)$ are isolated. Define $M = 2$ if $\vartheta_\ast > 0$ and $M = 4$ if $\vartheta_\ast = 0$. Since $i \vartheta_\ast$ is a root of multiplicity $M$ of $\widehat\Delta(\cdot, \vartheta_\ast)$ and $\widehat\Delta(\cdot, \vartheta_n)$ converges to $\widehat\Delta(\cdot, \vartheta_\ast)$ as $n \to +\infty$ uniformly on compact subsets of $\mathbb C$, it follows from Theorem~\ref{thm:Hurwitz} that there exists $N \in \mathbb N$ such that, for every $n \in \mathbb N$ with $n \geq N$, $\widehat\Delta(\cdot, \vartheta_n)$ admits exactly $M$ roots (counted with their multiplicities) in $B_{\mathbb C}(i \vartheta_\ast, \varepsilon)$. In the case $\vartheta_\ast > 0$, we have $M = 2$ but there are at least three roots of $\widehat\Delta(\cdot, \vartheta_n)$ (counted with their multiplicities) in $B_{\mathbb C}(i \vartheta_\ast, \varepsilon)$ for $n$ large enough: the double root $i \vartheta_n$ and the root $z_n$, yielding the desired contradiction. In the case $\vartheta_\ast = 0$, we have $M = 4$ but there are at least five roots of $\widehat\Delta(\cdot, \vartheta_n)$ (counted with their multiplicities) in $B_{\mathbb C}(i \vartheta_\ast, \varepsilon)$ for $n$ large enough: the double roots $i \vartheta_n$ and $-i\vartheta_n$ and the root $z_n$, yielding once again the desired contradiction.

We have thus shown by contradiction that the set in \eqref{eq:empty-set} is empty, which concludes the proof.
\end{proof}

\subsection{Numerical computation of the roots of \texorpdfstring{$\widehat\Delta$}{\^{}Delta}}
\label{sec:numerical-widehat-delta}

As discussed in Section~\ref{sec:MID-complex-auxiliary}, any quasipolynomial of the form \eqref{eq:Delta-DDE} with $n = 2$ and $m = 1$ admitting a pair of complex conjugate roots of multiplicity $2$ each (or, equivalently, any quasipolynomial of the form \eqref{eq:Delta-DDE} with $n = 2$ and $m = 1$ whose coefficients are given by \eqref{eq:coeffs-complex-dominance}) can be put into the form \eqref{eq:complex-widehat-Delta} up to a translation and a rescaling of the complex plane defined by \eqref{eq:complex-change-of-vars}. Theorem~\ref{thm:complex-dominance} guarantees that the multiple roots $i \vartheta$ and $-i\vartheta$ of $\widehat\Delta(\cdot, \vartheta)$ are strictly dominant, but says nothing about how the roots on the open left half-plane behave. In order to get a grasp on their behavior, we have performed numerical computations\footnote{Numerical computations were performed using Python's \texttt{cxroots} module \cite{Parini2018cxroots}, which implements numerical methods described in \cite{Kravanja2000Computing}.} of all roots of $\widehat\Delta(\cdot, \vartheta)$ on the region $\{s \in \mathbb C \suchthat -4.75 \leq \Real s \leq 0.25 \text{ and } -25 \leq \Imag s \leq 25\}$ for several values of $\vartheta \in [0, 8]$. The results are provided in Figure~\ref{FigRoots}, with different values of $\vartheta$ being represented with different colors.

\begin{figure}[ht]
\centering
\resizebox{0.667\textwidth}{!}{\input{Figures/complex-roots.pgf}}
\caption[{Roots of the quasipolynomial $\widehat\Delta(\cdot, \vartheta)$ from \eqref{eq:complex-widehat-Delta} for $\vartheta \in [0, 8]$}]{Roots of the quasipolynomial $\widehat\Delta(\cdot, \vartheta)$ from \eqref{eq:complex-widehat-Delta} for $\vartheta \in [0, 8]$, with a detailed view of the region $\{s \in \mathbb C \suchthat -1.80 \leq \Real s \leq -1.55 \text{ and } 10 \leq \Imag s \leq 11\}$}
\label{FigRoots}
\end{figure}

One observes in Figure~\ref{FigRoots} the movement of the dominant roots $i \vartheta$ and $-i\vartheta$ along the imaginary axis. Concerning the other roots, as $\vartheta$ increases, the imaginary parts of the nonreal roots represented in the figure increase in absolute value, while the real parts oscillate. Table~\ref{TabFirstRoot} synthesizes this oscillatory behavior for the first pair of nondominant complex conjugate roots in Figure~\ref{FigRoots} (in order of increasing absolute value of imaginary part) by presenting the values of $\vartheta$ and the corresponding roots at local extrema of their real part, as well as for the initial and final values $\vartheta = 0$ and $\vartheta = 8$ used in the numerical computations.

\begin{table}[ht]
\centering
\caption[First pair of nondominant complex conjugate roots in Figure~\ref{FigRoots}]{First pair of nondominant complex conjugate roots in Figure~\ref{FigRoots} at their initial and final positions and on local extrema of their real parts}
\label{TabFirstRoot}
\begin{tabular}{ccc}
\toprule
 & $\vartheta$ & Roots \\
\midrule
Initial & $0$ & $-1.731 \pm 10.16 i$ \\
First local maximum & $2.51$ & $-1.586 \pm 10.46i$ \\
First local minimum & $4.59$ & $-2.735 \pm 12.14i$ \\
Second local maximum & $6.19$ & $-1.764 \pm 13.74i$ \\
Second local minimum & $7.83$ & $-2.508 \pm 15.32i$ \\
Final & $8$ & $-2.466 \pm 15.66i$ \\
\bottomrule
\end{tabular}
\end{table}

\begin{figure}[ht]
\centering
\resizebox{\textwidth}{!}{\input{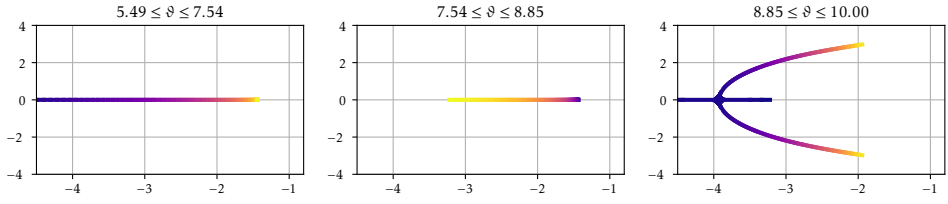}}
\caption[Detailed view of the roots of $\widehat\Delta(\cdot, \vartheta)$ around the real axis for different ranges of $\vartheta$]{Detailed view of the roots of $\widehat\Delta(\cdot, \vartheta)$ around the real axis for different ranges of $\vartheta$. In each figure, darker colors correspond to smaller values of $\vartheta$}
\label{FigRealRoots}
\end{figure}

We also notice in Figure~\ref{FigRoots} the presence of a real-valued root. Its detailed behavior obtained from numerical computations for $\vartheta \in [5.49, 10.00]$ is provided in Figure~\ref{FigRealRoots}, which is split in three different ranges for $\vartheta$ corresponding to different observed behaviors of the root. This root first appears in the domain under consideration for $\vartheta \approx 5.49$ and moves to the right, reaching a local maximum at $\vartheta \approx 7.54$, at which point its value is approximately $-1.437$. It then starts moving to the left for $\vartheta \in [7.54, 8.85]$. At $\vartheta \approx 8.85$, a second real root appears in the domain under consideration, moving to the right, and both roots meet, given rise, when $\vartheta \approx 8.88$, to a real root of multiplicity $2$ whose value is approximately $-3.927$. For $\vartheta \in [8.88, 10]$, these roots become a pair of complex conjugate roots which start moving to the right. As $\vartheta$ increases beyond $10$ (not represented in Figure~\ref{FigRealRoots}), one observes that this pair of roots oscillates like the other pairs of complex conjugate roots from Figure~\ref{FigRoots}.

\section{Systems with more than one delay}
\label{sec:MID-two-delays}

In Sections~\ref{sec:MID-maximal}, \ref{sec:MID-intermediate}, and \ref{sec:MID-complex}, all the presented results concern characteristic equations of delay-differential equations of the form \eqref{eq:DDE}, with a single delay $\tau$. A natural question is whether those results can be generalized to systems with more than one delay. A general answer to this questions is still missing in the literature, but some important insights were obtained for a particular first-order system with two delays in \cite{Fueyo2021Insights, Fueyo2023Pole}, whose results we now report.

\subsection{Statement of the main result}

We shall consider in the sequel the delay-differential equation
\begin{equation}
\label{eq:main-syst-two-delays}
y^\prime(t) + a_0 y(t) + a_1 y(t - \tau_1) + a_2 y(t - \tau_2) = 0,
\end{equation}
where $a_0, a_1, a_2$ are real coefficients and $\tau_1$ and $\tau_2$ are positive real numbers representing delays. Its characteristic quasipolynomial is the function $\Delta\colon \mathbb C \to \mathbb C$ defined for $s \in \mathbb C$ by
\begin{equation}
\label{eq:main-quasipoly-two-delays}
\Delta(s) = s + a_0 + a_1 e^{-s \tau_1} + a_2 e^{-s \tau_2}.
\end{equation}
Note that the degree of $\Delta$ is $3$, and thus, by Corollary~\ref{coro:Polya-Szego}, the maximal multiplicity that a root $s_0$ of $\Delta$ may have is $3$. The main result we will prove here is the following, characterizing the unique choice of coefficients $(a_0, a_1, a_2) \in \mathbb R^3$ such that a given real number $s_0$ is a root of multiplicity $3$ of $\Delta$ for given delays $\tau_1, \tau_2$ with $\tau_1 \neq \tau_2$, and proving that the MID property is satisfied in this configuration.

\begin{theorem}
\label{thm:main-two-delays}
Let $\tau_1 > 0$ and $\tau_2 > 0$ be such that $\tau_1 \neq \tau_2$, $a_0, a_1, a_2$ be real numbers, and consider the quasipolynomial $\Delta$ from \eqref{eq:main-quasipoly-two-delays}. Let $s_0 \in \mathbb R$. Then $s_0$ is a root of multiplicity $3$ of $\Delta$ if and only if
\begin{equation}
\label{eq:two-delays-coeffs}
a_0 = -s_0 - \frac{1}{\tau_1} - \frac{1}{\tau_2}, \qquad a_1 = \frac{\tau_2}{\tau_1 (\tau_2-\tau_1)}e^{s_0 \tau_1}, \qquad a_2 = -\frac{\tau_1}{\tau_2(\tau_2-\tau_1)}e^{s_0 \tau_2}.
\end{equation}
In addition, if \eqref{eq:two-delays-coeffs} is satisfied, then $\Real s < s_0$ for every $s \in \mathbb C \setminus \{s_0\}$ such that $\Delta(s) = 0$.
\end{theorem}

Similarly to Theorems~\ref{thm:main-MID} and \ref{thm:complex-dominance}, the first part of Theorem~\ref{thm:main-two-delays} consists in the characterization of the existence of a root $s_0$ of multiplicity $3$, while the second part asserts its dominance. We will prove Theorem~\ref{thm:main-two-delays} through a series of auxiliary results split into several sections. We first perform, in Section~\ref{sec:normalization}, a translation and a rescaling of the complex plane in order to normalize the root of multiplicity $3$ to $0$ and the largest delay to $1$. This allows us to provide a simple proof of the first part of Theorem~\ref{thm:main-two-delays}, which is the subject of Section~\ref{sec:conditions}. In order to prove the second part of Theorem~\ref{thm:main-two-delays}, on the strict dominance of the root at $0$, we need to study in details how the roots of $\Delta$ change as the ratio $\frac{\tau_1}{\tau_2}$ changes, and this is the subject of Section~\ref{sec:branches}. Finally, Theorem~\ref{thm:main-two-delays} is proved in Section~\ref{sec:proof-two-delays}, and some numerical experiments are provided in Sections~\ref{sec:roots-Q} and \ref{sec:three-delays}.

\subsection{Normalization}
\label{sec:normalization}

Our next result provides a normalization of the quasipolynomial $\Delta$ from \eqref{eq:main-quasipoly-two-delays}. Its proof can be obtained by straightforward computations and is omitted here.

\begin{lemma}
\label{lemm:transformation}
Let $\tau_1$ and $\tau_2$ be such that $0 < \tau_1 < \tau_2$, $a_0, a_1, a_2$ be real numbers, and consider the quasipolynomial $\Delta$ from \eqref{eq:main-quasipoly-two-delays}. Let $\lambda = \frac{\tau_1}{\tau_2} \in (0, 1)$, $s_0 \in \mathbb R$, and define the quasipolynomial $Q\colon \mathbb C \to \mathbb C$ by setting $Q(z) = \tau_2 \Delta(s_0 + \frac{z}{\tau_2})$ for $z\in \mathbb{C}$. Then
\begin{enumerate}
\item the coefficients $a_0, a_1, a_2$ of $\Delta$ are given by \eqref{eq:two-delays-coeffs} if and only if $Q$ is given by
\begin{equation}
\label{eq:def-Q}
Q(z) = z - \frac{\lambda + 1}{\lambda} + \frac{1}{\lambda (1 - \lambda)} e^{-\lambda z} - \frac{\lambda}{(1 - \lambda)} e^{-z}, \qquad \text{ for } z \in \mathbb C;
\end{equation}
\item $s_0$ is a root of multiplicity $3$ of $\Delta$ if and only if $0$ is a root of multiplicity $3$ of $Q$; and
\item $\Real s < s_0$ for every $s \in \mathbb C \setminus \{s_0\}$ such that $\Delta(s) = 0$ if and only if $\Real z < 0$ for every $z \in \mathbb C^\ast$ such that $Q(z) = 0$.
\end{enumerate}
\end{lemma}

\subsection{Conditions for maximal multiplicity}
\label{sec:conditions}

As a consequence of Lemma~\ref{lemm:transformation}, we can provide a proof of the first part of the statement of Theorem~\ref{thm:main-two-delays}.

\begin{proposition}
\label{prop:main-two-delays-conditions}
Let $\tau_1 > 0$ and $\tau_2 > 0$ be such that $\tau_1 \neq \tau_2$, $a_1, a_2, a_3$ be real numbers, and consider the quasipolynomial $\Delta$ from \eqref{eq:main-quasipoly-two-delays}. Let $s_0 \in \mathbb R$. Then $s_0$ is a root of multiplicity $3$ of $\Delta$ if and only if \eqref{eq:two-delays-coeffs} is satisfied.
\end{proposition}

\begin{proof}
Without loss of generality, we assume that $0 < \tau_1 < \tau_2$. From Lemma~\ref{lemm:transformation}, it suffices to show that a quasipolynomial of the form \eqref{eq:main-quasipoly-two-delays} with $\tau_1 = \lambda \in (0, 1)$ and $\tau_2 = 1$ admits $0$ as a root of multiplicity $3$ if and only if it is the quasipolynomial $Q$ defined in \eqref{eq:def-Q}.

Let $\widetilde a_0, \widetilde a_1, \widetilde a_2$ be real numbers and consider the quasipolynomial $\widetilde\Delta$ defined by
\[
\widetilde\Delta(z) = z + \widetilde a_0 + \widetilde a_1 e^{-\lambda z} + \widetilde a_2 e^{-z}.
\]
By Corollary~\ref{coro:Polya-Szego}, the maximal multiplicity of any root of $\widetilde\Delta$ is $3$, and thus $0$ is a root of multiplicity $3$ of $\widetilde\Delta$ if and only if $\widetilde\Delta(0) = \widetilde\Delta^\prime(0) = \widetilde\Delta^{\prime\prime}(0) = 0$. Hence, $0$ is a root of multiplicity $3$ of $\widetilde\Delta$ if and only if $\widetilde{a}_0$, $\widetilde{a}_1$, and $\widetilde{a}_2$ satisfy the linear system
\begin{equation*}
\begin{pmatrix}
1 & 1 & 1 \\
0 & -\lambda & -1 \\
0  & \lambda^2 & 1
\end{pmatrix} \begin{pmatrix}
\widetilde{a}_0 \\
\widetilde{a}_1 \\
\widetilde{a}_2
\end{pmatrix} = \begin{pmatrix}
0 \\
-1 \\
0 
\end{pmatrix}.
\end{equation*}
One can immediately verify that the above system admits a unique solution, given by
\begin{equation*}
\widetilde{a}_0=-\frac{\lambda+1}{\lambda},\qquad \widetilde{a}_1=\frac{1}{\lambda(1 - \lambda)}, \qquad \widetilde{a}_2 = -\frac{\lambda}{1 - \lambda},
\end{equation*}
and thus $0$ is a root of multiplicity $3$ of $\widetilde\Delta$ if and only if $\widetilde\Delta = Q$, as required.
\end{proof}

\subsection{Branches of roots}
\label{sec:branches}

The strategy of the proof for the second part of Theorem~\ref{thm:main-two-delays} that we use here bears several similarities with the proof of the corresponding part of Theorem~\ref{thm:complex-dominance}. We first remark that, thanks to Lemma~\ref{lemm:transformation}, the second part of the statement of Theorem~\ref{thm:main-two-delays} is proved if we show that, for every $\lambda \in (0, 1)$, the quasipolynomial $Q$ defined in \eqref{eq:def-Q} satisfies that $\Real z < 0$ for every $z \in \mathbb C^\ast$ such that $Q(z) = 0$. We then regard $Q$ as a function of two variables, $z$ and $\lambda$, i.e., we consider that $Q\colon \mathbb C \times (0, 1) \to \mathbb C$ is the function defined for $(z, \lambda) \in \mathbb C \times (0, 1)$ by
\begin{equation}
\label{eq:def-Q-lambda}
Q(z, \lambda) = z - \frac{\lambda + 1}{\lambda} + \frac{1}{\lambda (1 - \lambda)} e^{-\lambda z} - \frac{\lambda}{1 - \lambda} e^{-z}.
\end{equation}
As we did in Section~\ref{sec:MID-complex-proof} for the proof of Theorem~\ref{thm:complex-dominance}, our strategy will then consist in studying how the roots of the quasipolynomial $Q(\cdot, \lambda)$ behave as $\lambda$ changes. We will first prove that $0$ is a strictly dominant root of $Q(\cdot, \lambda)$ for $\lambda$ close enough to $1$, by studying the behavior of $Q(\cdot, \lambda)$ as $\lambda \to 1^-$. The result is then extended to any $\lambda \in (0, 1)$ by a contradiction argument.

Contrarily to Section~\ref{sec:MID-complex-proof}, the behavior of the roots of $Q(\cdot, \lambda)$ is more intricate than that of the roots of $\widehat\Delta(\cdot, \vartheta)$ from \eqref{eq:complex-widehat-Delta}, and we will need to study branches of roots of $Q(\cdot, \lambda)$, whose definition we now provide.

\begin{definition}
Let $Q$ be the function defined by \eqref{eq:def-Q-lambda}. A \emph{branch of roots} of $Q$ is a function $z\colon I\to\mathbb C$ such that $I \subset (0, 1)$ is an open interval and $Q(z(\lambda), \lambda) = 0$ for every $\lambda \in I$. If $z\colon (\underline\lambda, \overline\lambda)\to\mathbb C$ is a branch of roots of $Q$, an \emph{extension} $\widehat z$ of $z$ is a branch of roots of $Q$ defined on $(\underline\lambda, \widehat\lambda)$ for some $\widehat\lambda \in [\overline\lambda, 1]$ and such that $z(\lambda) = \widehat z(\lambda)$ for every $\lambda \in (\underline\lambda, \overline\lambda)$.
\end{definition}

Note that, here, we only consider extensions of branches to the right of the domain of definition, since these are the only kinds of extensions we shall need. When there is no ambiguity, we will also denote extensions of $z$ by the same letter $z$. We will also do the slight abuse of notation here of assuming that an extension always has the same regularity as the original branch: for instance, when we say ``an extension of an analytic branch of roots $z$'', we implicitly assume that the extension is also an analytic function of $\lambda$.

We start our auxiliary results by showing that any nontrivial root of $Q(\cdot, \lambda)$ in the complex right half-plane is necessarily simple and hence gives rise to an analytic branch of roots of $Q$.

\begin{lemma}
\label{lem:simple-root}
Let $Q$ be the function defined by \eqref{eq:def-Q-lambda}, $\lambda_0 \in (0, 1)$, and assume that there exists $z_0 \in \mathbb C^\ast$ such that $Q(z_0, \lambda_0) = 0$ and $\Real(z_0) \geq 0$. Then $z_0$ is a simple root of $Q(\cdot, \lambda_0)$, i.e., $\partial_z Q(z_0, \lambda_0) \neq 0$, and there exists an analytic branch $z\colon I\to\mathbb C$ of roots of $Q$ such that $\lambda_0\in I$ and $z(\lambda_0) = z_0$. Moreover, $z$ satisfies
\begin{equation}
\label{eq:derivative}
z^\prime(\lambda_0) = -\frac{(1 - 2 \lambda_0) z(\lambda_0) + 2 \lambda_0 - 2\lambda_0 e^{-z(\lambda_0)} - z(\lambda_0) e^{-\lambda_0 z(\lambda_0)}}{\lambda_0 (1 - \lambda_0) + \lambda_0^2 e^{-z(\lambda_0)} - \lambda_0 e^{-\lambda_0 z(\lambda_0)}}.
\end{equation}
\end{lemma}

\begin{proof}
Consider the function $\widetilde Q\colon \mathbb C \times (0, 1) \to \mathbb C$ defined for $(z, \lambda) \in \mathbb C \times (0, 1)$ by
\begin{equation*}
\widetilde{Q}(z,\lambda) = \lambda(1-\lambda)Q(z,\lambda) = \lambda(1-\lambda)z-1+\lambda^2-\lambda^2e^{-z}+e^{-\lambda z}.
\end{equation*}
For every $(z, \lambda) \in \mathbb C \times (0, 1)$, we have
\begin{subequations}
\begin{align}
\label{eq:lem_dp4}
\partial_z \widetilde{Q}(z,\lambda) & =\lambda(1-\lambda)+\lambda^2e^{-z}-\lambda e^{-\lambda z}, \\
\label{eq:lem_dp8}
\partial_{\lambda} \widetilde{Q}(z,\lambda) & = (1-2\lambda)z+2\lambda-2\lambda e^{-z}-z e^{-\lambda z}.
\end{align}
\end{subequations}
Since $\widetilde Q(z_0, \lambda_0) = 0$, we deduce that
\begin{equation}
\label{eq:lem_branch2}
\lambda_0^2 e^{-z_0} = \lambda_0 (1 - \lambda_0) z_0 - 1 + \lambda_0^2 + e^{-\lambda_0 z_0}.
\end{equation}
Inserting~\eqref{eq:lem_branch2} into \eqref{eq:lem_dp4}, one obtains
\begin{equation*}
\partial_z Q(z_0,\lambda_0)=- (1 - \lambda_0) (1-\lambda_0z_0-e^{-\lambda_0 z_0}).
\end{equation*}
In particular, we have $\Imag(\partial_z \widetilde Q(z_0, \lambda_0)) = (1 - \lambda_0)\allowbreak \left(\lambda_0 y_0 - e^{-\lambda_0 x_0} \sin(\lambda_0 y_0)\right)$, where $x_0 = \Real(z_0)$ and $y_0 = \Imag(z_0)$. Since $x_0 \geq 0$ and $(x_0, y_0) \neq (0, 0)$, we deduce that $\Imag(\partial_z \widetilde Q(z_0, \lambda_0)) \neq 0$ thanks to the classical properties of the sinc function, yielding that $z_0$ is a simple root of $Q(\cdot, \lambda_0)$. The existence of the branch of roots $z$ asserted in the statement is thus an immediate consequence of the simplicity of the nontrivial root $z$ and the implicit function theorem for analytic functions (see, e.g., \cite[Theorem~2.1.2]{Hormander1990Introduction}), and \eqref{eq:derivative} follows immediately from \eqref{eq:lem_dp4} and \eqref{eq:lem_dp8}.
\end{proof}

We next provide some properties of analytic branches of roots. We first notice that, as an immediate consequence of Hurwitz's theorem, since $0$ is a triple root of $Q(\cdot, \lambda)$ for every $\lambda \in (0, 1)$, no other branch of roots of $Q$ can pass through the origin.

\begin{lemma}
\label{lem:branch-not-zero}
Let $Q$ be the function defined by \eqref{eq:def-Q-lambda} and $z\colon I\to\mathbb C$ be a continuous branch of roots of $Q$. Assume that there exists $\lambda_0 \in I$ such that $z(\lambda_0) \neq 0$. Then $z(\lambda) \neq 0$ for every $\lambda \in I$.
\end{lemma}

\begin{proof}
Assume, to obtain a contradiction, that the set $\{\lambda \in I \suchthat z(\lambda) = 0\}$ is nonempty. Since $z$ is not constantly equal to $0$ over the interval $I$, there exists a sequence $(\lambda_n)_{n \in \mathbb N}$ in $I$ converging to some $\lambda_\ast \in I$ with $z(\lambda_n) \neq 0$ for every $n \in \mathbb N$ and $z(\lambda_\ast) = 0$. Since $0$ is a root of multiplicity three of $Q(\cdot, \lambda_\ast)$, by Theorem~\ref{thm:Hurwitz}, there exist neighborhoods $U \subset \mathbb C$ of $0$ and $V \subset (0, 1)$ of $\lambda_\ast$ such that, for every $\lambda \in V$, $Q(\cdot, \lambda)$ admits exactly three roots (counted with their multiplicity) in $U$. As $0$ is a triple root of $Q(\cdot, \lambda)$ for every $\lambda \in (0, 1)$, we then deduce that, for $\lambda \in V$, $0$ is the only root of $Q(\cdot, \lambda)$ in $U$. However, since $\lambda_n \to \lambda_\ast$ as $n \to +\infty$, we deduce that $\lambda_n \in V$ for $n$ large enough. Since $z$ is continuous, we also have $z(\lambda_n) \in U$ for $n$ large enough. We thus obtain the desired contradiction since $z(\lambda_n) \neq 0$ is a root of $Q(\cdot, \lambda_n)$.
\end{proof}

We now exploit the consequences of the P\'olya--Szeg\H{o} bound stated in Proposition~\ref{prop:Polya-Szego} to deduce that continuous branches of nontrivial roots never cross the real axis.

\begin{lemma}
\label{lem:cannot-cross-real-axis}
Let $Q$ be the function defined by \eqref{eq:def-Q-lambda} and $z\colon I \to \mathbb C$ be a continuous branch of nontrivial roots of $Q$. Then $\Imag(z(\lambda)) \neq 0$ for every $\lambda \in I$ and, in particular, the sign of $\Imag(z(\lambda))$ is constant.
\end{lemma}

\begin{proof}
Since $0$ is a root of multiplicity $3$ of $Q(\cdot, \lambda)$ and $Q(\cdot, \lambda)$ is a quasipolynomial of degree $3$, it follows from Proposition~\ref{prop:Polya-Szego} that $0$ is the unique real root of $Q(\cdot, \lambda)$, yielding the first part of the conclusion. The second part is a consequence of the continuity of $\lambda \mapsto \Imag(z(\lambda))$.
\end{proof}

We now show a boundedness property of branches of roots of $Q$ which holds away from the extremities $\lambda = 0$ and $\lambda = 1$.

\begin{lemma}
\label{lem:bound-branch}
Let $Q$ be the function defined by \eqref{eq:def-Q-lambda}, $\alpha \in \mathbb R$, and $K \subset (0, 1)$ be compact. Then there exists $C > 0$ such that, for every continuous branch $z\colon I\to\mathbb C$ of roots of $Q$ such that $\Real z(\lambda) \geq \alpha$ for every $\lambda \in I$, we have $\abs{z(\lambda)} \leq C$ for every $\lambda \in I \cap K$.
\end{lemma}

\begin{proof}
Let $\lambda_{\min} = \inf K$ and $\lambda_{\max} = \sup K$. Then, using the fact that $Q(z(\lambda), \lambda) = 0$ together with \eqref{eq:def-Q-lambda}, we obtain that, for every $\lambda \in I \cap K$, we have
\begin{align*}
\abs{z(\lambda)} & \le \frac{\lambda+1}{\lambda} + \frac{\lambda}{1 - \lambda} \abs*{e^{-z(\lambda)}} + \frac{1}{\lambda(1 - \lambda)} \abs*{e^{-\lambda z(\lambda)}} \\
& \leq \frac{\lambda_{\min} + 1}{\lambda_{\min}} + \frac{\lambda_{\max}}{1 - \lambda_{\max}} e^{-\alpha} + \frac{1}{\lambda_{\min} (1 - \lambda_{\max})} e^{\lambda_{\max} \abs{\alpha}}.
\end{align*}
The latter quantity only depends on $K$ and $\alpha$, and provides thus the required bound.
\end{proof}

Our next result shows that, if a branch of roots crosses the imaginary axis outside of the origin, then it necessarily crosses it from the left to the right.

\begin{lemma}
\label{lem:derivative-on-imaginary-axis}
Let $Q$ be the function defined by \eqref{eq:def-Q-lambda}, $z\colon I\to\mathbb C$ be an analytic branch of roots of $Q$, and $\lambda_0 \in I$ be a value such that $\Real(z(\lambda_0)) = 0$ and $z(\lambda_0) \neq 0$. Then
\begin{equation}
\label{eq:lem_dp1}
\Real(z^\prime(\lambda_0)) > 0.
\end{equation}
\end{lemma}

\begin{proof}
To simplify the notations in the proof, let $y = \Imag(z(\lambda_0))$. Note that, by Lemma~\ref{lem:simple-root}, $z(\lambda_0)$ is a simple root of $Q(\cdot, \lambda_0)$. By taking the real and imaginary parts, we have that $Q(i y, \lambda_0) = 0$ if and only if
\begin{subequations}
\begin{align}
 \label{eq:-2}
\cos(\lambda_0 y)-\lambda_0^2 \cos y+\lambda_0^2 -1& =0,\\
 \label{eq:-2bis}
\sin(\lambda_0 y)-\lambda_0^2 \sin y+\lambda_0^2 y-\lambda_0 y& =0.
\end{align}
\end{subequations}
Let $\widetilde Q$ be defined as in the proof of Lemma~\ref{lem:simple-root} and recall that \eqref{eq:derivative} is equivalent to
\begin{equation}
\label{eq:deriv-z-lambda0}
z^\prime(\lambda_0) = - \frac{\partial_\lambda\widetilde Q(i y, \lambda_0)}{\partial_z \widetilde Q(i y, \lambda_0)}.
\end{equation}
In the sequel, we study these two partial derivatives separately.

\step{1}{Computation of $\partial_z \widetilde{Q}(i y, \lambda_0)$}

From \eqref{eq:lem_dp4}, we have
\begin{equation*}
\begin{split}
\partial_z \widetilde{Q}(i y,\lambda_0) & = \lambda_0(1-\lambda_0)+\lambda_0^2e^{-i y}-\lambda_0 e^{-i \lambda_0 y} \\
& = \lambda_0(1-\lambda_0)+\lambda_0^2(\cos y-i \sin y) - \lambda_0 (\cos( \lambda_0 y)-i \sin(\lambda_0 y)).
\end{split}
\end{equation*}
Setting $R_1(y, \lambda_0) = \Real(\partial_z \widetilde Q(i y, \lambda_0))$ and using \eqref{eq:-2}, we compute
\begin{equation}
\label{eq:lem_dp6}
\begin{split}
R_1(y,\lambda_0) & = -\lambda_0 \cos(\lambda_0 y)+\lambda_0^2 \cos y-\lambda_0^2+\lambda_0 \\
& = -\lambda_0 (\lambda_0^2 \cos y-\lambda_0^2+1) +\lambda_0^2 \cos y-\lambda_0^2+\lambda_0 \\
& = \lambda_0^2(1-\lambda_0)(\cos y-1).
\end{split}
\end{equation}
Similarly, setting $I_1(y, \lambda_0) = \Imag(\partial_z \widetilde Q(i y, \lambda_0))$ and using \eqref{eq:-2bis}, we have
\begin{equation}
\label{eq:lem_dp7}
\begin{split}
I_1(y,\lambda_0) & = \lambda_0 \sin(\lambda_0 y)-\lambda_0^2 \sin y \\
& = \lambda_0(\lambda_0^2 \sin y-\lambda_0^2 y+\lambda_0 y)-\lambda_0^2 \sin y \\
& = \lambda_0^2(1-\lambda_0)(y - \sin y).
\end{split}
\end{equation}

\step{2}{Computation of $\partial_\lambda \widetilde{Q}(i y, \lambda_0)$}

From \eqref{eq:lem_dp8}, we have
\begin{equation*}
\begin{split}
\partial_{\lambda} \widetilde{Q}(i y,\lambda_0) & = (1-2\lambda_0)i y+2\lambda_0 -2\lambda_0 e^{-i y}-i y e^{-i \lambda_0 y} \\
& = (1-2\lambda_0)i y+2\lambda_0 -2\lambda_0 (\cos y-i \sin y) -i y (\cos(\lambda_0 y) - i \sin(\lambda_0 y)).
\end{split}
\end{equation*}
Setting $R_2(y, \lambda_0) = \Real(\partial_\lambda \widetilde Q(i y, \lambda_0))$ and using \eqref{eq:-2bis}, we compute
\begin{equation}
\label{eq:lem_dp10}
\begin{split}
R_2(y,\lambda_0) & = 2\lambda_0-2 \lambda_0 \cos y - y \sin(\lambda_0 y) \\
& = 2\lambda_0-2 \lambda_0 \cos y - y (\lambda_0^2 \sin y-\lambda_0^2 y+\lambda_0 y) \\
& = 2 \lambda_0-2 \lambda_0 \cos y - y \lambda_0^2 \sin y+\lambda_0^2 y^2-\lambda_0 y^2.
\end{split}
\end{equation}
Similarly, setting $I_2(y, \lambda_0) = \Imag(\partial_\lambda \widetilde Q(i y, \lambda_0))$ and using \eqref{eq:-2}, we have
\begin{equation}
\label{eq:lem_dp11}
\begin{split}
I_2(y,\lambda_0) & = (1-2\lambda_0) y+2 \lambda_0 \sin y- y \cos(\lambda_0 y) \\
& = (1-2\lambda_0) y+2 \lambda_0 \sin y - y (\lambda_0^2\cos y-\lambda_0^2+1) \\
& = 2 \lambda_0(\sin y-y)+\lambda_0^2 y(1-\cos y).
\end{split}
\end{equation}

\step{3}{Sign of $\Real(z^\prime(\lambda_0))$}

From \eqref{eq:deriv-z-lambda0} and the definitions of $R_1(y, \lambda_0)$, $I_1(y, \lambda_0)$, $R_2(y, \lambda_0)$, and $I_2(y, \lambda_0)$, we have
\begin{equation}
\label{eq:lem_dp12}
\Real(z^\prime(\lambda_0)) = -\frac{R_1(y,\lambda_0) R_2(y, \lambda_0)+I_1(y,\lambda_0) I_2(y,\lambda_0)}{\left(R_1(y,\lambda_0)^2 + I_1(y,\lambda_0)^2\right)}.
\end{equation}
In particular, $\Real(z^\prime(\lambda_0))$ has the same sign as the numerator of the right-hand side of \eqref{eq:lem_dp12}. Using \eqref{eq:lem_dp6}, \eqref{eq:lem_dp7}, \eqref{eq:lem_dp10}, and \eqref{eq:lem_dp11}, we have
\begin{equation}
\label{eq:lem_dp12bis}
\begin{split}
\MoveEqLeft[4] -R_1(y,\lambda_0) R_2(y, \lambda_0) - I_1(y,\lambda_0) I_2(y,\lambda_0) \\
& = \lambda_0^3 (1-\lambda_0) \bigl((1-\cos y) (2-2\cos y-y \lambda_0 \sin y+\lambda_0 y^2-y^2) \\
& \hphantom{{} = {}} + (\sin y-y)(2(\sin y-y) +\lambda_0 y (1-\cos y) \bigr) \\
& = \lambda_0^3 (1-\lambda_0)\bigl(2(1-\cos y)^2 +2(\sin y-y)^2-y^2(1-\cos y) \bigr) \\
& = 2\lambda_0^3 (1-\lambda_0)\left(y \cos\left(\tfrac{y}{2}\right)-2\sin\left(\tfrac{y}{2}\right)\right)^2.
\end{split}
\end{equation}
Hence $\Real(z^\prime(\lambda_0)) \geq 0$.

In order to prove that this inequality is strict, we use the identity $\cos^2(\lambda_0 y) + \sin^2(\lambda_0 y) = 1$ together with~\eqref{eq:-2} and \eqref{eq:-2bis} to obtain that
\begin{equation*}
2(\lambda_0+1)(1-\cos y)-(1-\lambda_0)y^2-2 y \lambda_0 \sin y = 0.
\end{equation*}
Isolating $\lambda_0$, we obtain
\begin{equation}
\label{eq:-4}
\lambda_0 = \frac{y^2 - 2(1-\cos y)}{2(1-\cos y)+y^2-2 y \sin y} = \frac{y^2+2(\cos y-1)}{(y-\sin y)^2+(1-\cos y)^2}.
\end{equation}
Let us now assume, to obtain a contradiction, that $\Real(z^\prime(\lambda_0)) = 0$. Hence, by \eqref{eq:lem_dp12bis}, we have $\tan\left(\tfrac{y}{2}\right) = \frac{y}{2}$ and
\[(1-\cos y)^2 + (\sin y-y)^2 = \frac{y^2}{2}(1-\cos y).\]
Inserting the above into \eqref{eq:-4}, we get
\begin{equation}
\label{eq:lem_dp13}
\lambda_0 = \frac{2(y^2+2(\cos y-1))}{y^2(1-\cos y)}.
\end{equation}
On the other hand, the equality $\tan\left(\tfrac{y}{2}\right) = \frac{y}{2}$ implies that
\[
\frac{y^2}{4} = \tan^2\left(\frac{y}{2}\right) = \frac{1 - \cos y}{1 + \cos y},
\]
i.e.,
\[
4(\cos y - 1) = -y^2(1 + \cos y).
\]
Inserting the above into the numerator of \eqref{eq:lem_dp13}, we deduce that $\lambda_0 = 1$, yielding the required contradiction. Hence, \eqref{eq:lem_dp1} is proved.
\end{proof}

Using the previous result, we can now show that, if a branch is in the right half-plane at some $\lambda_0$, it must remain there for all larger $\lambda$ and it will be defined until $\lambda = 1$.

\begin{lemma}
\label{lemm:extension}
Let $Q$ be the function defined by \eqref{eq:def-Q-lambda}, $z\colon (\underline\lambda, \overline\lambda)\to\mathbb C$ be an analytic branch of roots of $Q$, and $\lambda_0 \in (\underline\lambda, \overline\lambda)$ be such that $\Real(z(\lambda_0)) \geq 0$. Then $z$ admits a unique extension (still denoted by $z$) defined on $(\underline\lambda, 1)$. In addition, this extension satisfies $\Real(z(\lambda)) \geq 0$ for every $\lambda \in (\underline\lambda, 1)$. Moreover, either $z\equiv 0$ on $(\underline\lambda, 1)$ or $z(\lambda)\neq 0$ for every $\lambda \in (\underline\lambda, 1)$.
\end{lemma}

\begin{proof}
The case $z(\lambda_0) = 0$ is trivial since, in this case, $z$ can be extended to the function $z\colon (\underline\lambda, 1) \to \mathbb C$ defined by $z(\lambda) = 0$ for every $\lambda \in (\underline, 1)$, and this extension is unique by Hurwitz's theorem. Hence, in the sequel, we assume that $z(\lambda_0) \neq 0$.

By Lemma~\ref{lem:branch-not-zero}, we have that any extension of $z$ will never be equal to $0$. We now claim that any extension of $z$ will remain in the closed right-half plane for all $\lambda \geq \lambda_0$ in its domain of definition. Indeed, let us consider an extension $z\colon I \to \mathbb C$ to an interval $I$ with $(\underline\lambda, \overline\lambda) \subset I$. Assume, to obtain a contradiction, that there exists $\lambda_1 \in I$ with $\lambda_1 > \lambda_0$ such that $\Real(z(\lambda_1)) < 0$. As $\lambda \mapsto \Real(z(\lambda))$ is $\mathcal C^1$, we deduce that there exists $\lambda_\ast \in [\lambda_0, \lambda_1)$ such that $\Real(z(\lambda_\ast)) = 0$ and $\Real(z^\prime(\lambda_\ast)) \leq  0$. This, however, contradicts Lemma~\ref{lem:derivative-on-imaginary-axis}, proving the claim.

As a consequence of the above and Lemma~\ref{lem:simple-root}, for any extension of $z$, $z(\lambda)$ is a simple root of $Q(\cdot, \lambda)$ for every $\lambda \geq \lambda_0$ in the domain of the extension. Hence, arguing by Hurwitz's theorem, extensions of $z$ necessarily coincide on the intersection of their domains. In particular, if $z$ admits an extension to $(\underline\lambda, 1)$, then it is necessarily unique.

To prove that $z$ admits an extension to $(\underline\lambda, 1)$, let
\[
\lambda_{\max} = \sup\{\lambda \in [\lambda_0, 1] \suchthat z \text{ admits an extension to } (\underline\lambda, \lambda)\}.
\]
Assume, to obtain a contradiction, that $\lambda_{\max} < 1$, and consider an extension of $z$ to $(\underline\lambda, \lambda_{\max})$. Since $[\lambda_0, \lambda_{\max}]$ is compact and $z$ belongs to the closed right half-plane on $[\lambda_0, \lambda_{\max})$, we deduce by Lemma~\ref{lem:bound-branch} that there exists $C > 0$ such that $\abs{z(\lambda)} \leq C$ for every $\lambda \in [\lambda_0, \lambda_{\max})$. By Lemma~\ref{lem:simple-root}, we deduce that, up to increasing $C$, we have $\abs{z^\prime(\lambda)} \leq C$ for every $\lambda \in [\lambda_0, \lambda_{\max})$. In particular, $z(\lambda)$ is Lipschitz continuous on $[\lambda_0, \lambda_{\max})$, and thus $\lim_{\lambda \to \lambda_{\max}} z(\lambda)$ exists, and it is nonzero by Hurwitz's theorem. We can thus extend $z$ to a Lipschitz continuous function on $(\underline\lambda, \lambda_{\max}]$ with $Q(z(\lambda_{\max}), \lambda_{\max}) = 0$. By Lemma~\ref{lem:simple-root}, there exists a neighborhood $I$ of $\lambda_{\max}$ in $(0, 1)$ and an analytic branch $\widehat z\colon I \to \mathbb C$ of roots of $Q$ such that $\widehat z(\lambda_{\max}) = z(\lambda_{\max})$. Arguing by Hurwitz's theorem, one deduces that $z$ and $\widehat z$ coincide on $(\underline\lambda, \lambda_{\max}] \cap I$, and thus $z$ can be extended to $(\underline\lambda, \lambda_{\max}] \cup I$. This, however, contradicts the maximality of $\lambda_{\max}$, yielding the conclusion.
\end{proof}

To proceed, it is enough to rule out of the existence of a branch of nontrivial roots of $Q$ in the right half-plane. We start by proving that any such branch is uniformly bounded.

\begin{lemma}
\label{lem:bdd}
Let $Q$ be the function defined by \eqref{eq:def-Q-lambda}, $z\colon (\underline\lambda, 1) \to \mathbb C$ be an analytic branch of nontrivial roots of $Q$, and assume that there exists $\lambda_0 \in (\underline\lambda, 1)$ such that $\Real(z(\lambda_0)) \geq 0$. Then $z$ is bounded on $[\lambda_0, 1)$.
\end{lemma}

\begin{proof}
Notice first that, by Lemma~\ref{lemm:extension}, we have $\Real(z(\lambda)) \geq 0$ for every $\lambda \in [\lambda_0, 1)$. For $\lambda \in [\lambda_0, 1)$, define $x(\lambda) = \Real(z(\lambda))$ and $y(\lambda) = \Imag(z(\lambda))$ and note that $x(\lambda) \geq 0$. Since $Q(\cdot, \lambda)$ is a quasipolynomial with real coefficients, its nonreal roots appear in complex conjugate pairs and, combining this fact with Lemma~\ref{lem:cannot-cross-real-axis}, we assume, with no loss of generality, that $y(\lambda) > 0$ for every $\lambda \in [\lambda_0, 1)$.

We first claim that $x$ is bounded on $[\lambda_0, 1)$. Indeed, assume, to obtain a contradiction, that there exists a sequence $(\lambda_n)_{n \in \mathbb N}$ in $[\lambda_0, 1)$ such that $\lambda_n \to 1$ and $x(\lambda_n) \to +\infty$ as $n \to +\infty$. For every $\lambda \in [\lambda_0, 1)$, since $Q(z(\lambda),\lambda) = 0$, we have
\begin{equation}
\label{eq:lem:bor-1}
z(\lambda) - \frac{\lambda+1}{\lambda} = \frac{e^{-\lambda z(\lambda)}}{1-\lambda}\left(\lambda e^{-Z(\lambda)} - \frac{1}{\lambda} \right),
\end{equation}
where we have set $Z(\lambda) = (1 - \lambda) z(\lambda)$. After multiplying \eqref{eq:lem:bor-1} by $1-\lambda$ and using the fact that $e^{-\lambda_n z(\lambda_n)} \to 0$ as $n \to +\infty$, we deduce that $Z(\lambda_n) \to 0$ as $n \to +\infty$. Using a Taylor expansion of the term $e^{-Z(\lambda_n)}$, we obtain that
\begin{equation*}
\begin{split}
 z(\lambda_n) - \frac{\lambda_n + 1}{\lambda_n} & = \frac{e^{-\lambda_n z(\lambda_n)}}{1-\lambda_n} \left(\lambda_n \left(1 - Z(\lambda_n) + O\left(Z(\lambda_n)^2\right)\right) - \frac{1}{\lambda_n} \right) \\
  & = -e^{-\lambda_n z(\lambda_n)} \biggl(\frac{1}{\lambda_n} + 1 + \lambda_n z(\lambda_n) + z(\lambda_n) O(Z(\lambda_n)) \biggr).
\end{split}
\end{equation*}
By factoring in the above equation terms containing $z(\lambda_n)$, one gets
\begin{equation*}
z(\lambda_n) = \frac{ 1 - e^{-\lambda z(\lambda_n)}}{1 + \lambda_n e^{-\lambda_n z(\lambda_n)} + e^{ -\lambda_n z(\lambda_n)} O(Z(\lambda_n))}\left(1 + \frac{1}{\lambda_n}\right).
\end{equation*}
The right-hand side of the above equation converging to $2$ as $n \to +\infty$, we reach a contradiction and prove the claim on $x$.

Since $z$ is a branch of roots, by taking the real and imaginary parts of the identity $(1 - \lambda) Q(z(\lambda), \lambda) = 0$ and multiplying by $e^{x(\lambda)}$, we deduce that, for every $\lambda \in [\lambda_0, 1)$,
\begin{subequations}
\begin{align}
(1-\lambda) e^{x(\lambda)} x(\lambda)-\frac{1-\lambda^2}{\lambda} e^{x(\lambda)} & = \lambda \cos(y(\lambda))-\frac{e^{(1-\lambda) x(\lambda)}}{\lambda}\cos(\lambda y(\lambda)), \\
\label{eq:lem:bor3}
(1-\lambda) e^{x(\lambda)} y(\lambda) & =- \lambda \sin(y(\lambda)) + \frac{e^{(1-\lambda) x(\lambda)}}{\lambda} \sin(\lambda y(\lambda)).
\end{align}
\end{subequations}
Taking into account that $x$ is bounded, one can simplify \eqref{eq:lem:bor3} to obtain that, as $\lambda \to 1^-$,
\begin{equation*}
\sin(y(\lambda))-\frac{1}{\lambda^2}\sin(\lambda y(\lambda)) = -\frac{1-\lambda}{\lambda} e^{x(\lambda)} y(\lambda) + O(1-\lambda).
\end{equation*}
Using that $1/\lambda^2 = 1 + O(1 - \lambda)$ as $\lambda \to 1^-$, one obtains
\begin{equation}\label{eq:lem:5}
\sin(y(\lambda))-\sin(\lambda y(\lambda)) = -\frac{1-\lambda}{\lambda} e^{x(\lambda)} y(\lambda) + O(1-\lambda).
\end{equation}
Assuming that $z$ is not bounded means that its imaginary part is not bounded, i.e., the continuous function $\lambda\mapsto y(\lambda)$ is unbounded as $\lambda \to 1^-$. Then one can assume with no loss of generality that there exists a sequence $(\lambda_n)_{n \in \mathbb{N}}$ in $[\lambda_0, 1)$ so that $\lambda_n \to 1$ as $n \to +\infty$ and 
\[
y(\lambda_n) = \frac{\pi}2 + 2 n \pi
\]
for $n$ large enough. Considering \eqref{eq:lem:5} along the sequence 
$(\lambda_n)_{n \in \mathbb{N}}$ yields the equation
\[
\frac{1-\sin(\lambda_n y(\lambda_n))}{(1-\lambda_n)y(\lambda_n)}=
-\frac{e^{x(\lambda_n)}}{\lambda_n}+\frac{O(1)}{y(\lambda_n)}
\]
for $n$ large enough. The left-hand side of the above equation is clearly nonnegative while the limsup of the right-hand side, as $n$ tends to infinity, is upper bounded by $-1$. We have reached a contradiction and hence proved that $z$ is bounded on $[\lambda_0,1)$.
\end{proof}

The next step consists in proving that any analytical branch of nontrivial roots of $Q$ in the right half-plane defined in an interval of the form $(\underline\lambda, 1)$ can be extended as an analytic function in a right neighborhood of $1$. For that purpose, we extend the definition of $Q$ to $\mathbb C \times \mathbb R \setminus \{0, 1\}$ by using the same formula \eqref{eq:def-Q-lambda}. We must also identify the limit of $Q(\cdot, \lambda)$ as $\lambda \to 1$, and a Taylor expansion shows that
\begin{equation}
\label{eq:Q-lambda-1}
\lim_{\lambda \to 1} Q(z, \lambda) = z - 2 + e^{-z}(z + 2),
\end{equation}
uniformly for $z$ on compact subsets of $\mathbb C$. In the sequel, we denote the right-hand side of the above equality by $Q(z, 1)$. Note that $Q(\cdot, 1)$ is a quasipolynomial of the form \eqref{eq:Delta-DDE} with $n = m = 1$, of degree $3$, and admitting $0$ as a root of multiplicity $3$. Hence, by Theorem~\ref{thm:main-MID}, we have $\Real z = 0$ for every root $z$ of $Q(\cdot, 1)$. Moreover, Proposition~\ref{prop:spectrum-neutral} shows that the set of nonzero roots of $Q(\cdot, 1)$ coincides with the set of complex numbers of the form $2 \zeta i$, where $\zeta\in\mathbb R^\ast$ is any nontrivial solution of the equation $\tan(x)=x$ with $x\in \mathbb R$. Finally, an immediate computation shows that any nonzero root of $Q(\cdot, 1)$ is necessarily simple.

Using these facts, we can prove the following result.

\begin{lemma}
\label{lem:near_1}
Let $Q$ be the function defined by \eqref{eq:def-Q-lambda}, $z\colon (\underline\lambda, 1) \to \mathbb C$ be an analytic branch of nontrivial roots of $Q$, and assume that there exists $\lambda_0 \in (\underline\lambda, 1)$ such that $\Real(z(\lambda_0)) \geq 0$. Then there exist $\epsilon>0$ and an analytic function $\widehat z\colon (\underline\lambda, 1 + \epsilon) \to \mathbb C$ such that $Q(\widehat z(\lambda),\lambda) = 0$ for $\lambda\in (\underline\lambda, 1+\epsilon)$, $\widehat z(1) = 2 \zeta i$ for a nonzero real number $\zeta$, and $z(\lambda) = \widehat z(\lambda)$ for every $\lambda \in (\underline\lambda, 1)$.
\end{lemma}

\begin{proof}
By Lemma~\ref{lem:bdd}, $z$ is bounded on $[\lambda_0, 1)$, and thus, arguing by Hurwitz's theorem, $z(\lambda)$ must converge to a nontrivial root $2 \zeta i$ of $Q(\cdot, 1)$ as $\lambda$ tends to $1$.

Let $U$ be an open neighborhood of $(2 \zeta i,1)$ in $\mathbb{C}^2$ and consider the natural extension of $Q$ to $U$ (i.e., $Q(z, \lambda)$ is given by \eqref{eq:def-Q-lambda} if $\lambda \neq 1$ and $Q(z, 1) = z - 2 + e^{-z}(z + 2)$), still denoted by $Q$. A trivial computation from \eqref{eq:def-Q-lambda} shows that the partial derivatives $\partial_z Q(z,\lambda)$ and $\partial_{\lambda} Q(z,\lambda)$ exist and are analytic for $(z,\lambda) \in U \setminus(2 \zeta i,1)$. Moreover, $\partial_z Q(z,\lambda)$ and $\partial_{\lambda} Q(z,\lambda)$ admit limits as $(z,\lambda) \to (2 \zeta i,1)$, so that $\partial_z Q$ and $\partial_{\lambda} Q$ are continuous on $U$. We deduce that $Q(\cdot,\cdot)$ is analytic on $U$. In addition, we compute $\partial_z Q( 2 \zeta i,1)=1-e^{-2 \zeta i}(1-2 \zeta i)$, and thus $\partial_z Q( 2 \zeta i,1) \neq 0$ since otherwise $1=e^{-2 \zeta i}(1-2 \zeta i)$ and one would get that $\zeta=0$ after taking the modulus. The conclusion of the lemma follows by applying the analytic implicit function theorem.
\end{proof}

\subsection{Proof of Theorem~\ref{thm:main-two-delays}}
\label{sec:proof-two-delays}

We are finally in position to provide the proof of Theorem~\ref{thm:main-two-delays}.

\begin{proof}[Proof of Theorem~\ref{thm:main-two-delays}]
With no loss of generality, we assume here that $0 < \tau_1 < \tau_2$. The first part of the proof of Theorem~\ref{thm:main-two-delays} was shown in Proposition~\ref{prop:main-two-delays-conditions}, while, by Lemma~\ref{lemm:transformation}, the second part of the statement amounts to showing that $\Real z < 0$ for every $(z, \lambda) \in \mathbb C^\ast \times (0, 1)$ such that $Q(z, \lambda) = 0$, where $Q$ is the function defined in \eqref{eq:def-Q-lambda}.

Assume, to obtain a contradiction, that there exist $\lambda_0 \in (0, 1)$ and $z_0 \in \mathbb C^\ast$ such that $Q(z_0, \lambda_0) = 0$ and $\Real z_0 \geq 0$. Combining Lemmas~\ref{lem:simple-root}, \ref{lemm:extension}, and \ref{lem:near_1}, there exist $\underline\lambda \in (0, \lambda_0)$, $\overline\lambda > 1$, and an analytic function $z\colon (\underline\lambda, \overline\lambda) \to \mathbb C$ such that $z(\lambda_0) = z_0$, $Q(z(\lambda), \lambda) = 0$ for every $\lambda \in (\underline\lambda, \overline\lambda)$, $\Real(z(\lambda)) \geq 0$ for every $\lambda \in [\lambda_0, 1)$, and $z(1) = 2 \zeta i$ for some $\zeta \in \mathbb R^\ast$.

Let us compute $z^\prime(1)$ and $z^{\prime\prime}(1)$. Since $z$ is analytic, we have
\begin{equation}\label{eq:serie_entiere1}
z(\lambda)=\alpha_0+\alpha_1(1-\lambda)+\alpha_2(1-\lambda)^2+O\left((1-\lambda)^3\right)
\end{equation}
for $\lambda$ in a neighborhood of $1$, where $\alpha_0 = 2 \zeta i$, $\alpha_1 = -z^\prime(1)$, and $\alpha_2 = \frac{1}{2} z^{\prime\prime}(1)$. In particular, we compute
\begin{equation}
\label{eq:serie_entiere1b}
e^{-z(\lambda)} = e^{-\alpha_0} \Biggl(1 - \alpha_1 (1 - \lambda) + \left(\frac{\alpha_1^2}{2} - \alpha_2\right) (1 - \lambda)^2 + O\left((1 - \lambda)^3\right)\Biggr).
\end{equation}
We also expand \eqref{eq:def-Q-lambda} in a neighborhood of $\lambda = 1$, obtaining
\begin{equation}
\label{eq:serie_entiere2}
\begin{split}
Q(z, \lambda) = {} & z - 2 + e^{-z}(z + 2) + \left(\left(1 + z + \frac{z^2}{2}\right)e^{-z} - 1\right) (1 - \lambda) \\
&{} + \left(\left(1 + z + \frac{z^2}{2} + \frac{z^3}{6}\right)e^{-z} - 1\right) (1 - \lambda)^2 + O\left((1-\lambda)^3\right),
\end{split}
\end{equation}
where $O\left((1-\lambda)^3\right)$ is uniform with respect to $z$ on compact subsets of $\mathbb C$. Inserting \eqref{eq:serie_entiere1} and \eqref{eq:serie_entiere1b} in \eqref{eq:serie_entiere2}, one obtains that, as $\lambda \to 1$,
\begin{equation}
\label{eq:serie_entiere3}
\begin{split}
0=Q(z(\lambda),\lambda)= {} & (\alpha_0+2)e^{-\alpha_0}+\alpha_0-2 + \left( \tfrac{\alpha_0^2+2\alpha_0+2-2\alpha_1(1+\alpha_0)}{2}e^{-\alpha_0} +\alpha_1-1 \right) (1-\lambda) \\
& {} + \left( \tfrac{\alpha_0^3+ 3(1 - \alpha_1) \alpha_0^2 + 3(2 - 2\alpha_2 + \alpha_1^2)\alpha_0 + 6 (1 - \alpha_2)}{6}e^{-\alpha_0} +\alpha_2-1 \right) (1-\lambda)^2 \\
& {} + O\left((1-\lambda)^3\right).
\end{split}
\end{equation}
We deduce from \eqref{eq:serie_entiere3} that
\begin{equation}
\label{eq:serie_entiere4}
(\alpha_0+2)e^{-\alpha_0}+\alpha_0-2=0, \qquad \alpha_1=\frac{\alpha_0}{2}, \qquad \alpha_2=\frac{\alpha_0(\alpha_0+10)}{24}.
\end{equation}
As $\alpha_0 = 2 \zeta i$ with $\zeta \in \mathbb R^\ast$, $\alpha_1 = -z^\prime(1)$, and $\alpha_2 = \frac{1}{2} z^{\prime\prime}(1)$, it follows from \eqref{eq:serie_entiere4} that $\Real(z^\prime(1))=0$ and $ \Real(z^{\prime\prime}(1)) = -\frac{\zeta^2}{3} < 0$, yielding that $ \Real(z(\lambda))<0$ for $\lambda$ 
in a real neighborhood of $1$, in contradiction with the fact that $\Real(z(\lambda)) \geq 0$ for every $\lambda \in [\lambda_0, 1)$.
\end{proof}

\subsection{Numerical computation of the roots of \texorpdfstring{$Q$}{Q}}
\label{sec:roots-Q}

Similarly to what was done in Section~\ref{sec:numerical-widehat-delta}, we now provide a numerical computation\footnote{Numerical computations were performed using Python's \texttt{cxroots} module \cite{Parini2018cxroots}, searching for roots of $Q(\cdot, \lambda)$ in the rectangle $\{z \in \mathbb C \suchthat -10 \leq \Real z \leq 5,\, \abs{\Imag z} \leq 33\}$ and for 1001 evenly distributed values of $\lambda$ in $[0, 1]$. For $\lambda = 1$, we have computed the roots of the quasipolynomial $Q(\cdot, 1)$ defined, for $z \in \mathbb C$, as the right-hand side of \eqref{eq:Q-lambda-1}, while, for $\lambda = 0$, we have computed the roots of the quasipolynomial $z \mapsto \frac{z^2}{2} - z + 1 - e^{-z}$.} of the roots of the quasipolynomial $Q(\cdot, \lambda)$ from \eqref{eq:def-Q-lambda} for different values of $\lambda \in (0, 1)$. The result of this computation is provided in Figure~\ref{fig:roots-two-delays}, which represents the computed roots of $Q(\cdot, \lambda)$ as $\lambda$ varies in $(0, 1)$ through different colors, with blue shades representing values of $\lambda$ close to $0$ and yellow shades representing values of $\lambda$ close to $1$. The triple root at $0$ is represented by a black dot.

\begin{figure}[ht]
\centering
\resizebox{\textwidth}{!}{\import{Figures}{Two_Delays.pgf}}
\caption[{Roots of the quasipolynomial $Q(\cdot, \lambda)$ from \eqref{eq:def-Q-lambda} for $\lambda \in [0, 1]$}]{Roots of the quasipolynomial $Q(\cdot, \lambda)$ from \eqref{eq:def-Q-lambda} for $\lambda \in [0, 1]$, with a detailed view of the region $\{z \in \mathbb C \suchthat -6 \leq \Real z \leq 0,\, 27 \leq \Imag z \leq 32\}$}
\label{fig:roots-two-delays}
\end{figure}

We observe in Figure~\ref{fig:roots-two-delays} that, as shown in Theorem~\ref{thm:main-two-delays}, for every $\lambda \in (0, 1)$, all nontrivial roots of $Q(\cdot, \lambda)$ have strictly negative real part. One can also observe further information on the behavior of the roots obtained in the proof of Theorem~\ref{thm:main-two-delays}: for any nontrivial branch of roots $\lambda \mapsto z(\lambda)$, as $\lambda \to 1$, $z(\lambda)$ approaches a number of the form $2 \zeta i$ for some $\zeta \in \mathbb R^\ast$ satisfying $\tan \zeta = \zeta$, and this is done tangentially to the imaginary axis and from the left, since $\Real z(1) = \Real z^\prime(1) = 0$ and $\Real z^{\prime\prime}(1) < 0$.

Note that, as $\lambda \to 1$, the represented roots converge to the roots of the quasipolynomial $Q(\cdot, 1)$. As for the limit $\lambda \to 0$, we observe that
\[
\lim_{\lambda \to 0} \frac{Q(z, \lambda)}{\lambda} = \frac{z^2}{2} - z + 1 - e^{-z},
\]
uniformly for $z$ on compact subsets of $\mathbb C$. The quasipolynomial in the right-hand side of the above equality corresponds, up to a multiplicative factor of $\frac{1}{2}$, to the quasipolynomial of the form \eqref{eq:Delta-DDE} with $n = 2$, $m = 0$, $\tau = 1$, degree $D = 3$, and whose coefficients are given by \eqref{Coeffs} with $s_0 = 0$, admitting thus a root of maximal multiplicity at $0$, which is strictly dominant thanks to Theorem~\ref{thm:main-MID}.

\subsection{Insights on systems with more delays}
\label{sec:three-delays}

At the light of Theorem~\ref{thm:main-two-delays}, a natural question is whether a similar result would hold for systems with three (or more) delays. The answer to this problem is negative, as we now detail.

Consider the delay-differential equation
\begin{equation}
\label{eq:main-syst-three-delays}
y^\prime(t) + a_0 y(t) + a_1 y(t - \tau_1) + a_2 y(t - \tau_2) + a_3 y(t - \tau_3) = 0,
\end{equation}
where $a_0, a_1, a_2, a_3$ are real coefficients and $\tau_1, \tau_2, \tau_3$ are pairwise distinct positive real numbers representing delays. Its characteristic quasipolynomial is the function $\Delta\colon \mathbb C \to \mathbb C$ defined for $s \in \mathbb C$ by
\begin{equation*}
\Delta(s) = s + a_0 + a_1 e^{-s \tau_1} + a_2 e^{-s \tau_2} + a_3 e^{-s \tau_3}.
\end{equation*}
Note that the degree of $\Delta$ is $4$, and thus, by Corollary~\ref{coro:Polya-Szego}, the maximal multiplicity that a root $s_0$ of $\Delta$ may have is $4$. Straightforward computations, similar to those from the proof of Proposition~\ref{prop:main-two-delays-conditions}, show that $0$ is a root of multiplicity $4$ of $\Delta$ if and only if
\begin{align*}
a_0 & = - \frac{1}{\tau_{1}} - \frac{1}{\tau_{2}} - \frac{1}{\tau_{3}}, & a_1 & = \frac{\tau_{2} \tau_{3}}{\tau_{1} \left(\tau_{2} - \tau_{1}\right) \left(\tau_{3} - \tau_{1}\right)}, \\
a_2 & = \frac{\tau_{1} \tau_{3}}{\tau_{2} \left(\tau_{1} - \tau_{2}\right) \left(\tau_{3} - \tau_{2}\right)}, & a_3 & = \frac{\tau_{1} \tau_{2}}{\tau_{3} \left(\tau_{1} - \tau_{3}\right) \left(\tau_{2} - \tau_{3}\right)}.
\end{align*}
Let us then consider the function
\begin{multline}
\label{eq:def-Q-three-delays}
Q(s, \tau_1, \tau_2, \tau_3) = s - \frac{1}{\tau_{1}} - \frac{1}{\tau_{2}} - \frac{1}{\tau_{3}} + \frac{\tau_{2} \tau_{3}}{\tau_{1} \left(\tau_{2} - \tau_{1}\right) \left(\tau_{3} - \tau_{1}\right)} e^{-s \tau_1} \\ + \frac{\tau_{1} \tau_{3}}{\tau_{2} \left(\tau_{1} - \tau_{2}\right) \left(\tau_{3} - \tau_{2}\right)} e^{-s \tau_2} + \frac{\tau_{1} \tau_{2}}{\tau_{3} \left(\tau_{1} - \tau_{3}\right) \left(\tau_{2} - \tau_{3}\right)} e^{-s \tau_3},
\end{multline}
defined for $(s, \tau_1, \tau_2, \tau_3) \in \mathbb C \times (\mathbb R_+^\ast)^3$ such that $\tau_1 \neq \tau_2$, $\tau_2 \neq \tau_3$, and $\tau_1 \neq \tau_3$. One can compute that
\begin{equation}
\label{eq:three-delays-limit}
\lim_{(\tau_1, \tau_3) \to (1, 1)} Q(s, \tau_1, 1, \tau_3) = s - 3 + e^{- s} \left(\frac{s^{2}}{2} + 2 s + 3\right).
\end{equation}
The quasipolynomial on the right-hand side of the above expression corresponds to \eqref{eq:Delta-DDE} with $n = 1$, $m = 2$, $\tau = 1$, and with a root of multiplicity $4$ at $0$ (or, equivalently, with coefficients given by \eqref{Coeffs}). Hence, by Theorem~\ref{thm:main-MID}, we have $\Real s > 0$ for every nonzero root $s$ of the quasipolynomial on the right-hand side of \eqref{eq:three-delays-limit}, and so the quadruple root at $0$ is not dominant. Since the convergence in \eqref{eq:three-delays-limit} is uniform on compact subsets of $\mathbb C$, it follows from Hurwitz's theorem that the quadruple root at $0$ of $Q(\cdot, \tau_1, \tau_2, \tau_3)$ is not dominant if $\tau_2 = 1$ and $\tau_1$ and $\tau_3$ are close enough to $1$.

\begin{figure}[ht]
\centering
\resizebox{0.5\textwidth}{!}{\import{Figures}{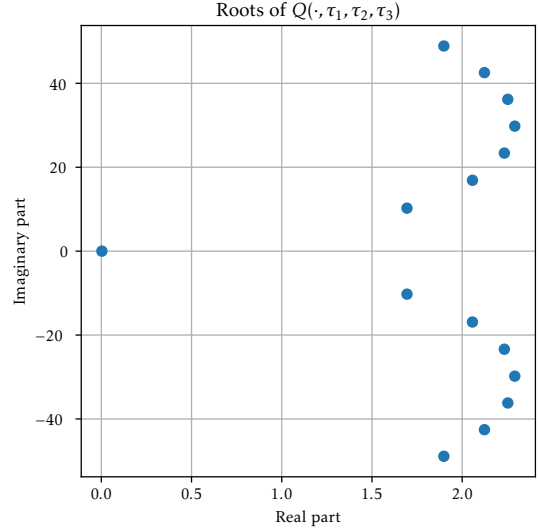}}
\caption[Roots of the quasipolynomial $Q(\cdot, \tau_1, \tau_2, \tau_3)$ from \eqref{eq:def-Q-three-delays} with $\tau_1 = 0.917686$, $\tau_2 = 1$, and $\tau_3 = 1.067836$]{Roots of the quasipolynomial $Q(\cdot, \tau_1, \tau_2, \tau_3)$ from \eqref{eq:def-Q-three-delays} with $\tau_1 = 0.917686$, $\tau_2 = 1$, and $\tau_3 = 1.067836$. The root at the origin has multiplicity $4$}
\label{fig:three-delays}
\end{figure}

Figure~\ref{fig:three-delays} represents the roots of $Q(\cdot, \tau_1, \tau_2, \tau_3)$ computed numerically\footnote{These roots were computed with Python's \texttt{cxroots} module \cite{Parini2018cxroots}, the search for roots being restricted to the rectangle $\{s \in \mathbb C \suchthat \abs{\Real s} \leq 10,\, \abs{\Imag s} \leq 50\}$.} for the values of delays $\tau_1 = 0.917686$, $\tau_2 = 1$, and $\tau_3 = 1.067836$, illustrating that the root of maximal multiplicity $4$ at $0$ is not dominant. Finding conditions on the delays $\tau_1, \tau_2, \tau_3$ ensuring that the MID property holds is currently an open problem.

\section{Examples and applications}
\label{sec:MID-applis}

We now present examples and applications of the results from Sections~\ref{sec:MID-maximal}--\ref{sec:MID-two-delays}, extracted from \cite{Mazanti2020Qualitative, Mazanti2021Effects, Benarab2020MID, Mazanti2021Multiplicity, Boussaada2022Generic, Boussaada2025Decay, Fueyo2023Pole}.

\subsection{Application of Theorem~\ref{thm:main-MID} to a third-order system}

As a first example, extracted from \cite{Mazanti2020Qualitative}, consider the delay-differential equation \eqref{eq:DDE} with $n = 3$, $m = 2$, and $\tau = 2.5$. The corresponding characteristic quasipolynomial $\Delta$ from \eqref{eq:Delta-DDE} is of degree $6$. We want to choose the real coefficients $a_0, a_1, a_2, \alpha_0, \alpha_1, \alpha_2$ of this equation in order for the origin to be an exponentially stable equilibrium point and, for that purpose, we impose that $\Delta$ should have a root at $s_0 = -0.5$ of maximal multiplicity $6$. By Theorem~\ref{thm:main-MID}, this is possible if and only if the coefficients $a_0, a_1, a_2, \alpha_0, \alpha_1, \alpha_2$ are given by
\begin{align*}
a_0 & = -1.735, & a_1 & = 2.91, & a_2 & = -2.1, \displaybreak[0] \\
\alpha_0 & \approx 1.736219, & \alpha_1 & \approx 1.443984, & \alpha_2 & \approx 0.3438058,
\end{align*}
and, in this case, we have $\Real s < s_0$ for every other root $s \in \mathbb C \setminus \{s_0\}$ of $\Delta$. In particular, the origin of \eqref{eq:DDE} is exponentially stable.

\begin{figure}[ht]
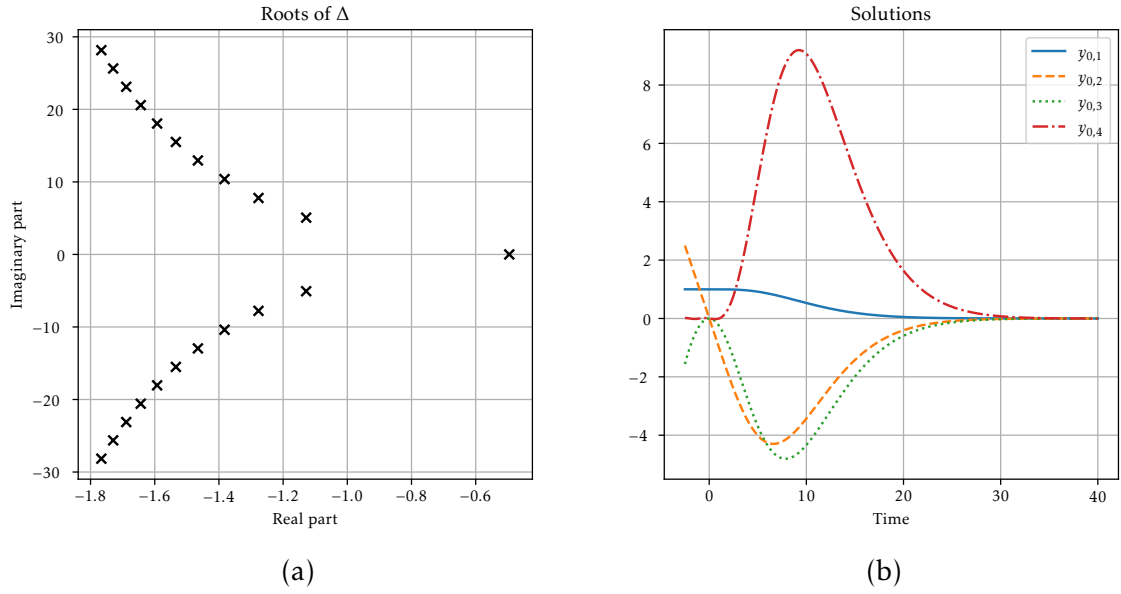

\centering
\begin{tabular}{@{} c @{} c @{}}
\resizebox{0.5\textwidth}{!}{\input{Figures/Expl1-zeros.pgf}} & \resizebox{0.5\textwidth}{!}{\input{Figures/Expl1-sols.pgf}} \tabularnewline
(a) & (b) \tabularnewline
\end{tabular}
\caption[{Roots of $\Delta$ from \eqref{eq:Delta-DDE} with $n = 3$, $m = 2$, $\tau = 2.5$, and a root $s_0 = -0.5$ of multiplicity $6$, and some solutions of the corresponding DDE}]{(a) Roots of the quasipolynomial $\Delta$ from \eqref{eq:Delta-DDE} with $n = 3$, $m = 2$, $\tau = 2.5$, and a root $s_0 = -0.5$ of multiplicity $6$. (b) Some solutions of the corresponding delay-differential equation \eqref{eq:DDE}}
\label{fig:expl-MID-1}
\end{figure}

Figure~\ref{fig:expl-MID-1}(a) shows the roots of $\Delta$ within the region $\{s \in \mathbb C \suchthat -5 \leq \Real s \leq 1,\, \abs{\Imag s} \leq 30\}$, computed numerically using Python's \texttt{cxroots} module \cite{Parini2018cxroots}. We observe in Figure~\ref{fig:expl-MID-1}(a) the presence of the multiple root $s_0 = -0.5$ as well as other roots dominated by $s_0$. Figure~\ref{fig:expl-MID-1}(b) illustrates the stability of the corresponding DDE \eqref{eq:DDE} with the above choice of coefficients by representing four particular solutions, obtained with the initial conditions $y_{0, i}$ for $i \in \{1, 2, 3, 4\}$ given, for $t \in [-\tau, 0]$, by
\begin{align*}
y_{0, 1}(t) & = 1, & y_{0, 2}(t) & = -t, \displaybreak[0] \\
y_{0, 3}(t) & = -\frac{t^2}{4}, & y_{0, 4}(t) & = -\frac{1}{6 \omega^2} \sin(\omega t),
\end{align*}
with $\omega = 2 \pi$.

\subsection{A delay differential-algebraic equation with scalar algebraic part}

This example, extracted from \cite{Mazanti2021Effects}, considers a DDAE of the form \eqref{eq:DDAE-MID} with $n = m = 1$, which we write as
\begin{equation}
\label{DiffAlgeb}
\left\{
\begin{aligned}
x^\prime(t) + a x(t) & = b k y(t - \tau), \\
y(t) & = c x(t) + d k y(t - \tau),
\end{aligned}
\right.
\end{equation}
where $x(t) \in \mathbb R$, $y(t) \in \mathbb R$, $\tau > 0$, and $a, b, c, d, k$ are real coefficients.

According to Proposition~\ref{prop:Delta-DDAE-simplified}, the characteristic quasipolynomial of \eqref{DiffAlgeb} is the function $\Delta\colon \mathbb C \to \mathbb C$ given by
\begin{equation*}
\Delta(s) = s + a - k e^{-s\tau} (d s + a d + b c), \qquad s \in \mathbb C.
\end{equation*}
Note that $\Delta$ is a quasipolynomial of degree $3$.

Applying Theorem~\ref{thm:main-MID}, Proposition~\ref{prop:spectrum-neutral}, and Corollary~\ref{coro:DDAE} to this problem, we deduce that a given $s_0 \in \mathbb R$ is a root of multiplicity $3$ of $\Delta$ if and only if the coefficients $a, b, c, d, k$, the number $s_0$, and the delay $\tau$ satisfy the relations
\begin{equation}
\label{DiffAlgebConditions}
a = -s_0 - \frac{2}{\tau}, \qquad d k = -e^{s_0 \tau}, \qquad b c k = -\frac{4}{\tau} e^{s_0 \tau}.
\end{equation}
In this case, we have $\Real s = s_0$ for every $s \in \mathbb C \setminus \{s_0\}$ such that $\Delta(s) = 0$, and, more precisely, the set of roots of $\Delta$ is $\left\{s_0 + i \frac{2}{\tau} \xi \suchthat \tan \xi = \xi\right\}$. Moreover, under \eqref{DiffAlgebConditions}, the trivial solution of \eqref{DiffAlgeb} is exponentially stable if and only if $s_0 < 0$ or, equivalently, $a > -\frac{2}{\tau}$.

\subsection{Stabilization of a delay differential-algebraic equation}
\label{sec:appli-DDAE}

We now consider another example from \cite{Mazanti2021Effects}, which consists in the stabilization of a DDAE of the form \eqref{eq:DDAE-MID} with $n = 1$ and $m = 2$. More precisely, consider the DDAE
\begin{equation}
\label{AppliSystem}
\left\{
\begin{aligned}
x^\prime(t) + a x(t) & = b k_1 y_1(t - \tau) + b k_2 y_2(t - \tau), \\
y_1(t) & = c_1 x(t) + d_1 k_1 y_1(t - \tau) + d_1 k_2 y_2(t - \tau), \\
y_2(t) & = c_2 x(t) + d_2 k_1 y_1(t - \tau) + d_2 k_2 y_2(t - \tau), \\
\end{aligned}
\right.
\end{equation}
where $x(t) \in \mathbb R$, $y_1(t) \in \mathbb R$, $y_2(t) \in \mathbb R$, $\tau > 0$, and $a, b, c_1, c_2, d_1, d_2, k_1, k_2$ are real coefficients. From Proposition~\ref{prop:Delta-DDAE-simplified}, its characteristic quasipolynomial is the function $\Delta\colon \mathbb C \to \mathbb C$ given by
\begin{equation*}
\Delta(s) = s + a - e^{-s\tau}\bigl((d_1 k_1 + d_2 k_2)s + k_1(a d_1 + b c_1) + k_2(a d_2 + b c_2)\bigr), \qquad s \in \mathbb C,
\end{equation*}
which is of degree $3$.

From Theorem~\ref{thm:main-MID}, Proposition~\ref{prop:spectrum-neutral}, and Corollary~\ref{coro:DDAE}, we deduce that a given $s_0 \in \mathbb R$ is a root of multiplicity $3$ of $\Delta$ if and only if the coefficients $a, b, c_1, c_2, d_1, d_2, k_1, k_2$, the number $s_0$, and the delay $\tau$ satisfy the relations
\begin{equation}
\label{AppliConditions}
\begin{aligned}
a & = -s_0 - \frac{2}{\tau}, \\
d_1 k_1 + d_2 k_2 & = - e^{s_0 \tau}, \\
(a d_1 + b c_1) k_1 + (a d_2 + b c_2) k_2 & = \left(s_0 - \frac{2}{\tau}\right) e^{s_0 \tau}.
\end{aligned}
\end{equation}
Under \eqref{AppliConditions}, the set of roots of $\Delta$ is $\left\{s_0 + i \frac{2}{\tau} \xi \suchthat \tan \xi = \xi\right\}$ and the trivial solution of \eqref{DiffAlgeb} is exponentially stable if and only if $s_0 < 0$ or, equivalently, $a > -\frac{2}{\tau}$.

Let us now assume that \eqref{AppliSystem} comes from a system of the form \eqref{eq:DDAE-open-loop} with $n = 1$, $m = 2$ and fixed coefficients $a, b, c_1, c_2, d_1, d_2$, and that the control goal is to select a feedback control $u$ of the form \eqref{eq:control-DDAE} (i.e., to select the parameters $k_1$ and $k_2$ in \eqref{eq:control-DDAE}) in order to render the trivial solution of the closed-loop system \eqref{AppliSystem} exponentially stable. Note that the first equation in \eqref{AppliConditions} can be rewritten as $s_0 = -a - \frac{2}{\tau}$, hence one may stabilize the system by placing a negative real root of maximal multiplicity $s_0$ only if $a \geq 0$ or $\tau < \frac{2}{\abs{a}}$. Setting $\delta_i = a d_i + b c_i$ for $i \in \{1, 2\}$, the second and third equations of \eqref{AppliConditions} can be rewritten as
\begin{equation}
\label{SystK}
\begin{pmatrix}d_1 & d_2 \\ \delta_1 & \delta_2 \\\end{pmatrix} \begin{pmatrix}k_1 \\ k_2 \\\end{pmatrix} = -\begin{pmatrix}1 \\ a + \frac{4}{\tau}\end{pmatrix} e^{-a \tau - 2}.
\end{equation}
This system admits at least one solution if and only if $d_1 \delta_2 \neq d_2 \delta_1$ or $\delta_i = \left(a + \frac{4}{\tau}\right) d_i$ for $i \in \{1, 2\}$, with exactly one solution in the first case and infinitely many solutions in the second case.

Hence, if $a \geq 0$ or $\tau < \frac{2}{\abs{a}}$, and assuming in addition that either $d_1 \delta_2 \neq d_2 \delta_1$ or $\delta_i = \left(a + \frac{4}{\tau}\right) d_i$ for $i \in \{1, 2\}$, there exist real parameters $k_1$, $k_2$ such that \eqref{AppliSystem} is exponentially stable, with exponential decay rate $-a - \frac{2}{\tau}$, the parameters $k_1$ and $k_2$ being solutions of the linear system \eqref{SystK}.

\begin{figure}[ht]
\centering
\resizebox{0.5\textwidth}{!}{\input{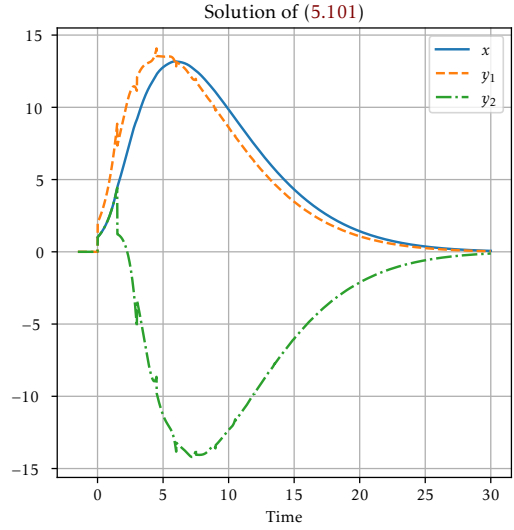}}
\caption{A solution of \eqref{AppliSystem} with the parameters from Section~\ref{sec:appli-DDAE}}
\label{AppliFigure}
\end{figure}

Let us consider, as an example, the case $a = -1$, $b = 1$, $c_1 = 2$, $c_2 = 1$, $d_1 = 1$, $d_2 = 2$, and $\tau = \frac{3}{2}$. Since $a = -1$, the corresponding open-loop system \eqref{eq:control-DDAE} is unstable when no control is applied. Note that the inequality $\tau < \frac{2}{\abs{a}}$ is satisfied and that the linear system \eqref{SystK} admits a unique solution, given by $k_1 \approx -0.87610$ and $k_2 \approx 0.13478$. With this choice of $k_1$ and $k_2$, the trivial solution of \eqref{AppliSystem} is exponentially stable, with exponential decay rate $-a - \frac{2}{\tau} = -\frac{1}{3}$. Figure~\ref{AppliFigure} presents a numerical simulation of the solution of \eqref{AppliSystem} with these parameters and with initial condition given by $x(0) = 1$ and $y_1(t) = y_2(t) = 0$ for $t < 0$. One observes, as expected, that solutions converge to $0$ exponentially.

\subsection{Improving the decay rate of a second-order control system without instantaneous velocity observation}
\label{sec:second-order-damped-oscillator}

We now present an application example of Theorem~\ref{thm:main-MID} extracted from \cite{Mazanti2021Multiplicity}. Consider the stable second-order control system
\begin{equation}
\label{Ex1SecondOrderSystem}
y^{\prime\prime}(t) + 2 \zeta \omega y^\prime(t) + \omega^2 y(t) = u(t),
\end{equation}
where $\zeta, \omega$ are positive real numbers and $u(t)$ is a control input. Under no control, the characteristic polynomial of this equation is $\Delta_0(s) = s^2 + 2 \zeta \omega s + \omega^2$, whose roots are $s_\pm = -\zeta\omega \pm i \omega \sqrt{1 - \zeta^2}$ if $0 < \zeta \leq 1$ or $s_{\pm} = -\zeta\omega \pm \omega \sqrt{\zeta^2 - 1}$ if $\zeta \geq 1$. Hence, the spectral abscissa $\gamma_0$ of $\Delta$ is given by
\[
\gamma_0 = \begin{dcases*}
-\zeta\omega, & if $0 < \zeta \leq 1$, \\
-\omega \left[\zeta - \sqrt{\zeta^2 - 1}\right], & if $\zeta \geq 1$.
\end{dcases*}
\]
In particular, $\gamma_0 < 0$ and the system is exponentially stable.

A classical problem in control theory is to choose the control $u$ in feedback form in order to improve the stability properties of \eqref{Ex1SecondOrderSystem}, and more precisely to decrease the value of the spectral abscissa $\gamma_0$. If one may choose $u(t)$ as a function of both instantaneous measures $y(t)$ and $y^\prime(t)$, then the linear feedback
\begin{equation}
\label{Ex1FullFeedback}
u(t) = -a_1 y^\prime(t) - a_0 y(t)
\end{equation}
for some real parameters $a_0, a_1$ yields the closed-loop equation
\begin{equation}
\label{Ex1ClosedLoopFullFeedback}
y^{\prime\prime}(t) + (a_1 + 2 \zeta \omega) y^\prime(t) + (\omega^2 + a_0) y(t) = 0.
\end{equation}
Hence, by choosing the coefficients $a_0$, $a_1$ appropriately, one may place the two roots of the characteristic polynomial of \eqref{Ex1ClosedLoopFullFeedback} anywhere in the complex plane (as long as they are real or form a complex conjugate pair), which allows one to obtain arbitrary values for the corresponding spectral abscissa, and hence arbitrary exponential decay rates of solutions.

In order to implement the control law \eqref{Ex1FullFeedback} for a practical system, one must be able to obtain instantaneous measures of $y(t)$ and $y^\prime(t)$ in order to use them in an instantaneous computation of the control $u(t)$ to be applied to the system, which is not always possible. In particular, in some situations, the measure of $y(t)$ may be available sufficiently fast in order to be considered approximately instantaneous, but good approximations of the velocity $y^\prime(t)$ may require extra time, due either to a velocity estimation procedure from the signal $y(t)$ or to slower measure processes. In the sequel, we explore two alternative feedback laws for \eqref{Ex1SecondOrderSystem} that one may use when not all instantaneous measures are possible.

Assuming that $y(t)$ can be measured instantaneously but not $y^\prime(t)$, one may propose the feedback law
\[
u(t) = - a_0 y(t)
\]
for \eqref{Ex1SecondOrderSystem}. In this case, the closed-loop system becomes
\[
y^{\prime\prime}(t) + 2 \zeta \omega y^\prime(t) + (\omega^2 + a_0) y(t) = 0,
\]
whose characteristic equation $\Delta_p(s) = s^2 + 2 \zeta\omega s + \omega^2 + a_0$ admits as roots $s_\pm = -\zeta\omega \pm i \omega \sqrt{1 + a_0 - \zeta^2}$ if $a_0 \geq \zeta^2 - 1$ and $s_{\pm} = -\zeta\omega \pm \omega \sqrt{\zeta^2 - 1 - a_0}$ if $a_0 \leq \zeta^2 - 1$. The spectral abscissa $\gamma_p$ of $\Delta_p$ is then
\[
\gamma_p = \begin{dcases*}
-\zeta\omega, & if $a_0 \geq \zeta^2 - 1$, \\
-\omega\left[\zeta - \sqrt{\zeta^2 - 1 - a_0}\right], & if $a_0 \leq \zeta^2 - 1$.
\end{dcases*}
\]
The spectral abscissa $\gamma_p$ can be minimized by choosing any $a_0 \in [\zeta^2 - 1, +\infty)$, in which case $\gamma_p = - \zeta\omega$. This coincides with $\gamma_0$ if $0 < \zeta \leq 1$ and improves $\gamma_0$ only in the case $\zeta > 1$.

The spectral abscissa $\gamma_p$ can be improved by considering a delayed feedback law when the velocity $y^\prime(t)$ is available for measure after some delay $\tau > 0$. Consider the control law
\begin{equation}
\label{Ex1Feedback}
u(t) = - a_0 y(t) - \alpha_1 y^\prime(t - \tau) - \alpha_0 y(t - \tau)
\end{equation}
in \eqref{Ex1SecondOrderSystem}, which yields the closed-loop system
\begin{equation}
\label{Ex1DelayedFeedback}
y^{\prime\prime}(t) + 2 \zeta \omega y^\prime(t) + (\omega^2 + a_0) y(t) + \alpha_1 y^\prime(t - \tau) + \alpha_0 y(t - \tau) = 0.
\end{equation}
Equation~\eqref{Ex1DelayedFeedback} is of the form \eqref{eq:Delta-DDE} and its characteristic quasipolynomial is
\[
\Delta(s) = s^2 + 2 \zeta \omega s + \omega^2 + a_0 + e^{- s \tau} (\alpha_1 s + \alpha_0).
\]
The conditions \eqref{Coeffs} on the coefficients of $\Delta$ are satisfied for some $s_0 \in \mathbb R$ if and only if
\begin{align*}
2 \zeta \omega & = -\frac{4}{\tau} - 2 s_0, & \omega^2 + a_0 & = \frac{6}{\tau^2} + \frac{4}{\tau}s_0 + s_0^2, \\
\alpha_1 & = -\frac{2}{\tau} e^{s_0 \tau}, & \alpha_0 & = \frac{2}{\tau} e^{s_0 \tau} \left(s_0 - \frac{3}{\tau}\right).
\end{align*}
Hence, if one chooses $s_0, a_0, \alpha_1, \alpha_0$ as
\begin{align*}
s_0 & = - \zeta \omega - \frac{2}{\tau}, & a_0 & = \frac{6}{\tau^2} + \frac{4}{\tau}s_0 + s_0^2 - \omega^2, \\
\alpha_1 & = -\frac{2}{\tau} e^{s_0 \tau}, & \alpha_0 & = \frac{2}{\tau} e^{s_0 \tau} \left(s_0 - \frac{3}{\tau}\right),
\end{align*}
then \eqref{Coeffs} is satisfied and, by Theorem~\ref{thm:main-MID}, $s_0$ is a strictly dominant root of $\Delta$ of multiplicity $4$. In particular, the spectral abscissa $\gamma$ of $\Delta$ is $\gamma = s_0 = - \zeta\omega - \frac{2}{\tau}$, which is strictly smaller than both $\gamma_0$ and $\gamma_p$. Hence, even when the instantaneous velocity $y^\prime(t)$ is not available for measure, one may design a control law depending on the delayed velocity $y^\prime(t - \tau)$ that stabilizes the closed-loop system with a spectral abscissa strictly better than both the open-loop one, $\gamma_0$, and the one obtained with only a position feedback, $\gamma_p$.

\begin{figure}[ht]
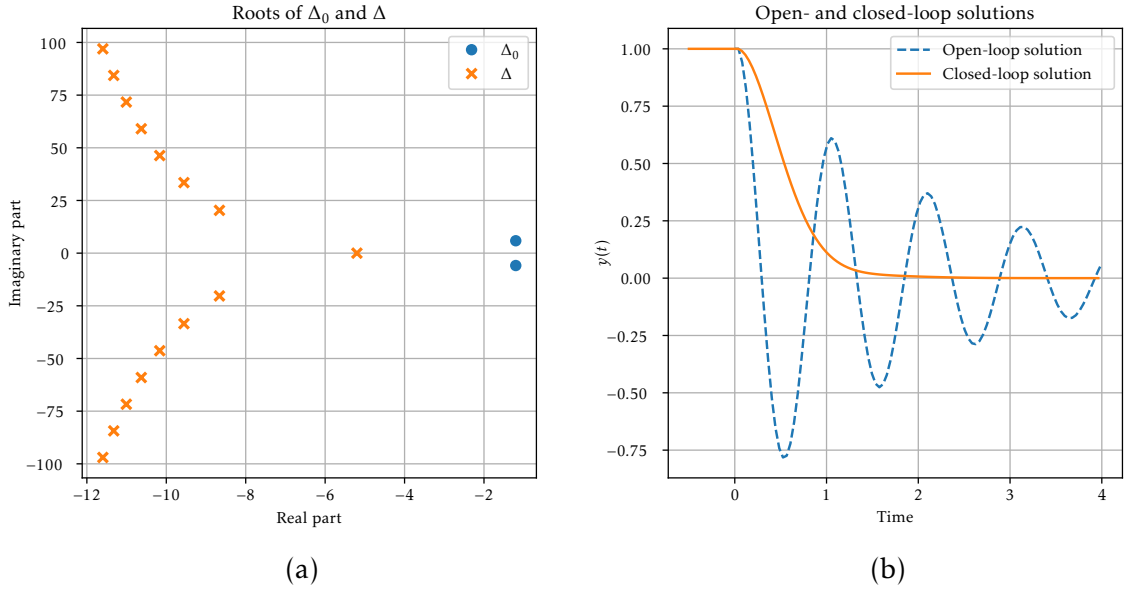

\centering
\begin{tabular}{@{} >{\centering} m{0.5\textwidth} @{} >{\centering} m{0.5\textwidth} @{}}
\resizebox{0.5\textwidth}{!}{\input{Figures/Expl4-roots.pgf}} & \resizebox{0.5\textwidth}{!}{\input{Figures/Expl4-sols.pgf}} \tabularnewline
(a) & (b) \tabularnewline
\end{tabular}
\caption[Comparison between the open-loop system \eqref{Ex1SecondOrderSystem} and the closed-loop system with delayed velocity \eqref{Ex1DelayedFeedback}]{Comparison between the open-loop system \eqref{Ex1SecondOrderSystem} with $u = 0$ and the closed-loop system with delayed velocity \eqref{Ex1DelayedFeedback} for $\zeta = 0.2$, $\omega = 6$, and a delay $\tau = 0.5$. (a) Roots of $\Delta_0$ (circles) and $\Delta$ (crosses). (b) Solution of \eqref{Ex1SecondOrderSystem} with $y(0) = 1$, $y^\prime(0) = 0$, $u(t) = 0$ (dashed line) and solution of \eqref{Ex1DelayedFeedback} with $y(t) = 1$, $y^\prime(t) = 0$ for $t \leq 0$ (solid line)}
\label{FigOLCL}
\end{figure}

A comparison between the open-loop system corresponding to \eqref{Ex1SecondOrderSystem} with $u(t) = 0$ and the closed-loop system \eqref{Ex1DelayedFeedback} with the feedback law \eqref{Ex1Feedback} is provided in Figure~\ref{FigOLCL} for the choice of parameters $\zeta = 0.2$, $\omega = 6$, and $\tau = 0.5$. In this case, the spectral abscissa $\gamma_0$ is $\gamma_0 = -1.2$, the root of multiplicity $4$ of $\Delta$ is $s_0 = -5.2$, the spectral abscissa $\gamma$ of $\Delta$ is also $\gamma = s_0 = -5.2$, and the coefficients of the feedback law \eqref{Ex1Feedback} ensuring that $s_0$ is a root of multiplicity $4$ are given by $a_0 = -26.56$, $\alpha_1 \approx -0.2971$, and $\alpha_0 \approx -3.327$. Figure~\ref{FigOLCL}(a) shows that, even though the spectrum of \eqref{Ex1DelayedFeedback} is infinite, its spectral abscissa is indeed smaller than that of \eqref{Ex1SecondOrderSystem}, which has a finite spectrum. The behavior of solutions is illustrated in Figure~\ref{FigOLCL}(b), in which the solution of \eqref{Ex1SecondOrderSystem} with initial condition $y(0) = 1$ and $y^\prime(0) = 0$ converges to $0$ much slower than the solution of the closed-loop system \eqref{Ex1DelayedFeedback} with initial condition defined by $y(t) = 1$ and $y^\prime(t) = 0$ for $t \leq 0$.

\subsection{Control of a transonic flow in a wind tunnel}

We consider in this section the application of Theorem~\ref{thm:main-MID} to the control of a transonic flow in a wind tunnel, which is presented in \cite{Mazanti2021Multiplicity}. The analysis of transonic flows is a challenging problem in compressible fluid dynamics, since a full model of the flow would involve considering the Navier--Stokes equations in a three-dimensional domain and boundary controls for temperature and pressure regulation. A first simplification presented in \cite{Tripp1983Development} assumes that the tunnel is a one-dimensional tube of varying cross-sectional area and that the flow is uniform across every cross section, which allows one to described the flow through a coupled system of nonlinear partial differential equations in one space dimension.

A further simplified model was presented in \cite{Armstrong1981Application} in order to analyze the response of the Mach number of the flow to changes in the guide vane angle. Instead of using partial differential equations, propagation phenomena are modeled in \cite{Armstrong1981Application} through a time delay, leading to the time-delay system
\begin{equation}
\label{EqApplication}
\left\{
\begin{aligned}
\kappa m^\prime(t) + m(t) & = k \vartheta(t - \tau_0), \\
\vartheta^{\prime\prime}(t) + 2 \zeta \omega \vartheta^\prime(t) + \omega^2 \vartheta(t) & = \omega^2 u(t),
\end{aligned}
\right.
\end{equation}
in which $m$, $\vartheta$, and $u$ represent perturbations of the Mach number of the flow, the guide vane angle, and the input of the guide vane actuator, respectively, with respect to steady-state values. The parameters $\kappa$ and $k$ depend on the steady-state operating point and are assumed to be constant as long as $m$, $\vartheta$, and $u$ remain small, and satisfy $\kappa > 0$ and $k < 0$. The parameters $\zeta \in (0, 1)$ and $\omega > 0$ come from the design of the guide vane angle actuator and are thus independent of the operating point. The time-delay $\tau_0$ is assumed to depend only on the temperature of the flow. In the absence of control (i.e., when $u(t) = 0$ for all $t$), the open-loop system \eqref{EqApplication} is exponentially stable.

The design of exponentially stabilizing feedback laws for \eqref{EqApplication} improving its stability properties has been considered, for instance, in \cite{Manitius1984Feedback}, in which the author designs controllers depending linearly on the state variables and on integrals of some of the state variables over an interval of length equal to the time-delay. This leads to a closed-loop system with finite spectrum. However, the practical implementation of these controllers is difficult due to the integral terms, which motivates the research for control laws with reduced implementation complexity. In general, one may not expect to obtain reduced-complexity controllers guaranteeing a closed-loop system with finite spectrum, but the multiplicity-induced dominance techniques developed in Section~\ref{sec:MID-maximal} are a good candidate for proposing reduced-complexity controllers while dealing with systems with infinite spectrum.

The goal of this section is to illustrate how Theorem~\ref{thm:main-MID} can be used to obtain a feedback controller for \eqref{EqApplication} improving its open-loop characteristics, with a reduced complexity with respect to the controller proposed in \cite{Manitius1984Feedback}. We assume that, up to redesigning the guide vane angle actuator, one may choose the values of its parameters $\zeta \in (0, 1)$ and $\omega > 0$. The control law we propose is
\begin{equation}
\label{ApplicationControlLaw}
u(t) = \beta_0 m(t - \tau_1) + \beta_1 m^\prime(t - \tau_1) + \beta_2 m^{\prime\prime}(t - \tau_1),
\end{equation}
where $\tau_1 > 0$ should be greater than or equal to the time-delay corresponding to measuring $m$ and its first two derivatives. Inserting \eqref{ApplicationControlLaw} into \eqref{EqApplication}, one obtains that the closed-loop characteristic quasipolynomial $\Delta$ is given by
\begin{multline}
\label{ApplicationDelta}
\frac{\Delta(s)}{\kappa} = s^{3} + \left(2 \omega \zeta + \frac{1}{\kappa}\right) s^{2} + \left(\omega^{2} + \frac{2 \omega \zeta}{\kappa}\right) s + \frac{\omega^{2}}{\kappa} \\ {} + \left(- \frac{\beta_{2} k \omega^{2}}{\kappa} s^{2} - \frac{\beta_{1} k \omega^{2}}{\kappa} s - \frac{\beta_{0} k \omega^{2}}{\kappa}\right) e^{- s \left(\tau_{0} + \tau_{1}\right)},
\end{multline}
where the division by $\kappa$ is performed in order to obtain a quasipolynomial under the form \eqref{eq:Delta-DDE}, for which the coefficient of the term of highest degree $s^3$ is $1$. As a consequence of Theorem~\ref{thm:main-MID}, one gets the following result.

\begin{proposition}
\label{PropExpl}
Let $r_0 = -3 - \sqrt[3]{9} + \sqrt[3]{3}$ and $s_0 \in \mathbb R$. Then $s_0$ is a root of maximal multiplicity $6$ of the quasipolynomial $\Delta$ from \eqref{ApplicationDelta} if and only if
\begin{subequations}
\label{ExplConditionsMaximalMultiplicity}
\begin{align}
s_0 & = \frac{r_0}{\tau} - \frac{1}{\kappa}, \label{PropS0} \displaybreak[0] \\
\omega & = \sqrt{- \kappa \left(s_{0}^{3} + \frac{9 s_{0}^{2}}{\tau} + \frac{36 s_{0}}{\tau^{2}} + \frac{60}{\tau^{3}}\right)}, \label{PropOmega} \displaybreak[0] \\
\zeta & = - \frac{3 s_{0}}{2 \omega} - \frac{9}{2 \omega \tau} - \frac{1}{2 \omega \kappa}, \label{PropZeta} \displaybreak[0] \\
\beta_{0} & = - \frac{3 \kappa \left(s_{0}^{2} \tau^{2} - 8 s_{0} \tau + 20\right) e^{s_{0} \tau}}{k \omega^{2} \tau^{3}}, \label{PropBeta0} \displaybreak[0] \\
\beta_{1} & = \frac{6 \kappa \left(s_{0} \tau - 4\right) e^{s_{0} \tau}}{k \omega^{2} \tau^{2}}, \displaybreak[0] \\
\beta_{2} & = - \frac{3 \kappa e^{s_{0} \tau}}{k \omega^{2} \tau}, \label{PropBeta2}
\end{align}
\end{subequations}
where $\tau = \tau_0 + \tau_1$. Moreover, if \eqref{ExplConditionsMaximalMultiplicity} is satisfied, then $s_0 < 0$, $s_0$ is a strictly dominant root of $\Delta$, and the trivial solution of \eqref{EqApplication} with the control law \eqref{ApplicationControlLaw} is exponentially stable. 
\end{proposition}

\begin{proof}
By Theorem~\ref{thm:main-MID}, a real number $s_0$ is a root of maximal multiplicity $6$ of $\Delta$ if and only if
\begin{subequations}
\begin{align}
2 \omega \zeta + \frac{1}{\kappa} & = - 3 s_{0} - \frac{9}{\tau}, \label{Syst1ZetaOmega} \displaybreak[0] \\
\omega^{2} + \frac{2 \omega \zeta}{\kappa} & = 3 s_{0}^{2} + \frac{18 s_{0}}{\tau} + \frac{36}{\tau^{2}}, \displaybreak[0] \\
\frac{\omega^{2}}{\kappa} & = - s_{0}^{3} - \frac{9 s_{0}^{2}}{\tau} - \frac{36 s_{0}}{\tau^{2}} - \frac{60}{\tau^{3}}, \label{Syst1OmegaSquare} \displaybreak[0] \\
- \frac{\beta_{2} k \omega^{2}}{\kappa} & = \frac{3 e^{s_{0} \tau}}{\tau}, \label{Syst1Beta2} \displaybreak[0] \\
- \frac{\beta_{1} k \omega^{2}}{\kappa} & = \frac{6 \left(- s_{0} \tau + 4\right) e^{s_{0} \tau}}{\tau^{2}}, \displaybreak[0] \\
- \frac{\beta_{0} k \omega^{2}}{\kappa} & = \frac{3 \left(s_{0}^{2} \tau^{2} - 8 s_{0} \tau + 20\right) e^{s_{0} \tau}}{\tau^{3}}. \label{Syst1Beta0}
\end{align}
\end{subequations}
Clearly, equations \eqref{Syst1Beta2}--\eqref{Syst1Beta0} are equivalent to \eqref{PropBeta0}--\eqref{PropBeta2}. Let us now prove that equations \eqref{Syst1ZetaOmega}--\eqref{Syst1OmegaSquare} are equivalent to \eqref{PropS0}--\eqref{PropZeta}.

Let $Q$ be the polynomial given by $Q(X) = X^3 + 9 X^2 + 36 X + 60$ and notice that $r_0$ is its unique real root. Equations~\eqref{Syst1ZetaOmega}--\eqref{Syst1OmegaSquare} can be rewritten as
\begin{subequations}
\label{ApplicationLast3}
\begin{align}
2 \omega \zeta \tau + \frac{\tau}{\kappa} & = - \frac{1}{2} Q^{\prime\prime}(s_0 \tau), \displaybreak[0] \\
\omega^{2} \tau^2 + \frac{2 \omega \zeta}{\kappa} \tau^2 & = Q^\prime(s_0 \tau), \displaybreak[0] \\
\frac{\omega^{2}}{\kappa} \tau^3 & = -Q(s_0 \tau),
\end{align}
\end{subequations}
which is equivalent to
\begin{gather*}
\begin{aligned}
\omega^{2} \tau^2 & = -\frac{\kappa}{\tau} Q(s_0 \tau), \\
2 \omega \zeta \tau & = \frac{\kappa^2}{\tau^2} Q(s_0 \tau) + \frac{\kappa}{\tau} Q^\prime(s_0 \tau),
\end{aligned} \displaybreak[0] \\
\frac{\kappa^3}{\tau^3} Q(s_0 \tau) + \frac{\kappa^2}{\tau^2} Q^\prime(s_0 \tau) + \frac{1}{2} \frac{\kappa}{\tau} Q^{\prime\prime}(s_0 \tau) + \frac{1}{6} Q^{\prime\prime\prime}(s_0 \tau) = 0.
\end{gather*}
Note that the last of the above equalities can be equivalently rewritten as
\begin{equation}
\label{InterestingEquality}
Q(s_0 \tau) + Q^\prime(s_0 \tau) \frac{\tau}{\kappa} + \frac{1}{2} Q^{\prime\prime}(s_0 \tau) \frac{\tau^2}{\kappa^2} + \frac{1}{6} Q^{\prime\prime\prime}(s_0 \tau) \frac{\tau^3}{\kappa^3} = 0.
\end{equation}
The right-hand side of the above equality corresponds to the Taylor expansion of order $3$ of $x \mapsto Q(s_0 \tau + x)$ around the origin and evaluated at $x = \frac{\tau}{\kappa}$. Since $Q$ is a polynomial of degree $3$, this Taylor expansion coincides with $Q\left(s_0 \tau + \frac{\tau}{\kappa}\right)$, and thus \eqref{InterestingEquality} is equivalent to $Q\left(s_0 \tau + \frac{\tau}{\kappa}\right) = 0$. Since the only real root of $Q$ is $r_0$, \eqref{InterestingEquality} is thus equivalent to
\[
s_0 = \frac{r_0}{\tau} - \frac{1}{\kappa}.
\]
Hence \eqref{ApplicationLast3} is equivalent to
\begin{subequations}
\begin{align}
s_0 & = \frac{r_0}{\tau} - \frac{1}{\kappa}, \label{s0} \\
\omega^{2} & = -\frac{\kappa}{\tau^3} Q(s_0 \tau), \label{OmegaSquare} \\
\zeta & = - \frac{3 s_{0}}{2 \omega} - \frac{9}{2 \omega \tau} - \frac{1}{2 \omega \kappa}.
\end{align}
\end{subequations}
One has $Q^\prime(X) = 3 (X^2 + 6 X + 12)$, which is positive for every $X \in \mathbb R$. Hence $Q$ is strictly increasing and, in particular, under \eqref{s0}, one has $Q(s_0 \tau) < Q\left(s_0 \tau + \frac{\tau}{\kappa}\right) = 0$, proving that the right-hand side of \eqref{OmegaSquare} is positive. Thus, if \eqref{s0} is satisfied, there exists a unique $\omega > 0$ such that \eqref{OmegaSquare} holds, and this $\omega$ is precisely the one from \eqref{PropOmega}. This concludes the proof that \eqref{Syst1ZetaOmega}--\eqref{Syst1OmegaSquare} is equivalent to \eqref{PropS0}--\eqref{PropZeta}.

Hence, $s_0$ is a root of maximal multiplicity of $\Delta$ if and only if \eqref{ExplConditionsMaximalMultiplicity} holds. The other assertions of the statement of the proposition follow immediately from \eqref{PropS0}, the fact that $r_0 < 0$, and Theorem~\ref{thm:main-MID}.
\end{proof}

\begin{remark}
We stress that the proof of Proposition~\ref{PropExpl} guarantees that, under \eqref{PropS0}, the argument of the square root in \eqref{PropOmega} is positive. Thus $\omega > 0$ is well-defined by \eqref{PropOmega}.
\end{remark}

\begin{remark}
For typical values of the parameters $k$, $\kappa$, $\tau_0$, $\zeta$, and $\omega$ (see \cite{Armstrong1981Application, Manitius1984Feedback}) and under no control, the dynamics of the deviations $\vartheta$ of the guide vane angle are considerably faster than those of the deviations $m$ of the Mach number, and thus one may expect \eqref{EqApplication} to be approximated by the first-order differential equation
\begin{equation}
\label{SimplifiedM}
\kappa m^\prime(t) + m(t) = 0.
\end{equation}
The characteristic polynomial of this equation admits a unique root at $s^*=-\frac{1}{\kappa}$. Since $r_0 < 0$, it follows from \eqref{PropS0} that the strictly dominant root $s_0$ of the quasipolynomial $\Delta$ from \eqref{ApplicationDelta} satisfies $s_0 < s^*$. Hence, the deviation of the Mach number $m$ has a faster decay rate in \eqref{EqApplication} with the control law \eqref{ApplicationControlLaw} than in the simplified open-loop system \eqref{SimplifiedM}.
\end{remark}

Thanks to Proposition~\ref{PropExpl}, given $k$, $\kappa$, $\tau_0$, and $\tau_1$, one may easily compute from \eqref{ExplConditionsMaximalMultiplicity} the unique possible real root of multiplicity $6$ of $\Delta$, $s_0$, as well as the parameters $\zeta$ and $\omega$ of the guide vane angle actuator and the parameters $\beta_0$, $\beta_1$, and $\beta_2$ of the control law \eqref{ApplicationControlLaw} ensuring that $s_0$ is indeed a root of multiplicity $6$ of the quasipolynomial $\Delta$.

As a numerical example, we consider the same system parameters as in \cite{Manitius1984Feedback}, given by $\kappa = \SI{1.964}{\second}$, $k = \SI{-0.67036}{\per\radian}$, and $\tau_0 = \SI{.33}{\second}$, corresponding to a linearization around the steady state with Mach number $\num{0.9}$ and air temperature $\SI{166}{\kelvin}$. The parameters $\zeta$, $\omega$, $\beta_2$, $\beta_1$, and $\beta_0$ computed according to \eqref{ExplConditionsMaximalMultiplicity} for $\tau_1 = \SI{0.33}{\second}$ and $\tau_1 = \SI{0.70}{\second}$ are provided in Table~\ref{TabParams}, together with the corresponding multiple root $s_0$.

\begin{table}[ht]
\centering
\caption[Parameters of the guide vane angle actuator and of the control law \eqref{ApplicationControlLaw}]{Parameters of the guide vane angle actuator and of the control law \eqref{ApplicationControlLaw} and location $s_0$ of the multiple root computed according to \eqref{ExplConditionsMaximalMultiplicity} for two different values of the delay $\tau_1$}
\label{TabParams}
\begin{tabular}{ccccccc}
\toprule
$\tau_1$ [\si{\second}] & $s_0$ [\si{\per\second}] & $\zeta$ & $\omega$ [\si{\radian\per\second}] & $\beta_0$ & $\beta_1$ [\si{\second}] & $\beta_2$ [\si{\square\second}] \tabularnewline
\midrule
$\num{0.33}$ & $\num{-6.021}$ & $\num{0.3902}$ & $\num{5.020}$ & $\num{1.542}$ & $\num{0.2401}$ & $\num{0.00994}$ \tabularnewline
$\num{0.70}$ & $\num{-4.041}$ & $\num{0.4368}$ & $\num{3.292}$ & $\num{0.8161}$ & $\num{0.1943}$ & $\num{0.01226}$ \tabularnewline
\bottomrule
\end{tabular}
\end{table}

\begin{figure}[ht]
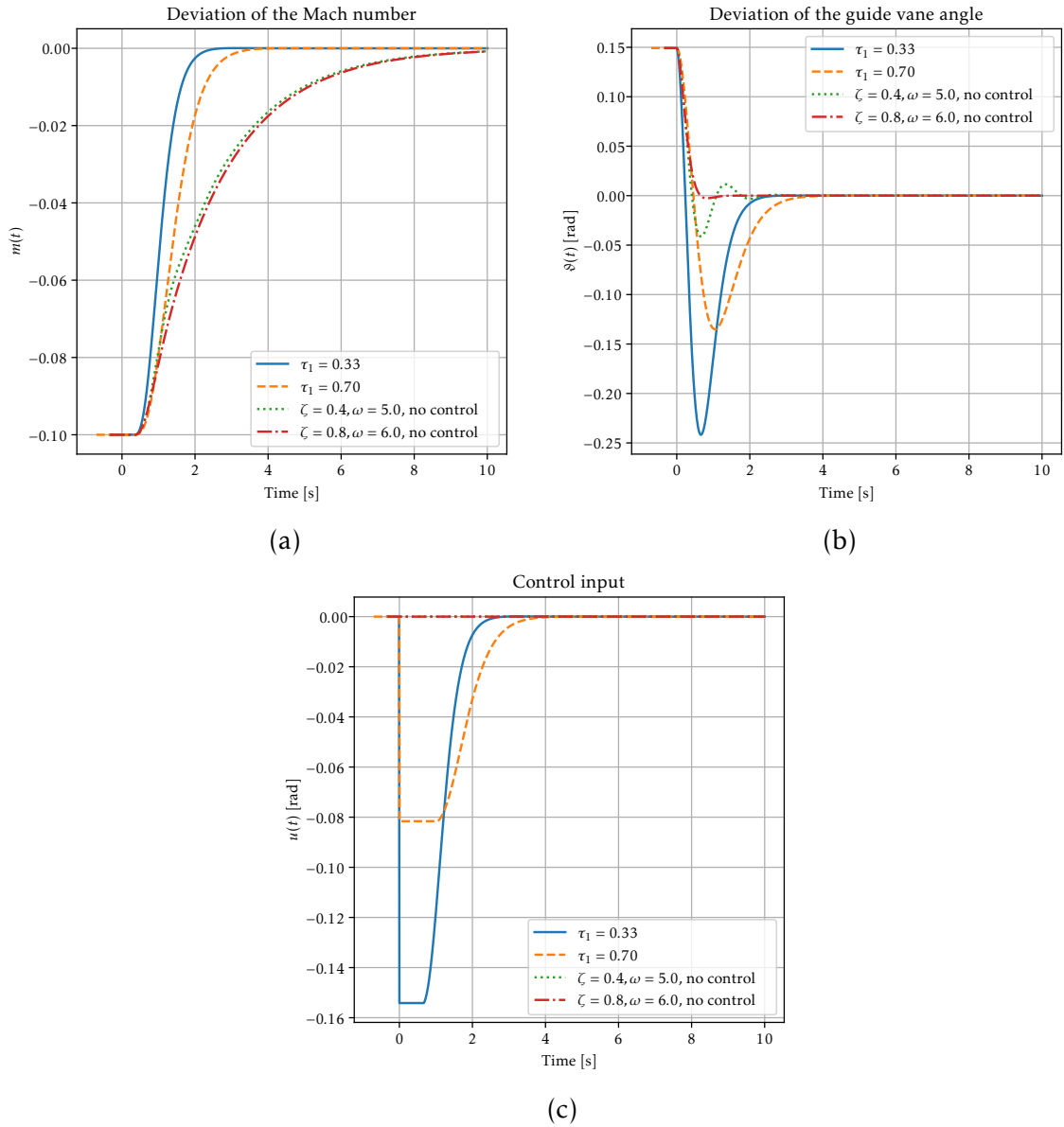

\centering

\begin{tabular}{@{} >{\centering} m{0.5\textwidth} @{} >{\centering} m{0.5\textwidth} @{}}
\resizebox{0.5\textwidth}{!}{\input{Figures/Expl5-mach.pgf}} & \resizebox{0.5\textwidth}{!}{\input{Figures/Expl5-angle.pgf}} \tabularnewline
(a) & (b) \tabularnewline
\end{tabular}

\resizebox{0.5\textwidth}{!}{\input{Figures/Expl5-control.pgf}}

(c)
\caption[Mach number and guide vane angle deviations and control input]{(a) Deviation $m(t)$ of the Mach number with respect to its steady-state value. (b) Deviation $\vartheta(t)$ of the guide vane angle with respect to its steady-state position. (c) Control input $u(t)$}
\label{FigApplication}
\end{figure}

Figure~\ref{FigApplication} provides the evolution of the deviation $m(t)$ of the Mach number with respect to its steady-state value, the deviation $\vartheta(t)$ of the guide vane with respect to its steady-state position, and the control input $u(t)$ in the case where the initial conditions are given by $m(t) = \num{-0.1}$, $\vartheta(t) = \frac{m(0)}{k} \approx \SI{0.14917}{\radian}$, and $\vartheta^\prime(t) = 0$ for $t \leq 0$. Notice that this choice of initial conditions, which is the same as in \cite{Manitius1984Feedback}, ensures that $m^\prime(t) = 0$ and $m(t) = \num{-0.1}$ for $t \in (0, \tau_0)$.

Together with the solutions of \eqref{EqApplication} under the control law \eqref{ApplicationControlLaw} with the parameters given in Table~\ref{TabParams} computed with $\tau_1 = \SI{0.33}{\second}$ (solid line) and $\tau_1 = \SI{0.70}{\second}$ (dashed line), we also represent in Figure~\ref{FigApplication}, for comparison, the solutions of \eqref{EqApplication} with no control for two sets of parameters $\zeta$ and $\omega$, $\zeta = \num{0.4}$ and $\omega = \SI{5}{\radian\per\second}$ (dotted line), and $\zeta = {0.8}$ and $\omega = \SI{6}{\radian\per\second}$ (dash-dotted line). The first set of parameters $\zeta$ and $\omega$ is an approximation of the ones from the first row of Table~\ref{TabParams}, while the second set of parameters $\zeta$ and $\omega$ is the one used in \cite{Manitius1984Feedback}.

Figure~\ref{FigApplication} shows that the control law \eqref{ApplicationControlLaw} acting on the guide vane angle dynamics from \eqref{EqApplication} allows one to improve the convergence rate of the deviations $m$ of the Mach number of the flow, as required. As previously remarked, one observes in Figure~\ref{FigApplication} that \eqref{EqApplication} is already exponentially stable under no control, but $m$ converges faster to $0$ with the control law \eqref{ApplicationControlLaw}. The convergence rate of $m$ increases as $\tau_1$ decreases since, as seen in \eqref{PropS0} and illustrated in Table~\ref{TabParams}, the dominant root $s_0$ of $\Delta$ moves to the left in the complex plane as $\tau_1$ decreases.

Since the control $u$ acts on the dynamics of the guide vane angle deviations $\vartheta$ in \eqref{EqApplication}, one might expect important differences in the dynamics of $\vartheta$ with respect to the situation with no control, which is indeed observed in Figure~\ref{FigApplication}(b). The magnitude of the deviations $\vartheta$ of the guide vane angle increases when the control $u$ from \eqref{ApplicationControlLaw} is applied, as one may see comparing the solid line in Figure~\ref{FigApplication}(b), corresponding to the control \eqref{ApplicationControlLaw} computed with $\tau_1 = \SI{0.33}{\second}$, and the dotted line in the same figure, corresponding to the case with no control and with the choice of parameters $\zeta = \num{0.4}$ and $\omega = \SI{5}{\radian\per\second}$, close to the ones obtained with $\tau_1 = \SI{0.33}{\second}$ in Table~\ref{TabParams}. We also observe, in Figure~\ref{FigApplication}(c), that, as $\tau_1$ decreases, the amplitude of the control input $u$ increases.

In conclusion, the control law \eqref{ApplicationControlLaw}, based on Theorem~\ref{thm:main-MID}, allows one to improve the convergence rate of deviations of the Mach number in the linearized model of a transonic flow. This control law only requires delayed measures of the Mach flow deviations $m$ and improves the convergence rate of $m$ even when the measurement delay is relatively large. However, one of its major drawbacks is that it also requires measurements of the derivatives $m^{\prime}$ and $m^{\prime\prime}$, which may introduce noise in practical applications.

An interesting open problem is whether similar techniques can be applied for other kinds of control laws, easier to implement in practice and avoiding the noisy derivative terms. A natural candidate would be to consider control laws based on measurements of $m$ with a certain number of different delays, which is motivated by discrete approximations of the derivatives in \eqref{ApplicationControlLaw}. However, the currently known theoretical results on the MID property for systems with several delays, such as the one from Section~\ref{sec:MID-two-delays}, still do not allow one to consider the feedback control of \eqref{EqApplication} with several delays in the feedback law.

\subsection{Exponential stabilization of the classical pendulum}

Following the example presented in \cite{Boussaada2022Generic}, we now consider the feedback stabilization of a friction-free classical pendulum. Let us consider the controlled classical pendulum described by the second-order differential equation
\begin{equation}
\label{InvPenNLin}
\theta^{\prime\prime}(t) + \frac{g}{L}\sin(\theta(t)) = u(t),
\end{equation}
where $\theta(t)$ stands for the angular displacement of the pendulum at time $t$ with respect to the stable equilibrium position, $L$ is the length of the pendulum, $g$ is the gravitational acceleration, and $u(t)$ represents the control law, which stems from an applied external torque.

When $u(t) = 0$ for all $t$, the trivial solution $\theta(t) = \theta^\prime(t) = 0$ for all $t$ is stable, but not exponentially stable. The problem we consider is to render this equilibrium locally asymptotically stable via the action of the external torque $u$. The linearization of \eqref{InvPenNLin} around this equilibrium is
\begin{equation}
\label{InvPenLin}
\theta^{\prime\prime}(t) + \frac{g}{L} \theta(t) = u(t).
\end{equation}
We focus in the sequel on the feedback exponential stabilization of the trivial solution of \eqref{InvPenLin}. Note that this system is of the form \eqref{Ex1SecondOrderSystem} with $\zeta = 0$ and $\omega^2 = \frac{g}{L}$.

Contrarily to what was done in Section~\ref{sec:second-order-damped-oscillator} for the second-order damped oscillator \eqref{Ex1SecondOrderSystem}, we will assume here that only delayed measures of $\theta(t)$ and $\theta^\prime(t)$ are available, but that the delay itself can be used as a control parameter. More precisely, we will consider a delayed proportional-derivative controller of the form
\begin{equation}
\label{DPD}
u(t) = -k_p \theta(t - \tau) - k_d \theta^\prime(t - \tau),
\end{equation}
where $k_p$ and $k_d$ are real parameters and $\tau > 0$ is a delay. The exponential stability of the linearized pendulum equation \eqref{InvPenLin} with the control law \eqref{DPD} can be studied through the location of the roots of the quasipolynomial 
\begin{equation}
\label{DeltaPendulum}
\Delta(s) = s^{2} + \frac{g}{L} + \left(k_{d} s + k_{p}\right) {e}^{-\tau  s},
\end{equation}
which is a quasipolynomial of the form \eqref{eq:Delta-DDE} of degree 4. We will consider in the sequel two possible strategies for stabilizing \eqref{InvPenLin}--\eqref{DPD}, according to whether one selects the parameters $k_d$, $k_p$, and $\tau$ of \eqref{DeltaPendulum} to place a root $s_0$ of maximal multiplicity $4$ or of multiplicity at least $3$.

\subsubsection{Root of maximal multiplicity \texorpdfstring{$4$}{4}}
\label{sec:pendulum-4}

By Theorem~\ref{thm:main-MID}, \eqref{DeltaPendulum} admits a root of multiplicity $4$ at some real number $s_0$ if and only if $s_0$, $g$, $L$, $k_d$, $k_p$, and $\tau$ satisfy
\begin{equation}
\label{IPQR}
\begin{aligned}
s_0 & = -\sqrt{\frac{2 g}{L}}, & \qquad \tau & = \sqrt{\frac{2 L}{g}}, \\
k_d & = - e^{-2} \sqrt{\frac{2 g}{L}}, & k_p & = - \frac{5 e^{-2} g}{L},
\end{aligned}
\end{equation}
and, in this case, $s_0$ is the spectral abscissa of the closed-loop system \eqref{InvPenLin}--\eqref{DPD}. Since $s_0 < 0$, this ensures the exponential stability of the trivial solution of \eqref{InvPenLin}--\eqref{DPD}. The main issue with the above method is that imposing that $s_0$ reaches its maximal multiplicity constrains $s_0$ and $\tau$ to have very specific values, described in \eqref{IPQR}.

\subsubsection{Root of multiplicity at least \texorpdfstring{$3$}{3}}

A possible strategy in order to allow additional freedom in the assignment of $s_0$ is to impose a lower multiplicity on $s_0$, in the spirit of Section~\ref{sec:MID-intermediate}. More precisely, for the quasipolynomial $\Delta$ given by \eqref{DeltaPendulum}, assuming that the free parameters are $k_d$ and $k_p$ and using the notation of Section~\ref{sec:iMID-EPF}, we select $\Lambda = 3$ and $I = \{3, 4\}$, and thus $\overline I = \{1, 2\}$ and we have
\begin{gather*}
\ell(s_0, \tau) = \begin{pmatrix}s_0 & 1 & s_0 e^{-s_0 \tau} & e^{-s_0 \tau}\end{pmatrix}, \displaybreak[0] \\
\ell_I(s_0, \tau) = \begin{pmatrix} (\tau^2 s_0 - 2 \tau) e^{-s_0 \tau} & \tau^2 e^{-s_0 \tau}\end{pmatrix}, \qquad
\ell_{\overline I}(s_0, \tau) = \begin{pmatrix}0 & 0 \end{pmatrix}, \displaybreak[0] \\
A(s_0, \tau) = \begin{pmatrix}
s_0 & 1 & s_0 e^{-s_0 \tau} & e^{-s_0 \tau} \\
1 & 0 & (-\tau s_0 + 1) e^{-s_0 \tau} & -\tau e^{-s_0 \tau} \\
\end{pmatrix}, \displaybreak[0] \\
A_I(s_0, \tau) = \begin{pmatrix}
s_0 e^{-s_0 \tau} & e^{-s_0 \tau} \\
(-\tau s_0 + 1) e^{-s_0 \tau} & -\tau e^{-s_0 \tau} \\
\end{pmatrix}, \qquad
A_{\overline I}(s_0, \tau) = \begin{pmatrix}
s_0 & 1 \\
1 & 0 \\
\end{pmatrix}, \displaybreak[0] \\
b(s_0) = \begin{pmatrix}
-s_0^2 \\
-2 s_0
\end{pmatrix}, \qquad
b_3 = -2, \qquad
X = \begin{pmatrix}0 \\ \frac{g}{L} \\ k_d \\ k_p \end{pmatrix}, \qquad
X_I  = \begin{pmatrix}k_d \\ k_p \end{pmatrix}, \qquad
X_{\overline I} = \begin{pmatrix}0 \\ \frac{g}{L} \end{pmatrix}.
\end{gather*}

By Proposition~\ref{prop:EPP-explicit}, the elimination-produced polynomial in this setting is
\[
\mathbb P(s_0, \tau) = \tau^2 s_0^2 + 4 \tau s_0 + \tau^2 \frac{g}{L} + 2,
\]
and thus, by Theorem~\ref{thm:charact-EPF}, $s_0$ is a root of multiplicity at least $3$ of $\Delta$ if and only if
\[
\begin{pmatrix}k_d \\ k_p \end{pmatrix} = -e^{2 s_0 \tau}\begin{pmatrix}
-\tau e^{-s_0 \tau} & -e^{-s_0 \tau} \\
(\tau s_0 - 1) e^{-s_0 \tau} & s_0 e^{-s_0 \tau} \\
\end{pmatrix} \left[\begin{pmatrix}
-s_0^2 \\
-2 s_0
\end{pmatrix} - \begin{pmatrix}
s_0 & 1 \\
1 & 0 \\
\end{pmatrix} \begin{pmatrix}0 \\ \frac{g}{L} \end{pmatrix}\right]
\]
and
\[
\tau^2 s_0^2 + 4 \tau s_0 + \tau^2 \frac{g}{L} + 2 = 0.
\]
Hence, $s_0$ is a root of multiplicity at least $3$ of $\Delta$ if and only if
\begin{subequations}
\begin{align}
k_d & = -e^{s_0 \tau} \left(\tau s_0^2 + 2 s_0 + \frac{g}{L} \tau\right), \label{eq:example-kd} \\
k_p & = e^{s_0 \tau} \left(\tau s_0^3 + s_0^2 + \tau \frac{g}{L} s_0 - \frac{g}{L}\right), \label{eq:example-kp}
\end{align}
\begin{equation}
\label{eq:example-EPP}
\tau^2 s_0^2 + 4 \tau s_0 + \tau^2 \frac{g}{L} + 2 = 0.
\end{equation}
\end{subequations}

Note that \eqref{eq:example-EPP} is satisfied if and only if
\begin{subequations}
\label{eq:conditions-mult-3}
\begin{equation}
\label{eq:example-EPP-2}
g \tau^2 \leq 2 L \quad \text{ and } \quad s_0 \in \left\{\frac{- 2 - \sqrt{2 - \frac{g \tau^2}{L}}}{\tau}, \frac{- 2 + \sqrt{2 - \frac{g \tau^2}{L}}}{\tau}\right\},
\end{equation}
and \eqref{eq:example-kd} and \eqref{eq:example-kp} can be further simplified to
\begin{equation}
\label{eq:example-kd-kp-2}
k_d = 2 e^{s_0 \tau} \left(s_0 + \frac{1}{\tau}\right), \qquad k_p = 2 e^{s_0 \tau} \left(\frac{g}{L} + \frac{5 s_0}{\tau} + \frac{3}{\tau^2}\right).
\end{equation}
\end{subequations}

While \eqref{eq:conditions-mult-3} provides a necessary and sufficient condition for $s_0$ to be a root of multiplicity at least $3$ of $\Delta$, the results provided in Section~\ref{sec:MID-intermediate} do not allow one to conclude on the dominance of $s_0$. This particular low-order example has been studied in further details in \cite[Theorem~4.2]{Boussaada2020Multiplicity} and that result shows, in particular, that, if $s_0$ is given by
\begin{equation}
\label{eq:example-EPP-3}
s_0 = \frac{- 2 + \sqrt{2 - \frac{g \tau^2}{L}}}{\tau},
\end{equation}
then $s_0$ is a dominant root of the quasipolynomial $\Delta$ from \eqref{DeltaPendulum}, while the other choice of $s_0$ yields a triple root at $s_0$ which is not dominant. In particular, choosing $s_0$ as in \eqref{eq:example-EPP-3} and $k_d$ and $k_p$ as in \eqref{eq:example-kd-kp-2}, the trivial solution of \eqref{InvPenLin}--\eqref{DPD} is exponentially stable.

\begin{figure}[ht]
\centering
\resizebox{0.5\textwidth}{!}{\import{Figures}{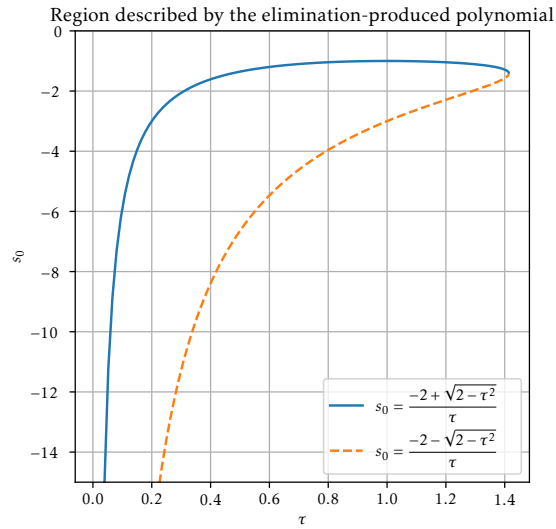}}
\caption[Region in the plane $(s_0, \tau)$ described by the elimination-produced polynomial \eqref{eq:example-EPP}]{Region in the plane $(s_0, \tau)$ described by the elimination-produced polynomial \eqref{eq:example-EPP} in the case $\frac{g}{L} = 1$. The two possible values for $s_0$ in \eqref{eq:example-EPP-2} are represented by a solid blue line and a dashed orange line}
\label{Pendulum-TRIPLE}
\end{figure}

Figure~\ref{Pendulum-TRIPLE} represents the region in the plane $(s_0, \tau)$ described by the polynomial equality \eqref{eq:example-EPP} in the case $\frac{g}{L} = 1$. Note that this region is restricted to the domain $\tau \in (0, \sqrt{2}]$, since, by \eqref{eq:example-EPP-2}, the quasipolynomial $\Delta$ from \eqref{DeltaPendulum} cannot admit a triple root if $\tau > \sqrt{2}$. The solid blue curve describes the choice of $s_0$ from \eqref{eq:example-EPP-3}, which is guaranteed to be a dominant root of \eqref{DeltaPendulum}, while the dashed orange curve describes the other choice of $s_0$ from \eqref{eq:example-EPP-2}, which is not dominant. Both curves meet at the point $(s_0, \tau) = (-\sqrt{2}, \sqrt{2})$, which corresponds to the configuration from Section~\ref{sec:pendulum-4} in which $s_0$ is a root of multiplicity $4$.

\subsection{Stabilization of a chain of \texorpdfstring{$n$}{n} integrators}
\label{sec:expl-integrator}

Our next example, extracted from \cite{Boussaada2025Decay}, deals with the case of the scalar control system given by a chain of $n$ integrators (or, equivalently, an $n$th order integrator), in which the control law is given by a delayed linear feedback involving derivatives of the state of order up to $n-1$. More precisely, we consider the control system
\begin{equation}
\label{eq:n-integrator}
y^{(n)}(t) = u(t)
\end{equation}
where $u(t) \in \mathbb R$ is the control input at time $t$, $y(t) \in \mathbb R$, and $n$ is a positive integer. We assume that the delayed feedback control law $u$ is chosen as
\begin{equation}
\label{eq:control-n-integrator}
u(t) = -\sum_{k=0}^{n-1} \alpha_k y^{(k)}(t - \tau)
\end{equation}
for some delay $\tau > 0$ and real coefficients $\alpha_0, \dotsc, \alpha_{n-1}$. System~\eqref{eq:n-integrator}--\eqref{eq:control-n-integrator} is a DDE of the form \eqref{eq:DDE} with $a_0 = \dotsb = a_{n-1} = 0$ and $m = n - 1$, and hence its characteristic function $\Delta$ is of the form \eqref{eq:Delta-DDE} and is given by
\begin{equation}
\label{eq:Delta-integrator}
\Delta(s) = s^n + e^{-s \tau} \sum_{k=0}^{n-1} \alpha_k s^k.
\end{equation}
Note that $\Delta$ is a quasipolynomial of degree $D = 2 n$.

The main aim of this example is to show that roots of $\Delta$ with the smallest overorder multiplicity $\Lambda = n + 1$ are necessarily dominant, exploiting for that the results from Section~\ref{sec:MID-intermediate}. We will consider here the configuration where the coefficients $\alpha_0, \dotsc, \alpha_{n-1}$ are all freely available for choice, which corresponds to selecting $I = \llbracket n + 1, 2n\rrbracket$ in the notation of Section~\ref{sec:iMID-EPF}. In this case, by Proposition~\ref{prop:EPP}, the elimination-produced function from Definition~\ref{def:EPF} is a polynomial, given by \eqref{eq:EPP}.

In order to analyze this example, we shall make use in the sequel of two classical results, which we now state. The first one is known as the Hermite--Poulain Theorem and its proof can be found, for instance, in \cite[Theorem~7.3.3]{Iliev1987Laguerre}.

\begin{lemma}[Hermite--Poulain Theorem]
\label{lem:HermitePoulain}
Let $n \in \mathbb N$, $c_0, \dotsc, c_n$ be real numbers, and $h$ be the polynomial given by $h(s) = c_0 + c_1 s + \dotsb + c_n s^n$ for $s \in \mathbb C$. Assume that all complex roots of $h$ are real numbers. Let $f$ be a polynomial with real coefficients and define the polynomial $F$ by
\[
F(s) = c_0 f(s) + c_1 f^\prime(s) + \dotsb + c_n f^{(n)}(s), \qquad \text{ for } s \in \mathbb C.
\]
Then $F$ has at least as many real roots as $f$ has (counted with multiplicity).
\end{lemma}

The second classical result we shall need concerns the distribution of roots of a function $F\colon \mathbb C \to \mathbb C$ of the form
\begin{equation}
\label{eq:F-Laplace-f}
F(s) = \int_0^1 f(t) e^{-s t} \diff t, \qquad \text{ for } s \in \mathbb C,
\end{equation}
where $f\colon [0, 1] \to \mathbb R$ is a given function. Note that integrals of this form appear in the factorizations \eqref{eq:Delta-max-mult} and \eqref{eq:Delta-EPP}. The function $F$ can be seen as the Laplace transform of a function $f$ whose support is included in $[0, 1]$, and the study of the distribution of its roots is related to a wide range of problems on Physics and Engineering and dates back to the pioneering works \cite{Hardy1905Zeroes, Polya1918Uber, Titchmarsh1926Zeros}. In particular, we have the following result, adapted from \cite{Polya1918Uber} (see also \cite[Part Five, Chapter~3, Problem~177]{Polya1998ProblemsII}).

\begin{lemma}
\label{lem:Polya}
Let $f\colon [0, 1] \to \mathbb R$ be integrable and assume that $f \in \mathcal C^1((0, 1))$, $f$ is positive on $(0, 1)$, and $f^\prime$ is negative on $(0, 1)$. Let $F$ be defined from $f$ by \eqref{eq:F-Laplace-f}. Then $\Real s < 0$ for every zero $s$ of $F$.
\end{lemma}

We are now in position to study the chain of $n$ integrators \eqref{eq:n-integrator} with the delayed feedback control \eqref{eq:control-n-integrator}. Our first result is that the corresponding elimination-produced polynomial $\mathbb P$ from \eqref{eq:EPP} only depends on $(s_0, \tau)$ through the product $s_0 \tau$.

\begin{lemma}
\label{lem:hat-P}
Let $n \in \mathbb N^\ast$, $\tau > 0$, $\alpha_0, \dotsc, \alpha_{n-1}$ be real numbers, and consider the quasipolynomial $\Delta$ from \eqref{eq:Delta-integrator}. Let $\mathbb P$ be the elimination-produced function from Definition~\ref{def:EPF} associated with the smallest overorder multiplicity $\Lambda = n + 1$ and in the case $I = \llbracket n + 1, 2n\rrbracket$. Then, for every $(s_0, \tau) \in \mathbb C \times \mathbb R_+$, we have
\begin{equation}
\label{eq:EPP-integrator}
\mathbb P(s_0, \tau) = \widehat{\mathbb P}(s_0 \tau),
\end{equation}
where $\widehat{\mathbb P}$ is the polynomial defined for $s \in \mathbb C$ by
\begin{equation}
\label{eq:EPP-integrator-hat}
\widehat{\mathbb P}(s) = \sum_{k = 0}^{n} k! \binom{n}{k}^2 s^{n-k}.
\end{equation}
In addition, all roots of $\widehat{\mathbb P}$ are negative real numbers.
\end{lemma}

\begin{proof}
Note that the polynomial $P_{F, \Lambda}$ from \eqref{eq:P_FM} is given by $P_{F, \Lambda}(s) = s^n$, and thus an immediate computation from \eqref{eq:EPP} yields \eqref{eq:EPP-integrator}--\eqref{eq:EPP-integrator-hat}. In addition, \eqref{eq:EPP-integrator-hat} can be rewritten in terms of $P_{F, \Lambda}$ as
\[
\widehat{\mathbb P}(s) = \sum_{k = 0}^{n} \binom{n}{k} P_{F, \Lambda}^{(k)}(s).
\]
Since the polynomial $s \mapsto (s + 1)^n = \sum_{k = 0}^n \binom{n}{k} s^k$ admits a unique root of multiplicity $n$ at the real number $-1$, Lemma~\ref{lem:HermitePoulain} shows that $\widehat{\mathbb P}$ has at least as many real roots as $P_{F, \Lambda}$ has (counted with multiplicity). Since $P_{F, \Lambda}$ has a single root of multiplicity $n$ at $0$, we deduce that $\widehat{\mathbb P}$ as at least $n$ real roots and, since the degree of $\widehat{\mathbb P}$ is $n$, it follows that all roots of $\widehat{\mathbb P}$ are real numbers. Since all coefficients of $\widehat{\mathbb P}$ are positive, we infer that all its real roots are necessarily negative, yielding the conclusion.
\end{proof}

We also have the following result on the polynomial $\widehat{\mathbb P}$ from \eqref{eq:EPP-integrator-hat}.

\begin{lemma}
\label{lem:widehatP-sign}
Let $n \in \mathbb N^\ast$, $\widehat{\mathbb P}$ be the polynomial defined in \eqref{eq:EPP-integrator-hat}, $c_1$ be the rightmost root of $\widehat{\mathbb P}$, and $q$ be defined by $q(t) = \widehat{\mathbb P}(c_1 t)$ for $t \in [0, 1]$. Then $q \in \mathcal C^1([0, 1])$ and, for every $t \in (0, 1)$, we have $q(t) > 0$ and $q^\prime(t) < 0$.
\end{lemma}

\begin{proof}
Notice first that, by Lemma~\ref{lem:hat-P}, all roots of $\widehat{\mathbb P}$ are negative real numbers, and thus we have $c_1 < 0$. Let us denote by $c_2, \dotsc, c_n$ the other roots of $\widehat{\mathbb P}$, repeated according to their multiplicities. Since $\widehat{\mathbb P}$ is a monic polynomial by \eqref{eq:EPP-integrator-hat}, it follows that $\widehat{\mathbb P}(s) = \prod_{k=1}^n (s - c_k)$, and an explicit computation from this formula shows that both $\widehat{\mathbb P}$ and $\widehat{\mathbb P}^\prime$ are positive in $(c_1, +\infty)$, and in particular in $(c_1, 0)$. As $c_1 < 0$, this entails that the polynomial $q$ is positive and its derivative is negative in $(0, 1)$, as required.
\end{proof}

Finally, combining Proposition~\ref{prop:EPP-factorization} with Lemmas~\ref{lem:Polya}, \ref{lem:hat-P}, and \ref{lem:widehatP-sign}, we deduce that, if the root $s_0 \in \mathbb R$ of multiplicity $n+1$ and the delay $\tau \in \mathbb R_+^\ast$ are chosen in such a way that $s_0 \tau = c_1$, where $c_1$ is as in the statement of Lemma~\ref{lem:widehatP-sign}, then necessarily $s_0$ is the rightmost root of the quasipolynomial $\Delta$ from \eqref{eq:Delta-integrator}, and we thus have the required dominance property. More precisely, we have the following result.

\begin{proposition}
\label{prop:n-integrator}
Let $n \in \mathbb N^\ast$, $\tau > 0$, $\alpha_0, \dotsc, \alpha_{n-1}$ be real numbers, and consider the quasipolynomial $\Delta$ from \eqref{eq:Delta-integrator}. Let $\widehat{\mathbb P}$ be the polynomial defined in \eqref{eq:EPP-integrator-hat} and take $s_0 \in \mathbb R$. Assume that $s_0$ is a root of multiplicity at least $n + 1$ of $\Delta$ and that $s_0 \tau$ is the rightmost root of $\widehat{\mathbb P}$. Then $\Real s < s_0$ for every $s \in \mathbb C \setminus \{s_0\}$ such that $\Delta(s) = 0$. As a consequence, the trivial solution of \eqref{eq:n-integrator}--\eqref{eq:control-n-integrator} is exponentially stable.
\end{proposition}

\begin{proof}
By Proposition~\ref{prop:EPP-factorization} and Lemma~\ref{lem:hat-P}, we have, for every $s \in \mathbb C$, that
\[
\Delta(s) = \frac{\tau (s - s_0)^{n+1}}{n!} \int_0^1 \widehat{\mathbb P}(s_0 \tau t) e^{-t \tau (s - s_0)} \diff t.
\]
Let $c_1$ be the rightmost root of $\widehat{\mathbb P}$ and recall that $c_1 < 0$. Let $q$ be the polynomial defined by $q(t) = \widehat{\mathbb P}(c_1 t)$. Combining Lemmas~\ref{lem:Polya} and \ref{lem:widehatP-sign} with the assumption that $s_0 \tau = c_1$, we obtain that all roots of $s \mapsto \int_0^1 q(t) e^{-t s} \diff t$ have negative real parts, and, consequently, all roots of $\Delta$ different from $s_0$ have real parts strictly less than $s_0$. In particular, the spectral abscissa of \eqref{eq:n-integrator}--\eqref{eq:control-n-integrator} is $s_0$ and, since $c_1 < 0$ and $\tau > 0$, the relation $s_0 \tau = c_1$ implies that $s_0 < 0$. Hence, the trivial solution of \eqref{eq:n-integrator}--\eqref{eq:control-n-integrator} is exponentially stable.
\end{proof}

To illustrate Proposition~\ref{prop:n-integrator}, let us consider the triple integrator, $n = 3$, in which case \eqref{eq:n-integrator}--\eqref{eq:control-n-integrator} is given by
\begin{equation}
\label{eq:triple-integrator}
y^{(3)}(t) + \alpha_2 y^{\prime\prime}(t - \tau) + \alpha_1 y^\prime(t - \tau) + \alpha_0 y(t - \tau) = 0
\end{equation}
and the corresponding characteristic function is
\begin{equation}
\label{eq:triple-integrator-Delta}
\Delta(s) = s^n + e^{-s \tau} (\alpha_2 s^2 + \alpha_1 s + \alpha_0).
\end{equation}
The function $\widehat{\mathbb P}$ from \eqref{eq:EPP-integrator-hat} is given by
\[
\widehat{\mathbb P}(s) = s^3 + 9 s^2 + 18 s + 6,
\]
and one can verify that its three roots are real and negative, and that its rightmost root is
\[c_1 = -3 - 2 \sqrt{3} \cos\left(\frac{\arctan\left(\sqrt{2}\right) + 2 \pi}{3}\right) \approx -0.4157746.\]

Assume that $\tau = \frac{1}{2}$ and let $s_0 \in \mathbb R$. We wish to select the parameters $\alpha_0, \alpha_1, \alpha_2$ in \eqref{eq:triple-integrator} so that $s_0$ is a root of multiplicity at least $4$ of the function $\Delta$ from \eqref{eq:triple-integrator-Delta}. Imposing that $\Delta(s_0) = \Delta^\prime(s_0) = \Delta^{\prime\prime}(s_0) = \Delta^{(3)}(s_0) = 0$, we obtain that
\[
\left\{
\begin{aligned}
\alpha_0 & = -\frac{1}{2} e^{s_{0} \tau} s_{0}^{3} \left(s_{0}^{2} \tau^{2} + 4 s_{0} \tau + 2\right), \\
\alpha_1 & = e^{s_{0} \tau} s_{0}^{2} \left(s_{0}^{2} \tau^{2} + 5 s_{0} \tau + 3\right), \\
\alpha_2 & = -\frac{1}{2} e^{s_{0} \tau} s_{0} \left(s_{0}^{2} \tau^{2} + 6 s_{0} \tau + 6\right),
\end{aligned}
\right.
\]
and
\[
\widehat{\mathbb P}(s_0 \tau) = 0.
\]

In accordance with Proposition~\ref{prop:n-integrator}, we select $s_0$ so that $s_0 \tau = c_1$, which yields
\begin{equation}
\label{eq:triple-integrator-parameters}
\left\{
\begin{aligned}
s_0 & \approx -0.8315491, \\
\alpha_0 & \approx 0.09670293,\\
\alpha_1 & \approx 0.4991405,\\
\alpha_2 & \approx 1.009084.
\end{aligned}
\right.
\end{equation}
Proposition~\ref{prop:n-integrator} then ensures that the spectral abscissa of \eqref{eq:triple-integrator} is $s_0 \approx -0.8315491$, implying that \eqref{eq:triple-integrator} is exponentially stable.

\begin{figure}[ht]
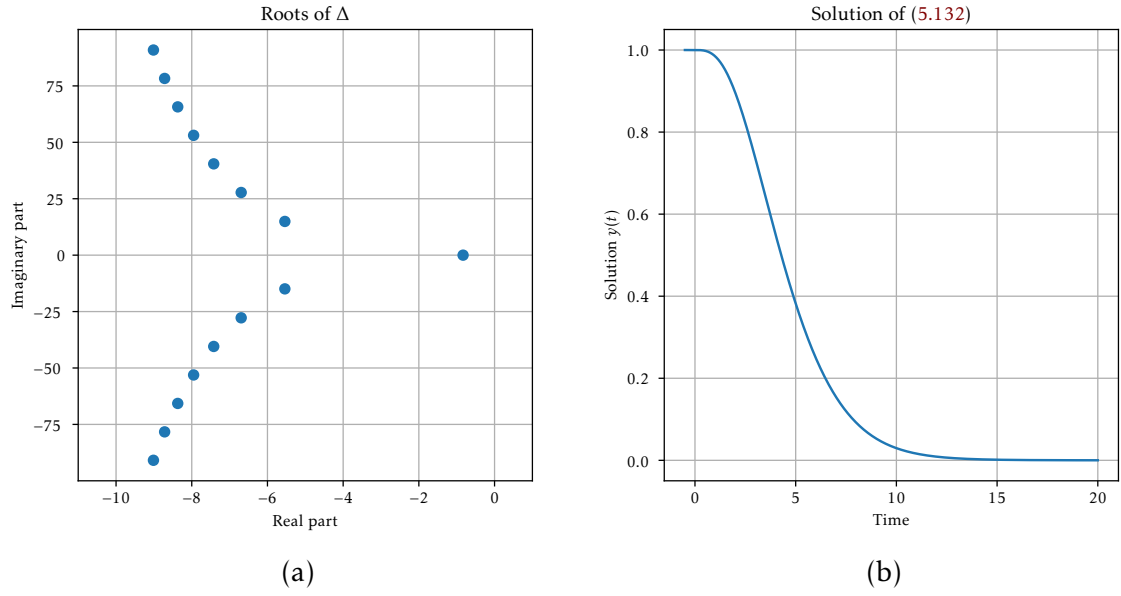

\centering
\begin{tabular}{@{} >{\centering} m{0.5\textwidth} @{} >{\centering} m{0.5\textwidth} @{}}
\resizebox{0.5\textwidth}{!}{\input{Figures/Expl7-TripleIntegrator-Spectrum.pgf}} & \resizebox{0.5\textwidth}{!}{\input{Figures/Expl7-TripleIntegrator-Solution.pgf}} \tabularnewline
(a) & (b) \tabularnewline
\end{tabular}
\caption[Spectrum and solution of \eqref{eq:triple-integrator} with parameters given by \eqref{eq:triple-integrator-parameters}]{(a) Spectrum and (b) solution of \eqref{eq:triple-integrator} for $\tau = \frac{1}{2}$ with parameters given by \eqref{eq:triple-integrator-parameters}}
\label{fig:triple-integrator}
\end{figure}

Figure~\ref{fig:triple-integrator} represents the spectrum and a solution of \eqref{eq:triple-integrator} with $\tau = \frac{1}{2}$ and parameters given by \eqref{eq:triple-integrator-parameters}. The roots of $\Delta$ from \eqref{eq:triple-integrator-Delta} were computed numerically using Python's \texttt{cxroots} module \cite{Parini2018cxroots}, and the numerical search for roots was restricted to the rectangle $\{s \in \mathbb C \suchthat -15 \leq \Real s \leq 10,\, \abs{\Imag s} \leq 100\}$. The solution represented in Figure~\ref{fig:triple-integrator}(b) corresponds to a constant initial condition equal to $1$ in the interval $[-\tau, 0]$, and it was computed by an explicit Euler scheme applied to the first-order system in $\mathbb R^3$ obtained by a standard state augmentation from \eqref{eq:triple-integrator}. We observe, in Figure~\ref{fig:triple-integrator}(a), the root of multiplicity $4$ at $s_0 \approx -0.8315491$, which corresponds to the spectral abscissa of the system. Figure~\ref{fig:triple-integrator}(b) shows that, as expected, the represented solution converges exponentially to $0$.

\subsection{Active vibration suppression}
\label{sec:avs}

\begin{figure}[ht]
\centering

%
%
\begin{tikzpicture}
\foreach \x in {-1.9, -1.8, ..., 2} {
  \draw[thin, gray] ({\x-0.1}, -0.1) -- (\x, 0);
}
\draw[thick] (-2, 0) -- (2, 0);

\draw (-1, 0) -- (-1, 0.45) (-1.2, 0.55) -- (-1.2, 0.45) -- (-0.8, 0.45) -- (-0.8, 0.55) (-1.15, 0.52) -- (-0.85, 0.52) (-1, 0.52) -- (-1, 1);
\draw (1, 0) -- (1, 0.2) -- (1.2, 0.25) -- (0.8, 0.35) -- (1.2, 0.45) -- (0.8, 0.55) -- (1.2, 0.65) -- (0.8, 0.75) -- (1, 0.8) -- (1, 1);

\node[draw, thick, rectangle, anchor=south, minimum width=4cm, minimum height=1.5cm] at (0, 1) {(P)};

\draw[red, ->] (1.9, 0.2) -- node[midway, right] {$f(t)$} (1.9, 1);

\draw (-0.5, 2.5) -- (-0.5, 2.95) (-0.7, 3.05) -- (-0.7, 2.95) -- (-0.3, 2.95) -- (-0.3, 3.05) (-0.65, 3.02) -- (-0.35, 3.02) (-0.5, 3.02) -- (-0.5, 3.5);
\draw (0.5, 2.5) -- (0.5, 2.7) -- (0.7, 2.75) -- (0.3, 2.85) -- (0.7, 2.95) -- (0.3, 3.05) -- (0.7, 3.15) -- (0.3, 3.25) -- (0.5, 3.3) -- (0.5, 3.5);

\node[draw, thick, rectangle, anchor=south, minimum width=3cm, minimum height=1.25cm] at (0, 3.5) {(A)};

\draw[blue, ->] (1.4, 2.7) -- node[midway, right] {$u(t)$} (1.4, 3.5);
\end{tikzpicture}

\caption[Active vibration suppression through an active vibration absorber above a primary structure]{Active vibration suppression through an active vibration absorber (A) above a primary structure (P) in order to suppress displacement induced by a harmonic disturbance force $f(t)$}
\label{fig:avs}
\end{figure}
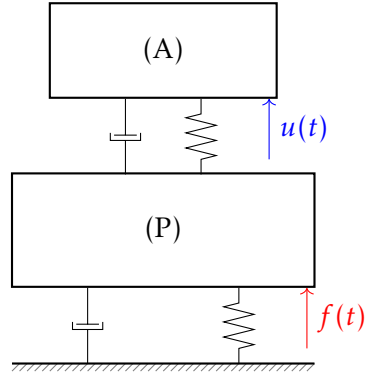

In this application, extracted from \cite{Mazanti2020Spectral}, we consider the situation from Figure~\ref{fig:avs}, in which a primary structure (P), modeled through a mass-spring-damper system, is subject to an external disturbance $f(t)$, which is assumed to be a trigonometric function with known angular frequency $\omega$. In order to reduce the effects of the disturbance, an active vibration absorber (A) is placed above the primary structure (P). The active vibration absorber is also modeled as a mass-spring-damper system, and it is actuated through a control input $u$. The control goal is then to choose the control $u$, under a feedback form depending only on the position $x_a$ of the vibration absorber (A), in order to suppress vibrations due to $f$ on the primary structure (P).

The dynamics of the vibration absorber can be described by
\begin{equation}
\label{eq:avs}
x_a^{\prime\prime}(t) + 2 \zeta \Omega x_a^\prime(t) + \Omega^2 x_a(t) = \frac{1}{m_a} u(t),
\end{equation}
where $m_a > 0$ is the mass of (A), $\zeta > 0$ and $\Omega > 0$ are the constants of the mass-spring-damper system corresponding to the absorber, $x_a$ is the position of the absorber, and $u$ is the control input.

The article \cite{Olgac1994Novel} introduced a new idea in vibration suppression by proposing to select the control $u(t)$ in \eqref{eq:avs} under a delayed feedback form, treating the delay as a design parameter. In this design, known as \emph{delayed resonator}, the parameters of the feedback law are selected in such a way as to ensure that the characteristic function of the closed-loop system admits $\pm i \omega$ as roots and that it has no roots with positive real part. We refer to \cite{Olgac1994Novel} for a more detailed description of the delayed resonator and its comparison to other vibration suppression techniques.

A further modification of the delayed resonator was proposed in \cite{Kure2018Spectral}, in order to impose that the roots $\pm i \omega$ of the closed-loop system are of multiplicity $2$ each, using for that purpose two different delays in the feedback law. In this section, similarly to \cite{Kure2018Spectral}, we aim at ensuring that $\pm i \omega$ are roots of multiplicity $2$ each of the characteristic function of the closed-loop system, but we use for that a feedback law with a single delay, namely
\begin{equation}
\label{eq:avs-feedback}
u(t) = -a_0 x_a(t) - a_1 x_a^\prime(t) - \alpha_0 x_a(t - \tau) - \alpha_1 x_a^\prime(t - \tau),
\end{equation}
where $a_0, a_1, \alpha_0, \alpha_1$ are real parameters to be selected and the delay $\tau > 0$ is also seen as a parameter to be selected.

The characteristic function of the closed-loop system \eqref{eq:avs}--\eqref{eq:avs-feedback} is
\begin{equation}
\label{eq:avs-Delta}
\Delta(s) = s^2 + \left(2 \zeta \Omega + \frac{a_1}{m_a}\right) s + \Omega^2 + \frac{a_0}{m_a} + e^{-s \tau} \left(\frac{\alpha_1}{m_a} s + \frac{\alpha_0}{m_a}\right).
\end{equation}
By Theorem~\ref{thm:complex-dominance}, the function $\Delta$ from \eqref{eq:avs-Delta} admits roots of multiplicity $2$ each at $\pm i \omega$ if and only if the coefficients $a_0, a_1, \alpha_0, \alpha_1$ are selected as
\begin{equation}
\label{eq:avs-coeffs}
\left\{
\begin{aligned}
a_1 & = -2 m_a \zeta \Omega - 2 m_a \omega \frac{\tau \omega - \sin(\tau \omega)\cos(\tau \omega)}{\tau^2 \omega^2 - \sin^2(\tau \omega)}, \\
a_0 & = -m_a \Omega^2 + m_a \omega^2 \frac{\tau^2 \omega^2 + \sin^2(\tau \omega)}{\tau^2 \omega^2 - \sin^2(\tau \omega)}, \\
\alpha_1 & = 2 m_a \omega \frac{\tau \omega \cos(\tau \omega) - \sin(\tau \omega)}{\tau^2 \omega^2 - \sin^2(\tau \omega)}, \\
\alpha_0 & = -2 m_a \omega \frac{\tau \omega^2 \sin(\tau \omega)}{\tau^2 \omega^2 - \sin^2(\tau \omega)},
\end{aligned}
\right.
\end{equation}
and in this case all other roots of $\Delta$ have negative real part. Note that, if the delay $\tau$ is selected as
\begin{equation}
\label{eq:avs-delay}
\tau = \frac{k \pi}{\omega} \qquad \text{ for some } k \in \mathbb N^\ast,
\end{equation}
then \eqref{eq:avs-coeffs} can be further simplified to
\begin{equation}
\label{eq:avs-coeffs-simple}
\left\{
\begin{aligned}
a_1 & = -2 m_a \left(\zeta \Omega + \frac{1}{\tau}\right), \\
a_0 & = -m_a (\Omega^2 - \omega^2), \\
\alpha_1 & = 2 m_a \frac{(-1)^k}{\tau}, \\
\alpha_0 & = 0.
\end{aligned}
\right.
\end{equation}
Hence, with the choices of parameters given by \eqref{eq:avs-delay}--\eqref{eq:avs-coeffs-simple}, the delayed feedback \eqref{eq:avs-feedback} ensures that the characteristic function of \eqref{eq:avs} has roots at $\pm i \omega$ of multiplicity $2$ each and no other roots with positive real part, ensuring thus the control goal of vibration suppression at angular frequency $\omega$.

\subsection{Flexible mode compensation}

We now consider another application extracted from \cite{Mazanti2020Spectral} in which, similarly to the active vibration suppression from Section~\ref{sec:avs}, the control goal can be achieved through the design of a characteristic function with imaginary roots $\pm i \omega$.

\begin{figure}[ht]
\centering
%
%
\begin{tikzpicture}

\node[draw, rectangle, minimum width = 1.5cm, minimum height = 1cm] (F) at (0, 0) {$F$};
\node[anchor=north, text width=1.8cm, align=center] at (F.south) {Flexible subsystem};

\node[draw, rectangle, minimum width = 1.5cm, minimum height = 1cm] (G) at (-4, 0) {$G$};
\node[anchor=north, text width=1.8cm, align=center] at (G.south) {Main body};

\node[draw, rectangle, minimum width = 1.5cm, minimum height = 1cm] (C) at (-7, 0) {$C$};
\node[anchor=north, text width=1.8cm, align=center] at (C.south) {Controller};

\node[draw, rectangle, minimum width = 1.5cm, minimum height = 1cm] (D) at (-5.5, -2) {$K$};
\node[anchor=north, text width=1.8cm, align=center] at (D.south) {Feedback};

\node[draw, circle, minimum size=0.4cm] (Sum1) at (-9, 0) {};
\node[draw, circle, minimum size=0.4cm] (Sum2) at (-2, 0) {};

\fill (-1.5, 0) circle[radius=0.05];

\draw[-Stealth] (F.east) -- ++(1, 0) node[above] {$y_s$};
\draw[-Stealth] (Sum2.east) --node[midway, above] {$y$} (F.west);
\draw[-Stealth] (Sum2.north) ++(0, 0.75) node[left] {$d$} -- node[at end, right] {$+$} (Sum2.north);
\draw[-Stealth] (G.east) -- node[at end, below] {$+$} (Sum2.west);
\draw[-Stealth] (C.east) -- (G.west);
\draw[-Stealth] (Sum1.east) -- (C.west);
\draw[-Stealth] (Sum1.west) ++(-0.75, 0) node[above] {$w$} -- node[at end, above] {$+$} (Sum1.west);
\draw[-Stealth] (-1.5, 0) -- (-1.5, -2) -- (D.east);
\draw[-Stealth] (D.west) -- node[very near start, above] {$v$} (Sum1 |- D) -- (Sum1.south) node[right] {$-$};
\end{tikzpicture}
\caption[Diagram for the control of a body attached to a flexible subsystem]{Block diagram for the control of a main body attached to a flexible subsystem}
\label{fig:fmc}
\end{figure}
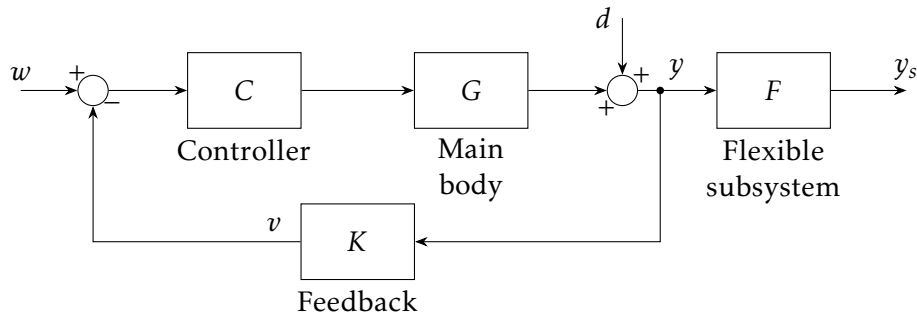

Consider a system made of a main body connected to a flexible subsystem, whose oscillations we want to compensate by acting on and observing the main body only. A typical application is the position control of a crane trolley (the main body) with a suspended payload (the flexible subsystem), the control goal being to compensate the oscillations of the suspended payload, i.e., the payload should not sway once the main body position $y$ reaches the set-point value $w$. For that purpose, we will consider the classical control scheme of Figure~\ref{fig:fmc}, with a controller $C$ and a feedback $K$ acting on the main body $G$, the output of the main body being the input of the flexible subsystem $F$. We also consider a possible additive disturbance $d$ between the main body and the flexible subsystem.

When the feedback $K$ is chosen under the form $K(s) = \frac{1}{\Delta(s)}$ for some function $\Delta$, the transfer function $T_{y_s, w}$ from the input $w$ to the output $y_s$ is
\[
T_{y_s, w}(s) = \frac{C(s) G(s) \Delta(s)}{\Delta(s) + C(s) G(s)} F(s),
\]
while the transfer function $T_{y_s, d}$ from the disturbance $d$ to the output $y_s$ is
\[
T_{y_s, d}(s) = \frac{\Delta(s)}{\Delta(s) + C(s) G(s)} F(s).
\]
In particular, if the oscillation of the flexible subsystem has a known angular frequency $\omega$, one can compensate for this oscillatory behavior by selecting $\pm i \omega$ to be zeros of the function $\Delta$. For instance, for the crane trolley application, the value of $\omega$ is known and given by $\omega = \sqrt{\frac{g}{L}}$, where $L$ is the length of the payload and $g$ is gravitational acceleration. According to Theorem~\ref{thm:complex-dominance}, a possible choice of $\Delta$ is thus
\begin{equation}
\label{eq:fmc-delta}
\Delta(s) = s^2 + a_1 s + a_0 + e^{-s \tau} (\alpha_1 s + \alpha_0)
\end{equation}
for some delay value $\tau > 0$ to be chosen and with $a_0, a_1, \alpha_0, \alpha_1$ being given by \eqref{eq:coeffs-complex-dominance} for $\sigma_0 = 0$ and $\theta_0 = \omega$. Let us note that if the mode to be compensated is damped, i.e., it is given by $-\zeta\omega \pm i \omega\sqrt{1-\zeta^2}$ for some $\zeta \in (0, 1)$, then the parameters of $\Delta(s)$ can be adapted according to \eqref{eq:coeffs-complex-dominance} by taking $\sigma_0 = -\zeta\omega$ and $\theta_0 = \omega\sqrt{1-\zeta^2}$.

The controller design from Figure~\ref{fig:fmc} in the case $K(s) = \frac{1}{\Delta(s)}$ corresponds to the \emph{inverse shaper} scheme proposed in \cite{Vyhlidal2016Feedback}, and the design of $\Delta$ through \eqref{eq:fmc-delta} is a generalization of the designs of $\Delta$ from \cite{Vyhlidal2016Feedback} with discrete or distributed delays.

\subsection{Constrained stabilization of a simple integrator}

We now consider the stabilization of a simple integrator with constraints on the coefficients of the feedback law, an example extracted from \cite{Fueyo2023Pole} in order to illustrate the interest of Theorem~\ref{thm:main-two-delays}. Consider the control system defined by the simple integrator, i.e., the control system
\begin{equation}
\label{eq:syst-integrator}
y^\prime(t) = u(t),
\end{equation}
where $y(t) \in \mathbb R$ is the state and $u(t) \in \mathbb R$ is the control. We consider the stabilization of \eqref{eq:syst-integrator} by linear feedback laws, with or without delays, with constraints in the corresponding controller gains.

In case of no delays, assume that 
\begin{equation}
\label{eq:feedback-0}
u(t) = -a y(t)
\end{equation}
with the constraint that $\abs{a} \leq 1$. In this case, the characteristic function $\Delta$ of closed-loop system $y^\prime(t) + a y(t) = 0$ is given by  $\Delta(s) = s + a$, with a single root at $s = -a$, also equal to the spectral abscissa of that closed-loop system. In particular, the best achievable exponential decay rate for \eqref{eq:syst-integrator}---i.e., the minimal spectral abscissa for \eqref{eq:syst-integrator}---with the feedback law \eqref{eq:feedback-0} and with the constraint $\abs{a} \leq 1$ is equal to $-1$, attained for $a = 1$.

Let us now consider, instead of \eqref{eq:feedback-0}, the linear feedback law with a single delay given by
\begin{equation}
\label{eq:feedback-1}
u(t) = -a y(t - \tau),
\end{equation}
where $\tau > 0$ is the delay, and still with the constraint $\abs{a} \leq 1$. In this case, the characteristic function $\Delta$ of the closed-loop system $y^\prime(t) + a y(t - \tau) = 0$ is the function $\Delta$ defined by $\Delta(s) = s + a e^{- s \tau}$. Using the results from \cite{Hayes1950Roots}, one can show that the minimal spectral abscissa for \eqref{eq:syst-integrator} with the feedback law \eqref{eq:feedback-1} and with the constraints $\abs{a} \leq 1$ and $\tau > 0$ is equal to $-e \approx -2.718$, attained for $a = 1$ and $\tau = e^{-1} \approx 0.3679$. Note that this configuration corresponds to the case where $\Delta$ admits a root of maximal multiplicity $2$, since one can easily verify that the relations \eqref{Coeffs} from Theorem~\ref{thm:main-MID} are satisfied in this case.

Using Theorem~\ref{thm:main-two-delays}, we can design a feedback control law with two delays and constraints on the gains yielding an even smaller spectral abscissa for the corresponding closed-loop system. Indeed, consider the linear feedback law with two delays given by
\begin{equation}
\label{eq:feedback-2}
u(t) = -a_1 y(t - \tau_1) - a_2 y(t - \tau_2),
\end{equation}
where $\tau_1$ and $\tau_2$ are positive delays with $\tau_2 > \tau_1$, and with the constraint $\abs{a_1} + \abs{a_2} \leq 1$. The closed-loop system is $y^\prime(t) + a_1 y(t - \tau_1) + a_2 y(t - \tau_2) = 0$ and the corresponding characteristic function is $\Delta(s) = s + a_1 e^{- s \tau_1} + a_2 e^{- s \tau_2}$. Let us select $a_1$ and $a_2$ as in Theorem~\ref{thm:main-two-delays}, i.e., given by \eqref{eq:two-delays-coeffs}. The constraint that $\abs{a_1} + \abs{a_2} \leq 1$ then reads
\[
\frac{\tau_2}{\tau_1 (\tau_2-\tau_1)}e^{s_0 \tau_1} + \frac{\tau_1}{\tau_2(\tau_2-\tau_1)}e^{s_0 \tau_2} \leq 1,
\]
and, by Theorem~\ref{thm:main-two-delays}, the spectral abscissa of the closed-loop system is $\gamma = -\frac{1}{\tau_1} - \frac{1}{\tau_2}$. A numerical constrained minimization algorithm\footnote{Computations were carried out using the function \texttt{minimize} from Python's \texttt{scipy.optimize} module \cite{Virtanen2020SciPy}, using the method of trust region with constraints \texttt{trust-const} from \cite{Conn2000Trust}.} shows that the smallest spectral abscissa that can be achieved with the feedback law \eqref{eq:feedback-2}, the constraints $\abs{a_1} + \abs{a_2} \leq 1$ and $\tau_2 > \tau_1 > 0$, and choosing $a_1$ and $a_2$ as in Theorem~\ref{thm:main-two-delays} is approximately $-3.353$, attained for $a_1 \approx 0.9882$, $a_2 \approx -0.01176$, $\tau_1 \approx 0.4063$, and $\tau_2 \approx 1.122$. In particular, this spectral abscissa is smaller than the smallest spectral abscissa $-e$ that can be achieved with the single-delay feedback law \eqref{eq:feedback-1} with the constraint $\abs{a} \leq 1$.

\begin{figure}[ht]
\centering
\resizebox{0.5\textwidth}{!}{\import{Figures}{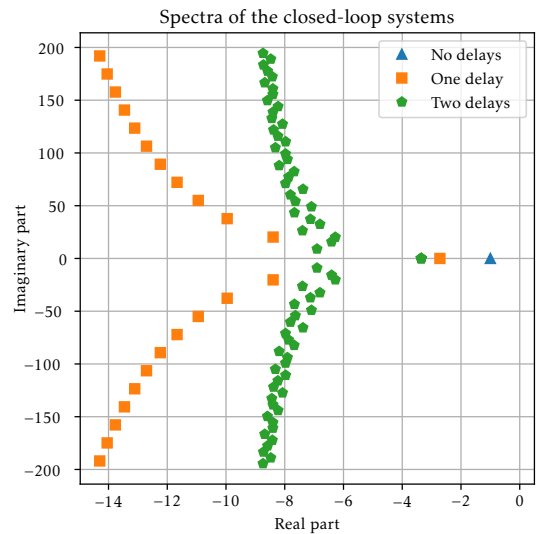}}
\caption[Spectra of \eqref{eq:syst-integrator} with feedback laws with zero, one, or two delays]{Spectra of \eqref{eq:syst-integrator} with the feedback law \eqref{eq:feedback-0} with no delays (blue triangles), the feedback law \eqref{eq:feedback-1} with one delay (orange squares), and the feedback law \eqref{eq:feedback-2} with two delays (green pentagons), when the parameters of these feedback laws are chosen as described in the text}
\label{fig:spectra}
\end{figure}

Figure~\ref{fig:spectra} shows the spectra\footnote{Spectrum computations were performed using Python's \texttt{cxroots} module \cite{Parini2018cxroots}. The search for roots was restricted to the rectangle $\{z \in \mathbb C \suchthat -30 \leq \Real(z) \leq 2,\, \abs{\Imag(z)} \leq 200\}$.} of \eqref{eq:syst-integrator} with the feedback laws \eqref{eq:feedback-0}, \eqref{eq:feedback-1}, and \eqref{eq:feedback-2} with the above choices of parameters. We remark that the delayed feedback law \eqref{eq:feedback-1} with a single delay allows one to improve the spectral abscissa with respect to the delay-free proportional feedback law \eqref{eq:feedback-0}, and that the spectral abscissa can be further improved with the feedback law with two delays \eqref{eq:feedback-2} by exploiting the result of Theorem~\ref{thm:main-two-delays}.

\begin{remark}
For the feedback laws \eqref{eq:feedback-0} and \eqref{eq:feedback-1}, the above choices of parameters can be shown to be those that minimize the spectral abscissa under the provided constraints on the coefficients of the feedback laws, either by simple arguments in the case of \eqref{eq:feedback-0} or by exploiting the results of \cite{Hayes1950Roots} in the case of \eqref{eq:feedback-1}. However, the above choice of parameters for the feedback law \eqref{eq:feedback-2} was selected under the additional constraint that these parameters satisfy \eqref{eq:two-delays-coeffs}, by performing a numerical minimization of the spectral abscissa under this constraint. Whether this coincides with the minimization of the spectral abscissa of \eqref{eq:syst-integrator} over all feedback laws \eqref{eq:feedback-2} satisfying $\abs{a_1} + \abs{a_2} \leq 1$ is an open question. Note, however, that the particular choice of parameters for \eqref{eq:feedback-2} above already performs better than any possible choice of parameters for \eqref{eq:feedback-0} and \eqref{eq:feedback-1} satisfying the corresponding constraints.
\end{remark}

{
\renewcommand{\thesection}{\thechapter.A}
\section{Appendix}
\label{sec:MID-app}

\subsection{An identity on the integral of a polynomial and an exponential}
\label{app:integralPExp}

This chapter makes use of the following result on the integral of the product of a polynomial and an exponential, which can be seen as the computation of the Laplace transform of a polynomial restricted to the interval $[0, 1]$.

\begin{proposition}
\label{PropIntegralPExp}
Let $d \in \mathbb N$ and $p$ be a polynomial of degree at most $d$. Then, for every $s \in \mathbb C \setminus \{0\}$,
\begin{equation}
\label{IntegralPExp}
\int_0^1 p(t) e^{-s t} \diff t = \sum_{k=0}^d \frac{p^{(k)}(0) - p^{(k)}(1) e^{-s}}{s^{k+1}}.
\end{equation}
\end{proposition}

\begin{proof}
The proof is done by induction on $d$. If $d = 0$, then $p$ is constant and one immediately verifies that \eqref{IntegralPExp} holds. Let now $d \in \mathbb N$ be such that \eqref{IntegralPExp} holds for every polynomial of degree at most $d$ and let $p$ be a polynomial of degree $d+1$. Integrating by parts and using the fact that $p^\prime$ is a polynomial of degree at most $d$, one gets
\begin{align*}
\int_0^1 p(t) e^{-st} \diff t & = \frac{p(0) - p(1) e^{-s}}{s} + \frac{1}{s} \int_0^1 p^\prime(s) e^{-st} \diff t \displaybreak[0] \\
& = \frac{p(0) - p(1) e^{-s}}{s} + \frac{1}{s} \sum_{k=0}^d \frac{p^{(k+1)}(0) - p^{(k+1)}(1) e^{-s}}{s^{k+1}} = \sum_{k=0}^{d+1} \frac{p^{(k)}(0) - p^{(k)}(1) e^{-s}}{s^{k+1}},
\end{align*}
establishing thus \eqref{IntegralPExp} by induction.
\end{proof}

\subsection{Some identities of binomial coefficients}
\label{app:binomial}

Let us present some properties of binomial coefficients used in this chapter. Even though the proofs of most of such properties are straightforward, they are provided here for sake of completeness. We recall that, with the convention that $\binom{k}{j} = 0$ whenever $k < 0$, $j < 0$, or $j > k$, binomial coefficients satisfy the relation $\binom{k}{j} = \binom{k-1}{j-1} + \binom{k-1}{j}$ for every $(k, j) \in \mathbb Z^2 \setminus \{(0, 0)\}$.

The first property on binomial coefficients needed in this chapter is the following shifting identity.

\begin{lemma}
\label{lemm:ProdBinom}
Let $j, k, \ell \in \mathbb N$ be such that $j \leq k \leq \ell$. Then
\[
\binom{k}{j} \binom{\ell}{k} = \binom{\ell}{j} \binom{\ell - j}{k - j} = \binom{\ell - k + j}{j} \binom{\ell}{\ell - k + j}.
\]
\end{lemma}

\begin{proof}
One immediately computes
\begin{align*}
\binom{k}{j} \binom{\ell}{k} & = \frac{k! \ell!}{(k-j)! (\ell-k)! j! k!} = \frac{\ell! (\ell - j)!}{(k-j)! (\ell-k)! j! (\ell - j)!} = \binom{\ell}{j} \binom{\ell - j}{k - j} \\
\intertext{and}
\binom{k}{j} \binom{\ell}{k} & = \frac{k! \ell!}{(k-j)! (\ell-k)! j! k!} = \frac{\ell! (\ell - k + j)!}{(k-j)! (\ell-k)! j! (\ell - k + j)!} = \binom{\ell - k + j}{j} \binom{\ell}{\ell - k + j}. \qedhere
\end{align*}
\end{proof}

Another auxiliary result we need is the following, concerning the sum of part of a row of binomial coefficients with alternating signs.

\begin{lemma}
\label{lemm:SumLine}
Let $k \in \mathbb N^\ast$ and $\ell \in \llbracket 0, k\rrbracket$. Then
\begin{equation}
\label{SumOfALineInThePascalTriangleWithAlternatingSigns}
\sum_{j=0}^{\ell} (-1)^j \binom{k}{j} = (-1)^{\ell} \binom{k-1}{\ell}.
\end{equation}
\end{lemma}

\begin{proof}
Let $k \in \mathbb N^\ast$. We prove \eqref{SumOfALineInThePascalTriangleWithAlternatingSigns} for $\ell \in \llbracket 0, k\rrbracket$ by induction on $\ell$. Clearly, \eqref{SumOfALineInThePascalTriangleWithAlternatingSigns} is satisfied for $\ell = 0$, since both left- and right-hand sides of \eqref{SumOfALineInThePascalTriangleWithAlternatingSigns} are equal to $1$ in this case. Assume now that $\ell \in \llbracket 0, k-1\rrbracket$ is such that \eqref{SumOfALineInThePascalTriangleWithAlternatingSigns} is satisfied. Then
\begin{equation*}
\sum_{j=0}^{\ell+1} (-1)^j \binom{k}{j} = (-1)^\ell \binom{k-1}{\ell} + (-1)^{\ell + 1} \binom{k}{\ell + 1} = (-1)^{\ell + 1} \left[\binom{k}{\ell + 1} - \binom{k-1}{\ell}\right] = (-1)^{\ell + 1} \binom{k-1}{\ell+1},
\end{equation*}
showing that \eqref{SumOfALineInThePascalTriangleWithAlternatingSigns} is satisfied with $\ell$ replaced by $\ell + 1$. Hence, by induction, \eqref{SumOfALineInThePascalTriangleWithAlternatingSigns} holds for every $\ell \in \llbracket 0, k\rrbracket$.
\end{proof}

The main identity used in this chapter is the following one.

\begin{proposition}
\label{PropHorribleSum}
For $j, k, n \in \mathbb N$, let
\begin{equation*}
S^n_{j, k} = \sum_{\ell = 0}^k (-1)^\ell \binom{n+k-j}{\ell} \binom{n+k-\ell}{n}.
\end{equation*}
Then, for every $k, n \in \mathbb N$ and $j \in \llbracket 0, n\rrbracket$, one has
\begin{equation*}
S^n_{j, k} = \binom{j}{k}.
\end{equation*}
\end{proposition}

\begin{proof}
If $j \neq n + k + 1$, then
\begin{align*}
S^n_{j, k} - S^n_{j-1, k} & = \sum_{\ell = 0}^k (-1)^\ell \left[\binom{n+k-j}{\ell} - \binom{n+k-j+1}{\ell}\right] \binom{n+k-\ell}{n} \displaybreak[0] \\
& = \sum_{\ell = 1}^k (-1)^{\ell-1} \binom{n+k-j}{\ell-1} \binom{n+k-\ell}{n} \displaybreak[0] \\
& = \sum_{\ell = 0}^{k-1} (-1)^{\ell} \binom{n+k-j}{\ell} \binom{n+k-\ell-1}{n} = S^n_{j-1, k-1}.
\end{align*}

For $j \geq 1$, using Lemma~\ref{lemm:ProdBinom}, we have
\[
S^n_{j, j} = \sum_{\ell = 0}^j (-1)^\ell \binom{n}{\ell} \binom{n+j-\ell}{n} = \sum_{\ell = 0}^j (-1)^\ell \binom{j}{\ell} \binom{n + j - \ell}{j}.
\]
For $\ell \in \llbracket 0, j\rrbracket$, let $a_\ell = (-1)^\ell \binom{j}{\ell}$, $b_\ell = \binom{n + j - \ell}{j}$, and $A_\ell = \sum_{m=0}^\ell a_m$. By Lemma~\ref{lemm:SumLine}, one has $A_\ell = (-1)^\ell \binom{j-1}{\ell}$, and one immediately computes $b_\ell - b_{\ell+1} = \binom{n+j-\ell-1}{j-1}$ for $\ell \in \llbracket 0, j-1\rrbracket$. Noticing that $S^n_{j, j} = \sum_{\ell = 0}^j a_\ell b_\ell$, one computes, using summation by parts (see, e.g., \cite[Theorem~3.41]{Rudin1976Principles}), that
\[
S^n_{j, j} = \sum_{\ell = 0}^{j-1} A_\ell (b_\ell - b_{\ell + 1}) = \sum_{\ell = 0}^{j-1} (-1)^\ell \binom{j-1}{\ell} \binom{n+j-\ell-1}{j-1} = S^n_{j-1, j-1}.
\]
We also compute, for $j \in \llbracket 0, n\rrbracket$, that
\[S^n_{j, 0} = \binom{n-j}{0} \binom{n}{n} = 1.\]

We have thus shown that
\begin{subequations}
\begin{align}
S^n_{j, k} & = S^n_{j-1, k-1} + S^n_{j-1, k}, & \qquad & \text{ for } j \neq n + k + 1, \label{SInduction} \\
S^n_{j, j} & = S^n_{j-1, j-1}, & & \text{ for } j \geq 1, \label{SDiagonal} \\
S^n_{j, 0} & = 1, & & \text{ for } j \in \llbracket 0, n\rrbracket. \label{SFirstCol}
\end{align}
\end{subequations}
In particular, from \eqref{SDiagonal} and \eqref{SFirstCol}, one obtains by an immediate induction that $S^n_{j, j} = 1$ for every $j, n \in \mathbb N$. Together with \eqref{SInduction}, one obtains that $S^n_{j, j+1} = 0$ for every $j, n \in \mathbb N$ and, using an immediate inductive argument and \eqref{SInduction}, one obtains that $S^n_{j, k} = 0$ for every $j, k, n \in \mathbb N$ with $k > j$. Moreover, it also follows from \eqref{SDiagonal} and \eqref{SFirstCol} that $S^n_{j, k} = \binom{j}{k}$ whenever $n \in \mathbb N$, $j \in \llbracket 0, n\rrbracket$, and $k \in \{0, j\}$, and using \eqref{SInduction} and an immediate inductive argument, one obtains that this equality also holds for $k \in \llbracket 0, j\rrbracket$.
\end{proof}
}

\chapter{Software for time-delay systems}
\label{chap:delay-software}

\epigraph{\foreignlanguage{brazil}{Software é a parte que você xinga, hardware é a parte que você chuta}}{\foreignlanguage{brazil}{Ditado popular da internet de autoria desconhecida}}

Motivated by the numerous applications of time-delay systems, in particular in Engineering, and by the rather difficult theoretical nature of some results, many researchers have tried to propose and develop toolboxes that facilitate and popularize the use of results in time-delay systems both in academic and industrial contexts. To list just a few of those toolboxes, let us mention the JiTCDDE Python module \cite{Ansmann2018Efficiently}, which uses methods from \cite{Shampine2001Solving, Bogacki1989Pair} to numerically approximate solutions of delay-differential equations; the parDMF Matlab toolbox \cite{Ramirez2021Scalable}, which computes the delay margin of linear neutral functional differential equations with a single delay; the DDE-BIFTOOL Matlab toolbox \cite{Engelborghs2002Numerical, Engelborghs2001DDEBIFTOOL} for numerical bifurcation and stability analysis of delay differential equations with several fixed discrete or state-dependent delays, whose implementation is based on \cite{Verheyden2008Efficient, Samaey2002Numerical, Luzyanina2002Computing, Engelborghs2002StabilityI, Engelborghs2002StabilityII, Luzyanina2001Numerical, Engelborghs2000Collocation}; the TRACE-DDE Matlab toolbox \cite{Breda2009TRACEDDE} for the computation of characteristic roots and stability charts of linear autonomous systems of delay differential equations with discrete and distributed delays, based on \cite{Breda2009Adaptive, Breda2005Pseudospectral, Breda2005Efficient}; and the QPmR Matlab toolbox and its later developments \cite{Vyhlidal2003Quasipolynomial, Vyhlidal2006Mapping, Vyhlidal2014QPmR}, which is a root-finding algorithm developed for computing the spectrum of both retarded and neutral time-delay systems based on \cite{Vyhlidal2009Mapping, Michiels2005Eigenvalue}. We refer the interested reader to \cite{Sipahi2019Mastering, Pekar2018Spectrum} for recent overviews on numerical methods for the stability analysis of time-delay systems and more details on some toolboxes.

In this short chapter, we will provide a brief description of two other recent toolboxes for time-delay systems, YALTAPy and P3$\delta$.

\section{YALTAPy}

This presentation of YALTAPy is extracted from \cite{Cavalera2022YALTAPy}. YAPTAPy builds upon the Matlab toolbox YALTA, which stands for ``Yet Another LTI TDS Algorithm''. YALTA is based on the theoretical work carried out in \cite{Fioravanti2012Numerical, Fioravanti2010Numerical} and the toolbox itself was first described in \cite{Avanessoff2014Hinfty, Avanessoff2013YALTA}. As the ``Py'' in its name suggests, YALTAPy is a Python module implementing the functionalities of YALTA.

\clearpage

\subsection{Scope}

YALTAPy considers time-delay systems whose characteristic function $\Delta \colon \mathbb C \to \mathbb C$ is of the form
\begin{equation}
\label{eq01}
\Delta(s) = P_0(s) + \sum_{k = 1}^N P_k(s) e^{-k \tau s},
\end{equation}
where $\tau > 0$ is a delay, $\alpha \in (0, 1]$ is the fractional exponent, and, for every $k \in \llbracket 0, N\rrbracket$, $s \mapsto P_k(s)$ is a polynomial in the variable $s^\alpha$. We assume that $\deg P_0 \geq \deg P_k$ for every $k \in \llbracket 1, N\rrbracket$, i.e., the system is either of retarded or neutral type.

A \emph{chain of roots} of $\Delta$ is a sequence $(s_n)_{n \in \mathbb N}$ of roots of $\Delta$ such that $\abs{s_n} \to +\infty$ as $n \to +\infty$. We say that a chain of roots is of \emph{retarded type} if $\Real s_n \to -\infty$ as $n \to +\infty$, and of \emph{neutral type} if the sequence $(\Real s_n)_{n \in \mathbb N}$ converges to some real number $c$, which is called the \emph{asymptotic axis} of the chain. As detailed in \cite{Bellman1963Differential, Fioravanti2010Numerical, Partington2004Linear},
\begin{itemize}
\item if $\deg P_0 > \deg P_k$ for all $k \in \llbracket 1, N\rrbracket$, then $\Delta$ only has chains of roots of retarded type; 
\item if $\deg P_0 = \deg P_N$, then $\Delta$ only has chains of roots of neutral type; and
\item if $\deg P_0 = \deg P_k > \deg P_N$ for some $k \in \llbracket 1, N-1\rrbracket$, then $\Delta$ has chains of both neutral and retarded types.
\end{itemize}

Given a characteristic function of the form \eqref{eq01}, YALTAPy can determine
\begin{itemize}
\item whether the corresponding time-delay system is retarded or neutral;
\item the position of asymptotic axes of neutral chains of roots for a neutral system;
\item whether the system is stable when $\tau = 0$;
\item given $\tau_\ast > 0$, intervals $T_1, \dotsc, T_K$ contained in $[0, \tau_\ast]$, called \emph{stability windows}, such that the system is exponentially stable if $\tau \in \bigcup_{k=1}^K T_k$ and not exponentially stable if $\tau \in [0, \tau_\ast] \setminus \bigcup_{k=1}^K T_k$, providing also the number of roots of $\Delta$ with nonnegative real part for each $\tau \in [0, \tau_\ast]$;
\item given $\tau_\ast > 0$, a graph showing all roots of $\Delta$ with nonnegative real part when $\tau = \tau_\ast$, and the movement of those roots in the complex plane as $\tau$ varies in the interval $[0, \tau_\ast]$, known as \emph{root locus}; and
\item a Pad\'e approximation of the system (only for nonfractional systems, i.e., $\alpha = 1$).
\end{itemize}

We refer to \cite{Cavalera2022YALTAPy} for technical details on how to use YALTAPy.

\subsection{Example}
\label{sec:expl-yaltapy}

To illustrate the use of YALTAPy, consider the time-delay system given by
\begin{equation*}
x^\prime(t) = A_0 x(t) + A_1 x(t - \tau),
\end{equation*}
with
\[
A_0 = \begin{pmatrix}
0 & 1 & -1 \\
1 & -\frac{1}{2} & -\frac{1}{4} \\
-1 & 0 & -\frac{4}{5}
\end{pmatrix}, \qquad A_1 = \begin{pmatrix}
-3 & -2 & \frac{1}{2} \\
-1 & 0 & \frac{1}{5} \\
\frac{2}{3} & \frac{1}{3} & 0
\end{pmatrix}, \qquad \tau = 3.
\]
The characteristic function of the above system is given by
\begin{equation*}
\begin{aligned}
\Delta(s) & = \det(s \Id - A_0 - A_1 e^{-s \tau}) \\
 & = s^3 + \frac{13}{10} s^2 - \frac{8}{5} s - \frac{31}{20} + \left(3 s^2 + \frac{163}{20} s + \frac{323}{60}\right) e^{-\tau s} - \left(\frac{12}{5}s + \frac{173}{60}\right) e^{-2 \tau s} + \frac{7}{30} e^{-3 \tau s}.
\end{aligned}
\end{equation*}

\subsubsection{Stability analysis}

A call to the stability analysis functionality of YALTAPy provides information on the stability of the system. YALTAPy correctly identifies the system as being of retarded type, and says that it has 4 roots with nonnegative real part, yielding the conclusion that, with the given parameters, the system is not exponentially stable. It also computes the roots of $\Delta$ when $\tau$ is set to $0$, i.e., the roots of the polynomial $s^3 + \frac{13}{10} s^2 - \frac{8}{5} s - \frac{31}{20}$, which are numerically computed as being $-3.07585975$ and $-0.61207012 \pm 0.10043131i$.

Since the system is of retarded type, it presents no chains of roots of neutral type, a fact that is recalled by YALTAPy. As the system admits 4 roots with nonnegative real part, YALTAPy provides information on how these roots behave as the delays vary from $0$ to the nominal value $3$. YALTAPy states that, at $\tau \approx 0.3820$, two roots cross the imaginary axis, at the values $\pm 3.3750 i$, and the direction of crossing is from the left to the right, and, at $\tau \approx 2.2437$, two other roots cross the imaginary axis, once again at the values $\pm 3.3750 i$ and from the left to the right.

\subsubsection{Stability windows}

\begin{figure}[ht]
\centering
\resizebox{0.75\textwidth}{!}{\input{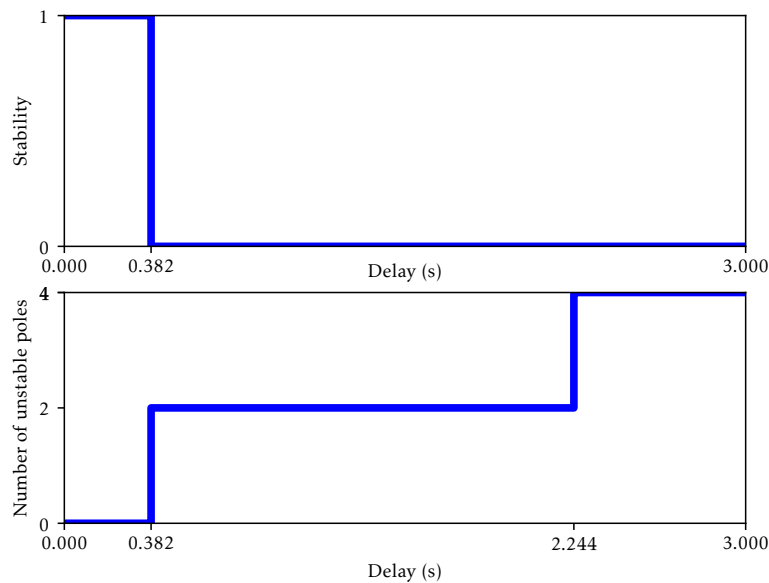}}
\caption[Stability windows for the example from Section~\ref{sec:expl-yaltapy}]{Stability windows computed by YALTAPy for the example from Section~\ref{sec:expl-yaltapy}}
\label{fig:expl1-stab-win}
\end{figure}

YALTAPy provides a function to obtain stability windows for the system under consideration. When called with the value $\tau_\ast = 3$ of maximal delay, it produces the graph from Figure~\ref{fig:expl1-stab-win}. This figure shows that, at $\tau = 0$, the system is exponentially stable (value $1$ at the top plot in Figure~\ref{fig:expl1-stab-win}) and has no roots in the right half-plane (value $0$ at the bottom plot in Figure~\ref{fig:expl1-stab-win}). At $\tau \approx 0.3820$, the system becomes unstable (value $0$ at the top plot of Figure~\ref{fig:expl1-stab-win}) and has $2$ roots in the right half-plane (value $2$ at the bottom plot of Figure~\ref{fig:expl1-stab-win}). At $\tau \approx 2.2437$, two additional roots appear in the complex right half-plane, and the system remains unstable and with $4$ roots in the right half-plane until the delay value $\tau = \tau_\ast = 3$. YALTAPy also returns the relevant values in a table format.

\subsubsection{Root locus}

\begin{figure}[ht]
\centering
\resizebox{0.75\textwidth}{!}{\import{Figures}{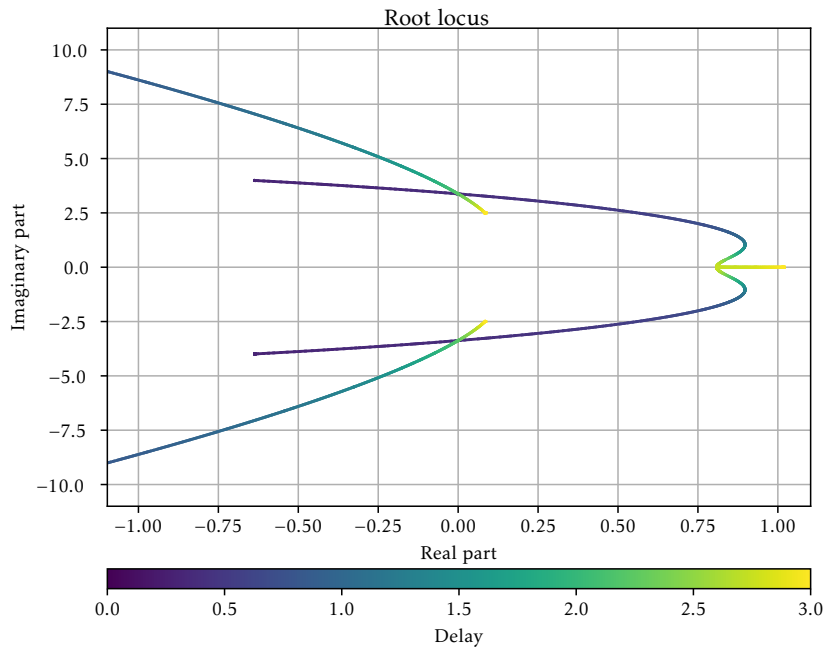}}
\caption[Root locus for the example from Section~\ref{sec:expl-yaltapy}]{Root locus computed by YALTAPy for the example from Section~\ref{sec:expl-yaltapy}}
\label{fig:expl1-root-locus}
\end{figure}

Given $\tau_\ast > 0$, YALTAPy provides a function that can compute the root locus of the roots of $\Delta$ which have nonnegative real part when $\tau = \tau_\ast$, showing how these roots vary in the interval $[0, \tau_\ast]$. For systems of neutral type, this function only runs if all chains of roots of neutral type have negative asymptotic axis. In the present example, with $\tau_\ast = 3$, we obtain the graph from Figure~\ref{fig:expl1-root-locus}, which allows us to observe the behavior of the four roots of the system that become unstable as the delay increases from $0$ to $3$. As for the stability windows, relevant information is also returned by YALTAPy.

\subsection{Further information}

YALTAPy is distributed as a free software under the GPL-3 license. The project website is \url{https://project.inria.fr/yalta/}, its public GitLab project page is \url{https://gitlab.inria.fr/disco/yaltapy-public}, and it can be easily installed in Python since it is available in the Python Package Index (PyPI) at \url{https://pypi.org/project/yaltapy}.

To facilitate the use of YALTAPy by researchers and engineers not familiar with Python, a graphical user interface, called YALTAPy\_Online, is available under the form of a Jupyter Notebook, also released under the GPL-3 license. YALTAPy\_Online can be downloaded from the project public GitLab page, \url{https://gitlab.inria.fr/disco/YALTAPy_Online}, and it can also be used online thanks to the Binder service \cite{ProjectJupyter2018Binder}. Readers of the PDF version of this document using Adobe Acrobat Reader can access YALTAPy\_Online by clicking on \href{https://mybinder.org/v2/git/https%3A%2F%2Fgitlab.inria.fr%2Fdisco%2FYALTAPy_Online/HEAD?urlpath=tree%2FYALTAPy_Online.ipynb}{\raisebox{-0.41em}{\small\color{white}
\begingroup%
  \makeatletter%
  \providecommand\color[2][]{%
    \errmessage{(Inkscape) Color is used for the text in Inkscape, but the package 'color.sty' is not loaded}%
    \renewcommand\color[2][]{}%
  }%
  \providecommand\transparent[1]{%
    \errmessage{(Inkscape) Transparency is used (non-zero) for the text in Inkscape, but the package 'transparent.sty' is not loaded}%
    \renewcommand\transparent[1]{}%
  }%
  \providecommand\rotatebox[2]{#2}%
  \newcommand*\fsize{\dimexpr\f@size pt\relax}%
  \newcommand*\lineheight[1]{\fontsize{\fsize}{#1\fsize}\selectfont}%
  \ifx\svgwidth\undefined%
    \setlength{\unitlength}{81.75bp}%
    \ifx\svgscale\undefined%
      \relax%
    \else%
      \setlength{\unitlength}{\unitlength * \real{\svgscale}}%
    \fi%
  \else%
    \setlength{\unitlength}{\svgwidth}%
  \fi%
  \global\let\svgwidth\undefined%
  \global\let\svgscale\undefined%
  \makeatother%
  \begin{picture}(1,0.18348624)%
    \lineheight{1}%
    \setlength\tabcolsep{0pt}%
    \put(0,0){\includegraphics[width=\unitlength,page=1]{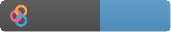}}%
    \put(0.38073394,0.04587156){\color[rgb]{0.00392157,0.00392157,0.00392157}\transparent{0.30000001}\makebox(0,0)[t]{\lineheight{1.25}\smash{\begin{tabular}[t]{c}launch\end{tabular}}}}%
    \put(0.38073394,0.05504587){\makebox(0,0)[t]{\lineheight{1.25}\smash{\begin{tabular}[t]{c}launch\end{tabular}}}}%
    \put(0.78440367,0.04587156){\color[rgb]{0.00392157,0.00392157,0.00392157}\transparent{0.30000001}\makebox(0,0)[t]{\lineheight{1.25}\smash{\begin{tabular}[t]{c}binder\end{tabular}}}}%
    \put(0.78440367,0.05504587){\makebox(0,0)[t]{\lineheight{1.25}\smash{\begin{tabular}[t]{c}binder\end{tabular}}}}%
  \end{picture}%
\endgroup%
}}. The direct link to access YALTAPy\_Online on Binder is also available on the project GitLab page.

\begin{figure}[ht]
\centering
\includegraphics[width=\textwidth]{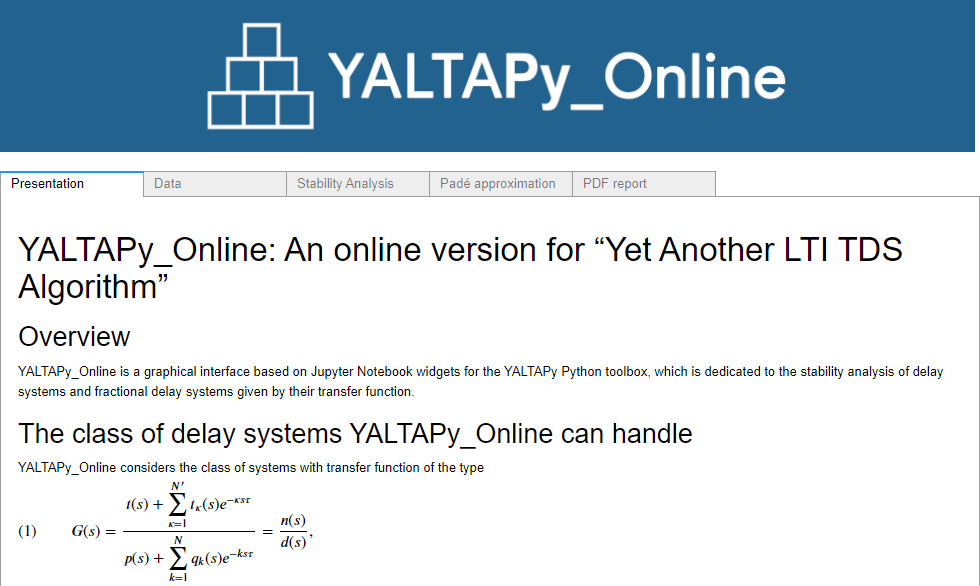}
\caption{Part of the main screen of YALTAPy\_Online}
\label{fig:online}
\end{figure}

YALTAPy\_Online starts by a presentation screen, represented in Figure~\ref{fig:online}, recalling the class of systems that it handles and its main features. It is divided into five tabs: the first one is the presentation screen, the second one is dedicated to the input of the data, the third one provides features of stability analysis, stability windows, and root locus, the fourth one provides a tool for computing Padé approximations, and the last one allows the user to obtain a PDF report with the obtained results.

\section{P3\texorpdfstring{$\delta$}{delta}}

This presentation of P3$\delta$ is based on \cite{Boussaada2022New, Boussaada2021New, Boussaada2020Partial}. P3$\delta$ stands for ``Partial pole placement via delay action'' and it is a Python toolbox implementing the major results on the MID property described in Chapter~\ref{chap:MID}, as well as the CRRID property studied in \cite{Bedouhene2019Pole, Amrane2020Some, Schmoderer2023Quasipolynomials, Schmoderer2024Insights, Amrane2018Qualitative, Bedouhene2020Real}.

\subsection{Scope}

P3$\delta$ is able to consider linear time-invariant differential equations with a single time delay under the form \eqref{eq:DDE}, whose characteristic function is \eqref{eq:Delta-DDE}.

P3$\delta$ has three different available modes, called ``Generic MID'', ``Control-oriented MID'', and ``Generic CRRID''. In the ``Generic MID'' mode, given nonnegative integers $n$ and $m$ with $m \leq n$, a delay $\tau \in \mathbb R_+^\ast$, and $s_0 \in \mathbb R$, P3$\delta$ computes the coefficients of the function $\Delta$ from \eqref{eq:Delta-DDE} so that $s_0$ is a root of maximal multiplicity $D = n + m + 1$ of $\Delta$, according to \eqref{Coeffs} from Theorem~\ref{thm:main-MID}. The software outputs the computed values of the coefficients, as well as the corresponding quasipolynomial $\Delta$ and its factorizations, both according to \eqref{eq:DeltaTilde-as-Kummer} and a similar one using the integral representation of Kummer functions from \eqref{eq:Kummer-integral}. It also shows a graph with the distribution of the spectrum, computed numerically using Python's \texttt{cxroots} module \cite{Parini2018cxroots}, and allows the user to simulate solutions of the corresponding DDE \eqref{eq:DDE} for some families of initial conditions (constant, exponential, polynomial, and trigonometric).

In the ``Control-oriented MID'' mode, P3$\delta$ implements the methods described in Section~\ref{sec:control-oriented}. Given nonnegative integers $n$ and $m$ with $m \leq n$ and the values of the known coefficients of the characteristic function $\Delta$ from \eqref{eq:Delta-DDE}, P3$\delta$ computes first the elimination-produced polynomial $\mathbb P$ from \eqref{eq:EPP}, seen as a function of $s_0 \in \mathbb R$ and $\tau \in \mathbb R_+^\ast$, and plots the curve $\{(s_0, \tau) \in \mathbb R \times \mathbb R_+^\ast \suchthat \mathbb P(s_0, \tau) = 0\}$. With this graph in hand, the user can choose the value of one among $s_0$ or $\tau$, and P3$\delta$ automatically computes the other value according to the elimination-produced polynomial, and then computes the unknown coefficients of $\Delta$ ensuring that $s_0$ is a root of the desired multiplicity. As in the ``Generic MID'' mode, it also produces the graph with the distribution of the spectrum and allows for the simulation of solutions of \eqref{eq:DDE}.

Finally, in the ``Generic CRRID'' mode, P3$\delta$ implements the results from \cite{Bedouhene2019Pole, Amrane2020Some, Schmoderer2023Quasipolynomials, Schmoderer2024Insights, Amrane2018Qualitative, Bedouhene2020Real}. Given nonnegative integers $n$ and $m$ with $m \leq n$, a delay $\tau \in \mathbb R_+^\ast$, $s_0 \in \mathbb C$, and $d \in \mathbb R_+^\ast$, P3$\delta$ computes the coefficients of the function $\Delta$ from \eqref{eq:Delta-DDE} so that $s_0, s_0 - d, \dotsc, s_0 - D + 1$ are real roots of $\Delta$, where $D = n + m + 1$ is the degree of $\Delta$. As in the other two modes, P3$\delta$ outputs the computed values of the coefficients, plots the distribution of the spectrum, and allows the user to simulate solutions of \eqref{eq:DDE}.

\subsection{Example}
\label{sec:expl-p3delta}

To illustrate the use of P3$\delta$, let us apply the ``Generic MID'' mode to the delay-differential equation
\begin{multline}
\label{eq:expl-p3delta-DDE}
y^{\prime\prime\prime\prime}(t) + a_3 y^{\prime\prime\prime}(t) + a_2 y^{\prime\prime}(t) + a_1 y^\prime(t) + a_0 y(t) \\ + \alpha_3 y^{\prime\prime\prime}(t - \tau) + \alpha_2 y^{\prime\prime}(t - \tau) + \alpha_1 y^\prime(t - \tau) + \alpha_0 y(t - \tau) = 0,
\end{multline}
whose characteristic quasipolynomial is
\begin{equation}
\label{eq:expl-p3delta-Delta}
\Delta(s) = s^4 + a_3 s^3 + a_2 s^2 + a_1 s + a_0 + e^{-s\tau} (\alpha_3 s^3 + \alpha_2 s^2 + \alpha_1 s + \alpha_0).
\end{equation}
Let us assume that $\tau = 1$ and that we wish to select the coefficients in order to place a root of multiplicity $D = 8$ at $s_0 = -5$. P3$\delta$ computes numerically that one must have
\begin{align*}
a_0 & = 65, & a_2 & = 30, & \alpha_0 & \approx -31.264074, & \alpha_2 & \approx -0.80855364, \\ 
a_1 & = 20, & a_3 & = 4, & \alpha_1 & \approx -8.4898132, & \alpha_3 & \approx -0.026951788, 
\end{align*}
and that in this case $\Delta$ can be factorized as both
\[
\Delta(s) = \frac{1}{6} (s + 5)^8 \int_0^1 t^3 (1 - t)^4 e^{-t (s + 5)} \diff t \quad \text{ and } \quad \Delta(s) = \frac{1}{1680} (s + 5)^8 M(4, 9, -(s + 5)),
\]
where $M$ is the Kummer function from \eqref{DefiConfluent}.

\begin{figure}[ht]
\centering
\resizebox{\textwidth}{!}{\import{Figures}{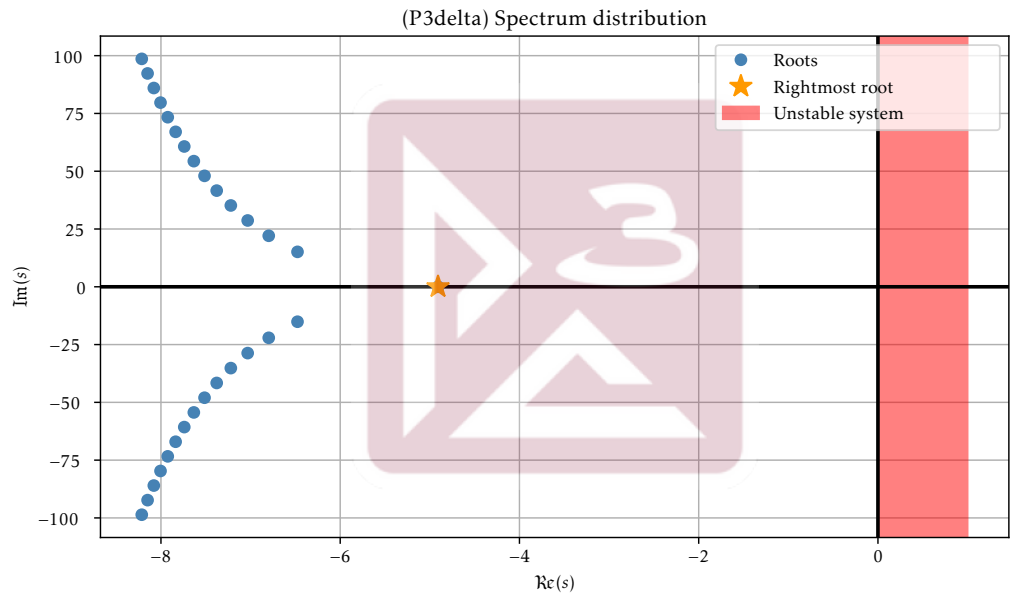}}
\caption[Spectrum distribution for the example from Section~\ref{sec:expl-p3delta}]{Spectrum distribution computed by P3$\delta$ for the example from Section~\ref{sec:expl-p3delta}}
\label{fig:p3delta-spectrum}
\end{figure}

\begin{figure}[ht]
\centering
\resizebox{\textwidth}{!}{\import{Figures}{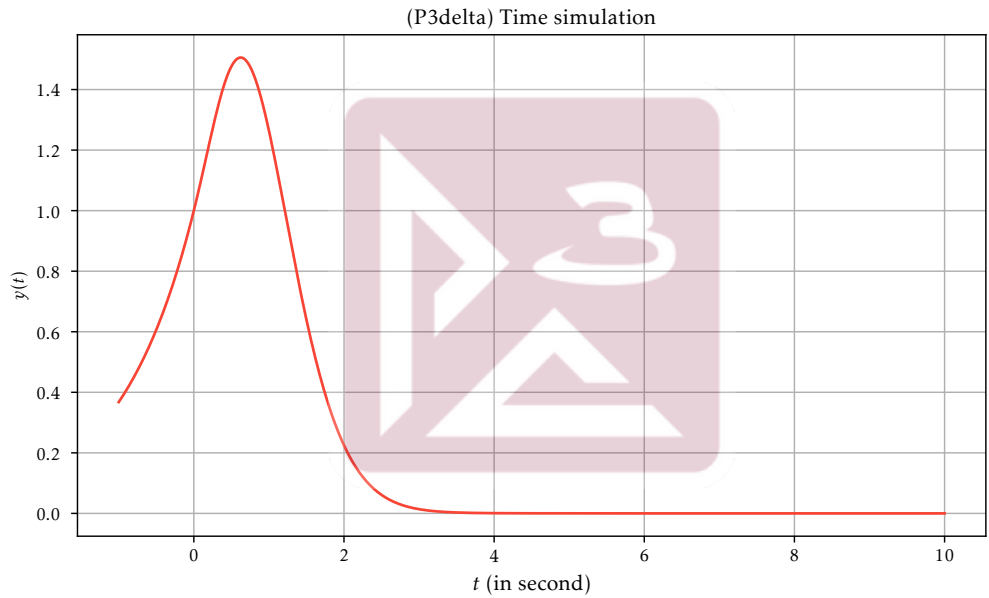}}
\caption[Solution for the example from Section~\ref{sec:expl-p3delta}]{Solution of \eqref{eq:expl-p3delta-DDE} with initial condition $y(t) = e^t$ for $t \in [-1, 0]$ computed by P3$\delta$ for the example from Section~\ref{sec:expl-p3delta}}
\label{fig:p3delta-solution}
\end{figure}

The distribution of the roots of the characteristic function $\Delta$ from \eqref{eq:expl-p3delta-Delta} computed by P3$\delta$ is represented in Figure~\ref{fig:p3delta-spectrum}, which represents the multiple root at $s_0$ by a star and all other roots of $\Delta$ by blue circles. Figure~\ref{fig:p3delta-solution} represents the solution of \eqref{eq:expl-p3delta-DDE} computed numerically by P3$\delta$ with initial condition $y(t) = e^{t}$ for $t \in [-1, 0]$.

\subsection{Further information}

\begin{figure}[ht]
\centering
\includegraphics[width=\textwidth]{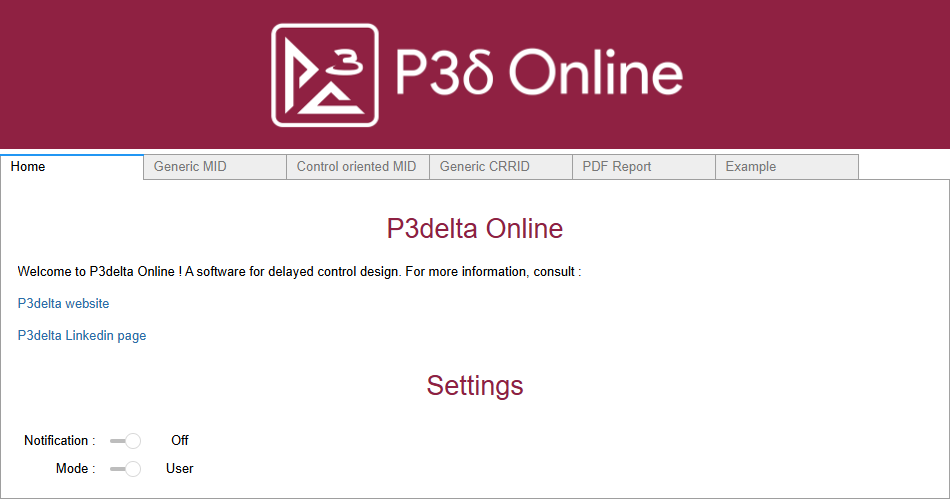}
\caption{Part of the main screen of P3$\delta$ Online}
\label{fig:p3delta-online}
\end{figure}

More information on P3$\delta$, including installation instructions, video demonstrations, the user guide, and download links, are available on the project website, \url{https://cutt.ly/p3delta}. P3$\delta$ has two available graphical user interfaces providing the same functionalities of the software: a standard one and a Jupyter Notebook-based one, the latter being called P3$\delta$ Online. The main screen of P3$\delta$ Online is represented in Figure~\ref{fig:p3delta-online}. This interface is divided into six tabs: a first ``Home'' tab with a welcome message, useful links, and some configurations; the three tabs ``Generic MID'', ``Control oriented MID'', and ``Generic CRRID'', containing the corresponding modes of P3$\delta$; a tab for the generation of a PDF report containing all results computed by P3$\delta$ so far; and a tab with some selected examples illustrating how to use P3$\delta$ and some applications in which it has been used.

\chapter{Open problems and perspectives}
\label{chap:delay-open-problems}

\epigraph{\foreignlanguage{brazil}{Quando alguém faz um discurso\\utiliza esse recurso\\se lhe fazem perguntas\\tchuin tchuin tchun clain}}{\foreignlanguage{brazil}{\emph{Chaves}}}

The results described in this monograph raise many interesting research questions on the stability and control of time-delay systems. In this concluding chapter, we present and discuss some of them.

\section{Minimal controllability time for difference equations}

Theorem~\ref{lem:3} in Chapter~\ref{chap:controllability} provides an upper bound for the minimal times of approximate and exact controllability, $T_{\min,\,\mathrm{ap},\,p}$ and $T_{\min,\,\mathrm{ex},\,p}$ from Definition~\ref{def:Tmin-ap-ex-p}, of the difference equation \eqref{eq:intro-linear-difference-equation}. As detailed in Corollary~\ref{Corollary_min_con}, these times are either $+\infty$ or belong to the set $[\Lambda_N, d \Lambda_N]$ and, as discussed after that corollary, it follows from the results of \cite{Chitour2020Approximate} that, in the case $N = d = 2$, $m = 1$, and $p = 2$, the only possible value for $T_{\min,\,\mathrm{ap},\,p}$ and $T_{\min,\,\mathrm{ex},\,p}$ in the interval $[\Lambda_N, d \Lambda_N]$ is $d \Lambda_N = 2 \Lambda_2$. This is also shown in \cite{Chitour2020Approximate} to be still true for \eqref{eq:intro-linear-difference-equation} in the single-input case, $m = 1$, under the additional assumption that the delays $\Lambda_1, \dotsc, \Lambda_N$ are commensurable. However, \cite{Chitour2020Approximate} also shows that, when $m \geq 2$, $T_{\min,\,\mathrm{ap},\,p}$ and $T_{\min,\,\mathrm{ex},\,p}$ may also take other values in the interval $[\Lambda_N, d \Lambda_N]$.

Several questions arise from these observations. The first one is whether the only possible value for $T_{\min,\,\mathrm{ap},\,p}$ and $T_{\min,\,\mathrm{ex},\,p}$ in the interval $[\Lambda_N, d \Lambda_N]$ is always $d \Lambda_N$ in the single-input case, $m = 1$. The proof provided in \cite{Chitour2020Approximate} for that result in the case of commensurable delays cannot be extended to more general delays, and a strategy of approximating more general delays by commensurable ones does not seem to give sufficient information to obtain a conclusion. A deeper study of the endpoint operator $E_p(T)$ from Definition~\ref{def_Upsilon_E}\ref{item_def_endpoint} seems to be required in order to address this problem.

A second question is on the equality between $T_{\min,\,\mathrm{ap},\,p}$ and $T_{\min,\,\mathrm{ex},\,p}$. Since exact controllability implies approximate controllability, we have the trivial inequality $T_{\min,\,\mathrm{ap},\,p} \leq T_{\min,\,\mathrm{ex},\,p}$. Strict inequality may occur when $T_{\min,\,\mathrm{ex},\,p} = +\infty$ and $T_{\min,\,\mathrm{ap},\,p} < +\infty$. For instance, consider the system with $N = d = 2$ and $m = 1$ given by
\[
x(t) = \begin{pmatrix}
0 & 1 \\
0 & 0
\end{pmatrix} x(t - \ell) + \begin{pmatrix}
1 & 1 \\
0 & 0
\end{pmatrix} x(t - 1) + \begin{pmatrix}
0 \\
1
\end{pmatrix} u(t),
\]
where $\ell \in (0, 1) \setminus \mathbb Q$. Using the results from Section~\ref{sec:N=D=2,m=1}, we deduce that $T_{\min,\,\mathrm{ap},\,p} \leq 2$ for every $p \in [1, +\infty)$, but $T_{\min,\,\mathrm{ex},\,1} = +\infty$ (and we also have $T_{\min,\,\mathrm{ex},\,2} = +\infty$ thanks to \cite[Theorem~4.1]{Chitour2020Approximate}). Thus, we have the strict inequality $T_{\min,\,\mathrm{ap},\,p} < T_{\min,\,\mathrm{ex},\,p}$ at least for $p \in \{1, 2\}$. The natural question is whether the only case where strict inequality may occur is when $T_{\min,\,\mathrm{ex},\,p} = +\infty$. In other words, assuming $T_{\min,\,\mathrm{ex},\,p} < +\infty$, is it true that $T_{\min,\,\mathrm{ap},\,p} = T_{\min,\,\mathrm{ex},\,p}$? The answer is true in the case of commensurable delays, since approximate and exact controllability are equivalent in that case (see, e.g., \cite[Proposition~3.3]{Chitour2020Approximate}), but an analysis of the general case is still missing.

Finally, one may wonder whether $T_{\min,\,\mathrm{ap},\,p}$ and $T_{\min,\,\mathrm{ex},\,p}$ actually depend on $p$, at least for $p \in [1, +\infty)$. The analysis of particular cases in \cite{Chitour2020Approximate} and in Section~\ref{sec:applications} has not produced any example in which these quantities may depend on $p$, and the particular structure of the endpoint operator from Definition~\ref{def_Upsilon_E}\ref{item_def_endpoint}, which depends on $p$ only through the choice of the domain and the codomain of the operator, may lead one to conjecture that approximate and exact controllability do not depend on $p$, and that $T_{\min,\,\mathrm{ap},\,p}$ and $T_{\min,\,\mathrm{ex},\,p}$ may thus be both constant on $p$, at least for $p \in [1, +\infty)$. Whether this is true remains an open problem.

\section{Exact controllability criteria for difference equations in \texorpdfstring{$L^p$}{Lp} for \texorpdfstring{$p \in (1, +\infty)$}{p in (1, infinity)}}

Theorem~\ref{main_result2} provides a necessary condition for the exact controllability of \eqref{eq:intro-linear-difference-equation} in $L^1$, and this condition was recently shown to be sufficient in \cite{Fueyo2024Corona}. However, the characterization of exact controllability in $L^p$ for $p \in (1, +\infty)$ remains an open problem. As detailed in Section~\ref{sec:N=D=2,m=1}, the analysis of \eqref{eq:intro-linear-difference-equation} in the case $N = d = 2$, $m = 1$, and $p = 2$ done in \cite{Chitour2020Approximate} provides the same criterion as Theorem~\ref{main_result2} in that particular case, showing that, at least in this case, exact controllability in time $2 \Lambda_2$ in $L^1$ is equivalent to exact controllability in the same time in $L^2$. Whether this is a more general fact is still an open problem. Note that, combining Corollaries~\ref{cor:distru_exac_equiv} and \ref{corollaire_L^1_equiv_radon}, one deduces that $L^1$ exact controllability implies $L^p$ exact controllability for every $p \in [1, +\infty)$, and hence the open question concerns the converse implication.

Equivalently, another way to view this question is to ask whether the solvability of Bézout's identity \eqref{eq:exact_cont0}, in the sense of Definition~\ref{def:solvable-Bezout}, is equivalent to exact controllability in $L^p$ for any $p \in [1, +\infty)$. That such a solvability property implies exact controllability in $L^p$ is a consequence of Proposition~\ref{prop_fond_exact_controllability} and Corollary~\ref{cor:distru_exac_equiv}, but the converse implication was proved only for $p = 1$ in Theorem~\ref{prop_fond_exact_controllability_L1}. As explained in Remark~\ref{remark1}, the proof of Theorem~\ref{prop_fond_exact_controllability_L1} uses the fact that $p = 1$ in a crucial way, and adapting this proof for any $p$ seems out of reach. New ideas would probably be required in order to obtain advances in this question.

\section{Controllability of difference equations in more general functional spaces}

We have considered, in Chapter~\ref{chap:controllability}, approximate and exact controllability of \eqref{eq:intro-linear-difference-equation} in $L^p$ for $p \in [1, +\infty)$. A natural question is whether one can also characterize these controllability properties in other functional spaces, such as in $L^\infty$, in $\mathcal C^k$, or in the Sobolev spaces $W^{k, p}$.

Regarding $L^\infty$, Theorem~\ref{lem:3} already applies to $L^\infty$, but Theorem~\ref{main_result1} requires a finite $p$, while Theorem~\ref{main_result2} only works for $p = 1$. Note that, in the proof of the equivalences of \ref{item:Yamamoto-approx-2}--\ref{item:Yamamoto-approx-Hautus} in Theorem~\ref{th_Yamamoto1}, the proof that \ref{item:Yamamoto-approx-p} implies \ref{item:Yamamoto-approx-Bezout} also holds true for $p = +\infty$. Indeed, roughly speaking, in this proof, one needs to transform the stronger notion of convergence in $L^p$ into the weaker notion of convergence in the sense of distributions for the approximate Bézout identity \eqref{eq:approx-Bezout}, and this is possible for any $p \in [1, +\infty]$. As a consequence, conditions \ref{assumption1bis} and \ref{assumption2bis} from Theorem~\ref{main_result1} are necessary for $L^\infty$ approximate controllability of \eqref{eq:intro-linear-difference-equation} in time $d \Lambda_N$. However, the proof that \ref{item:Yamamoto-approx-Bezout} implies \ref{item:Yamamoto-approx-p} in Theorem~\ref{th_Yamamoto1} uses the density of $X^{Q, p} \cap \mathcal C^\infty(\mathbb R_+, \mathbb R^d)$ in $X^{Q, p}$, which is not true for $p = +\infty$. Whether conditions \ref{assumption1bis} and \ref{assumption2bis} from Theorem~\ref{main_result1} are also sufficient for $L^\infty$ approximate controllability is an open problem.

Controllability notions in spaces with higher regularity, such as $\mathcal C^k$ and $W^{k, p}$, require compatibility conditions, in the sense of Definition~\ref{def:compatibility}, to be satisfied for the initial state $x_0$ and the target state $\phi$. While the approximate controllability problem may be addressed using density arguments similar to those used in the proof of Theorem~\ref{th_Yamamoto1}, exact controllability in such higher regularity spaces is a much more delicate question. A partial answer is provided, in the case of $\mathcal C^k$ spaces, in Theorem~\ref{thm:Ck}. However, this theorem does not provide a full response to the exact controllability problem since, in the one hand, controls are taken in $L^p$, and not necessarily in $\mathcal C^k$, and, in the other hand, the regularity index $k$ is not arbitrary. Whether one can characterize exact controllability in $\mathcal C^k$ or $W^{k, p}$ by Hautus-type conditions similar to \ref{assumption1-L1} and \ref{assumption2-L1} from Theorem~\ref{main_result2} is an open question. In particular, Theorem~\ref{thm:Ck} leads to the interesting conjecture that \ref{assumption1-L1} and \ref{assumption2-L1} from Theorem~\ref{main_result2} may also be necessary and sufficient conditions for exact controllability in $\mathcal C^k$.

\section{Controllability of more general systems using realization theory}

Chapter~\ref{chap:controllability} illustrated the use of Y.~Yamamoto's works on realization theory for obtaining controllability criteria for difference equations. The key aspect is that, contrarily to the classical approach of studying controllability of linear infinite-dimensional systems by duality through an observability inequality or a unique continuation principle (as described in details, for instance, in \cite[Section~2.3.2]{Coron2007Control}), Y.~Yamamoto's works allow one to characterize controllability in terms of Bézout's identities, as done in Theorems~\ref{th_Yamamoto1} and \ref{prop_fond_exact_controllability_L1} and Proposition~\ref{prop_fond_exact_controllability}. A natural question is whether these characterizations can be exploited to obtain controllability results for other classes of systems.

Similar controllability results can be obtained for a class of linear systems of hyperbolic PDEs in one space dimension without in-domain coupling between the components of the state vector, as detailed in \cite{Chitour2024Approximate}. Recent results have also extended the controllability criteria presented in this monograph to a class of difference equations containing also distributed delays \cite{Fueyo2025Left, Fueyo2025Lq}. There still remains the question of finding more general classes of systems for which the techniques of realization theory can be applied to obtain controllability results. For linear difference equations, a natural extension would be to consider systems of the form
\[
x(t) = \int_{-\tau}^0 \diff A(\theta) x(t + \theta) + \int_{-h}^0 \diff B(\theta) u(t + \theta),
\]
where $\tau$ and $h$ are positive delays and $A$ and $B$ are matrix-valued Radon measures defined over $[-\tau, 0]$ and $[-h, 0]$, respectively, and taking values in $\mathcal M_d(\mathbb R)$ and $\mathcal M_{d, m}(\mathbb R)$, respectively. The results from Chapter~\ref{chap:controllability} correspond to the case where $B$ is a Dirac mass at $0$ multiplied by some constant matrix from $\mathcal M_{d, m}(\mathbb R)$ and $A$ is a sum of $N$ Dirac masses, whereas those from \cite{Fueyo2025Left, Fueyo2025Lq} consider also the more general situation in which $A$ can be the sum of finitely many Dirac masses with an absolutely continuous measure with respect to the Lebesgue measure, still with the same structure for $B$. Other possible extensions could also concern linear hyperbolic systems of PDEs in one space dimension with in-domain coupling terms, and it would be interesting to investigate whether more general classes of PDEs could be fit into the framework of realization theory in order to obtain controllability criteria based on Bezout's identities.

\section{Spectrum-determined growth condition for difference equations}
\label{sec:open-SDG}

The main result of Chapter~\ref{chapter:strong-stability} in arbitrary dimension $d$ is Theorem~\ref{GHS-Theorem}, which considers strong stability in the class $\SDG$ from \eqref{eq_def-N-L-G}, namely, the class of functions $M \in \NBV([-1, 0],\allowbreak \mathcal M_d(\mathbb R))$ satisfying the spectrum-determined growth condition from Definition~\ref{def_big-def}\ref{item:def-sdg}. While it is known that the class $\SDG$ contains both the class of functions $\NoSingular$ from \eqref{eq_def-N-L-G} and the class of functions $M$ described in Remark~\ref{remk:big-remk}\ref{item_remk-sdg-condition-H}, a full characterization of $\SDG$ is still missing from the literature.

Some insight on this question can be obtained through the analysis of the proofs that $\NoSingular \subset \SDG$ in \cite[Chapter~12, Theorem~3.4]{Hale1977Theory} and of Remark~\ref{remk:big-remk}\ref{item_remk-sdg-condition-H} in \cite[Chapter~9, Theorem~4.4]{Hale1993Introduction}. Recall that the inequality $S_M \leq \omega_M$ is verified for all $M \in \WellPosed$ and, to prove the converse inequality, one should show that the exponential rate of a solution of \eqref{eq_NS} can be bounded in terms of the spectral abscissa $S_M$.

In \cite[Chapter~12]{Hale1977Theory}, the author first proves in \cite[Chapter~12, Lemma~3.3]{Hale1977Theory} the spec\-trum-determined growth condition for functions $M$ of the form \eqref{eq_NS} with $N = 0$, i.e., when $M$ is purely discontinuous. The main strategy of the proof is to remark first that any solution of \eqref{sys_general} can be written in terms of a fundamental matrix solution $Y$, and that an exponential rate on the latter yields the same exponential rate for solutions. It turns out that the Laplace transform of $Y$ can be written in terms of the characteristic function $\Delta_M$ of $M$ and, more precisely, the poles of the Laplace transform of $Y$ coincide (both in their locations and their multiplicities) with the zeros of $\Delta_M$, up to an additional pole at $0$. The analysis of the exponential rate of $Y$ can be done by using the inverse Laplace transform, and key facts in this analysis are that \begin{enumerate*}[label={(\roman*)}] \item $\frac{1}{\Delta_M(\cdot)}$ is analytic in $\{s \in \mathbb C \suchthat \Real s > S_M\}$ and \item\label{item:open-lower-bound} $\mathbb R \ni y \mapsto \frac{1}{\Delta_M(\alpha + i y)}$ is almost periodic for every $\alpha > S_M$ such that $\inf_{y \in \mathbb R} \abs{\Delta_M(\alpha + i y)} > 0$.\end{enumerate*} The conclusion then follows thanks to the theory of almost periodic functions (see, e.g., \cite{Bohr1947Almost, Corduneanu1968Almost} or the lecture notes \cite{Rodrigues2013Oscilacoes}), and the general case of $M$ of the form \eqref{eq_NS} is then treated by a perturbation argument, using the fact that the term corresponding to the absolutely continuous part $N$ of $M$ in the expression of the characteristic function $\Delta_M(s)$ tends to $0$ as $\abs{\Imag s} \to +\infty$. Such a perturbation technique would not work if $M$ had a singular part, since the corresponding term in $\Delta_M(s)$ might not tend to zero as $\abs{\Imag s} \to +\infty$ (see, e.g., \cite{Hu2002Fourier, Wiener1938Singular, Morrison1995Cosine}).

Note that the lower bound on $\inf_{y \in \mathbb R} \abs{\Delta_M(\alpha + i y)}$ from item~\ref{item:open-lower-bound} above is reminiscent of the lower bound condition from the definition of the class $\LowerBound$ in \eqref{eq_def-N-L-G}. At the light of this fact, a natural question is also whether there is some relation between the classes $\LowerBound$ and $\SDG$, for instance if one is included in the other or even if both turn out to be equal.

A similar perturbation argument based on almost periodic functions is also used in \cite[Chapter~9]{Hale1993Introduction}. More precisely, \cite[p.~271]{Hale1993Introduction} remarks that the structural assumption on $M$ from Remark~\ref{remk:big-remk}\ref{item_remk-sdg-condition-H} ensures that its characteristic function $\Delta_M$ is asymptotically almost periodic, i.e., it approaches an almost periodic function as $\abs{\Imag s} \to +\infty$. The link between the Laplace transform of the fundamental matrix solution and the characteristic function $\Delta_M$ is then once again exploited in order to obtain information on the exponential behavior of solutions from the information on the location of zeros of $\Delta_M$. An important estimate in this step is that of \cite[Chapter~9, Lemma~3.2]{Hale1993Introduction}, which provides upper and lower bounds for $\abs{\Delta_M(s)}$ when $s$ is not in a neighborhood of a zero of $\Delta_M$.

Combining the results from \cite[Chapter~12]{Hale1977Theory} and \cite[Chapter~9]{Hale1993Introduction}, one deduces that any $M \in \WellPosed \setminus \SDG$ should necessarily contain a singular part, and that the support of this singular part should contain $0$. It is not clear, however, how such a singular part concentrated around $0$ could prevent the spectrum-determined growth condition to hold true, and the analysis of simple examples, such as a singular part given by a Cantor-type function, does not seem to bring additional information on this question.

\section{Analysis of the impact of the singular part of the measure for a difference equation}

Following the discussion from Section~\ref{sec:open-SDG}, a natural question that arises is on the impact of the singular part of $M$ in the dynamics of the difference equation \eqref{sys_general}. It is know from \cite[Chapter~12]{Hale1977Theory} that the absence of a singular part ensures that $M$ satisfies the spectrum-determined growth condition, while, from Proposition~\ref{prop_positiveinf}, we also know that the absence of a singular part yields a lower bound on the modulus of the characteristic function on right half-planes to the right of the spectral abscissa. It is not clear, however, what happens to the spectrum-determined growth condition and to such a lower bound property in the presence of a singular part.

Following the discussion from Section~\ref{sec:open-SDG}, one expects that a singular part preventing the spectrum-determined growth condition to hold true should necessarily have $0$ in its support. In addition, the proof of Proposition~\ref{prop_positiveinf} makes use of the fact that the characteristic function of $M \in \NoSingular$ is asymptotically almost periodic, and this is also known to be the case if $M$ satisfies instead the assumption from Remark~\ref{remk:big-remk}\ref{item_remk-sdg-condition-H}, so that a function $M \in \WellPosed \setminus \LowerBound$ should also have $0$ in the support of its singular part.

One is the led to study functions $M \in \WellPosed$ whose singular part have $0$ in its support. An important step towards a greater understanding of the role of the singular part of $M$ in the dynamics of \eqref{sys_general} would be to exhibit examples of such $M$ that do not belong to $\LowerBound$ or to $\SDG$. Such examples would be of help in understanding what is the phenomenon created by the singular part of the measure that prevents $M$ from being in those classes.

\section{Difference equations in other functional spaces}

Difference equations can be written in many different functional settings, in particular according to the intended application of the theory. The classical theory of difference equations was developed mostly due to its relations with time-delay systems of neutral type (see, e.g., \cite{Melvin1974Stability, Cruz1970Stability, Henry1974Linear, Hale1977Theory, Hale1993Introduction}). As the latter are usually studied in spaces of continuous functions, many works have considered continuous solutions of difference equations, as it was done in Chapters~\ref{chapter:strong-stability} and \ref{chap:diff-eqn-varying-delay} of this monograph. On the other hand, difference equations were also studied due to their link with some partial differential equations of hyperbolic type in space dimension one, as recalled in Section~\ref{sec:intro-cont-time-diff-eqn}. Since such hyperbolic PDEs are typically studied in Sobolev or $L^p$ spaces, this has motivated the study of difference equations in these kinds of spaces \cite{Chitour2016Stability, Chitour2020Approximate, Auriol2019Explicit, Chitour2023Hautus, Mazanti2017Relative, Bastin2016Stability}, as done in Chapters~\ref{chap:controllability} and \ref{chap:diff-eqn-varying-delay} of this monograph. There is also an interest on difference equations, time-delay systems, and more general functional equations on function spaces of everywhere defined functions allowing for discontinuities, and a natural functional setting for that is to consider regulated solutions \cite{Gallegos2024Measure, Henriquez2024Existence, Honig1975Volterra, Prepelita1995Linear, McNabb1974Initial, Afonso2022Periodic, Slavik2013Measure}, as done in Chapter~\ref{chap:diff-eqn-varying-delay}.

As highlighted in Chapter~\ref{chap:diff-eqn-varying-delay} in the context of time-varying delays, assumptions for existence and uniqueness of solutions can be highly dependent on the functional space, and the notions of exponential stability do not necessarily coincide for different functional spaces, as highlighted in Example~\ref{expl:not-exponentially-stable-Lp}. However, for difference equations of the form \eqref{eq:intro-diff-eqn-sum}, exponential stability is equivalent in all $L^p$ spaces, $p \in [1, +\infty]$ (as shown in \cite{Chitour2016Stability} for the more general case of time-varying matrices), and it is also equivalent to exponential stability in the space of continuous solutions. Such a result was recently generalized in \cite{BraunSpectral} to difference equations of the form \eqref{sys_general} with $M$ belonging to the class $\NoSingular$ from \eqref{eq_def-N-L-G}, based on the links of such equations with Volterra integral equations and results for the latter from \cite{Gripenberg1990Volterra}. An interesting open problem is whether this equivalence of exponential stability also holds for everywhere-defined solutions with less regularity, such as regulated solutions or solutions of (locally) bounded variation.

Note that, when $M$ is of bounded variation, the integral in \eqref{sys_general} is well-defined in the Riemann--Stieltjes or the Lebesgue--Stieltjes senses for all $t \in \mathbb R_+$ if one assumes $x$ to be continuous, but this is not necessarily the case if $x$ is in $L^p$, $p \in [1, +\infty]$, regulated, or of bounded variation. A way to overcome this difficulty is to regard the map $t \mapsto \int_{-1}^0 \diff M(\theta) x(t + \theta)$ as a convolution, and then use results on convolution of measures of bounded variation and functions in $L^p$ or of bounded variation in order to deduce that this map is well-defined for almost every $t$. It is in this sense, for instance, that \cite{BraunSpectral} considers solutions of \eqref{sys_general}, which are thus required to satisfy this equation only almost everywhere, similarly to the setting in \cite{Gripenberg1990Volterra} for Volterra equations. Even though the norm of a function in $L^p$ does not change if the function is modified in a set of Lebesgue measure zero, this is not the case for the standard norm in the space of functions of bounded variation. Hence, in order to consider everywhere-defined discontinuous solutions of \eqref{sys_general} in the spaces of regulated functions or functions of bounded variation, one should provide a meaning to the integral in \eqref{sys_general} for all $t \in \mathbb R_+$. This can be done by considering more general integral notions than those of Riemann--Stieltjes or Lebesgue--Stieltjes. This is an important premilinary step before considering the exponential stability of \eqref{sys_general} in such spaces.

\section{Conjectures on the strong stability of linear difference equations}
\label{sec:conjectures-strong-stab}

While Chapter~\ref{chapter:strong-stability} provides many results on the strong stability of System~$\Sigma(M)$ from \eqref{sys_general}, many open problems are left unanswered, such as Conjectures~\ref{GHS-Conjecture} and \ref{GHS-Conjecture-SDG}.

Proving Conjecture~\ref{GHS-Conjecture-SDG} amounts to showing that $\SDG$-local (respectively, global) strong stability with uniform rate is equivalent to $\SDG$-local (respectively, global) strong stability. A direct proof of this fact would immediately establish Conjecture~\ref{GHS-Conjecture-SDG} as a consequence of Theorem~\ref{GHS-Theorem}, but such a direct proof does not seem easy. In other contexts, such as in switched systems, similar uniformity properties can be obtained by relying on results such as Fenichel's lemma (see \cite[Lemma~1.25]{Chitour2025Dynamics} or \cite[Lemma~5.2.7]{Colonius2000Dynamics}), but it is not clear whether these techniques can be adapted to the present setting. An alternative approach would be to show that $\SDG$-local strong stability implies $\rho_{\mathrm{HS}}(M) < 1$, as this would both complete the proof of Conjecture~\ref{GHS-Conjecture-SDG} and imply equivalence between strong stability with uniform rate or not. However, the only partial result that we presented in this direction, Corollary~\ref{coro:lss-implies-rho-leq-1}, only provides the weaker conclusion $\rho_{\mathrm{HS}}(M) \leq 1$, and it is not clear how its proof, or that of Lemma~\ref{lemma_roots}, could be reinforced in order to obtain the desired stronger conclusion.

As for Conjecture~\ref{GHS-Conjecture} (or even a version of it with the notions of strong stability being replaced by the corresponding ones with uniform rate), one of the main difficulties to establish it is that, by working on $\WellPosed$ instead of $\SDG$, spectral methods are not available, but the quantity $\rho_{\mathrm{HS}}(M)$ from \eqref{GHSC} is of spectral nature. It is not easy to imagine how to relate $\rho_{\mathrm{HS}}(M)$ directly with asymptotic properties of the system without going through spectral methods, but it has to be pointed out that this was done in \cite{Chitour2016Stability} in the particular case of \eqref{eq:intro-diff-eqn-sum}, in which case $\rho_{\mathrm{HS}}(M)$ reduces to \eqref{eq:rho-HS}. The strategy used in \cite{Chitour2016Stability} relies, however, on the fact that \eqref{eq:intro-diff-eqn-sum} contains only a finite sum, deducing then a representation formula for solutions of \eqref{eq:intro-diff-eqn-sum} that relies on products of the matrices $A_1, \dotsc, A_N$. These same products are retrieved from \eqref{eq:rho-HS} when using Gelfand's formula for the spectral radius (see \cite[Proposition~3.31]{Chitour2016Stability}), allowing one to relate \eqref{eq:rho-HS} directly with the representation formula for solutions of \eqref{eq:intro-diff-eqn-sum} without going through spectral methods. Generalizing these ideas to \eqref{sys_general} seems to be a challenging problem.

\begin{figure}[ht]
\centering
\begin{tikzpicture}
\node (rhoHS) at (0, 0) {\small $\rho_{\mathrm{HS}}(M) < 1$};
\node (V) at (0, -3) {\small $\V(M) < 1$};
\node[text width=3.5cm, align = center] (SS) at (5, 0) {\small $\Sigma(M)$ is $\WellPosed$-strongly stable};
\node[text width=3.5cm, align = center] (LSS) at (10, 0) {\small $\Sigma(M)$ is $\WellPosed$-locally strongly stable};
\node[text width=3.5cm, align = center] (SSUR) at (5, -3) {\small $\Sigma(M)$ is $\WellPosed$-strongly stable with uniform rate};
\node[text width=3.5cm, align = center] (LSSUR) at (10, -3) {\small $\Sigma(M)$ is $\WellPosed$-locally strongly stable with uniform rate};

\draw[double equal sign distance, double, Implies-Implies] (rhoHS) --node[midway, above] {\small Thm.~\ref{HS:scalar}} (SS);
\draw[double equal sign distance, double, Implies-Implies] (SS) --node[midway, above] {\small Thm.~\ref{HS:scalar}} (LSS);
\draw[double equal sign distance, double, -Implies] (SSUR) -- (SS);
\draw[double equal sign distance, double, -Implies] (LSSUR) -- (LSS);
\draw[double equal sign distance, double, -Implies] (SSUR) -- (LSSUR);
\draw[double equal sign distance, double, Implies-Implies] (V) --node[midway, above, sloped] {\small Rmk.~\ref{remk:rho-HS-TV}} (rhoHS);
\draw[double equal sign distance, double, -Implies] (V) --node[midway, above] {\small Rmk.~\ref{remk:TV-less-1-implies-global-strong-stab-unif-rate}} (SSUR);
\end{tikzpicture}
\caption{Diagram of the results of Chapter~\ref{chapter:strong-stability} for $M \in \WellPosed$ in dimension $d = 1$}
\label{fig:diagram-dim-1}
\end{figure}

We point out that, in the scalar case $d = 1$, $\WellPosed$-local (respectively, global) strong stability is equivalent to $\WellPosed$-local (respectively, global) strong stability with uniform rate. Indeed, the results of Chapter~\ref{chapter:strong-stability} give the implications from Figure~\ref{fig:diagram-dim-1} for any $M \in \WellPosed$, where implication arrows with no labels correspond to trivial implications. It is immediate to deduce the stated equivalence from this diagram.

\begin{figure}[ht]
\centering
\begin{tikzpicture}
\node (rhoHS) at (10, 0) {\small $\rho_{\mathrm{HS}}(M) < 1$};
\node[text width=3.5cm, align = center] (WPSS) at (0, 0) {\small $\Sigma(M)$ is $\WellPosed$-strongly stable};
\node[text width=3.5cm, align = center] (WPLSS) at (5, 0) {\small $\Sigma(M)$ is $\WellPosed$-locally strongly stable};
\node[text width=3.5cm, align = center] (KSS) at (0, -3) {\small $\Sigma(M)$ is $\mathsf K$-strongly stable};
\node[text width=3.5cm, align = center] (KLSS) at (5, -3) {\small $\Sigma(M)$ is $\mathsf K$-locally strongly stable};

\draw[double equal sign distance, double, Implies-Implies] (rhoHS) --node[midway, above] {\small Thm.~\ref{HS:scalar}} (WPLSS);
\draw[double equal sign distance, double, Implies-Implies] (WPSS) --node[midway, above] {\small Thm.~\ref{HS:scalar}} (WPLSS);
\draw[double equal sign distance, double, -Implies] (WPSS) -- (KSS);
\draw[double equal sign distance, double, -Implies] (WPLSS) -- (KLSS);
\draw[double equal sign distance, double, -Implies] (KSS) -- (KLSS);
\draw[double equal sign distance, double, -Implies] (KLSS) --node[midway, above, sloped] {\small $\star$} (rhoHS);
\end{tikzpicture}
\caption[Diagram of the results of Chapter~\ref{chapter:strong-stability} for $M \in \mathsf K$ in dimension $d = 1$]{Diagram of the results of Chapter~\ref{chapter:strong-stability} for $M \in \mathsf K$ in dimension $d = 1$, where $\mathsf K \subset \WellPosed$ contains all functions of the form \eqref{eq_M-finite-sum}}
\label{fig:diagram-K-SS}
\end{figure}

We also remark that, still in the scalar case $d = 1$, $\WellPosed$-local (respectively, global) strong stability is equivalent to $\mathsf K$-local (respectively, global) strong stability for every $\mathsf K \in \{\NoSingular,\allowbreak \LowerBound,\allowbreak \SDG\}$ (and, more generally, for any class $\mathsf K \subset \WellPosed$ containing all functions $M \in \NBV([-1, 0],\allowbreak \mathcal M_d(\mathbb R))$ of the form \eqref{eq_M-finite-sum} for $N \in \mathbb N^\ast$, $A_1, \dotsc, A_N$ in $\mathcal M_d(\mathbb R)$, and $\tau_1, \dotsc, \tau_N$ in $(0, 1)$). Indeed, in this case, we have the implications from Figure~\ref{fig:diagram-K-SS} where, once again, implication arrows with no labels correspond to trivial implications. As for the implication denoted by $\star$, one has to prove the implication \ref{item_main-loc} $\implies$ \ref{item_main-HS} from Theorem~\ref{HS:scalar} by replacing $\WellPosed$-local strong stability with $\mathsf K$-local strong stability. By inspecting the proof of this implication in Theorem~\ref{HS:scalar}, one notices that the function $\varphi \in \mathbf B$ constructed in that proof satisfies $\varphi_{\#} \mu \in \mathsf K$ for any class $\mathsf K$ as above, since the function of bounded variation corresponding to $\varphi_{\#} \mu$ is of the form \eqref{eq_M-finite-sum} by construction. An interesting open problem is whether this kind of equivalence is also true in higher dimension.

\section{Stability criteria for difference equations with time-varying delays}
\label{sec:open-diff-eqn-many-delays}

Chapter~\ref{chap:diff-eqn-varying-delay} has considered the stability analysis of the difference equation \eqref{eq:diff-eqn-N1-time} with a time-varying delay. One of the first main questions that can be asked is whether the analysis carried out in Chapter~\ref{chap:diff-eqn-varying-delay} can be extended to difference equations with finitely many time-varying delays, of the form
\begin{equation}
\label{eq:diff-eqn-time}
x(t) = \sum_{j = 1}^N A_j x(t - \tau_j(t)),
\end{equation}
where $N \in \mathbb N^\ast$ and, for every $j \in \llbracket 1, N\rrbracket$, $A_j \in \mathcal M_d(\mathbb R)$ and $\tau_j\colon \mathbb R_+ \to \mathbb R_+^\ast$. It is known, on the one hand, that the Hale--Silkowski condition from \eqref{eq:rho-HS} is not sufficient to ensure the exponential stability of \eqref{eq:diff-eqn-time}, as detailed in Example~\ref{expl:AvellarMarconato}, and one can easily show, by adapting the proof of Theorem~\ref{thm:exp-to-0}, that \eqref{eq:diff-eqn-time} is exponentially stable for continuous solutions if $\sum_{j=1}^N \abs{A_j} < 1$ for some norm $\abs{\cdot}$ in $\mathbb R^d$ and if the functions $\tau_1, \dotsc, \tau_N$ are bounded. A natural question is on the sharpness of the condition $\sum_{j=1}^N \abs{A_j} < 1$ for some norm $\abs{\cdot}$ in $\mathbb R^d$: can one find weaker conditions on $A_1, \dotsc, A_N$ that still ensure the exponential stability in the space of continuous functions, still under the assumption that $\tau_1, \dotsc, \tau_N$ are bounded functions? Another question of interest concerns conditions for the exponential stability of \eqref{eq:diff-eqn-time} in other functional spaces, including the spaces of regulated functions and of $L^p$ functions considered in Chapter~\ref{chap:diff-eqn-varying-delay}. Note that, by proceeding similarly to Section~\ref{sec:appli-varying-delay.2}, stability results for \eqref{eq:diff-eqn-time} would also have applications to transport equations of the form \eqref{eq:transport} in which each component of $u$ is transported at a different velocity, i.e., in which the real-valued function $\lambda$ is replaced by a diagonal matrix-valued function whose diagonal entries are positive.

More generally, one may also consider a variant of \eqref{eq:diff-eqn-time} in which the matrices $A_1, \dotsc, A_N$, instead of being constant, also depend on time. This question was addressed in \cite{Chitour2016Stability} in the case of constant delays, and it would be natural to study whether a mix of the techniques of Chapter~\ref{chap:diff-eqn-varying-delay} with those of \cite{Chitour2016Stability} could provide information on the behavior of this more general version of \eqref{eq:diff-eqn-time} with time-varying matrices. For simplicity, the case $N = 1$ could be considered as a starting point for this analysis.

\section{Difference equations with state-dependent delays}
\label{sec:open-diff-eqn-state-dependent}

As an application of the results from Section~\ref{sec:stability}, we have considered, in Section~\ref{sec:appli-varying-delay.1}, the stability analysis of the difference equation \eqref{equation-SD-1} with a state-dependent delay. As discussed in that section, we have not addressed questions on the existence and uniqueness of solutions of \eqref{equation-SD-1}, nor on sharper stability results. Interesting open problems in that direction involve the search for general conditions ensuring existence and uniqueness of solutions of \eqref{equation-SD-1} in different functional spaces, such as those of continuous functions, regulated functions, or $L^p$ functions, in an analysis similar to that carried out in Section~\ref{sec:well-posed-varying-delay} for the difference equation \eqref{eq:diff-eqn-N1-time} with a time-dependent delay. The dependence of $\tau$ on the history function $x_t$ should play an important role in such an analysis, in order to understand the behavior of the time-dependent delay $t \mapsto \tau(t, x_t)$ along a given trajectory. In particular, in view of the many delicate technical points that had to be addressed in Section~\ref{sec:well-posed-varying-delay} in order to consider solutions in $L^p$ for $p \in [1, +\infty)$, providing sufficient conditions that are not far from being sharp ensuring the existence and uniqueness of $L^p$ solutions for \eqref{equation-SD-1} for $p \in [1, +\infty)$ seems a particularly challenging question.

System~\eqref{equation-SD-1} is still the simplest form of difference equation with a state-dependent delay, and further research on the topic should focus on equations with more that just one delay, such as an equation of the form
\begin{equation}
\label{eq:diff-eqn-state}
x(t) = \sum_{j = 1}^N A_j x(t - \tau_j(t, x_t)),
\end{equation}
which is a generalization of \eqref{eq:diff-eqn-time} to state-dependent delays $\tau_1, \dotsc, \tau_N$. While results on existence and uniqueness of solutions for \eqref{equation-SD-1} would likely not be difficult to generalize to \eqref{eq:diff-eqn-state}, the analysis of the asymptotic behavior of solutions of \eqref{eq:diff-eqn-state} would certainly benefit from a preliminary study of the open problems described in Section~\ref{sec:open-diff-eqn-many-delays}.

\section{Stability of nonlinear difference equations with time-va\texorpdfstring{\-}{}rying or state-de\texorpdfstring{\-}{}pen\texorpdfstring{\-}{}dent delays}
\label{sec:open-diff-eqn-nonlinear}

The open problems described in Sections~\ref{sec:open-diff-eqn-many-delays} and \ref{sec:open-diff-eqn-state-dependent} all concern linear\footnote{Strictly speaking, difference equations with state-dependent delays are intrinsically nonlinear due to the dependence of the delay on the state: an expression of the form $x(t - \tau(t, x_t))$ is nonlinear in $x$, except in trivial cases. We make here the slight abuse of terminology of calling a difference equation with state-dependent delays \emph{linear} if the equation obtained by replacing the state-dependent delays by constants is linear.} difference equations. Generalizing the results of Chapter~\ref{chap:diff-eqn-varying-delay} to nonlinear equations is an open problem.

A first possible generalization would be to consider a nonlinear difference equation with a time-varying delay of the form
\begin{equation}
\label{eq:diff-eqn-nonlinear}
x(t) = f(t, x(t - \tau(t))),
\end{equation}
for some given function $f\colon \mathbb R_+ \times \mathbb R^d \to \mathbb R^d$. The study of such a class of equations can be motivated by its links to some hyperbolic PDEs: it was shown in \cite{Chitour2021One} that the study of a wave equation with a nonlinear (and even set-valued) boundary condition can be reduced to the study of a nonlinear difference equation of the form $x(t) = f(x(t - \tau))$ (or a set-valued generalization of it), for some constant delay $\tau > 0$. The fact that $\tau$ is constant comes from the fact that the analysis carried out in \cite{Chitour2021One} considers a wave equation with a constant velocity: a wave equation with a time-dependent velocity and time-dependent boundary conditions would lead, by the same arguments as in \cite{Chitour2021One}, to a nonlinear equation of the form \eqref{eq:diff-eqn-nonlinear}, and hence the analysis of the latter would carry out to the former. Generalizations of \eqref{eq:diff-eqn-nonlinear} involving more delays or state-dependent delays, in addition to their intrinsic interest, would also have implications for the analysis of more general systems of hyperbolic PDEs, including quasilinear systems.

\section{More general systems with varying delays}

Sections~\ref{sec:open-diff-eqn-many-delays}, \ref{sec:open-diff-eqn-state-dependent}, and \ref{sec:open-diff-eqn-nonlinear} describe some open problems for difference equations with varying delays. One may also wonder to which extent the techniques used in Chapter~\ref{chap:diff-eqn-varying-delay} can be adapted to treat the case of more general equations, such as linear neutral functional differential equations with time- or state-dependent delays, such as
\begin{equation}
\label{eq:neutral}
\frac{\df}{\df t} \left[x(t) - \sum_{j=1}^N A_j(t) x(t - \tau_j(t))\right] = \sum_{j=0}^M B_j(t) x(t - r_j(t)).
\end{equation}
Note that this can be seen as a generalization of \eqref{eq:diff-eqn-time}: if $B_j(\cdot) = 0$ for every $j \in \{0, \dotsc, M\}$, then solutions of \eqref{eq:diff-eqn-time} are particular solutions of \eqref{eq:neutral}. More generally, it has been shown in \cite{Henry1974Linear} by spectral methods that, if all matrices and delays are constant in time, then stability of \eqref{eq:diff-eqn-time} is a necessary condition for the stability of \eqref{eq:neutral}. A natural open question is whether this is still true with varying matrices and delays. The same questions can also be asked for the variant of \eqref{eq:neutral} in which the delays are state-dependent.

It should also be noted that similar systems also appear in the control literature, for instance in \cite{Bonnet2020L2}, which studies stability of a control system similar to \eqref{eq:neutral} but with constant matrices. The results of \cite{Bonnet2020L2} rely on a perturbation analysis with respect to the case of constant delays, and requires the assumption that $\sum_{j=1}^N \abs{A_j} < 1$ for some matrix norm $\abs{\cdot}$. An important open question, related to the similar question from Section~\ref{sec:open-diff-eqn-many-delays}, is whether such an assumption, which may seem quite restrictive at first, can be relaxed.

\section{The MID property for single-delay systems in higher dimensions}

One of the main drawbacks of the theory presented in Chapter~\ref{chap:MID} is that it only applies to delay-differential equations of the form \eqref{eq:DDE}, which are scalar equations with a single delay. As a first step towards generalizing those results, one might consider systems with a single delay, of the form
\begin{equation}
\label{eq:linear-one-delay}
x^\prime(t) = A_0 x(t) + A_1 x(t - \tau),
\end{equation}
where $A_0$ and $A_1$ are matrices in $\mathcal M_d(\mathbb R)$, $\tau > 0$, and the state $x(t)$ at time $t$ takes values in $\mathbb R^d$. This is a time-delay system of retarded type, and its characteristic function is the function $\Delta\colon \mathbb C \to \mathbb C$ defined for $s \in \mathbb C$ by
\[
\Delta(s) = \det(s \Id - A_0 - A_1 e^{-s \tau}).
\]
Note that, by expanding the above determinant, $\Delta$ can be rewritten as
\begin{equation}
\label{eq:Delta-higher-dimension}
\Delta(s) = \sum_{k = 0}^{d} P_{d - k}(s) e^{- k s \tau},
\end{equation}
where $P_d$ is a monic polynomial of degree $d$ and, for $k \in \llbracket 0, d-1 \rrbracket$, $P_k$ is a polynomial of degree at most $k$. The main difficulty with respect to the analysis carried out in Chapter~\ref{chap:MID} for the characteristic equation \eqref{eq:Delta-DDE} is that the quasipolynomial $\Delta$ from \eqref{eq:Delta-higher-dimension}, in addition to $e^{-s \tau}$, contains also terms corresponding to integer multiples of the delay $\tau$ up to $d$, $e^{-2 s \tau}, \dotsc, e^{-d s \tau}$. It is not clear how the results from Chapter~\ref{chap:MID} can be generalized to this setting, since the factorization-based approach taken in the proof of Theorem~\ref{thm:main-MID} and the links exploited in that proof with Padé approximations of the exponential function seem quite specific to the case of characteristic equations of the form \eqref{eq:Delta-DDE}. Other proofs of Theorem~\ref{thm:main-MID}, however, are also possible, based mostly on more direct and longer computations, such as the proofs previously provided in \cite{Mazanti2021Multiplicity, Boussaada2022Generic}, and they might be a better starting point towards a more general result for \eqref{eq:Delta-higher-dimension}.

Another difficulty in generalizing the results of Chapter~\ref{chap:MID} to \eqref{eq:linear-one-delay} and \eqref{eq:Delta-higher-dimension} is that, contrarily to what happens for \eqref{eq:DDE} and \eqref{eq:Delta-DDE}, the relations between the coefficients of the matrices $A_0$ and $A_1$ from \eqref{eq:linear-one-delay} and of the polynomials $P_0, \dotsc, P_d$ from \eqref{eq:Delta-higher-dimension} are more intricate, as the coefficients of $P_0, \dotsc, P_d$ will be sums of products of the coefficients of $A_0$ and $A_1$. Hence, even if one provides conditions on the coefficients of $P_0, \dotsc, P_d$ ensuring the existence of a dominant root of large enough or maximal multiplicity, the translation of those conditions in terms of the coefficients of $A_0$ and $A_1$ will involve the resolution of a nonlinear system of equations.

Results in these directions were recently obtained in \cite{Boussaada2025Toward} for systems in dimension $d = 2$. The main result of that reference makes use of Cauchy's argument principle and the counting formula from \cite{Hassard1997Counting, Stepan1989Retarded}, which require several explicit computations that are hard to generalize to higher dimension.

\section{The MID property for systems with several delays}

While the main result from Section~\ref{sec:MID-two-delays} shows that the MID property holds for the first-order system with two delays \eqref{eq:main-syst-two-delays}, a more general analysis of the multi-delay setting is still absent from the literature. Section~\ref{sec:three-delays} illustrates that Theorem~\ref{thm:main-two-delays} cannot be generalized to the first-order system with three delays \eqref{eq:main-syst-three-delays}. The main reason for that is that, as illustrated in \eqref{eq:three-delays-limit}, after imposing the existence of a root of maximal multiplicity, when all three delays converge to the same value, the limiting quasipolynomial will also have a root of maximal multiplicity, but degree of the polynomial corresponding to the delayed term is strictly greater than that of the polynomial corresponding to the nondelayed term, and thus we are in Case~\ref{item:m-greater-n} of Theorem~\ref{thm:main-MID}, in which the root of maximal multiplicity is not dominant.

Let us now consider the more general situation of a quasipolynomial of the form
\begin{equation}
\label{eq:Delta-many-delays}
\Delta(s) = P_0(s) + \sum_{k = 1}^{\ell} P_{k}(s) e^{- s \tau_k},
\end{equation}
where $\tau_1, \dotsc, \tau_\ell$ are positive delays and $P_0, \dotsc, P_\ell$ are polynomials of degrees $d_0, \dotsc, d_\ell$, respectively. Assume that the coefficients of $P_0, \dotsc, P_\ell$ are chosen in order to ensure the existence of a root $s_0 \in \mathbb R$ of maximal multiplicity $D = \ell + \sum_{k = 0}^\ell d_k$. Following the computations from Section~\ref{sec:MID-two-delays}, we expect that, as all delays $\tau_1, \dotsc, \tau_\ell$ converge to the same value $\tau_\ast > 0$, $\Delta(s)$ will converge (uniformly with respect to $s$ on compact subsets of $\mathbb C$) to a quasipolynomial of the form
\begin{equation}
\label{eq:Delta-ast-conjecture}
\Delta_\ast(s) = Q_0(s) + e^{-s \tau_\ast} Q_1(s),
\end{equation}
where $Q_0$ and $Q_1$ are given polynomials. Given the computations from Section~\ref{sec:MID-two-delays}, we conjecture that $\deg Q_0 = d_0$ and $\deg Q_1 = \ell - 1 + \sum_{k = 1}^{\ell} d_k$, and thus a necessary condition for the MID property to hold for \eqref{eq:Delta-many-delays} for any possible choice of delays $\tau_1, \dotsc, \tau_\ell$ would be that
\begin{equation}
\label{eq:conjecture-ineq}
d_0 \geq \ell - 1 + \sum_{k = 1}^{\ell} d_k.
\end{equation}
More precisely, we conjecture that the following result is true.

\begin{conjecture}
\label{conj:MID-many-delays}
Let $\ell \in \mathbb N^\ast$, $P_0, \dotsc, P_\ell$ be polynomials of degrees $d_0, \dotsc, d_\ell$, respectively, and $\tau_1, \dotsc, \tau_\ell$ be positive real numbers. Consider the function $\Delta \colon \mathbb C \to \mathbb C$ defined by \eqref{eq:Delta-many-delays} and assume that $\Delta$ admits a root $s_0 \in \mathbb R$ of multiplicity $D = \ell + \sum_{k = 0}^\ell d_k$ and that \eqref{eq:conjecture-ineq} is satisfied. Then $s_0$ is a dominant root of $\Delta$, i.e., $\Real s \leq s_0$ for every $s \in \mathbb C$ such that $\Delta(s) = 0$.
\end{conjecture}

A possible strategy to prove Conjecture~\ref{conj:MID-many-delays} would be to generalize the arguments from Section~\ref{sec:MID-two-delays}. This seems to be extremely technical, in particular because, instead of reducing the problem to the analysis of a quasipolynomial depending on a single parameter as in \eqref{eq:def-Q-lambda}, one would reduce the problem to the analysis of a quasipolynomial depending on $\ell - 1$ parameters $\lambda_1, \dotsc, \lambda_{\ell-1}$, all taking values in $(0, 1)$. The space of parameters is then $(0, 1)^{\ell - 1}$, and in order to prove the MID property for all parameters in this hypercube, one would first need to understand the limiting behavior of the quasipolynomial when the parameters converge to the boundary of this hypercube. Condition~\ref{eq:conjecture-ineq} should allow one to prove that the MID property holds for the limit quasipolynomial corresponding to the vertex $v_1 = (1, \dotsc, 1)$, and, at the light of the analysis carried out in Section~\ref{sec:MID-two-delays}, one expects that the limit quasipolynomial corresponding to the other vertices of the hypercube will be of the form \eqref{eq:Delta-ast-conjecture} with $\deg Q_0 \geq \deg Q_1$, allowing one to obtain the MID property from Theorem~\ref{thm:main-MID}. Once the MID property is established for the limiting polynomials corresponding to the vertices of the hypercube $(0, 1)^{\ell - 1}$, one may then proceed to prove that this property is preserved in each edge of this hypercube, by proceeding as in Section~\ref{sec:MID-two-delays}, and then proceed inductively to show that the MID property is preserved also in each $2$-dimensional face, $3$-dimensional face, etc., until obtaining it for the full hypercube. There are many technical challenges in order to make these ideas into a full proof of Conjecture~\ref{conj:MID-many-delays}.

Even if Conjecture~\ref{conj:MID-many-delays} is proved to be true, there still remain other important open questions for the MID property in the multi-delay case. Indeed, in the cases where \eqref{eq:conjecture-ineq} is not satisfied, such as for the first-order equation with three delays \eqref{eq:main-syst-three-delays} from Section~\ref{sec:three-delays}, a natural question is whether the MID property still holds true for some values of the delays $\tau_1, \dotsc, \tau_\ell$. This question seems more difficult to address, and it is not clear whether the techniques from Section~\ref{sec:MID-two-delays} can be of help. A good starting point would be to consider in more details the first-order equation with three delays from Section~\ref{sec:three-delays}: for which values of the delays $\tau_1, \tau_2, \tau_3$ is the root $0$ of multiplicity $4$ of the quasipolynomial $Q(\cdot, \tau_1, \tau_2, \tau_3)$ from \eqref{eq:def-Q-three-delays} a dominant root? A more detailed analysis of this case might be of help in order to treat the general case.

\section{The MID property for pairs of complex-conjugate roots}

Theorem~\ref{thm:complex-dominance} is the only known result on the MID property for pairs of complex-conjugate roots, and it only concerns the DDE \eqref{eq:DDE} in the case $n = 2$ and $m = 1$. It would be interesting to investigate whether Theorem~\ref{thm:complex-dominance} can be generalized to other values of $n$ and $m$.

Note that, for $n = 2$ and $m = 1$, the degree of the quasipolynomial \eqref{eq:Delta-DDE} is $4$, and Theorem~\ref{thm:complex-dominance} proceeds by placing a pair of complex-conjugate roots of multiplicity $2$ each, so that the total multiplicity of the chosen roots coincides with the degree of $\Delta$. It is not clear how Theorem~\ref{thm:complex-dominance} should be generalized to the cases where the degree of $\Delta$ is odd; a possible strategy would be to place a pair of complex-conjugate roots of multiplicity $\frac{D-1}{2}$ each, and then a simple real root equal to the real part of the chosen complex-conjugate pair, but this is not the only possible choice, and one may wonder which choice is more suitable in order to obtain the MID property.

In the cases $(n, m) \in \{(1, 0), (3, 2), (4, 3)\}$, numerical computations seem to suggest that Theorem~\ref{thm:complex-dominance} can be generalized to the case $(n, m) = (4, 3)$, but not to $(n, m) = (1, 0)$ or $(n, m) = (3, 2)$. Indeed, in the latter cases, whose degrees are $2$ and $6$, respectively, when placing a pair of complex-conjugate roots of multiplicity $1$ each or $3$ each, respectively, the corresponding quasipolynomial coefficients ensuring such multiplicities and generalizing \eqref{eq:coeffs-complex-dominance} have singularities for some values of the imaginary part $\theta_0$ of the chosen pair of complex-conjugate roots. Numerical computations suggest that the chosen pair of complex-conjugate roots is dominant for $\theta_0 \in (-\theta_\ast, \theta_\ast)$, where $\theta_\ast > 0$ is the first singularity of the coefficients (with $\theta_\ast = \frac{\pi}{\tau}$ in the case $(n, m) = (1, 0)$), but such dominance is lost for $\theta_0 > \theta_\ast$, with the apparition of a positive real root coming from complex infinity. More positive real roots appear when $\theta_0$ increases and crosses the other singularities of the coefficients, and the total number of roots with positive real part for some choice of $\theta_0 > 0$ seems to coincide with the number of singularities of the coefficients in the interval $[0, \theta_0]$. On the other hand, the quasipolynomial coefficients obtained in the case $(n, m) = (4, 3)$ have no singularities in $\theta_0$, and numerical computations suggest that the chosen pair of complex-conjugate roots is always dominant. We thus conjecture that Theorem~\ref{thm:complex-dominance} can be generalized to \eqref{eq:Delta-DDE} when $m = n - 1$ if $n$ is even, but that, when $n$ is odd, dominance will only be true if the imaginary part of the chosen root is smaller than the first singularity.

\section{Sufficient conditions for dominance of roots with large multiplicity}

Section~\ref{sec:MID-intermediate} provides several insights on the behavior of the spectrum of DDEs of the form \eqref{eq:DDE} when the characteristic function \eqref{eq:Delta-DDE} of \eqref{eq:DDE} admits a root of overorder multiplicity, i.e., of multiplicity strictly greater than the highest order of derivation $n$ of the nondelayed terms of \eqref{eq:DDE}. Even though these results can be used to prove dominance of such a multiple root in some cases, as illustrated in Section~\ref{sec:expl-integrator}, obtaining sufficient conditions for dominance of overorder roots in a general setting is still an important open problem.

For low-order systems or in particular configurations, dominance of an overorder root for \eqref{eq:DDE} was considered in \cite{Balogh2020Towards} in the case $m = n - 1$ for roots of multiplicity at least $n$ under the assumption that all roots of the polynomial $s^n + \sum_{k = 0}^{n-1} a_k s^k$ are real; in \cite{RojasRicca2022Multiplicity, RojasRicca2024Dominant} in the case $m = n - 1$ for the $n$th integrator (i.e., $a_0 = \dotsb = a_{n - 1} = 0$ in \eqref{eq:DDE}) for roots of multiplicity at least $n$; in \cite{Balogh2022Conditions, Sassi2022Observer} in the case $m = n - 1$ for roots of multiplicity at least $n$ or $n+1$, providing rather conservative sufficient conditions for dominance; in \cite{Boussaada2021Insights, Ma2020Multiplicity, Ma2022PID, Benarab2023MultiplicityI} in the case $n = m = 2$ for roots of multiplicities at least $3$, $4$ or $5$, with a more complete characterization of the cases of dominance; and in \cite{Benarab2023MultiplicityII} in the case $n = m = 1$ for roots of multiplicities at least $2$ or $3$. Sufficient conditions for dominance were obtained for the general case of any $n$ and $m$ with $n \geq m$ in \cite{Boussaada2025OverOrder, Boussaada2025Decay, Boussaada2022MID} based on the Green--Hille transform (see Section~\ref{sec:correction}), but these conditions are expected to be far from sharp, they are not easy to be verified, and the only known cases satisfying these conditions are trivial ones (for instance, when the multiplicity is maximal). We also mention that a methodology to establish proofs of dominance on a case-by-case basis was described in \cite{Benarab2023MultiplicityI}, based on a priori bounds on the imaginary parts of potential roots contradicting dominance, together with a suitable polynomial approximation of the exponential term in the characteristic quasipolynomial. However, the current techniques used in the analysis of the MID property seem to be insufficient to obtain sharper or easier-to-verify dominance conditions for overorder roots, and new ideas would probably be needed in order to obtain significant advances in these directions.

\section{More general partial pole placement methods guaranteeing dominance}

As described in Section~\ref{sec:stabilization-LTI-TDS}, \emph{partial pole placement} refers to techniques of stabilization of time-delay systems of the form \eqref{eq:neutral-control} through feedback laws of the form \eqref{eq:feedback}, consisting in choosing the parameters of the feedback law in order to impose some spectral values of the closed-loop system. The results from Chapter~\ref{chap:MID} can be interpreted at the light of partial pole placement: for instance, the equivalence between \ref{item:max-mult} and \ref{item:Coeffs} in Theorem~\ref{thm:main-MID} shows how one should choose the coefficients of \eqref{eq:DDE} in order to have a spectral value at a given $s_0 \in \mathbb R$ of multiplicity $D$. The main advantage of Theorem~\ref{thm:main-MID} is the additional information on the location of the other roots of the characteristic function $\Delta$ from \eqref{eq:Delta-DDE} given in \ref{item:m-less-n}--\ref{item:m-greater-n}, which gives dominance of $s_0$ in the case $m \leq n$.

An important limitation of the partial pole placement methods described in Chapter~\ref{chap:MID} is that they are all based on placing a root of ``large enough'' multiplicity. This may impose too many constraints on the coefficients of the system, such as in the case of maximal multiplicity from Theorem~\ref{thm:main-MID}, in which all coefficients are uniquely determined in terms of the desired root $s_0$ and the delay $\tau$, which reduces the applicability scope of the result. Roots with large multiplicity are often also not desirable in applications, as they are too sensitive with respect to variations in parameters: a variation of order $\varepsilon$ of the parameters of a system admitting a root of multiplicity $M$ may lead to a perturbation of the root of order up to $\varepsilon^{1 / M}$ \cite{Rellich1969Perturbation, Najman1986Remarks, Langer1989Remarks, Langer1992Remarks, Hryniv1999Perturbation, Michiels2017Explicit}. In addition, apart from Section~\ref{sec:MID-complex}, all the results in Chapter~\ref{chap:MID} concern the placement of real roots, while, in practical applications, complex conjugate roots are often preferred for performance reasons (see, e.g., \cite[Chapter~5]{Dorf2016Modern} for linear finite-dimensional systems without delays).

These limitations motivate the search for other partial pole placement techniques that, similarly to the MID property from Chapter~\ref{chap:MID}, ensure that, when placing a certain number of roots of a time-delay system through the choice of some coefficients, the rightmost placed root is dominant, but avoiding the limitations of MID. In order to avoid placing a multiple root and to reduce the amount of constraints on the coefficients, an alternative method investigated in recent years in the literature is to place several simple real roots (often assumed to be equidistant), based on the CRRID property described above in Section~\ref{sec:intro-MID} and developed in \cite{Bedouhene2019Pole, Amrane2020Some, Schmoderer2023Quasipolynomials, Amrane2018Qualitative, Bedouhene2020Real, Schmoderer2024Insights}. However, a general result similar to Theorem~\ref{thm:main-MID} is not yet available for the CRRID property and constitutes an important research line. Another important open question is to investigate whether one can develop techniques to place pairs of complex conjugate roots while guaranteeing dominance of the placed roots, in more general situations than those from Section~\ref{sec:MID-complex} and avoiding complex conjugate roots with large multiplicity.

\bibliographystyle{BibStyle}
\bibliography{Bibliography}

\includepdf[pages=2]{Cover/Cover.pdf}

\end{document}